\documentclass[reqno]{amsbook}

\newif\ifpdf
\ifx\pdfoutput\undefined
  \pdffalse
\else
  \pdfoutput=1
  \pdftrue
\fi

\usepackage{amssymb}
\usepackage{amscd}
\usepackage{verbatim}

\ifpdf
  \usepackage[pdftex]{graphicx}
  \usepackage[pdftex]{hyperref}
\else
  \usepackage{graphicx}
\fi
\numberwithin{section}{chapter}
\numberwithin{equation}{chapter}
\numberwithin{table}{chapter}
\numberwithin{figure}{chapter}

\theoremstyle{plain}
\newtheorem{theorem}{Theorem}[chapter]
\newtheorem{proposition}[theorem]{Proposition}
\newtheorem{lemma}[theorem]{Lemma}
\newtheorem{corollary}[theorem]{Corollary}
\newtheorem{Prop-def}[theorem]{Proposition-Definition}

\theoremstyle{definition}
\newtheorem{definition}[theorem]{Definition}
\newtheorem{example}[theorem]{Example}

\theoremstyle{remark}
\newtheorem{remark}[theorem]{Remark}
\newtheorem*{remark*}{Remark}
\newtheorem*{acknowledge}{Acknowledgment}

\newcommand{\clearemptydoublepage}{\newpage{\pagestyle{empty}\cleardoublepage}}

\newcommand{\C}{{\mathbf{C}}}
\newcommand{\N}{{\mathbf{N}}}
\renewcommand{\P}{{\mathbf{P}}}
\newcommand{\Q}{{\mathbf{Q}}}
\newcommand{\R}{{\mathbf{R}}}

\newcommand{\Z}{{\mathbf{Z}}}

\newcommand{\Nbar}{\overline{\N}}
\newcommand{\Rbar}{\overline{\R}_+}
\newcommand{\btilde}{{\tilde{\b}}}
\newcommand{\bbar}{{\overline{\b}}}

\newcommand{\bE}{{\bar{E}}}

\newcommand{\bk}{{\bar{k}}}
\newcommand{\bn}{{\bar{n}}}
\newcommand{\bp}{{\bar{p}}}

\newcommand{\bs}{{\bar{s}}}
\newcommand{\btheta}{{\bar{\theta}}}

\newcommand{\bcT}{{\bar{\mathcal{T}}}}

\newcommand{\cA}{{\mathcal{A}}}
\newcommand{\cB}{{\mathcal{B}}}
\newcommand{\cC}{{\mathcal{C}}}
\newcommand{\cD}{{\mathcal{D}}}
\newcommand{\cE}{{\mathcal{E}}}
\newcommand{\cF}{{\mathcal{F}}}

\newcommand{\cI}{{\mathcal{I}}}
\newcommand{\cM}{{\mathcal{M}}}
\newcommand{\cN}{{\mathcal{N}}}
\newcommand{\cO}{{\mathcal{O}}}
\newcommand{\cP}{{\mathcal{P}}}
\newcommand{\cS}{{\mathcal{S}}}
\newcommand{\cT}{{\mathcal{T}}}
\newcommand{\cV}{{\mathcal{V}}}

\newcommand{\hcV}{{\widehat{\mathcal{V}}}_x}

\newcommand{\hk}{{\hat{k}}}

\newcommand{\halpha}{{\hat{\alpha}}}
\newcommand{\hbeta}{{\hat{\beta}}}
\newcommand{\hmu}{{\hat{\mu}}}
\newcommand{\hnu}{{\hat{\nu}}}
\newcommand{\hphi}{{\hat{\phi}}}
\newcommand{\hpsi}{{\hat{\psi}}}

\newcommand{\tn}{{\tilde{n}}}

\newcommand{\ti}{{\tilde{\imath}}}

\newcommand{\tpi}{{\tilde{\pi}}}

\newcommand{\tpsi}{{\tilde{\psi}}}
\newcommand{\tE}{{\tilde{E}}}

\newcommand{\tcV}{{\tilde{\mathcal{V}}}}

\newcommand{\wi}{\widehat{\imath}}
\newcommand{\wtV}{\widetilde{\mathcal{W}}}

\newcommand{\fm}{{\mathfrak{m}}}
\newcommand{\fp}{{\mathfrak{p}}}
\newcommand{\fq}{{\mathfrak{q}}}
\newcommand{\fB}{{\mathfrak{B}}}

\newcommand{\fX}{{\mathfrak{X}}}

\newcommand{\bbA}{{\mathbb{A}}}
\newcommand{\bbH}{{\mathbb{H}}}
\newcommand{\bbP}{{\mathbb{P}}}

\newcommand{\bI}{{\overline{I}}}

\renewcommand{\=}{:=}
\renewcommand{\a}{\alpha}
\renewcommand{\b}{\beta}
\newcommand{\e}{\varepsilon}

\renewcommand{\l}{\lambda}
\newcommand{\half}{{\frac12}}
\newcommand{\s}{\sigma}

\newcommand{\lcm}{\operatorname{lcm}}
\renewcommand{\div}{\operatorname{div}}

\newcommand{\Gal}{\operatorname{Gal}}

\newcommand{\lgt}{\operatorname{length}}
\newcommand{\val}{\operatorname {val}}

\newcommand{\id}{\operatorname{id}}

\newcommand{\vol}{\operatorname{Vol}}

\newcommand{\con}{\operatorname{con}}
\newcommand{\gr}{\operatorname{gr}}
\newcommand{\rk}{\operatorname{rk}}
\newcommand{\ratrk}{\operatorname{rat.rk}}
\newcommand{\trdeg}{\operatorname{tr.deg}}
\newcommand{\krull}{\operatorname{krull}}
\newcommand{\supp}{\operatorname{supp}}
\newcommand{\mass}{\operatorname{mass}}
\newcommand{\Spec}{\operatorname{Spec}}
\newcommand{\diam}{\operatorname{diam}}

\newcommand{\self}{\circlearrowleft}

\newcommand{\ie}{i.e.\ }

\newcommand{\vv}{{\vec{v}}}
\newcommand{\ww}{{\vec{w}}}
\newcommand{\eg}{e.g.\ }
\newcommand{\qand}{{\quad\text{and}\quad}}
\newcommand{\qor}{{\quad\text{or}\quad}}

\newcommand{\cVqm}{{\mathcal{V}_\mathrm{qm}}}
\newcommand{\cVqmx}{{\mathcal{V}_{\mathrm{qm},x}}}
\newcommand{\cVdiv}{{\mathcal{V}_\mathrm{div}}}
\newcommand{\cVdivx}{{\mathcal{V}_{\mathrm{div},x}}}

\newcommand{\hcVfin}{{\widehat{\mathcal{V}}}_\mathrm{fin}}
\newcommand{\Gast}{\Gamma^\ast}

\makeindex
\begin{document}
%
%
%
%
%
%

\pagenumbering{roman}
\thispagestyle{empty}
\vspace*{2cm}
\begin{center}
  \Huge\textbf{THE VALUATIVE TREE}\\
  \vspace{2cm}
  \Large CHARLES FAVRE\\
  \vspace{.5cm}
  CNRS (Universit\'e Paris VII),
  Institut de Math\'ematiques\\
  Equipe G\'eom\'etrie et Dynamique\\
  F-75251 Paris Cedex 05, France\\
  \texttt{favre@math.jussieu.fr}\\
  \vspace{1.5cm}
  \Large MATTIAS JONSSON\footnote{Supported by NSF Grant No DMS-0200614.}\\
  \vspace{.5cm}
  Department of Mathematics,
  University of Michigan\\
  Ann Arbor, MI 48109-1109, 
  USA\\
  \texttt{mattiasj@umich.edu}\\
  \vspace{2cm}
  \today\\
\end{center}

\clearemptydoublepage
\vspace*{3cm}

\begin{center}
  \Large\textbf{Abstract}
\end{center}
\vspace{.5cm}

  We describe the set $\cV$ of all $\R\cup\{\infty\}$ valued
  valuations $\nu$ on the ring $\C[[x,y]]$ normalized by
  $\min\{\nu(x),\nu(y)\}=1$. It has a natural structure of an
  $\R$-tree, induced by the order relation $\nu_1\le\nu_2$ iff
  $\nu_1(\phi)\le\nu_2(\phi)$ for all $\phi$. It can also be metrized,
  endowing it with a metric tree structure.  From the algebraic point
  of view, these structures are obtained by taking a suitable quotient
  of the Riemann-Zariski variety of $\C[[x,y]]$, in order to force it
  to be a Hausdorff topological space.  The tree structure on $\cV$
  also provides an identification of valuations with balls of
  irreducible curves in a natural ultrametric.  We show that the dual
  graphs of all sequences of blow-ups patch together, yielding an
  $\R$-tree naturally isomorphic to $\cV$. Altogether, this gives many
  different approaches to the valuative tree $\cV$.  We then describe
  a natural Laplace operator on $\cV$. It associates to (special)
  functions of $\cV$ a complex Borel measure. Using this operator, we
  show how measures on the valuative tree can be used to encode
  naturally both integrally closed ideals in R and cohomology classes
  of the voute etoilee of $(\C^2,0)$.


\vfill


%
%
%
%
%
%
\clearemptydoublepage
\begin{center}
  \Large\textbf{The valuative tree}
\end{center}
\vspace{.5cm}
\begin{figure}[hl] 
  \includegraphics[width=\textwidth]{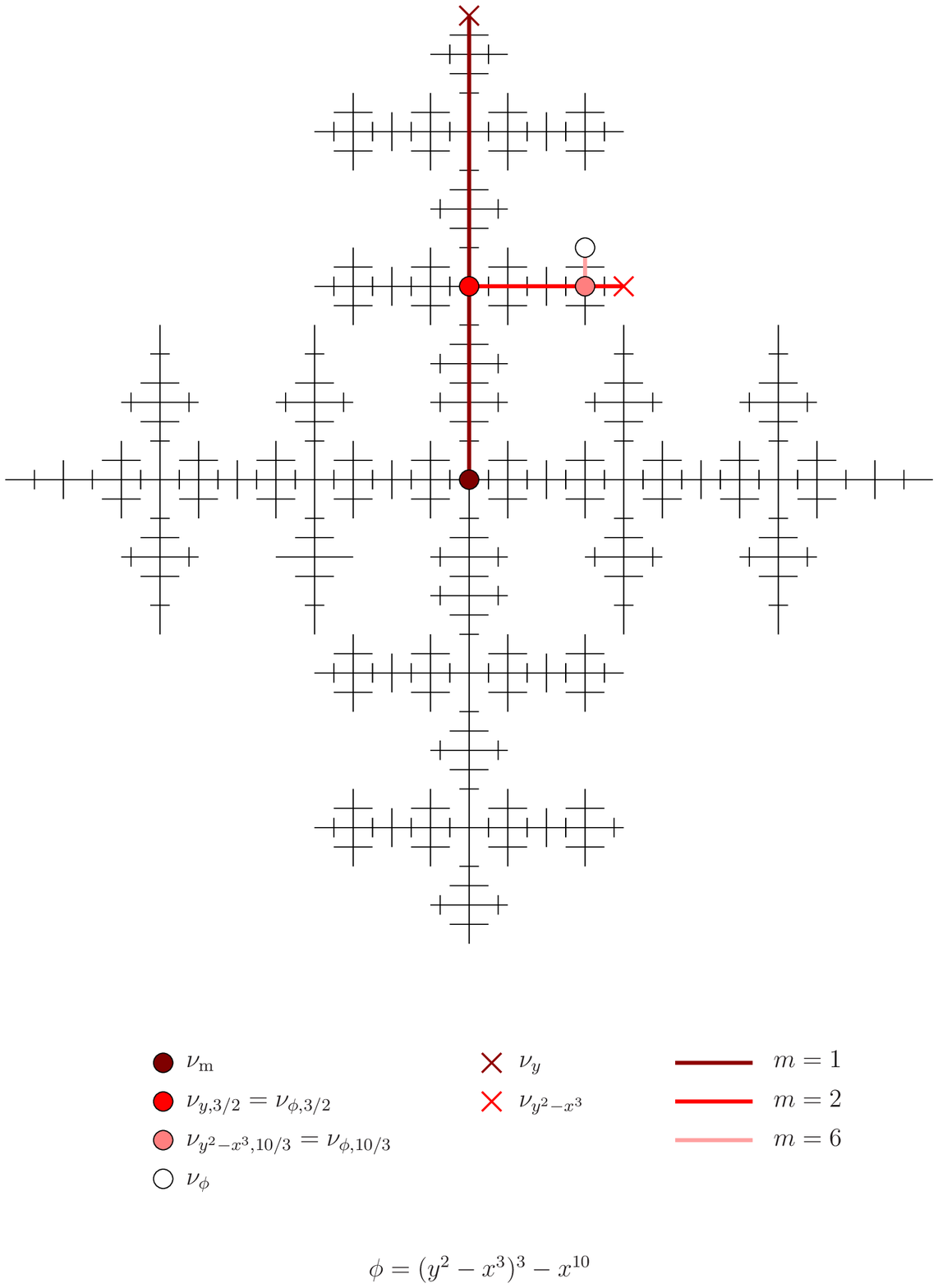}
\end{figure}
%
%
%
%
%
%
\clearemptydoublepage
\vspace*{3.5cm}
\begin{center}
  \Large\textbf{Structure of the memoir}
\end{center}
\vspace{.5cm}

In capital letters: chapter numbering; 
in Arabic numbers: section numbering.  
A plain arrow linking chapter A to chapter B indicates
that the understanding of B relies heavily on a previous lecture of A. 
A dashed arrow indicates a looser link between both chapters.
\vspace*{3cm}

\begin{figure}[hl] 
  \includegraphics[width=\textwidth]{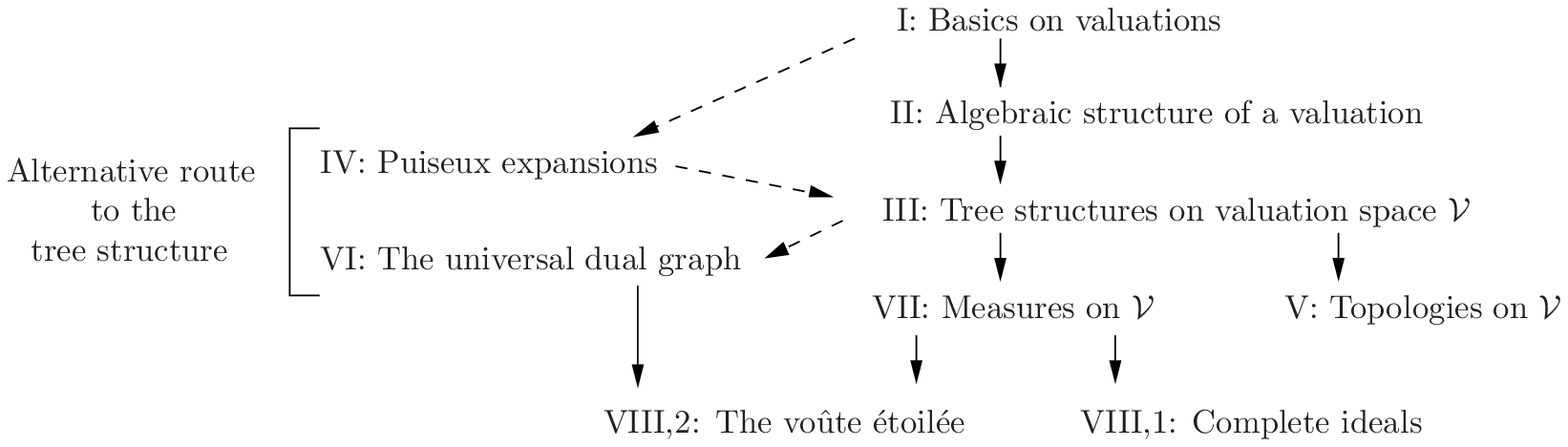}
\end{figure}
%
%
%
%
%
%
\clearemptydoublepage
\tableofcontents
%
%
%
%
%
%
\clearemptydoublepage
\pagenumbering{arabic}
\chapter*{Introduction}
The purpose of this monograph is to give a new approach to
singularities in a local, two-dimensional setting.  Our method enables
us to study curves, analytic ideals in $R= \C[[x,y]]$, and
plurisubharmonic functions in a unified way. It is also general and
powerful enough so that it can be applied to other situations, 
for instance when studying the dynamics of fixed point germs 
$f:(\C^2,0)\self$.

Curves and ideals can then be analyzed through the order of vanishing
of their pull-backs along irreducible components of the exceptional
divisor. These orders of vanishing define real-valued functions on the
ring R called divisorial valuations.\footnote{A valuation can be
evaluated at a psh function by replacing order of vanishing by generic
Lelong number.} It is a classical fact that the singularity of the
curves or ideals is completely determined by the values of 
\emph{all} divisorial valuations.


This naturally leads us to look at the set of all divisorial
valuations. Indeed our aim is to describe in detail the structure of
a slightly larger set $\cV$ that we call \emph{valuation space}. The
elements of $\cV$ are functions $\nu:R\to\Rbar=[0,\infty]$,
with $\nu(c) =0$ for all $c\in\C^*$, satisfying the standard axioms of
valuations: $\nu(\phi\psi)=\nu(\phi)+\nu(\psi)$; 
$\nu(\phi+\psi)\ge\min\{\nu(\phi),\nu(\psi)\}$ 
for all $\phi,\psi\in R$. 
We normalize them by $\nu(\fm)=\min\{\nu(x),\nu(y)\}=1$. 
The central theme of our work is that valuation space has a natural
structure of a \emph{tree} modeled on the real line, and that this
structure can be used to efficiently encode singularities of various
kinds.

\medskip
We distinguish between three different types of tree structures.
A \emph{nonmetric tree} is a poset
having a unique minimal element (its root) in which 
all sets of the form $\{\sigma\ ;\ \sigma\le\tau\}$
are isomorphic to real intervals, \ie there exists an
order-preserving bijection of each such set onto a real interval.
A \emph{parameterized tree} is a nonmetric tree in which 
all these bijections are chosen in a compatible way.
Finally, almost equivalent to parameterized trees are
\emph{metric trees}; these are 
metric spaces in which any two points are joined by a unique
path isometric to a real interval.\footnote{Metric trees are often 
called $\R$-trees in the literature.}

The nonmetric tree structure on $\cV$ arises as follows.
For $\nu,\mu\in \cV$, we declare $\nu\le\mu$ when 
$\nu(\phi)\le\mu(\phi)$ for all $\phi\in R$. 
Our normalization $\nu(\fm)=1$ implies that 
the \emph{multiplicity valuation} 
$\nu_\fm$ sending $\phi$ to its multiplicity $m(\phi)$
at the origin is dominated by any other valuation. 
This natural order defines a 
nonmetric tree structure on $\cV$, rooted at $\nu_\fm$ 
(Theorem~\ref{order}). 

As for the other two tree structures, any irreducible 
(formal local) curve $C$ defines a \emph{curve valuation}
$\nu_C\in\cV$: $\nu_C(\phi)$ is the normalized intersection
number between the curves $C$ and $\{\phi=0\}$.
A curve valuation is a maximal element under $\le$ and the segment 
$[\nu_\fm,\nu_C]$ is isomorphic, as a totally
ordered set, to the interval $[1,\infty]$.
We construct an increasing function $\a:\cV\to[1,\infty]$ 
that restricts to a bijection of 
$[\nu_\fm,\nu_C]$ onto $[1,\infty]$ for any $\phi$;
as a consequence $\a$ defines a parameterization of $\cV$.
The number $\a(\nu)$ is called the \emph{skewness}
of $\nu$.\footnote{The skewness is the inverse of the volume 
  of a valuation as defined in~\cite{ELS} (see Remark~\ref{rem-vol}).}
It is defined by the formula $\a(\nu)=\sup_\phi\nu(\phi)/m(\phi)$.

In addition to the partial ordering and skewness parameterization just
described, the valuative tree also carries an important
\emph{multiplicity} function.  The multiplicity of a valuation $\nu$
is equal to the infimum of the multiplicity of all curves whose
associated curve valuations dominate $\nu$ in the partial
ordering. Thus the multiplicity function is an increasing function on
$\cV$ with values in $\Nbar=\N\cup\{\infty\}$.  A second important
parameterization of $\cV$, \emph{thinness}, can be defined in terms of
skewness and multiplicity.\footnote{The thinness is also very precisely
  related to the Jacobian ideal of the valuation (see Remark~\ref{Rjaco}).}
We shall refer loosely to the combination of the partial ordering, the
parameterizations by skewness and thinness, and the multiplicity
function as the \emph{tree structure on $\cV$}.

\medskip
There are four types of inhabitants of the valuative tree $\cV$. 
The interior points, \ie the points that are not maximal in the
partial ordering, are valuations that become monomial 
(\ie determined by their values on a pair of local coordinates) 
after a finite sequence of blowups. 
We call them \emph{quasimonomial}.
They include all divisorial valuations but also all 
\emph{irrational} valuations such
as the monomial valuation defined by $\nu(x)=1$, $\nu(y)=\sqrt{2}$.
The other points in $\cV$, \ie the ends of the valuative tree,
are curve valuations and \emph{infinitely singular} valuations, 
which can be characterized
as the valuations with infinite multiplicity.  

There is in fact a fifth type of valuations. 
These valuations cannot
be defined as real-valued functions, but 
define functions on $R$ with values in $\R_+\times \R_+$ 
(endowed with the lexicographic order).
In fact they do have a natural place in the valuative tree,
as tree tangent vectors at points corresponding to divisorial
valuations (see Theorem~\ref{div-tangent}).  Geometrically they are
curve valuations where the curve is defined by an exceptional divisor.
We shall hence call them \emph{exceptional curve valuations}.

\medskip
The valuative tree is a beautiful object which may be
viewed in a number of different ways.
Each corresponds to a particular interpretation of a 
valuation, and each gives a new insight into it.
Some of them will hopefully lead to generalizations in a broader
context. Let us describe four such points of views.

The first way consists of identifying valuations with
balls of curves.
For any two irreducible curves $C_1$, $C_2$
set $d(C_1,C_2)=m(C_1)m(C_2)/C_1\cdot C_2$ 
where $m(C_i)$ is the multiplicity of $C_i$ and 
$C_1\cdot C_2$ is the intersection product
of $C_1$ and $C_2$. It is a nontrivial fact that 
$d$ defines an ultrametric
on the set $\cC$ of all irreducible formal curves (c.f.~\cite{GB}).  
This fact allows us to associate to $(\cC,d)$ 
a tree $\cT_\cC$ by declaring a point in $\cT_\cC$ 
to be a closed ball in $\cC$.
The tree structure on $\cT$ is given by
reverse inclusion of balls (partial ordering),
inverse radii of balls (parameterization) and
minimum multiplicity of curves in a ball (multiplicity).
Theorem~\ref{treeballs} states that the tree $\cT_\cC$
is isomorphic to the valuative tree $\cV$ 
with its ends removed
(\ie to the set $\cVqm$ of quasimonomial valuations).

A second way is through Puiseux series.  Just as irreducible curves
can be represented by Puiseux series, the elements in $\cV$ are
represented by valuations on the power series ring in one variable
with Puiseux series coefficients.  The set $\hcV$ of all such
(normalized) valuations has a natural tree structure and a suitably
defined restriction map from $\hcV$ to $\cV$ recovers the tree
structure on $\cV$.  In fact, $\cV$ is naturally the orbit space of
$\hcV$ under the action by the relevant Galois group (Theorem~\ref{T602}).
This approach can also be viewed in the context of Berkovich spaces
and Bruhat-Tits buildings. As a nonmetric tree, $\cV$ embeds as the
closure of a disk in the Berkovich projective line over the field of
Laurent series in one variable.  The metric on $\cVqm$ induced by
thinness then also arises from an identification of a subset of the
Berkovich projective line with the Bruhat-Tits building of
$\mathrm{PGL}_2$ (see Section~\ref{bruhat-tits}).

The third way is more algebraic in nature.  The earliest systematic
study of valuations in two dimensions was done in the fundamental work
by Zariski~\cite{Z0},~\cite{Z} who, among other things, identified the
set $\cV_K$ of (not necessarily real valued) valuations on $R$,
vanishing on $\C^*$ and positive on the maximal ideal $\fm=(x,y)$,
with sequences of infinitely nearby points.  The space $\cV_K$ carries
a natural topology (the Zariski topology) and is known as the
\emph{Riemann-Zariski variety}.  It is a non-Hausdorff quasi-compact
space.  The obstruction for $\cV_K$ being Hausdorff comes from the
fact that divisorial valuations do not define closed points. Namely,
their associated valuation rings strictly contain valuation rings
associated to exceptional curve valuations. One can then build a
quotient space by identifying all valuations in the closure of a
single divisorial one. This produces a compact Hausdorff
space. Theorem~\ref{weak-zariski} states that this space is
precisely $\cV$ (endowed with the topology of
pointwise convergence).

The last way uses Zariski's identification of valuations with
sequences of infinitely nearby points.  We let $\Gamma_\pi$ be the
dual graph of a finite composition of blow-ups $\pi$. It is a
simplicial tree whose set of vertices defines a poset $\Gast_\pi$.
When one sequence $\pi$ contains another $\pi'$, the poset $\Gast_\pi$
naturally contains $\Gast_{\pi'}$.  These posets therefore form an
injective system whose injective limit (or, informally, union) $\Gast$
is a poset with a natural tree structure modeled on the rational
numbers.  By filling in the irrational points and adding all the ends
to the tree we obtain a nonmetric tree $\Gamma$, the \emph{universal
dual graph}.  This nonmetric tree can in fact be equipped with with a
parameterization and multiplicity function.  These both derive from a
combinatorial procedure that to each element in $\Gast$ attaches a
vector $(a,b)\in(\N^*)^2$, the \emph{Farey weight}.\footnote{We
followed the terminology already used in~\cite{hubbard}.}  A fundamental
result (Theorem~\ref{thm-universal}) asserts that the universal dual
graph equipped with the Farey parameterization is canonically
isomorphic to the valuative tree with the thinness parameterization,
and that this isomorphism preserves multiplicity.

\medskip
As we mentioned above, singularities can be understood through
\emph{functions} on the valuative tree. It is a remarkable fact that the
information carried by these functions can also be described in terms of
\emph{complex measures} on $\cV$.
Let us be more precise. In the case of an
ideal $I\subset R$, the function on $\cV$ is given by
$g_I(\nu)\=\nu(I)$, and the measure $\rho _I$ is a positive
atomic measure supported on the Rees valuations of $I$.
This decomposition into atoms of $\rho_I$ corresponds
exactly to the Zariski decomposition of $I$ into simple
complete ideals.

In~\cite{pshfnts}, we shall show that a
plurisubharmonic function $u$ also determines a
function $g_u$ on $\cV$.
\index{tree transform!of a plurisubharmonic function}
The corresponding measure $\rho_u$, which is still
positive but not necessarily atomic,
captures essential information on the singularity of $u$.
In particular, we shall show in~\cite{criterion}
that $\rho_u$ determines the multiplier ideals of
all multiples of $u$.

The identifications $g_I\leftrightarrow\rho_I$ and
$g_u\leftrightarrow\rho_u$ are particular instances of a
general correspondence between measures on $\cV$ and certain
functions on $\cVqm$.
In fact, this correspondence, being
purely tree-theoretic in nature,
is even more general, and extends the equivalence between
positive measures and suitably normalized concave functions
on the real line.
By analogy, we thus write $\rho_I=\Delta g_I$,
$\rho_u=\Delta g_u$ and speak about the
\emph{Laplace operator} $\Delta$ on the valuative tree.

There is a second instance where complex measures on $\cV$ naturally
 appear, namely when we study the sheaf cohomology of the
 \emph{vo\^ute \'etoil\'ee} $\fX$. In our setting $\fX$ can be viewed
 as the total space of the set of all blowups above the
 origin. Elements of $ H^2(\fX, \C)$ naturally define functions on
 $\cV$, and their Laplacians are atomic measures supported on
 divisorial valuations.  The cup product on cohomology has a natural
 interpretation as an inner product on measures.  This inner product
 is a bilinear extension of an inner product on the valuative tree
 itself, and ultimately derives from intersections of curves.

\medskip
As mentioned above, the assumption that $R$ be the ring of
formal power series in two complex variables is 
unnecessarily restrictive. 
While we refer to Appendix~\ref{sec-notcomplete} for a precise
discussion, we mention here that our analysis goes through
in two important cases: the ring of holomorphic germs at
the origin in $\C^2$, and the local ring at a smooth (closed)
point of an algebraic surface over an algebraically closed field.

In order to make the monograph accessible to non-experts in
singularity theory, we have tried to keep the exposition
as self-contained as possible. This is in fact one of the
main reasons for always working in the context of complex formal
power series. 

We remark that a fair amount of the structure of the valuative tree
is implicitly contained in the analysis by Spivakovsky~\cite{spiv}.
However, our approach is quite different from his;
in particular we do not use continued fractions.
A tree structure was also described in a context similar 
to ours in~\cite{abh-an}, but without any explicit reference to 
valuations. 

Finally, the main applications of the tree structure 
on the valuative tree to analysis and dynamics will 
be explored in forthcoming papers: 
see~\cite{pshfnts},~\cite{criterion}, and~\cite{eigenvaluation}.

\medskip
We end this introduction by indicating the organization of
the monograph, which is divided into eight chapters and an appendix.

In the first chapter we give basic definitions, examples and 
results on valuations. In particular we describe the relationship 
between valuations and sequences of infinitely nearby points
(in dimension two).

Chapter~\ref{part2} is technical in nature.
We encode valuations by
finite or countable pieces of data that we call sequences of
key polynomials, or \emph{SKP's}. 
This encoding is an adaptation of a method
by MacLane~\cite{mac}, the possibility of which was
indicated to us by B.~Teissier. 
An SKP (or at least a subsequence of it) corresponds to generating
polynomials and approximate roots in the language of
Spivakovsky~\cite{spiv} and Abhyankar-Moh~\cite{abh-moh},
respectively. We are thus able to classify valuations on $R$. 
This classification is well-known to specialists 
(see~\cite{ZS,spiv} for instance) 
but we feel that our concrete approach is of independent interest. 
The representation of valuations by SKP's is
the key to the tree structure on the valuative tree.

The third chapter concerns trees. 
Our main goal is to visualize the encoding by SKP's in 
an elegant and coordinate free way.
We first discuss different definitions of trees and the
relations between them.
Using SKP's we then show that valuation space
$\cV$ does carry an intricate tree structure that we 
later in the chapter analyze in detail.

As an alternative to SKP's, Chapter~\ref{sec-puis} contains
an approach to the tree structure on $\cV$ through Puiseux series.
The results can be interpreted in the language of Berkovich. 
Specifically we indicate how the valuative tree
embeds inside the Berkovich projective line 
and Bruhat-Tits building of $\mathrm{PGL}_2$ over the local 
field of Laurent series in one variable. 
In fact, most of these results are at least implicitly contained 
in~\cite{Ber} but we felt it was worthwhile to write down the details.

In Chapter~\ref{part4} we analyze and compare different topologies on 
the valuative tree. 
The definition of a valuation as a function on $R$, as well as
the tree structure on $\cV$, gives rise to three types of 
topologies: the weak, the strong and the thin topology.
In addition, we have two topologies on $\cV_K$: the Zariski
topology and the Hausdorff-Zariski (or HZ) topology. As mentioned
above, the former gives rise to the weak topology on $\cV$
through the quotient construction. 
The HZ topology is in fact equivalent to the weak tree
topology induced by a natural discrete tree structure on $\cV_K$.

In Chapter~\ref{A3} we build the universal dual graph described above,
and show how to identify it with the valuative tree.  The fact that
valuations can be simultaneously viewed algebraically as functions on
$R$ and geometrically in terms of blowups is extremely powerful and we
spend a fair amount of time detailing some of the connections and
implications. In particular, we show that the valuation space has a
natural self-similar, or fractal, structure, see Figure~\ref{F16}.

Chapter~\ref{part-potent} is concerned with the relationship 
between measures on $\cV$ and certain classes of functions on $\cVqm$.
The analysis is purely tree-theoretic
and gives a connection between (complex) measures on a parameterized
tree and functions on the (interior of the) tree 
satisfying certain regularity properties. 
Apart from being of independent interest, this analysis 
is fundamental to many applications.

In Chapter~\ref{part-appli-analysis}, we describe two instances where
these measures appear naturally.  First we reinterpret in our context
Zariski's theory of simple complete ideals as explained
in~\cite{ZS}. This gives a new point of view on the decomposition
of any integrally closed ideal as a product of simple ideals.
We then construct Hironaka's ``vo\^ute \'etoil\'ee'' $\fX$ as
the projective limit of the total spaces of all 
sequences of blow-ups above the origin. 
We use measures on $\cV$ to understand the structure of the
sheaf cohomology group $H^2(\fX,\C)$.

Finally we conclude this monograph by an appendix containing 
a few results and discussions that did not find a natural home
elsewhere in the monograph.
Specifically, we discuss infinitely singular valuations; 
analyze the tangent space at a divisorial valuation; 
give tables summarizing the classification of valuations on
$R$ from different points of view;
present a short dictionary between our terminology 
and the more standard terminology from the theory of 
singular plane curves; 
and finally discuss what the essential assumptions are on the
ring $R$.

\medskip
\begin{acknowledge}
  The first author wishes to warmly thank Bernard Teis\-sier for his
  constant support and help, and Patrick Popescu-Pampu, Mark
  Spiva\-kovsky and Michel Vaqui\'e for fruitful discussions.  The second
  author extends his gratitude to Robert Lazarsfeld and
  Jean-Fran\c{c}ois Lafont for useful comments at an early stage of
  this project.
\end{acknowledge}


%
%
%
%
%
%
\clearemptydoublepage

\vspace*{2cm}
\begin{center}
  {\bf \large Different interpretations of valuations}
\end{center}

\begin{figure}[h] 
  \includegraphics[width=\textwidth]{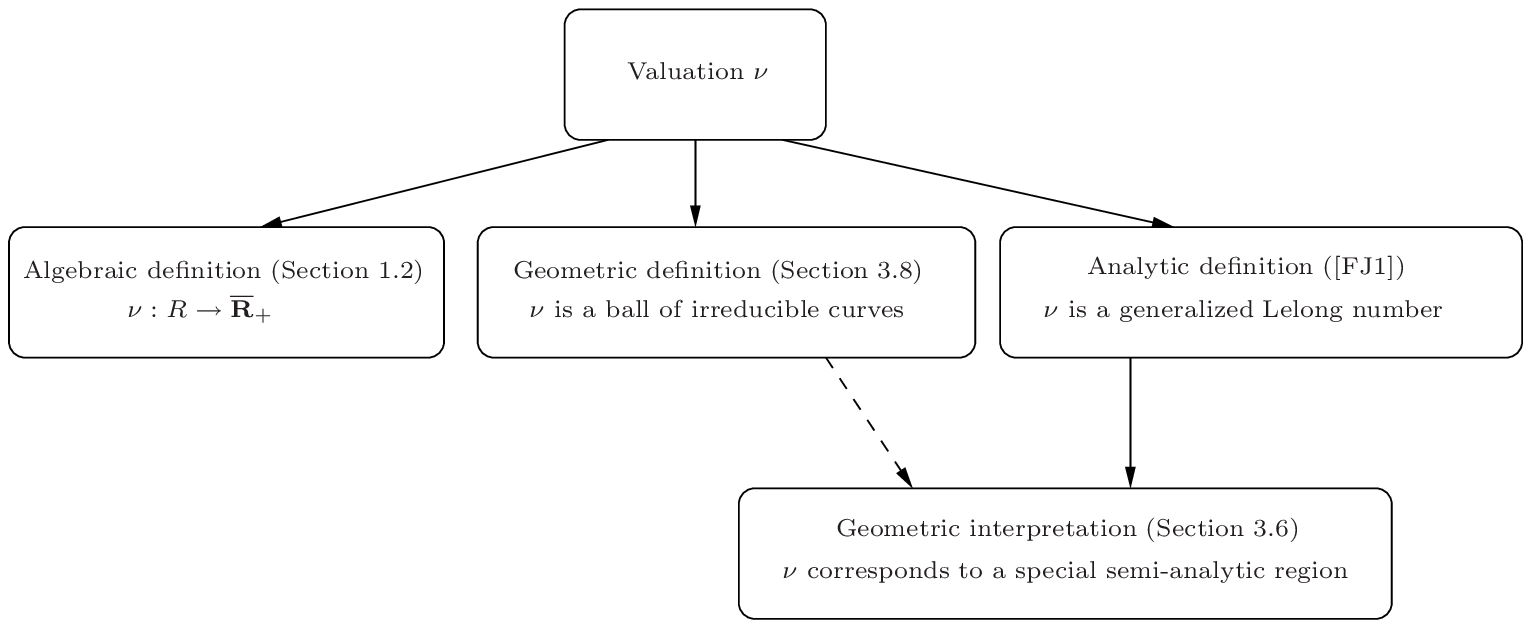}
\end{figure}

\vspace*{2cm}
\begin{center}
  {\bf \large Different approaches to the valuative tree}
\end{center}

\begin{figure}[h] 
  \includegraphics[width=\textwidth]{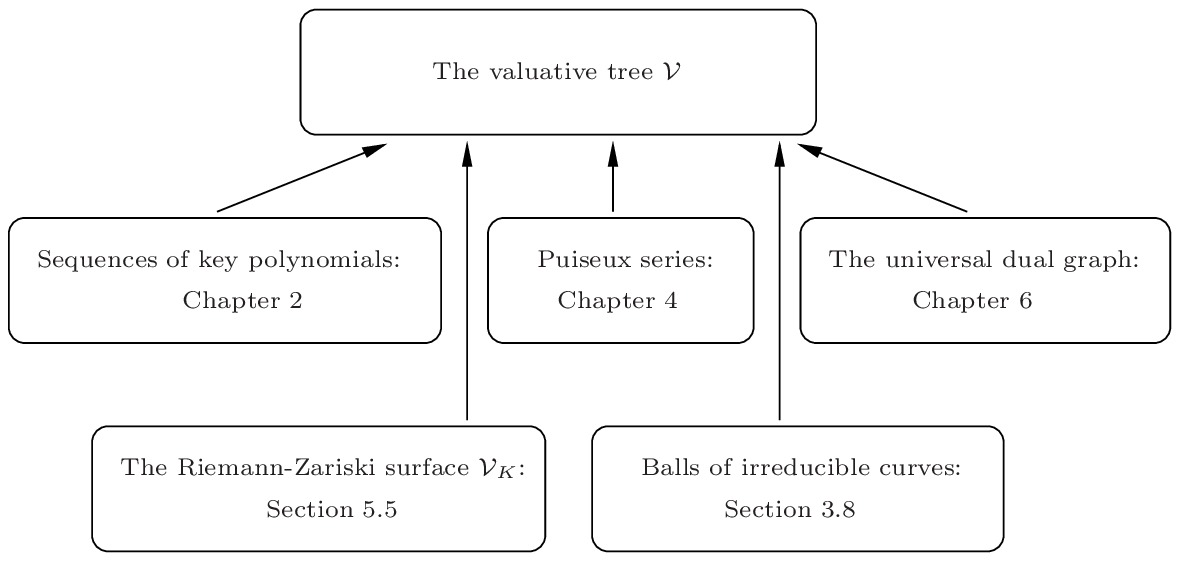}
\end{figure}
%
%
%
%
%
%
\chapter{Generalities}\label{part1}
In this chapter we give basic results on valuations. 
For definiteness we mostly 
restrict our attention to the ring of formal power series in 
two complex variables.

We distinguish between \emph{valuations} and \emph{Krull
valuations}\footnote{This is not standard terminology but will be used
throughout the monograph.}.  Valuations will always be $\R\cup\{\infty\}$
valued, whereas Krull valuations take values in an abstract totally
ordered group.  This distinction may seem unnatural but is useful for
our purposes.

We start in Sections~\ref{S25},~\ref{sec-val} and~\ref{krullval}
by giving precise definitions and noting basic facts on
valuations and Krull valuations. 
Most of this is standard and can be found 
in our main references which are~\cite{ZS} and~\cite{Vaquie}; 
see also~\cite{spiv}.
However,  we avoid using any deep facts from valuation 
theory and we have tried to make the exposition as self-contained
as possible.

As local formal curves play an important role in our study,
we spend some time discussing them in Section~\ref{S24}.

Then we try to give a ``road map'' to the valuative landscape
before embarking on the journey through it in subsequent chapters.
More precisely, we discuss in Section~\ref{ex}
six classes of Krull valuations that we shall later 
show form a complete list. 
The purpose of the discussion is to give the non-expert reader 
a feeling for the different valuations that will later 
appear in the monograph in various disguises.

We also try to give some analytic interpretations to valuations
whenever possible. For example, when applied to holomorphic germs,
valuations often give the the order of vanishing 
at the origin along a particular approach.
When applied to plurisubharmonic functions, many 
valuations can be viewed as Lelong numbers, (pushforward of) 
Kiselman numbers, or generalized Lelong numbers.
While these remarks are for the most part not necessary for
the purposes of the current monograph, they are important for
the applications.

In Section~\ref{versus} we make a precise comparison between 
valuations and Krull valuations.
Finally, Section~\ref{sec-equivalence}
contains a discussion of Zariski's correspondence between 
Krull valuations and sequences of infinitely nearby points.
%
%
%
%
\section{Setup}\label{S25}
Unless otherwise stated, $R$
\index{$R$ (ring)} 
will denote the ring of formal
power series in two complex variables.
This assumption on $R$ can be relaxed somewhat for the analysis
in this monograph; see Appendix~\ref{sec-notcomplete} for 
a discussion. Suffice it to say here that our method works for
the ring of holomorphic germs at the origin in $\C^2$, and for
the local ring at a smooth (closed) point of an algebraic surface 
over an algebraically closed field.
Certainly, most of the background results presented in 
this chapter hold in a quite general context. 

As is both convenient and customary (see~\cite{teissier-puis})
we sometimes ``ignore'' the fact that the power series in $R$ are not 
necessarily convergent. 
Thus, rather than talking about $\Spec R$ 
we most often write $(\C^2,0)$ and
view the elements of $R$ as functions on the latter space. 
In the same vein, we talk about local coordinates $(x,y)$ 
(rather than regular parameters) and local formal 
mappings $f:(\C^2,0)\to(\C^2,0)$ 
(rather than endomorphisms of $R$).

The ring $R$ is local, with unique maximal ideal $\fm$.
\index{$\fm$ (maximal ideal in $R$)}
If $(x,y)$ are local coordinates, then $R=\C[[x,y]]$ and $\fm=xR+yR$.
We let $K$
\index{$K$ (fraction field of $R$)} 
denote the fraction field of $R$.
%
%
%
%
\section{Valuations}\label{sec-val}
Denote the set of nonnegative real numbers by 
$\R_+$.
\index{$\R_+$ (nonnegative real numbers)}
Set $\Rbar\=\R_+\cup\{\infty\}$,
\index{$\Rbar$ ($\R_+\cup\{\infty\}$)}
extending the addition,
multiplication and order on $\R_+$ to $\Rbar$ in the usual way. 
\begin{definition}
  A \emph{valuation}\index{valuation} on $R$ is a nonconstant
  function $\nu:R\to\Rbar$ with: 
  \begin{enumerate}
  \item[(V1)]\label{p1} 
    $\nu(\phi\psi)=\nu(\phi)+\nu(\psi)$ for all
    $\phi,\psi\in R$; 
  \item[(V2)]\label{p2}
    $\nu(\phi+\psi)\ge\min\{\nu(\phi),\nu(\psi)\}$ 
    for all $\phi,\psi\in R$; 
  \item[(V3)]\label{p3}
    $\nu(1)=0$.
  \end{enumerate}
\end{definition}
\begin{remark}\label{R410}
  For applications to analysis and dynamics it is useful to
  think of a valuation as an \emph{order of vanishing} 
  at the origin along a particular approach direction. 
  See the examples in Section~\ref{ex} below.
\end{remark}
From the definition we see that $\nu(0)=\infty$ and $\nu|_{\C^*}=0$.
The set $\fp\=\{\phi\in R\ ;\ \nu(\phi)=\infty\}$ is a prime ideal
in $R$; we say $\nu$ is 
\emph{proper}\index{valuation!proper} 
if $\fp\subsetneq\fm$. 
It is
\emph{centered}\index{valuation!centered} if it is proper and if
$\nu(\fm)\=\min\{\nu(\phi)\ ;\ \phi\in\fm\}>0$.  Two centered
valuations $\nu_1$, $\nu_2$ are 
\emph{equivalent}\index{valuation!equivalent}, $\nu_1\sim\nu_2$,
if $\nu_1=C\nu_2$ for some constant $C>0$.  

We denote by $\tcV$\index{$\tcV$ (centered valuations)} 
the set of all centered valuations on $R$.
The primary object of study in this monograph will be 
quotient of $\tcV$ by the equivalence relation $\sim$.
We call it \emph{valuation space}\index{valuation space} 
for the time being, although we shall later give it the name
\emph{the valuative tree}\index{valuative tree} 
as it possesses a tree structure.

Notice that $\tcV$ is naturally a sort of bundle over 
$\tcV/\sim$ with 
fibers isomorphic to the real interval $(0,\infty)$.
For our purposes it will be important to pick a section of this bundle,
or, what amounts to the same thing, 
view $\tcV/\sim$ as a subset of $\tcV$ 
singled out by a normalizing condition.
There are several possible normalizations.
One is to pick any element $x\in\fm$ that defines a smooth curve 
$\{x=0\}$ and demand that $\nu(x)=1$.
This is useful in cases
where the curve $\{x=0\}$ is naturally given, \eg when it 
is an irreducible component of 
the exceptional divisor of some composition of blowups.
We thus define $\cV_x$
\index{$\cV_x$ (relative valuation space/valuative tree)}
to be the set of centered valuations $\nu$ on $R$ satisfying
$\nu(x)=1$.\footnote{In fact we shall later add to this space
  one valuation which is not centered: see Section~\ref{sec-relative}.}
We shall explore this space further in Section~\ref{sec-relative}.

However, the most important normalization in this monograph 
will be $\nu(\fm)=1$, \ie $\min\{\nu(x),\nu(y)\}=1$
for any local coordinates $(x,y)$.
We will then say that $\nu$ is 
\emph{normalized}\index{valuation!normalized}.
This normalization has the
advantage of being coordinate independent.
We define $\cV$
\index{$\cV$ (valuation space/valuative tree)}
as the set of centered valuations $\nu$ normalized by $\nu(\fm)=1$.

We endow $\tcV$, $\cV$ and $\cV_x$ with the 
weak 
topology\index{topology!weak (on $\cV$)}:
if $\nu_k,\nu$ are valuations in $\tcV$ ($\cV$, $\cV_x$), 
then $\nu_k\to\nu$ in $\tcV$ ($\cV$, $\cV_x$) 
iff $\nu_k(\phi)\to\nu(\phi)$ for all $\phi\in R$. 
As we will see later, the weak topologies 
on $\cV$ and $\cV_x$ are induced by tree structures.

Finally $\tcV$, $\cV$ and $\cV_x$ come with natural partial
orderings\index{partial ordering!on the valuative tree}: 
$\nu_1\le\nu_2$ iff $\nu_1(\phi)\le\nu_2(\phi)$ for all $\phi\in R$.
Again, the partial orderings on $\cV$ and $\cV_x$
arise from tree structures. 
%
%
%
%
\section{Krull valuations}\label{krullval}
We now introduce Krull valuations.\footnote{The reader mainly 
  interested in the analytic applications of  valuations may skip 
  this section on a first reading.}
They are defined in 
the same way as valuations, except that we replace $\Rbar$ with a
general totally ordered group. We shall analyze to what extent
this really makes a difference in Section~\ref{versus}. 
\begin{definition}
  Let $\Gamma$ be a totally ordered abelian group.  
  A \emph{Krull valuation}\index{valuation!Krull} on $R$
  is a function $\nu:R\setminus\{0\}\to\Gamma$ satisfying
  (V1)-(V3) above.
\end{definition}
A Krull valuation is \emph{centered}\index{valuation!centered Krull} 
if $\nu\ge0$ on $R$ and $\nu>0$ on $\fm$.
Two Krull valuations $\nu_1:R\to\Gamma_1$, $\nu_2:R\to\Gamma_2$ are
\emph{equivalent}\index{valuation!equivalent Krull} 
if $h\circ\nu_1=\nu_2$ for 
some strictly increasing homomorphism $h:\Gamma_1\to\Gamma_2$.
Any Krull valuation extends naturally to the fraction field $K$ 
of $R$ by $\nu(\phi/\psi)=\nu(\phi)-\nu(\psi)$.

A \emph{valuation ring}\index{ring!valuation} $S$ is a local ring with
fraction field $K$ such that $x\in K^*$ implies $x\in S$ or 
$x^{-1}\in S$.  
Say that $x,y\in S$ are equivalent iff $xy^{-1}$ is a unit in $S$. 
The quotient $\Gamma_S$ is endowed with a natural total order
given by $x\ge y$ iff $xy^{-1}\in\fm_S=\max S$. The projection
$S\to\Gamma_S$ then extends to a Krull valuation
$\nu_S:K\to\Gamma_S$.

Conversely, if $\nu$ is a centered Krull valuation, 
then $R_\nu\=\{\phi\in K\ ;\ \nu(\phi)\ge0\}$ 
is\index{$R_\nu$ (valuation ring of $\nu$)} 
a valuation ring with maximal 
ideal\index{$\fm_\nu$ (maximal ideal in $R_\nu$)}
$\fm_\nu\=\{\phi\ ; \nu(\phi)>0\}$.
One can show that $\nu$ is equivalent to $\nu_{R_\nu}$ defined above. 
In particular, two Krull valuations $\nu,\nu'$ are equivalent
iff $R_\nu=R_{\nu'}$.

We let $\tcV_K$\index{$\tcV_K$ (centered Krull valuations)} 
be the set of all centered Krull valuations,
and $\cV_K$\index{$\cV_K$ (centered Krull valuations modulo equivalence)} 
be the quotient of $\tcV_K$ by the equivalence
relation. Equivalently, $\cV_K$ is the set of valuation
rings $S$ in $K$ whose maximal ideal $\fm_S$ satisfies 
$R\cap\fm_S=\fm$.

The group $\nu(K)$ is called the 
\emph{value group}\index{valuation!value group/semigroup of} of $\nu$; 
similarly $\nu(R)$ is the \emph{value semigroup}.

A classical way to analyze a Krull valuation is through its
\emph{numerical invariants}.
\index{valuation!numerical invariants of}
\begin{definition}
  Let $\nu:R\to\Gamma$ be a centered Krull valuation.
  \begin{itemize}
  \item The \emph{rank}\index{valuation!rank of} 
    $\rk(\nu)$\index{$\rk$ (rank)} of $\nu$ is 
    the Krull dimension of the ring $R_\nu$.
  \item The \emph{rational rank}\index{valuation!rational rank of} 
    of $\nu$ is defined by\index{$\ratrk$ (rational rank)}.
    $\ratrk(\nu)\=\dim_{\Q}(\nu(K)\otimes_{\Z}\Q)$
  \item The \emph{transcendence degree}
    \index{valuation!transcendence degree of} 
    $\trdeg(\nu)$\index{$\trdeg$ (transcendence degree)} of $\nu$ is
    defined as follows. 
    Since $\nu$ is centered, we have a natural inclusion
    $\C=R/\fm\subset k_\nu\=R_\nu/\fm_\nu$. The field $k_\nu$ is 
    called the \emph{residue field} of $\nu$. 
    \index{valuation!residue field of}
    \index{$k_\nu$ (residue field)}
    We let $\trdeg(\nu)$ be the transcendence degree of this field
    extension.
  \end{itemize}
\end{definition}
One can show that $\rk(\nu)$ is the least integer $l$ so that $\nu(K)$
can be embedded as an ordered group into $(\R^l,+)$ endowed with the
lexicographic order. Hence both $\rk(\nu)$ and $\ratrk(\nu)$ 
depend only on the value group $\nu(K)$ of the valuation. 
We will later give a geometric interpretation of $\trdeg(\nu)$ 
(see Remark~\ref{residue}).

The main relation between these numerical invariants is given by
Abhyankar's inequalities,\footnote{Although we shall not prove this
formula, one may use results from Chapter~\ref{part2} to do so (see
Theorems~\ref{divis},~\ref{approx}).}\index{Abhyankar's inequalities}
which assert that:
\begin{equation}\label{AbIn}
  \rk(\nu)+\trdeg(\nu)\le\ratrk(\nu)+\trdeg(\nu)\le\dim R=2.
\end{equation}

Moreover, if $\ratrk(\nu)+\trdeg(\nu)=2$, 
then $\nu(K)$ is isomorphic (as a group) to $\Z^\e$ with $\e=\ratrk(\nu)$.
When $\rk(\nu)+\trdeg(\nu)=2$, $\nu(K)$ is isomorphic as
an \emph{ordered} group to $\Z^\e$ endowed with the lexicographic
order.  

To $\nu$ we associate a 
\emph{graded ring}\index{ring!graded}\index{$\gr_\nu R$ (graded ring)}
$\gr_\nu R=\oplus_{r\in\Gamma}\{\nu\ge r\}/\{\nu>r\}$.
It is in bijection with equivalence classes of $R$ under the 
equivalence relation $\phi=\psi$ \emph{modulo} $\nu$
iff $\nu(\phi-\psi)>\nu(\phi)$.
We have $\overline{\phi\cdot\psi}=\overline{\phi}\cdot\overline{\psi}$ 
in $\gr_\nu R$ for any $\phi,\psi\in R$.

Any formal mapping $f:(\C^2,0)\to(\C^2,0)$ induces a 
ring homomorphism $f^*:R\to R$. The latter in turn
induces actions on $\tcV$ and $\tcV_K$ given by 
$f_*\nu(\phi)=\nu(f^*\phi)$.\index{$f_*$ (pushforward on valuations)}
When $f$ is invertible
$f_*$ is a bijection and preserves
all three invariants $\rk$, $\ratrk$ and $\trdeg$ defined above.
It also restricts to a bijection on $\cV$ preserving the 
natural partial order: $f_*\nu\le f_*\nu'$ as soon as $\nu\le\nu'$.
%
%
%
%
\section{Plane curves}\label{S24}
Curves play a key role in our approach to valuations.
Indeed, we shall later see that ``most'' valuations can be seen
as balls of irreducible curves in a particular metric.
Here we recall a few classical results regarding plane curves: 
see~\cite{teissier-puis} for details.

Recall that $R$ is the ring of formal power series in two complex
variables. Thus a \emph{curve}\index{curve} $C$ for us is simply
an equivalence class of elements $\phi\in\fm$, where two elements
are equivalent if they differ by multiplication by a unit in $R$,
\ie by an element of $R\setminus\fm$. 
We write $C=\{\phi=0\}$ even though $\phi$ is not necessarily a
convergent power series (see Section~\ref{S25}) and say that the
curve is \emph{represented} by $\phi$.

The \emph{multiplicity}\index{curve!multiplicity of} 
\index{multiplicity! of a curve}
$m(C)$\index{$m(C)$ (multiplicity)} 
of a curve $C=\{\phi=0\}$ 
is defined by $m(C)=m(\phi)=\max\{n\ ;\ \phi\in\fm^n\}$; this does not
depend on the choice of $\phi$.\index{$m(\phi)$ (multiplicity)}
Recall that $R$ is a unique factorization domain. 
A curve is \emph{reduced}\index{curve!reduced} if it is represented by
an element of $\fm$ without repeated irreducible factors.
It is \emph{irreducible}\index{curve!irreducible} if it is
represented by an irreducible element of $\fm$.

Any irreducible curve $C$ admits a 
\emph{parameterization}\index{curve!parameterization of}
as follows. 
Pick local coordinates $(x,y)$ such that $C$ is transverse to 
the curve $\{x=0\}$, \ie if $C=\{\phi(x,y)=0\}$, then 
$\phi(0,y)$ is a formal power series in $y$ whose lowest order term
has degree $m=m(C)$. Then there exists a formal power series
$y(t)\in\C[[t]]$ such that $\phi(t^m,y(t))=0$. 

The \emph{intersection multiplicity}\index{curve!intersection
multiplicity of} of two curves $C=\{\phi=0\}$ and $D=\{\psi=0\}$ is
defined by\index{$C\cdot D$ (intersection multiplicity)} 
$C\cdot D=\dim_\C R/\langle\phi,\psi\rangle$, where
$\langle\phi,\psi\rangle$ denotes the ideal of $R$ generated by $\phi$
and $\psi$; this intersection multiplicity does not depend on the
choice of $\phi$ and $\psi$.  We sometimes write $\phi\cdot\psi$ for
$\{\phi=0\}\cdot\{\psi=0\}$.  Then $\phi\cdot\psi=\infty$ iff $\phi$
and $\psi$ have a common irreducible factor.

The intersection multiplicity can also be computed in terms of
parameterizations. Suppose $C$ is an irreducible curve, $(x,y)$ are
coordinates such that $C$ is transverse to $\{x=0\}$ and
$(x,y)=(t^m,y(t))$ is a parameterization of $C$.  Then for any curve
$D=\{\psi=0\}$ we have that $C\cdot D$ is the lowest order term in the
formal power series $\psi(t^m,y(t))\in\C[[t]]$.  In particular when
$D$ is smooth and transverse to $C$ then $C\cdot D = m(C)$.

%
%
%
%
\section{Examples of valuations}\label{ex}
It is now time to present examples of valuations and Krull
valuations. In fact, the list that we give below is complete, but that
will only be proved in Chapter~\ref{part2}, see also Appendix~\ref{sec-clas}.
The purpose of presenting the list below is not so much to
give precise definitions, but rather to introduce some notation and
give a road map to the valuative landscape that lies ahead.  In
particular we do not hesitate to mention features that are not obvious
from the definitions.

We also give some hints to how the valuations can be interpreted
analytically when applied to holomorphic germs or even to
plurisubharmonic functions.  This analytic point of view is not
strictly necessary for anything that we do in this monograph, but
hopefully serves to add to the reader's intuition and explain how
valuations can be used in applications such as
in~\cite{pshfnts}~\cite{criterion}.

%
%
\subsection{The multiplicity valuation}
The function
\begin{equation*}
  \nu_\fm(\phi )\=\max\{k\ ;\ \phi\in\fm^k\}, 
\end{equation*}
defines both a normalized valuation and a 
Krull valuation.\index{$\nu_\fm$ (multiplicity valuation)}
We call it the 
\emph{multiplicity valuation}.\index{multiplicity valuation}
Notice that $\nu_\fm(\phi)=m(\phi)$ is the multiplicity of the curve
$\{\phi=0\}$ as defined in Section~\ref{S24}.

The multiplicity valuation is the root of the valuative tree in the
sense that it is the minimal element under the partial ordering
defined above. See Proposition~\ref{P107}.

The numerical invariants are easy to compute 
(alternatively see Theorem~\ref{divis}). 
Clearly $\nu_\fm(K)=\Z$ so $\rk(\nu_\fm)=\ratrk(\nu_\fm)=1$.
It is not hard to see that the residue field is 
given by $k_\nu\simeq\C(y/x)$, where $(x,y)$ are local coordinates.
Hence $\trdeg(\nu_\fm)=1$.

Analytically, the multiplicity valuation can be interpreted as
the order of vanishing at the origin: 
if $\phi$ is a holomorphic germ, then
\begin{equation}\label{e438}
  \nu_\fm(\phi)=\lim_{r\to0}\frac1{\log r}\sup_{B(r)}\log|\phi|,
\end{equation}
where $B(r)$ is a ball of radius $r$ centered at the origin.  In fact,
this equation makes sense even when $\log|\phi|$ is replaced by a
plurisubharmonic function $u$ and recovers the 
\emph{Lelong number}~\cite{dem2}\index{Lelong number} 
of $u$ at the origin.
%
%
\subsection{Monomial valuations}\label{S31}
Fix local coordinates $(x,y)$ and $\a\ge1$.  Define a valuation
$\nu_{y,\a}$ as follows:\footnote{We will see later that $\nu_{y,\a}$
does not depend on the choice of the other coordinate $x$.}
$\nu_{y,\a}(x)=1$, $\nu_{y,\a}(y)=\a$ and
\begin{equation*}
  \nu_{y,\a}(\phi)\=\min\{i+\a j\ ;\ a_{ij}\ne0\},
\end{equation*}
for $\phi=\sum a_{ij}x^iy^j$.
Then $\nu_{y,\a}$ is a valuation and a Krull valuation: 
we say that $\nu_{y,\a}$ 
is \emph{monomial}\index{valuation!monomial} 
in the coordinates $(x,y)$.

In general, a valuation $\nu$ is monomial if it is
monomial in some local coordinates $(x,y)$. 
Notice that since $\a\ge1$, $\nu_{y,\a}$ is normalized
in the sense that $\nu_{y,\a}(\fm)=1$.
If $\a=1$, then $\nu_{y,1}=\nu_\fm$ is the multiplicity valuation.

A monomial valuation has numerical invariants given as follows 
(see Theorem~\ref{divis}).
We have $\rk(\nu_{y,\a})=1$ for any $\a$.  
When $\a\in\Q$, $\trdeg(\nu_{y,\a})=1$ and $\ratrk(\nu_{y,\a})=1$; 
when $\a\not\in\Q$, $\trdeg(\nu_{y,\a})=0$ and $\ratrk(\nu_{y,\a})=2$. 

Monomial valuations have the following analytic interpretation.
Assume $(x,y)$ are holomorphic coordinates (this is in fact not
a restriction as we may perturb $x$ and $y$ slightly without 
changing $\nu$). 
If $\phi$ is a holomorphic germ then
\begin{equation}\label{e439}
  \nu_{y,\a}(\phi)
  =\lim_{r\to0}\frac1{\log r}\sup_{|x|<r,|y|<r^\a}\log|\phi|,
\end{equation}
As before, this equation makes sense even when $\log|\phi|$ is
replaced by a plurisubharmonic function $u$ and recovers
the \emph{Kiselman number}~\cite{dem2}
\index{Kiselman number} 
\index{Lelong number!directional} 
of $u$ (with weights $(\a,1)$) at 
the origin.\footnote{Kiselman numbers are also sometimes called
  \emph{directional Lelong numbers}.}
%
%
\subsection{Divisorial valuations}\label{S28}
If $C=\{\psi=0\}$ is an irreducible curve, then
\begin{equation*}
  \div_C(\phi):=\max\{k\ ;\ \psi^k\mid\phi\}
\end{equation*}
defines a valuation 
on $R$.\index{$\div_C$ (order of vanishing along $C$)}
In other words, $\div_C$
is the order of vanishing along $C$.
Notice that such a valuation is not centered,
but it is the basis for a whole class of centered 
valuations defined as follows.

Consider a birational morphism $\pi:X\to(\C^2,0)$.  This means that
$X$ is a smooth complex surface, $\pi$ is holomorphic and proper, and
$\pi$ is a biholomorphism outside the exceptional divisor
$\pi^{-1}(0)$.  It is a classical theorem of Zariski that such a $\pi$
is a composition of point blowups.  Let $E$ be an irreducible
component of the exceptional divisor. In the sequel, we 
sometimes refer to such a divisor as an 
\emph{exceptional component}\index{exceptional component}. 
There is then a unique integer $b=b_E>0$ 
such that\index{$\nu_E$ (divisorial valuation)}
\begin{equation*}
  \nu_E(\phi):=b^{-1}\,\div_E(\pi^*\phi)
\end{equation*}
defines a normalized valuation (and a Krull valuation) on $R$, 
called a \emph{divisorial} valuation.\index{valuation!divisorial}
Thus $b_E\,\nu_E(\phi)$ is the order of vanishing of 
the pullback $\pi^*\phi$ along $E$.

The multiplicity valuation $\nu_\fm$ is divisorial, 
with $\pi:X\to(\C^2,0)$ being a single blowup of the origin,
$E=\pi^{-1}(0)$ and $b_E=1$.

The numerical invariants of $\nu_E$ are given as follows:
$\rk(\nu_E)=\ratrk(\nu_E)=1$ and $\trdeg(\nu_E)=1$. 
In fact the residue field of $\nu$ is canonically isomorphic 
to the field of rational functions on $E$: to the class
of $\phi\in K$ is associated the rational function
$\pi^*\phi|_E$.

Any valuation with the invariants above is in fact divisorial
(see~\cite[p.89]{ZS}).  We shall give an argument below: 
see Proposition~\ref{Pchardiv}. In particular a monomial valuation
$\nu_{y,\a}$ is divisorial iff $\a\in\Q$.

As we will show (Proposition~\ref{P107}), a divisorial valuation
$\nu_E$ constitutes a branch point in the valuative tree.  The tree
tangent space at $\nu_E$, \ie the set of branches emanating from
$\nu_E$, is naturally in bijection with the set of points on $E$
(Theorem~\ref{div-tangent}).

Analytically, divisorial valuations can be interpreted as follows.
If $C$ is a holomorphic curve and $\phi$ a holomorphic germ, 
then $\div_C(\phi)$ is the order of vanishing of $\phi$ at
a generic point of $C$.
Similarly, $b_E\,\nu_E(\phi)$ the order of vanishing of 
$\pi^*\phi$ at a generic point of $E$.

These interpretations generalize to plurisubharmonic functions $u$
(an example is $u=\log|\phi|$ for a holomorphic germ $\phi$).
If $C$ is a holomorphic curve, then $\div_C(u)$ is the 
Lelong number of $u$ at a generic point of $C$, or
equivalently the (normalized) mass of the current $dd^cu$ on $C$.
Similarly, $b_E\,\nu_E(u)$ is both the Lelong number of 
$\pi^*u$ at a generic point of $E$ and 
the normalized mass of the pullback $\pi^*dd^cu$ along $E$.

Divisorial valuations are special cases of quasimonomial valuations,
to be discussed next. This gives yet another analytic interpretation.
%
%
\subsection{Quasimonomial valuations}\label{S32}
As in the definition of a divisorial valuation above,
consider a birational morphism $\pi:X\to(\C^2,0)$ and
an exceptional component $E\subset\pi^{-1}(0)$.
Pick a point $p\in E$ and let $\mu$ be a monomial valuation at $p$. 
The image $\nu=\pi_*\mu$ defined by $\nu(\phi)=\mu(\pi^*\phi)$ 
is a \emph{quasimonomial} valuation.\index{valuation!quasimonomial}
Its numerical invariants are the same as those of $\mu$ 
(see Theorem~\ref{divis}).  
In particular, quasimonomial valuations for which the
the second inequality in~\eqref{AbIn} above is an equality. 
Such valuations were called rank one
Abhyankar valuations in~\cite{ELS}.\index{valuation!Abhyankar}
\index{Abhyankar valuation} 
Conversely, any rank one Abhyankar valuation is a 
quasimonomial valuation, see Theorem~\ref{divis}.

Beware that we will use an alternative definition of quasimonomial
valuations in Section~\ref{S26} and it will take until
Chapter~\ref{A3} before we see that the two definitions are
equivalent.

Any divisorial valuation is quasimonomial. 
A quasimonomial valuation that is not divisorial is called
\emph{irrational}\index{valuation!irrational}. Such a valuation
has numerical invariants $\rk(\nu)=1$, $\ratrk(\nu)=2$ and
$\trdeg(\nu)=0$, see Theorem~\ref{divis}.

We already observed that divisorial valuations are branch points
in the valuative tree. The irrational valuations are
regular points in the sense that there are exactly 
two branches emanating from each such point. 
See Proposition~\ref{P107}. In fact, the quasimonomial 
valuations are exactly the points in the valuative tree that
are not ends. The set of quasimonomial valuations therefore
constitute a subtree $\cVqm$ of $\cV$. 

Quasimonomial valuations are useful for applications as they can
be interpreted as the order of vanishing at the origin along a
particular approach direction. 
More precisely, we can associate 
to any quasimonomial valuation $\nu$ 
a \emph{characteristic region}\index{characteristic region}
of the form 
\begin{equation*}
  \Omega(r)
  =\left\{(x,y)\in\C^2\ ;\ |x|<r,\ |\phi(x,y)|<|x|^{tm}\right\}
  \subset\C^2.
\end{equation*}
Here $\phi\in\fm$ is an irreducible \emph{holomorphic} germ,
$(x,y)$ are local (holomorphic) coordinates 
such that the curves $\{x=0\}$ and $\{\phi=0\}$ are transverse,
$m=m(\phi)$ is the multiplicity of $\phi$ and $t>1$.
For $r>0$ small, the region $\Omega(r)$ 
is a small neighborhood of the curve $\{\phi=0\}$ with the origin removed. 
See Figure~\ref{F7} on page~\pageref{F7}.
Then $\nu(\psi)$ is given by 
\begin{equation}\label{e437}
  \nu(\psi)=\lim_{r\to0}\frac1{\log r}\sup_{\Omega(r)}\log|\psi|.
\end{equation}
In fact, this limit exists even when $\log|\psi|$ is replaced
by a general plurisubharmonic function, and this is the basis 
for the valuative study of singularities of plurisubharmonic functions.
The fact that a quasimonomial valuation is the pushforward of
a monomial valuation means that, analytically, a quasimonomial
valuation is the pushforward of a Kiselman number.
%
%
\subsection{Curve valuations}\label{S33}
Any irreducible curve defines a valuation by intersection.
More precisely, 
let $C$ be an irreducible curve of multiplicity $m(C)$ and 
define\index{$\nu_C$ (curve valuation)}
\begin{equation*}
  \nu_C(\psi)=\frac{C\cdot\{\psi=0\}}{m(C)},
\end{equation*}
using the intersection multiplicity between curves as defined
in Section~\ref{S24}. Then $\nu_C$ is a centered valuation
on $R$ and normalized in the sense that $\nu_C(\fm)=1$.
We call $\nu_C$ a \emph{curve valuation}.\index{valuation!curve} 
If $C=\{\phi=0\}$ then we also write 
$\nu_C=\nu_\phi$.\index{$\nu_\phi$ (curve valuation)}
Notice that $\nu_\phi(\psi)=\infty$ iff $\phi$ divides $\psi$.

Curve valuations will play a very important role in our approach.
Indeed, a quasimonomial valuation can be accurately understood
in terms of the curve valuations that dominate it.
This explains Spivakovsky's observation~\cite{spiv} that
the classification of valuations and curves 
is essentially the same, see Section~\ref{val-curve}.

If the curve $C$ is holomorphic, then the value of $\nu_C$ on
a holomorphic germ $\psi$ can be interpreted as the 
(normalized) order of vanishing of $\psi$ at the origin of
the restriction $\psi|_C$; compare~\eqref{e437}.

As a curve valuation $\nu_C=\nu_\phi$ can take infinite values, 
it is not strictly speaking
a Krull valuation, but it can be turned into one as follows.
For $\psi\in R$, write $\psi=\phi^k\hpsi$ 
with $k\in\N$, $\hpsi$ prime with $\phi$, and define
\begin{equation*}
  \nu'_C(\psi)\=(k,\nu_\phi(\hpsi))\in\Z\times\Q.
\end{equation*}
Then $\nu'_C=\nu'_\phi$ is a centered Krull valuation, 
where $\Z\times\Q$ is lexicographically ordered.
The numerical invariants of $\nu'_C$ are
$\rk(\nu'_C)=\ratrk(\nu'_C)=2$ and $\trdeg(\nu'_C)=0$
(use~\eqref{AbIn} or see Theorem~\ref{divis}).
%
%
\subsection{Exceptional curve valuations}
Consider a birational morphism $\pi:X\to(\C^2,0)$ as above
and a curve valuation $\mu$ at a point $p$ on
the exceptional divisor $\pi^{-1}(0)$. If the corresponding
curve is not contained in the exceptional divisor, then
the pushforward of $\mu$ under $\pi$ is proportional to the
curve valuation at the image curve. 

If the curve is instead an irreducible component $E$
of the exceptional divisor, then the pushforward of the
curve valuation is not proper. However, the pushforward of
the corresponding Krull valuation gives a new type
of Krull valuation.

Let $\div_E$ denote the order of vanishing along $E$ and
$\nu_E\in\cV$ the normalized divisorial valuation as 
in Section~\ref{S28}.
Pick local coordinates $(z,w)$ at $p$ such that $E=\{z=0\}$.
For $\phi\in R$ write $\pi^*\phi=z^{\div_E(\pi^*\phi)}\psi$,
where $\psi\in\C[[z,w]]$.
Thus $\psi(0,w)=w^k\hpsi(w)$, where $1\le k<\infty$ and $w\nmid\hpsi$.
Then\index{$\nu_{E,p}$ (exceptional curve valuation)}
\begin{equation*}
  \nu_{E,p}(\phi)\=(\nu_E(\phi),k)\in\Z\times\Z
\end{equation*}
defines a centered Krull valuation called an
\emph{exceptional curve valuation}.\index{valuation!exceptional curve}
Its numerical invariants are the same as those of a curve valuation,
\ie $\rk(\nu_{E,p})=\ratrk(\nu_{E,p})=2$ and $\trdeg(\nu_{E,p})=0$.

As we shall see in Lemma~\ref{L601}, 
the exceptional curve valuations are exactly the 
Krull valuations that are not valuations. They are therefore
not points in the valuative tree, but can be interpreted as
tangent vectors at divisorial valuations. 
See Appendix~\ref{sec-tangent}.
%
%
\subsection{Infinitely singular valuations}
The remaining centered valuations on $R$ are \emph{infinitely
singular} valuations.\index{valuation!infinitely singular} We shall
define them in Section~\ref{S26} but they can be characterized in a
number of equivalent ways. One way is through their numerical
invariants: $\ratrk=\rk=1$ and $\trdeg=0$ (see Theorem~\ref{divis}).
Another way is to say that their value (semi)groups are not finitely
generated.  Examples of valuations with this latter property are given
in~\cite[p.102-104]{ZS}.

Yet another, perhaps more illustrative, way is through Puiseux
series. Namely, any (normalized) infinitely singular valuation $\nu$
is represented by a (generalized) Puiseux series $\hphi$ as follows.
Pick local coordinates $(x,y)$ such that $\nu(x)=1$, $\nu(y)\ge1$.
Then $\hphi$ is a series of the form
$\hphi=\sum_{j=1}^\infty a_jx^{\hbeta_j}$,
where $a_j\in\C^*$ and $(\hbeta_j)_1^\infty$ is a strictly
increasing sequence of positive rational numbers with 
unbounded denominators.
If $\psi=\psi(x,y)\in\fm$, then 
$\psi(x,\hphi)$ defines a series of the same form and 
$\nu(\psi)$ is the minimum order of the terms in this series.

See Chapter~\ref{sec-puis} for details on the Puiseux
point of view and Appendix~\ref{sec-inf-sing} for even more 
characterizations of infinitely singular valuations.
In terms of the valuative tree, the infinitely singular 
valuations are ends, sharing the latter property with 
the curve valuations. In a way, the infinitely singular valuations
can be viewed as curve valuations associated to ``curves'' of
infinite multiplicity.\footnote{At least this is the reason for
their names.}

Infinitely singular valuations are nontrivial to interpret analytically,
but they can sometimes be understood as a 
generalized Lelong number
\index{Lelong number!generalized}
in the sense of Demailly, with respect to a plurisubharmonic
weight~\cite{dem2}.
%
%
%
%
\section{Valuations versus Krull valuations}\label{versus}
As we have made a point of distinguishing between valuations
and Krull valuations, it is reasonable to ask what the relation
is between the two concepts.\footnote{The 
  reader mainly interested in the analytic applications of 
  valuations may skip this section on a first reading.}

In order to understand the situation we will make use of the
following result, whose proof is postponed until the
end of the section.
\begin{lemma}\label{L601}
  Let $\mu,\nu$ be non-trivial Krull valuations on $R$ with
  $R\subset R_\mu\subsetneq R_\nu$. Then either 
  $\nu\sim\nu_E$ is divisorial and $\mu\sim\nu_{E,p}$ 
  is an exceptional curve valuation;
  or there exists a curve $C$ such that 
  $\mu\sim\nu'_C$ is a curve valuation and $\nu\sim\div_C$.
\end{lemma}

First consider a centered valuation $\nu:R\to\Rbar$.
In order to see whether or not 
it defines a Krull valuation we
consider the prime ideal $\fp=\{\nu=\infty\}$. 
As $\nu$ is proper, $\fp\subsetneq\fm$, so either
$\fp=(0)$ or $\fp=(\phi)$ for an irreducible $\phi\in\fm$.

When $\fp=(0)$, $\nu$ defines a Krull valuation whose value group is
included in $\R$. In particular $\rk(\nu)=1$.  When $\fp=(\phi)$ is
nontrivial, we may define a Krull valuation
$\krull[\nu]:R\to\Z\times\R$ as in Section~\ref{S33}.  
Namely, for $\psi\in R$, write $\psi=\phi^k\hat{\psi}$ with $k\in\N$, 
$\hat{\psi}$ prime with $\phi$, and set 
$\krull[\nu](\psi)=(k,\nu(\hat{\psi}))$.
Then Lemma~\ref{L601} applies to the pair
$\krull[\nu]$, $\div_\phi$. 
This shows that $\krull[\nu]$ is the associated 
Krull valuation $\nu'_\phi$.  
As a consequence, $\nu$ is equivalent to the 
curve valuation $\nu_\phi$.

It is straightforward to check that if two valuations give rise
to equivalent Krull valuations, then the valuations are themselves
equivalent. In particular we may define the numerical invariants
$\rk$, $\ratrk$ and $\trdeg$ for all valuations.

Conversely, pick a Krull valuation $\nu$. Let us see if it
arises from a valuation by the preceding procedure.

When $\rk(\nu)=1$ \ie when $\fp:=\{\nu=\infty\}=(0)$, 
the value group $\nu(K)$ of $\nu$ can be embedded 
(as an ordered group) in $(\R,+)$ so that $\nu$ defines a valuation.
When $\rk(\nu)=2$, two cases may appear.

(1) There are non-units $\phi_1,\phi_{\infty}\in R$ with
    $n\nu(\phi_1)<\nu(\phi_{\infty})$ for all $n\ge0$. 
    Then 
    \begin{equation*}
      \fp\=\{\phi\in R \ ;\ n\nu(\phi_1)<\nu(\phi)\ \text{for all}\ n\ge0\}
    \end{equation*}
    is a proper prime ideal generated by an irreducible element 
    we again denote by $\phi_\infty$. 
    We define a valuation $\nu_0:R\to\Rbar$ as
    follows. First set $\nu_0(\phi_1)=1$, and $\nu_0(\phi)=\infty$
    for $\phi\in\fp$. For any $\phi\in R\setminus\fp$, let
    \begin{equation*}
      n_k\=\max\{n\in\N\ ;\ n\nu(\phi_1)\le\nu(\phi^k)\}.
    \end{equation*}
    One checks that 
    $k^{-1}n_k\le (kl)^{-1}n_{kl}\le k^{-1}n_k+k^{-1}-(kl)^{-1}$ 
    for all $k,l\ge0$. 
    This implies that the sequence
    $n_k/k$ converges towards a real number we define to be 
    $\nu_0(\phi)$. The function $\nu_0$ is a valuation, 
    and $\krull[\nu_0]$ is isomorphic to $\nu$.

(2) For all non-units $\phi,\psi\in R$ there exists $n\ge0$ so that
    $n\nu(\psi)>\nu(\phi)$. In this case, there can be no valuation
    $\nu_0:R\to\Rbar$ so that $\krull[\nu_0]\sim\nu$.  Let us show
    that $\nu$ is an exceptional curve valuation.
    As $\rk\,\nu=2$, we may assume $\nu(R\setminus\{0\})\subset\R\times\R$, 
    and $\nu=(\nu_0,\nu_1)$. 
    The function $\nu_0:R\setminus\{0\}\to\R$ is also a Krull valuation, 
    as $\nu$ is and $\R\times\R$ is
    endowed with the lexicographic order. 
    Since $R_\nu\subsetneq R_{\nu_0}$, Lemma~\ref{L601}
    shows that $\nu_0$ is divisorial and $\nu$ an exceptional
    curve valuation. 
    
We have hence proved
\begin{proposition}\label{prop-versus}
  To any Krull valuation $\nu$ which is not an exceptional
  curve valuation is associated a unique (up to equivalence) 
  valuation $\nu_0:R\to\Rbar$ with $\krull[\nu_0]\sim\nu$. 
\end{proposition}
\begin{remark}\label{R411}
  In Section~\ref{sec-zar} we shall strengthen this result and
  show that $\cV$, endowed with the weak topology, 
  is homeomorphic to the quotient space resulting from
  $\cV_K$ with the Zariski topology
  after having identified all exceptional curve valuations with their
  associated divisorial valuations. See Theorem~\ref{weak-zariski}.
\end{remark}
\begin{proof}[Proof of Lemma~\ref{L601}]
  Whereas we have tried to keep the monograph as self-contained as
  possible, the following proof does make use of some (well-known)
  results from the literature. 
  Lemma~\ref{L601} is only used to explain the connections
  between Krull valuations and valuations and is not necessary for
  obtaining most of the structure on the valuative tree.
  
  The proof goes as follows. 
  As $R_\mu\subsetneq R_\nu$, the ideal $\fq\=\fm_\nu\cap R_\mu$ is
  prime, strictly included in $R_\mu$, and the quotient ring
  $R_\mu/\fq$ is a non-trivial valuation ring of the residue field
  $k_\nu=R_\nu/\fm_\nu$. Conversely, a valuation $\kappa$ on the
  residue field $k_\nu$ determines a unique valuation $\mu'$ on $K$
  with $R_{\mu'}\subsetneq R_\nu$, and $R_{\mu'}/\fm_\nu\cap
  R_{\mu'}=R_\kappa$. 
  One says that $\mu'$ is the composite valuation of $\nu$ and $\kappa$. 
  When the value group of $\nu$ is isomorphic
  to $\Z$, we have $\mu'(\psi)=(\nu(\psi),\kappa(\overline{\psi}))$,
  where $\overline{\psi}$ is the class of $\psi$ in $k_\nu$ (see
  \eg~\cite[p.554]{Vaquie}).

  Suppose $\nu$ is centered.  By Abhyankar's inequality~\eqref{AbIn},
  $\trdeg(\nu)$ is either 0 or 1.  In the first case, $k_\nu$ is
  isomorphic to $\C$, which admits no nontrivial nonnegative
  valuation.  Hence $\trdeg(\nu)=1$, $\nu\sim\nu_E$ is divisorial and
  the valuation $\kappa$ defined above has rank $1$.  We may hence
  assume $\kappa(K)\subset\R$ and $\mu=(\nu,\kappa)\in\R\times\R$.

  Performing a suitable sequence of blow-ups and lifting the
  situation, one can suppose $\nu=\nu_\fm$, the multiplicity valuation.
  Assume $\nu(y)\ge\nu(x)$, or else switch the roles of $x$ and
  $y$. Since $\mu\ne\nu_\fm$ there exists a unique $\theta\in\C$ for
  which $\mu(y-\theta x)>\mu(x)$.  Then $\mu$ is equivalent to the
  unique valuation sending $x$ to $(1,0)$ and $y-\theta x$ to $(1,1)$.
  This is the exceptional curve valuation attached to the exceptional
  curve obtained by blowing-up the origin once and at the point of
  intersection with the strict transform of $\{y-\theta x=0\}$.

  When $\nu$ is not centered, the prime ideal $R\cap\fm_\nu$ is
  generated by an irreducible element $\phi\in\fm$. Then 
  $\nu\sim\div_C$, where $C=\{\phi=0\}$.
  The residue field $k_\nu$ is isomorphic to $\C((t))$, on 
  which there exists a unique (up to equivalence) 
  nontrivial valuation vanishing on $\C^*$. 
  Therefore $\mu$ is equivalent to the (Krull) curve valuation
  $\nu'_C$.
\end{proof}
\begin{remark}\label{R413}
  Let $\nu$ be a divisorial valuation on $R$, and pick two Krull
  valuations $\mu_1,\mu_2$ such that $R_{\mu_i}\subsetneq R_\nu$ for
  $i=1,2$. Then $R_{\mu_1}=R_{\mu_2}$ iff one can find an element
  $\psi\in\fm_{\mu_1}\cap\fm_{\mu_2}$ such that $\nu(\psi)=0$.

  To see this, first note that the proof implies that
  $R_{\mu_1}=R_{\mu_2}$ iff the valuation rings 
  $S_i\=R_{\mu_i}/(\fm_\nu\cap R_{\mu_i})$
  coincide for $i=1,2$ inside $k_\nu$. But $k_\nu$ is isomorphic to
  $\C(T)$ so the two valuations rings coincide iff their maximal ideals
  intersect. This happens iff there exists 
  $\psi\in\fm_{\mu_1}\cap\fm_{\mu_2}$ such that $\nu(\psi)=0$.
\end{remark}
%
%
%
%
\section{Sequences of blow-ups and Krull valuations}
\label{sec-equivalence}
So far we have defined (Krull) valuations as functions on the ring $R$
satisfying certain conditions.  Let us now describe the geometric
point of view of valuations developed by Zariski
(see~\cite[p.122]{ZS},~\cite{Vaquie}).  This point of view will be
exploited systematically in Chapter~\ref{A3} as a geometric approach
to the valuative tree.

We start by a naive approach.
Let $\nu\in\tcV_K$ be a centered Krull valuation on $R$.
Fix local coordinates $(x,y)$ so that $R=\C[[x,y]]$.
Assume $\nu(y)\ge\nu(x)$.
Then $\nu\ge0$ on the larger ring $R[y/x]\subset K$ 
and $Z_\nu\=\{\phi\in R[y/x]\ ;\ \nu(\phi)>0\}$ 
is a prime ideal containing $x$.

When $Z_\nu=(x)$, the value $\nu(\phi)$ for $\phi\in R[y/x]$
is proportional to the order of vanishing of $\phi$ along $x$.
Thus $\nu=\nu_\fm$ is the multiplicity valuation.
Otherwise $Z_\nu$ is of the form $Z_\nu=(y/x-\theta,x)$ 
where $\theta\in\C$.
Define $x_1=x$ and $y_1=y/x-\theta$.
Then $\nu$ is a centered Krull valuation on 
the local ring $\C[[x_1,y_1]]\supset R$.
We have assumed $\nu(y)\ge\nu(x)$, but when 
$\nu(y)<\nu(x)$ we get the same conclusion
by exchanging the roles of $x$ and $y$. 
The idea is now to iterate the procedure and
organize the resulting information. 

A more geometric formulation of the construction above
requires some terminology and is usually formulated
in the language of schemes. 
If $(R_1,\fm_1)$, $(R_2,\fm_2)$ are local rings with
common fraction field $K$, then we say that $(R_1,\fm_1)$
\emph{dominates} $(R_2,\fm_2)$ if $R_1\supset R_2$ and
$\fm_2=R_2\cap\fm_1$.
\begin{definition}
  Let $X$ be a (complete) scheme whose field of rational functions 
  is $K$ and let $\nu$ be a Krull valuation on $K$ with valuation
  ring $R_\nu$.
  The \emph{center}\index{valuation!center of} 
  of $\nu$ in $X$ is the 
  unique (not necessarily closed) point $x\in X$ whose local ring
  $\cO_{X,x}\subset K$ is dominated by $R_\nu$.
\end{definition}
In our setting, the situation is quite concrete. 
The space $X$ is the total space of a finite
composition $\pi:X\to(\C^2,0)$ of (point) blowups above the origin.
There are then two cases:
\begin{itemize}
\item[(i)]
  the center of $\nu$ in $X$ is an irreducible component $E$ 
  of the exceptional divisor $\pi^{-1}(0)$;
  in this case $\nu$ is a divisorial valuation as in
  Section~\ref{S28}: $\nu(\phi)$ is proportional to the order of
  vanishing of $\pi^*\phi$ along $E$ for any $\phi\in R$;
\item[(ii)]
  the center of $\nu$ in $X$ is a (closed) point $p$ on the 
  exceptional divisor $\pi^{-1}(0)$; in this case the point $p$
  is characterized by the property that there exists
  a centered valuation $\mu$ on the local ring $\cO_p$ at $p$
  such that $\nu=\pi_*\mu$;
  in other words, there exist local coordinates $(z,w)$ at 
  $p$ and a valuation $\mu$ on the ring $\C[[z,w]]$ such that
  $\mu(z),\mu(w)>0$ and $\nu(\phi)=\mu(\pi^*\phi)$ for 
  every $\phi\in R$.
\end{itemize}
We shall take this more concrete point of view in the sequel.

Again consider a centered Krull valuation $\nu\in\tcV_K$. 
The center of $\nu$ in $(\C^2,0)$ is the origin, which we
denote by $p_0$ 
(in the language of schemes, the center of $\nu$ in $\Spec R$ 
is the maximal ideal $\fm$).
Blow-up the origin: $\tpi_0:X_0\to(\C^2,0)$, 
and let $E_0=\tpi_0^{-1}(0)$ be the exceptional divisor.
Consider the center of $\nu$ in $X_1$.  
Either this is the exceptional component $E_0$, 
in which case $\nu$ is
equivalent to the multiplicity valuation, $\nu\sim\nu_\fm$, 
or the center is a unique point $p_1\in E_0$. 

In the latter case, $\nu$ restricts to a centered Krull valuation
on the local ring $\cO_{p_1}$ at $p_1$ (if we were to follow~(ii),
then we should write $\nu=\pi_*\mu$, where $\mu$ is a 
centered valuation on $\cO_{p_1}$, but we shall not do this here).
Concretely, $\cO_{p_1}$ is of the form $\C[[x_1,y_1]]$ above.
The fraction field of $\cO_{p_1}$ is (isomorphic to) $K$, hence we 
can repeat the procedure.
That is, we blow-up $p_1$: $\tpi_1:X_1\to X_0$,
and consider the center of $\nu$ in $X_1$. 
This is either the exceptional divisor $E_1=\tpi_1^{-1}(p_1)$, 
in which case $\nu\sim\nu_{E_1}$ is divisorial, 
or a unique closed point $p_2\in E_1$.
Iterating this, we obtain a sequence 
$\Pi[\nu]=(p_j)_0^n$, $0\le n\le\infty$, of 
\emph{infinitely nearby points}\index{infinitely nearby points}
above the origin $p_0$.

Conversely, let $\bp=(p_j)_0^n$, $0\le n\le\infty$, be a sequence of
infinitely nearby points above the origin.  By definition, this means
$p_0$ is the origin ;  $\tpi_0:X_0\to(\C^2,0)$ is the blowup of the
origin with exceptional divisor $E_0$, and inductively, for $j<n$,
$\tpi_{j+1}:X_{j+1}\to X_j$ is the blowup of a point $p_{j+1}$ on the
exceptional divisor $E_j=\tpi_j^{-1}(p_j)$ of $\tpi$.  Write
$\pi_j\=\tpi_0\circ\dots\circ\tpi_j$.

Let us show how to associate a Krull valuation to $\bp$.
In fact, rather than defining a Krull valuation directly, 
we shall define a valuation ring $R_\bp\subset K$.

When $n<\infty$, we declare $\phi\in R_\bp$ iff 
$\pi_n^*\phi$ is a regular function at a generic point on $E_n$. 
If $n=\infty$, then $\phi\in R_\bp$ iff 
there exists $j\ge1$ such that 
$\pi_j^*\phi$ is regular at $p_{j+1}$ 
(which implies that
$\pi_k^*\phi$ is regular at $p_{k+1}$ for $k\ge j$).
A direct check shows that $R_\bp$ is a ring,
and since any curve can be desingularized by a finite sequence of
blow-ups, we get either $\phi\in R_\bp$ or $\phi^{-1}\in R_\bp$ 
for any $\phi\in K^*$. 
Thus $R_\bp$ is a valuation ring, and we denote 
by $\val[\bp]$ its associated Krull valuation.

It is clear that the map $\bp\to\val[\bp]$ is injective.  
Let us indicate how to show that $R_{\Pi[\nu]}=R_\nu$. 
This will prove that
the maps $\bp\to\val[\bp]$ and $\nu\to\Pi[\nu]$ are inverse one to
another. For sake of simplicity we suppose that $\Pi[\nu]=(p_j)$ is
an infinite sequence. Pick $\phi\in R_\nu$, \ie
$\nu(\phi)\ge0$. Write $\phi=\phi_1/\phi_2$ with $\phi_\e$ regular
functions. For $j$ large enough, the strict transform of 
$\{\phi_\e=0\}$, $\e=1,2$,
is smooth at $p_{j+1}$ and transverse to the exceptional
divisor, or does not contain $p_{j+1}$. 
In local coordinates $(z,w)$ at $p_{j+1}$, we therefore
have $\pi_j^*\phi=z^a w^b\times\,\text{unit}$ for large $j$.
When $a,b\ge0$, $\pi_j^*\phi$ is regular at 
$p_{j+1}$ and we are done. 
When $a,b \le0$ and $ab\ne0$, $\phi^{-1}$ 
lies in the maximal ideal of the ring $\cO_{p_j}$. 
By construction, this gives $\nu(\phi)=-\nu(\phi^{-1})<0$, 
a contradiction. 
Finally when $a>0$ and $b<0$, one goes one step further
considering the lift of $\phi$ to $p_{j+2}$. 
The coefficients $(a,b)$ are replaced by $(a,a+b)$ or $(b,a+b)$.  
After finitely many steps, we
are hence in one of the two preceding situations. 
This proves that $\pi_j^*\phi$ is regular at $p_{j+1}$ 
for $j$ sufficiently large, hence $\phi\in R_{\Pi[\nu]}$.
Conversely, pick $\phi\in R_{\Pi[\nu]}$. 
For $j$ large enough, $\pi_j^*\phi$ is a regular
function at $p_{j+1}$. 
By definition, the valuation ring associated to
$\nu$ at $p_{j+1}$ contains $\cO_{p_{j+1}}$.
Hence $\phi\in R_\nu$, which concludes the proof that
$R_{\Pi[\nu]}=R_\nu$.

Summing up, we have
\begin{theorem}\label{blw}
  The set $\cV_K$ of centered Krull valuations modulo equivalence 
  is in bijection with the set of sequences of infinitely nearby points.
  The bijection is given as follows: $\nu\in\cV_K$
  corresponds to $(p_j)$ iff $R_\nu$ dominates $\cO_{p_j}$
  for all $j$.
\end{theorem}  
\begin{remark}\label{residue}
  The correspondence above is a particular case of a general result:
  see~\cite[p.122]{ZS}. Moreover, the residue field 
  $k_\nu\cong R_\nu/\fm_\nu$ of a valuation is isomorphic to the 
  field of rational
  functions of the last exceptional divisor of the sequence
  $\Pi[\nu]$. In our case, either $\nu$ is non-divisorial and 
  $k_\nu\cong\C$; or $\nu$ is divisorial and $k_\nu \cong \C(T)$ 
  for some independent variable $T$. 
  \end{remark}
Theorem~\ref{blw} yields the following characterization of divisorial
valuations.
\begin{proposition}\label{Pchardiv}
  For any centered Krull valuation $\nu\in\cV_K$, 
  the following three conditions are equivalent:
  \begin{itemize}
  \item
    $\nu$ is divisorial in the sense
    of Section~\ref{S28};
  \item
    its sequence $\Pi[\nu]$ of infinitely nearby points is finite; 
  \item
    $\trdeg\nu=1$.
  \end{itemize}
\end{proposition}
\begin{proof}
  It is clear that if the sequence of infinitely nearby points
  associated to a Krull valuation is finite, then the valuation is
  divisorial.  We saw in Section~\ref{S28} that any divisorial valuation
  has transcendence degree $1$.  Finally, suppose the sequence
  $\Pi[\nu]=(p_j)_0^\infty$ of infinitely nearby points associated to
  $\nu$ is infinite. We follow the notation of the proof of
  Theorem~\ref{blw}.  Pick $\phi\in R_\nu\setminus\fm_\nu$.  Then
  $\pi_j^*\phi$ is a regular function near $p_{j+1}$ for $j$ large
  enough. As $\phi\not\in\fm_\nu$, it is invertible in $R_\nu$, so we
  may suppose that $\pi_j^*\phi^{-1}$ is also regular at $p_{j+1}$ by
  increasing $j$.  This implies that $\phi$ does not vanish at $p_j$,
  so $\phi=\phi(0)$ inside $k_\nu\=R_\nu/\fm_\nu$.  The residue field
  of $\nu$ is thus isomorphic to $\C$, and $\trdeg\nu=0$.
\end{proof}

In several cases, we may interpret geometrically the sequence of
infinitely nearby points associated to a Krull valuation.
\begin{example}\label{ex-curve}
  When $\nu=\nu_C$ is a curve valuation, then
  $\Pi[\nu]=(p_j)_0^\infty$ is the sequence of infinitely nearby
  points associated to $C$. This sequence can be defined as follows:
  $p_0$ is the origin and, inductively, $p_{j+1}$ is the intersection
  point of $E_j$ and the strict transform of $C$ under $\pi_j$.  See
  Section~\ref{S404}, and Table~\ref{table6} 
  in Appendix~\ref{sec-clas} for more
  information on this sequence of infinitely nearby points.
\end{example}
\begin{example}\label{E101}
  If $\nu=\nu_{E,p}$ is an exceptional curve valuation, then 
  $\Pi[\nu]=(p_j)_0^\infty$, where $\Pi[\nu_E]=(p_j)_0^{j_0}$,
  $p_{j_0+1}=p$ and, inductively for $j>j_0$, 
  $p_{j+1}$ is the intersection point
  of $E_j$ and the strict transform of $E=E_{j_0}$.
  See Appendix~\ref{sec-tangent}.
\end{example}
\begin{example}\label{ex-mono}
  Let $s\ge 1$, and $\nu=\nu_s$ the monomial valuation with
  $\nu(x)=1$, $\nu(y)=s$. Then $\Pi[\nu]$ is determined by
  the\index{continued fractions} 
  \emph{continued fractions expansion}
  \begin{equation*}
    s=a_1+\frac{1}{a_2+\dots}=[[a_1,a_2,\dots]],
  \end{equation*}
  and in particular $\Pi[\nu]$ is finite iff 
  the continued fractions expansion is finite, which happens
  iff $s$ is rational. See~\cite{spiv} for details.
\end{example}


%
%
%
%
%
%
\chapter{MacLane's method}\label{part2}
We now do the preparatory work for 
obtaining the tree structure on valuation space.
Namely, we show how to represent a valuation conveniently
by a finite or countable sequence of polynomials and real numbers.
In the next chapter we shall see how to visualize this encoding
through the language of trees.

Our approach is an adaptation of MacLane's method~\cite{mac}
and can be outlined as follows.
Fix local coordinates $(x,y)$.
Pick a centered valuation $\nu$ on $R=\C[[x,y]]$.
It corresponds uniquely to a valuation, still denoted $\nu$,
on the Euclidean domain $\C(x)[y]$. 
Write $U_0=x$ and $U_1=y$. 
The value $\btilde_0=\nu(U_0)$ determines $\nu$ on $\C(x)$
and the value $\btilde_1=\nu(U_1)$ determines $\nu$ on many 
polynomials in $y$. The idea is now to find 
a polynomial $U_2(x,y)$ of minimal degree in $y$ 
such that $\btilde_2=\nu(U_2)$ is not determined by 
the previous data. Inductively we construct polynomials $U_j$ 
and numbers $\btilde_j=\nu(U_j)$ which, as it turns out,
represent $\nu$ completely.

More precisely we proceed as follows. First we define 
what data $[(U_j);(\btilde_j)]$ to use, namely
sequences of  key polynomials, or SKP's.
These are introduced in Section~\ref{sec-stkp} where we also
show how to associate a canonical valuation to an SKP\@.
In Section~\ref{sec-stru} we study such a valuation in detail,
describing its graded ring and computing its numerical invariants. 
Then in Section~\ref{sec-appro} we show that 
in fact any valuation is associated to an SKP\@. 
Finally, in Section~\ref{some-compu} we compute $\nu(\phi)$ 
for $\phi\in\fm$ irreducible in terms of the
SKP's associated to $\nu$ and the curve valuation $\nu_\phi$. 
This computation is the key to parameterization of 
valuation space (by skewness) as described 
in Section~\ref{metric-tree-struc}.

We formulate almost all results for valuations rather than 
Krull valuations, but the method works for both classes.
We leave it to the reader to make the suitable adaptation in the
notation and statements.

This chapter is quite technical in nature. A number of proofs require a
fine combinatorial analysis and elaborate inductions.
We shall later exploit two other approaches to the valuative tree:
one based on Puiseux expansions (Chapter~\ref{sec-puis}) and one based
on the geometrical interpretation of valuations (Chapter~\ref{A3}).
The method by SKP's has, however, two advantages.
First, it is purely algebraic and may thus be generalized to
more general contexts, especially when $\C$ is replaced by a
field of positive characteristic.
Second, it is well adapted to global problems, for
instance when analyzing singularities of plane curves at infinity in
$\C^2$. These two points of view
already appear in the work of Abhyankar-Moh
on approximate roots, which are very closely related to SKP's.

%
%
%
%
\section{Sequences of  key polynomials}\label{sec-stkp}
%
%
\subsection{Key polynomials}
Let us define the data we use to represent valuations.
Fix local coordinates $(x,y)$ for the rest of this chapter.
\begin{definition}\label{def-key}
  A sequence of polynomials $(U_j)_{j=0}^k$, $1\le k\le\infty$, 
  in $\C[x,y]$ is called a 
  \emph{sequence of key-polynomials}\index{SKP (sequence of key polynomials)} 
  (SKP) if
  it satisfies:
  \begin{enumerate}
  \item[(P0)]
    $U_0=x$ and $U_1=y$;
  \item[(P1)] 
    to each polynomial $U_j$ is attached a number $\btilde_j\in\Rbar$ 
    (not all $\infty$) with
    \begin{equation}\label{key2}
      \btilde_{j+1}>n_j\btilde_j=\sum_{l=0}^{j-1}m_{j,l}\btilde_l
      \quad\text{for $1\le j<k$},
    \end{equation}
    where $n_j\in\N^*$ and $m_{j,l}\in\N$ satisfy,
    for $j>1$ and $1\le l<j$, 
    \begin{equation}\label{key3}
      n_j
      =\min\{l\in\Z\ ;\ l\,\btilde_j\in\Z\btilde_0+\dots+\Z\btilde_{j-1}\}
      \qand
      0\le m_{j,l}<n_l;
    \end{equation}
  \item[(P2)]
    for $1\le j<k$ there exists $\theta_j\in\C^*$ such that
    \begin{equation}\label{e501}
      U_{j+1}
      =U_j^{n_j}-\theta_j\cdot U_0^{m_{j,0}}\cdots U_{j-1}^{m_{j,j-1}}
    \end{equation}
  \end{enumerate}
\end{definition}
In the sequel, we call $k$ the \emph{length}\index{length of an SKP}
of the SKP.
\begin{remark}\label{R412}
  If an abstract semi-group $\Gamma$ is given, a sequence of polynomials
  satisfying (P0)-(P2) with $\btilde_j\in\Gamma$ will
  be called a $\Gamma$-SKP\@.
\end{remark}
\begin{remark}
  If $(U_j)_0^k$ is an SKP of length $k\ge 2$,
  then $(\btilde_j)_1^{k-1}$ are determined by
  $\btilde_0$ and $(U_j)_0^k$.
  Conversely, given a sequence $(\btilde_j)_1^{k+1}$ satisfying (P2) 
  and any sequence $(\theta_j)_1^k$ in $\C^*$ 
  there exists a unique associated SKP\@.
  Despite this redundancy we will typically write an SKP 
  as $[(U_j);(\btilde_j)]$.
\end{remark}
\begin{lemma}\label{lem:deg-irr}
  For $1\le j\le k$, the polynomial $U_j$ is irreducible and of
  Weierstrass form $U_j=y^{d_j}+a_1(x)y^{d_j-1}+\dots+a_{d_j}(x)$ with
  $a_l(0)=0$ for all $l$.  Moreover, $d_{j+1}=n_jd_j$ for $1\le j<k$.
\end{lemma}
\begin{proof}
  The proof of all assertions is by induction. The fact that $U_j$ is
  irreducible is not obvious and will follow from the proof of 
  Theorem~\ref{key}. We included this fact here for clarity.  
  Let us show by induction on $k$ that 
  $d_{j+1}=n_jd_j$ for $1\le j<k$.

  By~\eqref{key3} and the induction hypothesis, we have 
  $m_{j,l}<n_l=d_{l+1}/d_l$, hence $m_{j,l}+1\le d_{l+1}/d_l$ 
  as both sides are integers. We infer
  \begin{equation*}
    \sum_{l=1}^{l-1}m_{j,l}d_l
    \le\sum_{l=1}^{j-1}(\frac{d_{l+1}}{d_l}-1)d_l
    =d_j-1
    <n_jd_j,
  \end{equation*}
  hence $\deg_y(U_{j+1})=n_jd_j$.
  \end{proof}
Although we shall not need the following lemma in this chapter, we 
include it here for convenience. Recall that $m(\phi)$ denotes the
multiplicity of $\phi\in R$.
\begin{lemma}\label{lem:mult-irr}
  When $\btilde_1\ge \btilde_0$, $m(U_j)=d_j$ for $j>0$.
\end{lemma}
\begin{proof}
  The proof goes by induction on $j$. For $j=1$ the assertion
  is obvious as $U_1=y$ and $d_1=1$. 
  For $j=2$, $U_2=y^{n_1}-x^{m_{1,0}}$ 
  with $n_1\btilde_1=m_{1,0}\btilde_0$. 
  As $\btilde_1\ge\btilde_0$, $m(U_2)=n_1=d_1$. 

  Now assume $j\ge3$. 
  To prove $m(U_j)=n_jm(U_{j-1})$ it is sufficient to prove that 
  $n_jd_j<\sum_0^{j-1}m_{j,l}d_l$.
  We have $d_j = n_{j-1}d_{j-1}$ and $\btilde_j>n_{j-1}\btilde_{j-1}$
  hence $n_jd_j=n_jn_{j-1}d_{j-1}<n_j\btilde_jd_{j-1}/\btilde_{j-1}$.  
  Now $n_j\btilde_j=\sum_0^{j-1}m_{j,l}\btilde_l$ and 
  $\btilde_{j-1}>(n_{j-2}\dots n_l)\btilde_l$ for any $l<j$.
  Hence
  \begin{equation*}
    n_jd_j
    <\sum_{l=0}^{j-1}m_{j,l}\btilde_l
    \frac{d_{j-1}}{\btilde_{j-1}}
    < \sum_{l=0}^{j-1} m_{j,l}d_l.
  \end{equation*}
  This concludes the proof.
\end{proof}
\begin{remark}
  In the notation of Zariski and Spivakovsky the sequence
  $(\overline{\beta}_i)$ can be extracted from $(\btilde_j)$ as
  follows. 
  One has $\overline{\beta}_i=\btilde_{k_i}$ 
  where $k_i$ is defined inductively to be the smallest
  integer $k_i>k_{i-1}$ so that $n_{k_i}\ge 2$.
  In the Abhyankar-Moh terminology~\cite{abh-moh}, the $U_{k_i}$
  are the approximate roots of $U_k$.
  See Sections~\ref{sec-approx} and~\ref{S12} and Appendix~\ref{sec-clas}
  for more details.
\end{remark}
We include a result of arithmetic nature which will be used
repeatedly.
\begin{lemma}\label{lem:arith}
  Assume $\btilde_0,\dots,\btilde_{k+1}$ are given so that for all 
  $j=1,\dots,k$ one has $\btilde_{j+1}>n_j\btilde_j$ with
  $n_j\=\min\{n\in\N^*\ ;\ n\btilde_j\in\sum_0^{j-1}\Z\btilde_i\}$.
  Then for any $j=1,\dots,k$ there exists a unique decomposition
  $n_j\btilde_j=\sum_0^{j-1}m_{j,l}\btilde_l$
  where $0\le m_{j,l}<n_l$ for $l=1,\dots,j-1$.
\end{lemma}
\begin{proof}
  By assumption there exist $m_l\in\Z$ with 
  $n_j\btilde_j=\sum_0^{j-1}m_l\btilde_l$. 
  By Euclidean division,
  $m_{j-1}=q_{j-1}n_{j-1}+r_{j-1}$ with $0\le r_{j-1}<n_{j-1}$
  so by the fact that 
  $n_{j-1}\btilde_{j-1}\in\sum_0^{j-2}\Z\btilde_i$ 
  we can suppose that $0\le m_{j-1}<n_{j-1}$.
  Inductively we get $0\le m_l<n_i$ for $1\le l<j$. 
  It remains to prove that $m_0\ge 0$. 
  Set $S=\sum_1^{j-1} m_{l}\btilde_l$.
  We have $m_1\btilde_1<n_1\btilde_1<\btilde_2$ 
  and $m_2<n_2$ so $m_2+1\le n_2$ and 
  \begin{equation*}
    m_1\btilde_1+m_2\btilde_2
    <(1+ m_2)\btilde_2\le n_2\btilde_2
    <\btilde_3.
  \end{equation*}
  We infer that 
  $S<(1+ m_{j-1})\btilde_{j-1}\le n_{j-1}\btilde_{j-1}<\btilde_j$ 
  so that $m_0\btilde_0\ge\btilde_j-S\ge0$.
  \end{proof}
%
%
\subsection{From SKP's to valuations I}\label{S30}
We show how to associate a valuation $\nu$ to any finite SKP
in a canonical way.
\begin{theorem}\label{key}
  Let $[(U_j)_0^k;(\btilde_j)_0^k]$ be an SKP of length 
  $1\le k<\infty$. Then 
  there exists a unique centered valuation $\nu_k\in\tcV$ satisfying
  \begin{enumerate}
  \item[(Q1)]
    $\nu_k(U_j)=\btilde_j$ for $0\le j\le k$;
  \item[(Q2)] $\nu_k\le\nu$ for any $\nu\in\tcV$ satisfying (Q1).
  \end{enumerate}
  Further, if $k>1$ and $\nu_{k-1}$ is the valuation associated to 
  $[(U_j)_0^{k-1};(\btilde_j)_0^{k-1}]$, then
  \begin{enumerate}
  \item[(Q3)]
    $\nu_{k-1}\le\nu_k$;
  \item[(Q4)]
    $\nu_{k-1}(\phi)<\nu_k(\phi)$ iff
    $U_k$ divides $\phi$ in $\gr_{\nu_k}\C(x)[y]$.
  \end{enumerate}
\end{theorem}
\begin{remark}
  The valuation $\nu_k$ is normalized in the sense that $\nu_k(\fm)=1$
  iff $\min\{\btilde_0,\btilde_1\}=1$. 
  It is normalized in the sense that $\nu_k(x)=1$ iff $\btilde_0=1$.
\end{remark}
The proof of Theorem~\ref{key} occupies all of Section~\ref{pftheo}. 
It is a subtle induction on $k$ using divisibility properties 
in the graded ring $\gr_{\nu}\C(x)[y]$.

First, $\nu_1$ is defined to be the 
monomial valuation\index{valuation!monomial}\label{page1}
with $\nu_1(x)=\btilde_0,\nu_1(y)=\btilde_1$, \ie 
\begin{equation}\label{e440}
  \nu_1(\sum a_{ij}x^iy^j)=\min\{i\btilde_0+j\btilde_1\ ;\ a_{ij}\ne0\}.
\end{equation} 
Property~(Q2) clearly holds.

Now assume $k>1$, that the SKP $[(U_j)_0^k;(\btilde_j)_0^k]$ is given
and that $\nu_1,\dots,\nu_{k-1}$ have been defined. 
First consider a polynomial $\phi\in \C [x,y]$. As $U_k$ is unitary
in $y$, we can divide $\phi$ by $U_k$ in $\C[x][y]$:
$\phi=\phi_0+U_k\psi$ with $\deg_y(\phi_0)<d_k=\deg_y(U_k)$ and
$\psi\in \C[x,y]$. Iterating the procedure we get a unique
decomposition
\begin{equation}\label{divincrepower}
  \phi=\sum_j\phi_jU_k^j
\end{equation}
with $\phi_j\in \C[x,y]$ and $\deg_y(\phi_j)<d_k$. 
Define
\begin{equation}\label{min}
  \nu_k(\phi)
  \=\min_j\nu_k(\phi_jU_k^j)
  \=\min_j\{\nu_{k-1}(\phi_j)+j\btilde_k\}.
\end{equation}

The function $\nu_k$ is easily seen to satisfy~(V2) and~(V3). We will
show it also satisfies~(V1), hence defines a valuation on 
$\C[x,y]$, which automatically satisfies~(Q2). 
Theorem~\ref{key} then follows from
\begin{proposition}[\cite{spiv}]\label{extend}
  Any valuation $\nu:\C[x,y]\to\Rbar$ with $\nu(x),\nu(y)>0$
  has a unique extension to a centered valuation on $R=\C[[x,y]]$.
  This extension preserves the partial ordering:
  if $\nu_1(\phi)\le\nu_2(\phi)$ holds for polynomials 
  $\phi$, it also holds for formal power series.
\end{proposition}
\begin{proof}
  We may assume that $\nu$ is normalized in the sense that 
  $\nu(\fm)=1$.
  Pick $\phi\in R$, and write $\phi=\phi_k+\hphi_k$
  where $\phi_k$ is the truncation of $\phi$ at order $k$,
  \ie $m(\hphi_k)\ge k$. 
  Notice that this implies $\nu(\hphi_k)\ge k$.
  
  We claim that $\nu(\phi)\=\lim_k\nu(\phi_k)$ exists in $\Rbar$.
  Indeed, for any $k,l$, we have
  \begin{equation*}
    \nu(\phi_{k+l})
    \ge\min\{\nu(\phi_k),\nu(\phi_{k+l}-\phi_k )\}
    \ge\min\{\nu(\phi_k),k\}
  \end{equation*}
  Thus the sequence $\min\{\nu(\phi_k),k\}$ is 
  nondecreasing, hence converges. 

  Uniqueness can be proved as follows: 
  $\nu(\phi)\ge\min\{\nu(\phi_k),k\}$ for all $k$.  
  If $\nu(\phi_k)$ is unbounded, then $\nu(\phi)=\infty$. 
  Otherwise $\nu(\phi)=\nu(\phi_k)$ for $k$ large enough.
  By construction the extension
  preserves the partial ordering.
\end{proof}
The extension above needs not preserve numerical invariants.
For example, if $\nu$ is a curve valuation at 
a non-algebraic curve, then $\rk(\nu|_{\C[x,y]})=1$ 
but $\rk(\nu)=2$.
\begin{remark}\label{weier}
  The expansion \eqref{divincrepower} can be applied 
  for any $\phi\in R=\C[[x,y]]$. Indeed, 
  using Weierstrass' preparation theorem, we can assume
  that $\phi$ is a monic polynomial in $y$.
  Weierstrass' division theorem then yields~\eqref{divincrepower}
  with $\phi_j\in\C[[x]][y]$. The valuation $\nu_k$ can
  hence be defined by~\eqref{min} for any formal power series.
\end{remark}
%
%
\subsection{Proof of Theorem~\ref{key}}\label{pftheo}
By Proposition~\ref{extend}, we will restrict our attention to
centered valuations defined on the Euclidean ring $\C(x)[y]$.  
We first make the induction hypothesis precise.
\begin{itemize}
\item[$(H_k)$:]
  $\nu_k$ is a valuation satisfying (Q1)-(Q4).
\item[$(E_k)$:]
  the graded ring  $\gr_{\nu_k}\C(x)[y]$ is a Euclidean domain.
\item[$(I_k)$:]
  $U_k$ and $U_{k+1}$ are irreducible in $\gr_{\nu_k}\C(x)[y]$.
\item[$(\overline{I}_k)$:] 
  $U_j$ is irreducible both in $\C[x,y]$, 
  and $\C(x)[y]$ for $1\le j\le k$.
\end{itemize}
\medskip
For the precise definition of a Euclidean domain 
we refer to~\cite[p.23]{ZS0}.
One immediately checks that $H_1$ and $\overline{I}_1$ hold.
Our strategy is to show successively 
$\overline{I}_k\ \&\ H_k\Rightarrow E_k$; 
$\overline{I}_k\ \&\ H_k \Rightarrow I_k$; 
$E_k\ \&\ I_k\ \&\ H_k\Rightarrow H_{k+1}$; 
$I_k\ \&\ \overline{I}_k\ \&\ H_{k+1} \Rightarrow \overline{I}_{k+1}$.

\medskip
\noindent\textbf{Step 1}: $\overline{I}_k\ \&\ H_k\Rightarrow E_k$.
We want to show that $\gr_{\nu_k}\C(x)[y]$ is a Euclidean domain.
Let $\phi\in R$ and consider the expansion~\eqref{divincrepower}.
We define
\begin{equation*}
  \delta_k(\phi)\=\max\left\{j\ ;\ \nu_k(\phi_jU_k^j)=\nu_k(\phi)\right\}. 
\end{equation*}
By convention we let $\delta_k(0)=-\infty$. 
Note that if $\phi=\phi'$ modulo $\nu_k$, then 
$\delta_k(\phi)=\delta_k(\phi')$. 
Hence $\delta_k$ is well defined on $\gr_{\nu_k}\C(x)[y]$.
\begin{lemma}\label{L301}
  For all $\phi,\psi\in\C(x)[y]$, 
  $\delta_k(\phi\psi)=\delta_k(\phi)+\delta_k(\psi)$.
\end{lemma}
\begin{proof}
  First assume $\delta_k(\phi)=\delta_k(\psi)=0$. 
  Then $\phi=\phi_0$, $\psi=\psi_0$ in $\gr_{\nu_k}\C(x)[y]$, 
  so we can assume $\deg_y(\phi),\deg_y(\psi)<d_k$. 
  Thus $\deg_y(\phi\psi)<2d_k$ so
  $\phi\psi=\a_0+\a_1U_k$ with $\deg_y(\a_0),\deg_y(\a_1)<d_k$. 
  Suppose $\nu_k(\a_1U_k)=\nu_k(\phi\psi)$. 
  Then $\nu_{k-1}(\a_1U_k)<\nu_k(\a_1U_k)\le\nu_k(\a_0)=\nu_{k-1}(\a_0)$. 
  Thus
  \begin{align*}
    \nu_{k-1}(\a_1U_k) 
    =\nu_{k-1}(\phi\psi) 
    =\nu_{k-1}(\phi)+\nu_{k-1}(\psi)
    =\nu_k(\phi)+\nu_k(\psi)
    =\nu_k(\phi\psi).
  \end{align*}
  Hence $\nu_{k-1}(\a_1U_k)=\nu_k(\a_1U_k)$ contradicting~(Q4).
  So $\phi\psi=\a_0$ in $\gr_{\nu_k}\C(x)[y]$ and 
  $\delta_k(\phi\psi)=0$.

  In the general case, $\phi=\sum\phi_iU_k^i$, 
  $\psi=\sum\psi_jU_k^j$, we have $\phi\psi=\sum(\phi_i\psi_j)U_k^{i+j}$. 
  By what precedes $\delta_k(\phi_i\psi_j)=0$ for all $i,j$ 
  so $\delta_k(\phi\psi)=\delta_k(\phi)+\delta_k(\psi)$.
  \end{proof}
\begin{lemma}\label{unit}
  For $\phi\in\C(x)[y]$, $\delta_k(\phi)=0$ iff
  $\phi$ is a unit in $\gr_{\nu_k}\C(x)[y]$.
\end{lemma}
\begin{proof}
  If $\delta_k(\phi)=0$, then $\phi=\phi_0$ in $\gr_{\nu_k}\C(x)[y]$. 
  As $U_k$ is irreducible in $\C(x)[y]$ 
  and $\deg_y(U_k)>\deg_y(\phi_0)$, the polynomial 
  $U_k$ is prime with $\phi_0$. Hence we can find 
  $A,B\in\C(x)[y]$ with $\deg_yA,\deg_yB<d_k$ so that 
  $A\phi_0=1-BU_k$.
  Then, comparing with~\eqref{divincrepower}, we have
  $\nu_k(A\phi_0)=\nu_k(1)<\nu_k(BU_k)$. 
  Therefore $A\phi_0=1$ in $\gr_{\nu_k}\C(x)[y]$ so
  $\phi_0$, and hence $\phi$, is a unit in $\gr_{\nu_k}\C(x)[y]$.

  Conversely, if $\phi$ is a unit, say $A\phi=1$ in $\gr_{\nu_k}\C(x)[y]$ 
  for some $A\in\C(x)[y]$, then 
  $\delta_k(\phi)+\delta_k(A)=\delta_k(1)=0$, so $\delta_k(\phi)=0$.
  \end{proof}
\begin{lemma}\label{eucl}
  If $\phi,\psi\in\C(x)[y]$, then there exists $Q,R\in\C(x)[y]$ 
  such that $\phi=Q\psi+R$ in $\gr_{\nu_k}\C(x)[y]$ and
  $\delta_k(R)<\delta_k(\psi)$.
\end{lemma}
\begin{proof}
  Write
  $\psi=\sum_j\psi_jU_k^j$. 
  It suffices to prove the lemma when 
  $\psi_j=0$ for $j>M:=\delta_k(\psi)$
  and using Lemma~\ref{unit} we may assume $\psi_M=1$.
  As $\deg_y(\psi_j)<d_k$ for $j\le M$ we have 
  $\deg_y(\psi)=Md_k$.  
  Euclidean division in $\C(x)[y]$ yields 
  $Q,R^1\in\C(x)[y]$ with $\deg_y(R^1)<\deg_y(\psi)$
  so that $\phi=Q\psi+R^1$. 
  Write $R^1=\sum_iR_iU_k^i$ and set
  $N\=\delta_k(R^1)$, $R\=\sum_{i\le N}R_iU_k^i$.
  Then $\phi=Q\psi+R$ in $\gr_{\nu_k}\C(x)[y]$ and
  \begin{equation*}
    \deg_y(R)
    =\deg_y(R_N)+Nd_k
    <Md_k=\deg_y(\psi).
  \end{equation*}
  Hence $N<M$ and we are done.
  \end{proof}
This completes Step~1. The Euclidean property of 
$\gr_{\nu_k}\C(x)[y]$ makes
every ideal principal and supplies us with 
unique factorization and Gauss' lemma.

\medskip
\noindent\textbf{Step 2}: $\overline{I}_k\ \&\ H_k\Rightarrow I_k$.
We want to show that $U_k$ and $U_{k+1}$ are irreducible 
in $\gr_{\nu_k}\C(x)[y]$ and proceed in several steps.
\begin{lemma}
  $U_k$ is irreducible in $\gr_{\nu_k}\C(x)[y]$. 
\end{lemma}
\begin{proof}
  We trivially have $\delta_k(U_k)=1$ so
  if $\phi\psi=U_k$ in $\gr_{\nu_k}\C(x)[y]$, then 
  $\delta_k(\phi)=0$ or $\delta_k(\psi)=0$.
  Hence $\phi$ or $\psi$ is a unit in $\gr_{\nu_k}\C(x)[y]$
  (see Lemma~\ref{L301} and~\ref{unit}).
  \end{proof}
\begin{lemma}\label{L302}
  If $j<k$ then $U_j$ is a unit in $\gr_{\nu_k}\C(x)[y]$.
\end{lemma}
\begin{proof}
  By Lemma~\ref{unit} it suffices to show that $\delta_k(U_j)=0$.
  If $d_j<d_k$ then this is obvious. 
  If $d_j=d_k$, then $U_j=(U_j-U_k)+U_k$, where 
  $\deg_y(U_j-U_k)<d_k$. 
  Now $\nu_k(U_j)=\btilde_j<\btilde_k=\nu_k(U_k)$,
  so $\nu_k(U_j-U_k)<\nu_k(U_k)$ and $\delta_k(U_j)=0$.
  \end{proof}
\begin{lemma}\label{uni}
  If $\delta_k(\phi)< n_k$, then $\phi=\phi_iU_k^i$ 
  in $\gr_{\nu_k}\C(x)[y]$ for some $i<n_k$.
\end{lemma}
\begin{proof}
  Suppose $\nu_k(\phi_iU_k^i)=\nu_k(\phi_jU_k^j)=\nu_k(\phi)$,
  where $i\le j<n_k$. 
  Then 
  $(j-i)\btilde_k=\nu_{k-1}(\phi_i)-\nu_{k-1}(\phi_j)\in\sum_0^{k-1}\Z\btilde_j$.
  By~\eqref{key3} $n_k$ divides $j-i$ hence $i=j$.
  \end{proof}
\begin{lemma}\label{L303}
  $U_{k+1}$ is irreducible in $\gr_{\nu_k}\C(x)[y]$. 
\end{lemma}  
\begin{proof}
  We have $U_{k+1}=U_k^{n_k}-\tilde{U}_{k+1}$, where
  $\tilde{U}_{k+1}=\theta_k\prod_{j=0}^{k-1}U_j^{m_{k,j}}$.
  Assume $U_{k+1}=\phi\psi$ in $\gr_{\nu_k}\C(x)[y]$ with 
  $0<\delta_k(\phi),\delta_k(\psi)<n_k$. 
  By Lemma~\ref{uni}, we can write $\phi=\phi_iU_k^i$, 
  $\psi=\psi_jU_k^j$.
  Then $U_{k+1}=\phi_i\psi_jU_k^{n_k}$ so
  $(1-\phi_i\psi_j)U_k^{n_k}=\tilde{U}_{k+1}$. 
  As $U_k$ is irreducible and 
  $\tilde{U}_{k+1}$ is a unit, we have 
  $\phi_i\psi_j=1$ in $\gr_{\nu_k}\C(x)[y]$. 
  But then $\tilde{U}_{k+1}=0$
  in $\gr_{\nu_k}\C(x)[y]$, which is absurd.
  So we can assume $\delta_k(\phi)=n_k$ and $\delta_k(\psi)=0$. 
  Hence $\psi$ is a unit, completing the proof.
  \end{proof}

\medskip
\noindent\textbf{Step 3}: $E_k\ \&\ I_k\ \&\ H_k \Rightarrow H_{k+1}$.
We must show
$\nu_{k+1}(\phi\psi)=\nu_{k+1}(\phi)+\nu_{k+1}(\psi)$
for $\phi,\psi \in \C[x,y]$ where $\nu_{k+1}$ is defined 
as in~\eqref{min}. We first note
\begin{itemize}
\item[(i)]
  $\deg_y(\phi)<d_{k+1}$ implies $\nu_{k+1}(\phi)=\nu_k(\phi)$;
\item[(ii)]
  $\nu_k(U_{k+1})=n_k\btilde_k<\btilde_{k+1}=\nu_{k+1}(U_{k+1})$;
\item[(iii)]
  $\nu_{k+1}(U_{k+1}\phi)=\nu_{k+1}(U_{k+1})+\nu_{k+1}(\phi)$ 
  for all $\phi\in \C[x,y]$.
\end{itemize}
The second assertion follows from \eqref{key2}.
We now show
\begin{lemma}\label{st1}
  For $\phi\in \C[x,y]$ we have
  $\nu_{k+1}(\phi)\ge\nu_k(\phi)$. 
\end{lemma}
\begin{proof}
  We argue by induction on $\deg_y(\phi)$.  If $\deg_y(\phi)<d_{k+1}$,
  then we are done by~(i).  If $\deg_y(\phi)\ge d_{k+1}$, then we
  write $\phi=\sum_i\phi_iU_{k+1}^i$ as in~\eqref{divincrepower}.
  By the induction hypothesis and~(ii)-(iii) above we may assume that
  $\phi_0\not\equiv0$.  Write $\phi=\phi_0+\psi$ with 
  $\psi=\sum_{i\ge1}\phi_iU_{k+1}^i$. 
  When $\nu_k(\phi)=\min\{\nu_k(\phi_0),\nu_k(\psi)\}$, one has 
  \begin{equation*}
    \nu_{k+1}(\phi)
    =\min\{\nu_{k+1}(\phi_0),\nu_{k+1}(\psi)\}
    \ge\nu_k(\phi),
  \end{equation*}
  proving the lemma in this case. Otherwise, $\phi_0+\psi=0$ in
  $\gr_{\nu_k}\C(x)[y]$. This implies that $U_{k+1}$ divides $\phi_0$
  in this ring. But $\deg_y(\phi_0)<d_{k+1}$, hence 
  $\delta_k(\phi_0)<n_k$, so that $\phi_0$ is a power of 
  the irreducible $U_k$ times a unit in $\gr_{\nu_k}\C(x)[y]$ 
  by Lemma~\ref{uni}.  By Lemma~\ref{L303}, $U_{k+1}$ 
  is also irreducible, a contradiction.
  \end{proof}
Introduce $\fp\=\{\nu_{k+1}>\nu_k\}\subset\gr_{\nu_k}\C(x)[y]$. 
Using the preceding lemma, one easily verifies that 
$\fp$ is a proper ideal, which contains 
the irreducible element $U_{k+1}$ by~(ii). 
As $\gr_{\nu_k}\C(x)[y]$ is a Euclidean domain, $\fp$
is generated by $U_{k+1}$.

Fix $\phi,\psi\in R$. We want to show that
$\nu_{k+1}(\phi\psi)=\nu_{k+1}(\phi)+\nu_{k+1}(\psi)$.
First assume $\phi,\psi\not\in\fp$.
Then $\phi\psi\not\in\fp$. 
By $H_k$, $\nu_k$ is a valuation, so
\begin{equation*}
  \nu_{k+1}(\phi\psi)
  =\nu_{k}(\phi\psi)
  =\nu_k(\phi)+\nu_k(\psi)
  =\nu_{k+1}(\phi)+\nu_{k+1}(\psi).
\end{equation*}
In the general case, write $\phi=\widehat{\phi}\,U_{k+1}^n$,
$\psi=\widehat{\psi}\,U_{k+1}^m$ in $\gr_{\nu_k}\C(x)[y]$ where
$\widehat{\phi}$, $\widehat{\psi}$ are prime with $U_{k+1}$.
In particular they do not belong to $\fp$ so that
\begin{align*}
  \nu_{k+1}(\phi\psi) 
  &=\nu_{k+1}(\widehat{\phi}\,\widehat{\psi}\,U_{k+1}^{m+n})\\
  &=\nu_{k+1}(\widehat{\phi})+\nu_{k+1}(\widehat{\psi}) 
  +(m+n)\nu_{k+1} (U_{k+1}) 
  =\nu_{k+1}(\phi)+\nu_{k+1}(\psi).
\end{align*}
This completes Step 3.

\medskip
\noindent\textbf{Step 4}: 
$I_k\ \&\ \overline{I}_k\ \&\ H_{k+1}\Rightarrow\overline{I}_{k+1}$.
What we have to prove is 
\begin{lemma}\label{irr-k+1}
  $U_{k+1}$ is irreducible in 
  $\C[x,y]$ and $\C(x)[y]$.
\end{lemma}
\begin{proof}
  As $U_{k+1}$ is monic in $y$, irreducibility in $\C[x,y]$
  implies irreducibility in $\C(x)[y]$. 
  So suppose $\phi\phi'=U_{k+1}$ for $\phi,\phi'\in\C[x,y]$. 
  As $U_{k+1}$ is irreducible in $\gr_{\nu_k}\C(x)[y]$,
  we may assume that the image of $\phi$ in $\gr_{\nu_k}\C(x)[y]$ 
  is a unit \ie $\delta_k(\phi)=0$ and $\delta_k(\phi')=n_k$.
  On the other hand, one has 
  $\delta_k(\psi)d_k\le\deg_y(\psi)$ 
  for any $\psi\in\C[x,y]$. Using Lemma~\ref{lem:deg-irr} 
  (except the assertion of irreducibility which we want now to prove), 
  we infer
  \begin{equation*}
    d_k\delta_k(\phi')
    =d_kn_k
   =\deg_y(U_{k+1})
    =\deg_y(\phi)+\deg_y(\phi')
    \ge\deg_y(\phi')
    \ge d_k\delta_k(\phi'),
  \end{equation*}
  hence $\deg_y(\phi)=0$ and $\deg_y(\phi')=\deg_y(U_{k+1})$.
  Now $U_{k+1}$ is unitary in $y$ so $\phi\not\in\fm$ and is 
  hence a unit in $R$.
  \end{proof}
This completes the proof of Theorem~\ref{key}.
\begin{remark}\label{formal-irr}
  Since $U_{k+1}$ is monic in $y$ we may use
  Weierstrass' division theorem and 
  the argument above to show that
  $U_{k+1}$ is also irreducible in $R=\C[[x,y]]$.
\end{remark}
%
%
\subsection{From SKP's to valuations II}
We now turn to infinite SKP's.
\begin{theorem}\label{key-infinite}
  Let $[(U_j);(\btilde_j)]$ be an infinite SKP
  and let $\nu_k$ the valuation associated to 
  $[(U_j)_0^k;(\btilde_j)_0^k]$ for $k\ge 1$ 
  by Theorem~\ref{key}.
  \begin{itemize}
  \item[(i)]
    If $n_j\ge 2$ for infinitely many $j$, then for any 
    $\phi\in R$ there exists $k_0=k_0(\phi)$ such that 
    $\nu_k(\phi)=\nu_{k_0}(\phi)$ for all $k\ge k_0$. 
    In particular, $\nu_k$ converges to a valuation $\nu_\infty$.
  \item[(ii)]
    If $n_j=1$ for $j\gg1$, then 
    $U_k$ converges in $R$ to 
    an irreducible formal power series $U_\infty$
    and $\nu_k$ converges to a valuation $\nu_{\infty}$. 
    More precisely, for $\phi\in R$ prime to $U_{\infty}$
    we have $\nu_k(\phi)=\nu_{k_0}(\phi)$ for $k\ge k_0=k_0(\phi)$,
    and if $U_{\infty}$ divides $\phi$, then $\nu_k(\phi)\to\infty$.
  \end{itemize}
\end{theorem}
\begin{proof}[Proof of Theorem~\ref{key-infinite}]
  If $n_j\ge 2$ for infinitely many $j$, then
  $\deg_y(U_{k+1})=d_{k+1}=\prod_1^kn_j$ tends to infinity.  Pick
  $\phi\in R$. By Weierstrass' division theorem we may assume
  $\deg_y(\phi)<\infty$. By Remark~\ref{weier}, $\nu_k(\phi)$
  is defined using formula~\eqref{divincrepower}.  For $k\gg1$,
  $\deg_y(\phi)<\deg_y(U_{k+1})$, so that
  $\nu(\phi)=\nu(\phi)=\nu_k(\phi)=\nu_{k+l}(\phi)=\nu_{k+l}(\phi)$
  for all $l\ge0$.  Thus $\nu_k$ converges towards a valuation
  $\nu_\infty$.

  If  $n_k=1$ for $k\ge K$, then 
  set $d\=\max_j\deg_y(U_j)$. For $k\ge K$, 
  \begin{equation}\label{n=1}
    U_{k+1}=U_k-\theta_k\prod_0^{k-1}U_j^{m_{k,j}}.
  \end{equation} 
  As $\deg_y(U_{k+1})=\deg_y(U_k)=d$ 
  one has $m_{k,j}=0$ for any $j\ge K$ by~\eqref{key3}. 
  Then $\{\btilde_k\}_{k\ge K}$ is a
  strictly increasing sequence of real numbers belonging to the discrete
  lattice $\sum_0^{K}\Z\btilde_j$, so $\btilde_k\to\infty$.
  Write
  \begin{equation*}
    U_k=y^d+a_{d-1}^k(x)y^{d-1}+\dots+a_0^k(x),
  \end{equation*}
  the Weierstrass form of $U_k$. 
  As $m_{k,j}=0$ for any $j\ge K$ we get
  \begin{equation*}
    a_n^{k+1}(x) 
    =a_n^k(x)-\theta_k\sum_I a_{i_1}^{j_1}(x)\dots a_{i_l}^{j_l}(x)
  \end{equation*}
  where $\#\{j=\a\ ;\ j\in I\}=m_{k,\a}$, and $i_1+ \cdots + i_l =n$.
  Hence $a_n^{k+1}(x)-a_n^k(x)\in\fm^{\sum_0^{K}m_{k,i}}$. 
  But the sequence $\btilde_j$ is increasing so
  \begin{equation*}
    \sum_0^{K}m_{k,i}
    \ge\btilde_{K}^{-1}\sum_0^{K}\btilde_im_{k,i}
    =\frac{\btilde_k}{\btilde_{K}}\to\infty,
  \end{equation*}
  and $U_k$ converges towards a polynomial in $y$ 
  that we denote by $U_\infty$.  
  By~\eqref{n=1} we have $U_{k+1}=U_k$ modulo $\nu_{k-1}$ 
  and $U_{k+l}=U_k$ modulo $\nu_{k-1}$ for all $l \ge 0$. 
  So $U_{\infty}= U_k$ in
  $\nu_{k-1}$.  Therefore 
  $\nu_{\infty}(U_{\infty})\=
  \lim\nu_k(U_{\infty})=\lim\nu_k(U_k)=\lim\btilde_k=\infty$. 
  Let $\phi\in R$ and write $\phi=\phi_0+U_{\infty}\phi_1$ 
  with $\deg_y\phi_0<\deg_yU_{\infty}=d$. 
  When $\phi_0\equiv0$, then set 
  $\nu_{\infty}(\phi)\=\lim\nu_k(\phi)=\infty$. 
  Otherwise for $k\ge K$ large enough, for all $l\geq 0$, 
  $\nu_{k+l}(\phi)=\nu_{k+l}(\phi_0)=\nu_k(\phi_0)$. 
  Hence the sequence $\nu_k(\phi)$ is stationary for 
  $k\ge K$ and we can set $\nu_{\infty}(\phi)\=\lim\nu_k(\phi)$. 
  One easily checks that $\nu_{\infty}$ is a valuation. 
  As $\nu_\infty(\phi)=\infty$ iff $U_\infty|\phi$
  it follows that $U_{\infty}$ is irreducible in $R$.
\end{proof}
%
%
%
%
\section{Classification}\label{S26}
Having constructed a valuation associated to any SKP, we now
make a preliminary classification.
\begin{definition}\label{D101}
  Consider a centered valuation $\nu$ on $R$ given by
  an SKP: say $\nu\=\val[(U_j)_0^k;(\btilde_j)_0^k]$, where 
  $1\le k\le\infty$. 
  Assume that $\nu$ is normalized in the sense that $\nu(\fm)=1$.
  We then say that $\nu$ is
  \begin{itemize}
  \item[(i)]
    \emph{monomial}\index{valuation!monomial} 
    (in coordinates $(x,y)$)
    if $k=1$, $\btilde_0<\infty$ and $\btilde_1<\infty$;
  \item[(ii)]
    \emph{quasimonomial}\index{valuation!quasimonomial} 
    if $k<\infty$, $\btilde_0<\infty$ and $\btilde_k<\infty$;
  \item[(iii)]
    \emph{divisorial}\index{valuation!divisorial} 
    if $\nu$ is quasimonomial and $\btilde_k\in\Q$;
  \item[(iv)]
    \emph{irrational}\index{valuation!irrational} 
    if $\nu$ is quasimonomial but not divisorial;
  \item[(v)]
    \emph{infinitely singular}\index{valuation!infinitely singular} 
    if $k=\infty$ and $d_j\to\infty$, where $d_j=\deg_y(U_j)$;
  \item[(vi)]
    a \emph{curve valuation}\index{valuation!curve} 
    if $k=\infty$ and $d_j\not\to\infty$, or $k<\infty$ and 
    $\max\{\btilde_0,\btilde_k\}=\infty$.
  \end{itemize}
  If $\nu$ is not normalized, then the type of $\nu$ is defined to
  be the type of $\nu/\nu(\fm)$.
\end{definition}
Two remarks are in order.
First, this classification may seem unnatural as it
only applies to valuations associated to an SKP 
and since the technique of SKP's uses a fixed choice of local 
coordinates $(x,y)$.

However, as we shall see, every valuation is associated to an
SKP (Theorem~\ref{approx}). Moreover, the classification can
be rephrased in several equivalent ways, all of which are independent
on the choice of coordinates.

Second, we should compare the classification with the ``road map'' 
given in Section~\ref{ex}. The comparison goes as follows.

\begin{itemize}
\item[(i)]
  If $\nu$ is monomial in coordinates $(x,y)$ in the sense above, 
  then its SKP is of length 1. 
  By~\eqref{e440},
  this means that $\nu$ is monomial in the sense of 
  Section~\ref{S31}.
  Notice that by the definition above, any monomial valuation is also
  quasimonomial.
\item[(ii)]
  It is not obvious that the definition of quasimonomial, which is
  based purely on the definition of a valuation as a function on $R$,
  coincides with the more geometric point of view presented in
  Section~\ref{S32}. We shall prove in Chapter~\ref{A3} that the 
  two definitions are equivalent: see Proposition~\ref{monomialization}
  for a precise statement.
\item[(iii)] 
  We shall see at Theorem~\ref{divis} that $\nu$ is a
  divisorial valuation in the sense above, iff $\trdeg \nu =1$. 
  It follows from Proposition~\ref{Pchardiv} that the definition of
  divisorial valuation given above coincides with the one in 
  Section~\ref{S28}.
\item[(v)]
  As for infinitely singular valuations, we need to develop more
  machinery to see that the definition above agrees with the
  characterizations mentioned in Section~\ref{ex}.
  Informally speaking, however, the sequence $(U_k)_1^\infty$ 
  of key polynomials of increasing and unbounded 
  degrees correspond to the sequence 
  of truncated Puiseux expansions $\sum_1^kx^{\hbeta_j}$ where the 
  $\hbeta_j$'s are rational numbers with unbounded denominators.
\item[(vi)] Finally, consider a curve valuation $\nu$ in the sense
above.  First assume it is defined by a finite SKP, say
$\nu\=\val[(U_j)_0^k;(\btilde_j)_0^k]$, where $k<\infty$ and
$\btilde_k=\infty$. By construction, $\nu(U_k) = \infty$.  As $U_k$ is
irreducible, we conclude that $\nu=\nu_{U_k}$ (see the discussion
following Lemma~\ref{L601}). When $k=\infty$ and $\lim \btilde_k =
\infty$, by Theorem~\ref{key-infinite}, $U_k$ converges to an
irreducible formal power series $U_\infty$, and $\nu(U_\infty) =
\infty$. Again $\nu= \nu_{U_\infty}$.
\end{itemize}

%
%
%
%
\section{Graded rings and numerical invariants}\label{sec-stru}
The aim of this section is to give the structure of the valuation
associated to a finite or infinite SKP\@.  That is, we describe the
structure of the graded ring $\gr_\nu\C(x)[y]$, and compute the three
invariants $\rk(\nu)$, $\ratrk(\nu)$ and $\trdeg(\nu)$.  We refer
to~\cite{teissier} for general results in higher dimensions.

Given a finite or infinite SKP we denote by
$\nu\=\val[(U_j);(\btilde_j)]$
its associated valuation through 
Theorem~\ref{key} or~\ref{key-infinite}.
%
%
\subsection{Homogeneous decomposition I}
The following theorem gives the structure of the graded ring
of a valuation defined by a finite SKP.
\begin{theorem}\label{homo}
  Let $\nu\=\val[(U_j)_0^k;(\btilde_j)_0^k]$, where $1\le k<\infty$,
  define 
  $n_k\=\min\{n>0\ ;\ n\btilde_k\in\sum_0^{k-1}\Z\btilde_j\}$,
  and pick $0\ne\phi\in R$.
  \begin{itemize}
  \item[(i)]
    If $n_k=\infty$ or equivalently $\btilde_k\not\in\Q\btilde_0$, then
    \begin{equation} \label{equa1}
      \phi=\a\prod_{j=0}^{k}U_j^{i_j}\ \text{in}\ \gr_\nu R
      \ \text{and}\ \gr_\nu R_\nu
    \end{equation}
    with $\a\in\C^*$, $0\le i_j<n_j$ for $1\le j<k$
    and $i_0,i_k\ge0$.
  \item[(ii)]
    If $n_k<\infty$, write 
    $n_k\btilde_k=\sum_0^{k-1}m_{k,j}\btilde_j$
    with $0\le m_{k,j}<n_j$ for $1\le j<k$ and $m_0\ge0$
    as in Lemma~\ref{lem:arith}. Then
    \begin{equation}\label{equa2}
      \phi=p(T)\prod_{j=0}^{k}U_j^{i_j}
      \ \text{in}\ \gr_\nu R_\nu
    \end{equation}
    with $0\le i_j<n_j$ for $1\le j\le k$, $i_0,i_k\ge0$
    and where $p$ is a polynomial
    in $T\=U_k^{n_k}\prod_{j=0}^{k-1}U_j^{-m_{k,j}}$.
  \end{itemize}
  Both decompositions \eqref{equa1}, \eqref{equa2} are unique.
\end{theorem}
Recall that the ring $\gr_{\nu}\C(x)[y]$ is a Euclidean domain. 
Let us describe its irreducible elements.
\begin{corollary}\label{primes}
  Let $\nu$ be a valuation as above.
  \begin{itemize}
  \item[(i)]
    When $n_k=\infty$ the only irreducible element 
    of $\gr_{\nu}\C(x)[y]$ is $U_k$.
  \item[(ii)]
    When $n_k<\infty$ the irreducible elements of 
    $\gr_{\nu}\C(x)[y]$ consist of $U_k$ and all
    elements of the form 
    $U_k^{n_k}-\theta\prod_{j=0}^{k-1}U_j^{m_{k,j}}$ 
    for some $\theta\in\C^*$.
  \end{itemize}
\end{corollary}
\begin{proof}[Proof of Corollary~\ref{primes}]
  Assume $\phi\in\gr_\nu\C(x)[y]$ is irreducible.  If $n_k=\infty$,
  then $\phi=\a \prod_{j=0}^{k}U_j^{i_j}$ by~\eqref{equa1}. 
  But $U_j$ is a unit for $j<k$ by Lemma~\ref{L302}, so
  $U_k$ is the only irreducible element in $\gr_{\nu}\C(x)[y]$.
  
  When $n_k<\infty$, then we use~\eqref{equa2}
  and factorize 
  $p(T)=\prod(T-\theta_l)$.
  Modulo unit factors, we hence get
  \begin{equation}\label{fac}
    \phi=U_k^{i_k-Ln_k}\prod_l\left(
      U_k^{n_k}-\theta_l\prod_jU_j^{m_{k,j}}
    \right),
  \end{equation}
  where $L=\deg p$.
  On the other hand Lemma~\ref{L303} shows 
  that all elements of the form 
  $U_k^{n_k}-\theta_l\prod_jU_j^{m_{k,j}}$ are
  irreducible in $\gr_{\nu}\C(x)[y]$.  
  So~\eqref{fac} is the decomposition of $\phi$
  into prime factors in $\gr_{\nu}\C(x)[y]$.
  This concludes the proof.
\end{proof}
\begin{proof}[Proof of Theorem~\ref{homo}]
  First assume that $\deg_y\phi<d_k=\deg_yU_k$. 
  We will show that~\eqref{equa1} holds. 
  Write $\phi=\sum_0^{n_{k-1}-1}\phi_iU_{k-1}^i$ with 
  $\deg_y\phi_i<d_{k-1}$ for all $i$. 
  Iterating this procedure we get
  $\phi=\sum\a_IU_0^{i_0}\cdots U_{k-1}^{i_{k-1}}$ 
  with $\a_I\in\C$, $i_0\ge0$ and $0\le i_j<n_j$ for
  $j\ge 1$. 
  If 
  $\nu(U_0^{i_0}\dots U_{k-1}^{i_{k-1}})=\nu(U_0^{j_0}\dots U_{k-1}^{j_{k-1}})$
  then $\sum i_l\btilde_l=\sum j_l\btilde_l$ 
  so that $(i_{k-1}-j_{k-1})\btilde_{k-1}\in\sum_0^{k-2}\Z\btilde_j$. 
  Thus $n_{k-1}$ divides $i_{k-1}-j_{k-1}$ so 
  $i_{k-1}=j_{k-1}$ as $0\le i_{k-1},j_{k-1}<n_{k-1}$. 
  Iterating this argument we get~\eqref{equa1}.
  
  In the general case, we can assume $\deg_y\phi<\infty$
  by Weierstrass preparation theorem.
  Write $\phi=\sum\phi_iU_k^i$, with $\deg_y\phi_i<d_k$. 
  By what precedes, we have 
  \begin{equation*}
    \phi=\sum\a_IU_0^{i_0}
    \dots U_{k-1}^{i_{k-1}}U_k^{i_k}\ \text{in}\ \gr_\nu R 
  \end{equation*}
  with $\a_I\in\C^*$, $0\le i_j<n_j$ for $1\le j<k$, 
  $i_0,i_k\ge0$. We may assume 
  $\sum_0^k i_j\btilde_j=\nu(\phi)$
  for all $I$ with $\a_I \not=0$. 
  If~\eqref{equa1} is not valid, there exist
  $I\ne J$ with $\a_I,\a_J\ne0$.  
  Thus $(i_k-j_k)\btilde_k\in\sum_0^{k-1}\Z\btilde_j$ 
  implying $n_k<\infty$.
  This shows~\eqref{equa1} when $n_k=\infty$.
    
  Now assume $n_k<\infty$, and write
  $n_k\btilde_k=\sum_0^{k-1}m_{k,j}\btilde_j$ as in Lemma~\ref{lem:arith}.
  Let $I=(i_0,\dots,i_k)$ be any multiindex with $\a_I\ne0$. 
  Make the Euclidean division $i_k = r_k n_k + \wi _k$ with $0
  \leq \wi_k < n_k$, and write
  \begin{equation*}
    \prod_{j=0}^kU_j^{i_j}
    =U_k^{\wi_k}
    T^{r_k}
    \prod_{j=0}^{k-1}U_j^{a_j}
  \end{equation*}
  with $a_j\=i_j+r_km_{k,j}$ and
  $T\=U_k^{n_k}\prod_{j=0}^{k-1}U_j^{-m_{k,j}}$. 
  The key remark is now that
  $U_j^{n_j}=\theta_j\prod_0^{j-1}U_l^{m_{j,l}}$ in $\gr_\nu R$ 
  for $1\le j<k$. Making the Euclidean division 
  $a_{k-1}=r_{k-1}n_{k-1}+\wi_{k-1}$ with $0\le\wi_{k-1}<n_{k-1}$, 
  we get 
  \begin{equation*}
    \prod_{j=0}^{k-1}U_j^{a_j}
    =\theta_{k-1}U_k^{\wi_{k-1}}\prod_{j=0}^{k-2} U_j^{a'_j}
  \end{equation*}
  for some $a'_j\in\N$. We finally get by induction that
  \begin{equation*}
    \prod_{j=0}^kU_j^{i_j}
    =\theta_IT^{r_I}\prod_0^kU_j^{\wi_j}
  \end{equation*}
  in $\gr_\nu R_\nu$, with $\theta_I\in \C^*$, 
  $r_I(=r_k)\ge0$, $0\le\wi_j<n_j$ for $1\le j\le k$ 
  and $\wi_0\ge0$.
  
  Now $\sum_0^ki_j\btilde_j=\nu(\phi)$ and 
  $\nu(T)=0$, hence $\sum_0^k\wi_j\btilde_j=\nu(\phi)$. 
  Suppose $\sum_0^k\wi_j\btilde_j=\sum_0^k\ti_j\btilde_j$ 
  with $0\le\ti_j<n_j$ for $1\le j\le k$. 
  Since $|\wi_k-\ti_k |<n_k$, the definition of $n_k$
  gives $\wi_k=\ti_k$.
  By induction we get $\wi_j=\ti_j$ for all $j\ge0$.
  We have proved that $\phi$ can be written in the form~\eqref{equa2}.
  
  Uniqueness of both decompositions \eqref{equa1}, \eqref{equa2} comes
  from unique factorization in $\gr_{\nu}\C(x)[y]$.
  Indeed, when $n_k=\infty$, $i_k=\delta_k(\phi)$. 
  From $\sum i_j\btilde_j=\nu(\phi)$ and $i_j<n_j$ for $1\le j<k$
  we deduce uniqueness of the decomposition \eqref{equa1}.
  When $n_k<\infty$, $\nu(\phi)$ determines $i_j$ for all $j\ge0$, 
  and the polynomial $p(T)$ is determined 
  as $\phi$ admits a unique decomposition
  into prime factors (see Corollary~\ref{primes}). This concludes the
  proof of the theorem.
  \end{proof}
%
%
\subsection{Homogeneous decomposition II}
We now turn to infinite SKP's.
\begin{theorem}\label{graded-infinite}
  Let $\nu\=\val[(U_j)_0^\infty;(\btilde_j)_0^\infty]$.
  Pick $0\ne\phi\in R$.
  \begin{itemize}
  \item[(i)]
    If $n_j\ge2$ for infinitely many $j$, then 
    there exists $k=k(\phi)$ such that
    \begin{equation} \label{equa3}
      \phi=\a\prod_{j=0}^kU_j^{i_j}\ \text{in}\ \gr_\nu R
    \end{equation}
    with $\a\in\C^*$, $0\le i_j<n_j$ for $1\le j\le k$, 
    and $i_0\ge0$.
  \item[(ii)]
    If $n_j=1$ for $j\gg1$, then
    there exists $k=k(\phi)$ such that
    \begin{equation}\label{equa4}
      \phi=\a\cdot U_\infty^n\cdot\prod_{j=0}^kU_j^{i_j} 
      \ \text{in}\ \gr_{\nu}R
    \end{equation}
    with $\a\in\C$, $0\le i_j<n_j$ for $1\le j\le k$, $i_0\ge0$ 
    and $n\ge0$.
  \end{itemize}
  Both decompositions \eqref{equa3}, \eqref{equa4}  are unique.
\end{theorem}
The analogue of Corollary~\ref{primes} is
\begin{corollary}\label{primes-infinite}
  Assume $\nu$ is associated to an infinite SKP\@.
  \begin{itemize}
  \item[(i)]
    When $n_j\ge2$ for infinitely many $j$, 
    the ring $\gr_\nu\C(x)[y]$ is a field.
  \item[(ii)]
    When $n_j=1$ for $j\gg1$, the only irreducible element of
    $\gr_\nu\C(x)[y]$ is $U_\infty$.
  \end{itemize}
\end{corollary}
\begin{proof}[Proof of Theorem~\ref{graded-infinite}]
  Fix $\phi\in R$, and suppose it is in Weierstrass form, polynomial
  in $y$.  First assume that 
  $n_j\ge2$ for infinitely many $j$. 
  For $j\gg1$ $\deg_y\phi<\deg_y U_j$, hence 
  $\phi=\a\prod_{l=0}^{j-1}U_l^{i_l}$ 
  in $\gr_{\nu_j}R$  for some $\a\in \C^*$, and $i_l<n_l$  
  for $l\ge 1$ (see the beginning of the proof of Theorem~\ref{homo}).  
  As $\nu\ge\nu_j$, we get 
  $\phi=\a\prod_{l=0}^{j-1}U_l^{i_l}$ also in $\gr_\nu R$.
  
  When $n_j=1$ for $j\gg1$, write $\phi=U_{\infty}^n\phi'$ 
  with $\phi'$, $U_{\infty}$ prime. For $j\gg1$
  $\nu_{j+1}(\phi')=\nu_j(\phi')=\nu(\phi')$. 
  Hence $\phi'=\a\prod_{l=0}^{j-1}U_l^{i_l}$ in $\gr_{\nu_j}R$ 
  and $\gr_\nu R$ as above.
  \end{proof}
\begin{proof}[Proof of Corollary~\ref{primes-infinite}]
  The first assertion is immediate by Lemma~\ref{L302}.
  Assume $n_j=1 $ for $j\gg1$
  and pick $\phi\in R$. 
  If $U_\infty$ does not divide $\phi$, then~\eqref{equa3} holds.  
  For $k\gg1$ we get $\delta_k(\phi)=0$. 
  Thus $\phi$ is invertible in $\gr_{\nu_k}\C(x)[y]$, 
  hence in $\gr_\nu\C(x)[y]$.

  On the other hand $U_\infty=\phi\psi$ in $\gr_\nu\C(x)[y]$
  implies either $\nu(\phi)=\infty$ or $\nu(\psi)=\infty$. 
  Hence $U_\infty$ divides either $\phi$ or $\psi$ as formal power series
  hence in $\gr_\nu\C(x)[y]$. 
  We conclude noting that $U_\infty$ 
  cannot be a unit in $\gr_\nu\C(x)[y]$. 
  Indeed $\nu(U_\infty\phi)=\nu(1)$
  implies that $U_\infty\phi$ is prime with $U_\infty$.
  \end{proof}
%
%
\subsection{Value semigroups and numerical invariants}\label{numerical}
Knowing the structure of the graded rings allows us to compute the
value semigroups and numerical invariants introduced in
Section~\ref{krullval}. The first assumption
in the next theorem is in fact redundant as we shall see below in
Theorem~\ref{approx}.

\begin{theorem}\label{divis}
  Let $\nu$ be a valuation associated to a finite or infinite SKP.
  Then the value semi-group\footnote{Note that when $\nu$ is a curve
  valuation, it may  take the value $\infty$ on non-zero elements.}
  $\nu(R)$ is equal to $\sum\btilde_j\N$, Moreover, we have the
  following intrinsic characterization of the type of a valuation as
  introduced in Definition~\ref{D101} \index{valuation!numerical
  invariants of}.
  \begin{itemize} \item[(i)] $\nu$ is a curve valuation, iff
  $\rk(\nu)=\ratrk(\nu)=2$, $\trdeg(\nu)=0$.
  \item[(ii)]
    $\nu$ is divisorial, iff
    $\rk(\nu)=\ratrk(\nu)=1$ and $\trdeg(\nu)=1$.
  \item[(iii)]
    $\nu$ is irrational, iff
    $\rk(\nu)=1$, $\ratrk(\nu)=2$ and $\trdeg(\nu)=0$.
  \item[(iv)]
    $\nu$ is infinitely singular, iff
    $\rk(\nu)=\ratrk(\nu)=1$ and $\trdeg(\nu)=0$.
  \end{itemize}
\end{theorem}
\begin{proof}[Proof of Theorem~\ref{divis}]
  The computation of the value semigroup is a simple exercise
  that is left to the reader. As for the numerical
  invariants we treat the cases of finite and infinite SKP's
  separately.

  First assume that $\nu$ is associated to a finite SKP
  $[(U_j)_0^k;(\btilde_j)_0^k]$, $1\le k<\infty$.  The two invariants
  $\rk(\nu)$ and $\ratrk(\nu)$ can be read directly from the value
  semigroup $\nu(R)$ so their computation is straightforward thanks 
  to the first part of the proof.
  For the computation of $\trdeg(\nu)$ we rely
  on Theorem~\ref{homo} and proceed as follows.
  
  When $\btilde_0=\infty$, 
  $\btilde_k=\infty$ or $\btilde_k\not\in\Q^*$ 
  (so $\btilde_k\not\in\sum_0^{k-1}\Z\btilde_j$), 
  then~\eqref{equa1} applies. Pick $\phi,\psi\in R$ with
  $\nu(\phi)=\nu(\psi)$. 
  Then $\phi=\a\prod_{l=0}^kU_l^{i_l}$, 
  $\psi=\gamma\prod_{l=0}^k U_l^{j_l}$ modulo $\nu$, 
  with $\sum_0^ki_l\btilde_l=\sum_0^kj_l\btilde_l$.
  Hence $i_k=j_k$ since $n_k=\infty$. 
  Now $\sum_0^{k-1}i_l\btilde_l=\sum_0^{k-1}j_l\btilde_l$ and
  $0\le i_l,j_l<n_l$ for all $1\le j<k$ so $i_l=j_l$ for
  all $l$ by~\eqref{key3}. Hence $\phi/\psi=\alpha/\gamma\in\C^*$ 
  in $\gr_\nu R_\nu$.  
  We have shown that $k_\nu\=R_\nu/\fm_\nu\cong\C$ so
  $\trdeg(\nu)=0$.
  
  If $\btilde_0<\infty$ and $\btilde_k\in\Q^*$, then $n_k<\infty$. 
  For $\phi,\psi\in R$ with $\nu(\phi)=\nu(\psi)$ write
  \begin{equation*}
    \phi=p(T)\prod_{l=0}^kU_l^{i_l}
    \qand
    \psi=q(T)\prod_{l=0}^kU_l^{j_l}
    \quad\text{in $\gr_\nu R_\nu$}
  \end{equation*}
  as in~\eqref{equa2},
  with $0\le i_l,j_l<n_l$ for $1\le l\le k$, $i_0,j_0\ge0$ 
  and $p,q\in\C[T]$.
  As $\sum_0^ki_l\btilde_l=\sum_0^kj_l\btilde_l$, 
  one has $i_l=j_l$ for all $l$.  
  Hence $\phi/\psi=p(T)/q(T)$ in $\gr_\nu R_\nu$. 
  This shows $k_\nu\cong\C(T)$ so $\trdeg(\nu)=1$.

  Now assume $\nu$ is associated to an infinite SKP,
  $\nu\=\val[(U_j);(\btilde_j)]$.
  First assume $n_j\ge2$ for infinitely many $j$. 
  Set $\nu_k\=\val[(U_j)_0^k;(\btilde_j)_0^k]$. 
  Then the sequence $\nu_k(\phi)$ is eventually stationary
  for any $\phi\in R$.
  One easily checks that $\rk(\nu)=\ratrk(\nu)=1$.
  Let us compute $\trdeg(\nu)$. 
  Fix $\phi,\psi\in R$ with $\nu(\phi)=\nu(\psi)$.  
  For $k\gg 1$ we have 
  $\nu_k(\phi)=\nu(\phi)=\nu(\psi)=\nu_k(\psi)$. 
  Write $\phi=\a\prod_{l=0}^kU_l^{i_l}$ and
  $\psi=\gamma\prod_{l=0}^{k}U_l^{j_l}$
  where $\a,\gamma\in\C^*$, $0\le i_l,j_l<n_l$ 
  for $l\ge 1$ and apply the preceding arguments.  
  We get $\phi/\psi=\alpha/\gamma$ in 
  $\gr_{\nu_k}R_{\nu_k}$ hence in $\gr_\nu R_\nu$.
  This shows $k_\nu\cong\C$.

  The last case is when $n_j=1$ for $j\gg1$. We leave it to the reader to
  verify that $\rk(\nu)=\ratrk(\nu)=2$ and $\trdeg(\nu)=0$ in this case.
\end{proof}
%
%
%
%
\section{From valuations to SKP's}\label{sec-appro}
We now show that every centered valuation on $R$ is represented
by an SKP\@.
\begin{theorem}\label{approx}
  For any centered valuation $\nu$ on $R$, there exists 
  a unique SKP $[(U_j)_0^n;(\btilde_j)_0^n]$, 
  $1\le n\le\infty$, such that
  $\nu=\val[(U_j);(\btilde_j)]$.
  We have $\nu(U_j)=\btilde_j$ for all $j$.
  Further, if $k<n$ and $\nu_k:=\val[(U_j)_0^k;(\btilde_j)_0^k]$,
  then $\nu(\phi)\ge\nu_k(\phi)$ for all $\phi\in R$ and
  $\nu(\phi)>\nu_k(\phi)$ iff 
  $U_{k+1}$ divides $\phi$ in $\gr_{\nu_k}\C(x)[y]$.
\end{theorem}
\begin{remark}\label{R101}
  If $\nu$ is a Krull valuation, 
  the same result holds for a $\Gamma$-SKP\@.
\end{remark}
\begin{remark}
  In Spivakovsky's terminology~\cite{spiv}, the subset of $(U_j)$ for
  which $n_j>1$ forms a minimal generating sequence for the valuation
  $\nu$. The associated divisorial valuations $\nu_j$ will  play an
  important role in the sequel. They form what we call the 
  approximating sequence of $\nu$, see Section~\ref{sec-approx}.
\end{remark}
\begin{proof}[Proof of Theorem~\ref{approx}]
  We construct by induction on $k$ a valuation $\nu_k$ so that
  $\nu_k(U_j) = \btilde_j$ for $j\le k$.  Let $U_0=x$, $U_1=y$ and
  $\btilde_0=\nu(x)$, $\btilde_1=\nu(y)$.  Assume
  $\nu_k:=\val[(U_j)_0^k;(\btilde_j)_0^k]$ has been defined.  By
  Theorem~\ref{key}, $\nu(\phi)\ge\nu_k(\phi)$ for $\phi\in R$.  As
  $\nu(x)=\nu_k(x)$ this also holds for $\phi\in\C(x)[y]$.  If
  $\nu=\nu_k$, then we are done.  If not, set
  $\cD_k\=\{\phi\in\C(x)[y]\ ;\ \nu(\phi)>\nu_k(\phi)\}$.
  \begin{lemma} $\cD_k$ defines a prime ideal in
  $\gr_{\nu_k}\C(x)[y]$.  
  \end{lemma} 
  \begin{proof} Let us check that
  $\cD_k$ is well defined in $\gr_{\nu_k}\C(x)[y]$. Pick
  $\phi,\phi'\in\C(x)[y]$ with $\phi=\phi'$ modulo $\nu_k$.  If
  $\phi\in\cD_k$ but $\phi'\notin\cD_k$, then
  $\nu(\phi')=\nu_k(\phi')=\nu_k(\phi)<\nu(\phi)$.  So
  $\nu(\phi-\phi')=\nu(\phi')=\nu_k(\phi')<\nu_k(\phi-\phi')$, a
  contradiction.  To show that $\cD_k$ is a prime ideal is easy and
  left to the reader.  
  \end{proof}
  We continue the proof of the theorem.
  Recall that $\gr_{\nu_k}\C(x)[y]$ is a Euclidean domain.
  Hence Corollary~\ref{primes} and Lemma~\ref{lem:arith}
  show that $\cD_k$ is generated by 
  a unique irreducible element
  $U_{k+1}=U_k^{n_k}-\theta_k\prod_{j=0}^{k-1}U_j^{m_{k,j}}$
  with $\theta_k\in\C^*$, $0\le m_{k,j}<n_j$ 
  for $j\ge1$ and $m_{k,0}\ge0$.
  Define $\btilde_{k+1}\=\nu(U_{k+1})$. It is easy to check 
  that $[(U_j)_0^{k+1};(\btilde_j)_0^{k+1}]$ is an SKP\@.
  This completes the induction step.

  Either the induction terminates at some finite $k$, or 
  we get an infinite SKP $[(U_j);(\beta_j)]$. In the latter
  case we claim that $\nu=\val[(U_j);(\beta_j)]$. This 
  amounts to showing that if $\phi\in R$, then the
  increasing sequence $\nu_k(\phi)$ converges to 
  $\nu(\phi)$.
  
  If $\nu_k(\phi)=\nu(\phi)$ for $k\gg1$ then we are done,
  so assume $\nu_k(\phi)<\nu(\phi)$ for all $k$. 
  Then $\phi\in\cD_k$ so 
  $U_{k+1}$ divides $\phi$ in $\gr_{\nu_k}\C(x)[y]$. 
  In particular, $\nu_{k+1}(\phi)>\nu_k(\phi)$, and 
  $\deg_y(\phi)\ge\deg_y(U_{k+1})=d_{k+1}$. 
  Hence $d_k$ is bounded, and $n_k=1$ for $k\gg1$. 
  By Theorem~\ref{key-infinite}, 
  $U_k\to U_{\infty}$ in $R$. 
  Moreover $U_{\infty}$ divides $\phi$ in $R$
  or else $\nu_k(\phi)$ would be stationary for $k$ large. 
  But then $\nu(\phi)\ge\lim\nu_k(\phi)=\infty$ so 
  $\nu_k(\phi)\to\nu(\phi)$.
  This completes the proof of the theorem.
  \end{proof}
%
%
%
%
\section{A computation}\label{some-compu}
We now compute $\nu(\phi)$ for a normalized valuation 
$\nu$ and $\phi\in\fm$ irreducible
in terms of SKP's. This computation will be crucial
to describe the parameterization of valuation space; 
see Lemma~\ref{compuskew}.

Let $\nu_{\phi}$ be the curve valuation associated to $\phi$.
Write
$\nu=\val[(U_j);(\btilde_j)]$
and
$\nu_\phi=\val[(U_j^\phi);(\btilde_j^\phi)]$.
Assume $\nu\ne\nu_\phi$ and 
define the contact order of $\nu$ and $\nu_\phi$ by
\begin{equation}\label{contact}
  \con(\nu,\nu_\phi)=\max\{j\ ;\ U_j=U_j^\phi\}.
\end{equation}
Let $n_j^\phi$ be the integers defined by~\eqref{key3}
for $\nu_\phi$ 
and set $\gamma_k^\phi=\prod_{j\ge k}n_j^\phi$ for $k\ge1$.
These products are in fact finite as $n_j^\phi\,=1$ for $j\gg1$.
\begin{proposition}\label{vphi}
  If $\phi=x$ (up to a unit), then $\nu(\phi)=\btilde_0$. Otherwise
  \begin{equation}\label{e502}
    \nu(\phi)
    =\gamma^\phi_k\,\min\{\btilde_k,\btilde_k^\phi\}
    \min\{1,\btilde_0/\btilde_0^\phi\},
  \end{equation}
  where $k=\con(\nu,\nu_\phi)$.
\end{proposition}
Later on we will be interested in the quotient
$\nu(\phi)/m(\phi)$. 
Applying Proposition~\ref{vphi} to $\nu=\nu_\fm$ we obtain
$m(\phi)=\gamma^\phi_1/\btilde^\phi_0$
if $\phi\ne x$; this leads to
\begin{proposition}\label{vvm}
  If $\nu\in\cV$ and $\phi\in\fm$ is irreducible, then
  \begin{equation*}
    \frac{\nu(\phi)}{m(\phi)}
    =d_k^{-1}\min\{\btilde_k,\btilde_k^\phi\,\}
    \min\{\btilde_0,\btilde^\phi_0\},
  \end{equation*}
  where $k=\con(\nu,\nu_\phi)$ and $d_k=\deg_yU_k=\prod_1^{k-1}n_j$.
\end{proposition}
\begin{proof}[Proof of Proposition~\ref{vphi}]
  The case $\phi=x$ is trivial, so assume $\phi\ne x$, \ie
  $\btilde^\phi_0<\infty$.
  We may then assume $\phi=U^\phi_l$,
  where $l\=\lgt\nu_\phi\in[k,\infty]$.
  Note that 
  $\min\{\btilde_0,\btilde_1\}=\min\{\btilde^\phi_0,\btilde^\phi_1\}=1$
  since $\nu$ and $\nu_\phi$ are normalized.

  Let us first consider the case when
  $\btilde_0=\btilde_0^\phi$ 
  and
  $\btilde_1=\btilde_1^\phi$; this holds \eg if $k\ge2$.
  Then $\btilde_j=\btilde_j^\phi$ for $0\le j<k$.
  If $l=k$, then $\phi=U_k$ so $\nu(\phi)=\btilde_k$
  and $\btilde^\phi_k=\infty$, which implies~\eqref{e502}.
  Hence assume $l>k$.
  Define $\xi_j=\btilde^\phi_{k+j}$ and
  $\eta_j=\min\{\btilde_k,\btilde^\phi_k\}\gamma^\phi_k/\gamma^\phi_{k+j}$
  for $0\le j\le l-k$. 
  Then $\nu_\phi(U^\phi_{k+j})=\xi_j$ for $j\ge0$
  and $\nu(U^\phi_j)=\btilde_j$ for $0\le j\le k$.
  We will prove inductively that
  $\nu(U^\phi_{k+j})=\eta_j$ for $j\ge1$; when $j=l-k$ this
  gives~\eqref{e502}.
  The induction is based on~\eqref{e501}, which reads
  \begin{equation}\label{e503}
    U^\phi_{k+j+1}
    =(U^\phi_{k+j})^{n^\phi_{k+j}}
    -\theta^\phi_{k+j}\prod_{i=0}^{k+j-1}(U^\phi_i)^{m^\phi_{k+j,i}}
    =:A_j-B_j.
  \end{equation}

  Let us first show that $\nu(U^\phi_{k+1})=\eta_1$, using~\eqref{e503}
  for $j=0$. 
  There are three cases, depending on $\btilde_k$ and $\btilde^\phi_k$.
  The first case is when $\btilde_k>\btilde_k^\phi$. 
  Then
  $n_k^\phi\,\btilde_k>n_k^\phi\,\btilde_k^\phi=\sum_0^{k-1}m^\phi_{k,i}\btilde_i$, 
  so using~\eqref{e503} we obtain
  $\nu(U_{k+1}^\phi)=\btilde^\phi_kn_k^\phi=\eta_1$. 
  Similarly, in the second case, $\btilde_k<\btilde_k^\phi$, then 
  $\sum_0^{k-1}m^\phi_{k,i}\btilde_i=n_k^\phi\btilde_k^\phi>n_k^\phi\btilde_k$, 
  so that $\nu(U_{k+1}^\phi)=n_k^\phi\btilde_k=\eta_1$.
  The third case is when $\btilde_k=\btilde_k^\phi$.
  Set $\nu_k\=\val[(U_j)_0^k;(\btilde_j)_0^k]$. 
  Then $U_{k+1}^\phi$ is irreducible in 
  $\gr_{\nu_k}\C(x)[y]$ and
  $\nu_k(U_{k+1}^\phi)=\btilde_kn_k^\phi=\eta_1$. 
  If $\lgt\nu=k$ then $\nu=\nu_k$ and
  $\nu(U_{k+1}^\phi)=\eta_1$ so assume $\lgt\nu>k$.
  By assumption $U_{k+1}^\phi\not=U_{k+1}$, 
  hence $U_{k+1}^\phi$ does not belong to the ideal
  generated by $U_{k+1}$ in $\gr_{\nu_k}\C(x)[y]$. 
  But this ideal coincides with $\{\nu>\nu_k\}$ by Theorem~\ref{approx}, 
  so $\nu(U_{k+1}^\phi)=\nu_k(U_{k+1}^\phi)=\eta_1$.

  Now fix $1\le j<l-k$ and assume that 
  $\nu(U^\phi_{k+i})=\eta_i$ for $1\le i\le j$. 
  We will prove that $\nu(U^\phi_{k+j+1})=\eta_{j+1}$ using~\eqref{e503}.
  Write $a_i=m^\phi_{k+j,k+i}$ for $0\le i<j$ and
  $c=n^\phi_{k+j}$. 
  Then $\nu(A_j)=c\eta_j$ and
  \begin{multline*}
    \nu(B_j)
    =\sum_{i=1}^km^\phi_{k+j,i}\btilde_i
    +\sum_{i=1}^{j-1}a_i\eta_i
    =\sum_{i=1}^{k+j-1}m^\phi_{k+j,i}\btilde^\phi_i
    +a_0(\btilde_k-\xi_0)
    +\sum_{i=1}^{j-1}a_i(\eta_i-\xi_i)\\
    \ge\sum_{i=1}^{k+j-1}m^\phi_{k+j,i}\btilde^\phi_i
    +\sum_{i=0}^{j-1}a_i(\eta_i-\xi_i)
    =c\xi_j
    +\sum_{i=0}^{j-1}a_i(\eta_i-\xi_i).
  \end{multline*}
  To show that $\nu(U^\phi_{k+j+1})=\eta_{j+1}$ we
  only need to show $\nu(B_j)>c\eta_j$ or simply
  $\sum_0^{j-1}a_i(\xi_i-\eta_i)>c(\xi_j-\eta_j)$.
  Note that $\sum_0^{j-1}a_i\xi_i\le c\xi_j$ and that
  the sequence $(\eta_i/\xi_i)_0^j$ is strictly decreasing.
  Set $p_i=a_i\xi_i/\sum_0^{j-1}a_i\xi_i$. Then
  $\sum_0^{j-1}p_i=1$ so
  \begin{equation*}
    \sum_{i=0}^{j-1}a_i(\xi_i-\eta_i)
    =\left(\sum_{i=0}^{j-1}a_i\xi_i\right)
    \sum_{i=0}^{j-1}p_i\left(1-\frac{\eta_i}{\xi_i}\right)
    <c\xi_j\left(1-\frac{\eta_j}{\xi_j}\right)
    =c(\xi_j-\eta_j).
  \end{equation*}
  This completes the proof when
  $\btilde_0=\btilde_0^\phi$ and $\btilde_1=\btilde_1^\phi$.
  The remaining cases all have $k=1$ and are as follows:
  $\btilde_0\ge\btilde_1=1$ and $1=\btilde_0^\phi\le\btilde_1^\phi$; 
  $\btilde_0\ge\btilde_1=1$ and $\btilde_0^\phi\ge\btilde_1^\phi=1$; 
  $1=\btilde_0\le\btilde_1$ and $1=\btilde_0^\phi<\btilde_1^\phi$; 
  $1=\btilde_0\le\btilde_1$ and $\btilde_0^\phi\ge\btilde_1^\phi=1$.  
  These are handled in the same way as above. The details are
  left to the reader.
\end{proof}


%
%
%
%
%
%
\chapter{Tree structures}\label{part3}
In this third chapter we show that valuation space $\cV$
has the structure of a tree; we shall subsequently
refer to it as \emph{the valuative tree}.

Roughly speaking, 
a tree\footnote{We use the term ``tree'' instead
of ``$\R$-tree'' when the tree is modeled on the real line.} 
is a union of real intervals welded
together in such a way that no cycles appear. 
We will more precisely distinguish between three types of 
tree structures, see below.
The main result of the chapter is then 
that valuation space can be naturally
equipped with all three of these structures and even a
bit more. What we obtain is effectively a coordinate free
visualization of the encoding of valuations
by SKP's (which are defined using a fixed choice of
local coordinates). 
Indeed, SKP's play an instrumental role in
most proofs in this chapter.

It is good to keep in mind that 
the quite intricate structure of
the valuative tree that we are about to develop
all derives from the definition of 
an element of $\cV$ as a function on
$R$ satisfying certain axioms. 
In Chapter~\ref{A3} we shall arrive at the same 
tree structure using a purely geometric construction.

The organization of this chapter is as follows.
In Section~\ref{tree-subsec},
we discuss three different types of trees.
First we have \emph{nonmetric trees},
defined as partially ordered sets satisfying certain axioms. 
In particular, every ``full'', totally ordered subset can be
parameterized (in a noncanonical way) by a real interval.
This notion seems to be new, although related to the approach 
of Berkovich~\cite{Ber}. 
Then we have \emph{parameterized trees}. These are nonmetric
trees that come with a fixed parameterization.
Finally we have \emph{metric trees}
(typically called $\R$-trees in the literature). 
Metric trees are closely related to parameterized trees,
and are defined as metric spaces satisfying certain conditions.

As we show in Section~\ref{tree-struc}, 
the natural partial ordering on $\cV$
induces a nonmetric tree structure.
We analyze this structure in detail.
In particular, we show that the set $\cVqm$ of quasimonomial
valuations is exactly the tree $\cV$ with all ends removed,
hence in itself is a nonmetric tree. 

In Sections~\ref{metric-tree-struc},~\ref{sec-mult},~\ref{S10}, we
introduce two natural parameterizations of $\cV$ and $\cVqm$ and show
they can be used to exhibit these spaces as metric trees. These
parameterizations will play a fundamental role in
applications~\cite{pshfnts}~\cite{criterion},~\cite{eigenvaluation}.

We first define in Section~\ref{metric-tree-struc} a numerical
invariant of a valuation, its \emph{skewness}.  This gives the first
parameterization of both $\cV$ and $\cVqm$.  

We then define in Section~\ref{sec-mult} a discrete invariant,
the \emph{multiplicity} of a valuation. This invariant naturally
extends the notion of multiplicity of a curve.  To a divisorial
valuation is also associated a \emph{generic multiplicity}.

Multiplicity is an increasing function on $\cV$ with values in
$\Nbar$. By analyzing the points where it jumps, we 
can define a canonical \emph{approximating sequence} 
of a given valuation (Section~\ref{sec-approx}).

Using multiplicity and skewness, we define in Section~\ref{S10} a
third important invariant of a valuation: its \emph{thinness}.  It has
a particular geometric interpretation which makes it extremely useful
for applications. Thinness is obtained by integrating multiplicity
with respect to skewness and hence gives the second parameterization
of $\cV$ and $\cVqm$ as mentioned above.

In Section~\ref{S12} we use the approximating sequences
to compute the value group of a valuation, as well as
the generic multiplicity of a divisorial valuation.

After that, we present a different but intriguing approach 
to the valuative tree.
Namely, we show in Section~\ref{val-curve} that a quasimonomial
valuation can be identified with a ball of irreducible curves in a
particular (ultra-)metric. The partial ordering, skewness and
multiplicity on the valuative tree then have natural interpretations
as statements about balls of curves.

We conclude the chapter by a study of the 
\emph{relative valuative tree}. 
By definition, this is the (closure of) the set of centered 
valuations on $R$ normalized by the condition
$\nu(x)=1$, where $\{x=0\}$ is a smooth formal curve;
such a normalization is natural in several situations.
As we show, the relative valuative tree $\cV_x$ 
has a structure which is very similar to that of $\cV$.
%
%
%
%
\section{Trees}\label{tree-subsec}
In this section we discuss different types of trees.
%
%
\subsection{Rooted nonmetric trees}\label{S35}
Consider a partially ordered set, or \emph{poset}\index{poset} 
$(\cT,\le)$. Let us say that a totally ordered subset 
$\cS\subset\cT$ is \emph{full}\index{full set (in a tree)} 
if $\sigma,\sigma'\in\cS$, $\tau\in\cT$ and
$\sigma\le\tau\le\sigma'$ imply $\tau\in\cS$.
\begin{definition}\label{D402}
  A \emph{rooted nonmetric tree}\index{tree!rooted nonmetric} 
  is a poset $(\cT,\le)$ such that
  \begin{itemize}
  \item[(T1)]
    $\cT$ has a unique minimal element $\tau_0$, called the
    \emph{root}\index{root (of a tree)} of $\cT$;
  \item[(T2)]
    if $\tau\in\cT$, then the set 
    $\{\sigma\in\cT\ ;\ \sigma\le\tau\}$ is isomorphic to a real interval;
  \item[(T3)]
    every full, totally ordered subset of $\cT$ is isomorphic
    to a real interval.
  \end{itemize}
\end{definition}
Statements~(T2) and~(T3) assert that there exists an order preserving 
bijection from a real interval onto the corresponding set. However,
there may not be a canonical choice of bijection.
\begin{remark}
  Condition~(T3) may seem superfluous in view of~(T2) but is necessary
  to avoid a ``long half-line'', \ie a totally ordered set $(\cT,\le)$
  with a (unique) minimal element $\tau_0$ for which every set
  $\{\sigma\in\cT\ ;\ \tau_0\le\sigma\le\tau\}$ 
  is isomorphic to a real interval but the full set $\cT$ is not.
\end{remark}
\begin{remark}\label{R401}  
  In fact, it is useful---and not hard to see---that if (T1)-(T2) hold,
  then (T3) is equivalent to
  \begin{itemize}
  \item[(T3')]
    if $\cS$ is a totally ordered subset of $\cT$ without upper 
    bound in $\cT$, then there exists a countable increasing sequence 
    in $\cS$ without upper bound in $\cT$.
  \end{itemize}
\end{remark}
\begin{remark}
  More generally, if $\Lambda$ is a totally ordered set, then a
  \emph{rooted nonmetric $\Lambda$-tree}\index{tree!$\Lambda$-}
  is a partially ordered set such that~(T1)-(T3)
  hold, with the intervals in~(T2) and~(T3) being intervals
  in $\Lambda$. Besides $\Lambda=\R$\index{tree!$\R$-},
  interesting examples include $\Lambda=\N$,
  $\Lambda=\Nbar$ and $\Lambda=\Q$.
  \index{tree!$\N$-}\index{tree!$\Nbar$-}\index{tree!$\Q$-}
\end{remark}
It follows from the completeness of $\R$
that every subset $S\subset\cT$ admits an 
\emph{infimum}\index{infimum!in a tree},
denoted by $\wedge_{\tau\in S}\tau$.
Indeed, the set $\{\sigma\in\cT\ ;\ \sigma\le\tau\ \forall\tau\in S\}$
is isomorphic to the intersection of closed real intervals 
with common left endpoint.

If $\cT$ is a rooted, nonmetric tree and 
$\tau_1,\tau_2$ are two points in $\cT$, then we set
\begin{equation*}
  [\tau_1,\tau_2]:=\{\tau\in\cT\ ;\ 
  \tau_1\wedge\tau_2\le\tau\le\tau_1\qor\
  \tau_1\wedge\tau_2\le\tau\le\tau_2\}.
\end{equation*}
We call $[\tau_1,\tau_2]$ a \emph{segment}, and define
$[\tau_1,\tau_2[:=[\tau_1,\tau_2]\setminus\{\tau_2\}$.
The segments $]\tau_1,\tau_2]$ and $]\tau_1,\tau_2[$ are
defined similarly.

If $\cS$, $\cT$ are rooted nonmetric trees with roots
$\sigma_0$, $\tau_0$, then a mapping $\Phi:\cS\to\cT$ is a 
\emph{morphism of rooted nonmetric trees} 
\index{morphism!of rooted nonmetric trees} 
if for any $\sigma\in\cS$, $\Phi$ gives an 
order preserving bijection of $[\sigma_0,\sigma]$
onto $[\tau_0,\Phi(\sigma)]$.
If $\Phi$ is also bijective, then it is an 
\emph{isomorphism of rooted nonmetric trees}.
\index{isomorphism!of rooted nonmetric trees} 
and we say that $\cS$ and $\cT$ are \emph{isomorphic}.
\index{trees!isomorphic}

A \emph{subtree} of a rooted nonmetric tree $(\cT,\le)$ is a subset
$\cS$ such that $\sigma\in\cS$, $\tau\in\cT$ and $\tau\le\sigma$
implies $\tau\in\cS$. Clearly $\cS$ is then a rooted
nonmetric tree with root $\tau_0$.
We say that $\cS$ is a \emph{finite} (\emph{countable}) 
subtree\index{tree!finite}
\index{tree!countable}if $\cS$ has finitely (countably) many branch
points and the tangent space at each branch point is 
finite (countable).

A rooted nonmetric tree $\cT$ is \emph{complete}
\index{tree!complete}
if every increasing sequence $(\tau_i)_{i\ge1}$ in $\cT$ has a
majorant, \ie an element $\tau_\infty\in\cT$ with
$\tau_i\le\tau_\infty$ for every $i$.  
Thanks to (T3), any rooted
nonmetric tree $\cT$ has a \emph{completion}
\index{tree!completion of}
$\bcT$ obtained by adding points corresponding to unbounded increasing
\index{$\bcT$ (completion of a tree $\cT$)}
sequences $\tau_i$ in $\cT$. 
A maximal point in $\bcT$ (under $\le$) is called an \emph{end}
\index{end (of a tree)} 
\index{tree!end of} 
of $\cT$.
Hence all points in $\bcT\setminus\cT$ are ends.
We sometimes use the notation $\cT^o$ for the sets of elements
\index{$\cT^o$ (nonmaximal elements of a tree $\cT$)}
of $\cT$ that are not ends, \ie the set of nonmaximal elements
of $\cT$.
%
%
\subsection{Nonmetric trees}\label{tree-def}
In some situations, the choice of root in a rooted nonmetric tree
is not so important. By ``forgetting'' the root, we obtain
an object called a (nonrooted) \emph{nonmetric tree}.

Let us be more precise. 
Consider a rooted, nonmetric tree $(\cT,\le)$
with root $\tau_0$.  Pick any point $\tau'_0\in\cT$ and define a
new partial ordering $\le'$ on $\cT$ by declaring
$\tau_1\le'\tau_2$ iff $[\tau'_0,\tau_1]\subset[\tau'_0,\tau_2]$.
Then $(\cT,\le')$ is a nonmetric tree rooted at $\tau'_0$.  One
easily verifies that segments in $(\cT,\le')$ are the same
as segments in $(\cT,\le)$.

More generally, for any set $\cT$, let $P=P(\cT)$ be the 
(possibly empty) set of partial orderings $\le$ on $\cT$ 
such that $(\cT,\le)$ is a rooted, nonmetric tree. 
If $\le_1$ and $\le_2$ are partial orderings in $P$, 
then we say that $\le_1$ and $\le_2$ are \emph{equivalent}
if the segments in $(\cT,\le_1)$ and $(\cT,\le_2)$ are
the same. It is straightforward to verify that 
two partial orderings are equivalent iff one is obtained from
the other by changing the root as above.

This leads to the following definition.
\begin{definition}\label{D408}
  A \emph{nonmetric tree} \index{tree!nonmetric} 
  is a set $\cT$ with $P(\cT)\ne\emptyset$,
  together with a nonempty equivalence class in $P(\cT)$.
\end{definition}

If $\cT$ is a nonmetric tree, we obtain a canonical
rooted, nonmetric tree $(\cT,\le)$ by fixing a point 
(the root) in $\cT$.

Concepts on rooted, nonmetric trees that can be formulated
purely in terms of segments carry over to the nonrooted setting.
For instance, 
a nonmetric tree $\cT$ is \emph{complete}\index{tree!complete} 
if the rooted nonmetric tree $(\cT,\le)$ is complete 
for some choice of root.
This makes sense since $\cT$ is complete iff every open 
segment in $\cT$ is contained in a closed segment.

Given a nonmetric tree and a point $\tau\in\cT$ 
we define an equivalence relation on
$\cT\setminus\{\tau\}$ by declaring
$\sigma$, $\sigma'$ to be equivalent 
if the segments $]\tau,\sigma]$ and $]\tau,\sigma']$ intersect. 
An equivalence class is
called a \emph{tangent vector}
\index{tangent vector! in a tree} 
at $\tau$ and the set of tangent vectors is called the 
\emph{tangent space}\index{tangent space!in a tree} at $\tau$, 
denoted $T\tau$.  
We say that the point $\sigma$ \emph{represents} the tangent vector.
The tangent spaces should be thought of as projectivized. 
This is natural since there is no canonical parameterization
of segments by real intervals.

A point in $\cT$ is an \emph{end}\index{point (in a tree)!end} if its
tangent space has only one element.  It is a \emph{regular
point}\index{point (in a tree)!regular} if the tangent space has two
elements, and a \emph{branch point}\index{point (in a tree)!branch}
otherwise.

Notice that the concept of end differs slightly between rooted
and nonrooted nonmetric trees: the root of a rooted, nonmetric tree
$(\cT,\le)$ may be an end in the associated nonrooted tree 
but is never an end in $(\cT,\le)$ itself.

Two nonmetric trees $\cT_1$ and $\cT_2$ are
\emph{isomorphic}\index{trees!isomorphic} if there exists a bijection
$\Phi:\cT_1\to\cT_2$ and partial orderings $\le_1$ and $\le_2$ such
that $\Phi$ is an isomorphism of the rooted nonmetric trees
$(\cT_1,\le_1)$ and $(\cT_2,\le_2)$, \ie $\Phi$ is order-preserving.
Then $\Phi$ is called an 
\emph{isomorphism of nonmetric trees}.
\index{isomorphism!of nonmetric trees}
%
%
\subsection{Parameterized trees}
A \emph{parameterization}\index{parameterization!of a tree} of a
\index{tree!parameterized}
rooted, nonmetric tree $(\cT,\le)$ is an increasing (or decreasing)
mapping $\alpha:\cT\to[-\infty,+\infty]$ whose restriction to any
full, totally ordered subset of $\cT$ gives a bijection onto a real
interval. A rooted, nonmetric tree is
\emph{parameterizable}\index{tree!parameterizable} if it admits a
parameterization.

Note that by postcomposing with a suitable monotone function, 
we may require the parameterizations to be increasing with values in
$[0,\infty]$ or even in $[0,1]$. In the rest of this chapter, 
we shall always work with such parameterizations. 

Let us point out that in some situations (notably for the much of 
the analysis in Chapter~\ref{part-potent} and its applications
in Chapter~\ref{part-appli-analysis}) the choice of 
parameterization is important.
At any rate, on the valuative tree, all the parameterizations
that we shall be concerned with are increasing, 
with values in $\Rbar$.

Consider a rooted nonmetric tree $(\cT,\le)$ with root $\tau_0$
and a parameterization $\a:\cT\to[0,1]$.
Pick any $\tau_0'$ and consider the equivalent partial
ordering $\le'$ rooted in $\tau_0'$. 
Then the function $\a':\cT'\to[0,2]$ defined by
$\a'(\tau)=\a(\tau)+\a(\tau'_0)-2\a(\tau\wedge\tau'_0)$,
where $\wedge$ defines the minimum with respect to $\le$,
gives a parameterization of the rooted, nonmetric tree
$(\cT,\le')$.
We may therefore define a 
(nonrooted) nonmetric tree $\cT$ to be
\emph{parameterizable} if $(\cT,\le)$ is
parameterizable for any choice of root.

It follows easily that a nonmetric tree $\cT$ is 
parameterizable iff its completion $\bcT$ is. 
We do not know if there exists a non-parameterizable 
nonmetric tree.
All nonmetric trees we consider in this monograph will 
be parameterizable.

A \emph{morphism of parameterized trees} 
\index{morphism!of parameterized trees}
between $(\cT_1,\a_1)$ and $(\cT_2,\a_2)$ 
is a morphism $\Phi:\cT_1\to\cT_2$ of 
the underlying rooted, nonmetric trees, 
such that $\alpha_2\circ\Phi=\alpha_1$. 
If $\Phi$ is moreover a bijection then
it is an \emph{isomorphism of parameterized trees}
\index{isomorphism!of parameterized trees}
and we say that $(\cT_1,\a_1)$ and $(\cT_2,\a_2)$ are 
\emph{isomorphic}.
\index{trees!isomorphic}
Of course $\Phi$ is then also an isomorphism of the 
underlying (rooted) nonmetric trees $\cT_1$ and $\cT_2$.
\begin{example}\label{E1}
  Fix a set $X$ and set $\cT=X\times[0,\infty)/\sim$, where 
  $(x,0)\sim(y,0)$ for any $x,y\in X$. Then $\cT$ is a rooted,
  nonmetric tree under the partial ordering $\le$ defined by
  $(x,s)\le(y,t)$ iff either $s=t=0$, or $x=y$ and $s\le t$.

  Given $g:X\to(0,\infty)$ define a parameterization 
  $\a_g:\cT\to[0,\infty)$ on $\cT$ by $\a_g(x,s)=g(x)s$.
  The identity map $(\cT,\a_g)\to(\cT,\a_h)$ is then an 
  isomorphism of parameterized trees iff $g\equiv h$.
\end{example}
%
%
\subsection{The weak topology}\label{S302}
A nonmetric tree carries a natural 
\emph{weak topology}\index{topology!weak tree} defined as follows. 
If $\vv\in T\tau$ is a tangent vector at a point $\tau\in\cT$, 
\index{$\vv$ (tangent vector)}
set
\begin{equation*}
  U(\vv):=\{\sigma\in\cT\setminus\{\tau\}\ ;\ 
  \text{$\sigma$ represents $\vv$}\}.
\end{equation*}
Then the weak topology is generated by the sets $U(\vv)$ 
(\ie the open sets are unions of finite intersections of such sets). 

We shall study the weak topology in more detail in
Section~\ref{S27}, but let us summarize the main features.
The weak topology is Hausdorff. 
Any complete subtree $\cS$ of $\cT$ is weakly closed in $\cT$
and the injection $\cS\to\cT$ is an embedding.
In particular, any segment $\gamma=[\tau,\tau']$ in $\cT$
is closed, and the induced topology on $\gamma$
coincides with the standard topology on $[0,1]$ under the
identification of $\gamma$ with the latter set.
However, the branching of a nonmetric tree 
does play an important role for the weak topology:
\begin{example}
  Let $\cT=X\times[0.\infty[\,/\sim$ be as in Example~\ref{E1}.
  Assume that $X$ is infinite and pick a sequence $(x_n)_1^\infty$
  of distinct elements of $X$. Then $(x_n,1)$ 
  converges weakly to the root $(x,0)$ as $n\to\infty$.
\end{example}
See Proposition~\ref{P434} for a precise criterion of
weak sequential convergence.
In Section~\ref{S27} we shall prove that
\emph{any parameterizable, complete, nonmetric tree is weakly compact}.
Note that a nonmetric tree $\cT$ which is not complete cannot be 
weakly compact: a sequence in $\cT$ that increases (with respect
to some choice of partial ordering) to an element
in $\bcT\setminus\cT$ cannot have a convergent subsequence.
\begin{proposition}\label{P305}
  If $\a:\cT\to[0,\infty)$ is a parameterization of a rooted,
  nonmetric tree, then the function $\a$ is 
  weakly lower semicontinuous.
\end{proposition}
\begin{proof}
  We have to show that the superlevel set 
  $\cT_t=\{\tau\ ;\ \a(\tau)>t\}$ is open for every $t$. 
  Consider $t$ such that $\cT_t$ is nonempty and pick
  $\tau\in\cT_t$. Pick $\tau'<\tau$ with $\a(\tau')\ge t$ and let
  $\vv$ be the tree tangent vector at $\sigma'$ represented
  by $\tau$. Then $U(\vv)$ is an open neighborhood of $\tau$
  on which $\a>t$. This completes the proof.
\end{proof}
\begin{remark}
  The function $\a$ is \emph{not} weakly continuous in general 
  as can be seen from Example~\ref{E1} with $g\equiv1$
  and $X$ an infinite set.
  See also Propositions~\ref{P303} and~\ref{P302}.
\end{remark}
%
%
\subsection{Metric trees} 
Closely connected to parameterized trees are 
\emph{metric trees}\index{tree!metric}.\footnote{Metric trees are normally 
called $\R$-trees in the literature but this terminology
would not be precise enough for our purposes.}
These are metric spaces in which every two points is joined by a 
unique arc, or segment, and this segment is isometric to a real interval.
It is known~\cite{MayOve} that a metrizable topological
space is a metric tree (\ie admits a compatible metric 
under which it becomes a metric tree) iff it is 
uniquely pathwise connected and locally pathwise connected.

Roughly speaking, metric trees are to parameterized trees
what nonmetric trees are to rooted, nonmetric trees.
Let us make this precise.
First, a metric tree $(\cT,d)$ gives rise to a nonmetric tree.
Indeed, fix $\tau_0\in\cT$ (the root) 
and define a partial ordering $\le$ on $\cT$:
$\tau\le\tau'$ iff $\tau$ belongs to the segment between
$\tau_0$ and $\tau'$.  Clearly~(T1)-(T3) hold. 
Moreover, different choices of $\tau_0$ give rise,
by the very definition, to equivalent partial orderings
on $\cT$. Thus $\cT$ is naturally a (nonrooted) nonmetric tree.
We say that the metric $d$ on $\cT$ is \emph{compatible}
with the nonmetric tree structure.

Second, if $(\cT,d)$ is a metric tree, then 
given any choice of root $\tau_0\in\cT$,
the function $\a:\cT\to[0,\infty)$ 
defined by $\a(\tau)=d(\tau,\tau_0)$ gives a
parameterization of the rooted, nonmetric tree $(\cT,\le)$.

Conversely, consider a rooted, nonmetric tree $(\cT,\le)$ 
with a parameterization $\alpha:\cT\to[0,1]$. 
Define a function $d:\cT\times\cT\to[0,2]$ by
\begin{equation*}
  d(\sigma,\tau)=(\alpha(\sigma)-\alpha(\sigma\wedge\tau))
  +(\alpha(\tau)-\alpha(\sigma\wedge\tau)).
\end{equation*}
\begin{proposition}\label{P501}
  $(\cT,d)$ is a metric tree.
\end{proposition}
\begin{proof}
  It is straightforward to verify that $(\cT,d)$ is a metric space in
  which every two points is joined by an arc isometric to a real
  interval.  What remains to be seen is that there is a \emph{unique}
  arc between any two points in $\cT$.  For this, consider a
  continuous injection $\imath:I\to(\cT,d)$ of a real interval
  $I=[t_1,t_2]$ into $\cT$, and let $\tau_i=\imath(t_i)$. By
  declaring $\tau_1$ to be the root of $\cT$, we may suppose
  $\tau_1\le\tau_2$. Define $\pi(t)=\imath(t)\wedge\tau_2$, and
  suppose $\pi^{-1}\{t\}$ has an interior point for some $t\in I$. If
  $[a,b]$, $a<b$ denotes a non trivial connected component of
  $\pi^{-1}\{t\}$, the tree structure of $\cT$ implies
  $\imath(a)=\imath(b)$. This contradicts the injectivity of
  $\imath$. From the fact that the preimage by $\pi$ of any point in
  $[\tau_1,\tau_2]$ has empty interior, we infer that $\imath$
  maps $I$ into the segment $[\tau_1,\tau_2]$.  The injectivity of
  $\imath$ then gives that $\imath$ is an increasing homeomorphism of
  $I$ onto $[\tau_1,\tau_2]$.  This completes the proof.
\end{proof}  
If a metric tree $\cT$ is complete as a nonmetric tree, then $(\cT,d)$
is a complete metric space for every compatible tree metric $d$ on
$\cT$.  As a partial converse, if a nonmetric tree $\cT$ admits a
compatible, complete tree metric $d$ of finite diameter, then $\cT$ is
complete as a tree.  Notice, however, that $\R$ with its standard
metric is complete as a metric space but not as a nonmetric tree.

If a metric tree $\cT$ has finite diameter, then
the completion of $\cT$ as a metric space agrees with the completion
of $\cT$ as a nonmetric tree. Moreover, it is always possible to find 
an equivalent metric in which $\cT$ has finite diameter.
%
%
\subsection{Trees from ultrametric spaces}\label{def-ultra}
There is a natural way to construct trees from ultrametric spaces. 
\index{tree!from ultrametric space}
This procedure will be used to illustrate the connection 
between quasimonomial valuations and curves in Section~\ref{val-curve}.

Recall that a metric $d$ on a space $X$ is an \emph{ultrametric}
\index{ultrametric}
if it satisfies the stronger triangle inequality
$d(x,y)\le \max\{d(x,z),d(y,z)\}$
for any $x,y,z\in X$.

Let $(X,d_X)$ be an ultrametric space of diameter 1. 
Define an equivalence relation $\sim$ on $X\times(0,1)$ 
by declaring
$(x,s)\sim(y,t)$ iff $d(x,y)\le s=t$.  Note that
$(x,1)\sim(y,1)$ for any $x,y$. 
The set $\cT_X$ of equivalence classes is a nonmetric tree 
rooted at $(x,1)$ under the
partial order $(x,s)\le(y,t)$ iff $d(x,y)\le s\ge t$.
It is a metric tree under the
metric defined by $d((x,s),(x,t))=|s-t|$ and
$d(\sigma,\tau)=d(\sigma,\sigma\wedge\tau)+d(\tau,\sigma\wedge\tau)$
for general $\sigma,\tau\in\cT_X$. 
It can also be parameterized by declaring $\a(x,t)=t^{-1}$.

If the diameter of any ball of radius $r$ equals $r$ (!), then 
we may think of $(x,t)\in\cT_X$ as the closed ball of radius 
$t$ in $X$ centered at $x$.
%
%
\subsection{Trees from simplicial trees}\label{S303}
Next we show how the classical notion of simplicial tree fits into
our framework.

Recall that a \emph{simplicial tree}
\index{tree!simplicial}
is a set $V$ (vertices) together with a collection $E$
(edges) of subsets of $V$ of cardinality 2, such that the following
holds: for any two distinct vertices $\sigma,\tau\in V$ 
there exist a unique sequence 
$\sigma=\sigma_0,\sigma_1,\dots,\sigma_n=\tau$ of
distinct vertices such that $\{\sigma_{i-1},\sigma_i\}\in E$ for all 
$i=1,\dots,n$. We write $[\sigma,\tau]=\{\sigma_i\}_0^n$.
A \emph{rooted simplicial tree} is a triple
$(V,E,\tau_0)$ consisting of a simplicial tree $(V,E)$ 
together with a marked vertex $\tau_0\in E$ (the root).

If $(V,E,\tau_0)$ is a rooted, simplicial tree, then we can define
a partial ordering on $V$ by declaring $\tau\le\tau'$
iff $[\tau_0,\tau]\subset[\tau_0,\tau']$. It is straightforward
to verify that $(V,\le)$ is a rooted, nonmetric 
$\N$-tree in the sense of Section~\ref{S35}.
Conversely, to any rooted, nonmetric $\N$-tree
$(\cT,\le)$ we can associate a 
rooted, simplicial tree $(V,E,\tau_0)$ as follows:
$V=\cT$, $\tau_0$ is the root of $\cT$, and 
$\{\tau,\tau'\}\in E$ iff the segment $[\tau,\tau']$
in $\cT$ contains exactly two elements.
One easily verifies that $(V,E,\tau_0)$ is a rooted 
simplicial tree, and that the two operations just defined
are inverse to each other.
Moreover, changing the root of $(V,E,\tau_0)$ leads
to an equivalent partial ordering on $\cT=V$, and vice versa.
We thus conclude that
\emph{(rooted) simplicial trees can be 
  identified with (rooted) nonmetric $\N$-trees}.

Notice that any rooted, nonmetric $\N$-tree $(\cT,\le)$
has a unique parameterization $\a:\cT\to\N$ satisfying
$\a(\tau_0)=1$: this is given by
$\a(\tau)=\#\{\sigma\le\tau\}$.

To any rooted (nonmetric) $\N$-tree $(\cT_\N,\le)$ with root
$\tau_0$ we can associate a rooted, nonmetric $\R$-tree 
$(\cT_\R,\le)$ by ``adding the edges'' as follows: we set
\begin{equation*}
  \cT_\R=\{(\tau_0,0)\}\cup
  \left(\left(\cT_\N\setminus\{\tau_0\}\right)\times\,]-1,0]\right)
\end{equation*}
and define the partial ordering on $\cT_\R$ by lexicographic ordering:
$(\tau,s)\le(\tau',s')$ iff $\tau<\tau'$ or 
$\tau=\tau'$ and $s\le s'$.
It is straightforward to verify that $\cT_\R$ is a
nonmetric $\R$-tree rooted in $(\tau_0,0)$.
Notice that $\cT_\R$ naturally contains 
$\cT_\N\simeq\cT_\N\times\{0\}$ as a subset.
In fact, $\cT_\R$ is the minimal rooted, nonmetric $\R$-tree
with this property. We leave it to the reader to
make the last statement precise, as well as to verify that
equivalent partial orderings on $\cT_\N$ give rise
to equivalent partial orderings on $\cT_\R$.
Finally notice that the canonical parameterization
$\a:\cT_\N\to\N$ extends to a parameterization
$\a:\cT_\R\to\R_+$ by setting
$\a(\tau,s)=\a(\tau)+s$.
%
%
\subsection{Trees from $\Q$-trees}\label{sec-qtree}
As we show in this section, there is a natural way of passing from a
$\Q$-tree to an $\R$-tree by ``adding irrational points''. This
construction will be used in Section~\ref{sec-constr-univ}.

First, notice that definition of $\Q$-trees 
does not involve any arithmetic properties of the set $\Q$
but only its structure as a totally ordered set. 
As such, it is characterized by
\begin{lemma}\label{L421}
  A totally ordered set $\Lambda$ is isomorphic to an interval in 
  $\Q$ iff $\Lambda$ is countable and has no gaps in the sense that
  if $\lambda,\lambda'\in\Lambda$ and $\lambda<\lambda'$ then there 
  exists $\lambda''\in\Lambda$ with $\lambda<\lambda''<\lambda'$.
\end{lemma}
\begin{proof}
  Clearly any interval in $\Q$ has the stated property. 
  For the converse, fix a countable, totally ordered set $\Lambda$ 
  with no gaps. We wish to show it is isomorphic to an interval in
  $\Q$. There is no loss of generality in assuming that $\Lambda$
  contains its infimum and supremum.
  The proposition will follow if we can
  produce an isomorphism $\chi$ from $\Lambda$ onto the
  set $D$ of dyadic rational numbers in $[0,1]$, \ie 0,1 and
  all numbers of the form $i/2^j$, $j\ge 1$, $0<i<2^j$.
  Indeed, we can then post-compose with an isomorphism from
  $D$ onto $[0,1]\cap\Q$.

  Let $(\lambda_n)_0^\infty$ be an enumeration of $\Lambda$. 
  Assume that $\lambda_0=\min\Lambda$ and $\lambda_1=\max\Lambda$.
  Set $\chi(\lambda_0)=0$ and $\chi(\lambda_1)=1$. 
  Inductively, suppose $n\ge 2$ and that 
  we have defined $\chi(\lambda_m)$ for $m<n$ in an order-preserving way. 
  Define
  $\lambda'=\max\{\lambda_m\ ; \ m<n, \lambda_m<\lambda_n\}$
  and 
  $\lambda''=\min\{\lambda_m\ ; \ m<n, \lambda_m>\lambda_n\}$.
  Then pick $\chi(\lambda_n)$ to be a dyadic rational in 
  $]0,1[$ with minimal denominator such that
  $\chi(\lambda')<\chi(\lambda)<\chi(\lambda'')$.

  It is clear that this gives an order-preserving mapping of 
  $\Lambda$ into the set $D$. The fact that it is onto follows
  from the assumption that $\Lambda$ has no gaps.
\end{proof}
We now formulate the way in which a $\Q$-tree is canonically
embedded in an $\R$-tree.
\begin{proposition}\label{P420}
  Given a rooted, nonmetric $\Q$-tree $\cT_\Q$ there exists a
  rooted, nonmetric $\R$-tree $\cT_\R$ and an order-preserving 
  injection $\imath:\cT_\Q\to\cT_\R$ such that:
  \begin{itemize}
  \item[(i)]
    $\imath(\cT_\Q)$ is weakly dense in $\cT_\R$;
  \item[(ii)]
    every point in $\cT_\R\setminus\imath(\cT_\Q)$ is a
    regular point of $\cT_\R$;
  \item[(iii)]
    if $\imath':\cT_\Q\to\cT'_\R$ is an order-preserving injection
    into another rooted, nonmetric $\R$-tree with weakly dense image, 
    then there exists
    an injective morphism $\Phi:\cT_\R\to\cT'_\R$ of rooted,
    nonmetric trees such that $\Phi\circ\imath=\imath'$ and
    such that $\Phi$ extends to an isomorphism of rooted,
    nonmetric trees between the completions of $\cT_\R$
    and $\cT'_\R$.
  \end{itemize}
  Moreover, given a parameterization $\a_\Q:\cT_\Q\to\Q_+$ of $\cT_\Q$
  there exists a unique parameterization $\a_\R:\cT_\R\to\R_+$
  such that $\a_\R\circ\imath=\a_\Q$.
\end{proposition}
\begin{proof}
  We first show how to define $\cT_\R$ and $\imath$.
  Define $\cT_\R$ as the set of equivalence classes
  $[I_\bullet]$ of decreasing sequences 
  $I_\bullet=(I_n)_0^\infty$ of closed segments 
  in $\cT_\Q$ so that $\cap_{n\ge0} I_n$ contains at most one point.
  Here $I_\bullet$ and $J_\bullet$ are
  equivalent iff for every $n$ there exists $m$ such that
  $I_m\subset J_n$ and $J_m\subset I_n$. Further, 
  if $\sigma\in\cT_\Q$, then $\imath(\sigma)\in\cT_\R$ is
  defined to be the equivalence class containing 
  $I_\bullet=(I_n)_0^\infty$,
  where $I_n=[\sigma,\sigma]$ for all $n$.

  Let us show that $\cT_\R$ is naturally a rooted nonmetric
  $\R$-tree, \ie that it admits a natural partial ordering
  satisfying~(T1)-(T3).
  This partial ordering is defined as follows:
  $[I_\bullet]\le[J_\bullet]$ iff there exist increasing 
  sequences $(n_k)_1^\infty$, $(m_k)_1^\infty$
  and, for every $k$, elements $\sigma_k\in I_{n_k}$,
  $\tau_k\in J_{m_k}$ with $\sigma_k\le\tau_k$.
  We leave it to the reader to verify that this is well-defined
  and that $\imath:\cT_\Q\to\cT_\R$ is an order-preserving injection.
  If $\sigma_0$ is the root of $\cT_\Q$, then clearly
  $\imath(\sigma_0)$ is the unique minimal element of $\cT_\R$,
  so $\cT_\R$ satisfies~(T1). 
  We now consider~(T2). 
  We have to show that if $[I_\bullet]\in\cT_\R$, then 
  $\{[J_\bullet]\in\cT_\R\ ;\ [J_\bullet]\le[I_\bullet]\}$
  is a totally ordered set isomorphic to a real interval.
  First suppose $[I_\bullet]=\imath(\sigma)$ for some
  $\sigma\in\cT_\Q$. The segment $[\sigma_0,\sigma]$ is a 
  totally ordered set isomorphic to $[0,1]\cap\Q$. 
  But it is well-known that if we perform the 
  construction above on the nonmetric $\Q$-tree $[0,1]\cap\Q$,
  then we end up with $[0,1]\cap\R$. 
  Hence we are done in this case. 
  If $[I_\bullet]$ is not of the form $\imath(\sigma)$, then 
  there still exists a decreasing sequence $\sigma_n$ in
  $\cT_\Q$ such that $\imath(\sigma_n)$ decreases to $[I_\bullet]$. 
  Then $\{[J_\bullet]\le[I_\bullet]\}$
  is isomorphic to the intersection of real intervals 
  $[0,s_n]\subset\R$, where $(s_n)_1^\infty$ is a 
  decreasing sequence of real numbers with 
  $s_\infty=\lim s_n>0$. As $\{[J_\bullet]\le[I_\bullet]\}$
  has a maximal element, it must be isomorphic 
  to the real interval $[0,s_\infty]$.

  Thus~(T2) holds. Instead of proving~(T3) we prove
  the equivalent statement~(T3') in Remark~\ref{R401}.
  Consider a totally ordered subset $\cS_\R\subset\cT_\R$ 
  without upper bound in $\cT_\R$. 
  We may assume that $\cS_\R$ is full. 
  Then $\cS_\Q:=\imath^{-1}(\cS_\R)$ is a nonempty, totally
  ordered subset of $\cT_\Q$ without upper bound. 
  Since $\cT_\Q$ is a $\Q$-tree, there exists an increasing
  sequence $(\sigma_n)$ in $\cS_\Q$ without upper bound in
  $\cT_\Q$. 
  Then it is easy to see that $(\imath(\sigma_n))$
  is an increasing sequence in $\cS_\R$ without upper bound
  in $\cT_\R$.

  We conclude that $\cT_\R$ is an $\R$-tree. It follows from the
  construction that $\gamma\cap\imath(\cT_\Q)$ is weakly dense in
  $\gamma$ for every segment $\gamma\subset\cT_\R$, where $\gamma$ is
  equipped with the topology induced from $\R$.  Hence
  $\imath(\cT_\Q)$ is weakly dense in $\cT_\R$.  That all points of
  $\cT_\R\setminus\imath(\cT_\Q)$ are regular points of $\cT_\R$ is
  clear from the construction.  Thus we have proved~(i) and~(ii).

  Finally, suppose $\imath'$ is 
  an order-preserving injection of $\cT_\Q$ into
  an $\R$-tree $\cT'_\R$ with weakly dense image.
  Let us construct the mapping $\Phi:\cT_\R\to\cT'_\R$
  as in~(iii). 
  The construction is based on the following result,
  the proof of which is left to the reader:
  \begin{lemma}\label{L422}
    If $X$ and $Y$ are countable dense subsets of the 
    real interval $[0,1]$ and 
    $\chi:X\to Y$ is an order-preserving bijection, 
    then $\chi$ extends uniquely to an 
    order-preserving bijection (\ie an increasing
    homeomorphism) of $[0,1]$ onto itself.
  \end{lemma}
  We continue the proof of the proposition. 
  The fact that $\imath'(\cT_\Q)$ is weakly dense in
  $\cT'_\R$ implies that the root of $\cT'_\R$ is 
  $\imath'(\sigma_0)$, and that 
  $\imath'([\sigma_0,\sigma])$ is dense in
  the segment $[\imath'(\sigma_0),\imath(\sigma)]$ in $\cT'_\R$
  for any $\sigma\in\cT_\Q$ 
  (this uses the fact that $\cT_\Q$ is a $\Q$-tree
  and does not follow from the fact that 
  $\imath'(\cT_\Q)$ is a weakly dense subset of $\cT'_\R$).
  Lemma~\ref{L422} then implies that for any $\sigma\in\cT_\Q$
  there is a unique order-preserving bijection 
  $\Phi_\sigma:[\imath(\sigma_0),\imath(\sigma)]
  \to[\imath'(\sigma_0),\imath'(\sigma)]$ 
  such that $\Phi_\sigma\circ\imath=\imath'$ on $[\sigma_0,\sigma]$.
  These bijections patch together to form an order-preserving 
  injection $\Phi:\cT_\R\to\cT'_\R$ with 
  $\Phi\circ\imath=\imath'$ on $\cT_\Q$. 
  Since $\imath'(\cT_\Q)$ is dense in 
  $\cT'_\R$, $\Phi$ extends to an isomorphism of rooted trees
  between the completions of $\cT_\R$ and $\cT'_\R$, 
  completing the proof of~(iii).

  Finally, the fact that any parameterization of $\cT_\Q$
  extends uniquely to a parameterization of $\cT_\R$ is
  again a consequence of Lemma~\ref{L422}.
\end{proof}
%
%
%
%
\section{Nonmetric tree structure on $\cV$}\label{tree-struc}
We now show that the natural partial ordering on valuation space
$\cV$ turns it into a rooted, nonmetric tree, the \emph{valuative tree},
\index{valuative tree}
\index{valuation space}
\index{$\cV$ (valuation space/valuative tree)}
and describe its structure.
%
%
\subsection{Partial ordering.}
Recall that we have normalized the valuations in $\cV$ by 
$\nu(\fm)=1$ and that the partial ordering $\le$ on $\cV$
\index{partial ordering!on the valuative tree} 
is given by $\nu\le\mu$ iff 
$\nu(\phi)\le\mu(\phi)$ for all $\phi\in R$.
The multiplicity valuation is given by $\nu_\fm(\phi)=m(\phi)$.
\begin{theorem}\label{order}
  Valuation space $\cV$ is a complete nonmetric tree rooted at
  $\nu_\fm$\index{valuative tree}.
\end{theorem}
The general properties of nonmetric trees then give:
\begin{corollary}\label{mini}
  Any subset of $\cV$ admits an infimum. 
\end{corollary}
\index{$\wedge$ (infimum)}
\index{infimum!in $\cV$}
\begin{remark}\label{dim3}
  If $(\nu_i)$ are valuations in $\cV$, then the infimum
  $\nu:=\wedge \nu_i\in\cV$ can be constructed using the 
  following properties:
  $\nu(\phi)=\inf_i \nu_i(\phi)$ for all irreducible $\phi\in\fm$ 
  and 
  $\nu(\phi\psi)=\nu(\phi)+\nu(\psi)$ for all $\phi,\psi\in\fm$.

  This construction does not work on more general rings. For instance,
  consider monomial valuations $\nu_i$, $i=1,2,3$ on $\C[[x_1,x_2,x_3]]$
  defined by $\nu_i(x_j)=3$ if $i\ne j$ and $\nu_i(x_i)=1$. 
  Define $\nu=\min_i\nu_i$ by the construction above and consider
  $\phi=x_1x_2-x_2x_3+x_3x_1$,
  $\psi=x_1x_2+x_2x_3-x_3x_1$.
  Then $\nu(\phi)=\nu(\psi)=4$ but $\nu(\phi+\psi)=2$, so
  $\nu$ is not a valuation.

  A similar calculation shows that the natural partial 
  ordering on the
  set of normalized valuations on $\C[[x_1,x_2,x_3]]$ 
  does \emph{not} define a tree structure.
\end{remark}
Any invertible formal mapping $f:(\C^2,0)\to(\C^2,0)$ induces a
ring automorphism $f^*:R\to R$. Since $f^*\fm=\fm$, we have an induced
mapping $f_*:\cV\to\cV$ given by $f_*(\nu)(\phi)=\nu(f^*\phi)$. 
If $\mu\le\nu$, then clearly $f_*\mu\le f_*\nu$. Hence we get:
\begin{proposition}\label{P106}
  Any invertible formal mapping $f:(\C^2,0)\to(\C^2,0)$ 
  induces an isomorphism $f_*:\cV\to\cV$ of rooted, nonmetric trees.
\end{proposition}
We now turn to the proof of Theorem~\ref{order}. 
It is proved by describing the partial ordering 
on $\cV$ in terms of SKP's:
\begin{proposition}\label{pro-comp}
  Let $\nu$ and $\nu'$ be valuations in $\cV$, $\nu\ne\nu'$.
  Pick local coordinates $(x,y)$ such that
  $1=\nu(x)=\nu'(x)\le\min\{\nu(y),\nu'(y)\}$.
  Write $\nu=\val[(U_j);(\btilde_j)]$ and
  $\nu'=\val[(U'_j);(\btilde'_j)]$.  
  Then $\nu<\nu'$ iff 
  \begin{equation*}
    \lgt(\nu')\ge\lgt(\nu)=:k<\infty,
    \quad U_j=U'_j\ \text{for}\ 0\le j\le k 
    \qand \btilde'_k\ge\btilde_k.
  \end{equation*} 
\end{proposition}
As a direct consequence, the infimum of any family of normalized valuations
exists and is computable in terms of SKP's, at least as long as we
can choose the coordinates $(x,y)$ conveniently.
\begin{corollary}\label{C310}
  Let $(\nu^i)_{i\in I}$ be a family of valuations in $\cV$
  and suppose we can  find local coordinates $(x,y)$ such that
  $1=\nu^i(x)\le\nu^i(y)$ for all $i\in I$.
  Write $\nu^i=\val[(U^i_j);(\btilde^i_j)]$.
  Then the infimum of the $\nu^i$ exists and is
  given in terms of SKP's by
  \begin{equation*}
    \bigwedge_{i\in I}\nu^i
    =\val[(U_j)_0^k;(\btilde_j)_0^{k-1},\inf_{i\in I}\btilde^i_k],
  \end{equation*}
  where $1\le k<\infty$ is maximal such that 
  $U^i_j=U^{i'}_j=:U_j$ for $0\le j\le k$ and $i,i'\in I$, and
  $\btilde^i_j=\btilde^{i'}_j=:\btilde_j$ for $0\le j<k$ and $i,i'\in I$.
\end{corollary}
\begin{proof}[Proof of Proposition~\ref{pro-comp}]
  First suppose the three displayed conditions hold.  Then
  property~(Q2) of Theorem~\ref{key} implies
  $\nu'\ge\nu_0=\val[(U_j)_0^k;(\btilde_j)_0^{k-1},\btilde'_k]$.  But
  $\nu_0\ge\nu$, so $\nu'\ge\nu$.

  Conversely assume 
  $\val[(U'_j);(\btilde'_j)]=\nu'>\nu=\val[(U_j)_0^k;(\btilde_j)_0^k]$
  with $1\le k\le\infty$. 
  Let us show inductively that 
  $U_j'=U_j$ for $j\le k$. This is true by definition
  for $j = 0,1$.  Assume we proved it for $j<k$. 
  Define $\nu_{j-1}=\val[(U_l)_0^{j-1};(\btilde_l)_0^{j-1}]$.  
  Then $\nu'(U_j)\ge\nu(U_j)>\nu_{j-1}(U_j)$. 
  Hence Theorem~\ref{approx} implies $U'_j=U_j$. 
  Finally, $k<\infty$ (or else $\nu=\nu'$) 
  and $\btilde_k'=\nu'(U_k)\ge\nu(U_k)=\btilde_k$.
  \end{proof}
\begin{proof}[Proof of Theorem~\ref{order}]
  It is clear that $(\cV,\le)$ is a partially ordered set
  with unique minimal element $\nu_\fm$. Thus~(T1) holds.

  To prove~(T2), fix $\nu\in\cV$ with $\nu>\nu_\fm$. 
  We will show that the set 
  $I=\{\mu\ ;\ \nu_\fm\le\mu\le\nu\}$ is a totally 
  ordered set isomorphic to an interval in $\Rbar$. 

  We may pick local coordinates $(x,y)$ such that $1=\nu(x)\le\nu(y)$.
  Write $\nu=\val[(U_j)_0^k;(\btilde_j)_0^k]$.  First assume
  $k<\infty$.  Set $d_j=\deg_y(U_j)$ for $1\le j\le k$ and
  $d_0=\infty$ by convention.  Recall that the sequence
  $(\btilde_j/d_j)$ is strictly increasing (see
  Lemma~\ref{lem:deg-irr}).  We claim that $I$ is isomorphic to the
  interval $J=[1,\btilde_k/d_k]$.  To see this, pick $t\in J$. There
  exists a unique integer $l\in[1,k]$ such that
  $\btilde_{l-1}/d_{l-1}<t\le\btilde_l/d_l$.  Set
  $\nu_t=\val[(U_j)_0^l;(\btilde_j)_{j=0}^{l-1},td_l]$.
  Proposition~\ref{pro-comp} then shows that this gives an isomorphism
  from $J$ onto $I$.  The case $k=\infty$ is treated in a similar way.

  As for~(T3), it is easier to prove the equivalent statement~(T3') 
  given in Remark~\ref{R401}. Moreover~(T3') clearly follows if
  we can prove that every totally ordered subset of $\cV$
  has a majorant in $\cV$. This will in fact also prove that $\cV$
  is a complete tree. 
  Thus consider such a totally ordered subset $\cS\subset\cV$.
  We may pick local coordinates $(x,y)$ 
  such that $1=\nu(x)\le\nu(y)$ for every $\nu\in\cS$.
  By Proposition~\ref{pro-comp} the SKP
  defining $\nu\in\cS$ has a length that is a nondecreasing function
  of $\nu$. When $\lgt\nu\to\infty$,
  Theorem~\ref{key-infinite} shows that
  $\nu$ tends to a curve valuation 
  or an infinitely singular valuation
  dominating all the $\nu$'s. Otherwise
  $\lgt\nu$ is constant for large $\nu$, $\nu\in\cS$,
  and we can write 
  $\nu=\val[(U_j)_0^n;(\btilde_j)_0^{n-1},\btilde_n^{(\nu)}]$. 
  By Proposition~\ref{pro-comp}, $\btilde_n^{(\nu)}$ is an increasing
  function of $\nu$, hence is dominated by some $\btilde_n\in\Rbar$. 
  The valuation $\mu=\val[(U_j)_0^n;(\btilde_j)_0^n]$ 
  dominates all the $\nu$'s. Thus $\cV$ is a complete tree.
  \end{proof}
%
%
\subsection{Dendrology.}
We now undertake a more detailed study of the nonmetric tree structure 
on valuation space. 
Recall from Section~\ref{S26} that we have classified the valuations in 
$\cV$ into four categories: divisorial, irrational, infinitely singular 
and curve valuations. Moreover, we called a valuation quasimonomial
if it is either divisorial or irrational; 
see Definition~\ref{D101}.

This classification was defined in terms of SKP's 
and a priori depended on a choice of local coordinates. 
We saw in Theorem~\ref{divis} that the classification 
could also be formulated in terms of numerical invariants 
and in particular did not depend on coordinates. 
Here we show how the nonmetric tree structure interacts with
the classification.
\begin{proposition}\label{P107}
  The rooted, nonmetric tree structure on valuation space $\cV$ 
  has the following properties:
  \begin{itemize}
  \item[(i)]
    the root of $\cV$ is the multiplicity valuation $\nu_\fm$;
  \item[(ii)]
    the ends of $\cV$ are the infinitely singular and curve valuations;
  \item[(iii)]
    any tangent vector in $\cV$ is represented by
    a curve valuation as well as by an infinitely singular
    valuation;
  \item[(iv)]
    the regular points of $\cV$ are the irrational valuations;
  \item[(v)]
    the branch points of $\cV$ are the divisorial valuations. 
    Further, the tangent space at a divisorial valuation is in bijection 
    with $\P^1$.   
  \end{itemize}
\end{proposition}
\begin{remark}
  The bijection in~(v) will be made much more precise 
  in Appendix~\ref{sec-tangent}. 
  In particular we will show that each
  tangent vector has an interpretation both as a point on a rational
  curve and as a Krull valuation.
\end{remark}
\begin{remark}
  The nonmetric tree structure on $\cV$ does not allow us to 
  distinguish between a curve valuation and an infinitely singular
  valuation. 
  In Section~\ref{sec-mult} we will define the \emph{multiplicity} of 
  a valuation. This multiplicity is an increasing function on 
  $\cV$ with values in $\Nbar$ 
  and the infinitely singular valuations are the ones with 
  infinite multiplicity.
\end{remark} 
\begin{proof}[Proof of Proposition~\ref{P107}]
  We have already proved~(i).
  For~(ii), let $\nu\in\cV$ be quasimonomial
  and pick local coordinates $(x,y)$.
  Then $\nu=\val[(U_j)_0^k;(\btilde_j)_0^k]$, with 
  $k<\infty$ and $\btilde_k<\infty$.
  By Proposition~\ref{pro-comp}, $\nu$ is dominated by
  any valuation of the form
  $\val[(U_j)_0^k;(\btilde_j)_0^{k-1},\btilde'_k]$, where
  $\btilde'_k>\btilde_k$.
  Thus no quasimonomial valuation is an end in $\cV$. 
  The same proposition also shows that
  no valuation with an SKP of infinite length or with
  length $k<\infty$ and $\btilde_k=\infty$ can be dominated by
  another valuation. Thus curve valuations and infinitely singular
  valuations are ends in $\cV$, proving~(ii). 

  The proof of~(iii) is similar. 
  Indeed, consider $\nu\in\cV$ and
  a tangent vector $\vv$ at $\nu$. 
  First assume $\vv$ is not represented by $\nu_\fm$.
  By~(ii), $\nu$ is then quasimonomial and $\vv$ is 
  represented by another quasimonomial valuation $\mu$.
  Pick local coordinates $(x,y)$ and write 
  $\mu=\val[(U_j)_0^k;(\btilde_j)_0^k]$, where
  $k,\btilde_k<\infty$.
  The curve valuation 
  $\val[(U_j)_0^k;(\btilde_j)_0^{k-1},\infty]$
  then dominates $\mu$ and hence represents $\vv$.
  We can also construct an infinitely singular valuation
  dominating $\mu$: pick $\btilde'_k\in\Q$ with
  $\btilde'_k\ge\btilde_k$ and inductively extend
  the SKP $[(U_j)_0^k;(\btilde_j)_0^{k-1},\btilde'_k]$
  to an infinite SKP 
  $[(U_j)_0^\infty;(\btilde'_j)_0^\infty]$ where
  $\btilde'_j=\btilde_j$ for $0\le j<k$ and such
  that $d_j:=\deg_y(U_j)\to\infty$.
  The valuation
  $\val[(U_j)_0^\infty;(\btilde'_j)_0^\infty]$ is infinitely
  singular and dominates $\mu$; hence it represents $\vv$.

  If $\vv$ if represented by $\nu_\fm$, then $\nu\ne\nu_\fm$.
  Assume the local coordinates have been picked so that 
  $\nu(y)\ge\nu(x)=1$. 
  By Corollary~\ref{C310} the curve valuation $\nu_x$ 
  with SKP $[(x,y);(\infty,1)]$
  represents a tangent vector $\ww$ at $\nu_\fm$ 
  different from the one represented by $\nu$.
  From the preceding argument, $\ww$ is represented
  by both a curve valuation (\eg $\nu_y$) and by
  an infinitely singular valuation.
  These valuations then also represent the tangent vector $\vv$.
  Thus we have proved~(iii).

  Next consider a quasimonomial valuation $\nu$ and
  a tangent vector
  $\vv\in T\nu$ represented by a curve valuation $\nu_C$.
  First assume that $\nu_C\wedge\nu<\nu$. 
  Then the segment $]\nu,\nu_C]$ intersects
  $]\nu,\nu_\fm]$ at $\nu_C\wedge\nu$ so $\vv$ is the
  tangent vector represented by $\nu_\fm$.

  Otherwise, $\nu_C>\nu$. 
  First suppose $\nu$ is irrational. 
  We claim that all curve valuations $\nu_C$ with $\nu_C>\nu$
  represent the same tangent vector. 
  Indeed, in general, if $\nu_C$ and $\nu_D$ are distinct curve valuations,
  it follows from Corollary~\ref{C310} that
  $\nu_C\wedge\nu_D$ is divisorial. 
  Thus if $\nu_C,\nu_D>\nu$ and $\nu$ is irrational, then 
  $\nu_C\wedge\nu_D>\nu$ so $\nu_C$ and $\nu_D$ represent the
  same tangent vector at $\nu$. This proves~(iv).

  When $\nu$ is divisorial, write $\nu=\val[(U_j)_0^k,(\btilde_j)_0^k]$
  with $k<\infty$ and $\btilde_k\in\Q$.
  For $\mu >\nu$, either $\mu(U_k)>\nu(U_k)=\btilde_k$, 
  or $\mu(U_k)=\btilde_k$. In the former case, $\mu$ 
  represents the same vector as 
  $\nu_\infty=\val[(U_j)_0^k,(\btilde_j)_0^{k-1},\infty]$.
  In the latter case, Theorem~\ref{key}~(Q2) shows that there exists
  a polynomial 
  $U_{k+1}=U_k^{n_k}-\theta\prod_0^{k-1}U_j^{m_{kj}}$
  with $\theta=\theta(\mu)\in\C^*$, such that
  if $\phi\in R$ then 
  $\mu(\phi)>\nu(\phi)$ 
  iff $U_{k+1}$ divides $\phi$ in $\gr_\mu\C(x)[y]$. 
  Define $\nu_{\theta}=\nu_{U_{k+1}}$
  and $\btilde_{k+1}=\mu(U_{k+1})$. Then
  $\mu\wedge\nu_{\theta}\ge\val[(U_j)_0^{k+1},(\btilde_j)_0^{k+1}]>\nu$
  so that $\mu$ and $\nu_{\theta}$ define the same tangent vector. 
  Conversely $\nu_{\theta}\wedge\nu_{\theta'}=\nu$ 
  as soon as $\theta \not= \theta'$.
  The set of tangent vectors at $\nu$ is hence in bijection with
  $\{\nu_\fm,\nu_\infty,\nu_{\theta}\}_{\theta\in\C^*}$ which is
  in bijection with $\P^1$, proving~(v).
\end{proof}
\begin{figure}
  \begin{center}
    \includegraphics[width=\textwidth]{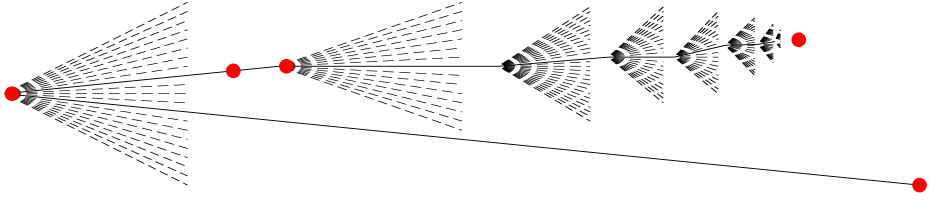}
  \end{center}
  \caption{The nonmetric tree structure on the valuative tree
    $\cV$. The valuations marked by dots are, from left to right: 
    the multiplicity valuation $\nu_\fm$,
    an irrational valuation, a divisorial valuation, 
    an infinitely singular valuation and a curve valuation.}\label{F101}
\end{figure}
%
%
\subsection{A model tree for $\cV$}\label{model}
The identification of valuations with SKP's gives an explicit model
for the tree structure on $\cV$. 
Namely, let $\cT$ be the set consisting of $(0)$ and of 
pairs $(\bs,\btheta)$ with
$\bs=(s_1,\dots,s_{m+1})$, $\btheta=(\theta_1,\dots,\theta_m)$, 
$0\le m\le\infty$, $s_j\in \Q_+^*$ for $1\le j<m+1$, 
$s_{m+1}\in(0,\infty]$ if $m<\infty$, $\theta_j\in\C^*$. 
Define a partial ordering on $\cT$ by 
$(\bs,\btheta)\le(\bs',\btheta')$ iff $m\le m'$, $s_j=s'_j$ and
$\theta_j=\theta'_j$ for $j<m$ and $s_m\le s'_m$, 
and by declaring $(0)$ to be the unique minimal element.
Then $(\cT,\le)$ is a complete nonmetric tree rooted in $(0)$.
  
Define $\imath:\cT\to\cV$ as follows. Fix local
coordinates $(x,y)$.
First set $\imath(0)=\nu_x$.
If $(\bs,\btheta)\in\cT$, then set $\btilde_0=1/\min\{1,s_1\}$,
$\btilde_1=s_1/\min\{1,s_1\}$, and define, inductively, 
$n_k=\min\{n\ ;\ n\btilde_k\in\sum_0^{k-1}\btilde_j\Z\}$ and
$\btilde_{k+1}\=n_k\btilde_k+s_{k+1}$.  Define $U_k$ by~(P2) and set
$\imath(\bs,\btheta)\=\val[(U_k);(\btilde_k)]$.
Proposition \ref{pro-comp}
shows that $\imath$ is an isomorphism of nonmetric trees.

The map $\imath$ is in fact an isomorphism of
\emph{rooted} nonmetric trees if we take
$\nu_x$ to be the root of $\cV$. 
We shall study this situation more
closely in Section~\ref{sec-relative}.

%
%
%
%
\section{Parameterization of $\cV$ by skewness}\label{metric-tree-struc}
While the nonmetric tree structure on $\cV$ induced by the 
partial ordering is quite appealing,
it only reflects some of the features of the valuative tree.
For applications it is crucial to \emph{parameterize} $\cV$.
As we will see, there are two canonical
parameterizations. Here we will discuss the parameterization
by skewness. The second one---thinness---will be introduced 
in Section~\ref{S10} after we have defined the concept of
the multiplicity of a valuation.
%
%
\subsection{Skewness}
We first introduce a new invariant of a valuation.
It measures in an intrinsic way how far the valuation
is from the multiplicity valuation $\nu_\fm$.
\begin{definition}
  For $\nu\in\cV$, define the \emph{skewness}\index{skewness}
 $\a(\nu)\in[1,\infty]$ by
  \begin{equation}\label{e504}
    \a(\nu)\=\sup\left\{\frac{\nu(\phi)}{m(\phi)}
      \ ;\ \phi\in\fm\right\}.
  \end{equation}
\end{definition}
This quantity is an invariant of a valuation, in the following sense:
\begin{proposition}\label{P502}
  If $f:(\C^2,0)\to(\C^2,0)$ is an invertible formal map, and
  $f_*:\cV\to\cV$ the induced map, 
  then $\alpha(f_*\nu)=\alpha(\nu)$ for any $\nu\in\cV$.
\end{proposition}
\begin{proof}
  This is immediate since
  \begin{equation*}
    \alpha(f_*\nu)
    =\sup_{\phi\in\fm}\frac{\nu(f^*\phi)}{m(f^*\phi)}
    =\sup_{\psi\in\fm}\frac{\nu(\psi)}{m(\psi)}
    =\alpha(\nu).
  \end{equation*}
\end{proof}
A main reason why skewness is useful is the following result that
will be used repeatedly in the sequel.
\begin{proposition}\label{P201}
  For any valuation $\nu\in\cV$ and any irreducible $\phi\in\fm$ we have
  \begin{equation*}
    \nu(\phi)=\alpha(\nu\wedge\nu_\phi)m(\phi).
  \end{equation*}
  In particular $\nu(\phi)\le\alpha(\nu)m(\phi)$ with equality
  iff $\nu_\phi\ge\nu$.
\end{proposition}
\begin{figure}
  \begin{center}
    \includegraphics[width=\textwidth]{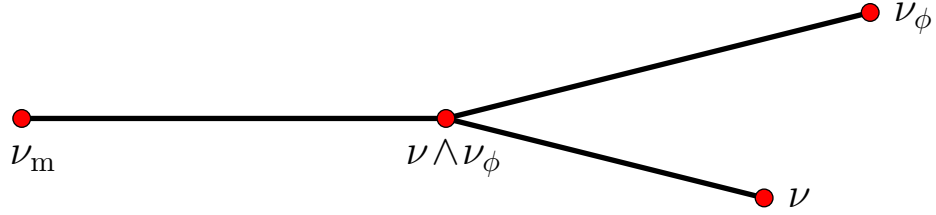}
  \end{center}
  \caption{The value $\nu(\phi)$ for $\phi\in\fm$ irreducible 
    depends on the multiplicity $m(\phi)$ and the relative 
    position of $\nu$ and the curve valuation $\nu_\phi$
    in the valuative tree $\cV$. 
    See Proposition~\ref{P201}.}\label{F4}
\end{figure}
%
%
\subsection{Parameterization}
\begin{figure}
  \begin{center}
    \includegraphics[width=\textwidth]{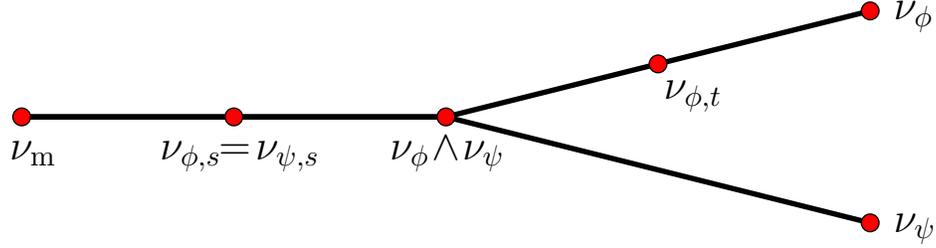}
  \end{center}
  \caption{Skewness is used to parameterize the segment
    $[\nu_\fm,\nu_\phi]$ by $\nu_{\phi,t}$, $1\le t\le\infty$.
    We have $\nu_{\phi,t}=\nu_{\psi,s}$ iff
    $s=t\le\alpha(\nu_\phi\wedge\nu_\psi)$.
    See Theorem~\ref{T401} and Definition~\ref{D401}.}\label{F5}
\end{figure}

As the following result asserts, skewness provides a good
parameterization of the valuative 
tree\index{parameterization!of the valuative tree (by skewness)}.
\begin{theorem}\label{T401}
  Skewness defines a parameterization $\a:\cV\to[1,\infty]$
  of the valuative tree $\cV$ rooted in the multiplicity
  valuation $\nu_\fm$. Moreover:
  \begin{itemize}
  \item[(i)]
    if $\nu$ is a divisorial valuation, then $\a(\nu)$ is rational;
  \item[(ii)]
    if $\nu$ is an irrational valuation, then $\a(\nu)$ is irrational;
  \item[(iii)]
    if $\nu$ is a curve valuation, then $\a(\nu)=\infty$;
  \item[(iv)]
    if $\nu$ is infinitely singular, then $\a(\nu)\in(1,\infty]$.
  \end{itemize}
\end{theorem}
It follows that skewness also defines a 
parameterization $\a:\cVqm\to[1,\infty)$.
We refer to Appendix~\ref{sec-inf-sing} for a construction of an
infinitely singular valuation with prescribed skewness
$t\in(1,\infty]$.
\begin{definition}\label{D401}
  Fix an irreducible curve $C$ and pick $t\in[1,\infty]$.
  We denote by $\nu_{C,t}$ the unique valuation 
  in the segment $[\nu_\fm,\nu_C]$ having skewness $\a(\nu)=t$.
  If $C=\{\phi=0\}$ for $\phi\in\fm$ irreducible, then we also 
  write $\nu_{\phi,t}=\nu_{C,t}$.
  See Figure~\ref{F5}.
\end{definition}
\begin{remark}
  The parameterization of $\cV$ by skewness does not allow us to 
  distinguish between curve valuations associated to smooth and 
  singular curves since both have infinite skewness, as do some 
  infinitely singular valuations. The distinction will be made in
  Section~\ref{sec-mult} by the multiplicity function on $\cV$.
\end{remark}
Before proving the theorem, let us note two additional results.
\begin{corollary}\label{C502}
  If $\nu\in\cV$ and $\phi\in\fm$ is irreducible, then 
  $\nu\ge\nu_{\phi,t}$, where $t=\nu(\phi)/m(\phi)$.
\end{corollary}
\begin{proof}
  By Proposition~\ref{P201} we have $t=\a(\nu\wedge\nu_\phi)$.
  Since $\nu\wedge\nu_\phi\in[\nu_\fm,\nu_\phi]$ we must have
  $\nu\wedge\nu_\phi=\nu_{\phi,t}$. Thus $\nu\ge\nu_{\phi,t}$.
\end{proof}
\begin{corollary}\label{C501}
  If $\nu\in\cV$, $\phi\in\fm$ is irreducible and
  $\nu(\phi)$ is irrational, then $\nu=\nu_{\phi,t}$ where
  $t=\nu(\phi)/m(\phi)$.
\end{corollary}
\begin{proof}
  Define $t=\nu(\phi)/m(\phi)$. By Corollary~\ref{C502} 
  we have $\nu\ge\nu_{\phi,t}$. Since $\nu_{\phi,t}$ is
  irrational, either equality holds, or else $\nu\ge\nu_{\phi,s}$
  for some $s>t$, in which case $\nu(\phi)\ge sm(\phi)>tm(\phi)$,
  a contradiction.
\end{proof}                                
\begin{proposition}\label{P303}
  Skewness $\a:\cV\to[1,\infty]$ is
  lower semicontinuous, but not continuous,
  in the weak tree topology on $\cV$.
\end{proposition}
\begin{proof}
  Lower semicontinuity follows from Proposition~\ref{P305} and the 
  fact that skewness defines a parameterization of $\cV$.
  Consider the sequence $\nu_n=\nu_{y-nx,2}$ in $\cV$. We have
  $\a(\nu_n)=2$ for all $n$, but $\nu_n\to\nu_\fm$ and $\a(\nu_\fm)=1$.
  Thus skewness is not weakly continuous.
\end{proof}
%
%
\subsection{Proofs}
The proofs of all the statements about skewness go by translating
them into properties of SKP's. Hence we first prove
\begin{lemma}\label{compuskew}
  Let $\nu\in\cV$ and 
  write $\nu=\val[(U_j)_0^k;(\btilde_j)_0^k]$
  in coordinates $(x,y)$.
  Then the following hold:
  \begin{itemize}
  \item[(i)]
    if $\nu$ is quasimonomial, then
    $\a(\nu)=d_k^{-1}\btilde_k\btilde_0$;
    where $d_k=\deg_yU_k$;
    in particular $\a(\nu)$ is rational iff $\nu$ is divisorial;
  \item[(ii)]
    if $\nu$ is a curve valuation, then $\a(\nu)=\infty$;
  \item[(iii)]
    if $\nu$ is infinitely singular, then
    $\a(\nu)=\lim_{j\to\infty}d_j^{-1}\btilde_j\btilde_0\in(1,\infty]$.
  \end{itemize}
\end{lemma}
\begin{proof}
  If $\nu=\nu_\phi$ is a curve valuation, then
  $\alpha(\nu)\ge\nu(\phi)/m(\phi)=\infty$, 
  which proves~(ii). 
  As for~(i) and~(iii) 
  note that it suffices to use $\phi$ irreducible 
  in~\eqref{e504}. Indeed, if $\nu(\phi)\le\alpha\ m(\phi)$
  for $\phi\in\fm$ irreducible and 
  $\phi=\prod\phi_i$ is reducible, then 
  \begin{equation*}
    \nu(\phi)
    =\sum_i\nu(\phi_i)
    \le\alpha\sum_i m(\phi_i)
    =\alpha\,m(\phi).
  \end{equation*}
  So assume $\phi$ irreducible and 
  write $\nu_\phi=\val [(U^\phi_j);(\btilde^\phi_j)]$.
  Let $l=\con(\nu,\nu_\phi)$ as defined in \eqref{contact}.
  Then $l\le k$. Recall that the sequence $(d_j^{-1}\btilde_j)_1^k$
  is increasing. We apply Proposition \ref{vvm}:
  \begin{equation}\label{e505}
    \frac{\nu(\phi)}{m(\phi)} 
    =d_l^{-1}\min\{\btilde_l,\btilde_l^\phi\}
    \min\{\btilde_0,\btilde^\phi_0\}
    \le\sup_{j\le k}d_j^{-1}\btilde_j\btilde_0.
  \end{equation}
  If $\nu$ is quasimonomial, \ie $k<\infty$, then
  equality holds in~\eqref{e505} if $l=k$ and 
  $\btilde^\phi_k=\btilde_k$, proving~(i).
  
  If $\nu$ is infinitely singular, so that $k<\infty$, then 
  for any $j$ we may pick $\phi$ with $l=j$ and
  $\btilde^\phi_l=\btilde_l$. Then 
  $\nu(\phi)/m(\phi)\ge d_j^{-1}\btilde_j\btilde_0$. 
  Letting $j\to\infty$ yields~(iii).
\end{proof}
\begin{proof}[Proof of Theorem~\ref{T401}]
  Assertions~(i)-(iv) follow from Lemma~\ref{compuskew}.  To show that
  $\a$ defines a parameterization it suffices to show that if $\nu$ is
  an end in $\cV$, \ie a curve or infinitely singular valuation, then
  $\a$ restricts to an increasing mapping of $[\nu_\fm,\nu[$ onto
  $[1,\a(\nu)[$.

  Pick local coordinates $(x,y)$ such that 
  $\nu(x)=1\le\nu(y)$.
  In terms of SKP's, write $\nu=\val[(U_j)_0^k;(\btilde_j)_0^k]$. 
  Here either $k=\infty$ or $\btilde_k=\infty$.
  Let $d_l=\deg_yU_l$. Recall that the sequence 
  $(\btilde_l/d_l)_1^k$ is strictly increasing.
  By Proposition~\ref{pro-comp} any valuation $\mu\in]\nu_\fm,\nu[$
  can be written (uniquely) in the form
  \begin{equation*}
    \mu=\val[(U_j)_0^l;(\btilde_j)_0^{l-1},\btilde]
  \end{equation*} 
  where $1\le l\le\min\{k,\infty\}$, 
  $\btilde<\infty$ if $l=k<\infty$, and where
  $d_{l-1}^{-1}d_l\btilde_{l-1}<\btilde\le\btilde_l$.
  By Lemma~\ref{compuskew} we have $\a(\mu)=\btilde/d_l$.
  This easily implies that $\a$ is an increasing bijection of 
  $[\nu_\fm,\nu[$ onto $[1,\a(\nu)[$.
\end{proof}
\begin{proof}[Proof of Proposition~\ref{P201}]
  Keep the notation of Lemma~\ref{compuskew} and its proof.
  If $\nu=\nu_\phi$, then $\nu\wedge\nu_\phi=\nu_\phi$ and 
  the result follows from~(ii).
  Otherwise $l\=\con(\nu_\phi,\nu)<\infty$ and
  $\nu\wedge\nu_\phi=\val[(U_j)_0^l;(\btilde_j)_0^{l-1},\btilde]$,
  where $\btilde=\min(\btilde_l,\btilde_l^\phi)$.
  Then 
  \begin{equation*}
    \alpha(\nu\wedge\nu_\phi)
    =d_l^{-1}\btilde\btilde_0
    =\nu(\phi)/m(\phi),
  \end{equation*}
  which completes the proof.
\end{proof}                                
\begin{remark}\label{rem-vol}
  Skewness is closely related to 
  \emph{volume}\index{valuation!volume of} 
  as defined in~\cite{ELS}. Indeed, if we define
  $\vol(\nu)= \limsup_{c\to\infty}2c^{-2}\dim_\C R/\{\nu\ge c\}$,
  then we have
  \begin{equation}\label{E-vol}
    \vol(\nu)=\a(\nu)^{-1}.
  \end{equation}
  We sketch the proof, referring to~\cite[Example~3.15]{ELS}
  for details. Write
  $\nu=\val[(U_j)_0^k;(\btilde_j)_0^k]$ and suppose for simplicity that
  $k$ is finite. For any $c>0$, a basis for the vector space 
  $R/\{\nu\ge c\}$ is given by the monomials $U_0^{r_0}\dots U_k^{r_k}$ 
  where $0\le r_j<n_j$ for $1\le j\le k-1$, and $\sum_0^kr_j\btilde_j<c$. 
  The number of such monomials is given (up to a bounded function) by
  $\prod_1^{k-1}n_j\cdot\mathrm{Area}
  \{r_0+\btilde_kr_k\le c\ ;\ r_0,r_k\ge0\}\simeq d_kc^2/2\btilde_k$.
  This gives~\eqref{E-vol}.
\end{remark}
%
%
\subsection{Tree metrics}\label{normal}
Skewness defines a parameterization of the valuative tree. As we
have seen, there is a close connection between parameterized trees
and metric trees. Here we use this relationship to define natural
tree metrics both on the valuative tree $\cV$ and on the subtree
$\cVqm$ consisting of quasimonomial valuations.

We start with the quasimonomial case. 
For $\nu,\mu\in\cVqm$ set
\begin{equation}\label{e112}
  d_\mathrm{qm}(\mu,\nu)
  =\left(\alpha(\mu)-\alpha(\mu\wedge\nu)\right)
  +\left(\alpha(\nu)-\alpha(\mu\wedge\nu)\right).
\end{equation}
This makes sense since any quasimonomial valuation has finite
skewness. As an immediate consequence of Proposition~\ref{P501}
we have
\begin{theorem}
  The metric $d_\mathrm{qm}$\index{$d_\mathrm{qm}$ (metric on $\cVqm$)}
  defines a metric tree structure
  on the set $\cVqm$ of quasimonomial valuations.
\end{theorem}
As a general valuation can have infinite skewness,~\eqref{e112} does
not define a metric on $\cV$, at least not in a standard sense. 
This problem can be resolved by postcomposing skewness by a
positive, monotone, bounded function on $[1,\infty]$. In view of
a later application (see Theorem~\ref{treeballs}) we use the 
function $\a\mapsto\a^{-1}$. Hence we define
\begin{equation}\label{e105}
  d(\mu,\nu)
  =\left(\frac1{\alpha(\mu\wedge\nu)}-\frac1{\alpha(\mu)}\right)
  +\left(\frac1{\alpha(\mu\wedge\nu)}-\frac1{\alpha(\nu)}\right).
\end{equation}
for any valuations $\nu$, $\mu$ in $\cV$. 
\begin{theorem}\label{met-tree-val}
  The metric $d$\index{$d$ (metric on $\cV$)}
  gives valuation space $\cV$ the structure
  of a metric tree. Further, $(\cV,d)$ is complete.
\end{theorem}
\begin{proof}
  A simple adaptation of the proof of Proposition~\ref{P501} 
  shows that $(\cV,d)$ is a metric tree.
  Completeness of $(\cV,d)$ follows from 
  completeness of $\cV$ as a nonmetric tree.
\end{proof}
\begin{remark}
  Proposition~\ref{P502} shows that 
  the metrics $d$ and $d_\mathrm{qm}$ are invariant
  metrics in the sense that if 
  $f:(\C^2,0)\to(\C^2,0)$ is an invertible formal map, then 
  the induced maps $f_*:\cV\to\cV$ and $f_*:\cVqm\to\cVqm$ 
  are isometries for $(\cV,d)$ and $(\cVqm,d_\mathrm{qm})$,
  respectively.
\end{remark}
%
%
%
%
\section{Multiplicities}\label{sec-mult}
Consider the maximal, increasing function $m:\cV\to\Nbar$ such that
$m(\nu_C)=m(C)$ for every irreducible formal curve $C$. Concretely,
this is given as follows: the
\emph{multiplicity}
\index{multiplicity!of a valuation} 
\index{multiplicity!in $\cV$}
\index{$m(\nu)$ (multiplicity)}
$m(\nu)$ of a
quasimonomial valuation $\nu\in\cVqm$ is defined by
\begin{equation*}
  m(\nu)=\min\{m(C)\ ;\ C\ \text{irreducible},\ \nu_C\ge\nu\}.
\end{equation*}
This multiplicity can be extended to arbitrary valuations in $\cV$
by observing that $\nu\mapsto m(\nu)$ is increasing.
We then have to allow for the possibility that $m=\infty$. 
\begin{proposition}\label{P402}
  If $\mu\le\nu$ then $m(\mu)$ divides $m(\nu)$. Further:
  \begin{itemize}
  \item[(i)]
    $m(\nu)=\infty$ iff $\nu$ is infinitely singular;
  \item[(ii)]
    $m(\nu)=1$ iff $\nu$ is monomial in some local coordinates $(x,y)$.
  \end{itemize}    
\end{proposition}
\begin{proposition}\label{P301}
  The multiplicity function $m:\cV\to\Nbar$ is
  lower semicontinuous, but not continuous, in the weak topology on $\cV$.
\end{proposition}
We postpone the proofs to the end of this section.

We also define $m(\vv)$ for a 
tangent vector
\index{$m(\vv)$ (multiplicity)}
\index{multiplicity!of a tangent vector} 
$\vv\in T\nu$ as follows: if
$\vv$ is represented by $\nu_\fm$, then $m(\vv):=m(\nu)$. 
Otherwise, $m(\vv)$ is the minimum of $m(C)$ over 
all $\nu_C$ representing $\vv$.
As there are uncountably many tangent vectors at a divisorial valuation,
it is useful to know that their multiplicities are well behaved:
\begin{proposition}\label{P111}
  Let $\nu$ be a divisorial valuation. Then there exists a
  positive integer $b(\nu)$ divisible by $m(\nu)$ such that 
  exactly one of the following holds:
  \begin{itemize}
  \item[(i)]
    all tangent vectors at $\nu$ share the same multiplicity $m(\nu)$;
    in this case we set $b(\nu)=m(\nu)$;
  \item[(ii)]
    among the tangent vectors at $\nu$, there are exactly two with
    multiplicity $m(\nu)$ and one of them is determined by $\nu_\fm$; 
    all other tangent vectors at $\nu$ share a common 
    multiplicity $b(\nu)>m(\nu)$.
  \end{itemize}
\end{proposition}
\begin{definition}
  We call $b(\nu)$ the 
  \emph{generic multiplicity}
  \index{multiplicity!generic} 
  \index{$b(\nu)$ (generic multiplicity)}
  of $\nu$.
\end{definition}
\begin{remark}\label{rem-gen-mult}
  As we shall show in Chapter~\ref{A3}, 
  the generic multiplicity $b=b(\nu)$ has the following 
  geometric interpretation.
  There exists a composition of blowups $\pi:X\to(\C^2,0)$
  and an irreducible component $E$ of the exceptional divisor
  $\pi^{-1}(0)$ such that $\pi_*\div_E=b\,\nu$.
  Moreover, if $C$ is an irreducible 
  curve whose strict transform intersects
  $E$ transversely at a point of $E$ which is smooth on $\pi^{-1}(0)$,
  then $C$ has multiplicity $m(C)=b$.
  The situation $b(\nu)=m(\nu)$ happens exactly when
  $\pi=\pi'\circ\tpi$ and $E$ is the exceptional divisor
  of the blowup $\tpi$ of a \emph{smooth} point on
  $(\pi')^{-1}(0)$.
\end{remark}
Let us now turn to the proofs of Propositions~\ref{P402},~\ref{P301}
and~\ref{P111}. They all rely on a translation to statements about
SKP's.
\begin{lemma}\label{L413}
  Consider a valuation $\nu\in\cV$ and pick local coordinates $(x,y)$
  such that $1=\nu(x)\le\nu(y)$. 
  Write $\nu=\val[(U_j)_0^k;(\btilde_j)_0^k]$, $1\le k\le\infty$.  
  Then the multiplicity $m(\nu)$
  is the maximum of the degree in $y$ of the polynomials $U_j$.
\end{lemma}
\begin{proof}
  By Lemma~\ref{lem:mult-irr} we know that $m(U_j)=d_j=\deg_yU_j$,
  which is an increasing function in $j$. The result then follows
  easily from the definition of $m(\nu)$ and from
  Proposition~\ref{pro-comp}. The details are left to the reader.
\end{proof}
\begin{proof}[Proof of Proposition~\ref{P402}]
  First pick $\mu\le\nu$.  Write $\nu=\val[(U_j)_0^k;(\btilde_j)_0^k]$,
  $1\le k\le\infty$. By Proposition~\ref{pro-comp}, we have
  $\nu=\val[(U_j)_0^{k'};(\btilde_j)_0^{k'-1}, \btilde'_{k'}]$, with
  $k'\le k$ and $\btilde'_{k'}\le\btilde_{k'}$. By Lemma~\ref{L413} and
  Lemma~\ref{lem:deg-irr}, we get $m(\mu) = \prod_0^{k'-1} n_j$ and
  $m(\nu) = \prod_0^{k-1} n_j$.  Whence $m(\mu)$ divides $m(\nu)$.

  Definition~\ref{D101} and Lemma~\ref{L413} immediately imply that
  $m(\nu)=\infty$ iff $\nu$ is infinitely singular. Thus~(i) holds.

  As for~(ii) it is clear that if $\nu$ is monomial in coordinates
  $(x,y)$, then $\nu\le\nu_x$ or $\nu\le\nu_y$, which implies
  $m(\nu)=1$. Conversely, suppose that $m(\nu)=1$.
  Pick an irreducible formal curve $C$ 
  such that $\nu_C\ge\nu$ and $m(C)=m(\nu)=1$. 
  We may then pick local coordinates $(x,y)$ with $C=\{y=0\}$.
  Thus $\nu_C=\nu_y$. 
  For each $t\in[1,\infty]$, let $\nu_t$ be the
  monomial valuation in coordinates $(x,y)$ with
  $\nu_t(x)=1$, $\nu_t(y)=t$ in the sense of~\eqref{e440}.
  It is then clear that $\nu_t\le\nu_y$ and that
  $\a(\nu_t)$ has skewness $t$. As $\cV$ is a tree,  $\nu_t=\nu_{y,t}$.
  This shows that all valuations in the segment
  $[\nu_\fm,\nu_y]$ are monomial. 
  In particular $\nu$ is monomial.
\end{proof}
\begin{proof}[Proof of Proposition~\ref{P301}]
  We first show that $m$ is lower semicontinuous on each segment of
  $\cV$. Pick an increasing sequence of valuations $\nu_n\to\nu$.
  We want to prove $m(\nu)= \lim m(\nu_n)$. To do so write
  $\nu=\val[(U_j)_0^k;(\btilde_j)_0^k]$. Then
  Proposition~\ref{pro-comp} imply that for $n$ large enough, $\nu_n =
  \val[(U_j)_0^k;(\btilde_j)_0^{k-1},\btilde_k^n]$ with $\btilde_k^n$
  increasing to $\btilde_k$. We conclude by Lemma~\ref{L413} above.

  Lower semicontinuity on all of $\cV$ is then 
  proved in much the same way as Proposition~\ref{P305}. 
  We have to show that the superlevel set
  $\{\nu\ ;\ m(\nu)>t\}$ is  weakly open for any real $t$.  Pick any
  $\nu_0\in\cV$ with $m(\nu_0)>t$.  If $\nu_0=\nu_\fm$, then $t<1$ and
  $m(\nu)>t$ for all $\nu\in\cV$. Hence suppose $\nu_0\ne\nu_\fm$.  Since
  $m$ is lower semicontinuous and increasing on the segment
  $[\nu_\fm,\nu_0]$, we may find $\mu<\nu$ such that
  $m(\mu)=m(\nu_0)>t$. Let $\vv$ be the tangent vector at $\mu$
  represented by $\nu_0$. Then $m>t$ on the open neighborhood $U(\vv)$
  of $\nu_0$.

  Thus $m$ is lower semicontinuous. 
  The example $\nu_n=\nu_{(y+nx)^2+x^3}$ in coordinates $(x,y)$
  shows that it is not continuous: $m(\nu_n)=2$, 
  $\nu_n\to\nu_\fm$ but $m(\nu_\fm)=1$.
\end{proof}

\begin{proof}[Proof of Proposition~\ref{P111}]
  We use the description of the tgt space at a divisorial valuation in
 terms of SKP's that is given in the proof of
 Proposition~\ref{pro-comp}.  Fix a divisorial valuation
 $\nu=\val[(U_j)_0^k;(\btilde_j)_0^k]$.  Let $n_k = \min \{ l\ ;\
 l\btilde_k\in \sum_0^{k-1} \Z \btilde_j \}$, and write $n_k
 \btilde_k= \sum m_{kj} \btilde_j$ with $0\le m_{kj}\le n_j$ as
 in~\eqref{key3}.  For any $\theta\in\C^*$ we define
 $U_{k+1}(\theta)=U_k^{n_k}-\theta\cdot\prod U_j^{m_{kj}}$.
  
  Then any tangent vector at $\nu$ which is not determined by $\nu_\fm$
  is represented by one and only one of the following curve valuations:
  $\nu_\phi$ with $\phi=U_k$, or $\phi=U_{k+1}(\theta)$ for some
  $\theta\in\C^*$. This representation gives a bijection between $T\nu$
  and $\P^1$.

  Let us now prove the proposition. Pick coordinates $(x,y)$
  such that $\nu(y)\ge\nu(x)$.
  Write $\nu=\val[(U_j)_0^k;(\btilde_j)_0^k]$ and pick $\vv\in T\nu$. 
  When $\vv$ is represented by $\nu_\fm$, its
  multiplicity is $m(\nu)$ by definition. 
  When $\vv$ is represented by
  the curve valuation associated to $U_k$, the multiplicity of $\vv$ is
  bounded from above by $m(U_k)=m(\nu)$ and hence equals $m(\nu)$.
  
  Otherwise, $\vv$ is represented by a curve valuation associated to
  $U_{k+1}=U_k^{n_k}-\theta\prod_0^{k-1}U_j^{m_{kj}}$, for some
  $\theta\in\C^*$.
  The multiplicity of $\vv$ is hence at most 
  $m(U_{k+1})=n_km(U_k)=n_km(\nu)$ (see Lemma~\ref{lem:mult-irr}). 
  On the other hand, pick any curve valuation
  $\nu_\phi$ representing $\vv$. As $\nu_\phi\ge\nu$ the SKP defining
  $\nu_\phi$ starts with $(U_j)_0^{k}$. 
  When $\nu_\phi(U_k)>\btilde_k$, $\nu_\phi$ determines the same 
  tangent vector as $U_k$. Therefore $\nu_\phi(U_k)=\btilde_k$ 
  and the next polynomial in the SKP of $\nu_\phi$ is of the form 
  $U_k^{n_k}-\theta'\prod_0^{k-1}U_j^{m_{kj}}$. 
  As $\nu_\phi$ and $U_{k+1}$ determine the same tangent vector, 
  $\theta'=\theta$. The multiplicity of $\phi$ is
  the supremum of the multiplicity of the polynomials appearing in the
  SKP of $\nu_\phi$, hence $m(\phi)\ge m(U_{k+1})=n_k m(\nu)$.  We
  conclude that $m(\vv)=n_k m(\nu)$.
  
  When $n_k=1$, we have proved that all tangent vectors have the same
  multiplicity. When $n_k\ge2$ all tangent vectors but two have the same
  multiplicity, and this multiplicity is a multiple of $m(\nu)$. 
  The remaining two tangent vectors have multiplicity
  $m(\nu)$. This completes the proof.
\end{proof}
Note that the proof gives
\begin{lemma}\label{L423}
  Consider a divisorial valuation $\nu\in\cV$ and pick local
  coordinates $(x,y)$ such that $1=\nu(x)\le\nu(y)$.  
  Write $\nu=\val[(U_j)_0^k;(\btilde_j)_0^k]$, $1\le k\le\infty$.  
  Then the generic multiplicity $b(\nu)$ is equal to 
  the product $n_kd_k$, where $d_k=\deg_y(U_k)$ and
  $n_k$ is defined by condition~(P2) in Definition~\ref{def-key}.
\end{lemma}
%
%
%
%
\section{Approximating sequences}\label{sec-approx}
\begin{figure}[h]
  \begin{center}
    \includegraphics{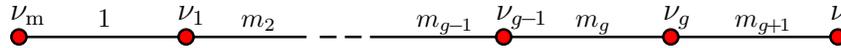}
  \end{center}
  \caption{The approximating sequence of a valuation of 
    finite multiplicity.}\label{F8}
\end{figure}

As the multiplicity function $m$ is increasing, it is piecewise
constant on any segment and the set of jump points is discrete. By
Propositions~\ref{P402},~\ref{P301} and~\ref{P111} this observation
immediately gives
\begin{proposition}\label{prop-approx}
  For any valuation $\nu\in\cV$ of finite multiplicity, there 
  exists a finite sequence of divisorial valuations $\nu_i$ and
  a strictly increasing sequence of integers $m_i$ such that:
  \begin{equation}\label{e705}
    \nu_\fm=\nu_0<\nu_1<\dots<\nu_g<\nu_{g+1}=\nu
  \end{equation}
  and such that $m(\mu)=m_i$ for $\mu\in]\nu_{i-1},\nu_i]$,
  $1\le i\le g+1$. Moreover $m_1=1$, $m_i$ divides $m_{i+1}$
  and $\nu_i$ has generic multiplicity $b(\nu_i)=m_{i+1}$.
  See Figure~\ref{F8}.
\end{proposition}
\begin{definition}
  We call the sequence $(\nu_i)_{i=0}^g$ the 
  \emph{approximating sequence}\index{approximating sequence} 
  associated to $\nu$. The concept of approximating sequences extends 
  naturally to infinitely singular valuations: 
  in this case they are infinite,
  of the form $(\nu_i)_0^\infty$.
\end{definition}
In Chapter~\ref{A3} we shall give a geometric interpretation of the
approximating sequence: see Remark~\ref{R306}.

As we show in Appendix~\ref{sec-termino}, the approximating
sequence of a curve valuation can be used to extract the classical
invariants of the associated irreducible curve.
%
%
%
%
\section{Thinness}\label{S10}
As we have indicated above, skewness is the first out of (at least)
two natural parameterizations of the valuative tree. Having defined
multiplicities we are now in position to introduce the second.
Namely, we define the \emph{thinness}\index{thinness}\footnote{We
normalize $A$ such that $A(\nu_\fm)=2$. This is natural in view
of Remarks~\ref{Rpict} and~\ref{Rjaco}.} of a valuation $\nu\in\cV$ by
\begin{equation}\label{e401}
  A(\nu)=2+\int_{\nu_\fm}^\nu m(\mu)\,d\alpha(\mu).
\end{equation}
Here the integral in~\eqref{e401} is defined by
$\int_{\nu_\fm}^\nu m(\nu)\,d\a(\mu)=\int_1^t m(\nu_{C,s})\,ds$,
if $\nu=\nu_{C,t}$ is quasimonomial or a curve valuation,
and as a suitable increasing limit if 
$\nu$ is infinitely singular.

Concretely, if $\nu$ has finite multiplicity and
$(\nu_i)_{i=0}^g$ is its approximating sequence, then 
\begin{equation}\label{e412}
  A(\nu)=2+\sum_{i=1}^{g+1}m_i(\alpha_i-\alpha_{i-1}),
\end{equation}
where $\alpha_i=\alpha(\nu_i)$. 
This formula extends naturally to infinitely singular valuations,
in which case $g=\infty$.
\begin{proposition}\label{P306}
  Thinness defines a parameterization $A:\cV\to[2,\infty]$
  \index{parameterization!of the valuative tree (by thinness)}
  of the valuative tree. Moreover:
  \begin{itemize}
  \item[(i)]
    if $\nu$ is divisorial, then $A(\nu)$ is rational;
  \item[(ii)]
    if $\nu$ is irrational, then $A(\nu)$ is irrational;
  \item[(iii)]
    if $\nu$ is a curve valuation, then $A(\nu)=\infty$;
  \item[(iv)]
    if $\nu$ is infinitely singular, then $A(\nu)\in(2,\infty]$.
  \end{itemize}  
\end{proposition}
It follows that thinness also defines a parameterization
$A:\cVqm\to[2,\infty)$ of the tree $\cVqm$ of
quasimonomial valuations.

We refer to Appendix~\ref{sec-inf-sing} for a construction of an
infinitely singular valuation with  prescribed thinness 
$t\in(2,\infty]$.
\begin{proof}
  All of this follows from~\eqref{e412} and Theorem~\ref{T401}.
\end{proof}

\begin{proposition}\label{P302}
  Thinness $A:\cV\to[2,\infty]$ is
  lower semicontinuous, but not continuous, 
  in the weak tree topology on $\cV$.
\end{proposition}
\begin{proof}
  Lower semicontinuity follows from Proposition~\ref{P305} and the 
  fact that thinness defines a parameterization of $\cV$.
  The same example $\nu_n=\nu_{y-nx,2}$ that was used for skewness
  (see Proposition~\ref{P303}) serves to show that thinness is
  not continuous in the weak topology. 
\end{proof}
Let us further analyze the relationship between skewness, thinness 
and multiplicity. When viewed as functions on $\cV$, any two of them
determines the third. Indeed,~\eqref{e401} gives thinness in terms
of skewness and multiplicity. 
We can then recover skewness from thinness and multiplicity:
\begin{equation}\label{e402}
  \a(\nu)=1+\int_{\nu_\fm}^\nu\frac1{m(\mu)}\,dA(\mu).
\end{equation}
Finally, multiplicity can easily be recovered from thinness and
skewness by differentiating either~\eqref{e401} or~\eqref{e402}.

\begin{figure}[h]
  \begin{center}
    \includegraphics{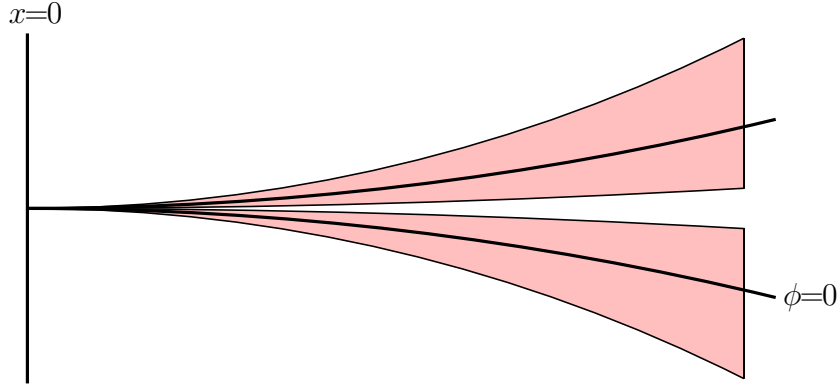}
  \end{center}
  \caption{A characteristic region}\label{F7}
\end{figure}

Alternatively, we can try to relate skewness, thinness and
multiplicity for a fixed valuation. There is no general formula,
but the following elementary bounds are useful
for applications.
\begin{proposition}\label{bounds}
  For any valuation $\nu\in\cV$ we have
  \begin{itemize}
  \item[(i)]
    $A(\nu)\ge1+\alpha(\nu)$ with equality iff $m(\nu)=1$, 
    \ie if $\nu$ is monomial in some local coordinates;
  \item[(ii)]
    if $m(\nu)>1$, then $A(\nu) < m(\nu)\a(\nu)$.
  \end{itemize}
\end{proposition}
\begin{proof}
  The first statement is immediate. As for the second
  we use the approximating sequence $(\nu_i)_0^g$. 
  Then $m(\nu)=m_{g+1}$, $\a(\nu)=\a_{g+1}$ and
  \begin{align*}
    A(\nu)-m(\nu)\a(\nu)
    &=1-(m_2-1)\a_1-(m_3-m_2)\a_2+\dots+(m_{g+1}-m_g)\a_g\\
    &\le 1 -\a_1 <0.
  \end{align*}
This concludes the proof.
\end{proof}
\begin{remark}\label{Rpict}
  The name ``thinness'' was chosen with the following picture
  in mind. Suppose $\nu$ is quasimonomial and write
  $\nu=\nu_{\phi,t}$ with $m(\phi)=m(\nu)$. Assume we have picked
  local coordinates $(x,y)$ with 
  $\nu_x\wedge\nu_\phi=\nu_\fm$. Then, for $r>0$ small, the region
  \begin{equation*}
    \Omega(r)
    =\left\{(x,y)\in\C^2\ ;\ |x|<r,\ |\phi(x,y)|<|x|^{tm(\phi)}\right\}
    \subset\C^2
  \end{equation*}
  is a small neighborhood of the curve $\phi=0$ (with the origin removed). 
  See Figure~\ref{F7}.
  A large value of $t$ (and hence $A(\nu)$) corresponds to a 
  ``thin'' neighborhood.
  In fact the  volume of $\Omega(r)$ is roughly 
   $r^{2A(\nu)}$. Regions of this type, called
   \emph{characteristic regions}\index{characteristic region},
   are important in~\cite{criterion}, where they are used
   to define the values of $\nu$ on plurisubharmonic functions.
   See also Section~\ref{S32}.
\end{remark}
\begin{remark}\label{Rjaco}
  In case $\nu$ is a divisorial valuation, thinness has the
  following algebraic characterization. 
  We shall freely use results from Chapter~\ref{A3}. 
  As the value group of $\nu$ is isomorphic to $\Z$,
  the maximal ideal $\fm_\nu$ of the valuation ring $R_\nu$ is
  principal, and any other ideal is a power of $\fm_\nu$.  
  We may define an ideal $J_{R_\nu/R}$ inside the valuation ring 
  $R_\nu$ called its Jacobian ideal.  
  By definition, this is the $0$-th Fitting ideal of
  the sheaf of differentials of $R_\nu$ over $R$.
  
  Geometrically the situation is as follows.  Let $\pi$ be a
  composition of blowups above the origin, and $E\subset\pi^{-1}(0)$
  an irreducible component of the exceptional divisor $\pi^{-1}(0)$
  such that $\pi_*\div_E=b\,\nu$ for some $b\in\N^*$.  Then an
  element $\phi\in R_\nu\subset K$ belongs to $\fm_\nu$ iff
  $\pi^*\phi$ defines a regular function at a generic point on $E$ and
  vanishes to order at least one. Similarly, $\phi\in J_{R_\nu/R}$
  iff $\pi^*\phi$ defines a regular function at a generic point on $E$
  and vanishes to at least the same order as the Jacobian determinant
  $J\pi$.  Write $a-1=\div_E(J\pi)$. Then, we have
  $\nu(\fm_\nu)=b^{-1}$, and $\nu(J_{R_\nu/R})=(a-1)/b$.  Thanks to
  Theorem~\ref{thm-universal}, the Farey parameter coincides with the
  thinness, \ie we have $A(\nu)=a/b$, whence
  $A(\nu)=\nu(J_{R_\nu/R}\,\fm_\nu)$.
\end{remark}
\begin{remark}
  Associated to the parameterization by thinness are
  tree metrics 
  $D_\mathrm{qm}$\index{$D_\mathrm{qm}$ (metric on $\cVqm$)} 
  and $D$\index{$D$ (metric on $\cV$)}
  on $\cVqm$ and $\cV$, 
  respectively. They are defined
  similarly to the tree metrics 
  $d_\mathrm{qm}$ and $d$ induced by skewness: 
  see~\eqref{e112} and~\eqref{e105}.
  More precisely, set 
  \begin{equation}\label{e512}
    D_\mathrm{qm}(\mu,\nu)
    =\left(A(\mu)-A(\mu\wedge\nu)\right)
    +\left(A(\nu)-A(\mu\wedge\nu)\right),
  \end{equation}
  for $\nu,\mu\in\cVqm$ and  
  \begin{equation}\label{e513}
    D(\mu,\nu)
    =(A^{-1}(\mu\wedge\nu)-A^{-1}(\mu))
    +(A^{-1}(\mu\wedge\nu)-A^{-1}(\nu))
  \end{equation}
  for $\nu,\mu\in\cV$.
  We shall compare these metrics with $d_\mathrm{qm}$ and $d$ in
  Chapter~\ref{part4}.
\end{remark}
%
%
%
%
\section{Value semigroups and approximating sequences}\label{S12}
The approximating sequence encapsulates most of the
structure of a valuation. 
We use it here to compute the value semigroup of a general 
valuation, and the generic multiplicity of 
a divisorial valuation.
\begin{proposition}\label{P556}
  Let $\nu\in\cV$ be a valuation and let $(\nu_i)_0^g$ be its
  approximating sequence ($g$ possibly infinite). Then the value
  semigroup $\nu(R)$ is given by
  \begin{equation}\label{e413}
    \nu(R)=\sum_{i=0}^{g+1}m_i\a_i\N,
  \end{equation}
  where $\nu_{g+1}=\nu$ if $g<\infty$, 
  $\a_i=\a(\nu_i)$ and $m_i=m(\nu_i)$.
\end{proposition}
\begin{proposition}\label{P403}
  Let $\nu\in\cV$ be a divisorial valuation with approximating sequence
  $(\nu_i)_0^g$. Write $\a_i=\a(\nu_i)$ for $0\le i\le g$.
  Then the generic multiplicity of $\nu$ can be characterized by the 
  equations:
  \begin{equation}\label{e707}
    b(\nu)
    =\inf\{b\in m\N\ ;\ b\a\in\sum_{i=0}^gm_i\a_i\N\}
    =\inf\{b\in m\N\ ;\ b\a\in\sum_{i=0}^gm_i\a_i\Z\},
  \end{equation}
  where $m=m(\nu)$, $\a=\a(\nu)$ and
  where we set $m_0=1$ by convention.
\end{proposition}
Notice that this implies that for $1\le j\le g$, 
$m_{j+1}$ is the smallest integer divisible by $m_j$ such that
$m_{j+1}\a_j\in\sum_{i=0}^{j-1}m_i\a_i\N$.
\begin{proposition}\label{P421}
  Consider a valuation $\nu\in\cV$. Then the following hold:
  \begin{itemize} 
  \item[(i)] if $\nu$ is a curve valuation with multiplicity $m=m(\nu)$,
    then $\nu(R)\subset m^{-1}\Nbar$;
    moreover, if $a\in\N$, then $a\,\nu(R)\subset\Nbar$
    iff $m$ divides $a$;
  \item[(ii)] if $\nu$ is divisorial with generic multiplicity
    $b=b(\nu)$, then $\nu(R)\subset b^{-1}\N$;
    moreover, if $a\in\N$, then $a\,\nu(R)\subset\N$
    iff $b$ divides $a$;
  \item[(iii)] if $\nu$ is irrational with multiplicity $m=m(\nu)$ 
    and skewness $\a=\a(\nu)$, then
    $\nu(R)\subset\Q+m\a\N$ but $\nu(R)\not\subset\Q$;
    moreover, if $a\in\N$, then $a\,\nu(R)\subset\N+\a\N$
    iff $m$ divides $a$;
  \item[(iv)] if $\nu$ is infinitely singular, then $\nu(R)\subset\Q$
    but $\nu(R)$ is not finitely generated.
  \end{itemize}
\end{proposition}

The proofs of these three propositions are based on the 
representation of a valuation by an SKP. Before embarking 
on the individual proofs, we note some common facts
regarding the translation of approximating sequences
into SKP's.

We are working with a fixed valuation $\nu\in\cV$.
Fix local coordinates $(x,y)$ such that $\nu(y)\ge\nu(x)=1$.
Write $\nu=\val[(U_j)_0^k ; (\btilde_j)_0^k]$ 
with $\btilde_1\ge\btilde_0=1$.
Also write $d_j=\deg_yU_j$.
Define $n_j$ as in~\eqref{key3}, \ie
$n_j\=\min\{n\in\N\ ;\ n\btilde_j\in\sum_{j'=0}^{j-1}\btilde_{j'}\Z\}$.
Then $d_k=\prod_0^{k-1}n_j$.
Define a sequence $(k_i)_0^{g+1}$ by $k_0=0$ and by the 
conditions $n_{k_i}\ge 2$ and $n_j=1$ if $k_i<j<k_{i+1}$.

By the proof of Proposition~\ref{prop-approx}
we have $\nu_i=\val[(U_j)_0^{k_i};(\btilde_j)_0^{k_i}]$
and $m_i:=m(\nu_i)=d_{k_i}$. 
Moreover, by Lemma~\ref{compuskew}, 
$\a_i:=\a(\nu_i)=\btilde_{k_i}/d_{k_i}$.
Also recall that $m_{i+1}$ is the generic multiplicity
of $\nu_i$ for $1\le i<g+1$. 
In particular, if $\nu$ is divisorial, then its generic
multiplicity is given by $b(\nu)=d_kn_k=\prod_1^kn_j$.
\begin{proof}[Proof of Proposition~\ref{P556}]
  By Theorem~\ref{divis}, $\nu(R)=\sum_0^k\btilde_j\N$.
  But if $k_i<j<k_{i+1}$, then $n_j=1$,
  \ie $\btilde_j\in\sum_{j'<j}\btilde_{j'}\N$.  
  This gives $\nu(R)=\sum_0^{g+1}\btilde_{k_i}\N$
  which amounts to $\nu(R)=\sum_0^{g+1}m_i\a_i\N$
  using the translation above.
\end{proof}
\begin{proof}[Proof of Proposition~\ref{P403}]
  We have $b(\nu)=d_kn_k$. 
  Here $n_k=\min\{a\in\N\ ;\ a\btilde_k\in\sum_{j=0}^{k-1}\btilde_j\Z\}$. 
  As in the previous proof, we have
  $\sum_0^{k-1}\btilde_j\Z=\sum_0^gm_i\a_i\Z$.
  Moreover, $\btilde_k=m\a$, so multiplying with $d_k=m$ we get
  \begin{equation*}
    b(\nu)=\min\{b\in m\N\ ;\ b\a\in\sum_0^gm_i\a_i\Z\}.
  \end{equation*}
  Finally we can replace $\Z$ by $\N$ by Lemma~\ref{lem:arith}.
\end{proof}
\begin{proof}[Proof of Proposition~\ref{P421}]
  We first consider~(ii), so assume $\nu$ is divisorial.  Then
  $k<\infty$.  Define $a=\min\{a\in\N\ ;\ a\nu(R)\subset\N\}$.
  Then~(ii) will follow if we can show that $a=b(\nu)$.  By
  Theorem~\ref{divis} we have $\nu(R)=\sum_0^k\btilde_j\N$.  Write
  $\btilde_j=r_j/s_j$, where $r_j$ and $s_j$ are integers without
  common factor and set $S_j=\lcm\{s_0,s_1,\dots,s_j\}$ for $0\le j\le
  k$.  Then $S_0=1$ and $a=S_k^{-1}$.  Moreover, it is an elementary
  arithmetic fact that $\sum_0^j\btilde_{j'}\Z=S_j^{-1}\Z$, although
  $\sum_0^j\btilde_{j'}\N\subsetneq S_j^{-1}\N$ in general.  Since
  $b(\nu)=\prod_1^k n_j$ we are done if we can show that
  $n_j=S_j/S_{j-1}$ for $1\le j\le k$.  But \begin{multline*}
  n_j=\min\{n\in\N\ ;\ n\btilde_j\in\sum_0^{j-1}\btilde_j\Z\}
  =\min\{n\in\N\ ;\ n\btilde_j\in S_{j-1}^{-1}\Z\}\\ =\min\{n\in\N\ ;\
  nr_jS_{j-1}\in s_j\Z\} =\min\{n\in\N\ ;\ nS_{j-1}\in s_j\Z\}\\
  =\lcm\{s_j,S_{j-1}\}/S_{j-1} =S_j/S_{j-1}.  \end{multline*} This
  completes the proof of~(ii).

  Next suppose that $\nu$ is a curve valuation.
  Let $(\nu_i)_0^g$ be its approximating sequence. 
  Here $0\le g<\infty$. 
  It follows from Proposition~\ref{P556}
  that $\nu(R)=\nu_g(R)\cup\{\infty\}$. 
  Moreover, $\nu_g$ is a divisorial valuation with
  generic multiplicity $b(\nu_g)=m(\nu)$. 
  Hence~(i) follows from(ii).

  As for~(iii), suppose $\nu$ is irrational,
  and let $(\nu_i)_0^g$ be its (finite) approximating sequence.
  Then~\eqref{e413} gives $\nu(R)=\nu_g(R)+m\a\N$.
  As $\nu_g$ is divisorial, this immediately yields
  $\nu(R)\subset\Q+m\a\N$ but $\nu(R)\not\subset\Q$.
  Further, if $a\in\N$, then 
  $a\,\nu(R)\subset\N+\a\N$ iff $a\,\nu_g(R)\subset\N$, which 
  by~(ii) is the case iff $b(\nu_g)=m(\nu)$ divides $a$.

  Finally consider an infinitely singular valuation
  $\nu$, with approximating sequence $(\nu_i)_0^\infty$.
  Then $\a(\nu_i)\in\Q$ so~\eqref{e413} yields $\nu(R)\subset\Q$.
  Suppose that $\nu(R)$ is a finitely generated semigroup.
  Then there exists $b\in\N$ such that $b\nu(R)\subset\N$.
  But by~\eqref{e413} this implies
  $b\nu_i(R)\subset\N$, which by~(i) gives that $b_i$ divides $b$.
  But this is impossible as $b_i\to\infty$ as $i\to\infty$.

  The proof is complete.  
\end{proof}
%
%
%
%
\section{Balls of curves}\label{val-curve}
Above we have shown that the valuative tree $\cV$ carries a
very rich structure, induced by the partial ordering,
skewness, multiplicity and thinness. 
The key to this structure was the detailed analysis of 
valuations as functions on the ring $R$.

However, it is an intriguing fact that there are many
paths leading to the valuative tree. 
Here we show how to identify quasimonomial valuations 
with balls of irreducible curves in a particular (ultra-)metric. 
There are then natural interpretations of 
the partial ordering, skewness and multiplicity on $\cVqm$
as statements about balls of curves.

The description of a quasimonomial
valuations as a ball of curves adds to the observation by 
Spivakovsky~\cite[p.109]{spiv}
that the classification of valuations (in dimension 2)
is essentially equivalent to the classification of plane curve 
singularities.
%
%
\subsection{Valuations through intersections.}
The starting point is the fact that a curve valuation acts by
intersection (see Section~\ref{S33}). 
Given formal curves $C,D$, let $C\cdot D$ denote the
intersection product between $C$ and $D$
(see Section~\ref{S24}).
If $C$ is irreducible, then for any $\psi\in\fm$:
\begin{equation}\label{valu}
  \nu_C(\psi)=\frac{C\cdot\{\psi=0\}}{m(C)}.
\end{equation}
Spivakovsky~\cite[Theorem~7.2]{spiv} showed that a 
divisorial valuation $\nu$ acts (up to normalization) 
by intersection with a $\nu$-generic curve (or $\nu$-curvette). 
See Section~\ref{sec-curvette} for a precise statement.

Here we give a different way of realizing a valuation as the intersection
with curves. In fact this works for general quasimonomial valuations.
\begin{proposition}\label{P101}
  If $\nu\in\cV$ is quasimonomial and $\psi\in\fm$ then 
  \begin{align*}
    \nu(\psi)
    =\min\{\nu_C(\psi)\ ;\ C\ \mathrm{irr.},\ \nu_C\ge\nu\}
    =\min\left\{\frac{C\cdot\{\psi=0\}}{m(C)}\ 
      ;\ C\ \mathrm{irr.},\ \nu_C\ge\nu\right\}. 
  \end{align*}
\end{proposition}
\begin{proof}
  By~\eqref{valu} we only have to show the first equality.
  Pick $\psi\in\fm$. 
  The inequality $\nu(\psi)\le\nu_C(\psi)$ 
  for $\nu_C\ge\nu$ is trivial. For the other inequality
  write $\psi=\psi_1\cdots\psi_n$, with $\psi_i\in\fm$ irreducible. 
  First assume that $\nu$ is
  divisorial. Pick a tangent vector $\vv$ at $\nu$ which is not
  represented by any $\nu_{\psi_i}$ nor $\nu_\fm$, and choose $C$
  irreducible such that $\nu_C$ represents $\vv$.
  Then $\nu_C\ge\nu$ and $\nu_C\wedge\nu_{\psi_i}=\nu$ for all $i$.  
  so by Proposition~\ref{P201} we get
  $\nu(\psi_i)=\nu_C(\psi_i)$ so that
  $\nu(\psi)=\nu_C(\psi)$.

  Thus Proposition~\ref{P101} holds for $\nu$ divisorial. 
  But if $\nu\in\cV$
  is irrational, then we may pick $\nu_n$ divisorial with 
  $\nu_n>\nu$ and $\nu_n$ decreasing to $\nu$. 
  Fix $\psi\in\fm$. 
  For each $n$ there exists an irreducible curve $C_n$
  with $\nu_{C_n}>\nu_n>\nu$ such that
  $\nu_{C_n}(\psi)=\nu_n(\psi)$. 
  Since $\nu_n(\psi)\to\nu(\psi)$ we obtain the desired equality.
\end{proof}
%
%
\subsection{Balls of curves}\label{sec-balls}
\index{balls of curves}
Proposition~\ref{P101} shows that a quasimonomial valuation $\nu$ is
determined by the irreducible curves $C$ satisfying
$\nu_C\ge\nu$. We proceed to show that this gives an isometry
between the subtree of quasimonomial valuations, and the set of balls
of irreducible curves in a particular (ultra)-metric. 
Namely, let  $\cC$
\index{$\cC$ (space of irreducible curves)} 
be the set of local irreducible curves. 
Let us define a metric on $\cC$ by
\begin{equation}\label{e302}
  d_\cC(C_1,C_2)=\frac{m(C_1)m(C_2)}{C_1\cdot C_2},
\end{equation}
where $C_1\cdot C_2$ denotes the intersection multiplicity 
between the curves $C_1$ and $C_2$.
\begin{lemma}{\cite[Corollary 1.2.3]{GB}.}\label{ballmetric}
  Equation~\eqref{e302} defines an ultrametric on $\cC$.
  For $C_1,C_2\in\cC$, we have 
  \begin{equation}\label{e103}
    d_\cC(C_1,C_2)
    =\half d(\nu_{C_1},\nu_{C_2})
    =\frac1{\alpha(\nu_{C_1}\wedge\nu_{C_2})}.
  \end{equation}
  Further, the ultrametric space $(\cC,d_\cC)$ has diameter 1.
\end{lemma}
The proof is given at the end of this section.

In Section~\ref{def-ultra} we constructed a parameterized tree as
well as a metric tree associated to any ultrametric space (of diameter 1). 
Write $\cT_\cC$ for the parameterized tree associated to $(\cC,d_\cC)$. 
Let us recall its definition. 
A point in $\cT_\cC$ is a closed ball in $\cC$ of positive radius. 
The partial order on $\cT_\cC$ is given by reverse inclusion. 
The parameterization on $\cT_\cC$ sends a ball of radius $r$ to
the real number $r^{-1}$, and the distance between two balls
$B_1$, $B_2$ in the metric tree
is given by $(r_1-r)+(r_2-r)$, where 
$r_1$, $r_2$ and $r$ are the radii of $B_1$, $B_2$ and $B_1\cap B_2$,
respectively. 

There is a natural multiplicity function on $\cT_\cC$.  Its value on a
ball is the minimum multiplicity of any curve in the ball.

The completion of $\cT_\cC$ can be taken in the sense of nonmetric
trees or, equivalently, in the sense of metric spaces. An element in
the completion is represented by a decreasing sequence of balls. The
multiplicity function and the parameterization on $\cT_\cC$ both
extend naturally to the completion. 
An element in the completion, represented by a decreasing sequence 
of balls, has finite multiplicity
iff the sequence has nonempty intersection in $\cC$.

Let us define a map from $\cVqm$ to $\cT_\cC$ by sending 
a quasimonomial valuation $\nu$ to the ball 
$B_\nu:=\{C\ ;\ \nu_C>\nu\}$ (we will see shortly that this
is indeed a ball).
Similarly, $B\mapsto\nu_B:=\wedge_{C\in B}\nu_C$ gives a map
from $\cT_\cC$ to $\cVqm$.
\begin{theorem}\label{treeballs}
  The mappings $\nu\mapsto B_\nu$ and $B\mapsto\nu_B$ preserve
  multiplicity and give inverse isomorphisms between the parameterized 
  trees $(\cVqm,\a)$ and $(\cT_\cC,r^{-1})$.

  Further, these mappings extend uniquely to isomorphisms 
  between $\cV$ and the completion of $\cT_\cC$ 
  and induce isometries between the corresponding metric trees.
\end{theorem}
\begin{proof}
  We claim that if $B$ is centered at $C$ 
  and of radius $r\in(0,1]$, then $\nu_B=\nu_{C,r^{-1}}$. 
  Indeed, if $D\in\cC$, then $\nu_{C,r^{-1}}$ and
  $\nu_C\wedge\nu_D$ both belong to the segment 
  $[\nu_\fm,\nu_C]$. Hence~\eqref{e302} implies
  \begin{equation*}
    d_\cC(C,D)\le r
    \Leftrightarrow\alpha(\nu_C\wedge\nu_D)\ge r^{-1}
    \Leftrightarrow\nu_C\wedge\nu_D\ge\nu_{C,r^{-1}}
    \Leftrightarrow\nu_D\ge\nu_{C,r^{-1}}.
  \end{equation*}    
  Since every quasimonomial valuation is of the form $\nu_{C,t}$,
  this shows that the mappings $B\mapsto\nu_B$ 
  and $\nu\mapsto B_\nu$ are well defined and each others
  inverse. They are clearly increasing, hence define
  isomorphisms between the rooted,
  nonmetric trees $\cVqm$ and $\cT_\cC$.
  As $\a(\nu_{C,r^{-1}})=r^{-1}$, they are in fact isomorphism
  of parameterized trees. Multiplicity is preserved by its very 
  definition. 
  
  The remaining facts, namely that we get an isomorphism
  between $\cV$ and the completion of $\cT_\cC$ as well as
  isometries between the corresponding metric trees, are easy
  consequences and left to the reader.
\end{proof}
\begin{remark}
  Similarly one may interpret a tree tangent vector $\vv$
  at a divisorial valuations $\nu\in\cV$ in terms of irreducible
  curves. 
  If $\vv$ is represented by $\nu_\fm$, then $\cC\cap U(\vv)$
  is the complement of the ball $B_\nu$ in $\cC$. 
  If $\vv$ is not represented by $\nu_\fm$, then 
  $\cC\cap U(\vv)$ is the \emph{open} ball of radius $\a(\nu)^{-1}$
  centered at any $\phi\in\cC$ such that $\nu_\phi$ represents $\vv$.
\end{remark}
\begin{proof}[Proof of Lemma~\ref{ballmetric}]
  Write $C_i=\{\phi_i=0\}$ for $i=1,2$.
  Let us first show \eqref{e103}.
  The second equality follows from~\eqref{e105}
  as $\nu_{C_1}$ and $\nu_{C_2}$ have infinite skewness.
  For the first equality we use Proposition~\ref{P201}. 
  Assume $\nu_{C_1}\ne\nu_{C_2}$
  and set $\nu=\nu_{C_1}\wedge\nu_{C_2}$. 
  This is divisorial and $\a(\nu)=\nu(\phi_1)/m(C_1)$
  as $\nu_{C_1}>\nu$.
  But $\nu(\phi_1)=\nu_{C_2}(\phi_1)=C_1\cdot C_2/m(C_2)$.
  Thus~\eqref{e103} holds. 

  The fact that $d_{\cC}$ is an ultrametric for which $\cC$ has
  diameter $1$ follows easily from the formula
  $d_\cC(C_1,C_2)=\alpha^{-1}(\nu_{C_1}\wedge\nu_{C_2})$
  and the tree structure of $\cV$. 
  The details are left to the reader.
\end{proof}
%
%
%
%
\section{The relative tree structure}\label{sec-relative}
We have normalized the valuations in $\cV$ by $\nu(\fm)=1$. 
This condition is natural when we study properties
invariant by the group of (formal) automorphisms of
$(\C^2,0)$. However, in some situations, other normalization 
may be more natural. We shall refrain from undertaking a 
general study of normalizations, but rather focus
on one particular, but important, case.

More precisely, we will describe the set of valuations $\cV_x$ 
normalized by the condition $\nu(x)=1$, where $x\in\fm$ and
$m(x)=1$, \ie $\{x=0\}$ defines a smooth formal curve at the origin in
$\C^2$. By going back and forth between the two normalizations
$\nu(x)=1$ and $\nu(\fm)=1$, we shall see that $\cV_x$ can be endowed
with a natural partial ordering $\le_x$ which induces a nonmetric tree
structure. We hence refer to $\cV_x$ as \emph{the relative valuative
tree}. It also has two natural parameterizations, 
the \emph{relative skewness} $\a_x$, 
and the \emph{relative thinness} $A_x$, as well as a
\emph{relative multiplicity} $m_x$. 
We shall see that the theory of the
relative valuative tree is much the same as that of its nonrelative
counterpart.

The relative valuative tree will appear in two different contexts in
this monograph. 
The first place is in Chapter~\ref{sec-puis}, where we describe
valuations in terms of Puiseux series. The relative valuative tree
appears as the quotient of the set of valuations 
on $\hk[[y]]$ where $\hat{k}$ denotes the set of Puiseux series in $x$
(Theorem~\ref{T602}); and also as a natural (metric) subtree of the
Bruhat-Tits building of $\mathrm{PGL}_2 (\C((x)))$ 
(see Section~\ref{bruhat-tits}).

In Chapter~\ref{A3}, the relative valuative tree
appears in the following situation. 
Let $\pi:X\to(\C^2,0)$ be a composition of blowups, 
and $p$ be a smooth point
on the exceptional divisor $\pi^{-1}(0)$. 
Then $\pi$ induces a map between the space of centered valuations 
on the local ring $\mathcal{O}_p$ at $p$, into 
the space of centered valuations on $R$.
If the exceptional divisor is given by $\{z=0\}$ at $p$,
then $\pi$ induces a natural isomorphism of parameterized trees
between the relative valuative tree $\cV_z$ and
its image in the valuative tree $\cV$: this image is the closure
of a set of the form $U(\vv)$ as in Section~\ref{S302}.
See Theorem~\ref{Trelnat} and Figure~\ref{F16}.
%
%
\subsection{The relative valuative tree}
Fix a smooth formal curve at the origin in $\C^2$, say $\{x=0\}$,
where $x\in\fm$ and $m(x)=1$.  We let $\cV_x$ be the set of centered
\index{relative!valuative tree}
\index{$\cV_x$ (relative valuation space/valuative tree)}
valuations $\nu:R\to\Rbar$ vanishing on $\C^*$ and such that 
$\nu(x)=1$, together with the valuation $\div_x$, defined by
$\div_x(\phi)=\max\{n\ ;\ x^n\mid\phi\}$.  
(The center of $\div_x$ is the curve $\{x=0\}$.)
Note that $\cV_x$ is compact for the weak topology.  

While it is perfectly feasible to study $\cV_x$ using SKP's normalized
by $\btilde_0=1$, we shall instead take advantage of the theory that
we have already developed for $\cV$.  Indeed, when $\nu$ is a
valuation normalized by $\nu(\fm)=1$ and different from $\nu_x$, then
$\nu/\nu(x)$ belongs to $\cV_x$. Conversely, $\nu\in
\cV_x\setminus\{\div_x\}$ implies $\nu/\nu(\fm)\in\cV$. Whence
\begin{lemma}\label{Lrel1}
  The map $N : \cV \to \cV_x$ sending $\nu$ to $\nu/\nu(x)$, and $\nu_x$
  to $\div_x$ is a homeomorphism for the weak topologies.
\end{lemma}
An immediate consequence is that, except for $\nu_x$ and $\div_x$, the
valuations $\nu$ and $N\nu$ define equivalent Krull valuations on $R$,
hence have the same invariants $\rk, \ratrk,
\trdeg$, and also the same type (see Theorem~\ref{divis}). 
We shall thus say a valuation in $\cV_x$ is
quasimonomial, divisorial, curve or infinitely singular, 
when its image in $\cV$ has the required property. 
In the sequel, we denote by $\cVqmx$ the set
of all quasimonomial valuations in $\cV_x$. 
By convention we consider
$\cVqmx$ to contain the valuation $\div_x$.

For $\nu_1,\nu_2\in\cV_x$, we write $\nu_1 \le_x\nu_2$ when
$\nu_1(\phi) \le \nu_2(\phi)$ for all $\phi\in R$. It is clear that
$\div_x$ is the unique minimal element of $\cV_x$.
\begin{lemma}\label{Lrel2}
  Let $\nu_1,\nu_2\in\cV_x\setminus\{\div_x\}$.
  Then $\nu_1\le_x\nu_2$ iff we have
  $[\nu_x , N^{-1}\nu_1]\subset [\nu_x, N^{-1}\nu_2]$ in the valuative
  tree $\cV$.
\end{lemma}
As an immediate consequence (see Section~\ref{tree-def}), we get
\begin{proposition}\label{Prel1}
  The poset $(\cV_x,\le_x)$ is a complete,
  nonmetric tree rooted at $\div_x$. 
  The map $N:\cV\to\cV_x$ is an isomorphism of (nonrooted) nonmetric
  trees.
\end{proposition}
We shall write $\nu\wedge_x\mu$ for the infimum of $\nu$ and $\mu$ 
with respect to the tree structure on $(\cV_x,\le_x)$.
\begin{remark}
  The ends (\ie the maximal elements) 
  of $\cV_x$ under $\le_x$ are exactly of
  the form $N\nu$, where $\nu\in\cV$ is not quasimonomial 
  and $\nu\ne\nu_x$. Hence $(\cVqmx,\le_x)$ is also
  a nonmetric tree rooted in $\div_x$. Recall our convention
  that $\div_x\in\cVqmx$.
\end{remark}  
\begin{remark}
  Both $\cV$ and $\cV_x$ can be viewed as subsets of the set
  $\tcV$ of centered valuations on $R$. It follows from
  Lemma~\ref{Lrel2} that
  \begin{equation*}
    \cV\cap\cV_x
    =\{\nu\in\cV\ ;\ N(\nu)=\nu\}
    =\{\nu\in\cV\ ;\ \nu\wedge\nu_x=\nu_\fm\}
    =\{\nu\in\cV_x\ ;\ \nu\ge_x\nu_\fm\}.
  \end{equation*}  
  Moreover, the partial orderings $\le$ and $\le_x$ 
  agree on $\cV\cap\cV_x$.
\end{remark}
\begin{proof}[Proof of Lemma~\ref{Lrel2}]
  Write $\mu_i=N^{-1}\nu_i$. These are valuations normalized by
  $\mu_i(\fm)=1$. 
  The condition $\nu_1\le_x\nu_2$ is equivalent to
  $\mu_1(\phi)/\mu_1(x)\le\mu_2(\phi)/\mu_2(x)$ 
  for all $\phi$. Here we may assume that $\phi\in\fm$ 
  is irreducible and apply Proposition~\ref{P201}:
  \begin{equation}\label{e446}
    \frac{\mu_1(\phi)/\mu_1(x)}{\mu_2(\phi)/\mu_2(x)}
    =\frac{\a(\mu_1\wedge\nu_\phi)\a(\mu_2\wedge\nu_x)}
    {\a(\mu_2\wedge\nu_\phi)\a(\mu_1\wedge\nu_x)}.
  \end{equation}
  
  First assume $\nu_1\le_x\nu_2$.
  We need to show that $[\nu_x,\mu_1]\subset[\nu_x,\mu_2]$. 
  By choosing $\phi$ such that $\nu_\phi\wedge\mu_i=\nu_\fm$
  for $i=1,2$ we get $\mu_1\wedge\nu_x\ge\mu_2\wedge\nu_x$ 
  from~\eqref{e446}. 
  If equality holds, then $\mu_1(x)=\mu_2(x)$, which 
  together with $\nu_1\le_x\nu_2$ gives $\mu_1\le\mu_2$.
  From this we easily conclude $[\nu_x,\mu_1]\subset[\nu_x,\mu_2]$. 
  If instead $\mu_1\wedge\nu_x>\mu_2\wedge\nu_x$, 
  then we claim that $\mu_1\in [\nu_x,\mu_2\wedge\nu_x]$. 
  Indeed, otherwise pick $\phi\in\fm$ irreducible
  such that $\nu_\phi\ge\mu_1$.
  Then $\mu_1\wedge\nu_\phi>\mu_1\wedge\nu_x$ and
  $\mu_2\wedge\nu_\phi=\mu_2\wedge\nu_x$, which
  by~\eqref{e446} contradicts $\nu_1\le_x\nu_2$.
  Thus $\mu_1\in [\nu_x,\mu_2\wedge\nu_x]$, which implies
  $[\nu_x,\mu_1]\subset[\nu_x,\mu_2]$.
  
  Conversely assume $[\nu_x,\mu_1]\subset [\nu_x,\mu_2]$. 
  We want to prove $N\mu_1\le_x N\mu_2$ \ie
  $\mu_1(\phi)/\mu_1(x)\le\mu_2(\phi)/\mu_2(x)$ for all $\phi$.  
  The assumption implies that $\mu_1\wedge\nu_x\ge\mu_2\wedge\nu_x$.
  If equality holds, then 
  $\mu_1(x)=\mu_2(x)$, and the assumption gives $\mu_1\le\mu_2$, 
  so that $N\mu_1\le_x N\mu_2$.
  If instead $\mu_1\wedge\nu_x>\mu_2\wedge\nu_x$ 
  \ie $\mu_1(x)>\mu_2(x)$, then
  the assumption implies that $\mu_1\in[\nu_\fm,\nu_x]$.
  Pick $\phi\in\fm$ irreducible.
  If $\mu_1\wedge\nu_\phi\le\mu_2\wedge\nu_\phi$, 
  then~\eqref{e446} immediately gives
  $\mu_1(\phi)/\mu_1(x)\le\mu_2(\phi)/\mu_2(x)$.
  If instead $\mu_1\wedge\nu_\phi>\mu_2\wedge\nu_\phi$, 
  then $\mu_1\wedge\nu_\phi\ge\mu_1\wedge\nu_x$,
  and $\mu_2\wedge\nu_\phi=\mu_2\wedge\nu_x$ 
  so again~\eqref{e446} results in
  $\mu_1(\phi)/\mu_1(x)\le\mu_2(\phi)/\mu_2(x)$.
  Thus $\nu_1\le_x\nu_2$, which concludes the proof.
\end{proof}
%
%
\subsection{Relative parameterizations}
We now introduce natural parameterizations on $\cV_x$, analogous to
skewness and thinness on $\cV$.  The multiplicity
$m(\phi)=\nu_\fm(\phi)$, which plays a key role on $\cV$, is here
replaced by the \emph{relative
multiplicity}\index{relative!multiplicity}\index{multiplicity!relative}
$m_x(\phi)=\nu_x(\phi)$.
By~\eqref{valu}, this number also equals the
intersection multiplicity between the curves $\{\phi=0\}$ and
$\{x=0\}$.  In particular $m_x(\phi)$ is an integer unless $x$ divides
$\phi$ in which case it is infinite.
\begin{definition}
  Pick $\nu\in\cV_x$. We define
  \index{parameterization!of the relative valuative tree}
  \begin{itemize}
  \item
    \emph{the relative skewness}:\index{relative!skewness}
    $\a_x(\nu):=\sup_{\phi\in\fm}\frac{\nu(\phi)}{m_x(\phi)}\in\Rbar$;
  \item
    \emph{the relative multiplicity}:
    \index{relative!multiplicity}
    \index{$m_x(\nu)$ (relative multiplicity)}
    if $\nu$ is infinitely singular, then $m_x(\nu):=\infty$;
    otherwise 
    $m_x(\nu):=\min\{m_x(\phi)\ ;\ \phi\in\fm\
    \text{irreducible},\ N\nu_\phi\ge_x\nu\}\in\N$;
  \item  
    \emph{the relative thinness}:\index{relative!thinness}
    $A_x(\nu):=1+\int_{\nu_x}^\nu m_x(\mu)\,d\a_x(\mu)\in\Rbar$.
  \end{itemize}
\end{definition}
Notice that $\a_x(\div_x)=0$, $A_x(\div_x)=1$ and $m_x(\div_x)=1$.
\begin{proposition}\label{Prel2}
  Pick $\nu\in\cV\setminus\{\nu_x\}$ and let $N\nu$ be its image
  in $\cV_x$. Then
  \begin{equation*}
   \a_x(N\nu)=\a(\nu)/\nu(x)^2
   \qand
   A_x(N\nu)=A(\nu)/\nu(x).
  \end{equation*}    
  The relative multiplicity is given by
  \begin{equation*}
    m_x(N\nu)
    =\begin{cases}
      1 & \text{ when } \nu\in[\nu_x,\nu_\fm],\\
      m(\nu)\,\nu(x) & \text{ otherwise}. 
    \end{cases}
  \end{equation*}
\end{proposition}
In particular, $\a_x$ and $A_x$ are finite on quasimonomial valuations, 
and $m_x(\nu)$ is infinite iff $\nu$ is infinitely singular. 
We leave it to the
reader to generalize Theorem~\ref{T401} and
Propositions~\ref{P402},~\ref{P306}. In particular we have
\begin{itemize}
\item[(i)]
  $A_x(\nu)\ge1+\a_x(\nu)$ with equality iff $m_x(\nu)=1$, 
  \ie if $\nu$ is monomial in some local coordinates $(x,y)$;
\item[(ii)]
  if $m_x(\nu)>1$, then $A_x(\nu)<m_x(\nu)\a_x(\nu)$.
\end{itemize}

The main consequence of the preceding proposition is
\begin{corollary}\label{Crel1}
  The relative skewness defines parameterizations
  \begin{equation*}
    \a_x:\cV_x\to[0,\infty]
    \qand
    \a_x:\cVqmx\to[0,\infty[\,.
  \end{equation*}
  Similarly, the relative thinness gives parameterizations
  \begin{equation*}
    A_x:\cV_x\to[1,\infty]
    \qand
    A_x:\cVqmx\to[1,\infty[\,.
  \end{equation*}
\end{corollary}
As in the nonrelative case, we can use
these parameterizations to define metrics
on $\cV_x$ and $\cVqmx$.
\begin{proof}[Proof of Proposition~\ref{Prel2}]
  We fix $\nu\in\cV\setminus\{\nu_x\}$. 
  First consider the relative skewness.
  The supremum defining
  $\a_x(N\nu)$ can be taken over irreducible elements $\phi\in\fm$. 
  For any such $\phi$, Proposition~\ref{P201} yields
  \begin{equation}\label{e445}
    \frac{(N\nu)(\phi)}{m_x(\phi)}
    =\frac{\nu(\phi)}{\nu(x)\nu_x(\phi)}
    =\frac{\a(\nu\wedge\nu_\phi)m(\phi)}
    {\nu(x)\a(\nu_\phi\wedge\nu_x)m(\phi)}
    =\frac{\a(\nu\wedge\nu_\phi)\a(\nu\wedge\nu_x)}
    {\a(\nu_\phi\wedge\nu_x)\nu(x)^2}.
  \end{equation}

  If $\nu$ is divisorial or $\nu\notin\,]\nu_x,\nu_\fm]$, then
  we can choose $\phi\in\fm$ irreducible with $\nu_\phi\ge\nu$
  and $\nu_\phi\wedge\nu_x=\nu\wedge\nu_x$.
  Then $\nu\wedge\nu_\phi=\nu$ so~\eqref{e445} immediately gives
  $(N\nu)(\phi)/m_x(\phi)\ge\a(\nu)/\nu(x)^2$, 
  yielding the lower bound $\a_x(N\nu)\ge\a(\nu)/\nu(x)^2$. 
  If $\nu\in\,]\nu_x,\nu_\fm[$ is irrational, then we get
  the same lower bound by a simple approximation argument.
  The corresponding upper bound also follows from~\eqref{e445}. 
  Indeed, we always have either 
  $\nu\wedge\nu_x\le\nu_\phi\wedge\nu_x$ or
  $\nu\wedge\nu_\phi\le\nu_\phi\wedge\nu_x$,
  so $(N\nu)(\phi)/m_x(\phi)\le\a(\nu)/\nu(x)^2$, implying
  $\a(N\nu)\le\a(\nu)/\nu(x)^2$.

  \smallskip
  Let us now relate multiplicity to relative multiplicity. 
  We may assume that $\nu$ is not infinitely singular, 
  since otherwise $m_x(N \nu)=m(\nu)=\infty$.
  
  If $\nu\in[\nu_x,\nu_\fm]$ then $N\nu\le_xN\nu_y$ 
  whenever $(x,y)$ are local coordinates.  
  Hence $m_x(N\nu)\le m_x(y)=\nu_x(y)=1$, so
  $m_x(N\nu)=1$.
  
  If $\nu\not\in[\nu_x,\nu_\fm]$ then it follows from
  Proposition~\ref{Prel1} that for $\phi\in\fm$ irreducible, the
  conditions $N\nu_\phi\ge_x N\nu$ and $\nu_\phi\ge\nu$ are equivalent.
  Moreover, they both imply $m_x(\phi)=m(\phi)\nu(x)$.
  Hence $m_x(N\nu)=m(\nu)\nu(x)$.
  
  \smallskip
  To relate thinness to relative thinness we make use of
  the preceding computations. Write $\nu'=\nu\wedge\nu_x$.
  Then 
  \begin{multline*}
    A_x(N\nu)
    =1+\int_{\div_x}^{N\nu'}m_x(\tilde{\mu})\,d\a_x(\tilde{\mu}) 
    +\int_{N\nu'}^{N\nu}m_x(\tilde{\mu})\,d\a_x(\tilde{\mu})\\
    =1+\a_x(N\nu')
    +\int_{\nu'}^{\nu}m(\mu)\mu(x)\,\frac{d\a(\mu)}{\mu(x)^2}\\
    =1+\a(\nu')^{-1}+\frac{A(\nu)-A(\nu')}{\nu'(x)}
    =\frac{A(\nu)}{\nu'(x)}=\frac{A(\nu)}{\nu(x)}.
  \end{multline*}
  This completes the proof.
\end{proof}
\begin{proof}[Proof of Corollary~\ref{Crel1}]
  It suffices to consider the full tree $\cV_x$; the statements
  about $\cVqmx$ will be direct consequences.

  If $t\ge 1$, then $\a_x(N\nu_{x,t})=t/t^2=t^{-1}$.  
  Hence the restriction of $\a_x$ to the segment 
  $[\div_x,\nu_\fm]$ gives an
  order-preserving bijection onto the real interval $[0,1]$.  
  Any end in $\cV_x$ different from $\div_x$ is of the form
  $N\nu$, where $\nu\ne\nu_x$ is an end in $\cV$.
  Let us prove that $\a_x$ defines
  an order-preserving bijection of $[\div_x,N\nu]$ onto
  $[0,\a_x(N\nu)]$.  Set $\nu'=\nu\wedge\nu_x$.
  By the above calculation,
  $\a_x$ gives an order-preserving bijection of $[\div_x,N\nu']$ onto
  $[0,\a_x(N\nu')]$. 
  On the other hand, if $\mu\in[\nu',\nu]$, then
  $\mu(x)=\nu'(x)$ so $\a_x(N\mu)=\a(\mu)/\nu'(x)^2$. 
  It hence follows from Theorem~\ref{T401} that $\a_x$ gives an 
  order-preserving mapping of $[N\nu',N\nu]$ onto 
  $[\a_x(N\nu'),\a_x(N\nu)]$.
  
  This shows that $\a_x$ defines a parameterization of the rooted
  nonmetric tree $(\cV_x,\le_x)$. 
  The fact that relative thinness gives a parameterization is proved
  in the same way. Note that $A_x(N\nu_{x,t})=1+t^{-1}$ for $t\ge1$.
\end{proof}
Finally we have the following calculation which will be
needed in Chapter~\ref{sec-puis}.
The proof is left to the reader.
\begin{lemma}\label{L437}
  If $\nu,\nu'\in\cV_x$, 
  $\nu'<_x\nu$ and $m_x(\nu)=m_x(\nu')$,
  then $A_x(\nu)-A_x(\nu')=\nu(\phi)-\nu'(\phi)$ 
  for any irreducible $\phi\in\fm$ such that
  $N\nu_\phi>_x\nu'$ and $m_x(\phi)=m_x(\nu')$.
\end{lemma}
%
%
\subsection{Balls of curves}
\index{balls of curves}
The relative valuative tree can also be understood in terms of balls
of irreducible curves, just as in Section~\ref{val-curve}.
\begin{definition}
  As before, let $\cC$ be the space of (local formal) 
  irreducible curves. 
  Define $\cC_x=\cC\setminus\{x=0\}$.
  \index{$\cC_x$ (space of irreducible curve, relative)}
  The \emph{relative multiplicity}
  \index{relative!multiplicity}
  \index{curve!relative multiplicity of}
  \index{$m_x(C)$ (relative multiplicity)}
  \index{$m_x(\phi)$ (relative multiplicity)}
  of a curve $C\in\cC_x$ is defined by
  $m_x(C)=C\cdot\{x=0\}$.
  Thus $m_x(\{\phi=0\})=m_x(\phi)$ for 
  $\phi\in\fm$ irreducible.
  Given $C_1,C_2\in\cC_x$ set 
  \begin{equation*}
    d_{\cC_x}(C_1,C_2)
    :=\frac{m_x(C_1)m_x(C_2)}{C_1\cdot C_2}.
  \end{equation*}
\end{definition}
We leave to the reader to check the analogue of Lemma~\ref{ballmetric}:
\begin{lemma}
  For any $C_1,C_2\in\cC_x$ we have
  \begin{equation}\label{e418}
    d_{\cC_x}(C_1,C_2)
    =\frac1{\alpha_x(N\nu_{C_1}\wedge_xN\nu_{C_2})}.
  \end{equation}
  Moreover, $(\cC_x,d_{\cC_x})$ is an
  ultrametric space of infinite diameter.
\end{lemma}
As in the nonrelative case we can associate a 
parameterized tree $\cT_{\cC_x}$ to the ultrametric space 
$(\cC_x,d_{\cC_x})$. 
Let us review this construction, pointing out the differences
that we encounter in the relative setting.

A point in $\cT_{\cC_x}$ is a closed ball in 
$\cC_x$ of positive radius. 
The partial ordering is given by reverse inclusion. 
In order to have a minimal element in $\cT_{\cC_x}$ we also
add the ball of infinite radius, \ie the whole space $\cC_x$. 
The parameterization on $\cT_{\cC_x}$ sends a ball of radius 
$r$ to the real number $r^{-1}$.

There is a natural multiplicity function on $\cT_{\cC_x}$.
Its value on a ball is the minimum relative
multiplicity of any curve in the ball.

The completion of $\cT_{\cC_x}$ is taken 
in the sense of nonmetric trees.
An element in the completion is represented by a decreasing 
sequence of balls. The multiplicity function
and the parameterization on $\cT_{\cC_x}$ 
both extend naturally to the completion.
An element in the completion, represented by a 
decreasing sequence of balls,
has finite multiplicity iff the sequence has nonempty 
intersection in $\cC_x$.

We can now define a map from $\cVqmx$ to $\cT_{\cC_x}$ by sending 
a quasimonomial valuation $\nu$ to the ball 
$B_\nu:=\{C\ ;\ N\nu_C>_x\nu\}$.
(That this is a ball follows from~\eqref{e418}.)
Notice that $\div_x$ is sent to the whole space $\cC_x$.

Similarly, $B\mapsto\nu_B:=\wedge_x\{N\nu_C\ ;\ C\in B\}$ 
gives a map from $\cT_{\cC_x}$ to $\cVqmx$.  
A straightforward adaptation of the proof of 
Theorem~\ref{treeballs} gives
\begin{theorem}
  The mappings $\nu\mapsto B_\nu$ and $B\mapsto\nu_B$ 
  preserve multiplicity and give inverse isomorphisms 
  between the parameterized 
  trees $(\cVqmx,\a_x)$ and $(\cT_{\cC_x},r^{-1})$.
  They also extend uniquely to isomorphisms 
  between $\cV_x$ and the completion of $\cT_{\cC_x}$.
\end{theorem}
%
%
%
\subsection{Homogeneity}
Recall that we defined $\tcV$ as the set of \emph{all} centered
valuations $\nu:R\to\Rbar$ without any normalization.
This space $\tcV$ is naturally a bundle with 
fibers that are rays of the form $\R_+^*\nu$, $\nu\in\tcV$.
The trees $\cV$ and $\cV_x$, consisting of elements 
$\nu\in\tcV$ satisfying $\nu(\fm)=1$ and $\nu(x)=1$,
respectively, can then be viewed as (the images of) sections
of this bundle.

Suppose we wish to define functions skewness and thinness
on the full space $\tcV$ in such a way that the 
restrictions to $\cV$ and $\cV_x$ recover skewness (thinness)
and relative skewness (relative thinness), respectively.
We also want these functions to be homogeneous
on each ray $\R_+^*\nu$.
Proposition~\ref{Prel2} indicates that we should
then impose $\a(t\nu)=t^2\a(\nu)$, and $A(t\nu)=tA(\nu)$ 
for $\nu\in\tcV$ and $t>0$. 
In other words, \emph{skewness is homogeneous of degree two,
and thinness is homogeneous of degree one}.
Let us give further evidence that these are the correct
degrees.

First, as noted in Remark~\ref{rem-vol}, the skewness $\a(\nu)$
of a valuation $\nu\in\cV$ is 
equal to the inverse of its volume 
$\vol(\nu)\=\lim_{c\to\infty}\frac{2}{c^2}\dim_\C R/\{\nu\ge c\}$. 
It is clear that $\vol(t\nu)=t^{-2}\vol(\nu)$. 
Hence keeping $\a=\vol^{-1}$ on all of $\tcV$
implies the required homogeneity property. 

Alternatively, we shall see in Section~\ref{S416}
that there is a natural intersection product on $\cV$
given by $\nu\cdot\nu'\=\a(\nu\wedge\nu')$.
In particular $\a(\nu)=\nu\cdot\nu$.
By extending this intersection product in a bilinear way,
restricting to $\tcV$, and still requiring 
$\a(\nu)=\nu\cdot\nu$, we again get
$\a(t\nu)=t\nu\cdot t\nu=t^2\a(\nu)$.

As for thinness, we rely on Remark~\ref{Rjaco}.
Consider $\nu$ divisorial. Let $\fm_\nu$ be the
maximal ideal of the valuation ring $R_\nu$, and let
$ J_{R_\nu/R}$ be the Jacobian ideal. 
If $\nu\in\cV$ is normalized by $\nu(\fm)=1$, we have
seen that $A(\nu)=\nu(J_{R_\nu/R}\,\fm_\nu)$.  Taking this as a
definition for thinness, we get $A(t\nu)=tA(\nu)$.  Indeed,
multiplying $\nu$ by a positive constant does not change the ideals
$\fm_\nu$ and $J_{R_\nu/R}$.


%
%
%
%
%
%
\chapter{Valuations through Puiseux series}\label{sec-puis}
We now present an alternative approach, based on Puiseux series, for
classifying valuations and obtaining the tree structure on the
valuative tree $\cV$, or more precisely the relative valuative tree
$\cV_x$ as defined in Section~\ref{sec-relative}.

Fix a smooth formal curve at the origin, say $\{x=0\}$.  Then any
local irreducible curve, save for $\{x=0\}$, is represented by a
(non-unique) Puiseux series in $x$. In the same spirit, we show that
every valuation in $\cV_x$ is represented by a valuation on the power
series ring in one variable with coefficients that are Puiseux series
in $x$.  Moreover, the set $\hcV$ of all such (normalized) valuations
has a natural tree structure and the restriction map from $\hcV$ to
$\cV_x$ is a tree morphism. In fact, $\cV_x$ is naturally the orbit
space of $\hcV$ under the action by the relevant Galois group.

This Puiseux approach is closely related to Berkovich's theory of
analytic spaces and we show how the valuative tree naturally embeds
(as a nonmetric tree) as the closure of a disk in the Berkovich
projective line $\bbP^1(k)$ over the local field $k=\C((x))$.
Moreover, $\cVqmx$ embeds as a subtree of the Bruhat-Tits
building of $\mathrm{PGL}_2$. The latter space has a natural metric
and we show that the induced metric on $\cVqmx$ is
exactly the relative thin metric.

As references for this chapter we point out~\cite{teissier-puis} for
Puiseux expansions, and~\cite{Ber} and~\cite{BT} for Berkovich spaces
and Bruhat-Tits buildings, respectively. In fact, many of the
results presented here are contained, at least implicitly, 
in~\cite{Ber} but are presented there in a different language and
with less details.
We also refer to~\cite{Rivera} for a concrete exposition in a context
similar to ours.
%
%
%
%
\section{Puiseux series and valuations}\label{sec-puiseux}
We fix an element of $R$ representing a smooth formal curve.
In other words, fix $x\in\fm$ with multiplicity $m(x)=1$.
Think of $x$ as a variable
and let $k=\C((x))$ be the fraction field of $\C[[x]]$, 
\ie the field of Laurent series in $x$.
Let $\hat{k}$ be the algebraic closure of $k$: its elements are
finite or infinite \emph{Puiseux series} of the form
\begin{equation}\label{puis-exp}
  \hat{\phi}=\sum_{j\ge1}a_jx^{\hbeta_j}
  \quad\text{with}\quad
  a_j\in\C^*,\ \hbeta_{j+1}>\hbeta_j\in\Q
\end{equation}
and where the rational numbers $\hbeta_j$ have bounded denominators,
\ie $m\hbeta_j\in\Z$ for all $j$ for some integer $m$. 
Notice that this implies that if the series is infinite, then
$\hbeta_j\to\infty$ as $j\to\infty$.

We endow $k$ and $\hk$ with the valuation $\nu_\star$ defined by
$\nu_\star|_{\C^*}=0$, and $\nu_\star(x)=1$, 
and we let $\bk$ be the completion of $\hk$
with respect to $\nu_\star$. 
The elements of $\bk$
are finite or infinite series of the form~\eqref{puis-exp},
with the restriction that $\hbeta_j\to\infty$ if the series is
infinite, but without the restriction on the denominators.
Shortly we will consider even more general series, dropping the 
condition that $\hbeta_j\to\infty$. Clearly $\nu_\star$ extends
to series of that type: we have $\nu(\phi)=\hbeta_1$.

We define $\hk_+$ to be the subset of $\hk$ consisting of 
Puiseux series $\hphi$ with $\nu_\star(\hphi)>0$, \ie 
$\hbeta_1>0$ in~\eqref{puis-exp}. Similarly we define
$k_+$ and $\bk_+$.

Fix $y\in\fm$ such that $(x,y)$ are local formal coordinates 
at $\C^2$, \ie $\nu_x\wedge\nu_y=\nu_\fm$. 
Think of $y$ as a variable and consider the ring
$\bk[[y]]$ of formal power series with 
Puiseux series coefficients.

We are interested in valuations $\hnu:\bk[[y]]\to\Rbar$
extending $\nu_\star$ on $\bk$ and satisfying $\hnu(y)>0$.
As in Proposition~\ref{extend}, any such valuation is
determined by its values on polynomials in $y$ and---since 
$\overline{k}$ is algebraically closed---in fact on
its values on linear polynomials $y-\hpsi$, $\hpsi\in\bk$.
\begin{proposition}\label{thm-puiseux}
  Any valuation $\hnu:\bk[[y]]\to\Rbar$ extending $\nu_\star$ on $\bk$ 
  and satisfying $\hnu(y)>0$ can be uniquely represented 
  by a number $\hbeta\in\Rbar$ and a finite or infinite series 
  $\hphi$ of the form~\eqref{puis-exp} 
  (with no restriction on $\lim\hbeta_j$ if the series is infinite) 
  such that $\hbeta>\hbeta_j>0$ for all $j$,
  and $\hbeta=\lim\hbeta_j$ if the series is infinite.
  More precisely, we have
  \begin{equation}\label{val-puis}
    \hnu(y-\hpsi)
    =\min\{\hat{\beta},\nu_\star(\hpsi-\hphi)\}
    \quad\text{for all $\hpsi\in\bk$}.
  \end{equation}

  Conversely, if $\hphi$ and $\hbeta$ satisfy these conditions
  then there exists a unique valuation $\hnu$ on $\bk[[y]]$
  satisfying~\eqref{val-puis} and extending $\nu_\star$ on $\bk$ 
\end{proposition}
\begin{definition}\label{def-puis}
  We write $\val[\hphi;\hbeta]$ for the valuation 
  $\hnu:\bk[[y]]\to\Rbar$ determined by the pair 
  $(\hphi,\hbeta)$ as in Proposition~\ref{thm-puiseux}.
  Moreover, we say that $\hnu$ is
  \begin{itemize}
  \item[(i)]
    of \emph{point type} 
    if $\hphi\in\bk$ and $\hbeta=\infty$;
    we then write $\hnu=\hnu_\hphi$;
  \item[(ii)]
    of \emph{finite type}
    if $\hphi\in\hk$ and $\hbeta<\infty$;
  \item[(iii)]
    \emph{rational}
    if $\hnu$ is of finite type and $\hbeta$ is rational;
  \item[(iv)]
    \emph{irrational}
    if $\hnu$ is of finite type and $\hbeta$ is irrational;
  \item[(v)]
    of \emph{special type}
    if $\hphi\notin\bk$.
  \end{itemize}
\end{definition}
\begin{remark}
  The proof below shows that if $\hnu=\val[\hphi;\hbeta]$
  is of point type or special type then
  \begin{equation}\label{e425}
    \hnu(y-\psi)=\nu_\star(\hpsi-\hphi)
    \quad\text{for all $\hpsi\in\bk$},
  \end{equation}
  exhibiting the fact that $\hbeta$ is determined by 
  $\hphi$ in this case.
\end{remark}
\begin{proof}[Proof of Proposition~\ref{thm-puiseux}]
  Pick a valuation $\hnu$ on $\bk[[y]]$ extending $\nu_\star$
  and satisfying $\hnu(y)>0$. Let us construct $\hphi$ and
  $\hbeta$ such that~\eqref{val-puis} holds.

  Set $\hbeta_1:=\hnu(y)>0$. 
  Then $\hnu(y-\hpsi)\ge\min\{\hbeta_1,\nu_\star(\hpsi)\}$
  for all $\hpsi\in\bk$.
  There are now two possibilities: 
  either equality holds for all $\hpsi$,
  or $\hbeta_1\in\Q$ and there exists a 
  unique $\theta_1\in\C^*$ such that $\hnu(y-\hphi_1)>\hbeta_1$,
  where $\hphi_1:=a_1x^\hbeta_1$.
  In the first case we stop, 
  setting $\hphi=\hphi_1$ and $\hbeta=\hbeta_1$.
  In the second case we let
  $\hbeta_2:=\hnu(y-\hphi_1)$ and continue.
  Inductively we construct
  $\hphi_j=\hphi_{j-1}+\theta_jx^{\hbeta_j}$ 
  with $\theta_j\in\C^*$ and 
  $\hbeta_{j+1}>\hbeta_{j}\in\Q$, such that 
  $\hnu(y-\hpsi)\ge\min\{\hbeta_{j+1},\nu_\star(\hpsi-\hphi_j)\}$.
  If this process stops at some finite step $j_0$, then we set
  $\hphi=\hphi_{j_0}=\sum_1^{j_0}\theta_jx^{\hbeta_j}$, 
  $\hbeta=\hbeta_{j_0+1}$.
  Otherwise $(\hbeta_j)_1^\infty$ is strictly 
  increasing with limit $\hbeta$ and 
  $\hphi_j$ converges to a series $\hphi$ of the 
  type~\eqref{puis-exp}, where $\lim\hbeta_j=\hbeta$ may or
  may not be finite.
  In fact $\hphi\in\bk$ iff $\hbeta=\infty$.

  We claim that~\eqref{val-puis} holds. 
  In the case when $\hphi$ is a finite series, \ie the procedure
  above stops after finitely many steps, 
  then this is clear by construction.
  If $\hphi$ is an infinite series, then for 
  $\hpsi\in\bk$ we have
  $\hnu(y-\hpsi)\ge\min\{\hbeta_{j+1},\nu_\star(\hpsi-\hphi_j)\}$
  for each $j$, with equality unless 
  $\nu_\star(\hpsi-\hphi_{j+1})>\hbeta_{j+1}$.
  Passing to the limit we obtain 
  $\hnu(y-\hpsi)\ge\min\{\hbeta,\nu_\star(\hpsi-\hphi)\}$
  and strict inequality would imply that 
  $\nu_\star(\hpsi-\hphi_j)>\hbeta_j$ for all $j$.
  If $\hphi\in\bk$, \ie $\hbeta=\infty$, then 
  this gives $\nu_\star(\hpsi-\hphi)=\infty$ so that
  $\hpsi=\hphi$ and~\eqref{val-puis} holds.
  If $\hphi\not\in\bk$, \ie $\hbeta<\infty$,
  then $\nu_\star(\hpsi-\hphi_j)>\hbeta_j$ for all $j$
  leads to a contradiction, since 
  $\hpsi\in\bk$ but $\hphi\not\in\bk$.
  Thus~\eqref{val-puis} holds also in this case.

  Notice that $\hbeta$ and $\hphi$ are unique. 
  Indeed, it follows from~\eqref{val-puis} that 
  the number $\hbeta$ satisfies
  $\hbeta=\sup\{\hnu(y-\hpsi)\ ;\ \hpsi\in\bk\}$
  (compare Definition~\ref{D404} below).
  Moreover, suppose that $\hphi=\sum a_jx^{\hbeta_j}$ 
  and $\hphi'=\sum a'_jx^{\hbeta'_j}$ 
  both satisfy~\eqref{val-puis} and that $\hbeta_j,\hbeta'_j<\hbeta$
  for all $j$. Then $\nu_\star(\hphi-\hphi')\ge\hbeta$, so
  $\hphi=\hphi'$.
  
  Finally we show that given $\hphi$ and $\hbeta$ satisfying the
  conditions in the proposition, there exists a unique valuation
  $\hnu$ satisfying~\eqref{val-puis}. 
  Uniqueness is immediate from~\eqref{val-puis} since
  $\hnu$ is determined on its values on linear polynomials
  $y-\hpsi$.
  As for existence, let us define $\hnu:\bk[y]\to\Rbar$ 
  using~\eqref{val-puis}, $\hnu|_\bk=\nu_\star$
  and $\hnu(\prod(y-\hpsi_i))=\sum\hnu(y-\hpsi_i)$.
  The nontrivial part is to show that $\hnu$ is 
  actually a valuation on $\bk[y]$, \ie that
  \begin{equation}\label{e422}
    \hnu(P+Q)\ge\min\{\hnu(P),\hnu(Q)\}
    \quad\text{for all $P,Q\in\bk[y]$}. 
  \end{equation}
  Once we have proved~\eqref{e422}, the fact that
  $\hnu$ extends uniquely to a valuation on the power series
  ring $\bk[[y]]$ is proved as in Proposition~\ref{extend}.
  
  To prove~\eqref{e422} let us rephrase the defining 
  property~\eqref{val-puis} somewhat. 
  First assume that $\hphi$ is a finite Puiseux series and
  that $\hbeta\in\Q$. 
  For $a\in\C$ define $\hphi_a=\hphi+ax^\hbeta$. 
  Then $\hphi_a\in\bk$ and from~\eqref{val-puis} it follows that
  there exists a subset $X(\psi)$ of $\C$ containing at most one
  element such that if $a\not\in X(\psi)$, then
  $\hnu(y-\hpsi)=\nu_\star(\hpsi-\hphi_a)$.
  As a consequence, given $P\in\bk[y]$ there exists
  a finite subset $X(P)$ such that
  $\hnu(P)=\nu_\star(P(\hphi_a))$ for $a\not\in X(P)$.
  Here $P(\hphi_a)$ is the element of $\bk$ obtained by
  substituting $\hphi_a$ for $y$ in the polynomial $P\in\bk[y]$.
  If $P,Q\in\bk[y]$ and $a\not\in X(P)\cup X(Q)\cup X(P+Q))$,
  we then get
  \begin{equation*}
    \hnu(P+Q)
    =\nu_\star(P(\hphi_a)+Q(\hphi_a))
    \ge\min\{\nu_\star(P(\hphi_a)),\nu_\star(Q(\hphi_a))\}
    =\min\{\hnu(P),\hnu(Q)\},
  \end{equation*}
  so that $\hnu$ is a valuation in this case.

  The same argument works also if $\hphi$ is a finite
  Puiseux series and $\hbeta$ is irrational: we simply
  replace $\hphi_a$ by $\hphi+x^\hbeta$ in the argument above.
  Of course $\hphi+x^\hbeta$ is not an element of $\bk$,
  but the computations still make sense. 
  If $\hphi$ is an infinite Puiseux series and
  $\hbeta=\infty$ (so that $\hphi\in\bk$), 
  then~\eqref{val-puis} instead implies 
  $\hnu(P)=\nu_\star(P(\hphi))$ for all $P\in\bk[y]$,
  which gives~\eqref{e422} by the same argument as above.
  This method in fact also works if $\hphi\not\in\bk$ and 
  $\hbeta<\infty$.

  Thus $\hnu$ is a valuation, completing the proof.
\end{proof}
%
%
%
%
\section{Tree structure}\label{S11}
Let $\hcV$ be the set of valuations $\hnu:\bk[[y]]\to\Rbar$ extending
$\nu_\star$ on $\bk$ and satisfying $\hnu(y)>0$.  To $\hcV$ we also
add the valuation $\hnu_\star$ defined by $\hnu_\star(y-\hphi)=0$ for
all $\hphi\in\bk$ with $\nu_\star(\phi)\ge0$. Note that by definition
this valuation is a rational valuation of finite type. Our objective
now is to show that $\hcV$ carries a tree structure similar to that of
the relative valuative tree $\cV_x$
discussed in Section~\ref{sec-relative}.

Recall that the main technical tool for obtaining the tree structure
on $\cV_x$ was the technique of SKP's. Here, in the case of $\hcV$,
the situation is a bit simpler and 
Proposition~\ref{thm-puiseux} provides the technical result that we
need.
%
%
\subsection{Nonmetric tree structure}
There is a natural partial ordering $\le$ on $\hcV$:
$\hnu\le\hmu$ iff $\hnu\le\hmu$ pointwise as functions on $\bk[y]$. 
As $\bk$ is algebraically closed, it suffices to check this condition
on linear polynomials in $y$.
Also notice that if $\hnu\in\cV_x$, $\hphi\in\bk$ and 
$\nu_\star(\hphi)\le0$, then $\hnu(y-\hphi)=\nu_\star(\hphi)$.
Hence
\begin{equation}\label{e424}
  \hnu\le\hmu
  \quad\text{iff}\quad
  \hnu(y-\hpsi)\le\hmu(y-\hpsi)
  \ \text{for every $\hpsi\in\bk_+$}. 
\end{equation}
\begin{proposition}\label{P428}
  The partial ordering $\le$ gives $\hcV$ the structure of a
  complete nonmetric tree rooted at $\hnu_\star$. 
  Its ends are the valuations of point type or special type.
  Its branch points are the rational valuations.
  Its regular points are the irrational valuations.
\end{proposition}
\begin{corollary}
  Write $\hcVfin$ for the set of valuations in $\hcV$ of finite type.
  Then $\hcVfin$ is a rooted subtree of $\hcV$ containing
  none of its ends (apart from the root $\hnu_\star$).
\end{corollary}
We denote by $\wedge$ the infimum in $\hcV$ with respect to
the partial ordering above.
\begin{proof}[Proof of Proposition~\ref{P428}]
  The key to the proof is the following characterization of
  the partial ordering using the structure result above.
  \begin{lemma}\label{L424}
    Consider $\hnu_1,\hnu_2\in\hcV$ and write
    $\hnu_i=\val[\hphi_i,\hbeta_i]$ as in 
    Definition~\ref{def-puis}. 
    Then $\hnu_1\le\hnu_2$ iff $\hbeta_2\ge\hbeta_1$ and
    $\nu_\star(\hphi_2-\hphi_1)\ge\hbeta_1$.
  \end{lemma}
  We postpone the proof of the lemma.
  Let us show that $(\hcV,\le)$ satisfies the axioms~(T1)-(T3) for 
  a nonmetric tree. As for~(T1) this is easy as $\hnu_\star$ is
  clearly the unique minimal element of $\hcV$.
  Now consider~(T2). 
  Fix $\hnu\in\hcV$, $\hnu\ne\hnu_\star$. 
  We have to show that $\{\hmu\in\hcV\ ;\ \hmu\le\hnu\}$ 
  is a totally ordered set isomorphic to a real interval.
  Write $\hnu=\val[\hphi,\hbeta]$ as in Definition~\ref{def-puis}.
  Thus $\hbeta>0$ and $\hphi=\sum_{j=1}^na_jx^{\hbeta_j}$,
  where $1\le n\le\infty$.
  Write $\hphi_i=\sum_{j=1}^ia_jx^{\hbeta_j}$ for $1\le i<n+1$.
  By Lemma~\ref{L424} we have $\hmu\le\hnu$ iff
  $\hmu=\val[\hphi_i;\halpha]$ for some $i$ and
  $\hbeta_i<\halpha\le\hbeta_{i+1}$.
  Here we adopt the convention that $\hbeta_0=0$ and that
  $\hbeta_{n+1}=\hbeta$ if $n<\infty$.
  Thus $\{\hmu\in\hcV\ ;\ \hnu_\star<\hmu\le\hnu\}$ is isomorphic 
  to the union of the intervals $]\hbeta_i,\hbeta_{i+1}]$, hence
  to the real interval $]0,\hbeta]$.
  This proves~(T2).

  As for~(T3), it is easier to prove the equivalent statement~(T3') 
  given in Remark~\ref{R401}. Moreover~(T3') clearly follows if
  we can prove that every totally ordered subset of $\hcV$
  has a majorant in $\hcV$. This will in fact also prove that $\hcV$
  is a complete tree. 
  Thus consider such a totally ordered subset $\cS\subset\hcV$.
  By Lemma~\ref{L424} the Puiseux series defining 
  $\hnu\in\cS$ has a length that is a nondecreasing function of $\nu$.
  The parameter $\hbeta$ is also an increasing function on $\cS$.
  It follows easily from this and from Lemma~\ref{L424}
  that there exists a valuation in $\hcV$ of point type or
  special type dominating all the $\hnu\in\cS$.
  Thus $\hcV$ is a complete tree.

  It is clear from Lemma~\ref{L424} that valuations in $\hcV$ of point 
  type or special type are maximal elements, hence ends. 
  Similarly, a valuation of finite type, say $\hnu=\val[\hphi,\hbeta]$
  with $\hphi\in\hk$ and $\hbeta<\infty$ cannot be an end as we
  obtain a valuation dominating it by increasing $\hbeta$.

  Finally, let us show that rational valuations are branch points and
  irrational valuations regular points in $\hcV$. 
  First consider $\hnu=\val[\hphi,\hbeta]$ rational, say
  $\hphi=\sum_1^{j_0}a_jx^{\hbeta_j}$ and $\Q\ni\hbeta>\hbeta_{j_0}$.
  For $a\in\C^*$ define $\phi_a=\phi+ax^\hbeta$ and 
  $\hnu_a=\val[\phi_a;\infty]$. Then $\hnu_a>\hnu$ for all $a$
  but if $a\ne b$, then $\hnu_a\wedge\hnu_b=\hnu$. 
  Thus $\hnu$ is a branch point.
  Now suppose $\hnu$ is irrational, with
  $\hphi=\sum_1^{j_0}a_jx^{\hbeta_j}$ and $\Q\not\ni\hbeta>\hbeta_{j_0}$.
  If $\hnu_1,\hnu_2>\hnu$, then $\hnu_i=\val[\phi+\phi_i;\halpha_i]$
  where $\halpha_i>\hbeta$ and $\hnu_\star(\phi_i)>\hbeta$, $i=1,2$.
  Pick $\halpha$ such that $\halpha>\hbeta$ but $\halpha<\halpha_i$,
  and define $\hmu=\val[\hphi;\halpha]$. 
  Then $\hnu<\hmu<\hnu_i$ for $i=1,2$ so that $\hnu_1$ and $\hnu_2$
  define the same tree tangent vector at $\hnu$.
\end{proof}
\begin{proof}[Proof of Lemma~\ref{L424}]
  First suppose that $\hbeta_2\ge\hbeta_1$ and 
  $\nu_\star(\hphi_2-\hphi_1)\ge\hbeta_1$.
  Let us show that $\hnu_1\le\hnu_2$.
  It suffices to show that
  $\hnu_1(y-\hpsi)\le\hnu_2(y-\hpsi)$ for every $\hpsi\in\bk$.
  Fix $\hpsi\in\bk$.
  By definition, 
  $\hnu_i(y-\hpsi)=\min\{\hbeta_i,\nu_\star(\hphi_i-\hpsi)\}$
  for $i=1,2$.
  The assumptions now give 
  $\hbeta_2\ge\min\{\hbeta_1,\nu_\star(\hphi_1-\hpsi)\}$ and
  \begin{equation*}
    \nu_\star(\hphi_2-\hpsi)
    \ge\min\{\nu_\star(\hphi_1-\hphi_2),\nu_\star(\hphi_1-\hpsi)\}
    \ge\min\{\hbeta_1,\nu_\star(\hphi_1-\hpsi)\}
  \end{equation*}
  so that $\hnu_1(y-\hpsi)\le\hnu_2(y-\hpsi)$ as claimed.
  
  For the converse, first suppose 
  $\hnu_\star(\hphi_2-\hphi_1)<\hbeta_1$.
  Pick $\hpsi=\hphi_1$. 
  Then $\hnu_1(y-\hpsi)=\hbeta_1$ but
  $\hnu_2(y-\hpsi)=\min\{\hbeta_2,\hnu_\star(\hphi_2-\hphi_1)\}<\hbeta_1$ 
  so $\hnu_1\not\le\hnu_2$.
  Next suppose $\hnu_\star(\hphi_2-\hphi_1)\ge\hbeta_1$
  but $\hbeta_1>\hbeta_2$.
  Pick $\hpsi=\hphi_2$.
  Then 
  $\hnu_1(y-\hpsi)=\min\{\hbeta_1,\hnu_\star(\hphi_2-\hphi_1)\}=\hbeta_1$
  but $\hnu_2(y-\hpsi)=\hbeta_2$ so again $\hnu_1\not\le\hnu_2$.
  This completes the proof.
\end{proof}
%
%
\subsection{Puiseux parameterization}
Next we parameterize $\hcV$.
\begin{definition}\label{D404}
  If $\hnu\in\hcV$, then define the \emph{Puiseux parameter} 
  of $\hnu$ by
  \begin{equation*}
    \hbeta(\hnu):=\sup\{\hnu(y-\hpsi)\ ;\ \hpsi\in\bk)\}.
  \end{equation*}
\end{definition}
\begin{proposition}\label{P429}
  The Puiseux parameter defines a parameterization 
  $\hbeta:\hcV\to[0,\infty]$ of the nonmetric 
  tree $\hcV$ rooted in the valuation $\hnu_\star$. Moreover:
  \begin{itemize}
  \item[(i)]
    if $\hnu$ is rational, then $\hbeta(\hnu)$ is rational;
  \item[(ii)]
    if $\hnu$ is irrational, then $\hbeta(\hnu)$ is irrational;
  \item[(iii)]
    if $\hnu$ is of point type, then $\hbeta(\hnu)=\infty$;
  \item[(iv)]
    if $\hnu$ is of special type, then $\hbeta(\hnu)\in(0,\infty]$.
  \end{itemize}
\end{proposition}
\begin{corollary}
  The Puiseux parameter also defines a 
  parameterization of the tree 
  $\hcVfin$ with values in $[0,\infty)$.
\end{corollary}
\begin{definition}\label{D410}
  Given $\hphi\in\bk_+$ and $\hbeta\ge0$ we denote by
  $\val[\hphi;\hbeta]$ the unique valuation 
  in the segment $[\hnu_\star,\hnu_\hphi]$ of Puiseux parameter
  $\hbeta$.
\end{definition}
\begin{proposition}\label{P404}
  If $\hnu=\val[\hphi;\hbeta]$ in the sense of 
  Definition~\ref{D410} then~\eqref{val-puis} holds. 
  Thus the notation $\hnu=\val[\hphi;\hbeta]$ is compatible 
  with that of Definition~\ref{def-puis}.
\end{proposition}
\begin{proposition}\label{P409}
  If $\hnu\in\hcV$ and $\hphi\in\bk_+$, then
  $\hnu_\hphi\ge\hnu$ iff 
  $\hnu(y-\hphi)=\hbeta(\hnu)$.
\end{proposition}
\begin{definition}\label{D411}
  The \emph{Puiseux metric} is the metric on the tree
  $\hcVfin$ induced by the 
  Puiseux parameterization, \ie the unique tree metric such
  that the distance between $\hnu_\star$ and $\hnu$ equals
  $\hbeta(\hnu)$.
\end{definition}
\begin{proof}[Proof of Proposition~\ref{P429}]
  Pick $\hnu\in\hcV$ and write 
  $\hnu=\val[\hphi;\hbeta]$ as in Definition~\ref{def-puis}.
  It follows from~\eqref{val-puis} that $\hbeta(\hnu)=\hbeta$.
  The proof of Proposition~\ref{P428} then immediately
  gives that the Puiseux parameter defines a 
  parameterization of $\hcV$.
  Assertions~(i)-(iv) are direct consequences of
  Definition~\ref{def-puis}.
\end{proof}
\begin{proof}[Proof of Proposition~\ref{P404}]
  Write $\hnu=\val[\phi';\hbeta']$ in the sense of
  Definition~\ref{def-puis}. 
  It follows from~\eqref{val-puis} that $\hbeta'$ 
  is the Puiseux parameter of $\hnu$, hence 
  $\hbeta'=\hbeta$.
  Moreover, we have $\hnu_\hphi=\val[\hphi;\infty]$ in the
  sense of Definition~\ref{def-puis} so since 
  $\hnu_\hphi\ge\hnu$, Lemma~\ref{L424} gives
  $\nu_\star(\hphi-\hphi')\ge\hbeta$.
  It is then easy to see that
  \begin{equation*}
    \min\{\hbeta,\nu_\star(\hpsi-\hphi)\}
    =\min\{\hbeta,\nu_\star(\hpsi-\hphi')\}
    =\hnu(y-\hpsi)
  \end{equation*}
  for any $\hpsi\in\hk$, which completes the proof.
\end{proof}
\begin{proof}[Proof of Proposition~\ref{P409}]
  The proof is similar to that of Proposition~\ref{P404}.
  Write $\hbeta=\hbeta(\hnu)$. 
  Then $\hnu=\val[\hphi';\hbeta]$ in the notation
  of Definition~\ref{def-puis} for some $\hphi'\in\bk_+$.
  Lemma~\ref{L424} gives $\hnu_\hphi\ge\hnu$ 
  iff $\nu_\star(\hphi-\hphi')\ge\hbeta$.
  Since $\hnu(y-\hphi)=\min\{\hbeta,\nu_\star(\hphi-\hphi')\}$,
  this is equivalent to $\hnu(y-\hphi)=\hbeta$.
\end{proof}
%
%
\subsection{Multiplicities}
Consider a Puiseux series $\hphi\in\bk$ written in the 
form~\eqref{puis-exp}. 
Write $\hbeta_j=r_j/s_j$ with $\gcd(r_j,s_j)=1$.
We define the \emph{multiplicity} of $\hphi$ by
$m(\hphi)=\lcm(s_j)$. 
Thus $m(\hphi)\in\overline{\N}$
and $m(\hphi)<\infty$ iff $\hphi\in\hk$.

As in the case of $\cVqm$ 
we can extend the notion of multiplicities from elements
in $\bk$ to valuations in $\hcV$ of finite type, using the tree
structure.
\begin{definition}
  If $\hnu\in\hcV$ is of point type or finite type,
  then the \emph{multiplicity} of $\hnu$ is 
  $m(\hnu)=\min\{m(\hphi)\ ;\ \hphi\in\bk_+, \hnu_\hphi\ge\hnu\}$.
  If $\hnu$ is of special type then we set $m(\hnu)=\infty$.
\end{definition}  
\begin{proposition}
  If $\hmu\le\hnu$ then $m(\hmu)$ divides $m(\hnu)$. 
  Further:
  \begin{itemize}
  \item[(i)]
    $m(\hnu)=\infty$ iff $\hnu$ is of special type or
    of point type associated to an element not in $\hk$;
  \item[(ii)]
    $m(\hnu)=1$ iff there exists a change of variables 
    $z=y-\hphi$ with $\hphi\in k_+$, such that
    have $\hnu(z-\hpsi)=\min\{\hbeta(\hnu),\nu_\star(\hpsi)\}$
    for every $\hpsi\in\bk$.
  \end{itemize}    
\end{proposition}
\begin{proof}
  We first claim that if $\hnu\in\hcV$ of finite type or point type 
  and $\hnu=\val[\hphi;\hbeta]$ as in Definition~\ref{def-puis},
  then $m(\hnu)=m(\hphi)$. 
  This is trivial if $\hnu$ is of point type.
  If $\hnu$ is of finite type and $\hnu_\psi>\hnu$, then 
  $\nu_\star(\hpsi-\hphi)\ge\hbeta$. 
  This implies that either $\hpsi=\hphi$, in which case 
  $m(\hpsi)=m(\hphi)$, or $\hpsi=\hphi$ plus higher order terms, 
  in which case $m(\hpsi)$ is a multiple of $m(\hphi)$.
  
  This claim together with Lemma~\ref{L424} implies that 
  $\hmu\le\hnu$ implies that $m(\hmu)$ divides $m(\hnu)$.
  It also immediately gives~(i). 

  Finally we prove~(ii).
  If the condition in~(ii) holds then
  $\hnu\le\hnu_\hphi$ by Proposition~\ref{P409}. 
  Since $\hphi\in k$, $m(\hphi)=1$ and so $m(\hnu)=1$.
  Conversely, if $m(\hnu)=1$, then we can find 
  $\phi\in\hk_+$ with $m(\hnu)$ such that 
  $\hnu\le\hnu_\hphi$. 
  Thus $\hphi\in k_+$ and $\hnu=\val[\hphi;\hbeta(\hnu)]$ so
  we are done by Proposition~\ref{P404}. 
\end{proof}
%
%
%
%
\section{Galois action}\label{Sgalois}
The Galois group $G:=\Gal(\hk/k)$ acts on $\hk$ and by duality on
valuations in $\hcV$.  Here we show that this action respects the tree
structure so that the orbit space $\hcV/G$ is in itself a tree. Later
we shall show that $\hcV/G$ is in fact isomorphic to the relative
valuative tree $\cV_x$.
%
%
\subsection{The Galois group}
As the field $\hk$ can be viewed as the injective limit of
field extensions, the Galois group $\Gal(\hk/k)$
is naturally a projective limit.
Let us be more precise.

For any $m\ge1$, set $k_m=\C((x^{1/m}))$. 
The map $x^{1/m}\mapsto(x^{1/mn})^n\in k_{mn}$
defines a field extension $k_m\to k_{mn}$ and
the fields $k_m$ form an injective system with
injective limit $\hk$.

The Galois group $G_m=\Gal(k_m/k)$ is isomorphic to the set of complex
numbers $\omega$ with $\omega^m=1$.  The action of $\omega\in G_m$ on
$k_m=\C((x^{1/m}))$ is denoted by $\omega^*$ (or $\omega_*$ which is
the same as $G$ is abelian) and satisfies $\omega^*(x^{1/m})=\omega
x^{1/m}$.  With this identification, there is a group homomorphism
from $G_{mn}$ to $G_m$ given by $\omega\mapsto\omega^n$.  The groups
$G_m$ form a projective system whose projective limit is the Galois
group $G=\Gal(\hk/k)$.

If $\omega\in G$, then we write $\omega_m$ for its image in $G_m$.
The action of $\omega\in G$ on a monomial $x^\beta\in\hk$ is then
given as follows. Write $\beta=p/q$ with $\gcd(p,q)=1$ and pick $m$
such that $q$ divides $m$, say $m=qr$.  Then
$\omega^*(x^\beta)=\omega_m^{pr}x^\beta$. This does not depend on the
choice of $m$.  It follows that $G$ also acts naturally on the
completion $\bk$ of $\hk$.

It is clear that $\nu_\star(\omega^*\hphi)=\nu_\star(\hphi)$
and $m(\omega^*\hphi)=m(\hphi)$ for every $\hphi\in\bk$. 
Moreover, if $m(\hphi)=1$,
\ie if $\hphi\in k$, then $\omega^*\hphi=\hphi$.

We extend the action of the Galois group to $\bk[y]$
by declaring $\omega^*y=y$ for every $\omega\in G$.
It follows that the restriction of $\omega^*$ to
$k[[y]]\supset R$ is the identity for every $\omega$.
%
%
\subsection{Action on $\hcV$}
By duality, the Galois group $G$ acts on valuations in $\hcV$. 
If $\omega\in G$, then the action of $\omega$ is 
denoted by $\omega_*$ and is given by
$(\omega_*\hnu)(\hpsi):=\hnu(\omega^*\hpsi)$.
We claim that $\omega_*$ maps $\hcV$ into itself.
Indeed, if $\hnu(y)>0$, then $(\omega_*\hnu)(y)=\hnu(y)>0$,
and if $\hnu$ extends $\nu_\star$, then so does
$\omega_*\hnu$ since 
$\nu_\star(\omega^*\hphi)=\nu_\star(\hphi)$ for every
$\hphi\in\bk$.
This calculation in conjunction with~\eqref{e425} also
shows that the action is well-behaved on valuations of point type:
$\omega_*\hnu_\hphi=\hnu_{\eta^*\hphi}$ for every $\hphi\in\bk_+$,
where $\eta=\omega^{-1}\in G$.
\begin{proposition}\label{P406}
  If $\omega\in G$, 
  then $\omega_*:(\hcV,\hbeta)\to(\hcV,\hbeta)$ is an 
  isomorphism of parameterized trees and preserves
  multiplicity.
\end{proposition}
\begin{proof}
  We have to show that $\omega_*$ is a bijection that 
  preserves the partial ordering, 
  the Puiseux parameter and the multiplicity.

  That $\omega_*$ is a bijection is obvious since
  $\omega^{-1}_*\omega_*=\omega_*\omega^{-1}_*=\id$.
  That it preserves the partial ordering, 
  and the Puiseux parameter are simple consequence
  of the equation 
  \begin{equation}\label{e423}
    (\omega_*\hnu)(y-\hpsi)=\hnu(y-\omega^*\hpsi)
    \quad\text{for $\hnu\in\hcV$ and $\hpsi\in\bk$}.
  \end{equation}
  Indeed,~\eqref{e424} and~\eqref{e423} immediately give 
  that $\hnu\le\hmu$ iff $\omega_*\hnu\le\omega_*\hmu$.
  Similarly, Definition~\ref{D404} and~\eqref{e423}
  yield $\hbeta(\omega_*\hnu)=\hbeta(\hnu)$ for
  every $\hnu\in\hcV$. 

  Finally let us show that $m(\omega_*\hnu)=m(\hnu)$ for 
  $\hnu\in\hcV$. 
  Consider $\hphi\in\hk$ such that $\hnu_\hphi>\hnu$ and 
  $m(\hphi)=m(\hnu)$. Write $\hpsi=(\omega^{-1})^*\hphi$.
  Then $m(\hpsi)=m(\hphi)=m(\hnu)$ and
  $\hnu_\hpsi=\omega_*\hnu_\hphi>\omega_*\hnu$. 
  Hence $m(\omega_*\hnu)\le m(\hnu)$. The reverse
  inequality follows by applying the same argument to
  $\omega^{-1}$ rather than $\omega$.
\end{proof}
%
%
\subsection{The orbit tree}\label{S17}
Proposition~\ref{P406} implies that the orbit space
$\hcV/G$ has a natural tree structure. Let us be more precise.
Declare two elements $\hnu,\hmu$ to be equivalent if
there exists $\omega\in G$ with $\hmu=\omega_*\hnu$.
Let $[\hnu]$ be denote the equivalence class containing $\hnu$
and let $\hcV/G$ be the set of equivalence classes.
We define a partial ordering on $\hcV/G$ by
$[\hnu]\le[\hmu]$ iff $\hnu\le\hmu$ for some representatives
$\hnu$, $\hmu$. Proposition~\ref{P406} implies that 
$\hcV/G$ is a complete, rooted nonmetric tree under this partial 
ordering and that $\hbeta$ defines a well-defined parameterization.
Moreover, the multiplicity $m$ is also well-defined.
%
%
%
%
\section{A tale of two trees}
We now wish to relate the tree $\hcV$ to the valuative tree $\cV_x$.
The elements of $\hcV$ are valuations on the ring
$\bk[[y]]\supset\C[[x]][[y]]=R$. Any element in $\hcV$ restricts to
$\nu_\star$ on $\C[[x]]$ so valuations in $\hcV$ restrict to
valuations on $R$ satisfying $\nu(x)=1$, and we end up with a mapping
$\Phi:\hcV\to\cV_x$.  Our first goal is to show that $\Phi$ respects
the two tree structures in the strongest possible sense.

Notice that $\Phi$ is not injective: the valuations 
$\val[y\pm x^{3/2};5/2]$ in $\hcV$
both map to the valuation $\nu_{y^2-x^3,2}$ in $\cV_x$.
Our second goal is to understand this lack of injectivity.
In fact, we shall see that $\Phi$ factors through the action
of the Galois group $G=\Gal(\hk/k)$, resulting in an isomorphism
between the orbit tree $\hcV/G$ defined in Section~\ref{S17}
and the relative valuative tree $\cV_x$.

%
%
\subsection{Minimal polynomials}
The main tool in the study of the restriction map
from $\hcV$ to $\cV_x$ is the relation between 
Puiseux series in $\hk$ and irreducible elements in $\fm$.
Recall that if $\phi\in\hk_+$, then the \emph{minimal polynomial}
of $\hat{\phi}$ over $k$ is given by 
\begin{equation}\label{e421}
  \phi(x,y)=\prod_{j=0}^{m-1}(y-\hphi_j)
\end{equation}
where $m=m(\hphi)$ is the multiplicity of $\hphi$, where
$\omega\in G$ is an element whose image in $G_m$ generates
$G_m$, and where $\hphi_j=(\omega^j)^*(\hphi)$.
Thus $\phi\in\fm$ is irreducible in $R$ and if
$\psi\in\fm\subset\bk[[y]]$ is irreducible in $R$ and 
$y-\hphi$ divides $\psi$ in $\bk[[y]]$, then 
$\psi=\phi\xi$, where $\xi$ is a unit in $R$.
Conversely, if $\psi\in\fm$ is irreducible in $R$, 
then there exists a unit $\xi\in R$ and $\hphi\in\hk_+$ 
such that $\psi=\phi\xi$, where $\phi$ is the minimal polynomial 
of $\hphi$.

Finally, two elements $\hphi,\hpsi\in\hk_+$ have the same minimal
polynomial iff there exists $\omega\in G$ such that
$\hpsi=\omega^*\hphi$.
%
%
\subsection{The morphism}
Recall that we are working with the relative tree $\cV_x$. 
Thus the elements of $\cV_x$ are normalized by $\nu(x)=1$,
and the root of $\cV_x$ is the valuation $\div_x$
defined by  $\phi\mapsto\max\{m\ ;\ x^m\mid\phi\}$.
Also recall that $\cV_x$ is equipped with the relative partial ordering 
$\le_x$, the parameterization $A_x$ by relative thinness,
and the relative multiplicity function $m_x$.
(The relative skewness will be of lesser importance here and
can anyway be recovered from the relative thinness and
multiplicity.)

By Section~\ref{S12}, $\hcV$ has a tree structure, comprised
of the partial ordering $\le$, 
the Puiseux parameter $\hbeta$ and the multiplicity function $m$.
Let us denote by $1+\hbeta$ the parameterization of 
$\hcV$ given by $(1+\hbeta)(\hnu)=1+\hbeta(\hnu)$ for all 
$\hnu\in\hcV$.

We define a mapping
$\Phi:\hcV\to\cV_x$ letting $\Phi(\hnu)$ be the
restriction of $\hnu$ from $\bk[[y]]$ to $R\subset k[[y]]$.
The fact that $\hnu$ restricts to $\hnu_\star$ on $\bk$
implies that $\Phi(\hnu)(x)=1$.
In particular $\Phi(\hnu_\star)=\div_x$ 
with the normalization convention above.
\begin{theorem}\label{T602}
  The mapping
  \begin{equation*}
    \Phi:(\hcV,1+\hbeta)\to(\cV_x,A_x)
  \end{equation*}
  is a surjective morphism of parameterized trees 
  and preserves multiplicity.
  Further, $\Phi$ factors through the projection
  $\hcV\to\hcV/G$, where $G=\Gal(\hk/k)$, and the induced mapping
  \begin{equation*}
    (\hcV/G,1+\hbeta)\to(\cV_x,A_x)
  \end{equation*}
  is a multiplicity-preserving isomorphism of parameterized trees.  
  
  Moreover, for $\hnu\in\hcV$ and $\nu=\Phi(\hnu)$ we have:
  \begin{itemize}
  \item[(i)]
    $\hnu$ is of finite type iff $\nu$ is quasimonomial;
  \item[(ii)]
    $\hnu$ is rational iff $\nu$ is divisorial;
  \item[(iii)]
    $\hnu$ is irrational iff $\nu$ is irrational;
  \item[(iv)]
    $\hnu$ is of point type and finite multiplicity iff
    $\nu$ is a curve valuation; in this case 
    $\hnu=\hnu_\hphi$ and $\nu=\nu_\phi$, where $\phi$ is
    the minimal polynomial of $\hphi$ over $k$;
  \item[(v)]
    $\hnu$ is of point type and infinite multiplicity iff
    $\nu$ is infinitely singular with infinite thinness;
  \item[(vi)]
    $\hnu$ is of special type iff
    $\nu$ is infinitely singular with finite thinness.
  \end{itemize}
\end{theorem}
Recall that the relative and nonrelative tree structures coincide on
the subtree $\{\nu\ge_x\nu_\fm\}=\{\nu\wedge\div_x=\nu_\fm\}$ of
$\cV_x$ by Proposition~\ref{Prel1}. Hence we have
\begin{corollary}
  $\Phi$ restricts to a surjective morphism 
  \begin{equation*}
    (\{\hnu\ge\val[0;1]\},1+\hbeta)\to(\{\nu\ge_x\nu_\fm\},A)
  \end{equation*}
  of parameterized trees and preserves multiplicity.
\end{corollary}
The fact that (relative) thinness can be computed in terms of Puiseux
series is very useful in practice and will play an important role in
Chapter~\ref{A3}.
\begin{figure}[ht]
  \begin{center}
    \includegraphics[width=\textwidth]{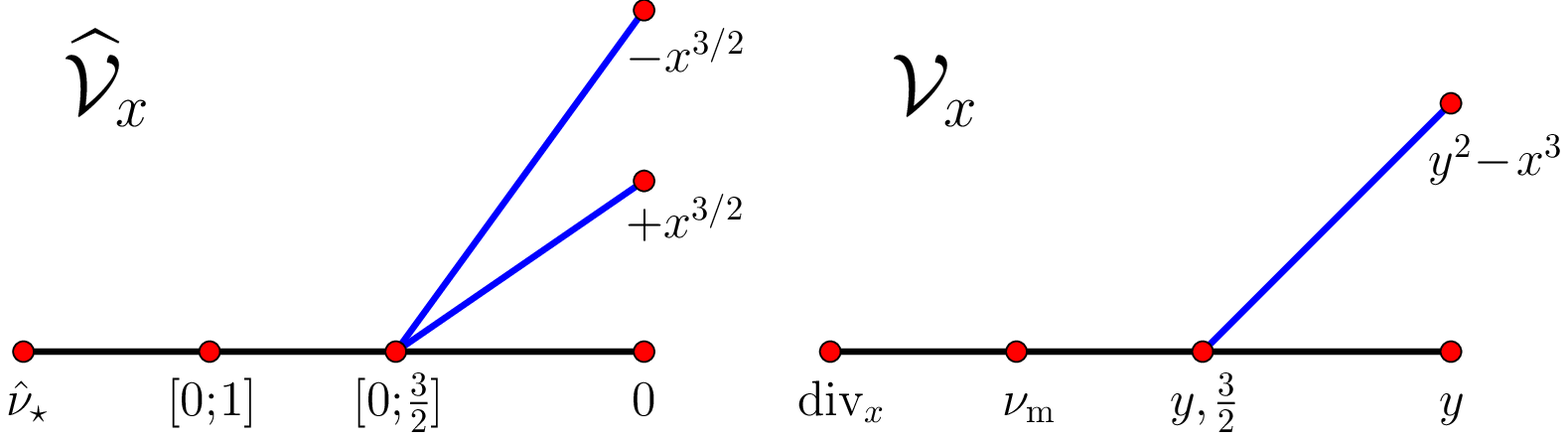}
  \end{center}
  \caption{The tree structures on the spaces $\hcV$ and $\cV_x$
    and their relation given by Theorem~\ref{T602}.
    The segment in $\hcV$ joining the root $\hnu_\star$ and
    the valuation $\hnu_0$ of point type is mapped onto
    the segment in $\cV_x$ joining the root $\div_x$ and
    the curve valuation $\nu_y$.
    The two segments in $\hcV$ joining $\hnu_\star$ and 
    the valuations $\nu_{\pm x^{3/2}}$ of point type are both 
    mapped onto
    the segment in $\cV_x$ joining the root $\div_x$ and
    the curve valuation $\nu_{y^2-x^3}$.}\label{F6}
\end{figure}
%
%
\subsection{Proof}
We now turn to the proof of Theorem~\ref{T602}. 
Almost the whole proof boils down to a study of the
factorization in~\eqref{e421}.

First notice that $\Phi$ is strictly increasing. 
Indeed, if $\hnu\le\hmu$ in $\hcV$, 
then $\hnu\le\hmu$ pointwise on $\bk[y]$ so 
by restriction $\Phi(\hnu)\le\Phi(\hmu)$ on $R$,
whence $\hnu\le_x\hmu$.
If moreover $\nu<\mu$, then there exists $\hphi\in\hk_+$
such that $\hnu(y-\hphi)<\hmu(y-\hphi)$. This implies
$\Phi(\hnu)(\phi)<\Phi(\hmu)(\phi)$,
where $\phi\in\hk[y]\subset\fm$ is the minimal polynomial
of $\hphi$ over $k$. Thus $\Phi(\nu)<_x\Phi(\mu)$.

Let us study the restriction of $\Phi$ to valuations of point type and
finite multiplicity.  Consider $\hphi\in\hk_+$ and let $\phi$ be the
minimal polynomial of $\hphi$ over $k$. It is given by~\eqref{e421}.
Then $\hnu_\hphi(y-\hphi)=\infty$, so $\Phi(\hnu_\hphi)(\phi)=\infty$,
which implies $\Phi(\hnu_\hphi)=\nu_\phi$ (normalized by
$\nu_\phi(x)=1)$). See the discussion following
Lemma~\ref{L601}. Essentially the same argument shows that if
$\hnu\in\hcV$ and $\Phi(\hnu)=\nu_\phi$ is a curve valuation, then
$\hnu=\hnu_\hphi$ as above. This proves~(iv).  Also notice
that~\eqref{e421} shows that the multiplicity of $\hphi$ equals the
relative multiplicity of $\phi$ (they are both equal to $m$).

The next step is to study the restriction of $\Phi$
to $\hcVfin$, the set of valuations of finite type. 
Generally speaking, 
a valuation of finite type can be understood through 
the valuations of point type (say of finite multiplicity)
that dominate it. 
Thus the previous paragraph, together with the fact that 
$\Phi$ is order-preserving,
allows us to understand the restriction of
$\Phi$ to $\hcVfin$. Of course we have to make this precise.

We claim that if $\hnu\in\hcVfin$, and $\phi\in\fm$ is irreducible,
then $\Phi(\hnu)<_x\nu_\phi$ iff there exists
$\hphi\in\hk$ such that $\phi$ is the minimal polynomial of
$\hphi$ over $k$ (up to a unit) and $\hnu<\hnu_\hphi$.
Indeed, write $\phi=\prod_j(y-\hphi_j)$
as in~\eqref{e421}, and write $\hnu_j=\hnu_{\hphi_j}$.
Suppose $\Phi(\hnu)<_x\nu_\phi$ but $\hnu\not<\hnu_j$ for all $j$.
Then we may pick $\hnu'<\hnu$ such that
$\hnu'\not<\hnu_j$ for all $j$.
Since $\Phi$ is order-preserving we get
$\Phi(\hnu')<_x\Phi(\hnu)<_x\nu_\phi$. 
But $\Phi(\hnu')(\phi)=\Phi(\hnu)(\phi)$ since
$\phi=\prod(y-\hphi_j)$, so this is a
contradiction. The other implication is easy.

Next we wish to show that the restriction of 
$\Phi$ to $\hcVfin$ preserves multiplicity
and transports the (adjusted) Puiseux parameterization 
to relative thinness. 

Let us first show that if $\hnu$ is
of finite type, then $m_x(\Phi(\hnu))=m(\hnu)$.
Consider $\hphi\in\hk_+$ with $\hnu_\hphi>\hnu$
and $m(\hphi)=m(\hnu)$.
Then $\nu_\phi>_x\Phi(\hnu)$, where $\phi$ is the minimal
polynomial of $\hphi$ over $k$. 
As noted above,~\eqref{e421} yields $m_x(\phi)=m(\hphi)$.
We infer $m_x(\Phi(\hnu))\le m(\hnu)$.
On the other hand, we saw above that if
$\phi\in\fm$ is irreducible and 
$\nu_\phi>\nu$, then $\phi$ is, up to a unit in $R$, 
the minimal polynomial over $k$ 
of some $\hphi\in\hk_+$ such that 
$\hnu_\hphi>\hnu$. This gives the reverse inequality.

To relate the Puiseux parameterization to relative thinness, 
we first consider $\hnu,\hnu'\in\hcV$ of finite type with
$\hnu'<\hnu$ such that the multiplicity is constant, 
say equal to $m$, on the segment $]\hnu',\hnu]$ in $\hcV$.
Pick $\hphi_0\in\hk_+$ with $\hnu_{\hphi_0}>\hnu$ and $m(\hphi_0)=m$. 
Write $\nu=\Phi(\hnu)$ and $\nu'=\Phi(\hnu')$. 
Let $\phi=\prod_j(y-\hphi_j)$ 
be the minimal polynomial of $\hphi_0$ over $k$.
Our assumptions imply
that $\nu_\star(\hphi_i-\hphi_j)<\hbeta(\hnu')$ for $i\ne j$.
Hence $\hnu(y-\hphi_j)=\hnu'(y-\hphi_j)$ for $j\ne0$
in view of~\eqref{val-puis}.
Moreover, $\hnu(y-\hphi_0)=\hbeta(\hnu)$ and
$\hnu'(y-\hphi_0)=\hbeta(\hnu')$.
This leads to 
\begin{equation*}
  A_x(\nu)-A_x(\nu')
  =\nu(\phi)-\nu'(\phi)
  =\sum_0^{m-1}\left(\hnu(y-\hphi_j)-\hnu'(y-\hphi_j)\right)
  =\hbeta(\hnu)-\hbeta(\hnu'),
\end{equation*} 
where we have used Lemma~\ref{L437}.
We also have 
$\Phi(\hnu_\star)=\div_x$ and $A_x(\div_x)=1=1+\hbeta(\hnu_\star)$.
Given $\hnu\in\hcVfin$ we can break up the segment 
$[\hnu_\star,\hnu]$ into finitely many pieces on each of
which the multiplicity is constant. 
By applying the equation above to each piece we 
get $A_x(\Phi(\hnu))=1+\hbeta(\hnu)$.

We have shown that $\Phi:\hcV\to\cV_x$ is order-preserving
and that $A_x\circ\Phi=1+\hbeta$ on $\hcVfin$. 
Since $\hcV$ is the completion of $\hcVfin$, this immediately
implies that $A_x\circ\Phi=1+\hbeta$ on all of $\hcV$, so
that $\Phi$ gives a morphism of parameterized trees.

We also know that $m_x(\Phi(\hnu))=m(\hnu)$ for $\hnu$ of
finite type and for $\hnu$ of point type of finite multiplicity.
If $\hnu$ is either of special type or of point type
with infinite multiplicity, then $m(\hnu)=\infty$. 
Thus $m(\hmu)\to\infty$ as $\hmu\to\hnu$ along the segment
$[\hnu_\star,\hnu[$. The valuations in this segment are
of finite type, so $m_x(\Phi(\hmu))\to\infty$. 
As $\Phi$ is order-preserving, $m_x(\Phi(\hmu))$ increases
to $m_x(\Phi(\hnu))$ so $m_x(\Phi(\hnu))=\infty$.
Thus $\Phi$ preserves multiplicity on all of $\hcV$.

Statements~(i)-(vi) now follow from the characterizations
of the elements in $\hcV$ and $\cV_x$ in terms of the tree
structure.

Finally, to show that $\Phi$ induces an isomorphism of
the orbit tree $\hcV/G$ onto $\cV_x$ it suffices to 
prove that if $\hnu,\hmu\in\hcV$, then 
$\Phi(\hnu)=\Phi(\hmu)$ iff there exists 
$\omega\in G$ such that $\hmu=\omega_*\hnu$.
One of these implications is trivial since
the restriction $\omega^*$ to $R$ is the identity for
any $\omega\in G$. 

For the other implication, suppose 
$\Phi(\hnu)=\Phi(\hmu)=:\nu$ for some $\hnu,\hmu\in\hcV$.
Then $\hnu$ and $\hmu$ share the same Puiseux parameter 
and multiplicity.
First assume the multiplicity is finite.
Then $m_x(\nu)<\infty$ so we can find
$\phi\in\fm$ irreducible with $\nu_\phi\ge\nu$.
As we showed above, this implies that there exist
$\hphi,\hpsi\in\hk$ both having $\phi$ as minimal
polynomial over $k$ and such that
$\hnu_\hphi\ge\hnu$ and $\hnu_\hpsi\ge\hmu$. 
Then $\eta^*\hpsi=\hphi$ for some $\eta\in G$,
hence $\omega_*\hnu_\hpsi=\hnu_\hphi$, where $\omega=\eta^{-1}\in G$.
By Proposition~\ref{P406} we infer that $\omega_*$ maps the
segment $[\hnu_\star,\hnu_\hpsi]$ onto $[\hnu_\star,\hnu_\hphi]$,
preserving the Puiseux parameter. 
Since $\hnu$ and $\hmu$ have the same Puiseux parameter, this
implies $\omega_*\hmu=\hnu$.

Now suppose $\Phi(\hnu)=\Phi(\hmu)=\nu$ and $m(\hnu)=m(\hmu)=\infty$.
Pick an increasing sequence $(\hnu_n)_1^\infty$ of rational
valuations in the segment $[\hnu_\star,\hnu[$ such that 
$m_n:=m(\hnu_n)$ is strictly increasing and
$\hnu_n\to\hnu$ as $n\to\infty$. 
Also pick $\hphi_n\in\hk$ such that $m(\hphi_n)=m_n$,
and $\hnu_{\hphi_n}\wedge\hnu=\hnu_n$.
Write $\nu_n=\Phi(\hnu_n)$ and let $\phi_n\in\fm$ be the minimal
polynomial of $\hphi_n$ over $k$.
Then $\nu_n$ increases to $\nu$, $m_x(\nu_n)=m_n$,
$m_x(\phi_n)=m_n$ and $\nu_{\phi_n}\wedge\nu=\nu_n$.
By what precedes there is a unique preimage $\hmu_n$ under $\Phi$
of $\nu_n$ in the segment $[\hnu_\star,\hmu[$ and we have
$m(\hmu_n)=m_n$. 
All preimages of $\nu_{\phi_n}$ under $\Phi$ are of the form
$\hnu_{\hpsi_n}$ where $\hpsi_n\in\hk$ and the 
minimal polynomial of $\hpsi_n$ over $k$ is $\phi_n$.
We automatically have $m(\hpsi_n)=m_n$ and 
we can pick $\hpsi_n$ such that $\hnu_{\hpsi_n}\wedge\hmu=\hmu_n$.

Since $\hpsi_n$ and $\hphi_n$ have the same minimal polynomial
over $k$ there exists $\eta_n\in G$ such that 
$\eta_n^*\hpsi_n=\hphi_n$. Notice while there are many such 
$\eta_n$, they all have the same image in $G_{m_n}$.
Thus there exists a unique $\eta\in G_{m_n}$ whose
image in $G_{m_n}$ agrees with that of $\eta_n$.
This $\eta$ satisfies $\eta^*\hpsi_n=\hphi_n$ for all $n$.
As a consequence $\omega_*\hnu_{\hpsi_n}=\hnu_{\phi_n}$
for all $n$, where $\omega=\eta^{-1}\in G$. 
Since $\omega_*$ preserves the Puiseux parameter, 
$\omega_*\hmu_n=\hnu_n$ for all $n$, so finally
$\omega_*\hmu=\hnu$ as desired.

This concludes the proof of Theorem~\ref{T602}.
%
%
%
%
\section{The Berkovich projective line}
We next indicate how the valuative spaces $\cV_x$ and $\hcV$ 
are naturally
embedded (as nonmetric trees) in the Berkovich projective lines 
over the local fields $k$ and $\bk$, respectively.

The Berkovich affine line $\bbA^1(k)$
\index{Berkovich!affine line} 
\index{$\bbA^1$ (Berkovich affine line)} 
is defined~\cite[p.19]{Ber} to be the set of valuations
$\nu:k[y]\to(-\infty,+\infty]$ extending $\nu_\star$. 
The line $\bbA^1(\bk)$ is defined analogously. 
We will write $\nu$ for elements
of $\bbA^1(k)$ and $\hnu$ for elements of $\bbA^1(\bk)$.

Both lines are endowed with the weak topology, defined in terms of
pointwise convergence.
By definition, $\hcV$ is a subset of $\bbA^1(\bk)$.
Berkovich defines~\cite[p.18]{Ber} the open
disk in $\bbA^1(\bk)$ with center $\hphi\in\bk$ and
radius $r>0$ to be 
$D(\hphi,r)=\{\hnu\ ;\ \hnu(y-\hphi)>\log r\}$.
It follows that 
\emph{$\hcV$ is the closure of the open disk in $\bbA^1(\bk)$ of 
  radius 1 centered at zero}. 
Disks in $\bbA^1(k)$ are defined as images of disks under 
the restriction map $\bbA^1(\bk)\to\bbA^1(k)$. 
It follows that the disk $D(0,1)$ in $\bbA^1(k)$ is the 
set of valuations $\nu:k[y]\to(-\infty,\infty]$ extending $\nu_\star$ and
satisfying $\nu(y)>0$. 
Any such valuation is an element of $\cV_x$, if we use the normalization
$\nu(x)=1$. Conversely, any element of $\cV_x$, save for $\div_x$, 
can be obtained in this way.
Hence \emph{$\cV_x$ is homeomorphic to the closure of 
  the open disk in $\bbA^1(\bk)$ of radius 1 centered at zero}. 

The Berkovich projective line $\bbP^1(k)$ is defined as
\index{Berkovich!projective line} 
\index{$\bbP^1$ (Berkovich projective line)} 
$\bbA^1(k)\cup\{\nu_\infty\}$ 
where $\nu_\infty:k[y]\to[-\infty,\infty]$ equals 
$\nu_\star$ on $k$ and $-\infty$ elsewhere.
The line decomposes as 
\begin{equation*}
  \bbP^1(k)=\{\hnu_\star\}\sqcup\bigsqcup_{\a\in\P^1(\C)}D_\alpha,
\end{equation*}
where $D_\alpha:=D(y-\alpha,1)$ for $\alpha\in\C$
and $D_\infty:=\{\nu\in\bbA^1(k)\ ;\ \nu(y)<0\}\cup\{\nu_\infty\}$
and $\hnu_\star\in\bbA^1(k)$ denotes the valuation
$\hnu_\star(\sum_ja_jy^j)=\min_j\nu_\star(a_j)$.

Any translation $y\mapsto y+\tau$, $\tau\in\C$,
induces an isomorphism of $\bbP^1(k)$ sending 
$D_\a$ homeomorphically to $D_{\a+\tau}$. 
Similarly, the inversion $y\mapsto y^{-1}$ induces a 
homeomorphism of $D_\infty$ onto $D_0$. 
The open disks $D_\a$ are hence all homeomorphic, 
and their closures are $\overline{D_\a}=D_\a\cup\{\hnu_\star\}$.

On any disk $\overline{D_\a}$ we can put the nonmetric tree structure
induced by the relative tree structure on $\cV_x$ 
(rooted at $\div_x$, see Section~\ref{sec-relative}).
For instance, on $\overline{D_0}\simeq\cV_x$, the tree structure 
is given by $\nu\le\mu$ iff $\nu(\phi)\le\mu(\phi)$ for all 
$\phi\in k[y]$. 
Patched together, these partial orderings
endow $\bbP^1(k)$ with a natural nonmetric tree structure rooted at 
$\nu_\star$ (see~\cite[Theorem~4.2.1]{Ber}).
As we shall see below (Theorem~\ref{wt}), the weak tree
topology on $\bbP^1(k)$ coincides with the weak topology.

This discussion goes through, essentially verbatim, 
with $k$ replaced by $\bk$.
%
%
%
%
\section{The Bruhat-Tits metric}\label{bruhat-tits}
The standard projective line $\P^1(\bk)$ embeds naturally 
in the Berkovich line $\bbP^1(\bk)$: $\hphi\in\bk$ 
corresponds to the valuation $\hnu_\hphi$ with 
$\hnu_\hphi(y-\hphi)=\infty$ and $\infty$ to $\hnu_\infty$. 
The points in $\P^1(\bk)$ are ends of the nonmetric tree $\bbP^1(\bk)$;
hence $\bbH:=\bbP^1(\bk)\setminus\P^1(\bk)$ has a nonmetric tree 
structure rooted at $\hnu_\star$. The group 
$\mathrm{PGL}_2(\bk)$ of M\"obius transformations acts 
on $\P^1(\bk)$ and this action extends to $\bbP^1(\bk)$ 
by $(M_*\hnu)(y-\hphi)=\hnu(My-\phi)$.

Berkovich noted that $\bbH$ is isomorphic (as a non-metric tree) to
the Bruhat-Tits building of $\mathrm{PGL}_2(\bk)$\index{Bruhat-Tits
building of $\mathrm{PGL}_2(\bk)$}.  We will not define this building
nor the isomorphism here.  Suffice it to say that the building is a
metric tree on which $\mathrm{PGL}_2(\bk)$ acts by \emph{isometries}.
This last condition in fact defines the tree metric on $\bbH$ up to a
constant.  To see this, consider the segment in $\bbH$ parameterized
by $\val[0;t]$, $0\le t<\infty$. Fix any rational $t_0>0$ and pick
$\hphi_0$ with $\hnu_\star(\hphi_0)=t_0$.  Then
$M\in\mathrm{PGL}_2(\bk)$ defined by $My=y\hphi_0$ gives
$M_*\val[0;t]=\val[0;t+t_0]$.  But the only translation invariant
metric on $\R$ is the Euclidean metric (up to a constant), so after
normalizing we get that $[0,\infty[\,\ni t\mapsto\val[0,t]\in\bbH$ is
an isometry.  Now if $\hphi\in\bk$, then $M(y):=y-\hphi$ induces an
isometry of the segment $[\hnu_\star,\hnu_\hphi[$ onto
$[\hnu_\star,\hnu_0[$.  Similarly, $M(y)=1/y$ gives an isometry of
$[\hnu_\star,\infty[$ onto $[\hnu_\star,\hnu_0[$.  Thus the metric on
$\bbH$ is uniquely defined.

In particular we see that if $\hphi\in\bk_+$, then
$[0,\infty[\,\ni t\mapsto\val[\hphi,t]\in\bbH$
is an isometry. 
Thus \emph{the Puiseux metric on $\hcVfin$ 
  is induced by the Bruhat-Tits metric on $\bbH\supset\hcV$}.
Passing to $\cV_x$ we conclude using Theorem~\ref{T602} than
\emph{the relative thin metric on $\cVqmx$ 
  is induced by the Bruhat-Tits metric}.
%
%
%
%
\section{Dictionary}
We end this chapter with a dictionary between valuations in $\cV_x$
\ie on $k[y]$ , their preimages in $\hcV$
(\ie on $\bk[y]$) under the restriction map (see Theorem~\ref{T602}),
and Berkovich's terminology~\cite[p.18]{Ber}.  See
Appendix~\ref{sec-clas} for a more extensive dictionary.

Note that except for $\nu=\div_x$, the terminology for a valuation
$\nu\in \cV_x$ and for its nonrelative counterpart
$\nu/\nu(\fm)\in\cV$ coincide. The valuation $\div_x$ in $\cV_x$ is
not centered at the origin, and is of divisorial type.  It corresponds
to $\nu_x$ which is a curve valuation in $\cV$. The preimage of
$\div_x$ in $\hcV$ is $\hnu_\star$ which is of rational type in our
terminology and of type 2 in Berkovich's terminology. With our
convention, $\div_x$ also fits into the table below.

\begin{table}[ht]
  \begin{center}
    \begin{tabular}{|c|c|c|c|c|c|}
      \hline
      \multicolumn{3}{|c|}{Valuations in $\cV_x$} & 
      \multicolumn{2}{|c|}{Valuations in $\hcV$} & 
      Berkovich\\
      \hline\hline

      & 
      \multicolumn{2}{|c|}{Divisorial} &
      & 
      Rational & 
      Type 2\\

      \cline{2-3}\cline{5-6}
      \raisebox{1.5ex}[0cm][0cm]{Quasim.} &
      \multicolumn{2}{|c|}{Irrational} &
      \raisebox{1.5ex}[0cm][0cm]{Finite type} & 
      Irrational & 
      Type 3\\ 
      \hline

      & 
      \multicolumn{2}{|c|}{Curve} &
      & 
      $m<\infty$ &
      \\
      \cline{2-3}\cline{5-5}
      \raisebox{2.5ex}[0cm][0cm]{Not quasi-}
      & 
      & 
      $A=\infty$ &
      \raisebox{1.5ex}[0cm][0cm]{Point type} & 
      $m=\infty$ & 
      \raisebox{1.5ex}[0cm][0cm]{Type 1}\\ 

      \cline{3-6}
      \raisebox{2.5ex}[0cm][0cm]{monomial} & 
      \raisebox{1.5ex}[0cm][0cm]{Inf sing} & 
      $A<\infty$ &
      \multicolumn{2}{|c|}{Special type} & 
      Type 4\\ 
      \hline
    \end{tabular}
  \end{center}
  \medskip
  \caption{Terminology}
  \label{table3}
\end{table}


%
%
%
%
%
%
\chapter{Topologies}\label{part4}
Our objective now is to introduce, analyze and compare 
different topologies on the valuative tree $\cV$.

There are at least three natural ways of defining topologies. 
First, we can exploit the tree structure of $\cV$. 
This leads to one topology defined using the nonmetric tree structure, 
and several others defined in terms of metric tree structures.

Second, recall that $\cV$ is by definition a collection of 
(normalized) functions from $R$ to $\Rbar$. 
Hence we can define topologies using pointwise or 
(suitably normalized) uniform convergence.

A third type of topologies can be defined on the set $\cV_K$ of all 
equivalence classes of centered Krull valuations, 
and then passing to $\cV$. 
As a Krull valuation takes values in an abstract ordered group, 
pointwise or uniform convergence a priori does not makes sense. 
However, we can always check whether the value is positive, 
zero or negative. 
Different topologies on the set $\{+,0,-\}$ then lead to  
different topologies on $\cV_K$. 
These can then be turned into topologies on the valuative tree.

The chapter is organized as follows.
We start in Section~\ref{sec-weak} by considering the \emph{weak topology}.
As we show, it can be defined equivalently either in terms of 
pointwise convergence or in terms of the nonmetric tree structure. 

Then in Sections~\ref{naturalmetric} and~\ref{S19}
we study the \emph{strong topology}.
More precisely, there is one strong topology on $\cV$ and one on
the subtree $\cVqm$ of quasimonomial valuations. They both arise
from the parameterization of $\cV$ by skewness. Alternatively,
they can be defined in terms of normalized uniform convergence.

As a slight variation, Section~\ref{S18} contains a brief discussion 
of the \emph{thin topology}, which also comes in two versions: 
one on $\cV$ and one on $\cVqm$. 
Both versions arise from the parameterization of $\cV$ by thinness.

After that we turn to the \emph{Zariski topology}. 
It is defined on the set $\cV_K$ of (equivalence classes of) 
centered Krull valuations. We analyze it in Section~\ref{sec-zar}. 
The main result is that if we start with $\cV_K$ endowed with 
the Zariski topology
and identify exceptional curve valuations with their 
associated divisorial valuations, then we recover the valuative tree
with the weak topology.

A refinement of the Zariski topology is the
\emph{Hausdorff-Zariski topology}, considered in Section~\ref{HZ}.
We show that it also has a natural interpretation as a weak
tree topology, but now given in terms of a $\Nbar$-tree
structure on $\cV_K$.

Finally we compare the different topologies on $\cV$ in
Section~\ref{topo-compare}.
%
%
%
%
\section{The weak topology}\label{sec-weak}
At this stage we have two candidates for the weak topology on the 
valuative tree $\cV$: the weak topology defined by pointwise convergence 
and the weak tree topology induced by the nonmetric tree structure on $\cV$ 
(see Section~\ref{S302}). 
In this section we show that these two topologies are in fact identical. 
We also study its main properties.
%
%
\subsection{The equivalence}
We defined the weak topology
\index{topology!weak (on $\cV$)} 
on $\cV$ by
describing converging sequences: $\nu_k\to\nu$ iff
$\nu_k(\phi)\to\nu(\phi)$ for every $\phi\in\fm$.  Equivalently, a
basis of open sets is given by $\{\nu\in\cV\ ;\ t<\nu(\phi)<t'\}$ over
$\phi\in\fm$ irreducible and $t'>t\ge1$.  The weak tree topology is
generated by $U(\vv)$, where $\vv$ runs over all tangent vectors
in $\cV$.
\begin{theorem}\label{wt}
  The weak topology on $\cV$ coincides with the weak tree topology.
\end{theorem}
\begin{proof}
  We use Proposition~\ref{P201} repeatedly.
  If $\vv\in T\nu$ is a tangent vector not represented by
  $\nu_\fm$ (true \eg if $\nu=\nu_\fm$), then
  $U(\vv)=\{\mu\ ;\ \mu(\phi)>\a\,m(\phi)\}$, where
  $\a=\a(\nu)$ and $\nu_\phi$ represents $\vv$. 
  If instead $\vv$ is represented by $\nu_\fm$, then
  $U(\vv)=\{\mu\in\cV\ ;\ \mu(\phi)<\a\,m(\phi)\}$,
  where $\nu_\phi\ge\nu$.
  In both cases $U(\vv)$ is open in the weak topology.

  Conversely, every nontrivial set of the form
  $\{\nu(\phi)>t\}$ or
  $\{\nu(\phi)<t\}$ with $t\ge1$ and 
  $\phi\in\fm$ irreducible is of the form $U(\vv)$. 
  This completes the proof.
  \end{proof}
%
%
\subsection{Properties}
We now investigate the weak topology further.
\begin{proposition}\label{Pweakcpt}
  The weak topology is compact but not metrizable.
\end{proposition}
\begin{proof}
  Compactness follows from Theorem~\ref{wt} and 
  the fact that any complete, parameterizable tree is weakly compact 
  (see Theorem~\ref{T401}).
  Alternatively, $\cV$ is naturally embedded as a 
  closed subspace of the compact product space $[0,\infty]^R$.

  To prove that $\cV$ is not metrizable 
  it suffices to show that $\nu_\fm\in\cV$ 
  has no countable basis of open neighborhoods. 
  Fix local coordinates $(x,y)$.
  For $\theta\in\C^*$, let 
  $\nu_\theta\=\nu_{y+\theta x,2}$ (see Definition~\ref{D401})
  and $U_\theta\=\{\nu\not\ge\nu_\theta\}$. 
  Any $U_\theta$ is an open neighborhood of $\nu_\fm$. 
  Suppose $\{V_k\}_{k\ge1}$ is a countable basis of open 
  neighborhoods of $\nu_\fm$. Then for any $\theta$ there exists
  $k=k(\theta)\ge0$ with $V_k\subset U_\theta$. 
  As $\C^*$ is uncountable, one can find $k$ and
  a sequence $\theta_i$ with $\theta_i\ne\theta_j$ for $i\ne j$ 
  such that $V_k\subset U_{\theta_i}$. 
  Now $\nu_{\theta_i}\to\nu_\fm$, so
  $\nu_{\theta_i}\in V_k$ for $i\gg1$.
  But $\nu_{\theta_i}\not\in U_{\theta_i}$. 
  This is a contradiction.
\end{proof}
\begin{proposition}\label{P110}
  The four subsets of $\cV$ consisting of
  divisorial, irrational, infinitely singular
  and curve valuations 
  are all weakly dense in $\cV$.
\end{proposition}
\begin{proof}
  Instead of proving all of this directly we will appeal to
  Proposition~\ref{P105} which asserts that the divisorial, irrational
  and infinitely singular valuations are all (individually) dense
  in $\cV$ in a stronger topology that the weak topology. Hence 
  it suffices to show that, say, any divisorial valuation can be 
  weakly approximated by curve valuations. 
  But if $\nu$ is divisorial,
  then by Proposition~\ref{P107} there exists
  a sequence $\phi_n$ of irreducible elements in $\fm$ 
  such that the associated curve valuations $\nu_{\phi_n}$ 
  represent distinct tangent vectors at $\nu$.
  Then $\nu_{\phi_n}\to\nu$ as $n\to\infty$. 
\end{proof}
\begin{remark}
  In Section~\ref{sec-relative} we considered the 
  relative valuative tree $\cV_x$ consisting of centered
  valuations $\nu$ on $R$ normalized by $\nu(x)=1$ 
  (and with the valuation $\div_x$ added).
  We could define two relative weak topologies on $\cV_x$: 
  one using the relative nonmetric tree structure
  and one using pointwise convergence.
  The proof of Theorem~\ref{wt} could easily be adjusted to show that
  these relative weak topologies are in fact the same. 
\end{remark}
%
%
%
%
\section{The strong topology on $\cV$}\label{naturalmetric}
As with the weak topology, we have two candidates for the
strong topology on $\cV$: one defined using uniform sequential
\index{topology!strong (on $\cV$)}
convergence (see~\eqref{e106} below) 
and one induced by the tree metric~\eqref{e105} on $\cV$.
We now show that these two coincide and analyze the properties of
the strong topology.
%
%
\subsection{Strong topology I}\label{S20}
A general metric tree admits two natural topologies, 
the weak topology as a nonmetric tree and the 
topology as a metric space.
We emphasize that the latter topology depends on the choice of metric.

The valuative tree carries several interesting metric tree structures.
For the moment we are mainly interested in the
tree metric $d$
\index{$d$ (metric on $\cV$)}
induced by skewness and defined in~\eqref{e105}:
\begin{equation}\label{e426}
  d(\mu,\nu)
  =\left(\a(\mu\wedge\nu)^{-1}-\a(\mu)^{-1}\right)
  +\left(\a(\mu\wedge\nu)^{-1}-\a(\nu)^{-1}\right).
\end{equation}
We shall refer to the induced topology as the 
\emph{strong tree topology} on $\cV$.
\index{topology!strong (on $\cV$)}

Before analyzing the strong topology in detail, let us note
the following general result.
\begin{proposition}\label{P430}
  The topology on a metric tree $(\cT,d)$ induced by the metric
  is at least as strong as the weak topology.
\end{proposition}
\begin{proof}
  Let $\sigma_n$ be a sequence of points in the metric 
  tree $(\cT,d)$ that converges to $\sigma_\infty\in\cT$ in the metric. 
  Use $\sigma_\infty$ as a root.
  If $\sigma_n\not\to\sigma_\infty$ weakly, then there would
  exist $\sigma>\sigma_\infty$ such that $\sigma_n\ge\sigma$ for 
  infinitely many $n$. 
  Thus $d(\sigma_n,\sigma_\infty)\ge d(\sigma,\sigma_\infty)>0$
  for these $n$, a contradiction.
  \end{proof}                                
\begin{remark}
  If $\cT$ is as in Example~\ref{E1} with $X$ infinite and 
  $g(x)\equiv1$, then the two topologies are not
  equivalent: if $(x_n)_{n\ge0}$ are distinct elements in $X$ then
  $(x_n,1)$ converges weakly to $(x,0)$, 
  but $d((x_n,1),(x,0))=1$.
\end{remark}
%
%
\subsection{Strong topology II}\label{S21}
Another strong topology on $\cV$ can be defined in a 
quite general setting (\ie for other
rings than $R$) in terms of the metric
\index{$d^\mathrm{str}$ (metric on $\cV$)}
\begin{equation}\label{e106}
  d^\mathrm{str}_\cV(\nu_1,\nu_2)
  =\sup_{\phi\in\fm\ \mathrm{irreducible}}
  \left|
    \frac{m(\phi)}{\nu_1(\phi)}-\frac{m(\phi)}{\nu_2(\phi)}
  \right|
\end{equation}
This topology is stronger than the weak topology and
any formal invertible mapping $f:(\C^2,0)\to(\C^2,0)$
induces an isometry $f_*:\cV\to\cV$ for $d^\mathrm{str}_\cV$.
%
%
\subsection{The equivalence}
We now show that the metric $d^\mathrm{str}_\cV$ is compatible with 
the tree metric $d$ given by~\eqref{e426}. As a consequence, the
two strong topologies on $\cV$ defined in Sections~\ref{S20} and~\ref{S21}
are really the same.
\begin{theorem}\label{T501}
  The strong topology on $\cV$ is identical to the strong tree topology.
  More precisely, if $d$ and $d^\mathrm{str}$ are the metrics on
  $\cV$ given by~\eqref{e426} and~\eqref{e106}, respectively, 
  then, for $\nu_1,\nu_2\in\cV$:
  \begin{equation}\label{e108}
    d^\mathrm{str}_\cV(\nu_1,\nu_2)
    \le d(\nu_1,\nu_2)
    \le 2d^\mathrm{str}_\cV(\nu_1,\nu_2)
  \end{equation}
\end{theorem}
\begin{proof}
  Let us consider $\phi\in\fm$ irreducible.  Set
  $\nu=\nu_1\wedge\nu_2$. We first show that
  \begin{equation}\label{e101}
    \left|\nu_1^{-1}(\phi)-\nu_2^{-1}(\phi)\right|
    =\max_{i=1,2}\left\{\nu^{-1}(\phi)-\nu_i^{-1}(\phi)\right\}.
  \end{equation} To see this, notice that Proposition~\ref{P201}
  implies $\mu(\phi)=(\mu\wedge\nu_\phi)(\phi)$ for any $\mu\in\cV$. 
  Thus we may replace $\nu_i$ by $\nu_i\wedge\nu_\phi$, so that
  $\nu\le\nu_i\le\nu_\phi$, $i=1,2$. But this means that
  $\nu=\nu_1\le\nu_2\le\nu_\phi$ or $\nu=\nu_2\le\nu_1\le\nu_\phi$ and
  then~\eqref{e101} is immediate.

  Multiplying~\eqref{e101} by $m(\phi)$ and
  taking the supremum over $\phi\in\fm$ we get
  $d^\mathrm{str}_\cV(\nu_1,\nu_2)=\max_i d(\nu_i,\nu)$,
  which implies~\eqref{e108}.
\end{proof}
%
%
\subsection{Properties}
We now investigate the strong topology further.
\begin{proposition}\label{P108}
  The strong topology on $\cV$ is strictly stronger than the 
  weak topology. It is not locally compact.
\end{proposition}
\begin{proof}
  The first assertion follows from Proposition~\ref{Pweakcpt}, or from
 the last assertion as $\cV$ is weakly compact Consider $\nu\in\cV$
 divisorial and pick $\phi_n\in\fm$ irreducible with
 $\nu_{\phi_n}>\nu$, such that $\nu_{\phi_n}$ represent distinct
 tangent vectors at $\nu$.  For fixed $\e>0$, set
 $\nu_n=\nu_{\phi_n,\a+\e}$, where $\a$ is the skewness of $\nu$.  Any
 $\nu_n$ is at distance $\e/\a(\a+\e)$ from $\nu$.  Further,
 $\nu_n\to\nu$ weakly, so $\nu_n$ has no strong accumulation point.
 If $\nu$ had a compact strong neighborhood, it would contain a ball
 of positive radius, say bigger than $\e$.  Then $\nu_n$ would have a
 strongly convergent subsequence.  This is impossible.
\end{proof}
\begin{proposition}\label{P105}
  The three subsets of $\cV$ consisting of
  divisorial, irrational and infinitely singular valuations 
  are all strongly dense in $\cV$.
  
  The strong closure of the set $\cC$ of curve valuations 
  is the set of valuations of infinite skewness.
  In particular, $\overline{\cC}\setminus\cC$ contains only 
  infinitely singular valuations.
\end{proposition}
\begin{proof}
  For the first assertion, we first prove that any valuation 
  $\nu\in\cV$ can be strongly approximated by divisorial and 
  irrational valuations. If $\nu$ is not infinitely singular,
  then $\nu=\nu_{\phi,t}$ for some irreducible $\phi\in\fm$
  and $t\in[1,\infty]$. Then $\mu:=\nu_{\phi,s}$ converges
  strongly to $\nu$ as $s\to t$ and $\mu$ is 
  divisorial (irrational) if $s$ is rational (irrational).
  If $\nu$ is infinitely singular, 
  then $\nu=\lim\nu_n$ for an increasing sequence $\nu_n$, 
  where $\nu_n$ can be chosen to
  be all divisorial or all irrational.

  The fact that any divisorial valuation can be approximated by an
  infinitely singular valuation follows from Lemma~\ref{L777} below and
  $\a(\mu) -\a(\nu) < A(\mu) - A(\nu)$ when $\mu> \nu$.
  This completes the proof of the first assertion. 
 
  For the second assertion, first notice that skewness defines a
  strongly continuous function $\a:\cV\to[1,\infty]$,
  so since every curve valuation has infinite skewness, so does 
  every valuation in $\overline{\cC}$. 
  Conversely, suppose $\nu$ is an infinitely singular valuation
  with infinite skewness and let $(\nu_n)_0^\infty$ be its
  approximating sequence. Then $\a(\nu_n)\to\infty$ as $n\to\infty$.
  For each $n$, let $\mu_n$ be a curve
  valuation with $\mu_n>\nu_n$ and 
  $m(\mu_n)=m(\nu_n)<b(\nu_n)$. 
  We get
  \begin{equation*}
    d(\nu_n,\nu)
    =\left(\frac1{\a(\nu\wedge\mu_n)}-\frac1{\a(\nu)}\right)
    +\left(\frac1{\a(\nu\wedge\mu_n)}-\frac1{\a(\mu_n)}\right)
    =\frac2{\a(\nu_n)}
    \to0,
  \end{equation*}
  so $\mu_n\to\nu$ as $n\to\infty$. This completes the proof. 
\end{proof}
\begin{remark}\label{R406}
  It is also possible to define a 
  \emph{relative strong topology} on $\cV_x$,
  based on the parameterization by relative skewness.
  Suitable versions of all the results above continue
  to hold in the relative setting. 
  In Proposition~\ref{P105} we should then use the convention that
  $\div_x$ is a divisorial valuation, and in particular not
  an element of $\cC$.
  Also notice that a valuation in $\cV_x\setminus\{\div_x\}$ 
  has infinite relative skewness iff the corresponding
  valuation in $\cV$ has infinite (nonrelative) skewness.
\end{remark}
%
%
%
%
\section{The strong topology on $\cVqm$}\label{S19}
In many instances it is natural and convenient to work on the
tree $\cVqm$ consisting of quasimonomial valuations.
As described in Section~\ref{normal} we can use the parameterization
by skewness to define a metric on $\cVqm$:
for $\nu,\mu\in\cVqm$ we have:
\index{$d_\mathrm{qm}$ (metric on $\cVqm$)}
\begin{equation}\label{e427}
  d_\mathrm{qm}(\mu,\nu)
  =(\a(\mu)-\a(\mu\wedge\nu))+(\a(\nu)-\a(\mu\wedge\nu)).
\end{equation}
We refer to the resulting topology as the 
\emph{strong tree topology} on $\cVqm$.

On the other hand, we could also consider the following strong metric:
\begin{equation*}
  d^\mathrm{str}_\mathrm{qm}(\mu,\nu)
  =\sup_{\phi\in\fm\ \mathrm{irreducible}}
  \left|
    \frac{\nu_1(\phi)}{m(\phi)}-\frac{\nu_2(\phi)}{m(\phi)}
  \right|
\end{equation*}
(It is in fact not hard to see that the supremum can be taken over 
all $\phi\in\fm$, not just $\phi$ irreducible.)
Following the proof of Theorem~\ref{T501}, we infer
\begin{theorem}
  The strong topology coincides with the 
  strong tree topology on $\cVqm$.
  More precisely, for $\nu_1,\nu_2\in\cVqm$ we have
  \begin{equation}\label{e109}
    d^\mathrm{str}_\mathrm{qm}(\nu_1,\nu_2)
    \le d_\mathrm{qm}(\nu_1,\nu_2)
    \le 2d^\mathrm{str}_\mathrm{qm}(\nu_1,\nu_2).
  \end{equation}
\end{theorem}  
\begin{proposition}\label{metrictor}
  The completion of $(\cVqm,d_\mathrm{qm})$ is a tree naturally
  isomorphic to the union of $\cVqm$ and all infinitely singular
  valuations with finite skewness. 
\end{proposition}
\begin{proof}
  The subset $\cV'$ of $\cV$ consisting of valuations with finite
  skewness is a strongly open subtree that contains all quasimonomial
  valuations. The metric $d_\mathrm{qm}$ extends naturally to $\cV'$ as
  a metric $d_{\cV'}$.  Let us show that $(\cV',d_{\cV'})$ is complete.
  Any $d_{\cV'}$-Cauchy sequence $\nu_n$ in $\cV'$ is $d_{\cV}$-Cauchy,
  hence $d_{\cV}(\nu_n,\nu)\to0$ for some $\nu\in\cV$ by completeness of
  $d_{\cV}$. But it is easy to see that $\nu_n$ must have uniformly
  bounded skewness, so $\a(\nu)<\infty$ and then
  $d_{\cV'}(\nu_n,\nu)\to0$ as well.
\end{proof}
\begin{remark}\label{R407}
  We can also equip the tree $\cVqmx$ 
  with a relative strong topology 
  defined in terms of the parameterization
  by relative skewness. See Section~\ref{sec-relative}
  and compare Remark~\ref{R406}.
\end{remark}
%
%
%
%
\section{Thin topologies}\label{S18}
Instead of skewness we can use thinness to parameterize the 
valuative tree.
As noted in Section~\ref{S10} we can use this parameterization
to define metrics $D$
\index{$D$ (metric on $\cV$)} 
on $\cV$ and $D_\mathrm{qm}$
\index{$D_\mathrm{qm}$ (metric on $\cVqm$)}:
\begin{align}
  D(\mu,\nu)
  &=\left(A(\mu\wedge\nu)^{-1}-A(\mu)^{-1}\right)
  +\left(A(\mu\wedge\nu)^{-1}-A(\nu)^{-1}\right)\label{e428}\\
  D_\mathrm{qm}(\mu,\nu)
  &=(A(\mu)-A(\mu\wedge\nu))+(A(\nu)-A(\mu\wedge\nu))\label{e429}.
\end{align}
We refer to $D$ and $D_\mathrm{qm}$ as 
the \emph{thin metric}\index{thin!metric} on $\cV$ and $\cVqm$,
respectively, and to the induced topologies as the
\emph{thin topologies}.
\index{topology!thin (on $\cV$)}

The thin metrics share a lot of similarities with the previously defined
metrics $d$ and $d_\mathrm{qm}$ defined in terms of skewness.
We summarize them in the following propositions.
\begin{proposition}\label{thin1}
  The metric space $(\cV,D)$ is complete, and not locally compact.  The
  three subsets consisting of divisorial, irrational and infinitely
  singular valuations are dense in $(\cV, D)$. 
  The closure of the set of curve valuations is the set 
  of valuations of infinite thinness.
\end{proposition}
\begin{proposition}\label{thin2}
  The completion of $(\cVqm,D_\mathrm{qm})$ is naturally
  isomorphic to the union of $\cVqm$ and all 
  infinitely singular valuations with finite thinness. 
\end{proposition}
\begin{proof}[Proof of Propositions~\ref{thin1} and~\ref{thin2}]
  The proofs of these propositions are completely analogous to those
  of Propositions~\ref{P108},~\ref{P105} and~\ref{metrictor}. The only
  point which is not clear is the fact that we may approximate any
  divisorial valuation by a sequence of infinitely singular
  valuations.  It is an immediate  consequence of the lemma below.
\end{proof}
\begin{lemma}\label{L777}
  For any divisorial valuation $\nu$ and $\e >0$, there exists an
  infinitely singular valuation $\mu > \nu$, with $A(\mu) < A(\nu) + \e$.
\end{lemma}
\begin{proof}
  Let us extend the approximating sequence $(\nu_i)_0^g$ 
  of $\nu$ to an infinite approximating sequence $(\nu_i)_0^\infty$
  as follows. First suppose $b(\nu)>m(\nu)$.
  By convention, $\nu_{g+1}=\nu$. 
  Pick a curve valuation $\mu_{g+2}>\nu$ with
  $m(\mu_{g+2})=b(\nu)$. Consider a divisorial valuation
  $\nu_{g+2}$ in the segment $]\nu,\mu_{g+2}[$. 
  It follows from~\eqref{e707} that
  if $\nu_{g+2}$ is close enough to $\nu$, then 
  $b(\nu_{g+2})>m(\nu)$. We may therefore pick $\nu_{g+2}$
  in this segment such that $b(\nu_{g+2})>m(\nu_{g+2})=b(\nu)$
  and $\a(\nu_{g+2})-\a(\nu)<\e/(2b(\nu))$.

  Inductively, given $(\nu_i)_1^{g+k}$ we construct $\nu_{g+k+1}$
  divisorial with $\nu_{g+k+1}>\nu_{g+k}$, 
  $b(\nu_{g+k+1})>m(\nu_{g+k+1})=b(\nu_{g+k})$ and
  $\a(\nu_{g+k+1})-\a(\nu_{g+k})<\e/(2^kb(\nu_{g+k}))$.
  Then $(\nu_i)_0^\infty$ defines an approximating sequence for
  an infinitely singular valuation $\mu$, satisfying
  $\mu>\nu$ and $A(\mu)<A(\nu)+\e$.

  If $b(\nu)=m(\nu)$ then we redefine $\nu_{g+1}$ to 
  be a divisorial valuation with $\nu_{g+1}>\nu$, 
  $b(\nu_{g+1})>m(\nu_{g+1})=m(\nu)$ and 
  $\a(\nu_{g+1})-\a(\nu)<\e/b(\nu)$. We may then continue as before
  and construct an infinite approximating sequence 
  $(\nu_i)_0^\infty$ of an infinitely singular valuation $\mu$
  with $\mu>\nu$ and $A(\mu)<A(\nu)+2\e$.
\end{proof}
\begin{remark}\label{R408}
  Using the relative thinness defined in 
  Section~\ref{sec-relative} we can define
  \index{thin!topology (relative case)}
  \index{thin!metric (relative case)}
  relative thin topologies on the trees $\cV_x$ and
  $\cVqmx$.
  Compare Remarks~\ref{R406} and~\ref{R407}.
\end{remark}
\begin{remark}\label{R409}
  Although we shall not pursue this further here, the analysis in 
  Chapter~\ref{sec-puis} implies that the (relative) thin topologies 
  can be defined in terms of (suitably normalized) 
  uniform convergence when the valuations are extended
  to the ring of power series in one variable with 
  coefficients that are Puiseux series.
\end{remark}
%
%
%
%
\section{The Zariski topology}\label{sec-zar}
We now turn to a classical, but quite different construction,
the Zariski topology. 
This topology, which is defined on the set $\cV_K$, 
of equivalence classes of Krull valuations, 
is a not  Hausdorff topology since divisorial valuations do
not define closed points.
We show how to make it Hausdorff by identifying a divisorial valuation 
with the valuations in its closure. 
The latter valuations are exactly the 
exceptional curve valuations, or,
equivalently, the elements of $\cV_K\setminus\cV$.
It is a remarkable fact that this procedure 
recovers $\cV$ endowed with the weak topology. 
%
%
\subsection{Definition}
The Zariski topology is defined on the set 
$\cV_K$ of equivalence classes of centered Krull valuations on $R$ 
(not necessarily $\Rbar$-valued, see Section~\ref{krullval}).
A Krull valuation $\nu$ is determined, up to equivalence,
by its valuation ring 
$R_\nu=\{\nu\ge0\}\subset K$, so an open set in 
$\cV_K$ can also be viewed as a
set of valuation rings satisfying certain conditions.
\begin{definition}
  A basis for the Zariski topology\index{topology!Zariski}
  on $\cV_K$ is given by
  \begin{equation*}
    V(z)
    \=\{\nu\ ;\ \nu(z)\ge 0\}
    =\{\nu\ ;\ z\in R_\nu\}
  \end{equation*}
  over $z\in K$.
  In other words, an arbitrary open set in the Zariski topology is
  a union of finite intersections of sets of the type $V(z)$.
\end{definition}
\begin{remark}
  The topological space $\cV_K$\index{Riemann-Zariski variety}
  endowed with the Zariski topology is called the 
  \emph{Riemann-Zariski variety} of $R$ 
  (see~\cite[p.110]{ZS},~\cite{Vaquie}).
\end{remark}
\begin{proposition}\label{closure-points}
  $\cV_K$ is quasi-compact but not Hausdorff. Moreover,
  \begin{itemize}
  \item[(i)]
    Any non-divisorial valuation is a closed point.
  \item[(ii)]
    The closure of a divisorial valuation associated to an 
    exceptional component $E$ 
    is the set of exceptional curve valuations $\nu_{E,p}$ whose
    centers lie on $E$.
  \end{itemize}
\end{proposition}
\begin{proof}
  That $\cV$ is not Hausdorff follows from~(ii). For the
  quasi-compactness of $\cV_K$, see Theorem~\ref{T601} below.
  Both~(i) and (ii) are consequences of Lemma~\ref{L601} since a
  valuation $\mu$ lies in the closure of another valuation $\nu\ne\mu$
  iff $R_\mu\subsetneq R_\nu$.  
\end{proof}
\begin{remark}
  $\cV\subset\cV_K$ is neither open nor closed for the 
  Zariski topology. 
  Indeed, Example~\ref{E403} shows that a sequence of divisorial
  valuations can converge to an exceptional curve valuation.
  Conversely, if $E$ is an irreducible component of 
  $\pi^{-1}(0)$ for some composition of blowups $\pi$, 
  and $(p_j)_1^\infty$ is a sequence of distinct points on
  $E$, then the exceptional curve valuations 
  $\nu_{E,p_j}$ converge to the divisorial valuation 
  $\nu_E$. See Appendix~\ref{sec-tangent}.
\end{remark}
\begin{remark}
  The Zariski topology can be geometrically 
  described as follows: a basis
  is given by $V(\pi,C)$, where $\pi$ ranges over
  compositions of finitely many blow-ups and
  $C$ over Zariski open subsets 
  of the exceptional divisor $\pi^{-1}(0)$. 
  Here $C$ is the complement in $\pi^{-1}(0)$ 
  of finitely many irreducible components and 
  finitely many points,
  and $V(\pi,C)$ is the set of all $\nu\in\cV_K$ 
  whose associated sequence of infinitely 
  nearby points $\Pi[\nu]=(p_j)$ are eventually in $C$.
\end{remark}
%
%
\subsection{Recovering $\cV$ from $\cV_K$}\label{recover}
We can try to turn $\cV_K$ into a Hausdorff space $\wtV_K$ by
identifying $\nu$ and $\nu'$ if they both belong to
the Zariski closure of the same Krull valuation $\mu$. 
By Proposition~\ref{closure-points} this amounts to identifying
exceptional curve valuations with their associated divisorial valuation.  
Let $p:\cV_K\to\wtV_K$ be the natural projection
and endow $\wtV_K$ with the quotient topology.  
Then $\wtV_K$ is quasi-compact, being the image of a quasi-compact
space by a continuous map.

Consider the natural injection
$\imath:\cV\to\cV_K$,
where $\cV$ carries the weak topology.
\begin{theorem}\label{weak-zariski}
  The composition $p\circ\imath:\cV\to\widetilde{\cV}_K$ 
  is a homeomorphism.
\end{theorem}
\begin{proof}
  Since $\cV$ contains all divisorial valuations but no
  exceptional curve valuations,
  injectivity and surjectivity of $p\circ\imath$ follow from
  Proposition~\ref{closure-points}.
  Since $\cV$ is Hausdorff and 
  $\widetilde{\cV_K}$ is quasi-compact, 
  continuity of $(p\circ\imath)^{-1}$ 
  will imply that $p\circ\imath$ is in fact a homeomorphism.

  Therefore, let us show that $(p\circ\imath)^{-1}$ is continuous. 
  A basis for the weak (tree) topology is given by the 
  open sets $U(\vv)$, over tangent vectors $\vv$.
  In fact, it suffices to take $\vv\in T\nu$ with $\nu$ divisorial.
  For $\phi,\psi\in\fm$ irreducible and $t>0$ define
  \begin{equation*}
    V(\phi,\psi,t)
    :=\left\{
      \mu\in\cV\ ;\ \frac{\mu(\phi)}{m(\phi)}<t\frac{\mu(\psi)}{m(\psi)}
    \right\}.
  \end{equation*}
  Proposition~\ref{P201} implies that if 
  $\vv\in T\nu$ is not represented by
  $\nu_\fm$ (true \eg if $\nu=\nu_\fm$), then
  $U(\vv)=V(\phi,\psi,1)$,
  where $\nu_\psi$ represents $\vv$ and $\nu_\phi\wedge\nu_\psi=\nu$. 
  If instead $\vv$ is represented by $\nu_\fm$, then
  $U(\vv)=V(\phi,\psi,\a(\nu))$,
  where $\nu_\phi\ge\nu$ and $\nu_\psi\wedge\nu_\phi=\nu_\fm$.
  
  Hence it suffices to show that 
  $(p\circ\imath)(V)$ is open in
  $\widetilde{\cV_K}$ for any $V=V(\phi,\psi,t)$. 
  This amounts to $p^{-1}(p\circ\imath)(V)$ 
  being open in $\cV_K$. Define
  \begin{equation*}
    W=W(\phi,\psi,t)
    :=\bigcup_{p/q>t}\left\{
      \nu\in\cV_K\ ;\ \nu\left(\psi^{pm(\phi)}/\phi^{qm(\psi)}\right)\ge0
    \right\}.
  \end{equation*}
  Then $W$ is open in $\cV_K$ and 
  $\imath(V)=\imath(\cV)\cap W$. 
  We claim that if $\nu\in\cV_K$ 
  is divisorial and $\nu'\in\overline{\{\nu\}}$
  an exceptional curve valuation, then
  $\nu\in W$ iff $\nu'\in W$. 
  This claim easily implies that $p^{-1}(p\circ\imath)(V)=W$ is open,
  completing the proof.
  
  As for the claim, $\nu\in W$ implies
  $t':=\nu(\psi)m(\phi)/\nu(\phi)m(\psi)>t$.  Pick $t''\in(t,t')$
  rational, $t''=p/q$.  Then $\nu$ belongs to the closed set 
  $\{\mu\ ;\ \mu(\psi^{pm(\phi)}/\phi^{qm(\psi)})>0\}\subset W$, 
  hence so does $\nu'$.  Conversely, if $\nu'\in W$ then $\nu\in W$ as
  $R_{\nu'}\subset R_\nu$.
\end{proof}
%
%
%
%
\section{The Hausdorff-Zariski topology}\label{HZ}
We now recall a natural refinement of the Zariski topology, 
the Hausdorff-Zariski (or simply HZ) topology.
It is still defined on the set $\cV_K$ of equivalence classes 
of Krull valuations.
As we show, the HZ topology is the weak tree topology 
for a natural $\Nbar$-tree structure on $\cV_K$. 
%
%
\subsection{Definition}
Let $\cF=\{0,+,-\}^K$ be the set of functions from $K$ to $\{0,+,-\}$.
Any Krull valuation defines an element in $\cF$ 
by setting $\nu(\phi)=+,-,0$ iff $\nu(\phi)>0$, $\nu(\phi)=0$ and
$\nu(\phi)<0$, respectively. 
Further, two Krull valuations are equivalent iff they define the same 
element of $\cF$. (This is essentially equivalent to the fact that 
equivalence classes of Krull valuations are in 1-1 correspondence with
valuation rings in $K$.)

Hence we can consider $\cV_K$ as a subset of $\cF$.
It is easy to see that Zariski topology on $\cV_K$ is 
exactly the topology induced from the product topology on $\cF$
associated to the topology on $\{0,+,-\}$ whose open sets 
are given by $\emptyset$, $\{0,+\}$ and $\{0,+,-\}$. 
We define the 
\emph{Hausdorff-Zariski topology}
\index{topology!Hausdorff-Zariski (or HZ)} 
or simply \emph{HZ topology}
on $\cV_K$ to be the topology induced by the product topology 
associated to the discrete one on $\{0,+,-\}$. 
The HZ topology is Hausdorff by construction and we have
\begin{lemma}\label{HZconv}
  A sequence $\nu_n\in \cV_K$ converges towards $\nu$ in the
  Hausdorff-Zariski topology iff for all $\phi\in K$ 
  with $\nu(\phi)>0$, $\nu(\phi)=0$ and $\nu(\phi)<0$ 
  one has $\nu_n(\phi)>0$, $\nu_n(\phi)$ and $\nu_n(\phi)<0$,
  respectively, for sufficiently large $n$.
\end{lemma}
From Tychonov's theorem, and the fact that $\cV_K$ is closed in
$\cF$ (see \cite[p.114]{ZS}) one deduces the following 
fundamental result.
\begin{theorem}[{\cite[p.113]{ZS}}]\label{T601}
  The space $\cV_K$ endowed with the Hausdorff-Zariski topology is
  compact. Thus $\cV_K$ is quasi-compact in the Zariski topology.
\end{theorem}
%
%
\subsection{The $\Nbar$-tree structure on $\cV_K$}\label{Z-tree}
We now introduce a tree structure on $\cV_K$ which plays the
same role for the HZ topology as the $\R$-tree structure 
plays for the weak topology on $\cV$.
The new tree will be modeled on the
totally ordered set $\Nbar=\N\cup\{\infty\}$ 
(see Section~\ref{tree-def}).
\begin{definition}
  Define a partial ordering $\trianglelefteq$\index{$\trianglerighteq$
(ordering on $\fB$)} on $\cV_K$ by
  \begin{equation*}
    \nu_1\trianglelefteq\nu_2 \text{ iff the sequence of blow-ups }
    \Pi[\nu_2] \text{ contains } \Pi[\nu_1]. 
  \end{equation*}
\end{definition}
\begin{proposition}\label{blwuporder}
  The space $(\cV_K,\trianglelefteq)$ is a $\Nbar$-tree rooted
  at $\nu_\fm$. Any divisorial valuation is a branch point, 
  with tangent space in bijection with $\mathbf{P}^1$.
  The ends of the tree are exactly the non-divisorial valuations.
\end{proposition}
\begin{proof}
  That $(\cV_K,\trianglelefteq)$ is a $\Nbar$-tree is a
  straightforward consequence of the definition, 
  as is the fact that the set of ends coincides with the set of
  non-divisorial valuations.

  A tangent vector at a divisorial valuation $\nu$ is given by a 
  point on the exceptional component $E$ defining $\nu$. 
  Hence the tangent space at $\nu$ is in bijection with
  $E\simeq\mathbf{P}^1$. This completes the proof.
  \end{proof}  
\begin{proposition}\label{P401}
  The Hausdorff-Zariski topology coincides with the weak tree
  topology induced by $\trianglelefteq$.
\end{proposition}
\begin{proof}
  Pick a tangent vector $\vv$ at a divisorial
  valuation $\nu_0$, and consider the weak open set $U(\vv)$.  
  If $\nu_0$ corresponds to the exceptional curve $E$, 
  and $\vv$ to the point $p\in E$, 
  $U(\vv)$ coincides with the set of valuations
  that are centered at $p$.  
  Choose $\phi,\phi'\in K$ defining two smooth transversal 
  curves at $p$ disjoint from $E$. Then
  $U(\vv)=\{\nu\ ;\ \nu(\phi),\nu(\phi')>0\}$ is a weak tree open set. 
  Hence the identity map from $\cV_K$ endowed with the HZ topology onto
  $\cV_K$ endowed with the weak tree topology associated to
  $\triangleleft$ is continuous. 
  As $\cV$ is compact in the HZ topology and Hausdorff in the
  weak tree topology, it follows that the identity map is
  a homeomorphism.
\end{proof}
In Chapter~\ref{A3} we shall use sequences of infinitely nearby points
in a different way and actually recover the valuative tree.
The analysis in that chapter can be used to 
explain the precise relation between the $\Nbar$-tree and $\R$-tree 
structures on $\cV$ associated to $\trianglelefteq$ and $\le$, respectively.
We will contend ourselves with the following illustrative example. 
Fix local coordinates $(x,y)$ and 
consider the segment $I=[\nu_\fm,\nu_y[$ in $\cV$, 
consisting of monomial valuations in these coordinates
satisfying $1=\nu(x)\le\nu(y)$. 
The restriction of $\le$ to $I$ coincides with the natural order 
on $[1,\infty[$. On the other hand,
$\trianglelefteq$ gives the lexicographic order on the continued
fractions expansions of elements in $[1,\infty[$.
%
%
%
%
\section{Comparison of topologies}\label{topo-compare}
We now have several topologies on the trees $\cV$ and $\cVqm$ and on
the space $\cV_K$. We present here some comparisons between them. Our
objective is not to prove, or even state, all possible results
relating the different topologies, but to illustrate a few of the
connections.

We split our analysis into two parts. 
In the first part we describe the relationships between the weak,
the strong, the thin and the HZ topology on $\cV$. 
In the second part we compare the strong and thin metrics on
$\cV$ and $\cVqm$.
%
%
\subsection{Topologies}
By the HZ topology on $\cV$ we mean the 
topology of $\cV$ as a subset of $(\cV_K,HZ)$.
\begin{theorem}\label{comparison}
  The strong, the thin and the HZ topology
  on $\cV$ are all stronger than the weak topology. Moreover:
  \begin{itemize}
  \item[(i)]
    if $\nu$ has infinite skewness and $\nu_n\to\nu$
    weakly, then $\nu_n\to\nu$ strongly;
  \item[(ii)]
    if $\nu$ has infinite thinness and $\nu_n\to\nu$
    weakly, then $\nu_n\to\nu$ thinly;
  \item[(iii)]
    if $\nu$ is non-divisorial and $\nu_n\to\nu$
    weakly, then $\nu_n\to\nu$ in the HZ topology.
  \end{itemize}      
\end{theorem}
Essentially no other implications hold as the following examples indicate.
\begin{example}\label{E403}
  Let $\nu_n\=\nu_{y,1+n^{-1}}$. 
  Then $\nu_n\to\nu_\fm$ strongly and thinly (hence weakly). 
  But $\nu_n(y/x)=1/n>0$ 
  and $\nu_\fm(y/x)=0$, so $\nu_n\not\to\nu$ in the
  HZ topology. In fact, it converges to the
  exceptional curve valuation
  $\val[(x,y);((1,0),(1,1))]$.
\end{example}
In this example, the HZ limit valuation is in the
closure of $\nu_\fm$ in $\cV_K$. 
This is a general fact for limits of sequences in
the weak and HZ topology.
\begin{example}
  Pick $p_n$, $q_n$ 
  relatively prime with 
  $p_n/q_n\to\sqrt{2}$. 
  Set $\phi_n=y^{q_n}-x^{p_n}$, 
  $\nu_n=\nu_{\phi_n}$ and $\nu=\nu_{y,\sqrt{2}}$.
  Then $\nu_n\to\nu$ weakly 
  (hence in the HZ topology since $\nu$ is irrational)
  but $d(\nu_n,\nu)\to1/\sqrt{2}$ and $D(\nu_n,\nu)\to1/(1+\sqrt{2})$
  so $\nu_n$ does not converge to $\nu$ neither strongly nor thinly.
\end{example}
\begin{example}\label{E402}
  Let $\nu$ be an infinitely singular valuation with
  finite skewness but infinite thinness. 
  See Remark~\ref{R405} for how to construct $\nu$.
  Pick a sequence $(\mu_n)_1^\infty$ of divisorial valuations
  increasing to $\nu$ and for each $n$ pick a curve valuation
  $\nu_n$ with $\nu_n\wedge\nu=\mu_n$.
  Then $D(\nu_n,\nu)=2A(\mu_n)^{-1}\to0$ so $\nu_n\to\nu$ thinly. 
  However, $d(\nu_n,\nu)\to\a(\nu)^{-1}>0$ so 
  $\nu_n\not\to\nu$ strongly.
\end{example}
\begin{example}\label{E401}
  It is possible to construct a sequence $(\nu_n)_1^\infty$
  of infinitely singular (or even divisorial) valuations such that 
  $A(\nu_n)\ge3$ but $\a(\nu_n)\to1$ as $n\to\infty$.
  See Remark~\ref{R404}. 
  Then $\nu_n\to\nu_\fm$ strongly but not thinly.
\end{example}
\begin{proof}[Proof of Theorem \ref{comparison}]
  We know from Proposition~\ref{P430} that the strong and thin
  topologies are both stronger than the weak topology.
  To compare the HZ and weak topologies we consider the 
  natural injection $\jmath:(\cV,\mathrm{HZ})\to(\cV_K,\mathrm{Z})$,
  which is continuous as the HZ topology is stronger than the Zariski
  topology. By Theorem~\ref{weak-zariski} there is a continuous
  mapping $q:(\cV_K,\mathrm{Z})\to(\cV,\mathrm{weak})$ which is the
  identity on $\cV\subset\cV_K$.  Since $q\circ\jmath=\id$, we see
  that $\id:(\cV,\mathrm{HZ})\to(\cV,\mathrm{weak})$ is continuous.
  Hence the HZ topology is stronger than the weak topology.

  For~(i), suppose $\nu_n\to\nu$ weakly.
  As $\a(\nu)=\infty$ we may find $\phi_k\in\fm$ such that 
  $\nu(\phi_k)\ge(k+1)m(\phi_k)$. 
  For $n\gg1$, $\nu_n(\phi_k)\ge k\,m(\phi_k)$.
  Thus $\nu,\nu_n\ge\nu_{\phi_k,k}$, so
  \begin{equation*}
    d(\nu_n,\nu)
    \le d(\nu_n,\nu_{\phi_k,k})
    +d(\nu_{\phi_k,k},\nu)
    \le k^{-1}+k^{-1}.
  \end{equation*}
  We let $k\to\infty$ and conclude that
  $\nu_n$ converges strongly towards $\nu$.

  The proof of~(ii) is essentially identical.
  As for~(iii), suppose $\nu_n\rightarrow\nu$ weakly. 
  Consider $\phi\in K$. If $\nu(\phi)>0$ ($<0$), then $\nu_n(\phi)>0$
  ($<0$) for $n\gg1$. If $\nu(\phi)=0$, then as the
  residue field $k_\nu$ is isomorphic to $\C$, we may write
  $\phi=\lambda+\psi$, where 
  $\lambda\in\C^*$ and $\nu(\psi)>0$.
  For $n\gg1$ $\nu_n(\psi)>0$ so that $\nu_n(\phi)=0$.
  Thus $\nu_n\to\nu$ in the HZ topology.
\end{proof}
%
%
\subsection{Metrics}
We have seen above that the strong and thin topologies on $\cV$
are distinct. Nevertheless, we can compare the thin and the strong metrics
whenever we have a bound on the multiplicity.
\begin{proposition}\label{comparison3}
  Fix $m\ge1$. Then the inequalities
  \begin{equation*}
    d_\mathrm{qm}\le D_\mathrm{qm}\le m\, d_\mathrm{qm}
  \end{equation*}
  hold on the subtree $\{\nu\in\cVqm\ ;\ m(\nu)\le m\}$ 
  of $\cVqm$ and these estimates are sharp.
  The identity map
  $(\cVqm,D_\mathrm{qm})\to(\cVqm,d_\mathrm{qm})$ is continuous,
  but its inverse is not.
\end{proposition}
\begin{proposition}\label{comparison2}
  Fix $m\ge 1$. Then the inequalities
  \begin{equation*}
    m^{-1}D\le d\le 2(m+1)\,D,
  \end{equation*}
  hold on the subtree $\{\nu\in\cV\ ;\ m(\nu)\le m\}$ of $\cV$.
  The identity map $(\cV,D)\to(\cV,d)$ is not continuous, and neither
  is its inverse.
\end{proposition}
\begin{proof}[Proof of Proposition~\ref{comparison3}]
  The estimate is an immediate consequence of the definition of 
  $d_\mathrm{qm}$ and $D_\mathrm{qm}$ and of~\eqref{e401}.
  That the identity map $(\cVqm,D_\mathrm{qm})\to(\cVqm,d_\mathrm{qm})$
  is continuous follows immediately.
  However, its inverse is not continuous in view of Example~\ref{E401}.
\end{proof}
\begin{proof}[Proof of Proposition~\ref{comparison2}]
  As $d$ and $D$ are both tree metrics, it suffices to consider
  $\nu_1\le\nu_2$ with $m(\nu_1)\le m$ 
  and compare $d(\nu_1,\nu_2)$ with $D(\nu_1,\nu_2)$.
  Write $A_i=A(\nu_i)$, $\a_i=\a(\nu_i)$ and $m_i=m(\nu_i)$ 
  for $i=1,2$. 
  Also write $d_{12}=d(\nu_1,\nu_2)$ and $D_{12}=D(\nu_1,\nu_2)$.
  Then $d_{12}=\a_1^{-1}-\a_2^{-1}=(\a_2-\a_1)/(\a_1\a_2)$
  and  $D_{12}=(A_2-A_1)/(A_1A_2)$.
  Moreover,~\eqref{e401} gives
  $m_1(\a_2-\a_1)\le A_2-A_1\le m(\a_2-\a_2)$
  and Proposition~\ref{bounds} implies
  $A_i\le m_i\a_i+1$. Hence
  \begin{equation*}
    D_{12}
    \le\frac{m(\a_2-\a_1)}{A_1A_2}
    \le m \frac{\a_2-\a_1}{\a_1\a_2}
    =m\,d_{12}
  \end{equation*}
  and
  \begin{equation*}
    D_{12}
    \ge\frac{m_1(\a_2-\a_1)}{(m_1\a_1+1)(m_2\a_2+1)}
    \ge\frac{m_1\a_1}{m_1\a_1+1}\frac{\a_2}{m_2\a_2+1}\,d_{12}
    \ge\frac{d_{12}}{2(m+1)}.
  \end{equation*}
  where we used $m_1,m_2\le m$ and $\a_1,\a_2\ge1$.

  That the identity maps $(\cV,D)\to(\cV,d)$
  and $(\cV,d)\to(\cV,D)$ are both discontinuous follows from
  Example~\ref{E402} and Example~\ref{E401}.
\end{proof}

%
%
%
%
%
%
\chapter{The universal dual graph}\label{A3}
We have already described several different approaches to the 
valuative tree.
They all fundamentally derive from the definition of a valuation as
a function on the ring $R$.
On the other hand, valuations can also be viewed geometrically
as sequences of infinitely nearby points.
It is therefore natural to ask whether the valuative tree can 
be recovered through a purely geometric construction.

In this chapter we show that this is indeed possible. 
The construction goes as follows. 
To any composition of (point) blowups we associate the dual graph 
of the exceptional divisor. 
The set of vertices of this graph is naturally a poset
and the collection of all such posets forms an injective system
whose limit is a nonmetric tree $\Gast$
modeled on the rational numbers. 
By filling in the irrational points and adding all the ends
we obtain a nonmetric tree $\Gamma$ modeled on the real line.
We call $\Gamma$ the \emph{universal dual graph}.
Its points are encoded by sequences of 
infinitely nearby points above the origin.
We show how to equip $\Gamma$ with a natural 
\emph{Farey parameterization} as well as an integer
valued \emph{multiplicity} function.

The main result is then that there exists a natural isomorphism
from the universal dual graph $\Gamma$ to the
valuative tree $\cV$.
Its inverse maps a valuation to its associated sequence of 
infinitely nearby points as defined in Section~\ref{sec-equivalence}.

The chapter is organized as follows. We start out by
defining $\Gamma$ as a nonmetric tree,
then equip it with the Farey parameterization and multiplicity 
function, respectively.
Along the way we show that sequences of infinitely nearby points 
correspond uniquely to either points in $\Gamma$, or tangent
vectors at branch points in $\Gamma$.
We then state and prove the main result, 
namely the isomorphism between $\Gamma$ and $\cV$.  
From this we deduce a number of applications,
and show how the dual graph of the minimal desingularization of
a reduced curve can be described inside the universal dual graph.
Finally we discuss a
relative version of the universal dual graph, analogous to the
relative valuative tree, and explain the
self-similar structure of these objects.
%
%
%
%
\section{Nonmetric tree structure}\label{sec-constr-univ}
We first define the universal dual graph as a nonmetric tree. 
Later we will equip it with a natural parameterization and
multiplicity function.
%
%
\subsection{Compositions of blowups}
Let us denote by $\fB$
\index{$\fB$ (modifications above the origin)} 
the set of all modifications above
the origin.  This means that each element of $\fB$ is a mapping
$\pi:X_\pi\to(\C^2,0)$ where $X_\pi$ is a smooth complex surface and
$\pi$ is a proper map which is a bijection outside the exceptional
divisor $\pi^{-1}(0)$.  It is well-known that each such $\pi$ is a
composition of (point) blowups 
(see~\cite[Theorem~5.7]{laufer} for instance). 
We emphasize that these blowups are not necessarily
associated to a sequence of infinitely nearby points.

The set $\fB$ has a natural partial ordering:
$\pi\trianglerighteq\pi'$\index{$\trianglerighteq$ (ordering on
$\fB$)} iff $\pi=\pi'\circ\tpi$ for some composition of blowups (see
Section~\ref{Z-tree}). We shall use repeatedly
\begin{lemma}\label{Lsup-inf}
  The set $\fB$ is naturally an inverse system: any finite subset of
  $\fB$ admits a supremum; and every nonempty subset of $\fB$ admits an
  infimum\index{infimum!in $\fB$}.
\end{lemma}
\begin{proof}
  To prove that any finite subset admits a supremum, it is sufficient to
  consider the case of two elements. If $\pi,\pi'$ belong to $\fB$, the
  identity map lifts to a bimeromorphic map 
  $\id:X_\pi\to X_{\pi'}$. 
  It is a classical fact that for any bimeromorphic map of
  surfaces there exist compositions of point blow-ups 
  $\varpi:Y\to X_\pi$, $\varpi':Y\to X_{\pi'}$ 
  such that $\pi\circ\varpi=\pi'\circ\varpi'$. 
  (Take a resolution of singularities of the graph of
  $\id$, and apply~\cite[Theorem 5.7]{laufer} for instance). 
  This proves that any finite set admits a supremum.
  
  In the discussion above, one can be more precise. By taking the
  minimal desingularization (see~\cite[Theorem 5.9]{laufer}) of the
  graph of $\id$, one can show that $Y$ is ``minimal'' in the sense that
  any other surface $Y'$ dominating both $X_\pi$ and $X_{\pi'}$
  dominates $Y$. We call the natural map sending $Y$ to the base
  $(\C^2,0)$ the \emph{join} of $\pi$ and $\pi'$.
  
  Now take any non-empty subset $B\subset\fB$, 
  and define 
  $M\=\{\pi\ ;\ \pi\trianglelefteq\varpi\ \text{for all}\ \varpi\in B\}$.
  The set $M$ contains the blow-up of the origin and is hence non-empty. 
  It is also a finite set, as any map $\pi\in \fB$ dominates at 
  most finitely many other elements in $\fB$. 
  Then the join of the collection of all $\pi\in M$ is
  dominated by all elements in $B$. 
  This is the infimum of $B$.
\end{proof}
%
%
\subsection{Dual graphs}\label{sec-dualgraphs}
To any $\pi\in\fB$ 
we may attach a \emph{dual graph}
\index{dual graph}
\index{$\Gamma_\pi$ (dual graph, simplicial tree)}
$\Gamma_\pi$: vertices are in bijection with irreducible
components of the exceptional divisor $\pi^{-1}(0)$
(we shall refer to these as 
\emph{exceptional components}\index{exceptional component}), 
and edges with intersection points between two components. 
By decomposing an element $\pi\in \fB$ into a composition of point
blowups, we see that the dual graph $\Gamma_\pi$ can be obtained
inductively as follows.

If $\pi$ is a single blowup of the origin, then $\Gamma_\pi$ is a
single point that we denote by $E_0$.  Otherwise we may write
$\pi=\pi'\circ\tpi$, where $\tpi$ is the blowup of a point $p$ on the
exceptional divisor $(\pi')^{-1}(0)$, resulting in a new exceptional
component $E_p$.  The dual graph of $\pi$ can then be obtained from
the dual graph $\Gamma_{\pi'}$ of $\pi'$ through an \emph{elementary
modification}\index{elementary modification}.  There are two kinds of
elementary modifications, depending on the location of $p$.  We say
that $p$ is a \emph{free}\index{point (on an exceptional
divisor)!free} (\emph{satellite}\index{point (on an exceptional
divisor)!satellite}) point if $p$ is a regular (singular) point on
$(\pi')^{-1}(0)$.

\begin{figure}[ht]
  \begin{center}
    \includegraphics[width=\textwidth]{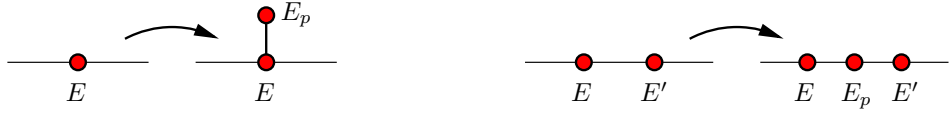}
  \end{center}
  \caption{Elementary modifications on the dual graph. 
    The first kind of modification is illustrated to the left, 
    the second kind to the right.}\label{F10}
\end{figure}

If $p$ is free, \ie $p\in E$ for a unique irreducible component
$E\subset(\pi')^{-1}(0)$, then the elementary modification consists
of adding to $\Gamma_{\pi'}$ a segment joining $E$ and the 
new vertex $E_p$. We say that this modification
is of the \emph{first kind}. See Figure~\ref{F10}.

If $p$ is a satellite point, it is the intersection
of two irreducible components $E$ and $E'$. 
The elementary modification now consists of adding a new vertex $E_p$ 
to the segment between $E$ and $E'$; we say that the modification
is of the \emph{second kind}. Again see Figure~\ref{F10}.

From this description, it follows that $\Gamma_\pi$ is a finite
simplicial tree. We shall denote by $E_0$ the vertex associated to the
proper transform in $\Gamma_\pi$ of the exceptional divisor obtained
by blowing up the origin. 
%
%
\subsection{The $\Q$-tree}\label{S310}
Following Section~\ref{S303} we can view the simplicial
trees $\Gamma_\pi$ as $\N$-trees. This goes as follows.
For $\pi\in\fB$, consider the set $\Gast_\pi$ of vertices of
\index{$\Gast_\pi$ (dual graph, $\N$-tree)}
$\Gamma_\pi$. It is a finite set, with a natural partial ordering
$\le_\pi$ derived from the simplicial tree structure on $\Gamma_\pi$,
rooted at $E_0$: $E_1\le_\pi E_2$ iff $[E_0, E_1]\subset [E_0,E_2]$ as
segments in $\Gamma_\pi$. The poset $\Gast_\pi$ is an $\N$-tree in
the sense of Chapter~\ref{part3}.

If $\pi,\pi'\in\fB$ and $\pi\trianglerighteq\pi'$, then the dual graph
$\Gast_\pi$ is obtained from $\Gast_{\pi'}$ by performing a sequence
of elementary modifications.  In particular, there is a natural
injective map $\imath_{\pi\pi'}:\Gast_{\pi'}\to\Gast_\pi$, and if
$\pi''\trianglelefteq\pi'\trianglelefteq\pi$, then
$\imath_{\pi\pi'}\circ\imath_{\pi'\pi''}=\imath_{\pi\pi''}$

We claim that each map $\imath_{\pi\pi'}$ is order-preserving. 
By induction it suffices to show this when $\pi=\pi'\circ\tpi$,
where $\tpi$ is the blowup of a point $p\in(\pi')^{-1}(0)$.
Then $\imath_{\pi\pi'}$ is given by an elementary modification
of the first or second kind: see Figure~\ref{F10}. 
If the elementary modification is of the first kind, \ie if
$E_p\in\Gast_{\pi}\setminus\imath_{\pi\pi'}(\Gast_{\pi'})$ 
is obtained by blowing up a free point $p$ on $E\in\Gast_{\pi'}$, 
then $E_p>\imath_{\pi\pi'}(E)$, $E_p$ is
an end in $\Gast_\pi$, and the segment $]\imath_{\pi\pi'}(E),E_p[$ in
$\Gast_\pi$ is empty.  If the elementary modification is of the second
kind, \ie if $E_p$ is instead obtained by blowing up the intersection
point of two elements $E,E'\in\Gast_{\pi'}$, say with $E<E'$, then
$\imath_{\pi\pi'}(E)<E_p<\imath_{\pi\pi'}(E')$ 
and the segments $]\imath_{\pi\pi'}(E),E_p[$ and
$]E_p,\imath_{\pi\pi'}(E')[$ in $\Gast_\pi$ are empty.

Thus $\imath_{\pi\pi'}$ is order-preserving
whenever $\pi'\trianglelefteq\pi$. 
Slightly abusively, we will consider 
$\Gast_{\pi'}$ to be a subset of $\Gast_\pi$
and ignore the map $\imath_{\pi\pi'}$.

The set $\fB$ defines an inverse system,  
and $(\Gast_\pi,\le_\pi)_{\pi\in\fB}$ forms an injective system.
We may therefore define the 
\emph{universal dual graph} $\Gast$
\index{dual graph!universal} 
\index{universal dual graph!$\Q$-tree} 
\index{$\Gast$ (universal dual graph, $\Q$-tree)} 
as the injective limit
\begin{equation*}
  (\Gast,\le)=\injlim_{\pi\in\fB}(\Gast_\pi,\le_\pi)
\end{equation*}
over all sequences of blow-ups $\pi$ above the origin.
Again we consider (slightly abusively) $\Gast_\pi$ as
a subset of $\Gast$ for all $\pi$.
Then $\Gast$ is the union of all the $\Gast_\pi$'s.
It is important to note that any finite subset of $\Gast$
is contained in $\Gast_\pi$ for some $\pi$.

The partial ordering on $\Gast$ is determined as follows. 
\index{partial ordering!on the universal dual graph}
If $E,E'\in\Gast$, then $E,E'\in\Gast_\pi$ for some $\pi\in\fB$ and 
$E\le E'$ in $\Gast$ iff $E\le E'$ in $\Gast_\pi$. 
\begin{proposition}\label{P407}
  The universal dual graph $(\Gast,\le)$ is a nonmetric $\Q$-tree rooted
  at $E_0$. All its points are branch points.
\end{proposition}
We shall describe the tangent spaces in 
Proposition~\ref{Ptangentgamma} below.
\begin{proof}
  Let us verify that the partial ordering on $\Gast$ 
  satisfies the axioms (T1)-(T3) for a rooted nonmetric
  $\Q$-tree on page~\pageref{D402}. 
  As for~(T1), this is clear: $E_0$ is the unique minimal element 
  of $\Gast$. 
 
  Next consider~(T2). Fix $E\in\Gast$.  In view of Lemma~\ref{L421} it
  suffices to show that $[E_0,E]:=\{F\in\Gast\ ;\ F\le E\}$ is totally
  ordered, countable, and has no gaps.  It is totally ordered as the
  intersection with any $\Gast_\pi$ is.  It cannot have any gaps,
  since if $E_1<E_2\le E$ and there was no $E'\in\Gast$ with
  $E_1<E'<E_2$, then $E_1$ and $E_2$ would be adjacent vertices in
  some $\Gamma_\pi$. We could then blow up the intersection point
  between $E_1$ and $E_2$ and obtain $E'\in\Gast$ with $E_1<E'<E_2$, a
  contradiction.  Thus the set $[E_0,E]$ has no gaps.  The hardest
  part is to show that it is countable.  Consider the minimal
  $\pi\in\fB$ such that $E\in\Gast_\pi$.  We shall describe the
  structure of $\Gast_\pi$ in Section~\ref{S404}, but for now it
  suffices to notice that any $F\in\Gast$ belonging to the segment
  $[E_0,E]$ can be obtained by performing finitely many satellite
  blowups at intersection points of elements of
  $[E_0,E]\cap\Gast_\pi$.  More precisely, we define a sequence of
  proper birational morphisms
  $\pi=\pi_0\trianglerighteq\pi_1\trianglerighteq\dots$ inductively
  through their associated posets $\Gast_n=\Gast_{\pi_n}$ as follows:
  the elements of $\Gast_{n+1}-\Gast_n$ are obtained by blowing up all
  intersection points of elements of $[E_0,E]\cap\Gast_n$. Then each
  $\Gast_n$ is a finite poset and
  $[E_0,E]=\bigcup_n[E_0,E]\cap\Gast_n$.  Thus $[E_0,E]$ is countable
  and~(T2) holds.

  Instead of proving~(T3) we prove the equivalent statement~(T3')
  in Remark~\ref{R401}. Thus consider an unbounded, 
  totally ordered subset $\cS$ of $\Gast$. We will prove that 
  $\cS$ admits an increasing sequence in $\cS$ without majorant 
  in $\Gast$. In doing so we will make use of the fact that
  if $E\in\Gast$, then there exists $n=n(E)<\infty$ such that 
  if $F\in\Gast$ with $F\le E$, then the minimal $\pi\in\fB$
  for which $F\in\Gast_\pi$ is a composition of blowups,
  at most $n$ of which are free (the number of satellite
  blowups can be arbitrarily high).
  Now consider $\cS$ as above. 
  Pick any $E_1\in\cS$, $E_1>E_0$. Let $\pi_1\in\fB$
  be minimal such that $E_1\in\Gast_1:=\Gast_{\pi_1}$.
  After replacing $E_1$ by $\max\Gast_1\cap\cS$ we may assume
  that $E_1=\max\Gast_1\cap\cS$. 
  Since $\cS$ is unbounded, there exists $E_2\in\cS$ with
  $E_2>E_1$. 
  Notice that the assumption $E_1=\max\Gast_1\cap\cS$
  implies that $E_2$ is obtained
  by blowing up a free point on $E_1$ followed by finitely many
  (free or satellite) blowups.
  As before we may assume that $E_2=\max\Gast_2\cap\cS$,
  where $\Gast_2=\Gast_{\pi_2}$ and $\pi_2\ge\pi_1$ is the
  minimal element of $\fB$ such that $E_2\in\Gast_2$.
  Inductively we construct an increasing sequence $(E_n)_1^\infty$
  in $\cS$ with the property that
  $E_{n+1}$ is obtained by blowing up a free point on $E_n$
  followed by finitely many (free or satellite) blowups.
  The remark above then implies that the sequence $(E_n)$ 
  is unbounded in $\Gast$, \ie there is no $E\in\Gast$
  with $E_n<E$ for all $n$.
  This proves~(T3') and hence~(T3).

  Thus $\Gast$ is a nonmetric $\Q$-tree rooted in $E_0$.
  Finally notice that any point $E\in\Gast$ is a branch point 
  since it becomes a branch point in some $\Gast_\pi$ 
  after blowing up two or more free points on $E$. 
  This concludes the proof.
\end{proof}
%
%
\subsection{Tangent spaces} 
The following geometric interpretation of the
tangent spaces in $\Gast$ will play an important role in
the sequel.
\begin{proposition}\label{Ptangentgamma}
  Let $E\in\Gast$ and pick $\pi\in\fB$ such that $E\in\Gast_\pi$.
  For $p\in E$, let $E_p\in\Gast$ be the exceptional component
  obtained by blowing up $p$, and denote by $\vv_p$ the tangent vector
  represented by $E_p$ at $E$.
  
  Then the map $p\to\vv_p$ induces a bijection between 
  the set of points in $E$ and the tangent space at $E$ in $\Gast$. 
\end{proposition}
\begin{proof}
  Pick a tangent vector $\vv$ at $E$. 
  By construction, $\vv$ is represented by some point in $\Gast$. 
  We may hence find $\pi'\in\fB$ such that 
  $\pi'\trianglerighteq\pi$, $E\in\Gast_{\pi'}$,
  and $\vv$ is represented by another component
  $F\in\Gast_{\pi'}$. 
  There is then a unique component $E'$
  intersecting $E$ and lying in the segment $[E,F]$.
  This component represents $\vv$ at $E$. 
  This is also the case for the exceptional
  component $E_p$ obtained by blowing up $p=E\cap E'$.  
  If $\pi'=\pi\circ\varpi$, $\vv$ is also represented by $E_q$ where
  $q\=\varpi(p)$. Whence $p\to\vv_p$ is surjective.
  
  To show that $p \to \vv_p$ is injective, pick $p\ne q$. 
  Then in the dual graph of the composition of $\pi$ 
  with the blow-ups at $p$ and $q$, 
  the component $E_p$ and $E_q$ represent different tangent vectors.
  This dual graph embeds in $\Gast$ by construction, 
  hence $\vv_p\ne\vv_q$.
\end{proof}
From the proof we deduce the following consequence:
\begin{corollary}\label{C309}
  Let $\pi$, $E$ and $p$ be as above, and consider a proper birational
  morphism $\pi':X_{\pi'}\to(\C^2,0)$ dominating $\pi$, \ie
  $\pi'=\pi\circ\varpi$. Assume that $\varpi$ is an isomorphism above
  $X_\pi\setminus\{p\}$.  Then all points in
  $\Gast_{\pi'}\setminus\Gast_\pi$ represent the tangent vector
  $\vv_p$ at $E$.
\end{corollary}
%
%
\subsection{The $\R$-tree}
In view of Section~\ref{sec-qtree} there is a canonical 
rooted, nonmetric $\R$-tree 
(or simply rooted, nonmetric tree) $\Gamma^o$ associated to the
rooted, nonmetric $\Q$-tree $\Gast$. 
The points in $\Gamma^o\setminus\Gast$
can be viewed as decreasing sequences of closed segments in $\Gast$
with empty intersection in $\Gast$.  We will somewhat abusively refer
also to $\Gamma^o$,
\index{$\Gamma^o$ (nonmaximal points in $\Gamma$)}
and even its completion $\Gamma$,
\index{$\Gamma$ (universal dual graph, $\R$-tree)} 
\index{universal dual graph!$\R$-tree}
as the universal dual graph.
In Section~\ref{S404} we shall interpret the elements of $\Gamma$ as
sequences of infinitely nearby points.
\begin{proposition}\label{P417}
  The universal dual graph $\Gamma$ is a complete nonmetric tree
  rooted in $E_0$ whose branch points are exactly the points in $\Gast$.
\end{proposition}
\begin{proof}
  Immediate consequence of Proposition~\ref{P420} and the definition
  of $\Gamma$ as the completion of $\Gamma^o$.
\end{proof}
We shall denote the infimum in $\Gamma$ by ``$\wedge$'', just as in 
\index{$\wedge$ (infimum)}
\index{infimum!in $\Gamma$}
the valuative tree.
\begin{remark}\label{R305}
  The tangent space in the $\R$-tree 
  $\Gamma$ at a branch point $E\in\Gast$
  is by construction canonically identified with the tangent 
  space at $E$ in the $\Q$-tree $\Gast$. 
  The latter tangent space 
  is described in Proposition~\ref{Ptangentgamma}: see
  also Theorem~\ref{div-tangent}.
\end{remark}
%
%
%
%
\section{Infinitely nearby points}\label{S404}
Next we show that the points in $\Gamma$ are encoded by sequences of
infinitely nearby points.  Recall that sequences of infinitely nearby
points are in bijection with Krull valuations centered at $\fm$ by
Theorem~\ref{blw}. We shall later see that $\Gamma$ is indeed
isomorphic to the valuative tree (Theorem~\ref{thm-universal}).
%
%
\subsection{Definitions and main results}
In general we can classify sequences of infinitely nearby points 
into five categories. The terminology below is essentially in accordance
\index{infinitely nearby points}
with that of Spivakovsky~\cite{spiv}.
\begin{definition}\label{D403}
  Let $\bp=(p_j)_0^n$, $0\le n\le\infty$
  be a finite or infinite sequence of infinitely
  nearby points. We say that $\bp$ is of
  \begin{itemize}
  \item\emph{Type~0}
    if $\bp$ is finite;
  \item\emph{Type~1}
    if $\bp$ is infinite and contains infinitely many free and
    infinitely many satellite points;
  \item\emph{Type~2}
    if $\bp$ is infinite, contains only finitely many free points,
    and is not of Type~3;
  \item\emph{Type~3}
    if $\bp$ is infinite, contains only finitely many free points,
    and has the following property: 
    there exists (a unique) $j_0\ge 1$ such that if $j>j_0$, 
    then $p_{j+1}$ is the satellite point defined by the intersection
    of $E_j$ and the strict transform of $E_{j_0}$;
  \item\emph{Type~4}
    if $\bp$ is infinite and contains only finitely many 
    satellite points.
  \end{itemize}
\end{definition}
We describe in Figures~\ref{F3} to~\ref{F17} below the structure of 
the dual graph appearing in the successive blow-ups associated 
to a sequence of infinitely nearby points. 
The notation is explained in Section~\ref{proof-infsg}.

\begin{definition}
  Fix a finite sequence $(p_j)_0^n$
  of infinitely nearby points (\ie of Type~0). 
  Write $\pi_\bp\in\fB$ for the composition of blow-ups
  at all points $p_0,\dots,p_n$, and define
  $\gamma(\bp)\in\Gast$ to be the exceptional divisor of the blow-up 
  at $p_n$.
\end{definition}
\begin{theorem}\label{Ttype0}
  The map $\bp\to\gamma(\bp)$ gives a bijection between 
  sequences of infinitely nearby points of Type~0, 
  and $\Gast$.
\end{theorem}
\begin{theorem}\label{Tnot3}
  Let $\bp = (p_j)_0^\infty$ be an infinite sequence of infinitely
  nearby points not of Type~3, and define the truncation
  $\bp_n=(p_j)_0^n$ for $n<\infty$.  Then the sequence
  $\gamma(\bp_n)\in\Gast$ converges weakly\footnote{We here use the
  weak topology induced by the tree structure on $\Gamma$.}  in
  $\Gamma$ to an element $\gamma(\bp)$. We have:
  \begin{itemize}
  \item
    $\gamma(\bp)$ is an end in $\Gamma$ if
    $\bp$ is of Type~1 or~4;
  \item
    $\gamma(\bp)$ is a regular point if
    $\bp$ is of Type~2.
  \end{itemize}
  The map $\bp\to\gamma(\bp)$ gives a bijection 
  between sequences of infinitely nearby points 
  not of Type~3, and $\Gamma$.
\end{theorem}
\begin{theorem}\label{Ttype3}
  Let $\bp=(p_j)_0^\infty$ be a sequence of infinitely
  nearby points of Type~3. 
  In the notation of the previous theorem, 
  the sequence $\gamma(\bp_n)\in\Gast$ converges weakly 
  in $\Gamma$ to an element $\gamma(\bp)\in\Gast$. 
  For $n$ large enough we have 
  $[\gamma(\bp),\gamma(\bp_{n+1})]\subset[\gamma(\bp),\gamma(\bp_n)]$. 
  In particular, the points $\gamma(\bp_n)$ define 
  the same tangent vector $\vv(\bp)$ at $\gamma(\bp)$.

  Moreover, the map $\bp\to\vv(\bp)$ gives a bijection between 
  infinitely nearby points of Type~3 and tangent vectors at 
  points in $\Gast$, \ie at branch points in $\Gamma$.
\end{theorem}
\begin{proposition}\label{P304}
  The sequence $\bp$ of infinitely nearby points associated to an 
  irreducible curve $C$ at the origin (see Example~\ref{ex-curve})
  is of Type~4. Conversely, every sequence of infinitely nearby
  points of Type~4 is associated to a unique irreducible curve.
\end{proposition}
In view of this proposition and Theorem~\ref{Tnot3} we see that
an irreducible curve $C$ naturally defines an end in the universal
dual graph $\Gamma$. We shall use this fact repeatedly in the sequel.

The following two results are consequences of Corollary~\ref{C309}
(see also Proposition~\ref{Ptangentgamma} and Remark~\ref{R305}):
\begin{corollary}\label{C311}
  If $\bp'=(p_j)_0^{n'}$ is a sequence of infinitely nearby points
  not of Type~3, $0\le n<n'\le\infty$ and $\bp=(p_j)_0^n$, 
  then the tangent vector at 
  $E:=\gamma(\bp)\in\Gast$ represented by $\gamma(\bp')\in\Gamma$ is
  given by the point $p_{n+1}$ on the exceptional component $E$.
 \end{corollary}
\begin{corollary}\label{C302}
  Consider any irreducible curve $C$ and view $C$ as an end in
  $\Gamma$. Let $\bp=(p_j)_0^\infty$ be the sequence of 
  infinitely nearby points (of Type~4) associated to $C$.
  Fix $n\ge0$, set $E_n:=\gamma((p_j)_0^n)\in\Gast$ and
  let $\pi_n$ be the composition of blowups at $p_0,\dots,p_n$.
  Then the tangent vector at $E_n$ represented by $C$
  is given by the intersection of $E_n$ with the strict transform 
  of $C$ under $\pi_n$.
\end{corollary}
%
%
\subsection{Proofs}\label{proof-infsg}
Let us now present the proofs of the results above, starting with
Theorem~\ref{Ttype0}.  If $\bp$ is of Type~0, then $\gamma(\bp)$ is
clearly an element in $\Gast$. Moreover, it is clear that different
sequences $\bp$ give rise to different exceptional components
$\gamma(\bp)$. To prove surjectivity, consider any $E\in\Gast$. By
Lemma~\ref{Lsup-inf} there is a minimal $\pi\in\fB$ for which
$E\in\Gast_\pi$. It is straightforward to verify that this $E$ is
associated to a sequence of infinitely nearby points. This completes
the proof of Theorem~\ref{Ttype0}.

\begin{figure}[ht]
  \begin{center}
    \includegraphics{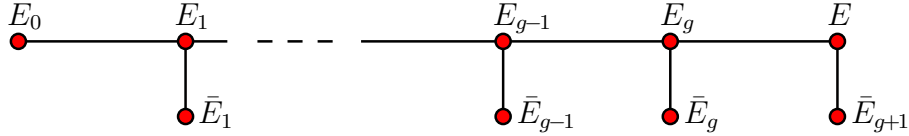}
  \end{center}
  \caption{The dual graph associated to a sequence of 
    infinitely nearby points of Type~0. 
    Only the branch points and ends are labeled.  
    Compare Figure~(5.7) in~\cite{spiv}.}\label{F3}
\end{figure}

\smallskip 
Before turning to the proof of the next theorem, we introduce some
notation which shall be used throughout the proofs below.  Consider a
finite sequence $(p_j)_0^n$ of infinitely nearby points.  Recall that
$p_0$ is the origin.  Therefore, $p_1$ is always free (assuming
$n\ge1$).  Define indices
\begin{equation}\label{e317}  
  0=n_0<\bn_1<n_1<\bn_2<n_2<\dots<\bn_g<n_g<\bn_{g+1}\le n
\end{equation}
as follows: $p_j$ is free for $n_i<j\le\bn_{i+1}$, $0\le i\le g$, and
satellite for $\bn_i<j\le n_i$, $1\le i\le g$. 
Let $\pi=\pi_\bp$ be the composition of blow-ups at
$p_0,\dots,p_n$ and let $E'_j$ be the strict
transform of the blow-up of $p_j$ in the total space 
$X_\pi$.\footnote{The notation differs marginally from 
  that of Section~\ref{sec-equivalence} where we would have written
  $E_j$ instead of $E'_j$.}
Identify $E'_j$ with its image in $\Gamma_\pi$. 
Write $E_i=E'_{n_i}$ for $0\le i\le g$
$\bE_i=E'_{\bn_i}$ for $1\le i\le g+1$.
Finally set $E=E'_n=\gamma(\bp)$.
The dual graph of $\pi_\bp$ then looks as in Figure~\ref{F3}. 

\smallskip
Let us now prove Theorems~\ref{Tnot3} and~\ref{Ttype3}.
Thus, pick an infinite sequence $\bp$ of infinitely nearby points 
(\ie of Type~1,~2,~3, or~4). We may define indices $n_i$, 
$0\le i\le g$ and $\bn_i$, $1\le i\le g+1$ as in~\eqref{e317},
allowing for the possibility that $g=\infty$.
Consider the truncations $\bp_n=(p_j)_0^n$ defined above.
Let us show that $\gamma(\bp_n)$ 
converges weakly to a point in $\Gamma$.

\begin{figure}[ht]
  \begin{center}
    \includegraphics{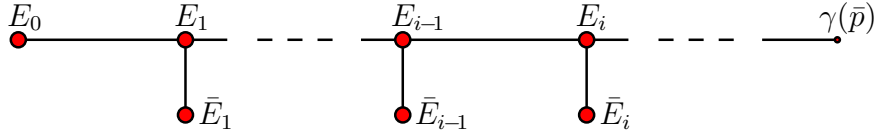}
    \caption{Infinitely nearby points of Type~1 corresponding to an 
     end in $\Gamma$ (of infinite multiplicity). 
     Compare Figure~(9.1) in~\cite{spiv}.}\label{F12}
  \end{center}
\end{figure}

If $\bp$ is of Type~1, then $g=\infty$ and the 
sequence $\gamma(\bp_{n_i})=E_i$ forms a strictly increasing 
sequence in $\Gast$. 
In fact, this sequence is unbounded in $\Gast$. 
This follows from the fact the sequence $\bp$ contains infinitely
many free blowups: see the proof of~(T3') in Section~\ref{S310}.
Thus $E_i:=\gamma(\bp_{n_i})$ converges to an end $\gamma(\bp)$
in $\Gamma$ as $i\to\infty$.  
As $\gamma(\bp_{n_i+k})\ge\gamma(\bp_{n_i})$ for all $k\ge0$, 
$\lim_{j\to\infty}\gamma(\bp_j)=\gamma(\bp)$.
See Figure~\ref{F12}. 

\begin{figure}[ht]
  \begin{center}
    \includegraphics{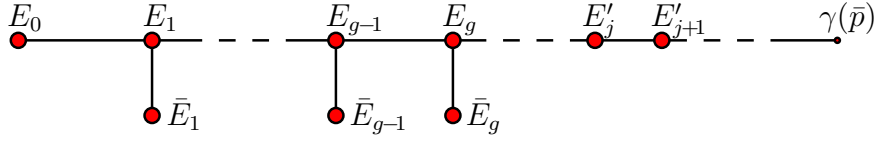}
    \caption{Infinitely nearby points of Type~4 corresponding to an 
     end in $\Gamma$ (of finite multiplicity). 
     Compare Figure~(9.6) in~\cite{spiv}.}\label{F13}
  \end{center}
\end{figure}

If $\bp$ is of Type~4, then $g<\infty$ and 
the sequence $\gamma(\bp_n)=E'_n$ for $n>n_g$ 
forms a strictly increasing sequence in $\Gast$. As before it
is unbounded, and therefore defines an end $E_\bp=\gamma(\bp)$
in $\Gamma$. See Figure~\ref{F13}. 

If $\bp$ is of Type~2 or ~3, then $g<\infty$ and 
for $j>\bn_{g+1}$, $p_{j+1}$ is a satellite point, 
intersection of two exceptional divisors, 
one of which is $E'_j$ and the other we denote by $F_j$. 
The segment $[E'_j,F_j]$ is then equal to either
$[E'_j,E'_{j-1}]$ or to $[E'_j, F_{j-1}]$, hence the sequence
$[E'_j,F_j]$ is decreasing. To conclude that $E_j$ converges in
$\Gamma$, we need to show that the intersection 
$\bigcap [E'_j,F_j]$ contains at most one point in $\Gast$.

Thus pick $E\in\Gast$ with $E\in [E'_j,F_j]$ for all $j$. 
If $E=E'_{j_0}$ for some $j_0$, then $\bp$ is of Type~3, and
$p_{j+1}$ is the intersection of $E$ with $E'_j$ for $j>j_0$.
Otherwise $E\ne E'_j$ for all $j$, in which case
we may find $j_0$ sufficiently large such that
$E$ can be obtained from $p_{j_0+1}=E'_{j_0}\cap F_{j_0}$ 
by a sequence of blow-ups of satellite points. 
Whence $E=E'_j$, or $E\not\in [E'_j,F_j]$ for some $j$. 
This is a contradiction, which shows that 
$\bigcap[E'_j,F_j]$ is empty for $\bp$ of Type~2, 
and is reduced to $E=E_{j_0}\in\Gast$ when $\bp$ is of Type~3.

\begin{figure}[ht]
  \begin{center}
    \includegraphics{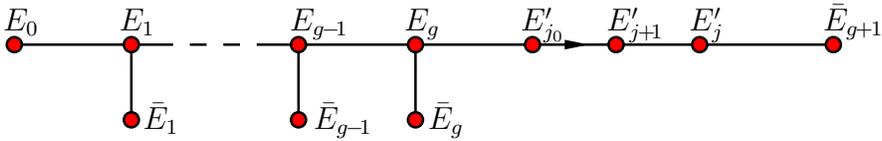}
    \caption{Infinitely nearby points of Type~3 corresponding 
      to a tangent vector at a branch point in $\Gamma$ 
      (here $E'_{j_0}$).
      Compare Figure~(9.4) in~\cite{spiv}.}\label{F15} 
  \end{center}
\end{figure}

\begin{figure}[ht]
  \begin{center}
    \includegraphics{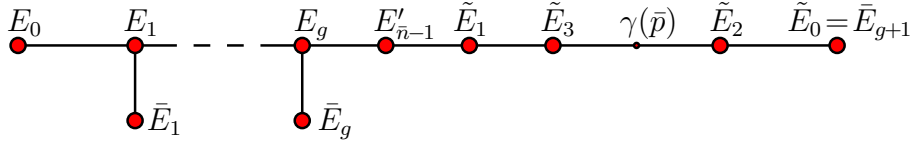}
    \caption{Infinitely nearby points of Type~2 corresponding to 
      a regular point in $\Gamma$. 
      Here $\tE_k=E_{\tn_k}$ for some increasing sequence 
      $(\tn_k)_0^\infty$ such that $\tE_{2k}$ ($\tE_{2k+1}$)
      decreases (increases) to $\gamma(\bp)$ as $k\to\infty$.
      Compare Figure~(9.2) in~\cite{spiv}.}\label{F17} 
  \end{center}
\end{figure}

This concludes the proof that $\bp\to\gamma(\bp)$ is
well-defined for infinite sequences $\bp$.
See Figures~\ref{F15} and~\ref{F17}.
Note that we also showed that $\gamma(\bp)$ belongs to
$\Gast$ if $\bp$ is of Type~3; is an end in $\Gamma$ if $\bp$
is of Type~1 or~4; and is regular point if $\bp$ is of Type~2.

To conclude the proof of Theorem~\ref{Tnot3} we 
only need to show that $\gamma$ is bijective from the set of infinitely
nearby points of Types~1,~2 and~4, to the set $\Gamma\setminus\Gast$.

First consider injectivity. 
Pick $\bp=(p_j)_0^\infty$ and $\bp'=(p'_j)_0^\infty$ 
of Types~1,~2 or~4 with $\gamma(\bp)=\gamma(\bp')$ and $\bp\ne\bp'$. 
Write $\gamma_n\=\gamma(\bp_n)$, 
$\gamma'_n\=\gamma(\bp'_n)$. 
Pick $n$ maximal such that $p_j=p'_j$ 
for all $j\le n$. Thus $p_{n+1}\ne p'_{n+1}$.

Suppose $p_{n+2}$ is a free point on
$\gamma_{n+1}$. Then $\gamma_{n+2}>\gamma_{n+1}$, 
and $[\gamma_{n+k},\gamma'_{n+1}]$ contains both 
$\gamma_{n+1}$ and $\gamma_n$ for all $k\ge2$. 
This is impossible, hence $p_{n+2}$ is a satellite point. 
Inductively, we see that
$p_{n+k+1}$ is necessarily the intersection of $\gamma_{n+k}$ and
$\gamma_n$ (interpreted as exceptional components), 
and $\gamma_{n+k}\to\gamma_n$. 
Thus $\bp$ is of Type~3, a contradiction. 
Hence $\bp\to\gamma(\bp)$ is injective on the set of
infinitely nearby points of Types~1,2 and~4.

Let us prove this map is also surjective. 
Pick $\gamma\in\Gamma\setminus\Gast$.
Let $p_0$ be the origin, and $E'_0$ the exceptional component obtained
by blowing up $p_0$.  The tangent vector $\vv_1$ at $E'_0$ represented
by $\gamma$ is associated to a unique point $p_1\in E'_0$ by
Proposition~\ref{Ptangentgamma}. We construct inductively an
infinite sequence $\bp=(p_j)_0^\infty$ 
of infinitely nearby points, so that $E'_j$
is the exceptional divisor of the blow-up at $p_j$, 
and $p_{j+1}\in E'_j$ is the point associated to the 
tangent vector represented by $\gamma$ at $E'_j$. 

When $\bp$ is of Type~1 or~4, $p_{j_k}$ is a free point 
for an increasing subsequence $j_k$. 
But $\gamma\ge E'_{j_k}$, and $\gamma(\bp)$ is an end, 
hence $\gamma=\gamma(\bp)$. 
When $\bp$ is of Type~2, there is a subsequence $(\tn_k)_0^\infty$
such that $\tE_{2k}:=E_{\tn_{2k}}$ ($\tE_{2k+1}:=E_{\tn_{2k+1}}$) 
decreases (increases) to $\gamma(\bp)$ as $k\to\infty$: 
see Figure~\ref{F17}.
It is not difficult to see that 
$\gamma\ge\tE_{2k+1}$ and $\gamma\wedge\tE_{2k}<\gamma(\bp)$ 
for all $k$. But $\gamma(\bp)$ is a regular
point in $\Gamma$, hence $\gamma=\gamma(\bp)$.

This concludes the proof of Theorem~\ref{Tnot3}. 

As for Theorem~\ref{Ttype3}, consider $\bp=(p_j)_0^\infty$ 
of Type 3 and pick $j_0$ as in the definition of Type~3.
Then 
$[\gamma(\bp),\gamma(\bp_{j+1})]\subset[\gamma(\bp),\gamma(\bp_j)]$
for all $j>j_0$. In fact, we saw above that the segments
$[\gamma(\bp),\gamma(\bp_j)]$ intersect only in $\gamma(\bp)$.
This proves the first few assertions in Theorem~\ref{Ttype3}:
in particular $\vv(\bp)$ is a well defined tangent vector
at $\gamma(\bp)$. Notice that the pair $(\gamma(\bp),\vv(\bp))$
is uniquely determined by the pair $((p_j)_0^{j_0},p_{j_0+1})$.
Using Proposition~\ref{Ptangentgamma} (see also Remark~\ref{R305}),
this easily implies that $\bp\to\vv(\bp)$ gives a bijection
between Type~3 infinitely nearby points and
tangent vectors at points in $\Gast$. 
The proof is complete.

\smallskip
Proposition~\ref{P304} is well-known so we shall only outline
a proof. First, to any irreducible curve we associate
a sequence of infinitely nearby points in terms of strict transforms
as in Example~\ref{ex-curve}. As the curve can be desingularized,
all but finitely many points in the sequence are free, \ie the sequence
is of Type~4. Conversely, consider a sequence $\bp=(p_j)_0^\infty$
of infinitely nearby points of Type~4. Let us show that $\bp$
is associated to a unique irreducible curve. By lifting the
situation, we may assume that all $p_j$ are free. Let us use
the notation $\tpi_j$ above.
Fix arbitrary preliminary coordinates $(z',w')$ at the origin.
As $p_1$ is free, we may write 
$\tpi_1(z'_0,w'_0)=(z'_0,z'_0(\theta_0+w'_0))$ 
for suitable coordinates $(z'_0.w'_0)$ at $p_1$ and $\theta_0\in\C^*$.
Inductively, we get coordinates $(z'_j,w'_j)$ at $p'_{j+1}$ 
such that $\tpi_j(z'_j,w'_j)=(z'_j,z'_j(\theta_j+w'_j))$
with $\theta_j\in\C^*$. 
Now define new coordinates $(z,w)$ at the origin and
$(z_j,w_j)$ at $p_{j+1}$ as follows:
$z=z'$, $z_j=z'_j$, $w=w'-\sum_0^\infty\theta_k(z')^{k+1}$
and $w_j=w'_j-\sum_1^\infty\theta_{k+j}(z'_j)^k$.
In these coordinates, $\tpi_j(z_j,w_j)=(z_j,z_jw_j)$ for $j\ge0$.
It is then straightforward to verify that $\{w=0\}$ is the unique
irreducible curve associated to the sequence $\bp$.

For the record we notice that all the maps $\tpi_j$ are monomial
in the coordinates that we constructed.
%
%
%
%
\section{Parameterization and multiplicities}
After having described the connection between $\Gamma$ and 
sequences of infinitely nearby points, we proceed to
equip $\Gamma$ with a parameterization and multiplicity function.
%
%
\subsection{Farey weights and parameters}\label{sec-farey}
\begin{figure}[ht]
  \begin{center}
    \includegraphics[width=\textwidth]{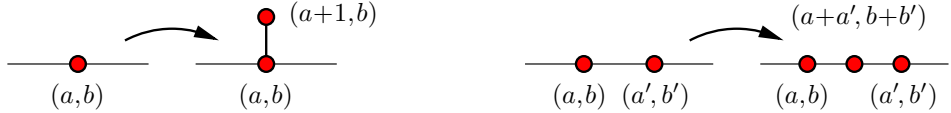}
  \end{center}
  \caption{Farey weights under elementary modifications. 
    The first kind of modification is illustrated to the left, 
    the second kind to the right.
    Compare Figure~\ref{F10}.}\label{F2}
\end{figure}

Let us associate to each element of $\Gast$ a vector
$(a,b)\in(\N^*)^2$, called its \emph{Farey
weight}\index{Farey!weight}.  It is defined inductively through the
following combinatorial procedure (see~\cite{hubbard} for the toric
case).  If $\pi$ is a single blowup of the origin, then $\Gamma_\pi$
is a single point $E_0$ whose weight is defined to be $(2,1)$.
Otherwise we may write $\pi=\pi'\circ\tpi$, where $\tpi$ is the blowup
of a point $p$ on the exceptional divisor $(\pi')^{-1}(0)$.  The Farey
weights of the vertices in $\Gamma_\pi$ that are strict transforms of
vertices in $\Gamma_{\pi'}$ inherit their weights from the latter
graph. The only other vertex in $\Gamma_\pi$ is the exceptional
divisor $E_p$ obtained by blowing up $p$.  The weight of $E_p$ is
determined as follows.  When $p$ is free, \ie $p\in E$ for a unique
exceptional divisor $E\subset(\pi')^{-1}(0)$, then the weight of $E_p$
is defined to be $(a+1,b)$, where $(a,b)$ is the weight of $E$.  When
$p$ is a satellite point, it is the intersection of two components $E$
and $E'$ whose weights are $(a,b)$ and $(a',b')$, respectively.  The
weight of $E_p$ is then $(a+a',b+b')$.  See Figure~\ref{F2} and
compare with Figure~\ref{F10}.

The \emph{Farey parameter}\index{Farey!parameter} of a point in
$\Gast$ with weight $(a,b)$ is defined to be the rational number
$A=a/b$. We shall later see that $A$ defines a parameterization of
$\Gast$. Although, this can be proved directly 
(in roughly the same way that we verified 
condition~(T2) in Section~\ref{S310}),
we shall only prove
\begin{lemma}\label{Lfareyincr}
  The Farey parameter $A$ is a strictly increasing function on $\Gast$.
\end{lemma}
\begin{proof}
  It is equivalent to prove that the restriction of $A$ to the set of
  vertices $\Gast_\pi$ for any $\pi\in\fB$ is strictly increasing.  We
  proceed by induction on the number of point blow-up necessary to
  decompose $\pi$.  When this number equals $1$, $\pi$ is the blow-up
  of the origin, $\Gast_\pi$ is reduced to one point, and the claim is
  obvious. Otherwise, pick $\pi\in\fB$, and consider elements
  $E,E'\in\Gast$ that are adjacent vertices in $\Gamma_\pi$.  Assume
  $E<E'$. There are two cases.  In the first case, $E'$ is obtained from
  $E$ by blowing up a free point and $a'/b'=(a+1)/b>a/b$.  In the second
  case, $E'$ is obtained by blowing up the intersection point between
  two irreducible components $E$ and $E''$. Then $E<E'<E''$ so the
  inductive assumption gives $a/b<a''/b''$, from which it is elementary
  to see that $a'/b'=(a+a'')/(b+b'')\in\,]a/b,a''/b''[$.  This concludes
  the proof.
\end{proof}
%
%
\subsection{Multiplicities}\label{S403}
We can also use the Farey weights to define \emph{multiplicities} in the
\index{multiplicity!in $\Gamma$}
\index{multiplicity!in $\Gast$}
universal dual graph. Namely, if $E\in\Gamma^o$, then we set
\begin{equation*}
  m(E):=\min\{b(F)\ ;\ F\in\Gast, F\ge E\},
\end{equation*}
\index{$m(E)$ (multiplicity)}
\index{multiplicity!of an exceptional component}
where $(a(F),b(F))$ denotes the Farey weight of $F$. 

\begin{figure}[ht]
  \begin{center}
    \includegraphics{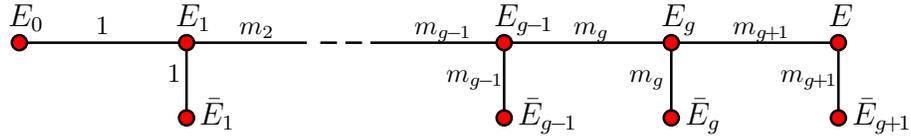}
  \end{center}
  \caption{Multiplicities.}\label{F9} 
\end{figure}

By definition, $m$ is integer-valued and increasing on $\Gamma^o$.
It hence has a minimal extension to a function on $\Gamma$
with values in $\overline{\N}$. 
Let us describe its basic structure.
The key to such a description is given by
\begin{lemma}\label{L415}
  Let $E\in\Gast$ and suppose $F$ is obtained from $E$ by blowing up a
  \emph{free} point. Then $m(F)=b(F)=b(E)$. Moreover, the multiplicity
  is constant equal to $b(E)$ in the segment $]E,F]$.
\end{lemma}
\begin{proof}
  Let $(a,b)$ be the Farey weight of $E$.  As $F$ is obtained by blowing
  a free point $p$ on $E$, its Farey weight is by definition
  $(a+1,b)$. Whence $b(E)=b(F)$.  We always have $m(F)\le b(F)$.  To
  prove the converse inequality, observe that if $F'\ge F$, then the
  sequence of infinitely nearby points associated to $F'$ starts with
  the one of $F$ followed by the blowup at a free point on $F$.  An
  easy induction shows that the Farey weight $(a',b')$ of $F'$ satisfies
  $b'\ge b(F)$. Thus $m(F)\ge b(F)$. This shows $m(F)=b(F)=b(E)$.

  Any $F'$ representing the same tangent vector as $F$ at $E$
  is obtained by a sequence of blowups starting with $p$,
  hence $b(F')\ge b=b(E)$. 
  In particular the multiplicity of any element in the
  segment $]E,F]$ is at least $b(E)$. On the other hand,
  this multiplicity cannot exceed $m(F)=b(E)$. 
  Thus $m\equiv b(E)$ on $]E,F]$.
\end{proof}
A direct consequence of this lemma is that, when $E\in\Gast$ and its
associated dual graph is as in Figure~\ref{F3},
$m(E)=b(E_g)=b(\bE_{g+1})$.  Define $m_i=b(E_{i-1})$ 
for $1\le i\le g+1$. Then $1=m_1<m_2<\dots<m_g<m_{g+1}$. 
By Lemma~\ref{L415} the
multiplicity equals $m_i$ on the segment $]E_i,\bE_{i+1}]$.
This is illustrated in Figure~\ref{F9}.
In analogy with the situation on the valuative tree, we shall refer to
$E_1$,\dots,$E_g$ as the 
\emph{approximating sequence}
\index{approximating sequence} 
of $E$.

Given $E\in\Gast$ and a tangent vector $\vv$ at $E$,
define the multiplicity\index{multiplicity!in $\Gast$} of $\vv$ by
$m(\vv)=m(E)$ if $\vv$ is represented by the root $E_0$, and
\begin{equation*}
  m(\vv)=\min\{m(F)\ ;\ F\ \text{represents $\vv$}\}
\end{equation*}
\index{multiplicity!of a tangent vector}
\index{$m(\vv)$ (multiplicity)}
otherwise. 

Let $E\in \Gast$ and $\pi\in\fB$ be minimal such that
$E\subset\pi^{-1}(0)$.  Its dual graph is given by Figure~\ref{F9}.
It is thus clear that the two tangent vectors represented by $E_0$ and
$\bE_{g+1}$ (if $\bE_{g+1}\ne E$)
satisfy $m(\vv) = m(E)$. Note that these tangent vectors
correspond exactly to the intersection points of $E$ with the other
components of the exceptional divisor.  Fix any other tangent vector
$\vv$, and $F\in\Gast$ representing $\vv$. Then $F$ is obtained by
first blowing up a free point $p$ on $E$. Thus $b(F)\ge b(E)$, \ie
$m(\vv)\ge b(E)$. On the other hand, Lemma~\ref{Ptangentgamma} shows
that $\vv$ is represented by the exceptional component $F$ obtained by
blowing up $p$, and Lemma~\ref{L415} gives $b(F) = b(E)$. Whence
$m(\vv)=b(F)$, and we obtain the following two corollaries.
\begin{corollary}\label{C306}
  We have $m(\vv)=b(E)$ for all but at most two tangent
  vectors $\vv$ at $E$. The exceptional tangent vectors
  have multiplicity $m(E)$.
\end{corollary}
We shall call $b(E)$ the \emph{generic multiplicity}
\index{multiplicity!generic} 
\index{$b(E)$ (generic multiplicity)} 
of $E$.
\begin{corollary}\label{C305}
  The multiplicity of a tangent vector $\vv$ at $E\in\Gast$ is equal
  to the generic multiplicity iff there exists $\pi\in\fB$ such that
  $E\in\fB$ and the point on $E$ defined by $\vv$ (see
  Proposition~\ref{Ptangentgamma}) is free.
\end{corollary}
We shall call such a tangent vector a \emph{generic tangent
vector}\index{tangent vector!generic}.

Note that when $m(\vv)=b(E)$, $\pi$ may be chosen to be the
composition of blow-up at the sequence of infinitely nearby points
associated to $E$.
Moreover, when blowing up at a free point, the generic multiplicity
does not increase. Therefore we have
\begin{corollary}\label{C432}
  An element of $\Gamma$ has infinite multiplicity iff it
  is an end associated to a sequence of infinitely nearby points 
  of Type~1.
\end{corollary}
Recall that any end of $\Gamma$ corresponding to a
sequence of infinitely nearby points of Type~4 can
be viewed as an irreducible curve $C$. 
Write $m_\Gamma(C)$ for the multiplicity of $C$ as
an element of $\Gamma$. (We shall later show that $m_\Gamma(C)$
coincides with the usual multiplicity $m(C)$ of $C$.)
\begin{corollary}\label{C307}
  If $C$ is an irreducible curve, $\pi\in\Gast$ and
  the strict transform $C'$ of $C$ by $\pi$ intersects 
  $\pi^{-1}(0)$ in a free point $p\in E$, then
  $m_\Gamma(C)\ge b(E)$ with equality iff
  $C'$ is smooth and transverse to $E$ at $p$.
\end{corollary}
\begin{proof}
  Let $(q_j)_0^\infty$ be the sequence of infinitely nearby points
  associated to $C'$ (note that $q_0=p$) and let $F_j\in\Gast$ be the
  exceptional divisor obtained by blowing up $q_j$.  Then
  $m_\Gamma(C)=\lim b(F_j)$. By induction, $b(F_j)\ge b(E)$ with
  equality iff $(q_k)_0^j$ are all free. Thus $m_\Gamma(C)\ge b(E)$
  with equality iff $q_j$ is free for all $j$.  The latter condition
  is equivalent to $C'$ being smooth and transverse to $E$ at $p$.
\end{proof}
\begin{corollary}\label{C304}
  For any $E\in\Gamma^o$ we have 
  $m(E)=\min\{m_\Gamma(C)\ ;\ C>E\}$.
\end{corollary}
\begin{proof}
  The inequality $m(E)\le\min\{m_\Gamma(C)\ ;\ C>E\}$
  is obvious. To prove the other inequality it suffices
  to show that any $E\in\Gast$ can be dominated by
  a curve $C$ with $m_\Gamma(C)=b(E)$. 
  For this, pick $\pi\in\fB$ such that $E\in\Gast_\pi$,
  pick any free point $p$, and let $C'$ be a smooth curve
  at $p$, transverse to $E$.
  By Corollary~\ref{C307} the choice $C=\pi(C')$ works.
\end{proof}
%
%
%
%
\section{The isomorphism}
We are now ready to formulate the main result of this chapter.
It states that the universal dual graph, whose construction is
purely geometric, is isomorphic---in the strongest possible sense---to 
the  valuative tree.

The definition of the isomorphism goes as follows. 
A point $E\in\Gast$
is an exceptional component of some $\pi\in\fB$. 
It hence defines a normalized divisorial valuation $\nu_E\in\cV$:
$\nu_E(\phi)$ is proportional to the order of vanishing of $\pi^*\phi$
along $E$. In other words, $\nu_E=b^{-1}\pi_*\div_E$ for some
constant $b>0$ (we will later see that $b$ is indeed the 
generic multiplicity of $E$). This defines a mapping
$\Phi:\Gast\to\cVdiv$.
\begin{theorem}\label{thm-universal}
  The map $\Phi:\Gast\to\cVdiv$ extends uniquely to 
  an isomorphism of parameterized trees
  $\Phi:(\Gamma,A)\to(\cV,A)$.
  Here $A$ denotes the Farey parameter on $\Gamma$ and
  thinness on $\cV$. Further, $\Phi$ preserves multiplicity.
\end{theorem}
Before starting the proof of Theorem~\ref{thm-universal}, let us 
state some immediate consequences and remarks. We refer to 
Section~\ref{Sappliuniv} below for more applications of the theorem.
The first consequence was in fact already implicitly used above.
\begin{corollary}\label{C406}
  The Farey parameter defines a parameterization of the rooted,
  nonmetric tree $\Gamma$. 
  The points in $\Gast$ are exactly the points in $\Gamma^o$
  having rational Farey parameter.
\end{corollary}
\begin{proof}
  The Farey parameter is indeed a parameterization, 
  as thinness gives a parameterization of $\cV$.  
  As $\cVdiv$ is the set of valuations in
  $\cVqm$ having rational thinness, the points in $\Gast$ must be
  exactly the ones in $\Gamma^o$ having rational Farey parameter.
\end{proof}
Our definition of the isomorphism $\Phi:\Gamma\to\cV$ is quite
indirect. In particular we do not explicitly define the value of
$\Phi$ on elements of $\Gamma\setminus\Gast$.  In fact, we have
\begin{corollary}\label{Cencoding}
  Let $\bp$ be a sequence of infinitely nearby points, 
  $\nu(\bp)$ be its associated Krull valuations 
  as in Section~\ref{versus}, and $\gamma(\bp)$
  be its associated point in $\Gamma$ as in 
  Section~\ref{S404}.
  
  When $\bp$ is not of Type~3, $\Phi(\gamma(\bp))=\nu(\bp)$.
  When $\bp$ is of Type~3, $\nu(\bp)$ is an exceptional curve
  valuation, and $\Phi(\gamma(\bp))=\nu'(\bp)$ where $\nu'(\bp)$ is the
  divisorial valuation associated to $\nu(\bp)$.
\end{corollary}
A direct consequence is:
\begin{corollary}\label{C303}
  For any irreducible curve $C$, viewed as an element of $\Gamma$,
  we have $\Phi(C)=\nu_C$, where $\nu_C$ is the curve valuation
  associated to $C$.
\end{corollary}
\begin{proof}[Proof of Corollary~\ref{Cencoding}]
  If $\bp$ is finite (\ie of Type~0), then the statement follows
  immediately from the definitions.
  Thus suppose $\bp=(p_j)_0^\infty$ is infinite and
  consider the truncations $\bp_n=(p_j)_0^n$.
  By the definition of $\gamma(\bp)$,
  $\gamma(\bp_n)\in\Gast$ tends to $\gamma(\bp)$ 
  in the weak tree topology of $\Gamma$. 
  On the other hand, it is also clear that $\nu(\bp_n)\in\cVdiv$ 
  converges to $\nu(\bp)$ in the
  Zariski topology on the set $\cV_K$ of Krull valuations (see
  Section~\ref{sec-zar}).
  Denote by $P(\nu(\bp))$ the valuation in $\cV$ corresponding to 
  $\nu(\bp)$.  We refer to the discussion in Section~\ref{versus}. 
  When $\bp$ is of Type~3, $P(\nu(\bp))$ is the
  divisorial valuation associated to $\nu(\bp)$.  
  By Theorem~\ref{weak-zariski}, $P(\nu(\bp_n))$ converges weakly to 
  $P(\nu(\bp))$, hence in the weak tree topology of $\cV$. 
  As $\Phi$ is a tree isomorphism, $P(\nu(\bp))=\Phi(\gamma(\bp))$.
\end{proof}
%
%
%
%
\section{Proof of the isomorphism}
To prove Theorem~\ref{thm-universal}, we shall consider
$\Phi$ as a mapping of $\Gast$ into $\cVdiv$ and proceed
in four steps: first, $\Phi$ is a bijection onto $\cVdiv$; 
second, the image of the Farey parameter under $\Phi$ is 
equal to the thinness;
third, $\Phi$ and its inverse are order preserving; 
and fourth, $\Phi$ preserves multiplicity.

This will show that $\Phi:\Gast\to\cVdiv$ is an isomorphism
of rooted, nonmetric $\Q$-trees. 
As thinness $A$ is a parameterization of $\cVdiv$, the Farey parameter
$A$ is then a parameterization of $\Gast$ and
$\Phi:(\Gast,A)\to(\cVdiv,A)$ is an isomorphism of parameterized 
$\Q$-trees. It follows immediately that $\Phi$ extends
to an isomorphism 
$\Phi:(\Gamma,A)\to(\cV,A)$ of parameterized $\R$-trees.
Finally, as $m$ is lower-semicontinuous on $\Gamma$ and $\cV$,
and $m=m\circ\Phi$ on $\Gast$, we see that
$m=m\circ\Phi$ on $\Gamma$, so that $m$ preserves
multiplicity.
%
%
\subsection{Step 1: $\Phi : \Gast \to \cVdiv$ is bijective} 
That $\Phi$ is surjective is automatic: any divisorial
valuation $\nu$ is equivalent to $\pi_*\div_E$ for some $\pi\in\fB$
and some exceptional component $E$ (see Proposition~\ref{Pchardiv}). 

Let us give two arguments to prove that $\Phi$ is injective.  
Pick $E,E'\in\Gast$ and $\pi\in\fB$ such that $E,E'\in\Gast_\pi$.  
Suppose $E\ne E'$. 
By construction the center of the valuation $\Phi(E)$ in
the total space $X_\pi$ is irreducible and contains $E$. It is hence
equal to $E$. Similarly, the center of $\Phi(E')$ is equal to
$E'$ so that $\Phi(E)\ne\Phi(E')$. We may also use a more geometric
argument.  Embed $X_\pi$ in some projective space
$\P^k$, $k\ge2$. The curves $E,E'$ are algebraic so we may
find a polynomial $P=P(z_0,\dots,z_k)$ homogeneous of degree $d$,
such that $E\subset\{P=0\}$, $E'\not\subset\{P=0\}$. 
Choose a coordinate axis $z_i$ such that 
$E\cup E'\not\subset\{z_i=0\}$. 
Then the rational function $h=P/z_i^d$ satisfies $\div_E(h)>0$
whereas $\div_{E'}(h)=0$. Whence $\Phi(E)\ne\Phi(E')$.
%
%
\subsection{Step 2: $A\circ\Phi= A$}
We first need a preliminary result that gives geometric 
interpretations of the components $a$ and $b$ of the Farey weight.
\begin{proposition}\label{Pgeo-inter}
  Let $\pi\in \fB$, and $E$ be an exceptional component with Farey
  weight $(a,b)$. Then
  \begin{align*}
    \div_E(\pi^*\fm)&\=\min_{\phi\in\fm} \div_E(\pi^*\phi)=b\\
    \div_E(J\pi)&=a-1.
  \end{align*}
  Here $\div_E$ denotes the order of vanishing along $E$, and $J\pi$ the
  Jacobian determinant of $\pi$. Note that the first equation shows that
  \begin{equation}\label{eq-univ1}
    \pi_*\div_E=b\,\nu_E.
  \end{equation}
\end{proposition}
This result allows us to give a local normal form for
the contraction map; the following lemma is a key ingredient in
the rest of the proof of Theorem~\ref{thm-universal}.
\begin{lemma}\label{L401}
  Let $E\in\Gast$ and $\pi\in\fB$ be as in the previous 
  proposition, and pick any free point $p\in E$.
  Then one can find local coordinates $(z,w)$ at $p$, 
  local coordinates $(x,y)$ at the origin, such that
  $E=\{z=0\}$ and such that the contraction map is given by 
  \begin{equation}\label{e403}
    \pi(z,w)=(z^b,z^{a-b}w+z^bh(z)),
  \end{equation}
  for some regular function $h$ with $h(0)\ne0$. 
  Moreover, given any smooth curve $V$ at $p$ intersecting
  $E$ transversely we may choose the coordinates such that
  $V=\{w=0\}$.
\end{lemma}
We are now able to prove the relation between thinness and the Farey
weight. Pick $\nu\in\cVdiv$, choose $\pi\in\fB$,
$E\subset\pi^{-1}(0)$ such that $\nu$ is equivalent to
$\pi_*\div_E$. Let $(a,b)$ be the Farey weight of $E$. 
We want to show $A(\nu)=a/b$. Pick a free point $p\in E$. 
By Proposition~\ref{Pgeo-inter}, $\pi_*\div_E=b\,\nu$.  
By Lemma~\ref{L401} we can choose coordinates $(z,w)$ at $p$ such that
$\pi(z,w)=(z^b,z^{a-b}w+z^bh(z))$, where $h(0)\ne0$.  Notice that
$\nu(x)=b^{-1}\pi_*\div_E(x)=1$. Recall from Chapter~\ref{sec-puis},
that we wrote $k = \C(x)$, and let $\hk$ be the field of Puiseux
series over $\C$ that is the algebraic closure of $k$.  Then the
valuation $\nu$ extends to a valuation on $k[y]$ of the form
$\hnu=\val[\hphi;\hbeta]$ for a Puiseux series $\hphi\in\hk$, and
where $1+\hbeta\in\Q_+$ is the thinness of $\nu$.  See
Theorem~\ref{T602}.  On the other hand, $\hbeta$ is also the maximum
of $\nu(y-\hpsi)$ when $\hpsi$ ranges over all Puiseux series.  But
\begin{equation*}
  \hnu(y-\hpsi)
  =b^{-1}(\pi_*\div_E)(y-\hpsi)
  =b^{-1}\div_E\left(wz^{a-b}+z^bh(z)-\hpsi(z^b)\right)
  \le\frac{a-b}{b},
\end{equation*} 
with equality for suitably chosen $\hpsi$. 
Hence $\hbeta=a/b-1$ and $A(\nu)=1+\hbeta=a/b$.

\begin{proof}[Proof of Proposition~\ref{Pgeo-inter}]
  First note that when $\pi$ is the blow-up of the origin, the
  proposition is clear as we have $E=\pi^{-1}(0)$, 
  $\div_E\pi^*\fm=1=b(E)$, and $\div_E(J\pi)=1=a(E)-1$. 
  We then proceed by induction as follows. 
  Pick $\pi\in\fB$, and suppose the proposition
  has been proved for all exceptional components
  $E\subset\pi^{-1}(0)$. We shall then prove it for the exceptional
  component obtained by blowing up an arbitrary point $p\in\pi^{-1}(0)$.
  
  Fix coordinates $(z,w)$ centered at $p$ such that the exceptional
  divisor is given by $E=\{z=0\}$ if $p$ is free, and $E=\{z=0\}$ ,
  $E'=\{w=0\}$ if $p$ is satellite. We let $(a,b)$, and $(a',b')$ be the
  Farey weights of $E$ and $E'$ respectively. The ideal $\pi^*\fm$ in
  the ring of regular functions at $p$ is generated by the pull-back
  of a generic smooth element $\phi\in\fm$. It is hence principal. By
  the inductive step, it is generated by $z^b$ when $p$ is free , and
  $z^bw^{b'}$ when $p$ is satellite. Similarly, the critical set of
  $\pi$ coincides with $\pi^{-1}(0)$, so that the determinant of the
  Jacobian of $\pi$ equals $J\pi = z^{a-1}\times\mathrm{units}$ when $p$
  is free, and $z^{a-1}w^{a'-1}\times\mathrm{units}$ when $p$ is
  satellite.

  Write $\tilde{\pi}$ for the blow-up of $p$, and 
  $F=\tilde{\pi}^{-1}(p)$. To compute the order of vanishing of
  $(\pi\circ\tilde{\pi})^*\fm$, and $J(\pi\circ\tilde{\pi})$ along
  $F$, we note that at a free point $p_1\in F$, we may choose coordinates
  $(z_1,w_1)$ such that $\tilde{\pi}(z_1,w_1)=(z_1,z_1(1+w_1))$. 
  Whence $(\pi\circ\tilde{\pi})^*\fm = \tilde{\pi}^*\pi^*\fm$ is a 
  principal ideal generated by
  \begin{equation*}
    \begin{cases}
      \tpi^*z^b=z_1^{b} & \text{when $p$ is free};\\
      \tpi^*z^bw^{b'}=z_1^{b+b'} &\text{when $p$ is satellite}. 
    \end{cases}
  \end{equation*}
  On the other hand, by the chain rule formula,
  \begin{multline*}
    J\pi\circ\tilde{\pi}
    =J\pi\circ\tilde{\pi}\ \times J\tilde{\pi}
    =\\
    \begin{cases}
      z_1^{a-1}\times\mathrm{units}\times z_1
      =z_1^a\times\mathrm{units}
      &\text{when $p$ is free};\\
      z_1^{a-1+a'-1}\times\mathrm{units}\times z_1
      =z_1^{a+a'-1}\times\mathrm{units}
      &\text{when $p$ is satellite}. 
    \end{cases}
  \end{multline*}
  The Farey parameter $(a_F,b_F)$ of $F$ equals $(a+1,b)$ when $p$ is
  free, and $(a+a',b+b')$ when $p$ is satellite. Thus
  $\div_F(\tilde{\pi}^*\pi^*\fm) = b_F$, and $\div_F 
  J(\pi\circ\tilde{\pi}) = a_F-1$. This concludes the proof.
\end{proof}
\begin{proof}[Proof of Lemma~\ref{L401}]
  Pick local coordinates $(x,y)$ at the origin such that
  the strict transforms of $\{x=0\}$ and $\{y=0\}$ 
  do not pass through $p$. Also pick local coordinates $(z,w)$
  at $p$ such that $E=\{z=0\}$ and $V=\{w=0\}$.
  By Proposition~\ref{Pgeo-inter}, 
  $\pi^*x$ and $\pi^*y$ vanish to order $b$ along $E$.
  After multiplying $z$ by a unit, we thus have
  $\pi(z,w)=(z^b,z^b\xi(z,w))$ for some unit $\xi$.
  Similarly, by Proposition~\ref{Pgeo-inter}, 
  the Jacobian determinant of $\pi$
  is given by $z^{a-1}\rho(z,w)$, where $\rho(0,0)\ne0$.
  Thus $b^{-1}z^{b-1}z^b\frac{\partial\xi}{\partial w}=z^{a-1}\rho(z,w)$,
  which gives $\xi=h(z)+dwz^{a-2b}\xi_1(z,w)$, where 
  $\xi_1+w\frac{\partial\xi_1}{\partial w}=\rho(z,w)$.
  Since $\xi(0,0)\ne0$ and $\rho(0,0)\ne0$ we have
  $h(0)\ne0$ and $\xi_1(0,0)\ne0$. 
  After multiplying $w$ by $d\xi_1(z,w)$ we arrive at~\eqref{e403}.
\end{proof}
%
%
\subsection{Step 3: $\Phi$ and $\Phi^{-1}$ are order preserving}
This is the key (and hardest) part of the argument. 
We start by stating
\begin{lemma}\label{Lkeylemma1}
  Pick $\pi\in\fB$ and $p\in\pi^{-1}(0)$.
  Let $F\in\Gast$ be the exceptional divisor of the blow-up of $p$. 
  Then:
  \begin{itemize}
  \item
    when $p$ is free, it belongs to a unique component
    $E\subset\pi^{-1}(0)$, and we have $\nu_F>\nu_E$;
  \item
    when $p$ is satellite, it lies at the intersection 
    of two components $E$ and $E'$ and $\nu_F\in[\nu_E,\nu_{E'}]$.
  \end{itemize}
\end{lemma}
In order to prove that both $\Phi$ and its inverse are order
preserving, it suffices to show that the map 
$\Phi:\Gast_\pi\to\cVdiv$ is a 
\emph{tree embedding} for any $\pi\in\fB$,
in the sense that $\Phi(E)<\Phi(E')$ iff $E<E'$.
We then proceed by induction on the number of blow-ups necessary to 
decompose $\pi$, \ie the cardinality of $\Gast_\pi$. 
When this number equals one, $\pi$ is the blow-up of the origin, and
$\Gast_\pi$ is reduced to one point, so that the statement is obvious.

Now suppose $\Phi:\Gast_\pi\to\cVdiv$ is a tree embedding
and pick $p\in\pi^{-1}(0)$. 
Let $\tpi$ denote the blow-up of $p$, and set $F\=\tpi^{-1}(p)$. 
When $p$ is a satellite point, it is the
intersection of two divisors $E$ and $E'$. 
The graph $\Gast_{\pi\circ\tpi}$ is obtained from 
$\Gast_\pi$ by adding one vertex in the segment $[E,E']$. 
By Lemma~\ref{Lkeylemma1},
$\nu_F\in[\nu_E,\nu_{E'}]$. This segment does not contain any other
valuation $\nu_{E''}$ for any exceptional components 
$E''\in\Gast_\pi$ by the inductive assumption.
Thus $\Phi:\Gast_{\pi\circ\tpi}\to\cVdiv$ is also a
tree embedding in this case.

When $p$ is free, it belongs to a unique component
$E\in\Gast_\pi$. The graph $\Gast_{\pi\circ\tpi}$ is
obtained from $\Gast_\pi$ by attaching a new vertex $F$ to $E$. 
By Lemma~\ref{Lkeylemma1}, $\nu_F>\nu_E$. 
This shows that $\Phi$ is order preserving. 
To conclude the induction step, we also need 
to show that $\nu_F$ does not define the same tangent vector as
$\nu_{E'}$ at $\nu_E$ for some other $E'\in\Gast_\pi$. 
Equivalently, we have to show that 
$\nu_F\wedge\nu_{E'}\le\nu_E$ for all $E'\in\Gast_\pi$. 
We may assume $\nu_{E'}\ge\nu_E$. 
As $\Phi:\Gast_\pi\to\cVdiv$ is a tree embedding, it is also
sufficient to consider the case when $E'$ intersects $E$. 
Pick an irreducible curve $C=\{\phi=0\}$, $\phi\in\fm$,
at the origin, whose strict transform $C'$
by $\pi$ is smooth and contains 
$p'\=E\cap E'$. 
Let $b$ and $b'$ be the generic multiplicities of $E$ and $E'$,
respectively; and let $F$ and $F'$ 
be the exceptional divisors of the blow-ups at $p$ and $p'$.
Introduce $\mu_p$ and $\mu_{p'}$, the multiplicity valuations 
at $p$ and $p'$, respectively.
From Proposition~\ref{Pgeo-inter} we have $\pi_*\div_E=b\,\nu_E$
and $\pi_*\div_{E'}=b'\,\nu_{E'}$, 
so that $\pi_*\mu_p=b\,\nu_F$ and $\pi_*\mu_{p'}=(b+b')\nu_{F'}$. 
Finally, let $\phi'$ be a defining function for $C'$:
this is a regular function vanishing at $p'$. 
We then have:
\begin{equation*}
  \nu_F(\phi)
  =b^{-1}(\pi_*\mu_p)(\phi)
  =b^{-1}\mu_p(\pi^*\phi)
  =b^{-1}\div_E(\pi^*\phi)
  =\nu_E(\phi).
\end{equation*}
On the other hand,
\begin{multline*}
  \nu_{F'}(\phi)
  =(b+b')^{-1}(\pi_*\mu_{p'})(\phi)
  =(b+b')^{-1}\mu_{p'}(\pi^*\phi)
  =\\
  =(b+b')^{-1}
  \left(\div_E(\pi^*\phi)+\div_{E'}(\pi^*\phi)+\mu_{p'}(\phi')\right)
  >\\  
  >(b+b')^{-1}\left(b\nu_E(\phi)+b'\nu_{E'}(\phi)\right)
  >\nu_E(\phi).
\end{multline*}
We conclude that $\nu_F(\phi)=\nu_E(\phi)<\nu_{F'}(\phi)$. 
But by Lemma~\ref{Lkeylemma1}, 
$\nu_{F'}\in[\nu_E,\nu_{E'}]$, so that $\nu_{F'}$
defines the same tangent vector as $\nu_{E'}$ at $\nu_E$. 
Thus $\nu_F\wedge\nu_{F'}=\nu_F\wedge\nu_{E'}=\nu_E$,
which completes the proof of Step 3.

Let us now turn to the proof of Lemma~\ref{Lkeylemma1}.
We rely on the following well-known result, which may
be viewed as the second basic ingredient in the proof of
the isomorphism, the first one being the normal form
in Lemma~\ref{L401}.
\begin{lemma}\label{Lkeylemma2}
  Fix a smooth curve $D$ at the origin and consider a finite sequence
  $(p_j)_0^n$ of infinitely nearby points above the origin such that
  $p_1,\dots,p_{n_0}$ are all free and $p_{n_0+1},\dots,p_n$ are all
  satellite.  Here $0\le n_0\le n$.  Let $\pi$ be the composition of
  blowups at the points $p_0,\dots,p_{n-1}$.  Then there exist
  coordinates $(x,y)$ at the origin and $(z,w)$ at $p_n$ such that
  $D=\{x=0\}$, $\{zw=0\}\subset\pi^{-1}(D)$ and such that $\pi$ is a
  monomial map in these coordinates.
\end{lemma}
Here, ``free/satellite point'' is to be understood as
``regular/singular point on the total transform of $D$''.
We postpone the proof of this lemma to the end of this section, and
conclude the proof of Lemma~\ref{Lkeylemma1}.
\begin{proof}[Proof of Lemma~\ref{Lkeylemma1}]
  When $p$ is free, it belongs to 
  a unique component $E$. By Lemma~\ref{L401}, there exist
  coordinates $(z,w)$ at $p$, $(x,y)$ at the origin,
  such that $\pi(z,w)=(z^b,z^{a-b}w+z^bh(z))$, for some regular
  function $h$ such that $h(0)\ne0$. 
  Here $E=\{z=0\}$, and the form of $\pi$ implies
  $\pi_*\div_z=b\,\nu_E$. 
  If $\mu_p$ denotes the multiplicity valuation at $p$
  (\ie the monomial valuation with weight $(1,1)$ on $(z,w)$),
  then we also have $\pi_*\mu_p=b\,\nu_F$. 
  As $\mu_p>\div_z$, we conclude that $\nu_F>\nu_E$.
  
  \smallskip
  Now assume that $p$ is satellite.
  Let us first reduce to the case when $\pi$ 
  is a composition of blow-ups at infinitely nearby points. 
  Consider the set of all $\pi'\in\fB$ such that 
  $\pi=\varpi\circ\pi'$ for some composition 
  of blow-ups $\varpi$ such that $\varpi$ is a 
  local biholomorphism at $p$. This set admits a minimum
  thanks to Lemma~\ref{Lsup-inf}. 
  It is not difficult to check that
  this minimum is necessarily a composition of blow-ups at infinitely
  nearby points.

  \smallskip
  We may thus suppose that 
  $\pi=\tpi_0\circ\dots\tpi_{n-1}$ is a composition of point 
  blow-ups at infinitely nearby points $p_0,\dots,p_{n-1}$
  and that $p=p_n$ is a satellite point, say $p=E\cap E'$.
  Define the numbers $n_g$ and $\bn_{g+1}$ as in~\eqref{e312}. 
  Thus $p_j$ is free for $n_g<j\le\bn_{g+1}$ and
  satellite for $\bn_{g+1}<j\le n$.
  Define $p'=p_{n_g+1}$,
  $\varpi=\tpi_{n_g+1}\circ\dots\circ\tpi_{n-1}$
  and $\varpi'=\tpi_0\circ\dots\circ\tpi_{n_g}$.
  Thus $\pi=\varpi'\circ\varpi$.

  We may apply Lemma~\ref{Lkeylemma2} to $\varpi$, with
  the origin being the point $p'$ and the curve $D$
  the exceptional divisor of $\varpi'$ at $p'$.
  This gives coordinates $(z,w)$ at $p$ and $(z',w')$ at $p'$
  in which $\varpi$ is a monomial map,
  say $\varpi(z,w)=(z^\a w^\b,z^{\gamma}w^{\delta})$ 
  for some $\a,\b,\gamma,\delta$. 
  Moreover, the exceptional divisor of $\varpi'$ is given
  by $\{z'=0\}$ at $p'$ and the exceptional divisor of
  $\pi$ by $\{zw=0\}$ at $p$.

  After permuting coordinates we may assume
  $E=\{z=0\}$ and $E'=\{w=0\}$ at $p$. Moreover, $F$ is
  the exceptional divisor of the blowup at $p$.
  Denote by $\mu_p$ the multiplicity valuation at $p$.
  Then $\nu_E$, $\nu_{E'}$ and $\nu_F$ are proportional
  to the pushforwards by $\pi$ of $\div_z$, $\div_w$ and
  $\mu_p$, respectively.

  Since $\varpi$ is monomial, 
  the three valuations $\varpi_*\div_z$, $\varpi_*\div_w$, and
  $\varpi_*\mu_p$ are monomial in the coordinates $(z',w')$.
  at $p'$. Their values on $(z',w')$ are given by 
  $(\a,\gamma)$, $(\b,\delta)$ and $(\a+\b,\gamma+\delta)$,
  respectively. 
  Write $\mu_{w',t}$ for the monomial valuation 
  sending $z'$ to $1$ and $w'$ to $t>0$. 
  This is an element of the relative valuative tree $\cV_{z'}$
  studied in Section~\ref{sec-relative}.
  The three valuations 
  $\varpi_*\div_z$, $\varpi_*\div_w$ and $\varpi_*\mu_p$ 
  are then equivalent to 
  $\mu_{w',\gamma/\a}$, $\mu_{w',\delta/\b}$ 
  and $\mu_{w',(\gamma+\delta)/(\a+\b)}$, respectively.  
  It is clear that $\mu_{w',(\gamma+\delta)/(\a+\b)}$ 
  belongs to the segment
  $[\mu_{w',\gamma/\a},\mu_{w',\delta/\b}]$, 
  hence $\varpi'_*\mu_{w',(\gamma+\delta)/(\a+\b)}\in
  [\varpi'_*\mu_{w',\gamma/\a},\varpi'_*\mu_{w',\delta/\b}]$.
  
  On the other hand, the point $p'$ 
  is a free point, lying on the component $E_g=E_{n_g}$,
  so the pullback under  $\varpi'$ of a smooth generic function 
  is given by $(z')^{b_g}\times\mathrm{units}$, 
  where $b_g$ is the generic multiplicity of $E_g$.
  See Proposition~\ref{Pgeo-inter}. 
  Thus
  \begin{equation*}
    \tilde{\pi}_*\mu_{w',(\gamma+\delta)/(\a+\b)}(\fm)
    =\tilde{\pi}_*\mu_{w',\gamma/\a}(\fm)
    =\tilde{\pi}_*\mu_{w',\delta/\b}(\fm)
    =b_g.
  \end{equation*} 
  This gives
  \begin{equation*}
    \nu_F=b_g^{-1}\tilde{\pi}_*\mu_{w',(\gamma+\delta)/(\a+\b)},
    \quad
    \nu_E=b_g^{-1}\tilde{\pi}_*\nu_{w',\gamma/\a}
    \qand
    \nu_{E'}=b_g^{-1}\tilde{\pi}_*\nu_{w',\delta/\b}. 
  \end{equation*}
  We conclude that $\nu_F\in[\nu_E,\nu_{E'}]$.
\end{proof}
\begin{proof}[Proof of Lemma~\ref{Lkeylemma2}]
  We denote by $\tpi_j$ the blow-up at $p_j$, 
  $E_j\=\tpi_j^{-1}(p_j)$ 
  (as well as its strict transform by further blow-ups).
  We shall construct coordinates $(x,y)$ at the origin and
  $(z_j,w_j)$ at each $p_{j+1}$ such that 
  $\tpi_j$ is a monomial map for all $j<n$. 
  This clearly implies the lemma.
  
  First suppose $n_0=n$ \ie all points are free. The existence of
  the coordinates above is the proved in exactly the same way as
  in Proposition~\ref{P304}; the minor adjustments necessary are
  left to the reader.
  
  Next suppose $n_0=0$, \ie $p_1,\dots,p_n$ are all
  satellite (in the sense that they are singular points
  on the total transform of $D$).
  Fix arbitrary coordinates $(x,y)$ at the origin such that
  $D=\{x=0\}$. 
  Our assumptions imply that $p_1$ is the intersection point of
  $E_0$ and the strict transform of $D$ (also denoted by $D$). 
  We may hence choose coordinates $(z_0,w_0)$ at $p_1$ 
  such that $\{z_0=0\}=E_0$, $\{w_0=0\}=D$ and
  $\tpi_0(z_0,w_0)=(z_0,z_0w_0)$. 
  As $p_2$ is a satellite point, 
  it is necessarily one of the two intersection
  points $E_1\cap D$ (case 1) or $E_1\cap E_0$ (case 2). 
  We may now choose coordinates $(z_1,w_1)$ at $p_2$ 
  such that $\{z_1=0\}=E_1$ and $\{w_1=0\}=D$ in case 1, 
  $=E_0$ in case 2; and
  $\tpi_1(z_1,w_1)=(z_1,z_1w_1)$ in case 1, 
  or $=(z_1w_1,z_1)$ in case 2.
  Inductively we obtain coordinates 
  $(z_j,w_j)$ at $p_{j+1}$ such that 
  $\{z_j=0\}=E_j$, $\{w_j=0\}$ represents the other exceptional 
  component containing $p_{j+1}$, and
  $\tpi_j(z_j,w_j)=(z_j,z_jw_j)$ or $=(z_jw_j,z_j)$. 
  Thus $\tpi_j$ is monomial for $0\le j<n$.
  
  Finally when $0<n_0<n$ we may compose the two constructions
  above. This completes the proof.
\end{proof}
%
%
\subsection{Step 4: $\Phi$ preserves multiplicity}
At this stage we may conclude that 
$\Phi:\Gast\to\cVdiv$ extends to an isomorphism
$\Phi:\Gamma\to\cV$ of rooted, nonmetric $\R$-trees.
Moreover, the Farey parameter defines a parameterization
of $\Gamma$ and $\Phi:(\Gamma,A)\to(\cV,A)$ is an
isomorphism of parameterized trees. 
Notice that the proofs of Corollaries~\ref{Cencoding} and~\ref{C303}
go through. Thus $\Phi$ maps any irreducible curve $C$, viewed
as an end in $\Gamma$, to the corresponding curve valuation $\nu_C$. 
\begin{lemma}\label{L311}
  If $C$ is any irreducible curve, then the multiplicity
  $m_\Gamma(C)$ of $C$ as an element of $\Gamma$ is the same as 
  its multiplicity $m(C)$ as a curve.
\end{lemma}
Now pick any $E\in\Gamma^o$. By the definition of the multiplicity 
in $\cV$ and by Corollary~\ref{C304} we have
\begin{multline*}
  m(\Phi(E))
  =\min\{m(C)\ ;\ \nu_C>\Phi(E)\}
  =\min\{m(C)\ ;\ C>E\}\\
  =\min\{m_\Gamma(C)\ ;\ C>E\}
  =m(E).
\end{multline*}
This immediately implies that $m\circ\Phi=m$ on all of 
$\Gamma$, completing Step~4 and thus the whole proof of
Theorem~\ref{thm-universal}.

\begin{proof}[Proof of Lemma~\ref{L311}]
  Let $(p_j)_0^\infty$ be the sequence of infinitely
  nearby points associated to $C$. 
  This is of Type~4 by Proposition~\ref{P304}.
  Pick $n$ minimal such that $p_j$ is free for $j>n$,
  let $\pi$ be the composition of blowups at
  $p_0,\dots,p_n$ and denote by $E$ the exceptional
  divisor obtained by blowing up $p_n$. 

  Let $C'$ be the strict transform of $C$ under $\pi$.
  Then $C'$ intersects $\pi^{-1}(0)$ only at $p_{n+1}$. 
  As the points $p_j$ are free for all $j>n$, $C'$
  must be smooth and $m_\Gamma(C)=b$, where $(a,b)$
  is the Farey weight of $E$. See Corollary~\ref{C307}.

  Lemma~\ref{L401} provides us with local
  coordinates $(z,w)$ at $p_{n+1}$ and $(x,y)$ at the
  origin, such that
  $E=\{z=0\}$, $C'=\{w=0\}$ and such that
  $\pi(z,w)=(z^b,z^{a-b}w+z^b h(z))$,
  where $h$ is a regular function with $h(0)\ne0$.
  A parameterization of $C$ is then given by
  $t\mapsto\pi(t,0)=(t^b,t^bh(t))$, showing that
  $m(C)=b=m_\Gamma(C)$.
\end{proof}
%
%
%
%
\section{Applications}\label{Sappliuniv}
The fact that $\Phi: \Gamma\to \cV$ is a tree isomorphism leads
to interpretations in more geometrical terms of several constructions
in the valuative tree. We shall give examples of this principle in
subsequent sections. From the proof that $\Phi$ is order preserving
(Step 3 of the proof), we also extract a monomialization procedure for
arbitrary quasimonomial valuations. 
This is explained in Section~\ref{Smonomial}.
%
%
\subsection{Curvettes}\label{sec-curvette}
Fix a divisorial valuation $\nu$, and an irreducible curve $C$.  
We say that $C$ defines a \emph{curvette}\index{curvette} 
for $\nu$ if there exist
$\pi\in\fB$ and $E\in\Gast_\pi$, such that $\nu=\nu_E$,
and the strict transform of $C$ is smooth and intersects 
$E$ transversely at a free point.
Before stating the proposition characterizing
curvettes inside $\cV$, let us introduce some terminology.  
For a fixed divisorial valuation $\nu$, 
a tangent vector $\vv$ at $\nu$ is
\emph{generic}\index{tangent vector!generic} if it is not represented
by $\nu_\fm$, and its multiplicity is the generic multiplicity of
$\nu$. By Proposition~\ref{P111}, any divisorial valuation admits at
most two tangent vectors which are not generic.
\begin{proposition}\label{P410}
  Pick $E\in\Gast$ and let $\nu=\nu_E$ be the associated divisorial
  valuation. Then if $C$ is an irreducible curve at the origin, 
  the following assertions are equivalent:
  \begin{itemize} 
  \item[(i)]
    $\nu_C>\nu$, $m(C)=b(\nu)$, the generic multiplicity of
    $\nu$, and $\nu_C$ represents a generic tangent vector at $\nu$; 
  \item[(ii)]
    $C>E$, $m(C)=b(E)$, the generic multiplicity of $E$, 
    and $C$ represents a generic tangent vector at $E$; 
  \item[(iii)]
    $C$ is a curvette for $\nu$.  
  \end{itemize}
\end{proposition}
It was observed by Spivakovsky, that a divisorial valuation acts by
intersection with a curvette: see~\cite[Theorem~7.2]{spiv}. 
In fact for any $\psi\in\fm$, $\nu(\psi)=\nu_C(\psi)$ 
as soon as $C$ is a curvette for $\nu$ and does 
not define the same tangent vector at $\nu$ as the curve valuation 
associated to any irreducible factor of $\psi$.
\begin{proof}
  The equivalence of~(i) and~(ii) is a direct consequence of
  the isomorphism between $\Gamma$ and $\cV$. Hence we need
  only show that~(ii) and~(iii) are equivalent. Notice
  that this is a statement purely inside $\Gamma$. 
  Also recall that the multiplicity of an irreducible
  curve $C$ coincides with its multiplicity in
  the universal dual graph $\Gamma$.

  First suppose $C$ is a curvette for $\nu$. 
  Pick $\pi\in\fB$, and an exceptional component 
  $E\in\Gast_\pi$ such that the strict transform of 
  $C$ by $\pi$ intersects $E$ at a free point $p$.
  By Corollary~\ref{C305}, we infer that $C>E$ and
  that $C$ represents a generic tangent vector at $E$.
  Further, Corollary~\ref{C307} gives $m(C)=b(E)$.
  Thus~(ii) holds.

  Conversely, suppose $C$ satisfies~(ii). 
  Consider the sequence of infinitely
  nearby points $(p_j)_0^n$ associated to $E$,
  and let $\pi\in\fB$ be the composition of blow-ups 
  at the points $p_0,\dots,p_n$.
  Then $E\in\Gast_\pi$ is obtained by blowing up the last point 
  $p_n$.
  The strict transform $C'$ of $C$ by $\pi$ intersects $E$ 
  at a point $p$. 
  As $C$ defines a generic tangent vector at $E$,
  $p$ is free by Corollary~\ref{C305}.
  As $C'$ is smooth and transverse to $E$ at $p$,
  $m(C)=b(E)$ by Corollary~\ref{C307}.
  Thus $C$ is a curvette for $\nu=\nu_E$,
  which completes the proof.
\end{proof}
%
%
\subsection{Centers of valuations}
Consider a valuation $\nu\in\cV$ and a proper birational morphism
$\pi\in\fB$, say $\pi:X\to(\C^2,0)$.  In Section~\ref{sec-equivalence}
we defined the center
\index{valuation!center of} 
of $\nu$ on $X$.
Here we shall compute the center in terms of tree data, using the
fundamental isomorphism between $\Gamma$ and $\cV$.

Recall more precisely that the center of $\nu$ on $X$ is
either an irreducible component $E$ of $\pi^{-1}(0)$ (\ie an element
of $\Gast_\pi$) or a (closed) point $p$ on $\pi^{-1}(0)$.  In the
first case, $\nu$ is the divisorial valuation $\nu_E$ associated to
$E$. In the second case, $\nu=\pi_*\mu$, where $\mu$ is a centered
valuation on the ring $R_p$ of formal power series at $p$, \ie there
exist local coordinates $(z,w)$ at $p$ such that $\mu(z),\mu(w)>0$.

Let us analyze the second case in more detail.
Consider the sequence $(q_j)_0^\infty$ of
infinitely nearby points associated to $\mu$. 
Thus $q_0=p$ and $0\le n\le\infty$. 
Let $\tpi_j$ be the blow-up at $q_j$, $F_j=\tpi_j^{-1}(q_j)$
and set $\varpi_j=\tpi_0\circ\dots\circ\tpi_j$.
For each $j$ we may apply Corollary~\ref{C309} to $\varpi_j$.
This shows that $F_j$ represents the tangent vector 
$\vv_p$ defined by $p$, at $E$.
In particular, the segment $[F,F_j]$ contains $F$ for any
$F\in\Gast_\pi$.

When $\nu$ is divisorial, $n<\infty$ and $\nu=\nu_{F_n}$.
As $\Phi$ is an isomorphism of rooted, nonmetric trees, we
conclude that $[\nu_F,\nu]$ contains $\nu_E$ for any
$F\in\Gast_\pi$. 
When $\nu$ is nondivisorial, Corollary~\ref{Cencoding} 
implies that $\nu_{F_j}\to\nu$ as $n\to\infty$.
Again we conclude that $[\nu_F,\nu]$ contains $\nu_E$ for any
$F\in\Gast_\pi$. 

We may summarize our result as follows.
\begin{proposition}\label{Pcenter-of-valuation}
  Pick $\pi\in\fB$ and $\nu\in\cV$.  
  Let  $\cE$ be the set consisting of divisorial valuations 
  $\nu_E$ with $E\subset\pi^{-1}(0)$ such that
  $[\nu_E,\nu]$ contains no other divisorial valuations 
  $\nu_F$ with $F\subset\pi^{-1}(0)$. 
  It consists of one or two valuations.
  \begin{itemize}
  \item[(i)]
    When $\cE=\{\nu_E\}$, either $\nu=\nu_E$ and the center of
    $\nu$ in $\pi^{-1}(0)$ equals $E$; or $\nu\ne\nu_E$ and the center
    of $\nu$ in $\pi^{-1}(0)$ is the (free) point on $E$ associated by
    Proposition~\ref{Ptangentgamma} 
    to the tangent vector $\vv$ at $\nu_E$ represented by $\nu$.
\item[(ii)]
    When $\cE = \{\nu_E,\nu_{E'}\}$, the center of $\nu$ is equal to
    the (satellite) point $E\cap E'$. 
    This point is on $E$ ($E'$) the point 
    associated to the tangent vector at $\nu_E$ ($\nu_{E'}$) 
    represented by $\nu$.
  \end{itemize}
\end{proposition}
\begin{remark}\label{R304}
  By applying this result to a curve valuation 
  $\nu=\nu_C$ we obtain 
  a description of the intersection point of the 
  strict transform of $C$ with the exceptional divisor 
  $\pi^{-1}(0)$.
  (As follows from the proof, this could have been achieved
  without passing to the valuative tree.)
\end{remark}
%
%
%
%
\subsection{Potpourri on divisorial valuations}
In this section we prove three results on divisorial valuations.

First, we describe which divisorial valuations are obtained by blowing
up a free point $p$ on some exceptional component.
\begin{proposition}
  Let $\nu$ be a divisorial valuation with 
  associated sequence $(p_0,\dots,p_n)$ of 
  infinitely nearby points.
  Then $m(\nu)=b(\nu)$ iff $p_n$ is free.
\end{proposition}
\begin{proof}
  If $p_n$ is free, then by Lemma~\ref{L415} the multiplicity of the
  exceptional divisor $F\in\Gast$ associated to the blow-up at $p_n$
  satisfies $m(F) = b(F)$. As $\Phi$ preserves multiplicity, 
  $m(\nu)=b(\nu)$.
  
  For the converse, suppose $p_n$ is satellite, say the 
  intersection of two components $E$ and $E'$ of generic 
  multiplicity $b$ and $b'$ respectively. Then the
  generic multiplicity of $F$ is $b+b'$, whereas 
  its multiplicity is bounded by
  $\max\{ m(E), m(E')\} \le \max\{ b, b'\}$. Hence $b(F) > m(F)$, and
  $b(\nu) > m(\nu)$.
\end{proof}

\smallskip
Second, let us summarize what happens in general when 
blowing up a point $p$ on an exceptional component. 
\begin{proposition}\label{prop-divisorial}
  Fix $\pi\in\fB$, pick a point $p\in\pi^{-1}(0)$,
  and define $\nu_F$ to be the divisorial valuation associated to the
  blow-up at $p$.
  \begin{itemize}
  \item[(i)]
    If $p$ is a free point, \ie $p$ 
    belongs to a unique exceptional component $E$ of $\pi$, then:
    \begin{itemize}
    \item[(a)]
      $\nu_F>\nu_E$ and $\nu_F$ does not define the same
      tangent vector at $\nu_E$ as any $\nu_{E'}$ for
      $E'\in\Gast_\pi\setminus\{E\}$;
    \item[(b)]
      the multiplicity $m(\nu_F)$ of $\nu_F$ is equal to its generic
      multiplicity $b(\nu_F)$, and both coincide with $b(\nu_E)$;
      moreover, the multiplicity is constant, equal to 
      $b(\nu_E)$, on the segment $]\nu_E,\nu_F]$;
    \item[(c)]
      $A(\nu_F)=A(\nu_E)+b(\nu_E)^{-1}$. 
    \end{itemize}
  \item[(ii)]
    If $p$ is a satellite point, \ie $p$ is the intersection point of
    two exceptional components $E$ and $E'$, then:
    \begin{itemize}
    \item[(a)]
      $\nu_{E'}>\nu_F>\nu_E$ or $\nu_E>\nu_F>\nu_E'$;
    \item[(b)]
      the multiplicity is constant, equal to 
      $\max\{m(\nu_E),m(\nu_E')\}$ on the segment 
      $]\nu_E,\nu_{E'}[\,$; moreover,
      the generic multiplicity of 
      $\nu(F)$ is given by $b(\nu_F)=b(\nu_E)+b(\nu_{E'})$;
    \item[(c)]
      $A(\nu_F)=(a(\nu_E)+a(\nu_{E'}))/(b(\nu_E)+b(\nu_{E'}))$,
      where $a(\nu_E)-1$ and $a(\nu_{E'})-1$ are the orders 
      of vanishing of the Jacobian of $\pi$ 
      along $E$ and $E'$, respectively.
    \end{itemize}
  \end{itemize}
\end{proposition}
\begin{proof}
  By the fundamental isomorphism $\Phi:\Gamma\to\cV$ it suffices
  to prove the corresponding statements in the universal dual graph.
  Most of them are then straightforward consequences of the 
  combinatorial definition of the partial ordering and Farey weights.
  Specifically, assertions~(a) and~(c) in both~(i) and~(ii) 
  are immediate; (see also Proposition~\ref{Pgeo-inter}) and
  assertion~(b) in~(i) follows from Lemma~\ref{L415}. 

  Let us prove the first statement in assertion~(b) in~(ii) 
  for completeness and as we shall use it below.
  For this, we show in general that if $\pi\in\fB$ and
  $E,E'$ are adjacent elements in $\Gast_\pi$ 
  (\ie $E$ and $E'$ intersect), then 
  either $E<E'$ or $E'<E$, 
  and the multiplicity is constant on the segment $]E,E'[$.

  For this, we may assume that $\pi$ is minimal such that
  $E,E'\in\Gast_\pi$. The proof now goes by induction on the 
  cardinality of $\Gast_\pi$. 
  Without loss of generality, $E'$ is obtained by blowing up a 
  point on $E$. 
  If this point is free, then $E'>E$ and the multiplicity is 
  constant equal to $b$ on $]E,E'[\,$, where $(a,b)$ is the Farey weight 
  of $E$. 
  If the point is not free, it is the intersection point between $E$ and
  another irreducible components $E''$. By induction we may assume
  $E''>E$ and that the multiplicity is constant on $]E,E''[$.
  As $E'\in]E,E''[\,$, this completes the proof.
\end{proof}

\smallskip
Third, let us define two divisorial valuations $\nu$, $\nu'$ to be
\emph{adjacent}\index{valuation!adjacent} if there exists a
composition of blowups $\pi\in\fB$ such that $\nu$ and $\nu'$ are
proportional to the pushforward of the order of vanishing along two
irreducible components $E$, $E'$ of $\pi^{-1}(0)$ with nonempty
intersection.  In other words, $\nu=\nu_E$, $\nu'=\nu_{E'}$ with
$E,E'$ being adjacent vertices in some $\Gamma_\pi$.

We wish to characterize when two divisorial valuations are adjacent,
purely in terms of quantities on the valuative tree $\cV$. 
It follows from Proposition~\ref{prop-divisorial} 
 that if $\nu$ and $\nu'$ are 
adjacent divisorial valuations, then either
$\nu<\nu'$ or $\nu'<\nu$. Further, the
multiplicity is constant, equal to $m=\max\{m(\nu),m(\nu')\}$,
on the segment $]\nu,\nu'[$. 
By lower semicontinuity we then conclude that the multiplicity is
constant on $]\nu,\nu']$ if $\nu<\nu'$ and
$]\nu',\nu]$ if $\nu'<\nu$.
\begin{proposition}\label{P411}
  Let $\nu$ and $\nu'$ be divisorial valuations with generic 
  multiplicities $b$ and $b'$, respectively. 
  Assume that $\nu<\nu'$ and that the multiplicity is
  constant, equal to $m=m(\nu')$ on the segment $]\nu,\nu']$. 
  Then
  \begin{equation}\label{e301}  
    A(\nu')-A(\nu)\ge\frac{m}{b\,b'},
  \end{equation}
  with equality iff $\nu$ and $\nu'$ are adjacent.
\end{proposition}
\begin{remark}\label{rem-compumul}
  As a consequence, if $\nu<\nu'$ are adjacent valuations, then 
  the multiplicity equals $|ab'-ba'|$ on the segment $]\nu,\nu']$.
  Here $(a,b)$ and $(a',b')$ denote the Farey weights
  of (the elements of $\Gast$ associated to) 
  $\nu$ and $\nu'$, respectively.
\end{remark}
\begin{proof}
  Let us first prove~\eqref{e301}. By Proposition~\ref{P111}, the
  multiplicity $m = m(\nu')$ divides $b' = b(\nu')$. As $\nu' > \nu$,
  the tangent vector $\vv$ represented by $\nu'$ at $\nu$ is not
  represented by $\nu_\fm$. By assumption $m(\vv) = m$. 
  Again by Proposition~\ref{P111}, either $m(\vv)=m(\nu)$
  or $m(\vv)=b(\nu)$. 
  In both cases, $m$ divides $b(\nu)$.  Write $\nu=\nu_{E}$ and
  $\nu'=\nu_{E'}$, for $E,E'\in\Gast$ and let $(a,b)$, $(a',b')$ be the
  Farey weights of $E$ and $E'$, respectively. Note that $b= b(\nu)$,
  $b'= b(\nu')$, hence $m$ divides $a'b-ab'$.
  Then~\eqref{e301} holds since
  \begin{equation*}
    A(\nu')-A(\nu)
    =\frac{a'}{b'}-\frac{a}{b}
    =\frac{m}{bb'}\frac{a'b-ab'}{m}
    \ge\frac{m}{bb'}.
  \end{equation*}

  Next we show that for any $\pi\in\fB$, and any adjacent vertices in
  $\Gast_\pi$, equality in~\eqref{e301} holds.  We proceed by induction on
  the number of blowups in $\pi$, writing $\pi=\pi'\circ\tpi$ with
  $\tpi$ being the blow-up at a point $p$. We let $E$ be the exceptional
  divisor of $\tpi$. Using the induction step, we only need to
  prove~\eqref{e301} for $E$ and a vertex adjacent to $E$. When $p$ is
  free, lying on a divisor $F$ with Farey weight $(a,b)$, the Farey
  weight of $E$ is $(a+1,b)$ by definition, and the multiplicity in 
  $]E,F[$ is constant (this is clear in $\Gast$), 
  equal to $m\=b=b(E)=b(F)$. 
  Whence $A(E)-A(F)=\frac{(a+1)b-ba}{b^2}=\frac{m}{b^2}$. 
  This proves the induction step in this case, as $E$ is
  adjacent to a unique vertex $F$.  When $p$ is satellite, intersection
  of two divisors $F_1$, $F_2$ with Farey weights $(a_1,b_1)$,
  $(a_2,b_2)$, $E$ belongs to the segment $[F_1,F_2]$ by
  construction and has Farey weight $(a_1+a_2,b_1+b_2)$ by
  definition. The multiplicity on the segment $]F_1,F_2[$ is constant
  (again look in $\Gast$), equal to $m=|a_1b_2-a_2b_1|$ by the
  induction step. This immediately implies~\eqref{e301} for both
  pairs $E,F_1$ and $E,F_2$. As $E$ is adjacent to either $F_1$ or
  $F_2$ this completes the induction step.

  Finally suppose $E<E'$, that the multiplicity is constant equal to $m$
  on $]E,E'[$, and that $|a'b-ba'|=m$ where $(a,b)$ and $(a',b')$ are
  the Farey weights of $E$ and $E'$, respectively. 
  We want to prove that $E$ and
  $E'$ are adjacent. Let $\pi\in\fB$ be minimal (for the order relation
  $\trianglerighteq$, see Lemma~\ref{Lsup-inf}) such that
  $E,E'\in\Gast_\pi$. Clearly $E$ and $E'$ are adjacent iff they are
  adjacent vertices in this dual graph $\Gamma_\pi$.  We can write
  $\pi=\pi'\circ\tpi$, where $\pi'\in\fB$, $\tpi$ is the blowup at
  $q\in(\pi')^{-1}(0)$ and $\tpi^{-1}(q)\in\{E,E'\}$.

  First suppose $q$ is a free point lying on a single component $F$. 
  Then $\tpi^{-1}(q)=E'$ or else $E\not<E'$ 
  by Proposition~\ref{prop-divisorial}.
  The Farey weight of $F$ equals $(a'-1,b')$. 
  If $F\ne E$, then $F\in\,]E,E'[$ 
  so that the multiplicity of $F$ is equal to $m$. 
  The preceding argument shows 
  $\frac{m}{bb'}\le\frac{a'-1}{b'}-\frac{a}{b}=\frac{m}{bb'}-\frac{1}{b'}$, 
  a contradiction. 
  Thus $F=E$ and $E'$ is adjacent to $E$.  
  When $q$ is satellite, it is the intersection point of $F_1< F_2$, 
  whose Farey weights are $(a_1,b_1)$
  and $(a_2,b_2)$, respectively. If $E\ne F_1,F_2$, then $E<F_1<E'<F_2$,
  thus $F_1$ and $F_2$ have multiplicity $m$. By what precedes, 
  \begin{equation*}
    A(E')-A(E)
    =A(E')-A(F_1)+A(F_1)-A(E)
    \ge\frac{m}{b'b_1}+\frac{m}{b_1b}
    =\frac{m}{bb'}\frac{b+b'}{b_1}.
  \end{equation*}
  But $b'=b_1+b_2>b_1$, thus $A(E')-A(E)>\frac{m}{bb'}$.
  This contradiction concludes the proof.
\end{proof}
%
%
%
%
\subsection{Monomialization}\label{Smonomial}
\index{monomialization (of a quasimonomial valuation)}
Let $\nu$ be a quasimonomial valuation.  By
Theorem~\ref{thm-universal}, $\Phi^{-1}(\nu)$ is either a branch point
or a regular point in $\Gamma$, and by Theorem~\ref{Tnot3} the
sequence of infinitely nearby points $\bp=(p_j)$ associated to $\nu$
is either of Type~0 or~2. Let $n$ be the smallest integer such that
for some $k\ge0$, $p_{n+1},\cdots, p_{n+k-1}$ are all free, and
$p_{n+k+l}$ are satellite for all $l\ge0$.  Thus $n=n_g$ in the
notation of Section~\ref{S404}.  Let $\pi\in\fB$ be the composition of
the point blow-ups at $p_0,\cdots,p_n$.  By Lemma~\ref{Lkeylemma2},
the valuation determined by $p_{n+1},\dots,p_{n+l}$ is monomial in
suitable coordinates at $p=p_{n+1}$, for all $l\ge1$.  When $\bp$ is
of Type~2, the valuation associated to the infinite sequence
$(p_j)_{n+1}^\infty$ is also monomial by continuity.  Therefore the
valuation $\nu$ is equivalent to $\pi_*\mu$ for some monomial
valuation $\mu$ at $p$. We have thus proved
\begin{proposition}\label{monomialization}
  Let $\nu\in\cVqm$ be a quasimonomial valuation which is not monomial. 
  Then there exists a proper birational morphism $\pi\in\fB$, 
  a smooth point $p\in\pi^{-1}(0)$ and a monomial valuation $\mu$ 
  centered at $p$, such that $\pi_*\mu$ is proportional to $\nu$. 
\end{proposition}
Note that in~\cite{ELS} it is proved that, in any dimension, 
an Abhyankar valuation (\ie giving equality in~\eqref{AbIn})
of rank one can be made monomial by a proper birational morphism. 
In our setting, the Abhyankar valuations of rank one 
are exactly the quasimonomial ones. 

Proposition~\ref{monomialization} above is, however, 
more precise than what can be extracted from~\cite{ELS}. 
First, the monomial valuation
is centered at a \emph{point}---something that cannot be guaranteed in
higher dimensions. Second, the birational morphism is explicitly
constructed. In fact, it can be detected in terms of data
in the valuative tree as follows:
\begin{corollary}\label{C301}
  The divisorial valuation associated to the sequence
  $p_0,\dots,p_n$ of infinitely nearby points
  above is exactly the last element $\nu_g$ in the
  approximating sequence of $\nu$.
\end{corollary}
\begin{proof}
  This is a consequence of the fact that the isomorphism
  $\Phi:\Gamma\to\cV$ preserves the partial ordering and
  multiplicity. 
  More precisely, we described in Section~\ref{S404} the dual
  graph associated to the sequence $\bp=(p_j)$ of infinitely
  nearby points of Type~0 and Type~2. 
  See Figures~\ref{F3} and~\ref{F17}. 
  In particular,
  this analysis shows that $E_g<\gamma(\bp)$, where 
  $E_g\in\Gast$ is the exceptional divisor of the blowup
  at $p_n$, $n=n_g$.
  In Section~\ref{S403} we described the restriction of 
  the multiplicity function to this dual graph. 
  See Figure~\ref{F9}. In particular,
  $E_g$ can be characterized as the 
  maximum element in the segment $[E_0,\gamma(\bp)]$ 
  having multiplicity strictly smaller
  than $m(\gamma(\bp))$. 
  On the other hand, by the definition of the 
  approximating sequence, $\nu_g$ is 
  the maximum element in the segment $[\nu_\fm,\nu]$
  having multiplicity strictly smaller than $m(\nu)$.
  As $\Phi$ preserves the partial ordering and multiplicity,
  and as $\Phi(\gamma(\bp))=\nu$, 
  this shows that $\Phi(E_g)=\nu_g$.
\end{proof}
%
%
%
%
\section{The dual graph of the minimal desingularization}\label{Sdualmini}
A \emph{desingularization}
\index{curve!desingularization of}
\index{desingularization!of a curve} 
of a (reduced, formal) curve
$C$ is a composition of point blow-ups $\pi\in\fB$ such that the total
transform $\pi^{-1}(C)$ has normal crossings.  By
Lemma~\ref{Lsup-inf}, the set of all desingularization maps admits a
minimal element, the \emph{minimal desingularization},
\index{desingularization!minimal}
\index{curve!minimal desingularization of}
\index{$\pi_C$ (minimal desingularization)}
which we denote by $\pi_C$.

Our aim in this section is two-fold. First we describe the embedding
of the dual graph $\Gamma_C$
\index{$\Gamma_C$ (dual graph of minimal desingularization, simplicial tree)} 
of $\pi_C$ inside the
universal dual graph. In particular, when $C$ is irreducible, we shall
see that the branch points of $\Gamma_C$ correspond to the
approximating sequence of $C$.  Then we explain how to recover
$\Gamma_C$ from the following finite set of data: the Farey
parameters of all elements in the approximating sequences of the
irreducible components of $C$, and the intersection multiplicities
between these components. It is well-known since the work of Zariski
that these data determine exactly the topology of the embedding of
$C$. In the literature it also is referred to as the 
\emph{equisingularity type} of $C$.
\index{curve!equisingularity type of} 
In Appendix~\ref{S501}, we
shall discuss a way of encoding these data in a tree called the Eggers
tree.

An algorithm describing the dual graph $\Gamma_C$ in terms of the
equisingularity type of $C$ is already described 
in~\cite[Section~1.4.3-1.4.5]{GB} 
using decompositions into continued fractions. 
In our algorithm this decomposition is included in the
recursive computation of the ``Farey weights'' of the exceptional
divisors.  The algorithm of~\cite{GB} was implemented in a
computer by the Spanish researchers A.~and
J.~Castellanos. Unfortunately, to our present knowledge these works
have not been published yet, and there seems to be no other precise
reference concerning this problem.  We hope that our approach
will lead the interested reader to read the excellent work of
E.~G.~Barroso.

\medskip
The dual graph may be viewed equivalently as a simplicial tree 
(\ie a collection of vertices and edges) $\Gamma_C$, or as an $\N$-tree 
(\ie a finite poset) $\Gast_C$.  Then $\Gast_C$ is the set
\index{$\Gast_C$ (dual graph of minimal desingularization, $\N$-tree)} 
of vertices in $\Gast_C$; see Section~\ref{S303}.
%
%
\subsection{The embedding of $\Gast_C$ in $\Gast$}
First suppose $C$ is irreducible. Then $C$ defines an
end in $\Gamma$. Recall that this end is defined in terms 
of the sequence $\bp=(p_j)_0^\infty$
of infinitely nearby points
associated to $\Gamma$: as $\bp$ is of Type~4, the points
$p_j$ are free for large $j$, and the components 
$E'_j\in\Gast$ increase to $C$ as $j\to\infty$.
Define indexes $n_i$, $0\le i\le g$ and $\bn_i$, $1\le i\le g$
as in~\eqref{e317}, \ie $p_j$ is free for $n_i<j\le\bn_{i+1}$
and satellite for $\bn_i<j\le n_i$. In particular, 
$p_j$ is free for $j>n_g$.

\begin{figure}[ht]
  \begin{center}
    \includegraphics[width=\textwidth]{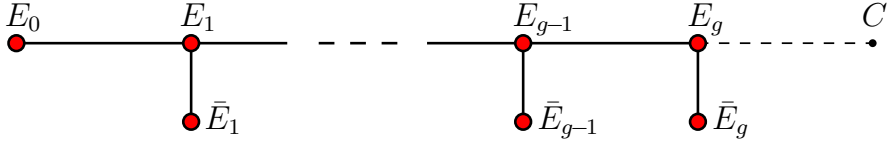}
  \end{center}
  \caption{The dual graph of the minimal desingularization
    of an irreducible curve.
    Not all points are marked.}\label{F19} 
\end{figure}

For $j$ large enough, the strict transform of $C$ is smooth and
transverse at $p_{j+1}$.  We denote this curve by $C_{j+1}$. Suppose
moreover that $j>n_g$. Then the contraction of the exceptional divisor
containing $p_{j+1}$ maps it onto a free point, namely $p_j$. Thus
$C_j$ is still smooth and transverse at $p_j$.  On the other hand,
when $j=n_g$, the contraction map sends $C_{j+1}$ to a curve passing
through $p_{n_g}$ which is a satellite point. Whence the total
transform of $C$ by the blowups at $p_0,\dots,p_{n_g-1}$ has not
normal crossing singularities.  This yields
\begin{proposition}\label{P310}
  The minimal desingularization $\pi_C$ of $C$ is given by the
  composition of blowups at the sequence of infinitely nearby points
  $p_0,\dots,p_{n_g}$.
\end{proposition}
Write $E'_j$ for the exceptional divisor obtained by blowing
up $p_j$, and set $E_i=E'_{n_i}$, $\bE_i=E'_{\bn_i}$.
Thus $E_1,\dots,E_g$ is the approximating sequence of $C$ and
the dual graph $\Gast_C$ looks as in Figure~\ref{F19}.

\begin{remark}\label{R306}
  This yields a geometric interpretation of the approximating 
  sequence of $C$ (and hence of $\nu_C$). 
  Let $\hat{\Gamma}_C$ be the simplicial graph (or tree) 
  whose vertices are irreducible components of the total transform
  $\pi^{-1}(C)$ and where two vertices are joined by an edge iff
  they intersect. 
  Then $E_1$,\dots, $E_g$ are exactly the branch points in 
  $\hat{\Gamma}_C$.
  Moreover, $E_g$ is the unique vertex with an edge joining the strict
  transform of $C$.
\end{remark}
Now suppose $C$ is reducible, with irreducible components $C_j$. The
minimal desingularization $\pi_C$ of $C$ clearly dominates the minimal
desingularizations $\pi_{C_j}$, hence also their join, which we denote
by $\pi'_C$.  The strict transform $C'_j$ of each $C_j$ by $\pi'_C$ is
smooth and intersects the exceptional divisor $(\pi')^{-1}(0)$ in a
free point. However, it may happen that $C'_j$ and $C'_k$ intersect
for $j\ne k$.  It is then clear that $\pi_C$ is obtained from $\pi'_C$
by blowing up these possible intersection points until the strict
transforms of $C_j$ and $C_k$ no longer intersect for $j\ne k$. Notice
that this involves only free blowups, and that the number of blow ups
needed to separate these components is given by the order of
contact of $C'_j$ with $C'_k$.

\begin{figure}[ht]
  \begin{center}
    \includegraphics[width=\textwidth]{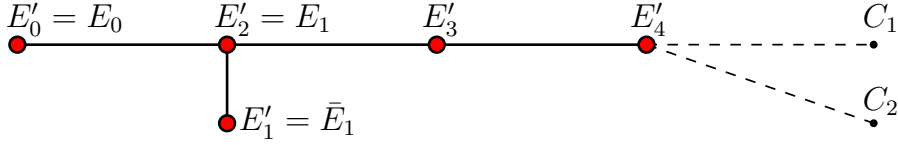}
  \end{center}
  \caption{The dual graph of the minimal desingularization
    of the curve whose irreducible components 
    $C_1$ and $C_2$ are parameterized by 
    $t\mapsto(t^2,t^3)$ 
    and $t\mapsto(t^2(1+t^2),t^3(1+t^2)^2)$, 
    respectively.}\label{F18} 
\end{figure}

This allows us to describe the dual graph $\Gast_C$
as follows. Consider the dual graph of $\pi'_C$.
This is simply the union of the dual graphs 
$\Gast_j:=\Gast_{C_j}$.
For each $j$, define the exceptional component 
$E_{j,g_j}\in\Gast_j$ as above.
If $j\ne k$ and $C'_j\cap C'_k\ne\emptyset$,
then we must have $E_{j,g_j}=E_{k,g_k}=:E_{jk}$ 
and $C'_j$ and $C'_k$ define the same tangent 
vector $\vv_{jk}$ at $E_{jk}$. Thus the exceptional component
$C_j\wedge C_k\in\Gast$ also represents $\vv_{jk}$.
This component is obtained by blowing up a
sequence of infinitely nearby (free) points, 
above $C'_j\cap C'_k$.
From these observations we conclude
\begin{proposition}
  The dual graph $\Gast_C$ of the minimal desingularization
  $\pi_C$ is the union of $\bigcup_j\Gast_{C_j}$ and
  all exceptional components $E\in\Gast$ such that
  $E_{jk}<E\le C_j\wedge C_k$ and $b(E)=b(E_{jk})$ for
  some $j\ne k$.
\end{proposition}
In the discussion of the Eggers tree in Appendix~\ref{S501}
we shall need the following consequence of the proposition
\begin{corollary}\label{C308}
  All branch points in $\Gast_C$ are dominated by some $C_j$.
\end{corollary}
In Figure~\ref{F18} we illustrate the 
dual graph $\Gast_C$, where $C$ has two 
tangential cusps as irreducible components.

%
%
\subsection{Construction of $\Gamma_C$ from the equisingularity type of $C$}
As before, we decompose $C$ into irreducible components $C_j$, and
view them as points in $\Gamma$.  We let $E_{ji}$ for $i=1,\dots,g_j$
be the approximating sequence of $C_j$.  Our aim is to present an
algorithm which has as input all Farey parameters $A(j,i):=A(E_{ji})$
and $A(C_j\wedge C_{j'})$, and as output the dual graph of the minimal
desingularization of $C$.

Note that we preferred working with the Farey parameter $A$ of
$C_j\wedge C_{j'}$ instead of the skewness $\a$ of the corresponding
valuation, because the Farey parameter appears more naturally in our
algorithm. However, the intersection multiplicity between $C_j$ and
$C_{j'}$ is more easily related to $\a$ than to $A$,
see~\eqref{e103}. But note that the knowledge of the Farey parameters
of the approximating sequence is essentially equivalent to the
knowledge of its skewness. Indeed, we may use results of
Section~\ref{S12} to compute the multiplicities, and~\eqref{e412} to
compute the Farey parameters (\ie the thinness thanks to
Theorem~\ref{thm-universal}), see also Table~\ref{table5}.

\medskip
For convenience we denote by $J$ the set of irreducible components of
$C$.  At each step our data consist of several elements: a finite
simplicial graph $\cT$; a function $w$ defined on the set of vertices
of $\cT$ with values in $(\N^*)^2$; to each vertex $E\in\cT$
(resp. each edge $e$ of $\cT$) a subset $J(E), J(e)\subset J$ such
that the collection $\{J(E), J(e)\}$ form a partition of $J$; finally
to each $j\in J$ an integer $I(j)\in \{1,\cdots, g_j\}$.

Before explaining the algorithm, let us interpret these data in
geometric terms.  Each step corresponds to the blow up of the
finitely many intersection points of the strict transform of $C$ with
the exceptional divisor. The graph $\cT$ is the dual graph of the
exceptional divisor at the current step.  In particular, a vertex in
$\cT$ corresponds to a divisor $E\in\Gast$. If $E$ is a vertex of
$\cT$, then $w(E)$ is exactly the Farey weight of $E$. The subset
$J(E)$ is the set of branches of $C$ which intersects the divisor $E$
at a free point. An edge $e$ of $\cT$ corresponds to the intersection
point of two divisors: $J(e)$ is the set of branches of $C$ which
intersects the exceptional divisor at this point.  Finally a branch
$C_j$ intersects one or two irreducible components of the exceptional
divisor which have the same multiplicity $m$.  The number $I(j)$ is
then the unique integer such that $m(E_{j,I(j)})=m$.

When $w=(a,b) \in (\N^*)^2$, we shall write $A(w)= a/b$. 
This is consistent with our previous notation: 
when $w$ is the Farey weight
of $E\in\Gast$, then $A(w)$ is the Farey parameter of $E$.


\smallskip
$\blacktriangle$ \textit{Initiation}: $\cT$ consists of one point $E_0$, 
$w(E_0)=(2,1)$, $J(E_0) = J$, and $I(j) =1$ for any $j\in J$.

\medskip
$\blacktriangle$ \textit{Loop}: we modify $\cT$ by adding a certain number of
vertices/edges to it, leaving the weight $w$ of the previous vertices
unchanged. Modifications of the graph are done exactly at the vertices
and edges for which $J(E)$ or $J(e)$ is non-empty.

-\textit{Modification at vertices}: to each vertex $E$ of $\cT$ for which
$J(E)\ne\emptyset$ apply the following procedure.
\begin{itemize}
\item Step 1. Define a partition of $J(E)$ using the equivalence
relation $j\sim j' \in J(E)$ iff $A(C_j\wedge C_{j'})> A(E)$.
\item Step 2. For each element of this partition $J'\subset J(E)$, add
a new vertex $F$ attached to $\cT$ by one edge $f$ at $E$. Set $w(F) =
w(E) + (1,0)$.
\begin{itemize}
\item[(a)] Define $J(F)$ to be the elements $j\in J'$ for which $A(F)<
A(j,I(j))$. For any $j\in J'$, leave $I(j)$ unchanged.
\item[(b)] Define $J(f)$ to be the elements $j\in J'$ for which $A(F)>
A(j,I(j))$. For any $j\in J'$, leave $I(j)$ unchanged.
\end{itemize}
\item Step 3. Set $J(E) =\emptyset$.
\end{itemize}

-\textit{Modification at edges}: to each edge $e$ of
$\cT$ for which $J(e)\ne\emptyset$ apply the following procedure. We
adopt the following notation. The boundary of $e$ consists of two
vertices $E_1$ and $E_2$, and we will assume $A(E_1)< A(E_2)$.
\begin{itemize}
\item Step 1. Replace the edge $e$ by a vertex $F$ and two edges,
$e_1$ joining $F$ to $E_1$, and $e_2$ joining $F$ to $E_2$. Define
$w(F) = w(E_1)+ w(E_2)$.
\item Step 2. Define $J(F)$ to be the collection of $j\in J(e)$ for which
$A(F) = A(j,I(j))$. For any $j\in J(F)$, replace $I(j)$ by $I(j)+1$.
\item Step 3. Define $J(e_1)\= \{ j\in J(e),\, A(F) > A(j,I(j))\}$,
 and $J(e_2)\= \{ j\in J(e),\, A(F) < A(j,I(j))\}$. For any elements
 in $J(e_1)\cup J(e_2)$, leave $I(j)$ unchanged.
\end{itemize}

$\blacktriangle$ \textit{End}: the algorithm stops when 
\begin{itemize}
\item[(i)] $J(e)$ is empty for all edges of $\cT$;
\item[(ii)] for any couple $j\ne j'\in J(E)$, $A(C_j\wedge C_{j'}) = A(E)$; 
\item[(iii)] for any $j$, $I(j)= g_j$.
\end{itemize}

\medskip

Let us explain a bit what is going on. We start with the dual graph of
the blowup at the origin. Suppose we have made the loop finitely
many times. This produces a finite graph $\cT$ which corresponds to a
certain modification $\pi\in\fB$. We let $C'$ be the strict transform
of $C$ by $\pi$.  

The set of vertices with $J(E) \ne \emptyset$ is the set of
exceptional component $E$ that $C'$ intersects at a free point.  The
step ``modification at a vertex'' corresponds to the blowup at all these
free points.  Pick any component $E$ that $C'$ intersects at a free
point. We first subdivide $J(E)$ into subsets consisting of branches
intersecting $E$ at the same point (Step 1). 
Then we add to $\cT$ a vertex for
each of the intersection points, and compute the Farey weight of the
exceptional components which are created (Step 2).  After these 
blowups some branches will intersect the exceptional divisor at a free
point (Case (a)) or at a satellite point (Case (b)).  In any case,
being created by a blowup at a free point, the new divisor cannot
belong to the approximating sequence of a branch of $C$. Thus $I(j)$
is left unchanged in this case.

The set of edges with $J(e) \ne \emptyset$ is the set of singular
points of the exceptional set that $C'$ contains. Choose one of these
points $p$. This is a satellite point, which is the intersection of
two components $E_1$ and $E_2$. Blowing up $p$ induces an elementary
modification of $\cT$ which consists of adding a vertex $F$ in
between $E_1$ and $E_2$, and the Farey weight of the new exceptional
divisor is the sum of the Farey weights of $E_1$ and $E_2$ (Step 1 of
``modification at edges''). Then several cases appear for the strict
transform of a given branch $C_j$ of $C$ containing $p$. First it
may intersect  $F$ at a free point. This happens precisely when $F$
lies in the approximating sequence of $C_j$ \ie when $A(F) =
A(j,I(j))$. We thus have ``superseded'' the $I(j)$-th element of the
approximating sequence of $C_j$, and thus we add $1$ to $I(j)$ (Step
2).  Or second, it still intersects the exceptional divisor at a satellite
point, $F\cap E_1$ when $A(F)> A(j,I(j))$, or $F\cap E_2$ when $A(F)<
A(j,I(j))$. In these two cases, $F$ can not be in the approximating
sequence of $C_j$. Whence $I(j)$ is left unchanged (Step 3).

Finally, we have achieved the desingularization of $C$ when three
conditions are satisfied. First, the strict transforms of all branches
of $C$ intersect the exceptional divisor at a free point (this is the
condition $J(e)=\emptyset$ for all edges). Second, no two branches
intersect the exceptional divisor at the same point (when $C_j$ and
$C_{kj}$ intersect $E$, this happens precisely when $C_j \wedge C_{j'}
= E$). Finally, the strict transform of each branch $C_j$ is smooth
and transverse to the exceptional divisor (this happens precisely when
the last element of the approximating sequence of each branch has
``appeared'' in the process of blow-up, \ie $I(j) = g_j$).

We leave to the reader to check that this shows that 
the algorithm stops
after finitely many steps and produces $\Gamma_C$.

%
%
%
%

\section{The relative tree structure}\label{S23}
We described in Section~\ref{sec-relative} the relative valuative
$\cV_x$ \ie the union of $\div_x$ and the set of centered valuations
normalized by $\nu(x)=1$, where $x\in\fm$ is an element determining a
smooth formal curve $\{x=0\}$ at the origin. We explain here how to
recover $\cV_x$ geometrically from the dual graphs of compositions of
blow-ups. More precisely we will construct a 
\emph{relative universal dual graph}
\index{dual graph!relative universal} 
\index{relative!universal dual graph}
$\Gamma_x$ and show,
similarly to Theorem~\ref{thm-universal}, that $\cV_x$ is isomorphic
to $\Gamma_x$.

We already saw in Chapter~\ref{sec-puis} that $\cV_x$ was useful for the
approach to the valuative tree through Puiseux series. 
In Section~\ref{S311} we shall describe another situation 
where $\cV_x$ appears naturally.
%
%
\subsection{The relative dual graph}
Let us first construct the relative universal dual graph as a 
nonmetric tree.
Fix $x\in\fm$ with $m(x)=1$ and denote the curve
$\{x=0\}$ by $E_x$. 

Consider the set $\fB$ of proper 
birational morphisms above the origin as
before. To each $\pi\in\fB$ we associate a \emph{relative dual graph}
$\Gamma_{x,\pi}$. This is a finite graph whose set of vertices
$\Gast_{x,\pi}$ is in bijection with the union of (the strict
transform of) $E_x$ and the set of all irreducible components of
$\pi^{-1}(0)$. Two vertices are joined by an edge iff they
intersect. Note that $\Gamma_{x,\pi}$ is always a simplicial tree,
that $\Gast_{x,\pi} = \Gast_\pi\cup\{E_x\}$, and that $E_x$ is joined
by exactly one other vertex.

We endow $\Gast_{x,\pi}$ with the natural partial ordering in which
$E_x$ is the unique minimal element.  As before, the set of all posets
$\Gast_{x,\pi}$ defines an injective system and we denote the limit by
$\Gast_x$.  It is naturally a nonmetric $\Q$-tree rooted in $E_x$.
All of its points are branch points, except $E_x$ which is an end.  In
fact, $\Gast_x$ is naturally isomorphic to the $\Q$-tree
$\Gast\cup\{E_x\}$ rooted in $E_x$.  We construct a canonical
nonmetric $\R$-tree $\Gamma^o_x$ from $\Gast_x$ using
Proposition~\ref{P420} and also consider its completion $\Gamma_x$.
Then $\Gamma^o_x$ is naturally isomorphic to the $\R$-tree
$\Gamma^o\cup\{E_x\}$ rooted in $E_x$, 
and $\Gamma_x$ is isomorphic to $\Gamma$, also rooted in $E_x$.

As in the case of $\Gamma$ (see Section~\ref{S404}), 
the points in $\Gamma_x$ can be uniquely
encoded by sequences of infinitely nearby points of Types~0,1,2 and~4.
The only difference is that $E_x$ is encoded by the empty sequence, 
and that the unique sequence
of Type~4 that encodes the end in $\Gamma$ corresponding to $E_x$,
now encodes the unique tangent vector at $E_x$ in $\Gamma_x$.

Identify $\Gamma_x$ with $\Gamma$ as sets but write
$\le_x$ and $\le$ for their partial orderings.
Also denote the associated infimum operators by 
$\wedge_x$ and $\wedge$.
As before, let $E_0$ be the exceptional divisor obtained by blowing
up the origin (which lies on $E_x$).
Then we have
\begin{lemma}\label{L426}
  If $I$ is a segment in $\Gamma$, 
  then the partial orderings $\le$ and $\le_x$
  on $I$ coincide iff the intersection of $I$ with the 
  segment $[E_x,E_0]$ contains at most one point. 
  On the other hand, on 
  $[E_x,E_0]$ we have $E'\le_x E''$ iff $E'\ge E''$.
\end{lemma}
%
%
\subsection{Weights, parameterization and multiplicities}
Next we define \emph{relative Farey weights}
$(a_x,b_x)$.  
\index{Farey!weight (relative case)}
\index{relative!Farey weight}
These are defined using elementary
modifications just like their nonrelative counterparts, with the
following exception: $E_x$ is set to have relative Farey weight
$(1,1)$.  As $E_0$ is obtained by a free blowup of $E_x$, it has
relative Farey weight $(2,1)$, which is also its nonrelative Farey
weight. This easily gives that the relative and nonrelative Farey
weights agree on the tree $\{E\ge_x E_0\}$.  In general we have the
following relation:
\begin{lemma}\label{L425}
  The two Farey weights of $E\in\Gast$ are related by
  \begin{equation}\label{e435}
    (a_x,b_x)=\left(a,b\left(\frac{c}{d}-1\right)\right),
  \end{equation}
  where $(c,d)$ is the nonrelative Farey weight of
  $E\wedge_xE_0$.
\end{lemma}
\begin{proof}
  We have to prove that~\eqref{e435} holds 
  for any $\pi\in\fB$, and any $E\in\Gast_\pi$.
  The proof goes by induction on the cardinality of $\Gast_\pi$
  and is left to the reader.
\end{proof}

The \emph{relative Farey parameter}
\index{Farey!parameter (relative case)}
\index{relative!Farey parameter}
of $E\in\Gast_x$ is defined by $A_x(E)=a_x/b_x$, where
$(a_x,b_x)$ is the relative Farey weight of $E$. It follows from
Lemma~\ref{L425} that
\begin{equation}\label{e434}
  A_x(E)=A(E)/(A(E\wedge_x E_0)-1).
\end{equation}
This formula easily implies that $A_x$ defines a parameterization
of $\Gamma_x$. Moreover, $A_x=A$ on the subtree $E\ge_x E_0$.

We can also use the relative Farey weights to define a 
\emph{relative multiplicity} function $m_x$ on $\Gamma^o_x$.  
\index{multiplicity!relative}
\index{relative!multiplicity}
\index{$m_x(E)$ (relative multiplicity)}
Namely, we set
\begin{equation}\label{e433}
  m_x(E)=\min\{b_x(F)\ ;\ F\in\Gast_x, F\ge_x E\}.
\end{equation}
This multiplicity function naturally extends to 
the completion $\Gamma_x$.
Lemma~\ref{L426} and Lemma~\ref{L425} easily imply
\begin{lemma}\label{L427}
  For any $E\in\Gamma$ we have $m_x(E)=1$ if
  $E\in[E_x,E_0]$ and
  $m_x(E)=m(E)(A(E\wedge_x E_0)-1)$ otherwise.
\end{lemma}
%
%
\subsection{The isomorphism}
We can now define a relative version of the fundamental isomorphism
$\Phi:\Gamma\to\cV$ in Theorem~\ref{thm-universal}. 
Let us define a map $\Phi_x$ from $\Gast_x$ into $\cV_x$ by 
setting $\Phi_x(E_x)\=\div_x$ and
$\nu_{x,E}=b_x^{-1}\pi_*\div_E(\phi)$
for $E\in\Gast_x\setminus\{E_x\}$.
In view of Lemma~\ref{L425} it is clear that
\begin{equation}\label{e432}
  \Phi_x(E)=\Phi(E)\left(A(E\wedge_x E_0)-1\right)
\end{equation}
Using Theorem~\ref{thm-universal} we have
$A(E\wedge_x E_0)=A(\nu_E\wedge_x\nu_\fm)=1+\nu_E(x)$.
Thus~\eqref{e432} implies 
$\nu_{x,E}=\nu_E/\nu_E(x)$, \ie $\nu_{x,E}$ is the 
centered valuation on $R$ equivalent to $\nu_E$ and 
normalized by $\nu_{x,E}(x)=1$.
\begin{theorem}\label{T402}
  The map $\Phi_x:\Gast_x\to\cV_x$ extends uniquely to 
  an isomorphism of parameterized trees
  $\Phi_x:(\Gamma_x,A_x)\to(\cV_x,A_x)$
  Here $A_x$ denotes the relative Farey parameter on $\Gamma_x$ 
  and relative thinness on $\cV_x$. 
  Further, $\Phi_x$ preserves relative multiplicity.
\end{theorem}
\begin{proof}
  Let $\cVdivx$ denote the set of divisorial valuations
  in $\cV_x$, with the convention that 
  $\div_x\in\cVdivx$.
  Clearly $\Phi_x$ maps $\Gast_x$ into $\cVdiv_x$.
  It suffices to prove that $\Phi_x:\Gast_x\to\cVdivx$ 
  is an isomorphism of rooted, nonmetric $\Q$-trees, 
  and that $\Phi_x$ preserves 
  the parameterization and multiplicity.

  As in Section~\ref{sec-relative} we denote by 
  $N:\cV\to\cV_x$ the map sending $\nu_x$ to $\div_x$
  and $\nu\ne\nu_x$ to $\nu/\nu(x)$. 
  Then $N$ restricts to an isomorphism of nonrooted, nonmetric
  $\Q$-trees $\cVdiv\cup\{\nu_x\}\to\cVdivx$.

  Above we identified $\Gast_x$ and $\Gast\cup\{E_x\}$ as sets.
  Let us for the purpose of this proof think of this identification
  as an isomorphism $N:\Gast\cup\{E_x\}\to\Gast_x$ of nonrooted, 
  nonmetric $\Q$-trees. We saw above that $A_x\circ N=A$.
  
  It is then clear from the construction that 
  $N\circ\Phi=\Phi_x\circ N$ on $\Gast\cup\{E_x\}$.
  Since by Theorem~\ref{thm-universal}, $\Phi$ is an isomorphism
  of nonrooted, nonmetric $\Q$-trees, so is $\Phi_x$.
  As $\Phi_x$ maps $E_x$ to $\div_x$, it is an isomorphism
  of rooted, nonmetric $\Q$-trees.
  Since $A_x\circ N\circ\Phi=A$ and $A_x\circ N=A$ 
  on $\Gast\cup\{E_x\}$ we get $A_x\circ\Phi_x=A_x$
  on $\Gast_x$, so that $\Phi_x$ preserves the parameterization.
  The same reasoning, in conjunction with
  Lemma~\ref{L427} and Proposition~\ref{Prel2},
  shows that $\Phi_x$ preserves multiplicity.
  This completes the proof.
\end{proof}
%
%
\subsection{The contraction map at a free point}\label{S311}
Let us end this section by exhibiting a situation where the 
relative point of view appears naturally.

Pick $\pi\in\fB$, an exceptional component $E\in\Gast_\pi$, 
and a free point $p\in E$. 
Let $(a,b)$ be the Farey weight of $E$.
Pick a regular function $z$ at $p$ such that $E=\{z=0\}$,
and let $\Gamma_z$ and $\cV_z$ be the corresponding relative
universal dual graph and valuative tree, respectively.

The point $p$ defines a tangent vector in $\Gamma$ at $E$.  Write
$U_\Gamma$ for the corresponding (weak) open set, \ie the set of
elements in $\Gamma\setminus\{E\}$ representing this tangent vector.
The closure of $U_\Gamma$ is given by
$\overline{U_\Gamma}=U_\Gamma\cup\{E\}$ (see Lemma~\ref{L310} below).
Notice that $\overline{U_\Gamma}$ is naturally a tree rooted in $E$.

Let $\nu_E\in\cV$ be the divisorial valuation defined by $E$.
The point $p$ also defines a tangent vector at $\nu_E$ in $\cV$.
Write $U_\cV$ for the corresponding (weak) open set, \ie the
set of valuations $\nu\in\cV\setminus\{\nu_E\}$ representing the
same tangent vector at $\nu_E$ as the divisorial valuation
obtained by blowing up $p$.
Again, $\overline{U_\cV}=U_\cV\cup\{\nu_E\}$ and
$\overline{U_\cV}$ is naturally a tree rooted in $\nu_E$.

Define a map $\varpi_\Gamma:\Gamma_z\to\Gamma$ in terms of sequences 
of infinitely nearby points.
Let $(p_j)_0^n$ be the (finite) sequence of infinitely
nearby points associated to $E$. 
Then $\varpi_\Gamma$ maps a sequence $(q_k)_0^l$, $0\le l\le\infty$,
not of Type 3, with $q_0=p$, to the concatenated sequence
$p_0,\dots,p_n,q_0,\dots,q_l$. In particular, the empty sequence is
mapped to $(p_0)_0^n$. Notice that the concatenated sequence is never
of Type~3.

We also define a map $\varpi_\cV:\cV_z\to\cV$ by 
$\varpi_\cV(\mu)=b^{-1}\pi_*\mu$. 
As $\mu(z)=1$, it follows from Lemma~\ref{L401}
that $(\pi_*\mu)(\fm)=b$, so that $\varpi_\cV(\mu)\in\cV$.

\begin{figure}[ht]
  \begin{center}
    \includegraphics[width=\textwidth]{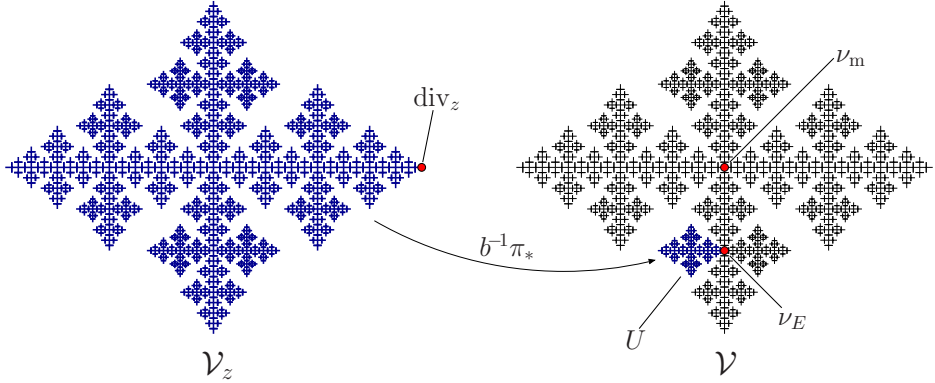}
  \end{center}
  \caption{The contraction map $\pi$ 
    at a free point $p$ on an exceptional
    component $E$ induces an isomorphism between the 
    relative valuative tree $\cV_z$
    and a subset $\overline{U}$ of the valuative tree $\cV$.
    See Theorem~\ref{Trelnat}}\label{F16} 
\end{figure}

\begin{theorem}\label{Trelnat}
  The following properties hold:
  \begin{itemize}
  \item[(i)]
    the map $\varpi_\Gamma$ gives an isomorphism
    $\varpi_\Gamma:\Gamma_z\to\overline{U_\Gamma}$ of rooted,
    nonmetric trees. Moreover, for any $F\in\Gamma_z$ we have:
    \begin{align}
      m(\varpi_\Gamma(F))&=b\,m_z(F)\label{e313}\\
      A(\varpi_\Gamma(F))&=A(E)+b^{-1}(A_z(F)-1)\label{e312}
    \end{align}
  \item[(ii)]
    the map $\varpi_\cV$ gives an isomorphism
    $\varpi_\cV:\cV_z\to\overline{U_\cV}$ of rooted,
    nonmetric trees. Moreover, for any $\mu\in\cV_z$ we have:
    \begin{align}
      m(\varpi_\cV(\mu))&=b\,m_z(\mu)\label{e315}\\
      A(\varpi_\cV(\mu))&=A(\nu_E)+b^{-1}(A_z(\mu)-1)\label{e314}\\
      \a(\varpi_\cV(\mu))&=\a(\nu_E)+b^{-2}\a_z(\mu)\label{e316}
    \end{align}
  \item[(iii)]
    the isomorphisms in~(i) and~(ii) respect the identification of the 
    universal dual graph and the valuative tree: 
    if $\Phi:\Gamma\to\cV$ and $\Phi_z:\Gamma_z\to\cV_z$ are as
    in Theorem~\ref{thm-universal} and Theorem~\ref{T402},
    then $\varpi_\cV\circ\Phi_z=\Phi\circ\varpi_\Gamma$.
  \end{itemize}
\end{theorem}
The theorem shows that the valuative tree (and the universal dual
graph) has a self-similar, or fractal, structure. 
\index{valuation space!self-similar structure}
\index{valuative tree!self-similar structure}
This is illustrated in Figure~\ref{F16}.
\begin{proof}
  It is immediate from the definition of $\varpi_\Gamma$ 
  in terms of concatenations of sequences of 
  infinitely nearby points, 
  from the definition of
  the partial orderings on $\Gast_z$ and $\Gast$, 
  and from Corollary~\ref{C311},
  that $\varpi_\Gamma$ gives an isomorphism of rooted, nonmetric 
  $\Q$-trees between $\Gast_z$ and $\overline{U_\Gamma}\cap\Gast$.
  Thus it also gives an isomorphism between the corresponding rooted,
  nonmetric $\R$-trees.
  To prove~\eqref{e313} and~\eqref{e312}
  we have to control the Farey weights.
  The relative Farey weight of $E_z$ is $(1,1)$.
  Write $(a,b)$ for the Farey weight of $E$.
  A simple induction shows that 
  if $F$ has relative Farey weight $(c,d)$, then 
  $\varpi_\Gamma(F)$ has relative Farey weight
  $(c+(a-1),bd)$. This easily implies 
  both~\eqref{e313} and~\eqref{e312} and completes
  the proof of~(i).

  \smallskip
  Both~(ii) and~(iii) now follow from 
  Theorem~\ref{thm-universal} and Theorem~\ref{T402}
  together with Corollary~\ref{Cencoding} and its corresponding
  relative version. Notice that~\eqref{e316} is a consequence
  of~\eqref{e314} and~\eqref{e315}.

  \smallskip
  We can also give a direct proof of~(ii), not using the
  universal dual graph.
  As we already noticed, Lemma~\ref{L401} implies that
  $\varpi_\cV$ is a well-defined map from $\cV_z$ into $\cV$.
  It is clearly order-preserving and (weakly) continuous.
  It is also injective as $\pi$ is a birational map
  and hence induces an isomorphism between the fraction
  fields of the rings of formal power series at $p$ and at the origin.
  If $C$ is an irreducible curve at the origin, then it follows from
  Proposition~\ref{Pcenter-of-valuation}
  that $\nu_C\in U_\cV$ 
  iff the strict transform of $C$ contains $p$. 
  This immediately implies that $\varpi_\cV$ restricts to
  a bijection between curve valuations in $\cV_z$ and curve
  valuations in $U_\cV$. By continuity, this shows 
  that $\varpi_\cV$ is a surjective map of $\cV_z$ onto
  $\overline{U_\cV}$.
  Equations~\eqref{e315}-\eqref{e316} can also be proved
  directly, using Lemma~\ref{L401}. The details
  are left to the reader.
\end{proof}


%
%
%
%
%
%
\chapter{Tree measures}\label{part-potent}
One of the main motivations behind the present monograph is the fact
that the valuative tree $\cV$ provides an efficient means of encoding
singularities. We shall see in the next chapter that this encoding can
be nicely given in terms of measures on $\cV$.
These measures are obtained indirectly by passing through functions
defined on the tree $\cVqm$ of quasimonomial valuations, 
a procedure analogous to the identification of 
positive measures on the real line with suitably normalized 
concave functions.

Here we develop a general methodology for
identifying certain classes of functions on $\cVqm$ 
with measures on $\cV$. 
In fact, the analysis is purely tree-theoretic and does not
depend on the fact that the elements of the valuative tree are
valuations. We shall therefore work on a general complete, rooted
nonmetric tree, equipped with an increasing parameterization 
taking value 1 at the root.
The main example we have in mind is $\cV$
with the skewness parameterization.

An outline of the theory, as well as a more precise comparison
to the situation on the real line, is given in 
Section~\ref{S304}. We refer to that section for the organization 
of the remainder of the chapter.

While the main results in this chapter are of fundamental importance
for applications,
the techniques are quite different from those of the rest of the
monograph. Hence many of the details can be skipped on a first
reading.
%
%
%
%
\section{Outline}\label{S304}
As the analysis that we are about to undertake is
somewhat delicate,
we start in this section 
by giving an outline of what we are trying to accomplish.
%
%
\subsection{The unbranched case}\label{S305}
It is a standard result that finite positive measures 
on $\R$ can be identified with either decreasing functions
or increasing, concave functions (up to certain normalizations). 
More precisely we have 
\begin{proposition}
  The following objects are naturally in 1-1 correspondence
  with each other:
  \begin{itemize}
    \item[(i)]
      finite positive measures $\rho$ on $\R$;
    \item[(ii)]
      left continuous decreasing functions
      $f:\R\to\R$ with $f(\infty)=0$;
    \item[(iii)]
      increasing concave functions
      $g:\R\to\R$ with $g(0)=0$ and $g'(\infty)=0$.
  \end{itemize}
  This correspondence goes as follows:
  \begin{itemize}
  \item
    to a measure $\rho$ is associated the function
    $f(x)=\rho\,[x,\infty)$;
  \item
    to a function $f$ is associated the function
    $g(x)=\int_0^xf(y)\,dy$;
  \item
    to a function $g$ is associated 
    the measure $\rho=-d^2g/dx^2$.
  \end{itemize}
\end{proposition}

We shall extend this correspondence in two ways.
First, we shall work with complex measures in~(i).
The corresponding functions in~(ii) are then known as
(normalized) functions of bounded variation, and
the functions in the analogue of~(iii) are
characterized by having a complex measure as
second derivative (in the sense of distributions);
again subject to a normalization.

The more difficult extension consists of replacing
$\R$ with a general, complete, parameterized tree.
This presents two challenges. 
First, as opposed to a complete tree, the real line has neither a 
maximal nor a minimal element. 
Essentially, this means that 
the normalizations of the functions corresponding to the
classes~(ii) and~(iii) above have to be done with great care.
Second, and more importantly, a tree typically has branch points.
This forces us to redefine notions such as ``increasing'',
``concave'' and ``bounded variation''.
Notice that in general there is no bound on the 
amount of branching that a tree can exhibit. 
Already in the case of the valuative tree, the tangent space at a 
divisorial valuation has the cardinality of the continuum.
Luckily, we shall not have to worry about the amount of
branching, and in any case it turns out that 
the measures and functions that we shall consider
all ``live'' on completions of 
countable unions of finite subtrees.
%
%
\subsection{The general case}\label{S306}
We now outline the analysis that we will undertake in
the remainder of the chapter. Our objective is
to provide the reader with an intuitive picture 
rather than stating all the precise definitions. 

Consider a complete, parameterizable, rooted
nonmetric tree $\cT$. Let $\cT^o$ denote the set of non-ends
in $\cT$. By convention, the root $\tau_0$ of $\cT$ belongs to
$\cT^o$.
Also fix an increasing parameterization $\a$ of $\cT$
with $\a(\tau_0)=1$. The choice of parameterization is
important for some, but not all, of the constructions below.

We will work exclusively with the \emph{weak topology} on $\cT$.
Any parameterizable complete tree is weakly compact.
Denote by $\cM$ the set of \emph{complex Borel measures} on $\cT$.
We think of the elements of $\cM$ as complex valued functions
on the Borel $\sigma$-algebra generated by the weak open sets,
but $\cM$ can also be viewed as the dual of $C(\cT)$. 
The latter point of view provides $\cM$ with a norm,
turning it into a Banach space.

\smallskip
Our main objective is to identify $\cM$ with two classes of
functions on $\cT^o$. This goes as follows. To any
complex Borel measure $\rho\in\cM$ we associate two complex-valued
functions $f_\rho$ and $g_\rho$ on $\cT^o$. These are defined by
\begin{align}
  f_\rho(\tau)&=\rho\,\{\sigma\ge\tau\}\label{e308}\\
  g_\rho(\tau)&=\int_\cT\a(\sigma\wedge\tau)\,d\rho(\sigma)\label{e309}.
\end{align}
We let $\cN$ and $\cP$ be the set of all functions of the form
$f_\rho$ and $g_\rho$, respectively.

The key point is now that the spaces $\cN$ and $\cP$ can be
characterized intrinsically, without any reference to measures. 
We shall refrain from giving 
the precise definitions here, but the functions in $\cN$ are 
\emph{left continuous} and of \emph{bounded variation}.
The \emph{total variation} $TV()$ provides $\cN$ with a norm.
As for the set $\cP$, its elements are called 
\emph{complex tree potentials}. 
They can be characterized by the fact that they have a well-defined
\emph{left derivative} at every point and this left-derivative
defines a function in $\cN$. 
This defines a natural isomorphism $\delta:\cP\to\cN$, which
enables us to define a norm on $\cP$ by $\|g\|=TV(\delta g)$.

The mappings $\rho\mapsto f_\rho$ and $\rho\mapsto g_\rho$
defined in~\eqref{e308}-\eqref{e309} can then be thought of as 
Banach space isomorphisms $\cM\to\cN$ and $\cM\to\cP$.
We denote their inverses by $d:\cN\to\cM$ 
and $\Delta:\cP\to\cM$. The mapping $\Delta$ is of particular
interest: we call it the \emph{Laplace operator} as it naturally
generalizes the Laplace operator both on the real line and on a 
simplicial tree.

\smallskip
Inside the three isometric Banach spaces $\cM$, $\cN$ and $\cP$ are
located three natural (positive) convex cones $\cM^+$, $\cN^+$ and
$\cP^+$.  The cone $\cM^+$ consists of all (finite) \emph{positive
Borel measures} on $\cT$, and $\cN^+$ of all \emph{left-continuous,
nonnegative, strongly decreasing} functions on $\cT^o$.  As before,
these properties have to be interpreted in the tree sense.  The cone
$\cP^+$ is of particular interest. It consists of all \emph{positive
tree potentials}, and can be characterized intrinsically, without
passing to $\cN$.
The positive tree potentials 
should be thought of as nonnegative functions that are concave in
a (strong) way that reflects the branching of the tree.
The isomorphisms $d:\cN\to\cM$ and $\Delta:\cM\to\cP$ then
restrict to bijections $d:\cN^+\to\cM^+$ and 
$\Delta:\cM^+\to\cP^+$. 

\smallskip
In addition to the topologies induced from the norms, the cones above
also carry \emph{weak topologies}, of great importance for
applications.  In the case of $\cM^+$, the weak topology is given by
the usual weak (or vague) topology on measures.  On $\cN^+$, it is
given by pointwise convergence outside a countable subset
of $\cT^o\setminus\{\tau_0\}$.
On the cone $\cP^+$ of positive tree potentials,
it is given by pointwise convergence on $\cT^o$.

The bijections $d:\cN^+\to\cM^+$ and $\Delta:\cM^+\to\cP^+$
then become homeomorphisms in the weak topologies. As the
subset of probability measures in $\cM^+$ is weakly compact,
the cone $\cP^+$ possesses strong compactness properties as well,
reminiscent of those of concave functions on the real line or
superharmonic functions on the unit disk. 
The analogy to the
latter two spaces is reinforced by the fact that the infimum
of any family of positive tree potentials remains a 
positive tree potential.

\smallskip
We also consider \emph{inner products} on the cones above.
On positive measures, the inner product is defined as the
bilinear extension of an inner product on Dirac masses, defined
as follows:
\begin{equation}\label{e310}
  \tau\cdot\tau':=\a(\tau\wedge\tau'),
  \quad\tau,\tau'\in\cT
\end{equation}
On the cone $\cN^+$, it is given by
\begin{equation}\label{e311}
  \langle f,f'\rangle 
  \= f(\tau_0)f'(\tau_0)
  +\int_{\cT^o}ff'\,d\l,
  \quad
  f,f'\in\cN^+,
\end{equation}
where $\l$ denotes the one-dimensional \emph{Hausdorff measure}
on $\cT^o$, induced by the parameterization $\a$.
Finally, the inner product is defined on the cone $\cP^+$ of
positive tree potentials by declaring 
$\langle g,g'\rangle=\langle\delta g,\delta g'\rangle$.
These inner products are preserved under the isomorphisms
of $\cM^+$, $\cN^+$ and $\cP^+$.

Integrability problems prevent us from extending these 
inner products to the full Banach spaces 
$\cM$, $\cN$ and $\cP$, but as we show, they are well defined 
on suitable subspaces $\cM_0$, $\cN_0$ and $\cP_0$.

\smallskip
Finally, let us comment on the role of the
parameterization $\a$. It does not appear in the definition
of the space $\cM$ of complex Borel measures, 
nor in the space $\cN$.
It does, however, appear in the definition
of complex tree potentials, \ie the space $\cP$.

Moreover, the parameterization is needed to define 
the inner products in all three cases (the Hausdorff measure
depends on the parameterization).

\smallskip
Throughout the chapter we shall assume that the parameterization is
increasing and takes value 1 at the root $\tau_0$. This is motivated
by the applications in Chapter~\ref{part-appli-analysis} where we work
with the valuative tree parameterized by skewness. 
Suitable adaptations should be made when working 
with other types of parameterizations. 
%
%
\subsection{Organization}\label{S307}
Let us end this section by giving a guide to the organization
of the remainder of the chapter. See the introduction to
each section for further details.

The weak topology on a tree, which was introduced in
Section~\ref{S302}, is analyzed in more detail in
Section~\ref{S27}. In particular we characterize connected open
sets and prove that a complete, parameterizable, nonmetric tree is
weakly compact.

We analyze (weak) Borel measures in Section~\ref{measure}.  Among
other things we show that every Borel measure is Radon (\ie regular),
supported on the completion of a countable union of finite subtrees,
and determined by its values on appropriate generalizations of
half-open intervals on the real line.

In Section~\ref{S410} we turn to functions of bounded variation.
They are defined and analyzed in much the same way as on the real
line (see \eg~\cite{folland}), 
but considerable care has to be taken because of the branching
in the tree. Here we define the space $\cN$, which is the first space
of functions to be identified with complex Borel measures.
This identification is proved in Section~\ref{S411}.

Then, in Section~\ref{S412}, we turn to complex tree potentials.
These constitute the second---and for applications most 
important---class of functions to be identified with complex
Borel measures. 
They are constructed by integrating functions in $\cN$ along segments.
The space $\cP$ of complex tree potentials is by definition
isomorphic to $\cN$, hence to $\cM$. The isomorphism between $\cP$ and
$\cM$ defines the Laplacian of a complex tree potential and is spelled
out in Section~\ref{S413}.

Atomic measures play an important role in the two applications
in Chapter~\ref{part-appli-analysis}. The corresponding functions
in $\cN$ and $\cP$ are easy to characterize; this is done
in Section~\ref{S414}.

Another important class of measures is given by the positive cone
$\cM^+\subset\cM$ of positive Borel measures. 
In Section~\ref{S415} we describe its preimages $\cN^+\subset\cN$
and $\cP^+\subset\cP$ under the isomorphism of $\cN$ and $\cP$
with $\cM$. Moreover, we show that the real-valued elements
of $\cM$, $\cN$ and $\cP$ admit canonical 
Jordan decompositions into differences of elements in 
$\cM^+$, $\cN^+$ and $\cP^+$, respectively.

The cone $\cP^+$ is very important in applications and
consists of the positive tree potentials. As opposed to
general complex tree potentials, these
can be easily characterized directly as functions on $\cT^o$,

On the cones $\cM^+$, $\cN^+$ and $\cP^+$ it is natural to
consider weak topologies, in addition to the restriction of
the norm topologies. The weak topologies are defined and analyzed
in Section~\ref{S308}. We show that the three cones are
homeomorphic in the weak topology, a fact that yields
compactness properties in $\cP^+$.

In Section~\ref{S409} we show that the three cones behave well
when passing to and from a complete subtree $\cS\subset\cT$.

Finally, in Section~\ref{S416}, we analyze inner products,
first on the cones $\cM^+$, $\cN^+$ and $\cP^+$, then on suitable 
complex subspaces of $\cM$, $\cN$ and $\cP$.
We shall use them in the next chapter.

\smallskip
The exact assumptions on the tree will vary somewhat from section to
section, but are stated in the respective introductions.
%
%
%
%
\section{More on the  weak topology}\label{S27}
Before studying measures and functions on trees, we
review the weak topology defined in Chapter~\ref{part3}
in more detail.

In this section, $\cT$  denotes a
nonmetric tree. Except if we mention it otherwise,
$\cT$ is supposed nonrooted.
%
%
\subsection{Definition}
Recall the definition of the 
weak topology\index{topology!weak tree} on $\cT$.
If $\vv\in T\tau$ is a tangent vector at a point 
$\tau\in\cT$, we define\index{$U(\vv)$ (weak open set)}
\begin{equation}\label{e450}
  U(\vv)=\{\sigma\in\cT\setminus\{\tau\}\ ;\ 
  \text{$\sigma$ represents $\vv$}\}.
\end{equation}
The weak topology is generated by the sets $U(\vv)$ in the sense
that the open sets are arbitrary unions of finite intersections of 
sets of the form $U(\vv)$.
%
%
\subsection{Basic properties}
First we state and prove some basic results.
\begin{lemma}
  The weak topology is Hausdorff.
\end{lemma}
\begin{proof}
  Suppose $\tau,\tau'$ are distinct points in $\cT$. 
  Pick a point $\tau''\in\,]\tau,\tau'[$ and let 
  $\vv$ and $\vv'$ be the tangent vectors at $\tau''$ 
  represented by $\tau$ and $\tau'$, respectively.
  Then $U(\vv)$ and $U(\vv')$ are disjoint open neighborhoods
  of $\tau$ and $\tau'$.
\end{proof}
\begin{lemma}\label{L310}
  If $\vv\in T\tau$, then the weak closure of 
  $U(\vv)$ equals $U(\vv)\cup\{\tau\}$.
\end{lemma}
\begin{proof}
  Any open set of the form $\bigcap_1^nU(\vv_i)$
  containing $\tau$ must intersect $U(\vv)$. Thus
  $\tau\in\overline{U(\vv)}$.
  On the other hand, if 
  $\sigma\not\in U(\vv)\cup\{\tau\}$, let $\ww$
  be the tangent vector at $\tau$ represented by $\sigma$.
  Then $\sigma\in U(\ww)$ and $U(\vv)\cap U(\ww)=\emptyset$
  so $\sigma\not\in\overline{U(\vv)}$.
\end{proof}
\begin{lemma}\label{L447}
  If $\cT$ has no branch points, then $\cT$ is homeomorphic
  to a real interval.
\end{lemma}
\begin{proof}
  Suppose for simplicity that $\cT$ is complete.
  Then $\cT$ is isomorphic as a nonmetric tree to the real 
  interval $I=[0,1]$. Moreover, the sets $U(\vv)$ correspond
  to intervals of the form $(x,1]$ or $[0,x)$ with 
  $0\le x\le 1$. Such intervals generate the topology
  on $I$ so $\cT$ and $I$ are homeomorphic.
\end{proof}
\begin{proposition}\label{P434}
  If $(\tau_k)_1^\infty$ is a sequence of points in $\cT$, then
  $\tau_k\to\tau\in\cT$ iff 
for all subsequence $(k_j)_1^\infty$, the segments
  $]\tau,\tau_{k_j}]$ have empty intersection.  In particular, if the
  points $\tau_k$ all represent distinct tangent vectors at $\tau$,
  then $\tau_k\to\tau$.
\end{proposition}
\begin{proof}
  If there is a subsequence $(k_j)$ such that
  the segments $]\tau,\tau_{k_j}]$ have a point 
  $\sigma\in\cT$ in common, let $\vv$ be the tangent vector at 
  $\sigma$ represented by $\tau$.
  Then $U(\vv)$ is an open neighborhood of $\tau$
  not containing any of the points $\tau_{k_j}$, so
  $\tau_{k_j}\not\to\tau$, implying $\tau_k\not\to\tau$.

  Conversely, suppose that there is no such subsequence.
  It suffices to show that for any open set of the form
  $U(\vv)$ containing $\tau$ we have
  $\tau_k\in U(\vv)$ for large $k$. 
  Here $\vv$ is a tangent vector 
  at some point $\sigma\ne\tau$ and represented by $\tau$.
  But this is clear, since 
  otherwise we could find infinitely many $k$
  with $\sigma\in[\tau,\tau_k]$.
  The proof is complete.
\end{proof}
%
%
\subsection{Subtrees} 
As we show next, the weak topology is well behaved with respect
to subtrees.
\begin{lemma}\label{L308}
  If $\cS$ is a subtree of $\cT$, then the weak
  topology on $\cS$ coincides with the topology
  on $\cS$ induced from the weak topology on $\cT$.
  In other words, the inclusion map 
  $\imath:\cS\to\cT$ is an embedding.
\end{lemma}
\begin{proof}
  The weak topology on $\cS$ is generated by subsets of 
  $\cS$ of the form $U_\cS(\vv)$ defined as in~\eqref{e450}, but 
  where $\vv$ is a tangent vector in $\cS$. 
  The induced topology is generated by sets of the form
  $U(\vv)\cap\cS$, where $\vv$ ranges over tangent vectors
  in $\cT$.
 
  On the one hand, every tangent vector in $\cS$ also
  defines a tangent vector in $\cT$ and 
  $U_\cS(\vv)=U(\vv)\cap\cS$. On the other hand, 
  if $\vv$ is a tangent vector at some point $\tau\in\cT$
  such that $U(\vv)\cap\cS\ne\emptyset$, then either 
  $\tau\in\cS$ and $\vv$ is a tangent vector in $\cS$;
  or $\tau\not\in\cS$ and $U(\vv)\cap\cS=\cS$.
\end{proof}
\begin{lemma}\label{L306}
  The relative closure of a subtree $\cS\subset\cT$ consists of
  $\cS$ and all points in $\cT$ that are ends in $\cS$.
  In particular, if $\cT$ is complete, then the closure of 
  $\cS$ in $\cT$ is equal to the completion of $\cS$.
\end{lemma}
\begin{proof}
  Consider a point $\tau\in\cT\setminus\cS$, and let $\ww$
  be the unique tangent vector $\ww$ at $\tau$ 
  such that $\cS\subset U(\ww)$.
  If $\tau$ is not an end in $\cS$, then
  we can find a point $\tau'\in U(\ww)$ such that
  $U(\vv)\cap\cS=\emptyset$, where $\vv$ is the tangent vector
  at $\tau'$ represented by $\tau$. 
  Thus $\tau$ does not belong to the closure of $\cS$ in $\cT$.
  On the other hand, if $\tau$ is an end in $\cS$, then it is easy to
  see that if $\vv$ is a tangent vector in $\cT$ and 
  $\tau\in U(\vv)$, then $U(\vv)\cap\cS\ne\emptyset$.
  This concludes the proof.
\end{proof}
If $\cS$ is a complete subtree of $\cT$, then we can define a
natural mapping $p_\cS:\cT\to\cS$ as follows. Pick a root
$\tau_0\in\cS$ of $\cT$ and set
$p_\cS(\tau)=\max([\tau_0,\tau]\cap\cS$); this does not
depend on the choice of $\tau_0$. Clearly $p_\cS=\id$ on $\cS$.
\begin{lemma}\label{L307}
  The mapping $p_\cS:\cT\to\cS$ is continuous, hence defines 
  a retraction of $\cT$ onto $\cS$.
\end{lemma}
\begin{proof}
  Any tangent vector $\vv$ in $\cS$ defines both an open set
  $U_\cS(\vv)$ in $\cS$ and an open set $U_\cT(\vv)$ in $\cT$.
  Clearly $U_\cT(\vv)=p_\cS^{-1}(U_\cS(\vv))$. 
  Thus $p_\cS$ is continuous.
\end{proof}
\begin{corollary}\label{C398}
  If $(\cT,\le)$ is a rooted, nonmetric tree, then the map  
  $\sigma\mapsto\tau\wedge\s$ is (weakly) continuous for any 
  $\tau\in\cT$.
\end{corollary}
\begin{proof}
  Apply Lemma~\ref{L307} to $\cS=[\tau_0,\tau]$,
  where $\tau_0$ is the root of $\cT$.
\end{proof}
%
%
\subsection{Connectedness}\label{SproofL446}
On the real line, the connected open sets are very easy to describe,
simply being the open intervals. On a general tree, the situation is
more complicated, yet understandable.  The following result will play
an important role later in the chapter. 
\begin{proposition}\label{L446}
  On a nonmetric tree $\cT$, any connected open subset
  $U\subsetneq\cT$ is a countable increasing union of (open) 
  sets of the form
  \begin{equation}\label{e400'}
    \bigcap_{i=1}^n U(\vv_i),\
    \vv_i\ \text{tangent vectors}. 
  \end{equation}
  When $\cT$ is rooted at $\tau_0\not\in U$, one can write
  the previous equation in the form
  \begin{equation}\label{e400}
    U(\vv)\setminus\bigcup_{i=1}^n\{\s\ge\tau_i\}; 
  \end{equation}
  where $\vv$ is a tangent vector, not represented by $\tau_0$, 
  at some $\tau\in U$, and $\tau_i\in U(\vv)$ for all $i$.
\end{proposition}
\begin{remark}\label{R446} 
  Note that any open set of the form $U(\vv)$ is also a countable
  increasing union of sets of the form $\{\s \ge \tau\}$ when $\vv$ is
  not represented by the root.
  Proposition~\ref{L446} can hence be rephrased by saying that
  any connected open set is a countable increasing union of sets of
  the form
  \begin{equation*}
    \{ \s \ge \tau , \,  \s \not\ge \tau_i , \, 1\le i \le n \},
  \end{equation*}    
  where $\tau, \tau_i$ are in $\cT$. These sets play a role 
  similar to that of half-open intervals on the real line 
  and will be explored systematically in Section~\ref{S407}.
\end{remark}
\begin{proof}[Proof of Proposition~\ref{L446}]
  As $U\ne\cT$, we may pick an element $\tau_0\in\cT \setminus U$, 
  and choose it as the root of $\cT$. 
  Define $\tau\=\inf U$. As $U$ is open, 
  $U\cap [\tau_0,\tau]$ is an open set in the segment $[\tau_0,\tau]$. 
  It is hence empty, and thus $\tau\notin U$.
  For simplicity of notation, let us change the root if necessary,
  so that $\tau=\tau_0$.

  Thus $U$ is a subset of $\cT\setminus\{\tau_0\}$, 
  which is a disjoint union of weak open sets $U(\vv)$ 
  where $\vv$ ranges over all tangent at $\tau_0$.
  As $U$ is connected,
  there exists a unique tangent vector $\vv$ such that 
  $U\subset U(\vv)$.
  As $U$ is open, $U(\vv)\setminus U$ is closed in $U(\vv)$, 

  Let $F$ be the set of all minimal elements in $U(\vv)\setminus U$. 
  More precisely, $\tau\in F$ iff 
  $\tau\in U(\vv)\setminus U$ and $]\tau_0,\tau[\,\subset U$.
  We claim that
  \begin{equation}\label{e410}
    U=U(\vv)\setminus\bigcup_{\tau\in F}\{\s\ge\tau\}. 
  \end{equation}
  Indeed, on the one hand, for any $\tau\in F$ we may cover
  $U$ by the disjoint open sets $U\cap\{\s>\tau\}$ and
  $U\cap U(\vv_0)$, where $\vv_0$ is the tangent vector at
  $\tau$ represented by $\tau_0$. As $U$ is connected, we conclude
  that $U\subset U(\vv)\setminus\bigcup_{\tau\in F}\{\s\ge\tau\}$.
  On the other hand, if $\s\in U(\vv)\setminus U$, then
  set $\tau=\inf\,]\tau_0,\sigma]\setminus U$. 
  As $U$ is open and nonempty we have $\tau\not\in U$. 
  Thus $\tau\in F$ and we obtain 
  $U\supset U(\vv)\setminus\bigcup_{\tau\in F}\{\s\ge\tau\}$.
  Hence~\eqref{e410} holds.

  Notice that $\inf F=\inf(U(\vv)\setminus U)$.
  Let us prove that $U$ can be written in the form~\eqref{e400} 
  under the assumption $\inf F>\tau_0$. 
  We shall explain later how to handle the other case.
  Set $\tau_1=\inf F$. 
  Then either $\tau_1\in F$ in which case 
  $U=U(\vv)\setminus\{\s\ge\tau_1\}$; 
  or $F\subset\{\s>\tau_1\}$. 
  In the latter case, we define $V$ to be the set of tangent vectors at
  $\tau_1$ which are represented by some element in $F$. 
  This set is finite, since otherwise we could find an infinite 
  sequence of elements in $F\subset U(\vv)\setminus U$ 
  representing distinct tangent vectors
  at $\tau_1$. By Proposition~\ref{P434}, and the fact that 
  $U(\vv)\setminus U$ is closed this would imply that 
  $\tau_1\not\in U$, a contradiction.
  For any $\ww \in V$, we define $\tau(\ww)\=\inf(F\cap U(\ww))$.
  Since $U(\vv)\setminus U$ is closed in $U(\vv)$ 
  and $\tau_1\in U(\vv)$, we have $\tau(\ww)>\tau_1$.
  Set $F_1\=\bigcup_{\ww\in V}\tau(\ww)$.
  
  Inductively we construct finite sets $F_k$ with the
  following properties: 
  \begin{itemize}
  \item[(i)]
    any element in $F_{k+1}$ dominates some element in $F_k$;
  \item[(ii)]
    any element in $F_k$ is dominated by some element in $F$;
  \item[(iii)]
    any element in $F$ dominates some element in $F_k$;
  \item[(iv)]
    any increasing sequence $(\tau_k)_1^\infty$ with 
    $\tau_k\in F_k$ for all $k$ converges to an element in $F$.
  \end{itemize}
  The set $F_1$ was already defined. The construction of $F_{k+1}$ in
  terms of $F_k$ is done as follows. Fix $\tau \in F_k$, and suppose
  $\tau\not\in F$. Let $V(\tau)$ be the set of tangent vectors at $\tau$
  represented by at least one point in $F$. By the same argument as
  before, $V(\tau)$ is a finite set, and we let
  \begin{equation*}
    F_{k+1}\=(F_k\cap F)\cup\bigcup_{\tau\in F_k\setminus F}
    \bigcup_{\ww\in V(\tau)}\inf(F\cap U(\ww)).
  \end{equation*}
  Clearly~(i)-(iii) hold.
  As for~(iv), pick an increasing sequence 
  $\tau_{k+1} \ge \tau _k \in F_k$.
  By construction, if $\tau_k=\tau_{k+1}$, then $\tau_k\in F$. We may
  hence suppose $\tau_{k+1} >\tau_k$ for all $k$.  As $\cT$ is complete,
  the sequence converges to a point $\tau\in\cT$.  
  We need to prove $\tau\in F$. 
  By construction, any $\tau_k$ is dominated by some $\s_k\in F$ 
  which does not represent the same tangent vector as $\tau$ at
  $\tau_k$. This implies that $\s_k\to\tau$ by Proposition~\ref{P434},
  hence $\tau\in U(\vv)\setminus U$. 
  But $\tau_k$ belongs to $U$ for all $k$ and the sequence 
  increases to $\tau$, so $\tau\in F$ and~(iv) holds.
  
  Define $U_k\=U(\vv)\setminus\bigcup_{\tau \in F_k}\{\s\ge\tau\}$. 
  This is a set of the form~\eqref{e400}.
  From~(i) it is clear that $U_{k+1}\supset U_k$. 
  By (ii) and~\eqref{e410}, $U_k\subset U$.
  We leave it to the reader to check that~(iii) and~(iv)
  imply $\bigcup U_k=U$.
  Hence we are done in the case $\inf F>\tau_0$.
  
  When $\inf F=\tau_0$, we fix a strictly decreasing sequence
  $(\tau_j)_1^\infty$ in $U$ converging to $\tau_0$,
  and apply the preceding result to 
  $U_j\=U(\vv_j)$, where $\vv_j$ is the tangent vector at $\tau_j$
  represented by $\tau_{j+1}$.
  This is a connected open set with 
  $\inf(F\cap U_j)>\inf U_j=\tau_j$.
  
  The proof is complete.
\end{proof}

%
%
\subsection{Compactness}
We end this section by proving 
\begin{proposition}\label{P104}
  Any parameterizable, complete, nonmetric tree $\cT$ is weakly
  compact.
\end{proposition}
\begin{proof}
  We will embed $\cT$ as a closed subspace of a product space.  Pick a
  root $\tau_0$ of $\cT$ and let $\le$ be the partial ordering rooted
  at $\tau_0$.  Also pick a parameterization $\alpha:\cT\to[0,1]$.
  Let $E$ be the set of functions $\chi:\cT\to[0,1]$ endowed with the
  product topology. This is a compact space, and convergence in this
  topology is given by pointwise convergence.
  Define a mapping $\jmath:\cT\to E$ by
  \begin{equation*}
    \jmath(\tau)(\sigma)=\alpha(\sigma\wedge\tau).
  \end{equation*}
  It is clear that $\jmath$ is injective:
  if $\tau\ne\tau'$, then without loss of generality
  $\tau\wedge\tau'<\tau=\tau\wedge\tau$,  and so 
  $\jmath(\tau)(\tau)\ne\jmath(\tau)(\tau')$.
  Let us show that $\jmath(\cT)$ is closed in $E$.
  For this, notice that if
  $\tau\in\cT$, then $\chi:=\jmath(\tau):\cT\to[0,1]$ 
  has the following two properties:
  \begin{itemize}
  \item[(i)]
    for any end $\sigma$ in $\cT$, the restriction of 
    $\chi$ to $[\tau_0,\sigma]$ is of the form
    \begin{equation*}
      \chi(\sigma')=\max\{\a(\sigma'),s\}
    \end{equation*} 
    for some $s\in[0,\a(\sigma)]$
    (we have $s=\a(\sigma\wedge\tau)$);
  \item[(ii)]
    $\chi(\sigma\wedge\sigma')=\inf\{\chi(\sigma),\chi(\sigma')\}$
    for any $\sigma,\sigma'\in\cT$.
  \end{itemize}
  Conversely, we claim that any function 
  $\chi:\cT\to[0,1]$ satisfying~(i) and~(ii) is of the form
  $\chi=\jmath(\tau)$ for some $\tau\in\cT$. 
  Indeed, let $\cS=\cS_\chi\subset\cT$ be the set of points 
  $\sigma$ such that $\chi$ is strictly increasing on 
  the segment $[\tau_0,\sigma]$. 
  If $\chi$ is nonconstant, then $\cS$ is nonempty by~(i) 
  and is a totally ordered set by~(ii). 
  Since $\cT$ is complete, $\cS$ has a maximal
  element $\tau$. It is then easy to check from~(i) and~(ii)
  that $\chi=\jmath(\tau)$. If $\chi$ is constant, 
  then $\chi=\jmath(\tau_0)$. This proves the claim.

 It is then clear that in $E$ both conditions~(i) and~(ii) define a
 closed set in the product topology, hence $\jmath(\cT)$ is closed.

  Thus $\jmath$ gives a bijection of $\cT$ onto the closed subset
  $\jmath(\cT)$ of $E$.
  We claim that $\jmath$ is a homeomorphism onto its image.
  As $\cT$ is Hausdorff and $E$ compact, it suffices to
  show that $\jmath$ is an open map onto its image. 
  By the definition of the weak topology it suffices to show that
  $\jmath(U(\vv))$ is relatively open in $\jmath(\cT)$ for every
  tangent vector $\vv$ at any point $\tau$ in $\cT$.
  If $\vv$ is represented by $\tau_0$, then
  \begin{equation*}
    U(\vv)
    =\{\sigma\in\cT\ ;\ \a(\sigma\wedge\tau)<\a(\tau)\}
    =\{\sigma\in\cT\ ;\ \jmath(\sigma)(\tau)<\a(\tau)\}
  \end{equation*}    
  and if $\vv$ is represented by some $\tau'>\tau$, then
  \begin{equation*}
    U(\vv)
    =\{\sigma\in\cT\ ;\ \a(\sigma\wedge\tau')>\a(\sigma\wedge\tau)\}
    =\{\sigma\in\cT\ ;\ \jmath(\sigma)(\tau)<\jmath(\sigma)(\tau')\}.
  \end{equation*}    
  In both cases it follows, after unwinding definitions, that
  $\jmath(U(\vv))$ is the intersection of an open set in $E$ 
  and $\jmath(\cT)$, hence is relatively open in $\jmath(\cT)$.
  This concludes the proof.
\end{proof}
%
%
%
%
\section{Borel measures}\label{measure}
In this section we let $\cT$ be a complete, nonmetric tree.
We also assume that $\cT$ is weakly compact: as we have seen,
this is the case if $\cT$ is parameterizable.
We shall study complex Borel measures on $\cT$.

Since a tree can be a quite ``large'' space in the sense that there is
no a priori bound on the amount of branching, we shall be very
detailed in our analysis. In particular we shall (temporarily)
distinguish between Borel measures and Radon measures.  Our general
reference is~\cite{folland}.
%
%
\subsection{Basic properties}
Let $\cB$\index{$\cB$ (Borel $\sigma$-algebra)} 
be the Borel $\sigma$-algebra on $\cT$, \ie the smallest
$\sigma$-algebra containing all (weakly) open sets in $\cT$.  A
\emph{complex Borel measure}\index{measure!complex Borel} is a
function $\rho:\cB\to\C$ satisfying $\rho(\emptyset)=0$ and
$\rho(\bigcup_1^\infty E_j)=\sum_1^\infty\rho(E_j)$ whenever
$(E_j)_1^\infty$ is a sequence of disjoint Borel sets; 
\emph{real} (or \emph{signed}) 
{Borel measures}\index{measure!real Borel}\index{measure!signed Borel}
are defined in the same way, with $\C$ replaced by 
$\R$.\footnote{We do not allow sets of measure $\pm\infty$.}  
A real Borel measure $\rho$ is
\emph{positive}\index{measure!positive Borel}
if $\rho(E)\ge0$ for all $E\in\cB$.  Obviously a
function $\rho:\cB\to\C$ is a complex Borel measure iff
$\mathrm{Re}\,\rho$ and $\mathrm{Im}\,\rho$ are real Borel measures.
Less trivial is the fact that a real Borel measure $\rho$ has a unique 
\emph{Jordan decomposition}\index{Jordan decomposition!of a measure} 
$\rho=\rho_+-\rho_-$, where $\rho_\pm$ are mutually singular
positive measures: there exist disjoint Borel sets $E_\pm$ with
$\rho_+(E_-)=\rho_-(E_+)=0$ and $\rho_+(\cT\setminus
E_+)=\rho_-(\cT\setminus E_-)=0$.  The measures $\rho_+$ and $\rho_-$
are called the \emph{positive and negative parts} of $\rho$.  These
decompositions allow us to reduce many questions about complex
measures to positive measures.
%
%
\subsection{Radon measures}
A positive measure $\rho$ is a\index{measure!Radon} 
\emph{Radon measure}\footnote{Some authors would call Radon measures 
\emph{regular Borel measures}.}
if
\begin{equation*}
  \rho(E)=\inf\{\, \rho(U)\ ;\ U\supset E,\, U\ \text{open}\ \}
\end{equation*}
for any Borel set $E$ and
\begin{equation*}
  \rho(U)=\sup\{\, \rho(K)\ ;\ K\subset U,\, K\ \text{compact}\ \}
\end{equation*}
for any open set $U$; the latter equality then
holds when $U$ is replaced by a general 
Borel set $E$~\cite[Proposition~7.5]{folland}.
A real Borel measure is said to be Radon iff
its positive and negative parts are Radon, and 
a complex measure is Radon iff its real and imaginary
parts are.

To any complex Borel measure $\rho$ is associated its
\emph{total variation measure} $|\rho|$. 
\index{total variation!measure} 
This is a positive measure which can be defined by
$|\rho|(E)=\sup\sum_1^n|\rho(E_i)|$, where the supremum is taken 
over collections of finitely many disjoint
Borel sets $E_1,\dots,E_n$ such that $E = \cup E_i$. 
We call $\Vert\rho\Vert:=|\rho|(\cT)\in\R_+$ the 
\emph{total variation}
\index{total variation!of a measure} of $\rho$.
If $\rho$ is a real measure, then $|\rho|=\rho_++\rho_-$,
where $\rho=\rho_+-\rho_-$ is the Jordan decomposition of $\rho$.
%
%
\subsection{Spaces of measures}
Let $C(\cT)$ be the space of continuous functions on $\cT$ 
with the topology of uniform convergence.
By the Riesz representation theorem~\cite[Theorem~7.17]{folland}, 
complex Radon measures
can be identified with continuous linear functionals 
on $C(\cT)$: to $\rho$ is associated the operator
$\varphi\mapsto\int\varphi\,d\rho$.
(In fact this operator is well-defined for any 
Borel measure---the point is that there is a canonical way of
associating a unique Radon measure to each continuous
linear functional.)
The total variation of $\rho$ is equal to its norm as a 
linear functional on $C(\cT)$.

On many topological spaces every finite complex Borel measure is Radon. 
This is the case for $\R^n$ and, more generally, 
any locally compact Hausdorff space in which every open set is 
a countable union of compact subsets~\cite[Theorem~7.8]{folland}.
While a tree may not fall into the latter category 
(if, for instance, the tangent space at some point is uncountable)
we still have
\begin{proposition}\label{L434}
  Every complex Borel measure on a weakly compact
  nonmetric tree is a Radon measure.
\end{proposition}
The proof is given in Section~\ref{S417} below. 
We shall consequently omit the adjective
``Radon'' in what follows, and we shall identify Borel measures
with continuous linear functionals on $C(\cT)$.

We let $\cM$\index{$\cM$ (complex Borel measures)} 
be the set of complex Borel measures on $\cT$.
Let us recall a few facts about $\cM$; they 
follow from the fact that $\cT$ is a
compact Hausdorff space (and that every Borel measure
is Radon), but do not use the tree structure of $\cT$.

There are two important topologies on $\cM$.  First, we have the
\emph{strong topology} induced by
the norm $\Vert\rho\Vert:=\sup\{|\int_\cT\varphi\,d\rho|\ ;\
\sup_\cT|\varphi|=1\}$.  
\index{topology!strong (on $\cM$)} 
Under this norm, $\cM$ is a Banach space.  
We shall also consider the 
\emph{weak topology},
defined in terms of convergence by $\rho_k\to\rho$ iff
$\int\varphi\,d\rho_k\to\int\varphi\,d\rho$ for all $\varphi\in
C(\cT)$.  This topology is Hausdorff.  The unit ball $\{\rho\in\cM\ ;\
\Vert\rho\Vert\le1\}$ is weakly compact.

We denote by 
$\cM^+$\index{$\cM^+$ (positive Borel measures)} 
the set of positive Borel measures on $\cV$. 
This is a closed subset of $\cM$ in both the
weak and strong topologies.
The subset of $\cM^+$ consisting of probability measures 
(\ie positive measure of total mass one) is weakly compact.

Atomic measures will play an important role in our study.  To
alleviate notation we shall identify $\tau\in\cT$ with the
corresponding point mass $\delta_\tau\in\cM^+$.  Indeed, the mapping
$\tau\mapsto\delta_\tau$ gives an embedding of $\cT$ as a weakly
closed subset of $\cM^+$ (see~\cite[III, \S 1, n.9, p.59]{Bou1} for
instance).
%
%
\subsection{The support of a measure}
A general tree can exhibit quite wild branching, 
but as the following result asserts, Borel measures 
live on reasonably well-behaved subsets.
\index{support!of a measure}
\begin{lemma}\label{supp-count}
  The support of any complex Borel measure $\rho$ on $\cT$ 
  is contained in the completion of a countable union of 
  finite subtrees.
  Moreover, the weak topology on the support is metrizable.
\end{lemma}
\begin{proof}
  We may suppose that $\rho$ is a positive measure of mass $1$.
  Consider the decreasing function $f:\cT\to[0,1]$ defined by
  $f(\tau)=\rho\{\sigma\ge\tau\}$.  The support of $\rho$ is included
  in the completion of the tree $\cS:=\{f>0\}$
  (see Lemma~\ref{L306}). 
  Let $\cS_n=\{f\ge n^{-1}\}$, $n\ge1$. 
  By construction the tree $\cS_n$ has at most $n+1$ ends, 
  hence at most $n+2$ branch points, so $\cS_n$ is finite.

  According to the Urysohn Metrization Theorem 
  (see~\cite[p.139]{folland}), any normal, 
  second countable space is metrizable. 
  Now the support $\supp\rho$ is closed, hence compact, 
  and therefore normal: see~\cite[Proposition 4.25]{folland}.  
  To conclude we need to show
  that it is also second countable, \ie that it admits a countable
  basis for its topology. To do so, write $\supp\rho$ as the 
  completion of an increasing union of finite trees 
  $\cS_n$ as before. Take a countable dense set $F$ in 
  $\bigcup\cS_n$. Then the collection of all open sets
  $\{\s>\tau\}$ over $\tau\in F$,
  and $\{\s>\tau,\ \s\not\ge\tau'\}$ over $\tau,\tau'\in F$,
  form a countable basis.
\end{proof}
\begin{remark}
  In general, a Borel measure is not necessarily supported on a
  countable union of finite trees. An example can be constructed as
  follows.
  
  Consider the Cantor set $X\=\{0,1\}^{\N}$ endowed with the
  (ultra-)metric $d(x,x')=\max 2^{-i}|x_i-x'_i|$ 
  where $x=(x_i)_{i\ge0}$, $x'=(x'_i)_{i\ge0}$. 
  Following Section~\ref{def-ultra}, define $\cT$ 
  to be the quotient of $X\times[0,1]$ by the
  equivalence relation $(x,t)\sim(x',t')$ iff $t=t'<d(x,x')$.

  On the Cantor set $X$ there exists a natural probability 
  measure $\rho_X$ giving measure $2^{-n}$ to any cylinder
  $\{x\ ;\ x_1,\dots,x_n\ \text{are fixed}\}$. 
  The pushforward of $\rho_X$ by the inclusion map
  $X\ni x\mapsto(x,1)\in\cT$ is a probability measure on $\cT$ 
  which is not supported on a countable union of finite trees.
\end{remark}
%
%
\subsection{A generating algebra}\label{S407}
Next we generalize the fact that on the real line, 
complex Borel measures are determined by their
values on half-open intervals.

Define collections $\cE$ and $\cA$ of subsets of $\cT$ as follows.
First, $\cE$\index{$\cE$ (elementary family)} consists of $\cT$ and
all sets of the form $\{\sigma\in\cT\ ;\ \sigma\ge\tau,
\sigma\not\ge\tau_i, 1\le i\le n\}$ where $\tau,\tau_i\in\cT$.  It is
clear that the complement of an element in $\cE$ is a finite disjoint
union of elements in $\cE$.  Hence $\cE$ is an \emph{elementary
family}\index{elementary family} in the sense
of~\cite[p.~22]{folland}.  Second, $\cA$\index{$\cA$ (algebra)} is the
set of finite disjoint unions of elements in $\cE$.  Then $\cA$ is an
\emph{algebra}\index{algebra (of Borel sets)} in the sense that $\cA$
is closed under complements and finite unions and intersections
(see~\cite[Proposition~1.7]{folland}).

Still following~\cite{folland}, we define a 
\emph{premeasure}\index{premeasure} (on $\cA$) to be a function 
$\rho:\cA\to[0,\infty[$ such that 
$\rho(\emptyset)=0$ and 
$\rho(\bigcup_1^\infty E_j)=\sum_1^\infty\rho(E_j)$
whenever $(E_j)_1^\infty$ is a sequence of disjoint
sets in $\cA$ such that $\bigcup_1^\infty E_j\in\cA$.
\begin{lemma}\label{L442}
  Any premeasure on $\cA$ has a unique extension to a (positive)
  Borel measure on $\cB$.
\end{lemma}
\begin{proof}
  This follows from Theorem~1.14 in~\cite{folland}.
\end{proof}
Recall from Remark~\ref{R446} that any connected open sets is a countable
increasing union of elements in $\cE$.  Another reason why the
elementary family $\cE$ is important lies in the following fact. Any
set $E\in \cE$ is a countable intersection of sets of the form
$U(\vv)\setminus\bigcup_1^nU(\vv_i)$, where the $\vv_i$'s are tangent
vectors, and $U(\vv_i)$'s are disjoint subsets of $U(\vv)$.
Conversely, any set $U(\vv)\setminus\bigcup_1^nU(\vv_i)$ as above, is
a countable increasing union of elements in $\cE$.

From this and Lemma~\ref{L442}, one deduces
\begin{lemma}\label{L448}
Let $\rho,\rho'$ be any two Borel measures on $\cT$.  If $\rho (U(\vv))
= \rho'(U(\vv))$ for any tangent vector $\vv$, then $\rho = \rho'$.
\end{lemma}

%
%
\subsection{Every complex Borel measure is Radon}\label{S417}
Finally we prove that every complex Borel measure is Radon:
\begin{proof}[Proof of Proposition~\ref{L434}]
  It suffices to show that every positive Borel measure $\rho$ on 
  $\cT$ is Radon.
  Following the proof of Theorem~7.8 in~\cite{folland} we
  can associate a positive Radon measure $\rho'$ to $\rho$ as follows:
  \begin{equation*}
    \rho'(U)
    :=\sup\left\{\int_\cT\varphi\,d\rho\ ;\, 
    \varphi\in C(\cT),\,  0\le\varphi\le1,\,  \supp\varphi\subset U\right\},
  \end{equation*}  
  for $U$ open and 
  $\rho'(E):=\inf\{\rho'(U) ;\, E\subset U,\, U\, \text{open}\}$ for
  any Borel set $E$.
  We claim that $\rho'=\rho$. 
  Consider a tangent vector $\vv$ at a point $\tau\in\cT$.
  First assume that $\vv$ is not represented by the root $\tau_0$.
  Pick a sequence $(\tau_k)_1^\infty$ in $\cT$, decreasing 
  to $\tau$, and such that $\tau_k$ represents $\vv$ for all $k$. 
  Choose nondecreasing continuous functions 
  $\chi_k:[\tau,\tau_k]\to[0,1]$ such that $\chi_k(\tau)=0$ and
  $\chi_k(\tau_k)=1$. Then define $\varphi_k:\cT\to[0,1]$ by
  $\varphi_k(\sigma) \= \chi_k (\sigma \wedge \tau_k)$ if $\s\ge \tau$, 
  and $\varphi_k(\s) =0$ otherwise. 
  Using Corollary~\ref{C398}, we see that
  $\varphi_k$ is weakly continuous,
  Moreover, $0\le\varphi_k\le1$, $\supp\varphi_k\subset U$,
  and $\varphi_k$ increases pointwise to the characteristic
  function of $U(\vv)$ as $k\to\infty$.
  Monotone convergence gives
  \begin{equation*}
    \rho(U(\vv))
    =\lim_{k\to\infty}\int_\cT\varphi_k\,d\rho
    =\rho'(U(\vv)).
  \end{equation*}
  A similar argument shows that $\rho(U(\vv))=\rho'(U(\vv))$
  also when $\vv$ is represented by $\tau_0$.
  In view of Lemma~\ref{L448} we conclude that $\rho=\rho'$.
  Hence $\rho$ is a Radon measure.
\end{proof}
%
%
%
%
\section{Functions of bounded variation}\label{S410}
Our objective now is to define a space $\cN$
of functions that can be naturally identified with the space
$\cM$ of complex Borel measures.
The elements in $\cN$ are functions of bounded variation. 
These functions are defined and studied in a way similarly to
the classical case of the real line, but some details turn
out to be trickier because of the branching in the tree.

Unless otherwise stated, 
$(\cT,\le)$ will denote a rooted, nonmetric tree 
(not necessarily complete) and $\tau_0$ its root.
Recall that, by convention, the ends of $\cT$ are the
maximal elements of $\cT$. 
In particular, the root $\tau_0$ is not an end.
%
%
\subsection{Definitions}
Consider a function $f:\cT\to\C$. 
Let us define what it means for $f$
to be of bounded variation, taking into account the
(possible) branching in $\cT$.
If $F$ is a finite subset of $\cT$, let 
$F_{\max}$ be the set of maximal elements of $F$
and write $\sigma\succ\tau$ when $\sigma,\tau\in F$ and
$\sigma$ is an immediate successor of $\tau$ in $F$, \ie
$\sigma>\tau$ and there is no $\sigma'\in F$ with 
$\tau<\sigma'<\sigma$.
Given a function $f:\cT\to\C$ and a finite set $F$, define
the \emph{variation of $f$ on $F$} 
by\index{$V(f;F)$ (variation of $f$ on $F$)}
\begin{equation}\label{e443}
  V(f;F)
  :=\sum_{\tau\in F\setminus F_{\max}}
  \left|f(\tau)-\sum_{\sigma\succ\tau}f(\sigma)\right|;
\end{equation}
if $F_{\max}=F$ then we declare $V(f;F)=0$.
\begin{remark}\label{R414}
  In contrast to the unbranched case, it can happen 
  $V(f;F)>V(f;F')$ when $F\subset F'$.
  Indeed, let $\cT$ be a rooted, nonmetric tree
  having exactly two ends $\tau_\pm$.
  Let $F=\{\tau_0,\tau_+\}$ and $F'=\{\tau_0,\tau_+,\tau_-\}$.
  Define $f:\cT\to\C$ by $f=0$ on $\cT\setminus\{\tau_+,\tau_-\}$ and
  $f(\tau_\pm)=\pm1$.
  Then $V(f;F')=0$ but $V(f;F)=1$.
\end{remark}
However, we have 
\begin{lemma}\label{L436}
  Let $(\cT,\le)$ be a nonmetric tree and let
  $F\subset F'$ be finite subsets of $\cT$ such that
  all the maximal points in $F'$ belong to $F$. 
  Then $V(f;F)\le V(f;F')$.
\end{lemma}
\begin{proof}
  By induction it suffices to consider the case when 
  $F'=F\cup\{\tau'\}$, where $\tau'$ is 
  not a maximal element in $F'$.
  First assume $\tau'$ is not a minimal element of $F'$. 
  Then $\tau'$ has an immediate predecessor $\tau\in F$.
  Let $E'$ ($E$) be the set of immediate successors 
  $\sigma$ of $\tau$ in $F$ such that 
  $\tau'\not<\sigma$ ($\tau'<\sigma$).
  Then we have
  \begin{equation*}
    V(f;F')-V(f;F)
    =\left|f(\tau)-f(\tau')-\sum_{E'}f\right|
    +\left|f(\tau')-\sum_{E}f\right|
    -\left|f(\tau)-\sum_{E\cup E'}f\right|
  \end{equation*}
  which is nonnegative by the triangle inequality.
  The case when $\tau'$ is a minimal element of 
  $F'$ is similar, but easier, and left to the reader.
\end{proof}
\begin{definition}
  A function $f:\cT\to\C$ is of \emph{bounded
  variation}\index{variation!bounded} if
  \begin{equation}\label{e442}
    TV(f):=\sup_F V(f;F)<\infty.
  \end{equation}
  We call $TV(f)$ the \emph{total variation}\index{variation!total} of
  $f$ (on $\cT$).
\end{definition}
Clearly a function $f:\cT\to\C$ is of bounded variation iff 
$\mathrm{Re}\,f$ and $\mathrm{Im}\,f$ are.
Notice that a constant function is \emph{not} necessarily of bounded
variation. In fact, the constant function $1$ is of bounded variation
iff $\cT$ is finite.

If $f:\cT\to\C$ has bounded variation, 
then so has the restriction of $f$ to any subtree 
$\cS$ and $TV(f|_\cS)\le TV(f)$. 

If $\cT$ is totally ordered (\ie has no branch points),
then we have $TV(f)=\sup\sum_1^{n-1}|f(\tau_{i+1})-f(\tau_i)|$,
where the supremum is taken over all finite, strictly
increasing sequences $(\tau_i)_1^n$ in $\cT$.
In particular, if $\cT$ is a real interval, 
then the definition above 
coincides with the classical definition.
%
%
\subsection{Decomposition}
On the real line, any bounded decreasing function is of bounded
variation, and it is a standard result~\cite[Theorem~3.27]{folland}
that any real-valued function of bounded variation is in fact a
difference of two nonnegative, decreasing
functions.\footnote{Folland~\cite{folland} works with
\emph{increasing} functions but the theory is completely analogous.}

In order for these results to generalize we have to
interpret ``bounded'' and ``decreasing'' in a way 
that reflects the branching.
Let us declare a finite subset $E$ of $\cT$
to be \emph{without relations} if 
$\sigma\not<\tau$ whenever $\sigma,\tau\in E$.
\begin{definition}\label{D412}
  We say that 
  \begin{itemize}
  \item[(i)] a function $f:\cT\to\C$ is 
    \emph{strongly bounded}\index{strongly!bounded function} 
    if $\sup_E|\sum_Ef|<\infty$,
    where the supremum is taken over finite sets $E$ without
    relations;
  \item[(ii)] a function $f:\cT\to\R$ is 
    \emph{strongly decreasing}\index{strongly!decreasing function} 
    if $f(\tau)\ge\sum_Ef$ whenever $E$ is a finite set 
    without relations and $\tau<\sigma$ for all $\sigma\in E$.
  \end{itemize}
\end{definition}
\begin{remark}\label{R417}
  The condition in~(i) is equivalent to
  $\sup_E\sum_E|f|<\infty$: see~\cite[Lemma~6.3]{rudin}.
  Also, any nonnegative strongly decreasing function is
  clearly strongly bounded.
\end{remark}
\begin{lemma}\label{L431}
  The following properties hold:
  \begin{itemize}
  \item[(i)]
    any function $f:\cT\to\C$ of bounded variation is
    strongly bounded;
  \item[(ii)]
    any strongly bounded and strongly decreasing function 
    $f:\cT\to\R$ is of bounded variation.
  \end{itemize}
\end{lemma}  
The proof is given in Section~\ref{S36} below.
\begin{corollary}\label{C408}
  The support of a function of bounded variation 
  (\ie the locus where it is nonzero) is
  contained in the completion of a countable union 
  of finite subtrees.
\end{corollary}
\begin{proof}
  If $f$ is has bounded variation, Lemma~\ref{L431}
  and Remark~\ref{R417} give us $C>0$ such that
  $\sum_E|f|<C$ for any finite set $E$ without relations.
  Thus the smallest subtree of $\cT$
  containing the set $\{|f|>n^{-1}\}$ has at most
  $Cn$ ends.
\end{proof}

Next we introduce an auxiliary function that plays the same role
as the total variation measure associated to a complex
Borel measure.
\begin{proposition}\label{L440}
  Suppose $\cT$ is complete, and consider a function
  $f:\cT\to\R$ of bounded variation that
  vanishes on $\cT\setminus\cT^o$. 
  Define the 
  \emph{total variation function}
  \index{total variation!function}
  $T_f:\cT\to[0,\infty[$\index{$T_f$ (total variation function)} 
  by
  \begin{equation}\label{e454}
    T_f(\tau)=TV(f|_{\{\sigma\ge\tau\}}).
  \end{equation}
  Then the three functions $T_f$, $T_f+f$ and $T_f-f$ are
  nonnegative, strongly decreasing functions on $\cT$ and 
  vanish identically on $\cT\setminus\cT^o$.
\end{proposition}
As $f=\frac12 (T_f+f)-\frac12 (T_f -f)$, we obtain
the following decomposition result: any function of
bounded variation is the difference of two functions in the convex
cone of nonnegative, strongly decreasing functions on $\cT$, vanishing
on $\cT\setminus\cT^o$.

\begin{proof}[Proof of Proposition~\ref{L440}]
  Let us start with the function $T_f+f$.  It is clearly nonnegative
  and nonincreasing.  Our conventions imply that $T_f+f=0$ on
  $\cT\setminus\cT^o$.  Let us show that $T_f+f$ is strongly
  decreasing.  Consider any finite set $E\subset\cT$ without relations
  and $\tau\in\cT$ such that $\tau<\sigma$ for every $\sigma\in E$.
  Let $\varepsilon>0$ and for every $\sigma\in E$ pick a finite subset
  $F_\sigma\subset\{\sigma'\ge\sigma\}$ such that $T_f(\sigma)\le
  V(f;F_\sigma)+\e/n$, where $n$ is the number of elements of $E$.
  Set $F=\{\tau\}\cup\bigcup_{\sigma\in E}F_\sigma$.  Then we have
  \begin{equation*}
    T_f(\tau)+ f(\tau) -\sum_E (T_f +f)
    \ge V(f;F)-\sum_{\sigma\in E}V(f;F_\sigma)-\e
    =-\e.
  \end{equation*}
  Letting $\e\to0$ we conclude that the function
  $T_f+f$ is strongly decreasing.
  
  As $T_{-f}=T_f$, this also shows that $T_f-f$ is strongly decreasing
  and non-negative. Finally $2T_f = (T_f +f) + (T_f-f)$, thus $T_f$ is
  also strongly decreasing and non-negative.  This completes the
  proof.
\end{proof}

%
%
\subsection{Limits and continuity}
If $\vv$ is a tangent vector at some point $\tau\in\cT$ then we
say that a sequence $(\tau_k)_1^\infty$ 
\emph{converges to $\tau$ along $\vv$} if $\tau_k$ represent
$\vv$ for all $k$ and $]\tau,\tau_k]$ form a decreasing sequence
of segments with empty intersection.

If $f$ is a complex-valued function on $\cT$ and $\vv$ is a tangent
vector then we say that $f$ has a 
\emph{limit along $\vv$}
\index{tangent vector!limit along}
if there exists a number $a\in\C$ such that
$f(\tau_k)\to a$ for all sequences $(\tau_k)_1^\infty$ converging to
$\tau$ along $\vv$; we then write $f(\vv)=a$.  If $\tau\ne\tau_0$,
then we define $f(\tau-)=f(\vv)$, where $\vv$ is the tangent vector at
$\tau$ represented by $\tau_0$ (assuming $f(\vv)$ exists).  
For convenience we also define $f(\tau_0-)=f(\tau_0)$.  Let us say that a
function $f:\cT\to\C$ is 
\emph{left continuous} at $\tau\ne\tau_0$ if
$f(\tau-)=f(\tau)$.  
A function $f$ is 
\emph{left continuous}\index{left continuity (in a tree)} 
if it is left continuous at every point $\tau\ne\tau_0$.
\begin{lemma}\label{L433}
  If $f:\cT\to\C$ has bounded variation, then 
  $f(\vv)$ exists for every tangent vector $\vv$.
\end{lemma}
\begin{proof}
  If $\tau\in\cT$, then the restriction of $f$ to the segment
  $[\tau_0,\tau]$ has bounded variation, allowing us to
  invoke the corresponding result on the real 
  line, see~\cite[Theorem~3.27]{folland}.
\end{proof}

Next we wish to introduce a number measuring roughly the 
discontinuity at a point 
of a function of bounded variation, in a way that reflects the tree
structure. In the case of the real line, this number equals the
absolute value of the difference between the left and right 
limits at a point, and can be used to prove that a function of 
bounded variation
has at most countably many discontinuities.
Note that on a general tree, a function of bounded variation may have 
uncountably many points of discontinuity in the usual sense.
\begin{Prop-def}\label{P444}
  Suppose $f:\cT\to\C$ has bounded variation. Then for each
  $\tau\in\cT$, the series $\sum_{\vv\in T\tau}f(\vv)$ is absolutely
  convergent. In particular, there are at most countably many tangent
  vectors $\vv$ at $\tau$ for which $f(\vv)\ne0$.
  
  We can hence set
  \index{$d(\tau;f)$ (atom of $df$ on $\tau$)}
  \begin{equation}\label{e444}
    d(\tau;f):=f(\tau-)-\sum_\vv f(\vv)
  \end{equation}
  where the sum is over all tangent vectors $\vv$ at $\tau$
  not represented by $\tau_0$.
\end{Prop-def}
Recall that we defined $f(\tau_0-)=f(\tau_0)$.
\begin{proof}
  Fix $\tau\in\cT$ and pick an arbitrary finite set $V$ of tangent
  vectors at $\tau$, none represented by the root $\tau_0$. Fix $\e
  >0$. For each of these vectors, choose a representing element $\s > \tau$,
  sufficiently close to $\tau$ so that $| f(\vv) - f(\s)| \le \e$.  We
  have $\sum | f(\vv) | \le \sum |f(\s)| + \e \#V $.  The collection of
  $\{ \s\}$ defines a set without relation, and $f$ is strongly bounded,
  thus $\sum | f(\vv) | \le B + \e \#V $, where $B>0$ is independent on
  the choices of tangent vectors (and on $\tau$ also).  By letting
  $\e$ tend to zero, we conclude $\sum | f(\vv) | \le B$.  This implies
  that $\sum_{T\tau} | f(\vv) |$ is absolutely convergent. 
\end{proof}
\begin{lemma}\label{L432}
  Suppose $f:\cT\to\C$ has bounded variation. Then, we have
  $d(\tau;f)=0$ for all but countably many points $\tau$, and the
  series $\sum_{\tau\in\cT}d(\tau;f)$ is absolutely convergent; in fact
  $\sum_{\tau\in\cT}|d(\tau;f)|\le TV(f)$.
\end{lemma}
\begin{proof}
  Fix $\e>0$ and 
  consider any finite subset $Z\subset\cT$ such that
  $d(\tau;f)\ne0$ for all $\tau\in Z$. 
  For each $\tau\in Z$, pick finitely many points
  $\tau_1,\dots,\tau_n$ close to $\tau$
  with the following properties:
  $\tau_1\le\tau$ and $\tau_1<\tau$ unless $\tau=\tau_0$;
  $\tau_i>\tau$ for $2\le i\le n$ and these $\tau_i$'s 
  represent distinct tangent vectors at $\tau$;
  $|f(\tau_1)-\sum_2^nf(\tau_i)-d(\tau;f)|<\e|d(\tau;f)|$.
  We may assume that the sets $\{\tau_1,\dots,\tau_n\}$
  are disjoint as $\tau$ varies over $Z$.
  Let $F$ be their union: then
  $\sum_{\tau\in Z}|d(\tau;f)|<(1-\e)^{-1}V(f;F)\le(1-\e)^{-1}TV(f)$.
  We conclude by letting $\e\to0$.
\end{proof}

%
%
\subsection{The space $\cN$}\label{S418}
We are now in position to define the first space of 
functions that can be identified with complex Borel
measures.
Assume that the rooted, nonmetric tree $(\cT,\le)$ 
is \emph{complete}.
Let $\cT^o$ be the set of nonmaximal elements of $\cT$
(\ie the points of $\cT$ that are not ends).
\begin{definition}
  We let $\cN$\index{$\cN$ (functions of bounded variation)} 
  be the set of functions 
  $f:\cT\to\C$ of bounded variation such that
  $f$ is left continuous at any point in $\cT^o$
  and $f=0$ on $\cT\setminus\cT^o$.
\end{definition}
\begin{definition}
  We let 
  $\cN^+$\index{$\cN^+$ (nonnegative strongly decreasing functions)} 
  be the set of nonnegative,
  strongly decreasing functions on $\cT$, which are
  left continuous on $\cT^o$ and vanish on $\cT\setminus\cT^o$.
\end{definition}
Notice that $\cN^+\subset\cN$, thanks to Lemma~\ref{L431}, and that
$\cN^+$ is a convex cone.  We shall prove the following crucial
decomposition result (see Proposition~\ref{L440}).
\begin{proposition}\label{P441}
  For any $f\in\cN$, the function 
  $T_f (\tau)\=TV(f|_{\{\s\ge\tau\}})$ 
  is left continuous, and the functions 
  $T_f$, $T_f+f$ and $T_f-f$ belong to $\cN^+$.

  We can thus write any real-valued function in $\cN$ as the difference
  of two functions in 
  $\cN^+$, $f=\half(T_f+f)-\half(T_f-f)$.
  This decomposition is called the 
  \emph{Jordan decomposition}\index{Jordan decomposition!of a function} 
  of $f$.
\end{proposition}
We shall later see that this decomposition $f=f_+-f_-$,
$f_\pm\in\cN^+$ is unique if one imposes $TV(f)=TV(f_+)+TV(f_-)$
(see Proposition~\ref{P419}).

We can equivalently view $\cN$ as the space of 
left continuous functions of bounded variation 
\emph{on the subtree $\cT^o$}.
In fact, this point of view is natural in many situations
when $\cT$ is the valuative tree $\cV$ and $\cT^o$ the
subtree $\cVqm$ of quasimonomial valuations.
The equivalence follows from
\begin{lemma}
  Let $f:\cT^o\to\C$ be a function of bounded variation.
  Extend $f$ to $\cT$ by declaring $f|_{\cT\setminus\cT^o}=0$.
  Then $f:\cT\to\C$ also has bounded variation.
\end{lemma}
\begin{proof}
  Pick any finite subset $F'\subset\cT$.
  Set $F:=F'\cap\cT^o$ 
  and let $E$ be the set of points in $F$ that are
  not maximal elements of $F'$ but have no
  immediate successors in $F$.
  Then $V(f;F')=V(f;F)+\sum_E|f|$. 
  Here the first term is bounded by $TV(f|_{\cT^o})$ 
  and the second term is uniformly bounded in view
  of Lemma~\ref{L431}~(i).
\end{proof}
In the sequel we shall often consider the
functions in $\cN$ as defined only on $\cT^o$.
However, it is important to notice that the norm 
$TV(f)$ of $f\in\cN$ is defined as the total variation
of $f$ as a function on the complete tree $\cT$.

Notice that the functions in $\cN$ are not necessarily left continuous 
at the ends of $\cT$ but they do have left limits there in view of
Lemma~\ref{L433}.
\begin{proposition}
  $(\cN,TV)$ is a normed, complex vector space\index{topology!strong
  (on $\cN$)}.
\end{proposition}

\begin{proof}
  That $\cN$ is a complex vector space is obvious, as are the facts
  that $TV(\lambda f)=\lambda TV(f)$ for $\lambda\in\C$, $f\in\cN$,
  and $TV(f)=0$ iff $f=0$. Finally, if $f_1,f_2\in\cN$, then 
  for any finite subset $F\subset\cT$ we have 
  $V(f_1+f_2;F)\le V(f_1;F)+V(f_2;F)$ by the triangle inequality.
  Hence $TV(f_1+f_2)\le TV(f_1)+TV(f_2)$.
\end{proof}
We shall later see that the norm $TV(\cdot)$ is complete, 
so that $\cN$ is in fact a Banach space (see Corollary~\ref{C409}).

%
%
\subsection{Finite trees}
The total variation $TV(f)$ of a function $f:\cT\to\C$
captures both the variation of $f$ along segments
of $\cT$ and the discontinuities of $f$ at branch points.
This assertion can be made precise as long as we deal with
finite, complete trees and functions in the space $\cN$:
\begin{proposition}\label{L435}
  Assume that $\cT$ is complete and finite.
  Let $B\subset\cT^o$ be any finite subset containing all the
  branch points of $\cT$ and let $\cI$ be the set of
  connected components of $\cT\setminus B$.
  Then for any $f\in\cN$ we have
  \begin{equation}\label{e447}
    TV(f)
    =\sum_{\tau\in B}|d(\tau;f)|
    +\sum_{I\in\cI}TV(f|_I).
  \end{equation}
\end{proposition}
\begin{remark}
  It is important that $f$ be zero on the ends of $\cT$.
  Indeed, let $\cT$ be a rooted, nonmetric tree
  having two ends $\tau_\pm$ and one branch point $\tau'$.
  Define $f:\cT\to\R$ by $f|_{[\tau_0,\tau']}=0$,
  and $f|_{]\tau',\tau_\pm]}=\pm1$.
  Then $TV(f)=1$, but if $B=\{\tau'\}$, then 
  the right hand side of~\eqref{e447} is zero.
\end{remark}
In view of this remark, the proof of this proposition,
which is given below, 
is trickier than might be expected. 
%
%
\subsection{Proofs}\label{S36}
We end this section by supplying the proofs of various assertions
above. 
\begin{proof}[Proof of Lemma~\ref{L431}]
  First suppose $f:\cT\to\C$ has bounded variation.
  Pick any finite subset $E$ without relations and
  not containing the root $\tau_0$, and set $F=E\cup\{\tau_0\}$.
  Then $|\sum_Ef|\le V(f;F)+|f(\tau_0)|\le TV(f)+|f(\tau_0)|<\infty$
  so $f$ is strongly bounded, proving~(i).

  As for~(ii), suppose $f:\cT\to\R$ is strongly bounded and 
  strongly decreasing. 
  Set $C=\inf_E\sum_Ef$, where $E$ runs over all
  finite subsets without relations.
  Then $C>-\infty$ as $f$ is strongly bounded.
  We claim that $TV(f)=f(\tau_0)-C$.
  To see this, first pick any $\e>0$ and choose $E$ finite
  without relations such that $\sum_Ef\le C+\e$.
  As $f$ is strongly decreasing, we may assume that
  $\tau_0\not\in E$. Set $F=E\cup\{\tau_0\}$.
  Then $TV(f)\ge V(f;F)\ge f(\tau_0)-C-\e$, hence 
  $TV(f)\ge f(\tau_0)-C$ when $\e\to0$.
  For the converse inequality,   
  consider any finite subset $F\subset\cT$.
  Let us show that $V(f;F)\le f(\tau_0)-C$.
  We may assume that $\tau_0\in F$ as 
  $V(f;F)\le V(f;F\cup\{\tau_0\})$. 
  Let $E$ be the set of maximal elements in $F$.
  Then $E$ is without relations and 
  \begin{align*}
    \sum_{\tau\in F\setminus F_{\max}}
    \left|f(\tau)-\sum_{F_\tau}f\right|
    =\sum_{\tau\in F\setminus F_{\max}}
    \left(f(\tau)-\sum_{F_\tau}f\right)
    =f(\tau_0)-\sum_Ef
    \le f(\tau_0)-C,
  \end{align*}
  where $F_\tau$ stands for the set of immediate successor of $\tau$ in $F$.
  This completes the proof.
\end{proof}
\begin{proof}[Proof of Proposition~\ref{P441}]
  In view of Proposition~\ref{L440}, 
  it suffices to show that
  the function $T_f$ defined in~\eqref{e454}
  is left continuous.

  For this, pick any $\tau\in\cT$, $\tau\ne\tau_0$.
  Suppose $T_f(\tau-)>T_f(\tau)$ and fix $\e>0$
  with $0<3\e<T_f(\tau-)-T_f(\tau)$.
  Pick $\tau'<\tau$ such that 
  $|f(\sigma)-f(\tau)|\le\e$ and
  $T_f(\sigma)\ge T_f(\tau)+3\e$ for all $\sigma\in[\tau',\tau[$.

  We shall construct a strictly increasing sequence $\tau_n \in
  [\tau',\tau[$ converging to $\tau$, and finite sets $F_n \subset\{
  \sigma\ ;\ \sigma>\tau_n, \sigma\not>\tau_{n+1}\}$, such that
  $V(f;F_n)\ge\e$ and $\tau_{n+1}\in F_n$.
  This implies $V(f ; \bigcup_1^n F_i) \ge \sum_1^n V(f ; F_i)$, thus
  $TV(f)\ge V(f;\bigcup_1^nF_i)\ge\sum_1^nV(f;F_i)\ge n\e$,
  which gives a contradiction as $n\to\infty$.
  
  We define $\tau_n$ and $F_n$ by induction. First set $\tau_1 \= \tau'$.
  Given $\tau_n \in [\tau',\tau)$, let us show how to construct $F_n$.
  First pick a finite subset $F'_n$ of $\{\sigma > \tau_n\}$ such that
  $V(f;F'_n)\ge T_f(\tau)+2\e$.  Pick any end $\bar{\tau}$ in $\cT$ with
  $\bar{\tau}>\tau$.  As $f$ vanishes on $\cT\setminus\cT^o$, we have
  $f(\bar{\tau})=0$, so that we may assume that $\bar{\tau}\in F'_n$.
  By Lemma~\ref{L436} we may then also assume that $\tau\in F'_n$.  In
  addition, we may assume that $F'_n$ contains a point $\tau'_n$ such
  that $\tau_n<\tau'_n<\tau$ and
  $F'_n\cap\{\sigma>\tau'_n\}=F'_n\cap\{\sigma\ge\tau\}$.  Let
  $F_n=F'_n\setminus\{\sigma\ge\tau\}$ and
  $F''_n=F'_n\cap\{\sigma\ge\tau\}$.  Then
  $V(f;F'_n)=V(f;F_n)+|f(\tau'_n)-f(\tau)|+V(f;F''_n)$.  Since
  $V(f;F''_n)\le T_f(\tau)$ and $|f(\tau'_n)-f(\tau)|<\e$, this gives
  $V(f;F_n)\ge\e$. 
  
  Finally we set $\tau_{n+1} \= \tau'_n (= \max F_n \cap [\tau',\tau[)$.
  It is clear that $F_n \subset \{ \sigma\ ;\ \sigma>\tau_n,
  \sigma\not>\tau_{n+1}\}$, and that $\tau_{n+1}$ belongs to $F_n$.
  This completes the inductive construction of $\tau_n$ and $F_n$ and
  ends the proof.
\end{proof}
Finally we address the formula for the total variation on a finite tree.
\begin{proof}[Proof of Proposition~\ref{L435}]
  Fix $\e>0$.  For each open segment $I\in\cI$ pick a finite set
  $F_I\subset I$ such that $V(f;F_I)\ge TV(f|_I)-\e$. 
  
  For a point $\tau \in B$, denote by $F_\tau$ the set of points
  $\s\ge\tau$ in $\bigcup_I F_I$ such that $[\tau,\s[\cap F_I$ is empty.
  This set consists of $\min F_I$ for the segments $I$ for which $\min I
  = \tau_0$. As any segment $I$ is totally ordered, one can add to $F_I$
  points arbitrarily closed to $\min I$. One can therefore assume
  $|d(\tau ; f) - V(f ; F_\tau\cup\{\tau\}) | \le \e$.  Declare $F$ to
  be the union of $B$ and all the $F_I$'s, and define $\underline{\tau}
  \= \max \{ \s \in F_\tau , \, \s < \tau \}$ for any $\tau\in B$. We
  then have
  \begin{multline*}
    TV(f)
    \ge V(f;F)
    =\sum_{\tau\in B}V(f;F_\tau)
    +\sum_{I\in\cI}V(f;F_I)
    + \sum_{\tau\in B} |f(\underline{\tau}) - f(\tau)| 
    \ge\\
    \ge\sum_{\tau\in B}|d(\tau;f)|
    +\sum_{I\in\cI}TV(f|_I)
    -(2|B|+|\cI|)\e.
  \end{multline*}
  Letting $\e\to0$, we get $TV(f) \ge \sum_{\tau\in B}|d(\tau;f)|
  +\sum_{I\in\cI}TV(f|_I)$.
  
  \medskip

  Conversely, pick a finite set $F$ such that $V( f ; F) \ge TV (f)
  -\e$. As $f$ vanishes at the ends of $\cT$, we can suppose $F$
  contains all ends of $\cT$.  For each point $\tau \in B$, and each
  segment $I$ containing $\tau$ in its boundary, pick a point $\tau_I
  \in I$ closed enough to $\tau$ such that $|f(\tau_I) - f(\vv_I)|\le
  \e$. Here $\vv_I$ denotes the tangent vector represented by $I$ at
  $\tau$.  When $\vv_I$ is represented by the root, we let as above
  $\underline{\tau} \= \tau_I$. As $f$ is left continuous, one can
  assume $|f(\underline{\tau}) - f(\tau)|\le \e$, for all $\tau\in B$.
  Otherwise, we denote by $F_\tau$ the union of all $\tau_I\ge\tau$.
  
  By Lemma~\ref{L436}, we can add all points $\tau\in B$ and $\tau_I$ to
  $F$, this may only increase $V( f; F)$.  Define $F_I \= F\cap I$. Then
  \begin{multline*}
    TV(f)
    \le \e + V(f;F) 
    =\e + \sum_{\tau\in B}V(f;F_\tau)
    +\sum_{I\in\cI}V(f;F_I)
    + \sum_{\tau\in B} |f(\underline{\tau}) - f(\tau)| 
    \le\\
    \le\sum_{\tau\in B}|d(\tau;f)|
    +\sum_{I\in\cI}TV(f|_I)
    +(1+M)\e,
  \end{multline*}
  where $M$ is the sum of the  number of branches at all points in $B$.
  We conclude the proof by letting $\e\to 0$.
\end{proof}                               
                                
%
%
%
%
\section{Representation Theorem I}\label{S411}
We are now ready to relate complex Borel measures and
normalized functions of bounded variation.
\begin{theorem}\label{T403}
  Let $(\cT,\le)$ be a weakly compact, complete, rooted nonmetric 
  tree.\footnote{For instance a complete, parameterizable, rooted nonmetric
    tree.}
  Then for any complex measure, the function 
  $I\rho=f_\rho:\cT^o\to\C$ defined by
  \index{$f_\rho$ (function of bounded variation associated to $\rho$)}
  \begin{equation}\label{e307}
    f_\rho(\tau)\=\rho\{\sigma\ge\tau\}
  \end{equation}
  belongs to $\cN$.
  Moreover, the map
  \begin{equation*}
    I:(\cM,\Vert\cdot\Vert)\to(\cN, TV), 
  \end{equation*}
  is an isometry and 
  restricts to a bijection between the set of positive measures
  $\cM^+$, and $\cN^+$.
  
  When $f\in\cN$, we shall denote by $df$
  \index{$df$ (measure associated to $f$)} 
  the unique complex Borel measure such that $I(df)=f$.
\end{theorem}
\begin{remark}
  If $\cT=[1,\infty]$, with the natural parameterization $\a(x)=x$,
  then $df=-df/dx$.
\end{remark}
If $\vv$ is a tangent vector not represented by $\tau_0$, then
$U(\vv)$ is a countable decreasing intersection of sets of the form
$\{\sigma\ge\tau\}$; and if $\vv$ is represented by $\tau_0$, then
$U(\vv)$ is the complement of a set $\{\sigma\ge\tau\}$. By regularity
of Borel measures, we immediately obtain
\begin{proposition}\label{P472}
  Pick $f\in \cN$, and $\rho\in\cM$ such that $\rho=df$. 
  Then 
  \begin{itemize}
  \item[(i)]
    $\rho\{\sigma\ge\tau\}=f(\tau)$ for every $\tau\in\cT^o$; 
  \item[(ii)]
    $\rho\,U(\vv)=f(\vv)$ for every tangent vector $\vv$
    not represented by $\tau_0$;
  \item[(iii)]
    $\rho\,U(\vv)=f(\tau_0)-f(\vv)$ for every tangent vector $\vv$
    represented by $\tau_0$;
  \item[(iv)]
    $\rho\{\tau\}=d(\tau;f)$ for every $\tau\in\cT^o$.
  \end{itemize}
\end{proposition}
Since $\cM$ is complete in the norm $\Vert\cdot\Vert$ we 
also infer
\begin{corollary}\label{C409}
  $(\cN,TV)$ is a Banach space.
\end{corollary}
We split the proof of the theorem into three parts.
%
%
\subsection{First step}
We show that $\rho \in \cM$ (resp.\ in $\cM^+$) implies $I\rho\in \cN$
(resp.\ in $\cN^+$) for any complex Borel measure; and $I: \cM \to \cN$
is injective.

\par
Write $f=f_\rho$, \ie $f(\tau)=\rho\,\{\sigma\ge\tau\}$ 
for any $\tau\in\cT^o$.  
Let us show that $f$ is left continuous and of bounded variation.  After
decomposing $\rho$ into real and imaginary parts, and further into
positive and negative parts, we may assume that $\rho$ is a positive
measure.  We need to show that $f$ is left continuous and strongly
decreasing.  That $f$ is left continuous is easy to prove: if $\tau_k$
increases to $\tau$ then $\rho\{\sigma\ge\tau_k\}$ decreases to
$\rho\{\sigma\ge\tau\}$.  To see that $f$ is strongly decreasing,
consider a finite subset $E$ of $\cT$ without relations and
$\tau\in\cT^o$ with $\tau<\sigma$ for all $\sigma\in E$.  Then the
subsets $\{\sigma'\ge\sigma\}_{\sigma\in E}$ are mutually disjoint and
contained in $\{\sigma\ge\tau'\}$, implying that $f(\tau)\ge\sum_Ef$,
so that $f$ is strongly decreasing.  This shows that $I$ sends $\cM$,
$\cM^+$ to $\cN$, $\cN^+$ respectively.

To prove the injectivity of $I$, suppose $I\rho = I\rho'$ for some
Borel measures $\rho,\rho'$.  Then $\rho \{ \s \ge \tau\} =\rho' \{ \s
\ge \tau\}$ for any $\tau \in \cT$, thus $\rho (E) = \rho'(E)$ for any
element of the elementary family $\cE$ defined in Section~\ref{S407}.
And $\rho = \rho'$ by Lemma~\ref{L442}. Thus $I$ is injective.
%
%
\subsection{Second step: from functions to measures}
Given $f\in\cN^+$ we shall find a \emph{positive} Borel measure
$\rho\in\cM$ such that $f(\tau)=\rho\,\{\sigma\ge\tau\}$ for any $\tau
\in \cT^o$.  This shows that $I : \cM^+ \to \cN^+$ is surjective,
hence bijective. Note that this implies $I : \cM \to \cN$ to be
surjective too.  Indeed, if $f$ belongs to $\cN$, we may write 
$f=\mathrm{Re}(f)+i\,\mathrm{Im}(f)$.
In view of Proposition~\ref{P441},
$\mathrm{Re}(f)$ and $\mathrm{Im}(f)$ are both differences 
of elements in $\cN^+$, thus lie in the image of $I$. Hence so does $f$.

So pick $f\in\cN^+$. Recall that this means that $f$ 
is nonnegative, left continuous, strongly decreasing, 
and vanishing on $\cT\setminus\cT^o$. 
To construct the positive measure $\rho$, we proceed as on
the real line (see~\cite[Section 1.5]{folland}), with suitable
adaptations to our setting, using
the elementary family $\cE$ and
the algebra $\cA$ from Section~\ref{S407}.
First define $\rho$ as a function on $\cE$ by
\begin{equation*}
  \rho\{\sigma\ge\tau, \sigma\not\ge\tau_k\}
  =f(\tau)-\sum_kf(\tau_k).
\end{equation*}
As $f$ is strongly decreasing, $\rho$ is nonnegative.
Recall that an element $E$ of the algebra $\cA$ is 
a finite disjoint union $\bigcup E_i$ of elements in $\cE$.
We can therefore try to extend $\rho$ to $\cA$ by declaring
$\rho(E)=\sum\rho(E_i)$. A priori, this is not well-defined,
as the decomposition of $E$ into elements of $\cE$ is
not unique. 
However, as in the proof of Proposition~1.15 in~\cite{folland},
if $E=\bigcup E_i=\bigcup F_j$ are two decompositions, then it
is easy to see that
$\sum_i\rho(E_i)=\sum_{i,j}\rho(E_i\cap F_j)=\sum_j\rho(F_j)$.

Hence $\rho$ is well defined on $\cA$. We claim that it defines
a premeasure on $\cA$. Clearly $\rho$ is finitely additive. 
It remains to show that if $(E_i)_1^\infty$ is a disjoint sequence
of elements in $\cA$ such that $E=\bigcup E_i\in\cA$, then
$\rho(E)=\sum_i\rho(E_i)$.
By finite additivity we may assume that $E\in\cE$, and we have
\begin{equation*}
  \rho(E)
  =\rho\left(\bigcup_1^nE_i\right)
  +\rho\left(E\setminus\bigcup_{n+1}^\infty E_i\right)
  \ge\rho\left(\bigcup_1^nE_i\right)
  =\sum_1^n\rho(E_i).
\end{equation*}
Letting $n\to\infty$ yields $\rho(E)\ge\sum_1^\infty\rho(E_i)$.  
For the converse inequality, fix $\e>0$ and write 
$E=\{\sigma\ge\tau,\sigma\not\ge\tau_k\}$ and 
$E_i=\{\sigma\ge\tau_i,\sigma\not\ge\tau_{ij}\}$.  
For each $k$, one can pick
$\tilde\tau_k\in\,]\tau,\tau_k[$ such that 
$0\le f(\tilde\tau_k)-f(\tau_k)<\e2^{-k}$ as $f$ is left continuous.
The set $K\=\{ \s\ge\tau,\, \s\not>\tilde{\tau}_k\}$
is then a compact set included in $E$, and as $f$ is decreasing
\begin{equation*}
  \rho(K)\ge\rho(E)-\sum (f(\tilde{\tau}_k)-f(\tau_k))\ge\rho(E)-\e. 
\end{equation*}
Similarly, for each $i$ pick $\tilde\tau_i<\tau_i$ such that 
$0\le f(\tilde\tau_i)-f(\tau_i)<\e2^{-i}$.  
Again this is possible, at least
as long as $\tau_i$ is not the root $\tau_0$. If $\tau_i=\tau_0$ for
some (unique) $i$, then we set $\tilde\tau_i=\tau_i=\tau_0$.  
The set $V_i\=\{\s>\tilde{\tau}_i,\,\s\not\ge\tau_{ij}\}$ is then
an open set containing $E_i$, and 
\begin{equation*}
  \rho(V_i)\le\rho(E_i)+(f(\tilde{\tau}_i)-f(\tau_i))\le\rho(E_i)+\e2^{-i}. 
\end{equation*}
The compact set $K$ can now be covered by finitely many of the open
sets $V_i$, say by $V_1,\dots,V_N$. This gives
\begin{equation*}
  \rho(E) \le \rho(K) + \e \le \sum_1^{N} \rho(V_i)+\e 
  \le \sum_1^\infty \rho(E_i ) + 2\e.
\end{equation*}
Letting $\e\to0$ we conclude that $\rho(E)\le\sum\rho(E_i)$,
implying that equality holds.

Thus $\rho$ is a premeasure on $\cA$. 
By Lemma~\ref{L442} it extends uniquely to a positive 
Borel measure $\rho$ which obviously satisfies $I\rho = f$.

%
%
\subsection{Total variation}
To complete the proof of Theorem~\ref{T403} we now
show that $d:\cN\to\cM$ preserves the norm.
Thus pick any $f\in\cN$ and write $\rho=df\in\cM$.
We must show that the total variation of $f$ equals
the total mass of $\rho$.

Let $T_f:\cT^o\to[0,\infty)$ be the function defined in~\eqref{e454}, 
\ie $T_f(\tau)=TV(f|_{\{\sigma\ge\tau\}})$.
By Proposition~\ref{P441}, $T_f$ belongs to $\cN^+$.
Therefore there exists a positive measure $\rho':=d(T_f)$ such that 
$T_f\{\tau\}=\rho'\{\sigma\ge\tau\}$ for all $\tau\in\cT^o$.
We shall prove that 
$\rho'=\vert\rho\vert$, which implies 
$TV(f)=\rho'(\cT)=\rho(\cT)=\Vert\rho\Vert$.

Consider $\tau\in\cT^o$, fix $\e>0$, and pick a finite set 
$F\subset\{\sigma\ge\tau\}$ such that
$T_f(\tau)\le V(f;F)+\e$. Following the notation in~\eqref{e443}
we have
\begin{equation*}
  V(f;F)
  =\sum_{\tau\in F\setminus F_{\max}}
  \left|f(\tau)-\sum_{\sigma\succ\tau}f(\sigma)\right|
  =\sum_{\tau\in F\setminus F_{\max}}
  \left|\rho\{\sigma'\ge\tau, \sigma'\not\ge\sigma\}\right|
  \le|\rho|\{\sigma'\ge\tau\}.
\end{equation*}
After letting $\e\to0$ we conclude that $\rho' \{ \s'\ge \tau \} = T_f
(\tau ) \le |\rho| \{ \s' \ge \tau\}$. This implies that $\rho'(E)\le
|\rho(E)|$ for all $E\in\cE$, so that 
$\rho'\le |\rho|$ by Lemma~\ref{L442}.

Let us now show that $|\rho|\le\rho'$. 
First consider a 
set $E$ in the elementary family $\cE$, \ie
$E=\{\sigma\ge\tau, \sigma\not\ge\tau_i\}$.
It is then clear that
\begin{equation*}
  |\rho (E)|
  =|f(\tau)-\sum_if(\tau_i)|
  \le|T_f(\tau)-\sum_iT_f(\tau_i)|
  =\rho'(E).
\end{equation*}
Next consider a connected open set $U$.  By Remark~\ref{R446}
there exists a countable increasing sequence $(E_i)_1^\infty$ such that
$E_i\in\cE$ and $\bigcup E_i=U$.  Thus
\begin{equation*}
  |\rho(U)|
  =\lim|\rho(E_i)|
  \le\limsup\rho'(E_i)
  =\rho'(U).
\end{equation*}
This easily implies that $|\rho(U)|\le\rho'(U)$ for all
(not necessarily connected) open sets $U$.
Finally consider an arbitrary Borel set $E$.
Since $\rho$ and $\rho'$ are Radon measures we may find
a decreasing sequence $(U_j)_1^\infty$ of open neighborhoods
of $E$ such that 
$\lim\rho(U_j)=\rho(E)$ and $\lim\rho'(U_j)=\rho'(E)$.
This gives $|\rho(E)|\le\rho'(E)$. 
By the definition of the total variation measure, we conclude that
$|\rho|\le\rho'$.

This ends the proof of Theorem~\ref{T403}.

For further reference we note that we proved
\begin{proposition}\label{P418}
For any function $f\in\cN$, the total variation measure of $df$ can be
computed as follows:
\begin{equation}\label{e927}
|df| = d \left( T_f \right).
\end{equation}
\end{proposition}
%
%
%
%
\section{Complex tree potentials}\label{S412}
We now turn to the second type of functions that will be
identified with complex Borel measures. They are slightly more
complicated to define, and their relationship to measures is 
less direct, but they are the functions that appear most naturally 
in applications.
Their analogues on the real line are (normalized) antiderivatives 
of functions of bounded variation, and we shall define them
accordingly also in the tree setting.

In this section we work with a complete, rooted, nonmetric tree
$(\cT,\le)$ and a fixed increasing parameterization
$\a:\cT\to[1,\infty]$ of $\cT$.  While the parameterization
was not important when studying functions of bounded variation, here
we must fix a choice. 
We require that $\a$ be increasing and $\a(\tau_0)=1$, 
where $\tau_0$ is the root of $\cT$.  
As before, denote by $\cT^o$ the ends of $\cT$,
\ie the set of nonmaximal points in $\cT$.  
%
%
\subsection{Definition}
Recall the Banach space $\cN$ defined in Section~\ref{S418}:
its elements are complex-valued, left continuous functions on
$\cT^o$ of bounded variation.
The norm $TV(f)$ of $f\in\cN$ is the total variation of
the extension of $f$ to $\cT$ obtained by setting $f=0$ on
$\cT\setminus\cT^o$. 

As in Section~\ref{S10} we need to integrate with respect to the 
parameterization $\a$.
Fix $\s_0<\s_1\in\cT$ and write $\a_i=\a(\s_i)$ for $i=0,1$. 
For $t\in[\a_0,\a_1]$, define $\s_t$ to be the unique element in 
$[\s_0,\s_1]$ such that $\a(\s_t)=t$. 
For any measurable function $f$ on $\cT$, set
$\int_{\s_0}^{\s_1}f(\sigma)\,d\a(\sigma)\=\int_{\a_0}^{\a_1}f(\s_t)\,dt$.
\begin{definition}\label{D414}
  A function $g:\cT^o\to\C$ is a 
  \emph{complex tree potential}\index{tree potential!complex} 
  if there exists a
  function $f:\cT^o\to\C$ such that $f\in\cN$ and
  \begin{equation}\label{e451}
    g(\tau)=f(\tau_0)+\int_{\tau_0}^\tau f(\sigma)\,d\a(\sigma)
  \end{equation}
  for any $\tau\in\cT^o$.  We shall write $f=\delta g$
  \index{$\delta g$ (left derivative of $g$)}, 
  and $g=If$, when $f,g$ are related by~\eqref{e451}.  
  We denote by $\cP$\index{$\cP$ (complex tree potentials)} 
  the set of all complex tree potentials.
\end{definition}
Notice that $f$ appears both in the constant term and in the
integrand. Also note that $f=\delta g$ is uniquely determined by $g$:
\begin{equation}\label{e452}
  f(\tau_0)=g(\tau_0)
  \qand
  f(\tau)
  =\frac{dg}{d\a}(\tau)
  :=\lim_{\sigma\to\tau-}
  \frac{g(\tau)-g(\sigma)}{\a(\tau)-\a(\sigma)}
  \quad\text{for $\tau\ne\tau_0$}.
\end{equation}
\begin{definition}
  We refer to $\frac{dg}{d\a}$ as the 
  \emph{left derivative}\index{left derivative (in a tree)} 
  of $g$.
\end{definition}
It is clear that $\cP$ is a complex vector space as $\cN$ is.
Moreover, $I:\cN\to\cP$ is a bijection. 
If $g\in\cP$, we may thus define a norm on $\cP$
by $\Vert g\Vert:=TV(\delta g)$.
\index{topology!strong (on $\cP$)}
It is clear from the definitions and Corollary~\ref{C409} that
\begin{proposition}
  The set $(\cP,\Vert\cdot\Vert)$ is a Banach space and the mappings
  $I:\cN\to\cP$ and $\delta:\cP\to\cN$ 
  defined by~\eqref{e451} and~\eqref{e452} 
  are isometric isomorphisms.
\end{proposition}
Define the \emph{support}
\index{support!of a complex tree potential} 
of a complex tree potential $g$
to be the smallest subtree $\cS$ of $\cT$ such that
$g$ is constant on any segment in $\cT$ disjoint from $\cS$.
We obtain from Corollary~\ref{C408}:
\begin{corollary}\label{C410}
  The support of any complex tree potential is contained
  in the completion of a countable union of finite subtrees.
\end{corollary}
\begin{remark}
  The definition of a complex tree potential is a bit 
  indirect. It is possible to give an
  equivalent, more direct definition by saying that $g\in\cP$ iff:
  \begin{itemize}
 \item[(i)]
    the restriction of $g$ to any segment $[\tau_0,\tau]$ 
    admits a left-derivative at all points;
  \item[(ii)]
    the function $f$ given by~\eqref{e452} is left continuous
    and has bounded variation.
  \end{itemize}
\end{remark}
%
%
\subsection{Directional derivatives}
Consider a function $g:\cT^o\to\C$. 
and pick a tangent vector $\vv$ at a point $\tau\in\cT^o$.
We define the \emph{derivative}\index{directional derivative} 
of $g$ along $\vv$ (when it exists) by
\index{$D_\vv g$ (derivative along $\vv$)}
\begin{equation*}
  D_\vv g
  =\lim_{k\to\infty}\frac{g(\tau_k)-g(\tau)}{|\a(\tau_k)-\a(\tau)|}
\end{equation*}
for any sequence $(\tau_k)_1^\infty$ converging
to $\tau$ along $\vv$.
Note that if $\vv$ is represented by $\tau_0$, then 
$D_\vv g=-\frac{dg}{d\a}$, where $\frac{dg}{d\a}$
is the left derivative as in~\eqref{e452}.
It is then clear from Lemma~\ref{L433} that
\begin{lemma}\label{L443}
  If $g\in\cP$ is a complex tree potential, 
  then $D_\vv g$ exists for every tangent vector $\vv$;
  in fact $D_\vv g=(\delta g)(\vv)$.
\end{lemma}
Next we define the analogue of the quantity $d(\tau;f)$ for $f\in\cN$
defined in~\eqref{e444}. 
Namely, if $g$ is a complex tree potential, 
then we set
\index{$\Delta(\tau;f)$ (atom of $\Delta f$ on $\tau$)}
\begin{equation}
  \Delta(\tau_0;g)=g(\tau_0)-\sum_{\vv\in T\tau_0}D_\vv g
  \qand
  \Delta(\tau;g)=-\sum_{\vv\in T\tau}D_\vv g\ \text{for $\tau\ne\tau_0$}.
\end{equation}
Note that this makes sense as the series $\sum_{\vv\in T\tau}D_\vv g$
for any $\tau\in\cT^o$ is absolutely convergent thanks
to Proposition~\ref{P444}.  
Lemma~\ref{L432} immediately implies 
\begin{lemma}
  Let $g:\cT^o\to\C$ be a complex tree potential.  Then
  $\Delta(\tau;f)=0$ for all but countably many points $\tau$, and the
  series $\sum_{\tau\in\cT}\Delta(\tau;g)$ is absolutely convergent;
  in fact $\sum_{\tau\in\cT}|\Delta(\tau;g)|\le\Vert g\Vert$.
\end{lemma}
%
%
\section{Representation Theorem II}\label{S413}
We can now state the relation between complex tree potentials and
complex Borel measures.
\begin{theorem}\label{T404}
  Let $(\cT,\le)$ be a complete, rooted nonmetric tree equipped with
  an increasing parameterization $\a:\cT\to[1,\infty]$ such that
  $\a(\tau_0)=1$.  Then for any complex Borel measure, 
  the function
  \index{$g_\rho$ (tree potential of a measure)}
  $g_\rho:\cT^o\to\C$ defined by
  \begin{equation}\label{e457}
    g_\rho(\tau)=\int_\cT\a(\sigma\wedge\tau)\,d\rho(\sigma)
  \end{equation}
  is a complex tree potential, and
  \begin{equation*}
    (\cM,\Vert\cdot\Vert)\ni\rho\mapsto g_\rho\in(\cP,\Vert\cdot\Vert)
  \end{equation*}
  is an isometry. When $g\in\cP$, we shall denote by $\Delta g$ the
  unique complex Borel measure such that $g_{\Delta g} = g$.
\end{theorem}
\begin{definition}
  We call $\Delta$\index{$\Delta$ (Laplace operator)} the 
  \emph{Laplace operator}\index{Laplace operator};
  if $g$ is a complex tree potential, then
  $\Delta g$ is the
  \emph{Laplacian}\index{Laplacian} of $g$.
  Given a complex Borel measure $\rho\in\cM$, the
  function $g_\rho$ in~\eqref{e457} is called
  \emph{the potential of $\rho$}
  \index{tree potential!of a measure}.
\end{definition}
The following formulas relate a complex tree potential and its Laplacian.
\begin{proposition}
Let $g\in \cP$, and $\rho\in\cM$ such that $\rho = \Delta g$. 
Then 
\begin{itemize}
  \item[(i)]
    $\rho\{\sigma\ge\tau\}=\frac{dg}{d\a}(\tau)$ for every $\tau\in\cT^o$; 
  \item[(ii)]
    $\rho(U(\vv))=D_\vv g$ for every tangent vector $\vv$
    not represented by $\tau_0$;
  \item[(iii)]
    $\rho(U(\vv))=D_\vv g+g(\tau_0)$
    for every tangent vector $\vv$ represented by $\tau_0$;
  \item[(iv)]
    $\rho\{\tau\}=\Delta(\tau;g)$ for every $\tau\in\cT^o$.
  \end{itemize}
\end{proposition}
\begin{proof}[Proof of Theorem~\ref{T404}]
  For any complex measure $\rho$ it is easy to check 
  from~\eqref{e307} and~\eqref{e457} that 
  $g_\rho=If_\rho$ and $f_\rho=\delta g_\rho$.
  As $\delta:(\cP,\Vert\cdot\Vert)\to(\cN,TV)$
  is an isometry, Theorem~\ref{T404} follows from Theorem~\ref{T403}.
\end{proof}
\begin{remark}\label{R303}
  In the absence of branching, $\Delta$ reduces to 
  the usual Laplacian on the real line:
  if $\cT=[1,\infty]$ is parameterized by $\a(x)=x$,
  then $\Delta=-\frac{d^2}{dx^2}$.
\end{remark}
\begin{remark}
  Our definition of $\Delta$ also generalizes 
  the Laplace operator on (rooted) simplicial trees 
  as presented in \eg~\cite{cartier}.

  To see this, we make use of the discussion in Section~\ref{S303}.
  We hence view a rooted simplicial tree
  as a rooted (nonmetric) $\N$-tree $\cT_\N$ and equip it with
  its canonical parameterization $\a:\cT_\N\to\N$ with 
  $\a(\tau_0)=1$. Let $\cT_\R$ be the associated
  rooted nonmetric $\R$-tree obtained by ``adding edges''.
  Finally let $\cT$ be the completion of $\cT_\R$
  (if $\cT_\N$ is finite, then $\cT=\cT_\R$).
  We view $\cT_\N$ as a subset of $\cT$ and equip $\cT$
  with a parameterization $\a:\cT\to[1,\infty]$ extending the
  one on $\cT_\N$.

  Any complex valued function $g:\cT_\N\to\C$ extends 
  uniquely as a function $g:\cT^o\to\C$ which is
  affine (in the parameterization $\a$) 
  on the edges of $\cT^o$. 

  If $\cT_\N$ is finite, then 
  $g$ is always a complex tree potential whose Laplacian
  $\Delta g$ is atomic, supported on $\cT_\N$, and
  \begin{align}
    \Delta g\{\tau\}
    &=-\sum_{\sigma\sim\tau}(g(\sigma)-g(\tau))\
    \text{if $\tau\ne\tau_0$}\label{e305}\\
    \Delta g\{\tau_0\}
    &=g(\tau_0)-\sum_{\sigma\sim\tau_0}(g(\sigma)-g(\tau_0))\label{e306}
  \end{align}
  Here the sums are over all $\sigma\in\cT_\N$ adjacent to
  $\tau$ and $\tau_0$, respectively.

  If $\cT_\N$ is infinite, the situation is more complicated.
  Suffice it to say that $g$ is a complex tree potential
  on $\cT$ if the quantities in the right hand sides 
  of~\eqref{e305} and~\eqref{e306} are all nonnegative.
  In that case the Laplacian $\Delta g$ is a positive measure
  on $\cT$ whose restriction to $\cT^o$ is a sum of (at most)
  countably many atoms at points in $\cT_\N$. The masses
  at these points are still given by~\eqref{e305}-\eqref{e306}.
  It is quite possible, however, for $\Delta g$ to be nonzero
  and even nonatomic on $\cT\setminus\cT^o$.
\end{remark}
%
%
%
%
\section{Atomic measures}\label{S414}
A measure $\rho\in\cM$ is \emph{atomic}\index{measure!atomic} 
if $\rho=\sum_{j=1}^nc_j\tau_j$, where $n<\infty$, 
$c_j\in\C$ and $\tau_j\in\cT$ (recall that we identify
a point in $\cT$ with the corresponding point mass in $\cM^+$).

The representation Theorems~\ref{T403} and~\ref{T404} 
immediately yield
\begin{proposition}\label{Patom1}
Pick $f \in \cN$. Then $df \in \cM$ is atomic iff there exists a
  finite subset $F\subset\cT$ containing the root $\tau_0$, such that:
\begin{itemize}
  \item[(i)]
    $f$ is zero outside the finite tree 
    $\bigcup_{\tau\in F}[\tau_0,\tau]$;
  \item[(ii)]
    $f$ is constant on any segment in $\cT$ not containing any point 
    in $F$.
  \end{itemize}
\end{proposition}
\begin{proposition}\label{Patom2}
Pick $g\in\cP$. Then  $\Delta g\in\cM$ is atomic iff
  there exists a finite subset $F\subset\cT$ containing the root
  $\tau_0$, such that:
  \begin{itemize}
  \item[(i)]
    $g$ is constant on any segment that intersects the finite subtree
    $\bigcup_{\tau\in F}[\tau_0,\tau]$ in at most one point;
  \item[(ii)]
    $g$ is an affine function of the parameterization $\a$ 
    on any segment in $\cT$ not containing any point in $F$
    in its interior.
  \end{itemize}
\end{proposition}

%
%
%
%
\section{Positive tree potentials}\label{S415}
Our next goal is to describe the image of the set of positive measures
respectively in the set of functions of bounded variation, 
and in the set of complex tree potentials (Theorem~\ref{T406} below).

In fact, we have already identified functions $f\in \cN$ giving 
rise to positive measures $df \in \cM^+$ as belonging to the cone $\cN^+$ of
nonnegative, left continuous, strongly decreasing functions on
$\cT^o$ (Theorem~\ref{T403}).
%
%
\subsection{Definition}
We now wish to describe the preimage of $\cM^+$ in $\cP$
under the Laplace operator.
\begin{definition}\label{treepot}
  A function $g:\cT^o\to\R$ is called a 
  \emph{positive tree potential}\index{tree potential!positive},
  or simply \emph{tree potential} (on $\cT$),
  if the following conditions are satisfied:
  \begin{itemize}
  \item[(P1)] 
    $g$ is nonnegative, increasing, and concave along totally ordered 
    segments;
   \item[(P2)]
    if $\tau\ne\tau_0$, then
    $\sum_{\vv\in T\tau}D_\vv g\le0$;
  \item[(P3)]
    $\sum_{\vv\in T\tau_0}D_\vv g\le g(\tau_0)$.
  \end{itemize}
  We denote by $\cP^+$\index{$\cP^+$ (positive tree potentials)}
  the set of all positive tree potentials on $\cT$.
\end{definition}

Note that~(P1) implies that
the directional derivative $D_\vv g$ is well defined for every
tangent vector $\vv$ at any point $\tau\in\cT^o$.
Moreover, $D_\vv g\ge0$ except if 
$\vv$ is represented by $\tau_0$, in which case
$D_\vv=-\frac{dg}{d\a}\le0$.
Hence the series in~(P2) and~(P3) are well-defined.

If $\tau\in\cT^o$, then~(P2) and~(P3) imply that
$D_\vv g=0$ for all but countably many tangent 
vectors $\vv\in T\tau$. 
Moreover, if $\vv$ is not represented
by $\tau_0$ and $D_\vv g=0$, then~(P1) implies that $g$
is constant in the open set $U(\vv)$.
Hence conditions~(P1)-(P3) are quite strong.
\begin{proposition}\label{P435}
  Every positive tree potential $g\in \cP^+$ can be written
  $g=If$ for some $f\in \cN^+$ (see~\eqref{e451}).
  As a consequence, every positive tree potential 
  is a complex tree potential, \ie $\cP^+\subset\cP$.
\end{proposition}
\begin{proof}
  Consider a positive tree potential $g\in\cP^+$.
  It follows from~(P1) that the left derivative
  of $g$ is defined at any point. Thus we may define
  $f:\cT^o\to\R$ by~\eqref{e442}, \ie
  $f(\tau_0)=g(\tau_0)$ and
  $f(\tau)=\frac{dg}{d\a}(\tau)$ for $\tau\ne\tau_0$.
  It then follows from~(P1) that $f$ is nonnegative
  and left continuous.
  We will show that $f$ is also strongly decreasing so
  that $f\in\cN^+\subset\cN$. By definition of $\cP$,
  this will show that $g$ is a complex tree potential.
  
  Hence consider a finite, nonempty set $E\subset\cT^o$ 
  without relations
  and $\tau\in\cT^o$ such that $\tau<\sigma$ for all $\sigma\in E$.
  We have to show that $f(\tau)\ge\sum_Ef$.

  First assume that $\sigma_1\wedge\sigma_2=\tau$ for
  any two distinct elements $\sigma_1,\sigma_2\in E$.
  (This is true if $E$ has only one element!)
  Let $V$ be the set of tangent vectors at $\tau$
  represented by the elements of $E$.
  Then~(P1) implies that $\sum_Ef\le\sum_VD_\vv g$.
  If $\tau\ne\tau_0$, then we conclude from~(P2) that
  \begin{equation*}
    \sum_Ef-f(\tau)
    \le\sum_VD_\vv g-\frac{dg}{d\a}(\tau)
    \le\sum_{T\tau}D_\vv g
    \le0.
  \end{equation*}
  If instead $\tau=\tau_0$, then~(P3) gives
  \begin{equation*}
    \sum_Ef-f(\tau)
    \le\sum_VD_\vv g-g(\tau_0)
    \le\sum_{T\tau}D_\vv g-g(\tau_0)
    \le0.
  \end{equation*}
  
  In the general case we proceed by induction on the
  number of elements in $E$. 
  By the previous step we may assume that there exists
  $\tau'>\tau$ and a decomposition 
  $E=E'\cup E''$, where $E'$ and $E''$ are disjoint, 
  $E'$ has at least two elements ($E''$ could be empty), 
  $\sigma'_1\wedge\sigma'_2=\tau'$ for all distinct
  elements $\sigma'_1,\sigma'_2\in E'$, and
  $\sigma'\wedge\sigma''<\tau'$ whenever
  $\sigma'\in E'$ and $\sigma''\in E''$.
  Then $E''\cup\{\tau'\}$ has no relations and
  by the inductive hypothesis and the first step 
  we obtain
  \begin{equation*}
    f(\tau)
    \ge f(\tau')+\sum_{E''}f
    \ge\sum_{E'}f+\sum_{E''}f
    =\sum_Ef.
  \end{equation*}
This concludes the proof.
\end{proof}
Let us mention here some continuity properties of positive tree
potentials with respect to the weak topology.
\begin{lemma}\label{L702}
  Any positive tree potential on $\cT$ is (weakly)
  lower semicontinuous and restricts
  to a continuous tree potential on any finite subtree.
\end{lemma}
\begin{proof}
Let $f\= \delta g\in\cN^+$. For each $n\in\N^*$, define $\cS_n \=\{f >
1/n\}$. This is a finite tree. Set $f_n \= \mathbf{1}_{\cS_n}\, f$, and
$g_n \= If_n$. Then~(P1) implies that the potential $g_n$ is
continuous on the finite tree $\cS_n$. As it is locally constant
outside $\cS_n$, it is weakly continuous on $\cT$.  The sequence $f_n$
increases pointwise towards $f$.  By integration $g_n$ increases
pointwise to $g$ as $n\to\infty$. Thus $g$ is lower semicontinuous.
\end{proof}                                
However, positive tree potentials are not necessarily 
continuous as the following example on the valuative tree shows.
\begin{example}\label{E701}
We work in the valuative tree $\cV$.
Fix local coordinates $(x,y)$,
set $\rho=\sum_{n\ge1}n^{-2}\nu_{y+nx}$ and let
  $g=g_\rho$ be the associated positive tree potential
  given by~\eqref{e451}.
  If $\nu_n=\nu_{y+nx,n^3}$, then $\nu_n\to\nu_\fm$ weakly but
  $g(\nu_n)>n\to\infty>g(\nu_\fm)$. 
\end{example}
%
%
\subsection{Jordan decompositions}
We now show that the three positive cones 
$\cM^+$, $\cN^+$ and $\cP^+$ are isomorphic and deduce the
existence of Jordan decompositions of real-valued elements
in $\cM$, $\cN$ and $\cP$. 
\begin{theorem}\label{T406}
  The isomorphisms
  \begin{equation*}
    d:\cN\to\cM
    \qand
    \Delta:\cP\to\cM
  \end{equation*}
  restrict to bijections
  \begin{equation*}
    d:\cN^+\to\cM^+
    \qand
    \Delta:\cP^+\to\cM^+.
  \end{equation*}
\end{theorem}
\begin{proof}
  The fact that $d:\cN^+\to\cM^+$ is a bijection follows from
  Theorem~\ref{T403}.  To prove that $\Delta:\cP^+\to\cM^+$ is also a
  bijection, we prove that $g \in \cP^+$ iff $f$ belongs to $\cN^+$
  where $g$ and $f$ are related by~\eqref{e451}.  That $f$ belongs to
  $\cN^+$ if $g \in \cP^+$ follows from Proposition~\ref{P435}.
  Conversely, pick $f \in \cN^+$, and define $g (\tau ) \= f(\tau_0) +
  \int_{\tau_0}^\tau f(\s) d\a(\s)$. As $f$ is non-decreasing and
  non-negative it is clear that $g$ satisfies~(P1). By
  Proposition~\ref{P472} (iv), we have $ d( \tau ; f ) = df \{\tau \}$
  which is non-negative as $df$ is a positive measure. 
  Recall from~\eqref{e444} that 
  $d(\tau ; f) = -\sum_{\vv\in T\tau}f(\vv)$ when
  $\tau \ne \tau_0$, and $= f(\tau_0)- \sum_{\vv\in T\tau}f(\vv)$ when
  $\tau = \tau_0$, and that $D_\vv g = f(\vv)$. 
  Thus~(P2) and~(P3) are satisfied, so 
  $g\in\cP^+$, completing the proof.
\end{proof}
Using the Jordan decomposition of a real measure into positive 
measures and the isomorphisms in 
Theorems~\ref{T403},~\ref{T404},~\ref{T406},
we infer from the result above:
\begin{Prop-def}{~}\label{P419}
  \begin{itemize}
  \item
    Any real Borel measure $\rho$ is the difference 
    of two positive measures $\rho=\rho_+-\rho_-$.
    These measures are 
    uniquely determined by the condition 
    $\|\rho\|=\|\rho_+\|+\|\rho_-\|$.
  \item
    Any real-valued function $f\in\cN$ is the difference 
    of two functions $f=f_+-f_-$,
    where $f_\pm\in\cN^+$. These functions
    are uniquely determined by the condition $TV(f)=TV(f_+)+TV(f_-)$.
  \item
    Any real-valued complex tree potential $g \in \cP$ 
    is the difference of two positive tree potentials 
    $g=g_+-g_-$, with $g_\pm\in\cP^+$. 
    This decomposition is uniquely determined by the condition 
    $\|g\|=\|g_+\|+\|g_-\|$.
  \end{itemize}
  Any of these decompositions above is called the 
  \emph{Jordan decomposition}\index{Jordan decomposition}
  of $\rho$, $f$ or $g$ respectively.
\end{Prop-def}
Recall that in the case of measures, 
$\rho_+$ and $\rho_-$ are also characterized 
by the fact that their support are disjoint in 
the sense there are Borel
sets $E_\pm$ such that $\rho_+(E_+)=\rho_-(E_-)=1$, 
and $\rho_+(E_-)=\rho_-(E_+)=0$.
%
%
\section{Weak topologies and compactness}\label{S308}
We have so far considered $\cM$, $\cN$ and $\cP$ with 
topologies induced by natural norms on these three spaces. 
It is important in applications to consider weaker topologies, 
in which these spaces, or at least subspaces of them,
are compact.
Here we shall show how to accomplish this in the cones 
$\cM^+$, $\cN^+$ and $\cP^+$.

Recall that the weak topology on $\cM^+$
\index{topology!weak (on $\cM^+$)}
is defined in terms of convergence: $\rho_k\to
\rho$ iff $\int \varphi \ d\rho_k\to \int \varphi\ d\rho$ for any
(weakly) continuous $\varphi$ on $\cT$.

We define the weak topology on $\cN^+$
\index{topology!weak (on $\cN^+$)} 
in terms of convergent sequences as follows:
$f_k\to f$ iff $f_k(\tau_0) \to f(\tau_0)$, and $f_k\to f$ pointwise
on $\cT\setminus\{\tau_0\}$ except at (at most) countably many points.
Note that this is a well-defined Hausdorff topology as $f_1 = f_2$
outside countably many points implies $f_1 = f_2$ on
$\cT\setminus\{\tau_0\}$. A weak limit in $\cN^+$ is hence uniquely
determined.

Finally the weak topology on $\cP^+$
\index{topology!weak (on $\cP^+$)} 
is defined in terms of pointwise convergence:
$g_k\to g$ iff $g_k(\tau)\to g(\tau)$ for any $\tau\in\cT^o$.
\begin{theorem}\label{T407}
  The maps
  \begin{equation*}
    \cN^+ 
    \mathop{\longrightarrow}\limits^{d}
    \cM^+ 
    \qand
    \cP^+
    \mathop{\longrightarrow}\limits^{\Delta}
    \cM^+
  \end{equation*}
  are homeomorphisms in the weak topology.
\end{theorem}
As the set of positive measures of mass $1$ is compact, so are
its images in $\cN^+$ and $\cP^+$. This remark has ramifications
for the structure of the cones $\cN^+$ and $\cP^+$. 
Here we only mention an application to positive tree potentials.
\begin{corollary}\label{C209}
  The space of positive tree potentials normalized by 
  $g(\tau_0)=1$ is compact in the topology of pointwise convergence. 
  Moreover:
  \begin{itemize}
  \item[(i)]
    from any sequence $(g_n)_1^\infty$ of positive tree potentials
    such that $\sup g_n(\tau_0)<+\infty$, one can extract a 
    subsequence $g_{n_k}$ converging pointwise;
  \item[(ii)]
    if $(g_i)_{i\in I}$ is any family of positive tree potentials,
    then $g=\inf_ig_i$ is also a positive tree potential;
  \item[(iii)]
    if $(g_n)_1^\infty$ is an increasing sequence of 
    positive tree potentials such that 
    $\sup g_n(\tau_0)<+\infty$, then $g=\sup_ng_n$ is also
    a positive tree potential.
  \end{itemize} 
\end{corollary}
\begin{remark}
  As the properties above indicate, 
  positive tree potentials play a role 
  similar to that of concave functions on the real line, 
  or superharmonic functions on the unit disk in the complex plane.
\end{remark}
\begin{proof}
  The first assertion is a consequence of the compactness of the
  set of positive measures with bounded mass. The second statement
  is proved using~(P1)-(P3) in the same way as the fact that
  the family of nonnegative concave functions on $[0,\infty[$ 
  are closed under infima.
  The same is true of~(iii). 
  Notice that $\sup g_n(\tau_0)<\infty$ implies
  $\sup_ng_n(\tau)\le\a(\tau)\sup_ng(\tau_0)<\infty$ for all
  $\tau\in\cT^o$. The details are left to the reader.
\end{proof}
\begin{proof}[Proof of Theorem~\ref{T407}]
  Thanks to Theorems~\ref{T403},~\ref{T404} and~\ref{T406}, both maps
  $d,\Delta$ are bijective.  To complete the proof, we need to show
  that these maps and their inverses are weakly continuous.
  For sake of simplicity we shall only prove that they are
  sequentially continuous. We leave to the reader to check that they
  are indeed continuous, using the language of nets, or of filters.
 
  Let us first prove that $d$ is a homeomorphism.  By
  Theorem~\ref{T403}, $d$ is an isometry. It is hence sufficient to
  prove that its restriction to $\cM^+(1)$ the set of positive measures
  of mass $1$, induces a homeomorphism onto its image $\cN^+(1)$,
  which consists of functions in $\cN^+$ with $f(\tau_0)=1$.  The set
  $\cM^+(1)$ is weakly compact, and $\cN^+(1)$ is Hausdorff
  (see above). We thus
  only need to prove that $I=d^{-1}$ is weakly continuous.
  
  So consider a sequence of positive measures $\rho_n$ of mass $1$
  converging weakly to $\rho$. Write $f_n=I\rho_n$, $f=I\rho$.
  It is clear that $f_n(\tau_0)=\rho_n\{\cT\}=1\to f(\tau_0)$.
  The following result generalizes Proposition~7.19 in~\cite{folland}
  and exemplifies the idea that the quantity $d(\tau;f)$ 
  measures the discontinuity of $f$ at $\tau$:
  \begin{lemma}\label{L-elementary}
    We have $f_n(\tau)\to f(\tau)$ for all $\tau\in\cT^o$
    with $d(\tau;f)=0$.
  \end{lemma}
  In view of Lemma~\ref{L432}, this result implies that
  $f_n\to f$ on $\cT^o$ except on a countable subset.
  Thus $f_n\to f$ weakly, which 
  completes the proof that $d$ is a homeomorphism.
  
  \medskip
  To prove that $\Delta$ is a homeomorphism, we proceed in the same way.  
  Note that 
  $\Delta^{-1}\rho=g_\rho=\int\a(\sigma\wedge\cdot)\,d\rho(\sigma)$ 
  for any measure $\rho$.  
  It is sufficient to prove that $g_{\rho_n} \to g_\rho$ 
  pointwise when $\rho_n\to\rho$ weakly.  
  But for any $\tau\in\cT^o$, then
  \begin{equation*}
    g_{\rho_n}(\tau)
    =\int\a(\sigma\wedge\tau)\,d\rho_n(\sigma)
    \to\int\a(\sigma\wedge\tau)\,d\rho(\sigma)
    =g_\rho(\tau)
  \end{equation*}
  since the function 
  $\sigma\mapsto\a(\sigma\wedge\tau)$ on $\cT^o$ is 
  weakly continuous (Corollary~\ref{C398}).
  This completes the proof.
\end{proof}
\begin{proof}[Proof of Lemma~\ref{L-elementary}]
  Fix $\tau\in\cT^0$, $\tau\ne\tau_0$. 
  First pick an increasing sequence 
  $\tau_k<\tau$ converging towards $\tau$. 
  Let $\varphi_k$ be a continuous
  increasing function, with values in $[0,1]$, $1$ on 
  $\{\s\ge\tau\}$, and $0$ on $\{\s\not\ge\tau_k\}$. Then
  \begin{equation*}
    \limsup_n f_n (\tau)
    \le\limsup_n\int\varphi_k\,d\rho_n 
    =\int\varphi_k\,d\rho
    \le f(\tau_k).
  \end{equation*}
  As $f$ is left continuous, we infer $\limsup f_n(\tau)\le f(\tau)$.

  Now assume $d(\tau;f)=0$. Fix $\e>0$.
  Pick finitely many tangent vectors $\{\vv_j\}$ at $\tau$, 
  not representing $\tau_0$, such that $\sum f(\vv_j)\ge f(\tau)-\e$.
  Pick sequences 
  $(\tau_{jk})_{k=1}^\infty$ decreasing to $\tau$ such that
  $\tau_{jk}$ represents $\vv_j$ for all $j,k$.
  Also pick continuous increasing functions 
  $\varphi_{jk}$ with values in $[0,1]$, 
  $1$ on $\{\s\ge\tau_{jk}\}$ and 
  vanishing outside $U(\vv_j)$.
  Set $\varphi_k=\sum_j\varphi_{jk}$.
  Then 
  \begin{equation*}
    \liminf_n f_n(\tau) 
    \ge\liminf_n\int\varphi_k\,d\rho_n
    =\int\varphi_k\,d\rho
    \ge\sum_jf(\tau_{jk}).
  \end{equation*}
  Letting $k\to\infty$ yields
  $\liminf_n f_n(\tau)\ge\sum_jf(\vv_j)\ge f(\tau)-\e$, so
  as $\e\to0$ we obtain
  $\liminf_n f_n(\tau)\ge f(\tau)$.
  Hence $\lim f_n(\tau)=f(\tau)$ and we are done.
\end{proof}
%
%
%
%
\section{Restrictions to subtrees}\label{S409}
It is often important in applications to consider the 
restriction of functions and measures to subtrees, as well
as extensions from a subtree to the larger tree.

Let $\cS$ be a complete subtree of $\cT$ 
(as $\cT$ is a rooted tree, 
we assume that the root $\tau_0$ is contained in $\cS$). 
Denote by $p=p_\cS:\cT\to\cS$ the retraction map defined by
$p_\cS(\tau):=\max[\tau_0,\tau]\cap\cS$ and
by $\imath=\imath_\cS:\cS\to\cT$ the inclusion map.
By Lemma~\ref{L308} and Lemma~\ref{L307}, these are
both continuous.

When $f:\cT^o\to[0,\infty)$ belongs to $\cN^+$, 
its restriction to $\cS^o$ is an element of $\cN^+(\cS)$.
(Here it is important to regard the elements of 
$\cN^+(\cS)$ as functions on $\cS^o$ and extend them
by zero on $\cS\setminus\cS^o$.)
Conversely, if $f\in\cN^+(\cS)$, we may extend it to 
a function in $\cN^+$ by declaring it to be zero outside 
$\cS^o$.
Clearly the composition 
$\cN^+(\cS)\to\cN^+(\cT)\to\cN^+(\cS)$
is the identity.

In the case of tree potentials the situation is only 
slightly more complicated. 
Let $g:\cT^o\to[0,\infty)$ be a positive tree potential.
It is straightforward to verify 
from~(P1)-(P3) that the restriction $\imath^*g$
of $g$ to $\cS^o$ is a positive tree potential on $\cS$.
Conversely, if $g$ is a positive tree potential on $\cS$, 
then the 
\emph{minimal extension}
\index{minimal extension (of a tree potential)} 
\index{tree potential!minimal extension of} 
of $g$ to $\cT^o$ given by $p^*g=g\circ p$
is a positive tree potential on $\cT$.
Again the composition 
$\cP^+(\cS)\to\cP^+(\cT)\to\cP^+(\cS)$
equals the identity.

Finally we consider measures.
If $\rho$ is a positive measure on $\cS$,
then its pushforward $\imath_*\rho$,
\ie the extension of $\rho$ by zero,
is a positive measure on $\cT$.
Conversely, if $\rho$ is a positive measure on $\cT$, 
then the pushforward $p_*\rho$ is a positive 
measure on $\cS$.
This time, too, the composition
$\cM^+(\cS)\to\cM^+(\cT)\to\cM^+(\cS)$ is the identity.
\begin{proposition}
  All the mappings above are continuous and respect
  the homeomorphisms between $\cN^+$, $\cP^+$ and $\cM^+$
  given by Theorem~\ref{T407}.
\end{proposition}
In particular, $\cN^+(\cS)$, $\cP^+(\cS)$ and $\cM^+(\cS))$
are retracts of $\cN^+(\cT)$, $\cP^+(\cT)$ and $\cM^+(\cT))$,
respectively.
\begin{proof}
  That the mappings respect the isomorphisms is straightforward
  to verify and is left to the reader. 
  It is immediate from the definition of the weak topology on
  $\cN^+$ that the two mappings between $\cN^+(\cT)$ and 
  $\cN^+(\cS)$ are continuous. Hence continuity
  holds also in the cones $\cP^+$ and $\cM^+$.
\end{proof}
\begin{remark}
  The image measure $\rho_S\in\cM^+(\cS)$ 
  of a measure $\rho\in\cM^+(\cT)$ can be written as
  \begin{equation*}
    \rho_\cS
    =\rho|_\cS  +\sum_{\tau\in\cS}
    \left( \sum_{\vv\in T\tau\setminus T_\cS\tau}\rho(U(\vv))\right) \tau,
  \end{equation*}
  where $T\tau$ and $T_\cS\tau$ denote the
  tangent spaces of $\tau$ in $\cT$ and $\cS$, respectively,
  and where $U(\vv)$ as usual denotes the open subset of
  $\cT$ consisting of points in $\cT\setminus\{\tau\}$
  represented by $\vv$.
\end{remark}
%
%
%
%
\section{Inner products}\label{S416}
We have constructed three isomorphic Banach spaces
$\cM$, $\cN$ and $\cP$ associated to a given
complete, parameterized tree.
In this section we wish to equip these spaces with inner products.
As we show in Chapter~\ref{part-appli-analysis}, there are several
interesting interpretations of these inner products 
when working on the valuative tree $\cV$.
First, general complex atomic measures on $\cV$ can be viewed as
cohomology classes on the vo{\^u}te {\'e}toil{\'e}e, and the intersection
product on $\cM$ agrees with the cup product on cohomology.
Second, positive atomic measures $\rho$ on $\cV$ with integer 
coefficients correspond to integrally closed ideals $\rho_I$ 
in the ring $R$, and
the intersection product $\rho_I\cdot\rho_J$ gives the mixed 
multiplicity $e(I,J)$. An analytic version of the latter result
is studied in~\cite{pshfnts}, where the intersection product is
instead interpreted as a mixed Monge-Amp{\`e}re mass at the
origin.

We start by showing how to define the inner products on the positive
cones $\cM^+$, $\cN^+$ and $\cP^+$. The definition relies on a
parameterization of the tree even in the case of measures and
functions of bounded variation.  In order to extend this inner product
in the complex case, we need to impose some integrability conditions.
The upshot is that the inner products are well defined on subspaces
$\cM_0$, $\cN_0$ and $\cP_0$, and turn these into isometric
pre-Hilbert spaces.  Finally, we compare the topologies induced by the
inner products with the strong and weak topologies and show that the
pre-Hilbert spaces are not complete.
%
%
\subsection{Hausdorff measure}
The inner products on (subspaces) of $\cN$ and $\cP$
are defined by suitable integrals over $\cT^o$ with
respect to (one-dimensional) 
\emph{Hausdorff measure}\index{measure!Hausdorff},
whose definition we now recall.

We fix here a complete tree $\cT$ with an increasing parameterization
$\a:\cT\to[1,\infty]$ with $\a(\tau_0) =1$. The latter restricts to a
parameterization $\a:\cT^o\to[1,\infty[$ of the subtree $\cT^o$.
Consider the metric on $\cT^o$ associated to the parameterization as
in Section~\ref{tree-subsec}, \ie
$d(\tau,\tau')=\a(\tau)+\a(\tau')-2\a(\tau\wedge\tau')$, 
and let $\l$\index{$\l$ (Hausdorff measure)} be
the corresponding one-dimensional Hausdorff measure:
\begin{equation}\label{e461}
  \l\{A\} \= \lim_{\delta\to0}
  \inf\left\{
    \sum_1^\infty\diam(E_i)\ ;\ 
    A\subset\bigcup_1^\infty E_i,\ \diam(E_i)\le\delta
  \right\}
\end{equation}
for any subset $A\subset\cT^o$.  Here we use the convention that
$\inf\emptyset=\infty$ so that $\l(A)=\infty$ if $A$ cannot be covered
by countably many sets of diameter $\le\delta$ for any $\delta>0$.
From the remarks following Proposition~10.20 in~\cite{folland} it
follows that $\l$ restricts to a measure on the Borel $\sigma$-algebra
generated by the open sets in the topology associated to the metric
$d$.  As the latter topology is at least as strong as the weak
topology (see Proposition~\ref{P430}) we conclude that $\l$ restricts
to a (positive, weak) Borel measure on $\cT^o$.  We can therefore
integrate Borel measurable functions against $\l$.  In particular we
can integrate the functions $\cN^+$ and $\cP^+$.

Notice that, in general, the mass of $\l$ is infinite. That $\l \{ \cT
\} < + \infty$ implies in particular that $\cT$ has at most countably
many ends and is bounded for the metric $d$ introduced above. 
Note, moreover, that if $\cT$ has no branch points, 
then $\l$ is isomorphic to
Lebesgue measure on the interval $\a(\cT)\subset\R$.
%
%
\subsection{The positive case}
The inner product is defined first on the set of Dirac masses,
and then extended to $\cM^+$ by bilinearity.
\begin{definition}\label{def-innerprod}
  If $\tau,\tau'\in\cT$, then we define
  $\tau\cdot\tau':=\a(\tau\wedge\tau')\in[1,\infty]$.
\end{definition}
Note that this definition depends on the choice of the
parameterization.
\begin{remark}
This definition has a natural
  interpretation in the case when $\cT$ is the
  valuative tree $\cV$: then
  $\nu_{C}\cdot\nu_{C'}=\frac{C\cdot C'}{m(C)m(C')}$
  for any two irreducible curves $C$, $C'$. 
  See Section~\ref{val-curve}.
\end{remark}
\begin{definition}
  If $\rho,\rho'\in\cM^+$, then 
  we define $\rho\cdot\rho'\in[1,\infty]$ by
  \begin{equation}\label{e459}
    \rho\cdot\rho'
    :=\iint\limits_{\cT\times\cT}\tau\cdot\tau'\,
    d\rho(\tau)d\rho'(\tau').
  \end{equation}
  \index{inner product!on $\cM^+$}
\end{definition}
Notice that by Fubini's theorem 
and by~\eqref{e457} we have
\begin{equation}\label{e460}
  \rho\cdot\rho'
  =\int_\cT g\,d\rho'
  =\int_\cT g'\,d\rho,
\end{equation}
where $g,g'\in\cP^+$ are the positive tree potentials associated to
$\rho$ and $\rho'$, respectively. Here it is important to consider
the natural extension of these potentials from $\cT^o$ to $\cT$
for the formula to work also when $\rho$ or $\rho'$ 
charges an end.

\medskip
The inner product is defined on $\cN^+$, the cone of
nonnegative, left continuous, strongly decreasing functions on
$\cT^o$, by
\begin{equation}\label{e836}
  \langle f,f'\rangle 
  \= f(\tau_0)f'(\tau_0)
  +\int_{\cT^o}ff'\,d\l,
\end{equation}
\index{inner product!on $\cN^+$}
where $\l$ is 1-dimensional Hausdorff measure as above.  Finally we
define an inner product on $\cP^+$, the cone of positive tree potentials,
by passing to $\cN^+$ (see~\eqref{e451} and~\eqref{e452}):
\index{inner product!on $\cP^+$}
\begin{equation}
  \langle g,g'\rangle 
  \=g(\tau_0)g'(\tau_0)
  +\int_{\cT^o}\frac{dg}{d\a}\,\frac{dg'}{d\a}\,d\l.
\end{equation}
\begin{theorem}\label{T408}
  The isomorphisms
  \begin{equation*}
    d:\cN^+\to\cM^+
    \qand
    \Delta:\cP^+\to\cM^+
  \end{equation*}
  given in Theorem~\ref{T406} preserve the inner products
  defined above.
\end{theorem}
The proof relies on the following lemma of independent interest.
\begin{lemma}\label{L444}
  Let $\rho\in\cM^+$ be any positive Borel measure, 
  and let $f=f_\rho$
  be its associated functions in $\cN^+$. Then one can find a sequence
  of positive atomic measures $\rho_n$ tending weakly to $\rho$, 
  such that $f_n\=f_{\rho_n}$ \emph{increases} pointwise
  (and in fact uniformly) to $f$. 
\end{lemma}
\begin{proof}[Proof of Theorem~\ref{T408}]
  Pick positive measures $\rho,\rho'\in\cM^+$ and let 
  $f,f'\in\cN^+$ and $g,g'\in\cP^+$ be their preimages under
  $d$ and $\Delta$, respectively.
  It is clear from the definition that
  $\langle g,g'\rangle=\langle f,f'\rangle$.
  Hence it suffices to show that 
  $\rho\cdot\rho'=\langle f,f'\rangle$.

  First suppose $\rho,\rho'$ are Dirac masses at 
  $\tau$ and $\tau'$ respectively. Then 
  $\rho\cdot\rho'=\tau\cdot\tau'=\a(\tau\wedge\tau')$.
  Moreover, $f$ and $f'$ are the characteristic functions
  of the segments $[\tau_0,\tau]\cap\cT^o$ and
  $[\tau_0,\tau']\cap\cT^o$, respectively, so
  \begin{equation*}
    \langle f,f'\rangle
    =f(\tau_0)f'(\tau_0)
    +\int_{\cT^o}ff'\,d\lambda
    =1+\int_{\tau_0}^{\tau\wedge\tau'}d\a
    =\a(\tau\wedge\tau').
  \end{equation*}
  By bilinearity we conclude that
  $\rho\cdot\rho'=\langle f,f'\rangle$ holds when
  $\rho$ and $\rho'$ are positive atomic measures.
  
  In the general case, we rely on Lemma~\ref{L444}, and find positive
  atomic measures $\rho_n, \rho_m'$ whose corresponding functions
  $f_n, f_m'\in\cN^+$ increase to $f$ and $f'$ respectively. Note that
  by integration it follows that $g_n$, $g'_m$ increase to $g$ and
  $g'$ respectively.
  By monotone  convergence, it is clear that
  \begin{equation*}
    \lim_{m,n \to \infty} \langle f_n,f_m'\rangle
    =
    \lim f_n(\tau_0)f_m'(\tau_0)+\lim\int f_nf_m'\,d\l
    =
    \langle f,f'\rangle.
  \end{equation*}
  On the other hand, monotone convergence applied to $g_n$ and $g_m'$
  gives
  \begin{multline*}
    \rho_n\cdot\rho_m'
    =\int g_n\,d\rho_m'
    \mathop{\longrightarrow}\limits^{n\to \infty}    
    \int g\,d\rho_m'
    =
    \int g'_m\,d\rho
    \mathop{\longrightarrow}\limits^{m\to \infty}
    \int g' \, d\rho = \rho \cdot \rho'.
  \end{multline*}
  By what precedes, $\rho_n\cdot \rho_m' = \langle f_n,f_m'\rangle$ for
  all $n,m$. 
  Hence $\rho\cdot\rho'=\langle f,f'\rangle$, 
  which completes the proof.
\end{proof}
\begin{proof}[Proof of Lemma~\ref{L444}]
  Suppose first that the support of $f$ is included in a finite tree
  $\cS$.  For each $n$, pick a finite subset $B_n\subset\cS$
  containing the root, all branch points, and all ends of $\cS$, and
  such that $\rho(I)\le1/n$ for any connected component $I$ of
  $\cS\setminus B_n$.  For $\tau\in B_n$ set
  $m(\tau)=\rho\{\tau\}+\sum_I\rho(I)$, where the sum is over all
  $I\in\cI$ having $\tau$ as a \emph{left} endpoint.  Now define the
  positive atomic measure $\rho_n$ by $\rho_n=\sum_{\tau\in
    B}m(\tau)\tau$.  In other words, we slide the mass on each segment
  $I$ to its left endpoint.  Set $f_n=I\rho_n$.  Clearly
  $\rho_n$ has the same mass as $\rho$, hence $f_n(\tau_0) =
  f(\tau_0)$.  From the construction we have, for $\tau\in\cS$:
  \begin{equation*}
    f(\tau)-f_n(\tau)
    =\rho\{\sigma\ge\tau\}
    -\rho_n\{\sigma\ge\tau\}
    \in[0,1/n].
  \end{equation*}
  By choosing $B_n\subset B_{n+1}$ for all $n$, we get that $f_n$
  increases uniformly to $f$.
  
  In the general case, the set $\cS_n \= \{ f > 1/n \}$ is a finite
  tree. Define $f_n \= f \times \mathbf{1}_{\cS_n}$. This is a sequence of
  functions in $\cN^+$ increasing uniformly to $f$, and
  supported on a finite tree. We conclude by a diagonal argument. 
\end{proof}
%
%
\subsection{Properties}\label{S424}
It is clear that the inner product on $\cN^+$ (and on $\cP^+$) defined
by~\eqref{e836} satisfies the Cauchy-Schwartz inequality. 
An immediate consequence of Theorem~\ref{T408} above is
\begin{corollary}\label{C411}
  The inner product defined on $\cM^+$ satisfies the 
  Cauchy-Schwartz inequality:\index{Cauchy-Schwartz inequality}
  \begin{equation*}
    (\rho\cdot\rho')^2\le(\rho\cdot\rho)(\rho'\cdot\rho').
  \end{equation*}
\end{corollary}
For further reference, we note that in general the intersection
product is not (weakly) continuous on $\cM^+\times\cM^+$, and not even on the
set of Dirac masses $\cT\times\cT$. An example on the valuative tree
$\cT=\cV$ is given by $\nu_n=\nu_{y-nx,2}$ in local coordinates
$(x,y)$.  Here $\nu_n\cdot\nu_n=2$ but $\nu_n\to\nu_\fm$ and
$\nu_\fm\cdot\nu_\fm=1$.  

It is also clear that the inner product is not continuous for the
strong topology induced by the mass norm on $\cM^+$. 
The self-intersection of a Dirac mass at an end with infinite
parameter is infinite; see also Example~\ref{E411}.  
Note however that $\rho\cdot\rho'\ge\mass\rho\cdot\mass\rho'$.
\begin{proposition}\label{P436}
  The intersection product on $\cM^+\times\cM^+$ 
  is (weakly) lower semicontinuous:
  if $\rho_n\to\rho$, $\rho'_n\to\rho'$ weakly,
  then $\liminf\rho_n\cdot\rho'_n\ge\rho\cdot\rho'$.
\end{proposition}
The same results are also true in $\cP^+$ and $\cN^+$
\begin{proof}
  We may suppose the mass of all measures $\rho,\rho',\rho_n,\rho'_n$
  is equal to $1$.  Write $f = I \rho$, $f' = I\rho'$. The functions
  $f,f'\in\cN^+$ are uniformly bounded from above by $1$. 
  Note that 
  $\rho\cdot\rho'=1+\int ff'\,d\l=1+\int_0^1\l\{ff'>t\}\,dt$
  by Theorem~\ref{T408}. 
  Fix $\e>0$ arbitrarily small, and $k$ sufficiently
  large such that 
  $1+\int_{1/k}^1\l\{ff'>t\}\,dt\ge\rho\cdot\rho'-\e$ 
  if $\rho\cdot\rho'<+\infty$, or $\ge\e^{-1}$
  if $\rho \cdot \rho' = +\infty$. As $\{ ff' > 1/k\}\subset\{ f>
  1/k\}\cup \{f' >1/k\}$, we may find a finite tree $\cS_k$ containing
  $\{ ff' > 1/k\}$ (see Corollary~\ref{C408}).
  
  Define $f_n \= I\rho_n$, $f'_n \= I\rho'_n$. As $\rho_n$ tends
  weakly to $\rho$, $f_n$ tends to $f$ $\l$-almost everywhere. These
  functions are bounded uniformly from above by $1=\sup \| \rho_n \|$,
  hence $f_n \to f$ in $L^2 (\cS_k)$ also.  The same is true for
  $f'_n$.  We hence have
  \begin{equation*}
    \rho_n \cdot \rho'_n
    \ge
    1+ \int_{\cS_k} f_n f'_n\, d\l
    \mathop{\longrightarrow}\limits^{n\to \infty}
    1+\int_{\cS_k} f f' \, d\l.
  \end{equation*}
  When $\rho\cdot \rho'$ is finite, the right term is greater than
  $\rho\cdot\rho' -\e$. When $\rho\cdot\rho' = \infty$, it is greater
  than $\e^{-1}$. In any case, we conclude that $\liminf
  \rho_n\cdot\rho'_n \ge \rho \cdot \rho'$ by letting $\e\to0$.
\end{proof}
%
%
\subsection{The complex case}\label{S419}
We now wish to extend the previous definitions and results to
the complex case. In order to do this, we have to impose 
suitable integrability restrictions.
\begin{definition}{~}\label{D628}
  \begin{itemize}
  \item
    A complex measure $\rho$ belongs to $\cM_0$
    \index{$\cM_0$ (Hilbert subspace of $\cM$)} 
    iff its total 
    variation measure $|\rho|$ satisfies $|\rho|\cdot|\rho|<\infty$.
  \item A function $f \in \cN$ belongs to $\cN_0$
    \index{$\cN_0$ (Hilbert subspace of $\cN$)}  
    iff the function $T_f$, defined in~\eqref{e454},
    satisfies $\langle T_f,T_f\rangle<+\infty$.
  \item A function $g \in \cP$ belongs to $\cP_0$
    \index{$\cP_0$ (Hilbert subspace of $\cP$)}  
    iff $\delta g$ belongs to $\cN_0$.
\end{itemize}
\end{definition}
\begin{remark}\label{R428}
  Note that any atomic measure supported on the set $\{ \a < \infty\}$
  lies in $\cM_0$. In particular, in the valuative tree $\cV$
  (parameterized by skewness), any atomic measure supported on 
  quasimonomial valuations lies in $\cM_0$.
\end{remark}
Let us now describe how to define inner products on the sets
$\cM_0,\cN_0,\cP_0$.
\begin{proposition}\label{P474}
  For any measures $\rho,\rho'\in\cM_0$ the function
  $(\tau,\tau')\mapsto\tau\cdot\tau'$ lies in $L^1
  (\rho\otimes\overline{\rho'})$. 
  One can thus define\index{inner product!on $\cM_0$}
  \begin{equation}\label{e458}
    \rho\cdot\rho'
    :=\iint\limits_{\cT\times\cT}\tau\cdot\tau'\,
    d(\rho\otimes\overline{\rho'})(\tau,\tau').
  \end{equation}
\end{proposition}
\begin{proof}
  We have $\tau\cdot\tau' \in L^1(\rho\otimes\overline{\rho'})$ iff
  $\tau\cdot\tau' \in L^1(|\rho|\otimes|\rho'|)$. 
  By the Cauchy-Schwartz inequality (Corollary~\ref{C411}), 
  the latter condition is satisfied when
  $\tau\cdot\tau'\in L^1(|\rho|\otimes|\rho|)\cap 
  L^1(|\rho'|\otimes |\rho'|)$, \ie when $\rho,\rho'\in\cM_0$.
\end{proof}
\begin{proposition}\label{P473}
  For any  $f,f'\in\cN_0$, one has $f\overline{f}'\in L^1(\l)$. 
  One can thus define\index{inner product!on $\cN_0$}
  \begin{equation}\label{e453}
    \langle f,f'\rangle 
    \= f(\tau_0)\,\overline{f'(\tau_0)}
    +\int_{\cT^o} f\,\overline{f'}\,
    d\l.
  \end{equation}
  As a consequence, when $g,g'\in\cP_0$, 
  we can set\index{inner product!on $\cP_0$} 
  \begin{equation}\label{e456}
    \langle g,g'\rangle 
    \=g(\tau_0)\,\overline{g'(\tau_0)}
    +\int_{\cT^o}\frac{dg}{d\a}\,\overline{\frac{dg'}{d\a}}\,
    d\l.
  \end{equation}
\end{proposition}
\begin{proof}
  Pick $f,f'\in \cN_0$. Then $\langle T_f, T_f \rangle < \infty$,
  $\langle T_{f'}, T_{f'} \rangle < \infty$, so that
  $\langle T_f,T_{f'}\rangle<\infty$ by Corollary~\ref{C411}. 
  But $|f\overline{f}'|\le T_f T_{f'}$, implying
  $f\overline{f}'\in L^1(\l)$.
\end{proof}
\begin{remark}
  The space $\cM_0$ can be characterized in a slightly different way.
  Namely, one can show that $\rho \in \cM_0$ iff $\tau \cdot \tau' \in
  L^1 ( \rho \otimes \overline{\rho})$. In a similar way, a
  real-valued function $f\in\cN_0$ iff 
  $\langle f_1,f_1\rangle<\infty$, 
  $\langle f_2,f_2\rangle<\infty$, 
  where $f=f_1-f_2$ is the Jordan decomposition of $f$ 
  (see Proposition~\ref{P419}).
  
  On the other hand, a function $f\in\cN$ may be in $L^2(\l)$ but
  not in $\cN_0$, as the following example shows.
\end{remark}
\begin{example}
  Pick an end $\tau_\infty\in\cT$ with $\a(\tau_\infty)=\infty$
  (assuming such an end exists) and define increasing sequences
  $(\tau_k)_1^\infty$ and $(\tau'_n)_1^\infty$ by
  $\tau_k,\tau'_n<\tau_\infty$, $\a(\tau'_n)=2^n$ and
  $\a(\tau_k)=2^n+2^{-n}$.  Set $\rho_1=\sum_1^\infty 2^{-n/2}\tau_k$
  and $\rho_2=\sum_1^\infty 2^{-n/2}\tau'_n$.  Finally let
  $f_1=d^{-1}(\rho_1)$, $f_2=d^{-1}(\rho_2)$ and $f=f_1-f_2$.  Then it
  is straightforward to see that $f\in L^2(\l)$ but
  $f_1,f_2\notin L^2(\l)$, so that $f\not\in\cN_0$.
\end{example}
\begin{theorem}\label{T405}
  The three spaces $\cM_0$, $\cN_0$ and $\cP_0$ are vector spaces, 
  and the natural inner product defined 
  in~\eqref{e458},~\eqref{e453},~\eqref{e456} endow them with a 
  pre-Hilbert space structure.
  Further, the maps $d:\cN\to\cM$, $\Delta:\cP\to\cM$ 
  restrict to bijections
  \begin{equation*}
    d:\cN_0\to\cM_0
    \qand
    \Delta:\cP_0\to\cM_0
  \end{equation*}
  which preserve the inner products.
\end{theorem}
\begin{proof}
  It is clear by definition that 
  $\delta$ maps $\cP_0$ bijectively onto $\cN_0$ and
  preserves the inner product. Hence we need only
  consider $\cN_0$ and $\cM_0$.

  First note that $d$ maps $\cN_0$ onto $\cM_0$. Indeed:
  \begin{multline*}
    f\in \cN_0
    \Leftrightarrow 
    \langle T_f, T_f \rangle < \infty
    \mathop{\Longleftrightarrow}\limits^{\text{thm~\ref{T408}}}
    d(T_f) \cdot d(T_f) < \infty
    \\
    \mathop{\Leftrightarrow}\limits^{\text{prop~\ref{P418}}}
    |df| \cdot |df| < \infty
    \Leftrightarrow  df \in \cM_0.
  \end{multline*}
  It is also clear that $\cN_0$ is a vector space because 
  $T_{f+f'}\le T_f+T_{f'}$ for any functions 
  $f,f'\in \cN$. 
  Hence $\cM_0$ is also a vector space.
  
  We are thus reduced to proving that 
  $\rho\cdot\rho'=\langle I\rho,I\rho'\rangle$ 
  for any $\rho,\rho'\in\cM_0$. Note that by
  linearity, we need only check that 
  $\rho\cdot\rho=\langle I\rho,I\rho\rangle$ for all $\rho\in\cM_0$. 

  First assume that $\rho$ is a real measure. 
  Let $\rho=\rho_+-\rho_-$ be the Jordan decomposition of $\rho$.
  Thus $\rho_\pm$ are positive measures 
  and $|\rho|=\rho_++\rho_-$.
  It is clear that $\rho\in\cM_0$ implies 
  $\rho_\pm\cdot\rho_\pm<\infty$.
  By Theorem~\ref{T408}, we infer that
  \begin{equation*}
    \rho\cdot\rho
    =\rho_+\cdot\rho_++\rho_-\cdot\rho_--2\rho_+\cdot\rho_-
    =\langle f_+,f_+\rangle+\langle f_-,f_-\rangle-2\langle f_+,f_-\rangle
    =\langle f,f\rangle.
  \end{equation*}

  If $\rho=\rho_1+i\rho_2$ is a complex Borel measure in $\cM_0$,
  with $\rho_i$ real Borel measures,
  then $|\rho_j|\le|\rho|$, implying $\rho_j\in\cM_0$
  for $j=1,2$.
  By what precedes, $\rho_j\cdot\rho_j=\langle f_j,f_j\rangle$,
  for $j=1,2$, where $f_j=f_{\rho_j}$. 
  This gives 
  \begin{equation*}
    \rho\cdot\rho
    =\rho_1\cdot\rho_1+\rho_2\cdot\rho_2
    =\langle f_1,f_1\rangle+\langle f_2,f_2\rangle
    =\langle f,f\rangle,
  \end{equation*}
  and completes the proof. 
\end{proof}
%
%
\subsection{Topologies and completeness}
We end this section by briefly discussing the topology on
the three pre-Hilbert spaces $\cM_0$, $\cN_0$ and $\cP_0$
induced by the associated norms. Namely, we give examples
showing that the spaces are not complete with respect to these 
norms, and that the topologies are not comparable to the strong
topologies in general.

\begin{example}\label{E411}
  Assume that $\a$ is unbounded on $\cT$ and pick a sequence
  $(\tau_n)_1^\infty$ in $\cT^o$ with $\a(\tau_n)=n^2$.
  Then $\rho_n:=n^{-1}\tau_n$ tends to zero strongly
  but $\rho_n\cdot\rho_n=n$.
  On the other hand, if $\a$ is bounded on $\cT$, say 
  $\a\le C$, then it is clear that
  $\rho\cdot\rho\le C\Vert\rho\Vert^2$.
\end{example}
\begin{example}\label{E412}
  Pick any strictly increasing sequence
  $(\tau_n)_1^\infty$ such that $\a(\tau_n)$ is bounded
  and set $\rho_n=\sum_{j=1}^{2n}(-1)^j\tau_{2n+j}$. 
  Then clearly $\rho_n\in\cM_0$ and
  \begin{equation*}
    \rho_n\cdot\rho_n
    =\sum_{j=1}^{n}(\a(\tau_{2n+2j})-\a(\tau_{2n+2j-1}))
    <\a(\tau_{4n})-\a(\tau_{2n+1})
    \to0
  \end{equation*}
  as $n\to\infty$.
  On the other hand, $\rho_n$ has total variation $2n$ so
  $\rho_n\not\to0$ strongly.
\end{example}
Notice that the fact that $\Vert\rho_n\Vert$ is unbounded
implies that $\rho_n$ does not tend to zero weakly either.

A small modification of the same example shows that 
$\cM_0$ (and hence $\cN_0$ and $\cP_0$) is not complete. 
\begin{example}
  Let $(\tau_n)_1^\infty$ be as in Example~\ref{E412} and set
  $\rho_n=\sum_1^{2n}(-1)^{j-1}\tau_j$. 
  Then $\rho_n\in\cM_0$ and
  \begin{equation*}
    (\rho_{n+m}-\rho_n)\cdot(\rho_{n+m}-\rho_n)
    =\sum_{j=1}^{2m}(-1)^j\a(\tau_{2n+j})
    \le\a(\tau_{2n+2m})-\a(\tau_{2n+1})
    \to0
  \end{equation*}
  as $n,m\to\infty$.
  Hence $(\rho_n)_1^\infty$ is a Cauchy sequence,
  but there is no $\rho\in\cM$ such that 
  $(\rho_n-\rho)\cdot(\rho_n-\rho)\to0$.
  (Notice that $\rho_n$ has total variation $2n$.)
\end{example}
\begin{remark}
  In fact, the completion of $\cN_0$ 
  is naturally isomorphic to $\C\oplus L^2(\cT^o)$.
\end{remark}


%
%
%
%
%
%
\chapter{Applications of the tree analysis}
\label{part-appli-analysis}
This chapter is devoted to applications of the tree analysis
developed in the previous chapter. We shall use measures on the
valuative tree $\cV$ to describe singularities of ideals in
Section~\ref{sec-zar-ideal}, and cohomology classes of the 
vo\^ute \'etoil\'ee in Section~\ref{sec-voute}. 
Further applications, to singularities of plurisubharmonic functions
and to the dynamics of fixed point germs $f:(\C^2,0)\self$, 
will be explored in forthcoming papers:
see~\cite{pshfnts},~\cite{criterion},~\cite{eigenvaluation}.

Let us describe in more detail the content of this chapter.

We first attach to any ideal $I\subset R$ a 
\emph{tree transform} $g_I:\cVqm\to\R_+$,
\index{tree transform!of an ideal}
by setting $g_I(\nu)=\min\{\nu(\phi) \ ; \ \phi\in R\}$.  
This is a positive tree
potential in the sense of Section~\ref{S415}, hence the Laplacian of
$g_I$ is a well-defined positive measure $\rho_I=\Delta g_I$, called
the \emph{tree measure} of $I$.  We characterize the positive measures
on $\cV$ that are tree measures of ideals (Theorem~\ref{Tcharideal}).
Any tree measure $\rho_I$ for an ideal $I\subset R$ is atomic, and its
support coincides with the set of Rees valuations of $I$ when $I$ is
primary.  This gives a tree-theoretic approach to the Rees valuations
of a primary ideal and to Zariski's factorization of integrally closed
ideals (complete ideals in Zariski's terminology).

We next introduce the \emph{vo\^ute \'etoil\'ee} $\fX$. 
This space was first defined by Hironaka~\cite{hironaka} 
in a quite general context. 
In our setting, $\fX$ has a simpler description
than in the general case, and can be viewed as the total space of the
set of all blow-ups above the origin. Our aim is to 
describe the cohomology of this space, that is, the sheaf cohomology
$H^2(\fX,\C)$. In doing so, we were much inspired by the
monograph of Hubbard-Papadopol~\cite{hubbard}, where the toric case
was described in detail.

Let us summarize our approach.
The cohomology $H^2(\fX,\C)$ is a natural complex vector
space, endowed with an intersection form coming from the cup product.
We describe in Section~\ref{sec-ve-ident} a natural map 
sending a cohomology class $\omega\in H^2(\fX,\C)$ 
to a function $g_\omega$ defined on $\cVqm$, the set of
quasimonomial valuations.
As we show, $g_\omega$ is always a complex tree potential
in the sense of Section~\ref{S412}, and its Laplacian
$\rho_\omega=\Delta g_\omega$ is a complex atomic measure 
supported $\cVdiv$, the set of divisorial valuations.
Thus $\rho_\omega$ belongs to the subspace $\cM_0$ of
complex measures on which
we defined a inner product in Section~\ref{S416}.
We show that the mapping
$\omega\mapsto\rho_\omega$ gives an isometric embedding
of $H^2(\fX,\C)$ into $\cM_0$ (Theorem~\ref{T347}).
Finally we identify the images inside $\cM_0$ 
of the subspaces $H^2(\fX,\C)$, $H^2(\fX,\Z)$ and $H^2(\fX,\R)$, 
as well as two natural positive cones in $H^2(\fX,\R)$
(Theorem~\ref{T286}).
%
%
%
%
\section{Zariski's theory of complete ideals}\label{sec-zar-ideal}
%
%
\subsection{Basic properties}
We define the \emph{tree transform} 
\index{tree transform!of an ideal}
of an ideal $I\subset R$ as the function $g_I:\cVqm\to\R$
given by
\begin{equation}\label{e702-ideal}
  g_I(\nu)=\nu(I)=\min\{\nu(\phi)\ ;\ \phi\in I\}.
\end{equation}
We will show that $g_I$ is a positive tree potential 
in the sense of Definition~\ref{treepot}, where the valuative
tree $\cV$ is parameterized by \emph{skewness}.
Moreover, we will characterize all positive tree potentials 
that are of this form.

Let us start with the case of a principal ideal. 
\index{ideal!principal}
For $\phi\in\fm$, set $g_\phi(\nu)=\nu(\phi)$.
\begin{lemma}\label{Lideal-tree}
  The function $g_\phi:\cVqm\to[1,\infty)$ is a positive 
  tree potential
  with the property that if $\vv$ is any tangent vector in $\cV$,
  then $D_\vv g_\phi$ is an integer divisible by $m(\vv)$.
  Further, if $\phi$ is irreducible, 
  then $\Delta g_\phi=m(\phi)\nu_\phi$.
\end{lemma}
\begin{proof}
  By additivity in $\cP^+$ and unique factorization in $R$
  we may assume that $\phi$ is irreducible.
  Then Proposition~\ref{P201} and Definition~\ref{def-innerprod} 
  imply that
  \begin{equation*}
    g_\phi(\nu)
    =m(\phi)\,\a(\nu\wedge\nu_\phi)
    =m(\phi)\,\nu\cdot\nu_\phi,
  \end{equation*}
  hence~\eqref{e457} shows that $g_\phi=g_\rho$, where
  $\rho=m(\phi)\nu_\phi$. If $\vv\in T\nu$
  is a tangent vector, then $D_\vv g_\phi=0$ unless 
  $\nu\in[\nu_\fm,\nu_\phi]$ and $\vv$ is represented
  by either $\nu_\fm$; 
  or $\nu_\phi$, in which case 
  $D_\vv g_\phi=\pm m(\phi)$ and $m(\vv)$ is a factor
  of $m(\phi)$.
\end{proof}
We let $\cM^+_\cI$ be the set of positive measures $\rho\in\cM^+$ 
of the form
\begin{equation}\label{e703-ideal}
  \rho=\sum_{i=1}^sn_ib_i\nu_i.
\end{equation}
Here $1\le s<\infty$, $n_i$ are positive integers, 
$\nu_i$ is a divisorial or curve valuation, and
$b_i=b(\nu_i)$ is the generic multiplicity of $\nu_i$
if $\nu_i$ is divisorial and $b_i=m(\nu_i)$ is the
multiplicity of $\nu_i$ if $\nu_i$ is a curve
valuation.

To any positive measure $\rho\in\cM^+$ we associate an ideal 
$I_\rho\subset R$ by
\begin{equation}\label{e704-ideal}
  I_\rho=\{\phi\in R\ ;\ g_\phi\ge g_\rho\}.
\end{equation}
If $\rho=b\nu$, \ie $i=1$ and $n_1=1$ in~\eqref{e703-ideal},
then we write $I_\nu=I_\rho$.
\begin{theorem}\label{Tcharideal}
  The tree transform $g_I$ of any ideal $I\subset R$
  is a positive tree potential whose Laplacian $\rho_I=\Delta g_I$
  is a positive measure in $\cM^+_\cI$ of mass $m(I)$.
  Conversely, if $\rho\in\cM^+_\cI$, then the ideal
  $I=I_\rho$ has associated measure $\rho_I=\rho$.
\end{theorem}
\begin{definition}
  We call the positive measure $\rho_I$ the
  \emph{tree measure} of $I$.
\end{definition}
\begin{remark}
  If $R'$ is any subring of $R= \C[[x,y]]$ whose completion equals $R$,
  and $I$ is any ideal in $R'$, the curve valuations
  in~\eqref{e703-ideal} (if any) are associated to elements in $R'$. In
  particular, when $R'$ is the ring of convergent power series, any
  such curve valuation is analytic.
\end{remark}
\begin{remark}
  If $I$, $J$ are ideals, then $g_{IJ}=g_I+g_J$, and
  $g_{I+J}=\min\{g_I,g_J\}$. Moreover $g_{I\cap J}$ is the smallest 
  positive tree potential dominating $\max\{g_I,g_J\}$.
  
  One can rephrase the first two of these properties nicely in the
  terminology of semi-rings.  The set of ideals in $R$ defines a
  semi-ring, with addition $I+J$, and multiplication $I\cdot J$.  We
  may endow $\R_+$ with its tropical\index{tropical semi-ring}
  semi-ring structure: addition is given by $\min\{ a,b\}$, and
  multiplication by $a+b$.  This induces a semi-ring structure on the
  set of functions $\cVqm\to R_+$, and in particular to the set
  $\cP^+$ of positive tree potentials.  The properties above assert
  that the map $I\mapsto g_I$ is a semi-ring homomorphism.
\end{remark}
\begin{remark}
  Fix a composition of blow-ups $\pi\in\fB$, 
  and let $E_i$ be its exceptional components. 
  It is a classical problem to characterize the collections
  of integers $r_i$ which appear as multiplicities of some ideal $I$
  \ie such that $r_i=\div_{E_i}(\pi^*I)$ for all $i$.  
  A necessary and
  sufficient condition is that the $r_i$'s satisfy certain
  \emph{proximity inequalities}\index{proximity inequalities} (or relations). 
  For principal ideals
  this goes back to Enriques~\cite{enriques} 
  (see~\cite{casas} for a more recent
  presentation). The case of general ideals was treated by
  Lejeune-Jalabert~\cite{Lej91} and Lipman in~\cite{Lip94}. 
  Theorem~\ref{Tcharideal} can be viewed as a translation of these
  proximity relations into the tree language.
\end{remark}
\begin{proof}[Proof of Theorem~\ref{Tcharideal}]
  That $g_I=\inf_{\phi\in I}g_\phi$ is a positive 
  tree potential is a
  consequence of Lemma~\ref{Lideal-tree} and the fact that $\cP^+$ 
  is closed under infima (Corollary~\ref{C209}).
  Clearly $\rho_I$ has mass $g_I(\nu_\fm)=m(I)$.
  
  Let $S\subset I$ be a finite set of generators for $I$.
  Then $g_I=\min_{\phi\in S}g_\phi$. This implies
  that $g_I$ is supported on the smallest subtree of 
  $\cV$ containing $\nu_\fm$ and any $\nu_\psi$, where $\psi$
  ranges over the irreducible factors of the elements of $S$.
  This is a finite subtree $\cS$. Moreover, it follows from 
  Lemma~\ref{Lideal-tree} that on any segment in $\cS$ parameterized
  by skewness, $g_I$ is a piecewise affine function with integer
  slopes. Thus $\rho_I$ is an atomic measure supported on 
  valuations that are either ends or branch points in $\cS$,
  or regular points in $\cS$ where $g_I$ fails to be locally
  affine: see Proposition~\ref{Patom2}.
  From the integer slope property we
  conclude that $\rho_I=\sum_{i=1}^r\tn_i\nu_i$, where 
  $\nu_i$ are divisorial (\ie have rational skewness)
  or curve valuations and $\tn_i$ are positive integers. 
  
  We have to show that $b_i$ divides $\tn_i$.
  For this, it suffices to show that if
  $\nu$ is a curve or divisorial valuation, then 
  $\rho_I\{\nu\}$, which we now know is an integer, 
  is divisible
  by $b(\nu)$ in the case of a divisorial valuation, and
  by $m(\nu)$ in the case of a curve valuation.
  
  If $\nu$ is a curve valuation, then on $[\nu_\fm,\nu[$
  $\mu\mapsto\mu(\phi)$ is an affine function of skewness 
  with slope in $m(\nu)\N$ as $\mu\to\nu$. 
  Thus so is $g_I$, which implies that $\rho_I\{\nu\}\in m(\nu)\N$.

  If $\nu$ is divisorial we have to work a bit harder. 
  We may assume that $b(\nu)>1$, so that in particular
  $\nu\ne\nu_\fm$. The proof relies on the following lemma.
  \begin{lemma}\label{L707-ideal}
    Let $\nu\ne\nu_\fm$ be a divisorial valuation with approximating
    sequence $\nu_\fm=\nu_0<\nu_1<\dots<\nu_g<\nu_{g+1}=\nu$
    as in Proposition~\ref{prop-approx} 
    and let $\phi\in\cC$. Assume that $\nu_\phi$ and $\nu_\fm$
    represent the same tangent vector at $\nu$.
    Then $\nu(\phi)\in\sum_{i=1}^g\N m_i\a_i$.
  \end{lemma}
  We continue the proof of the theorem. 
  Recall that $\rho_I\{\nu\}=-\sum_{\vv\in T\nu}D_\vv g_I$.
  Lemma~\ref{Lideal-tree} implies that $D_\vv g_I\in m(\vv)\N$ for every
  $\vv\in T\nu$. If $b(\nu)=m(\nu)$ then $m(\vv)=b(\nu)$
  for every $\vv$ and we are done, so suppose
  that $b(\nu)>m(\nu)$. 
  
  Let $\vv_-\in T\nu$ be the tangent vector
  represented by $\nu_\fm$. There is then a unique tangent vector
  $\vv_+\in T\nu$, $\vv\ne\vv_-$ such that $m(\vv)=m(\nu)$; for all
  other $\vv$ we have $m(\vv)=b(\nu)$. It hence suffices to show that
  $D_{\vv_+}g_I+D_{\vv_-}g_I\in b(\nu)\N$. Moreover, we may find
  $\psi_\pm\in S$ such that $g_I(\mu)=\mu(\psi_\pm)$ as 
  $\mu\to\nu$, $\mu\in U(\vv_\pm)$. Then $\nu(\psi_+)=\nu(\psi_-)$
  and we have to show that
  $D_{\vv_+}g_{\psi_+}+D_{\vv_-}g_{\psi_-}\in b(\nu)\N$.

  Write $\psi_\pm=\psi_\pm'\psi_\pm''\psi_\pm'''$, where 
  $\psi_\pm'$ ($\psi_\pm''$) is
  the product of all irreducible factors representing 
  $\vv_-$ ($\vv_+$). 
  Then
  \begin{equation*}
    0
    =\nu(\psi_+)-\nu(\psi_-)
    =\nu(\psi_+')-\nu(\psi_-')
    +\alpha(\nu)(m(\psi_+''\psi_+''')-m(\psi_-''\psi_-''')).
  \end{equation*}
  By Lemma~\ref{L707-ideal} we get that
  $\a(\nu)(m(\psi_+''\psi_+''')-m(\psi_-''\psi_-'''))
  \in\sum_{i=1}^g\N m_i\a_i$.
  This implies that 
  $m(\psi_+''\psi_+''')-m(\psi''_-\psi_-''')\in b(\nu)\N$
   by Proposition~\ref{P403}.
  But we always have $m(\psi_\pm''')\in b(\nu)\N$ so we conclude that
  $m(\psi_+'')-m(\psi_-'')\in b(\nu)\N$.
  Finally this gives
  \begin{equation*}
    D_{\vv_+}g_{\psi_+}+D_{\vv_-}g_{\psi_-}
    =m(\psi_+'')-m(\psi_-'')-m(\psi_-''')
    \in b(\nu)\N,
  \end{equation*}
  which completes the proof that $\rho_I\in\cM^+_\cI$.
  
  Conversely, if $\rho\in\cM^+_\cI$, define $I=I_\rho$ by~\eqref{e704-ideal}
  and let $\cS$ be the support of $g_\rho$ 
  (thus $\cS$ is a finite tree: see Section~\ref{S412}).
  Clearly $g_I=\inf_{\phi\in I}g_\phi\ge g_\rho$. 
  For the reverse inequality we pick
  irreducible elements $\psi_{ij}\in\fm$, $1\le i\le s$, $1\le j\le n_i$
  as follows. Write $\mu_{ij}=\nu_{\psi_{ij}}$.
  If $\nu_i$ is a curve valuation, then $\mu_{ij}=\nu_i$ for all $j$.
  If $\nu_i$ is divisorial, then $\mu_{ij}>\nu_i$, $m(\mu_{ij})=b_i$
  and $\mu_{ij}$ represent distinct tangent vectors at $\nu_i$,
  none of which is in $T_\cS\nu_i$.

  Write $\psi=\prod_{i,j}\psi_{ij}$. 
  It is then straightforward to verify that $g_\psi=g_\rho$ on $\cS$.
  Hence $\psi\in I$. Now consider $\nu\not\in\cS$
  and let $\nu_0=\max\cS\cap[\nu_\fm,\nu]$. If $\nu_0\ne\nu_i$ for
  all $i$, then $g_\psi(\nu)=g_\psi(\nu_0)=g_\rho(\nu_0)=g_\rho(\nu)$
  for any choice of $\psi$. If $\nu_0=\nu_i$, then we pick $\psi_{ij}$
  such that no $\mu_{ij}$ represent the same tangent vector as $\nu_0$
  at $\nu_i$. Again $g_\psi(\nu)=g_\rho(\nu)$. 
  Hence $g_I=g_\rho$, which implies $\rho_I=\rho$.
  The proof is complete.
\end{proof}
\begin{proof}[Proof of Lemma~\ref{L707-ideal}]
  Set $\mu=\nu_\phi\wedge\nu$ and write $\mu\in[\nu_j,\nu_{j+1}[$ for
  some $0\le j\le g$.
  Then $b(\mu)$ divides $m(\phi)$ and
  $b(\mu)\alpha(\mu)\in\sum_{i\le j}\N m_i\a_i\subset\sum_{i=1}^g\N m_i\a_i$.
  This proves the lemma as 
  $\nu(\phi)=m(\phi)\alpha(\mu)$.
\end{proof}
%
%
\subsection{Normalized blow-up}
The proof of Theorem~\ref{Tcharideal} was based on tree arguments.  
We now follow~\cite{teissier-mult}, and use the classical normalized
blow-up of an ideal to describe the geometric structure of an ideal of
the form $I_\rho$.

Recall that an ideal $I$ is \emph{primary}\index{ideal!primary}
iff $I\supset\fm^n$ for some $n\ge1$. 
When $I$ is primary, denote by $\pi:X\to (\C^2,0)$ the
normalization of the blowup along $I$. The exceptional components
$E_i$ of $\pi$ are associated to divisorial valuations called the 
\emph{Rees valuations}\index{valuation!Rees} of $I$.
\begin{proposition}\label{P0-ideal}
  Let $\rho=\sum n_ib_i\nu_i\in\cM_{\cI}$. Then
  $I_\rho=\prod I_{\nu_i}^{n_i}$.
  Moreover, $I_\rho$ is primary iff all $\nu_i$ are divisorial, in
  which case the latter coincide with the Rees valuations of $I$.
\end{proposition}
\begin{proof}
  Pick $\rho\in\cM^+_\cI$, and write $I=I_\rho$. Assume first that $I$
  is primary. Then $I\supset\fm^n$ for some $n$ so the tree
  transform of $I$ is bounded by $n$. Hence all the valuations
  $\nu_i$ are divisorial. Let us show that
  they coincide with the Rees valuations of $I$, and that $I=\prod
  I_{\nu_i}^{n_i}$.

  Let $\pi:X\to (\C^2,0)$ be the normalization of the blowup along $I$,
  and $E_i$ be the set of exceptional components of $\pi$.  Let
  $(\psi_0,\cdots, \psi_N)$ be a finite set of generators for $I$.  For
  $a\in \C^{N+1}$ write $\psi_a=\sum
  a_i\psi_i$, and let $V_a$ be the strict transform of $\psi_a^{-1}(0)$
  by $\pi$.
  
  The structure of $V_a$ is described in~\cite[p.\@ 332]{teissier-mult}.
  We may find integers $n'_i\ge1$ such that outside a proper closed
  Zariski subset $Z\subset\P^N$, $V_a$ is a union of smooth curves
  $V^{ij}_a$, $1\le j\le n_i'$. Moreover each $V^{ij}_a$
  intersects $E_i$ transversely at a smooth point 
  and $V_a\cap V_b=\emptyset$ for $a\ne b\in Z$.
  
  Let $\mu_i$ be the divisorial valuation associated to $E_i$,
  $\nu_{ij}^a$ the valuation associated to the irreducible curve
  $\pi(V^{ij}_a)$, and $\psi_{ij}\in\fm$ the irreducible element
  attached to $\pi(V^{ij}_a)$.  Pick $\phi\in\fm$ irreducible. Then
  $g_{\mu_i}(\phi)$ is equal to $b(\mu_i)$ times the order of
  vanishing of $\pi^*\phi$ along $E_i$. On the other hand,
  $\nu_{ij}^a(\phi)$ is equal to $m(V^{ij}_a)^{-1}$ times the order of
  vanishing of $\phi \circ h (t)$ where $ t\to h(t)\in \C^2$ is a
  parameterization of $\pi(V^{ij}_a)$. As $V^{ij}_a$ is smooth and
  transverse to $E_i$, $\pi(V^{ij}_a)$ is a curvette for $\nu_{ij}^a$
  in the sense of Section~\ref{sec-curvette}. By
  Proposition~\ref{P410}, its multiplicity equals $b(\mu_i)$.  The
  parameterization $h$ can be obtained by composing a parameterization
  of $V^{ij}_a$ with $\pi$,  hence  $\nu_{ij}^a(\phi)$ is equal 
  $b(\mu_i)$ times the intersection product of $\pi^*\phi^{-1}(0)$
  with $V^{ij}_a$. Whence $\nu_{ij}^a(\phi)\ge
  g_{\mu_i}(\phi)$, with equality when the strict transform of
  $\phi^{-1}(0)$ does not intersect $V^{ij}_a$.  In particular,
  $\psi^{ij}_a\in I_{\mu_i}$ for all $j$ and all $a\not\in Z$. When
  $a\ne b$, $V^{ij}_a\cap V^{ij}_b=\emptyset$ hence
  $g_{\mu_i}(\phi)=\min\{\nu_{ij}^a(\phi),\nu_{ij}^b(\phi)\}$, so that
  $g_{\mu_i}=\min\{\nu_{ij}^a\ ;\ a\ne Z,\ 1\le j\le n_i'\}$.  
  We infer that the tree transform $\min\{g_{\psi_a}\ ;\ a\not\in Z\}$
  coincides with the positive tree potential 
  associated to the measure $\sum n_i'b(\mu_i)\mu_i$.
  
  But $I$ is generated by a finite number of sufficiently generic
  elements $\psi_{a_i}$, $a_i\not\in Z$ so we conclude that 
  $\sum n_ib_i\nu_i=\sum n_i'b(\mu_i)\mu_i$.  
  In particular the $\nu_i$'s are divisorial and are 
  exactly the Rees valuation of $I$.
  
  We choose $a_i\ne Z$ such that $\psi_{a_i}$ belongs to $\prod
  I_{\mu_i}^{n'_i} =\prod I_{\nu_i}^{n_i}$, hence $I\subset \prod
  I_{\nu_i}^{n_i}$.  Conversely $\psi\in \prod I_{\nu_i}^{n_i}$ implies
  $\nu(\psi)\ge \sum n_i g_{I_{\nu_i}}(\nu) =g_\rho(\nu)$, hence $I=
  \prod I_{\nu_i}^{n_i}$.

  Suppose now $I$ is not primary. We will show that
  $I=\prod I_{\nu_i}^{n_i}$, and 
  that some $\nu_i$ are curve valuations.
  Write $\rho_I=\sum n_ib_i\nu_i\in\cM^+_\cI$.
  Suppose $\nu_i$ is divisorial for $i\le r$ and 
  a curve valuation, $\nu_i=\nu_{\phi_i}$ for $i>r$.
  Define $\rho'=\sum_{i\le r}n_ib_i\nu_i$,
  and $I'=I_{\rho'}$. 
  We have proved that $I'$ is a primary ideal equal to
  $\prod_{i\le r}I_{\nu_i}^{n_i}$.
  
  Pick $\psi\in I'\cdot\prod_{i>r}(\phi_i)^{n_i}$. The tree transform of
  $\nu_\psi$ clearly dominates $g_{I}$, hence $\psi\in I$.  Conversely
  suppose $g_\psi(\nu)\ge g_I(\nu)$ for all $\nu$. 
  Then letting $\nu\to\nu_{\phi_i}$ we infer that the 
  mass of $\Delta g_\psi$ at
  $\nu_{\phi_i}$ is greater than $n_i$ hence $\phi_i^{n_i}$ divides
  $\psi$. We may then write $\psi=\psi'\prod_{i>r}\phi_i^{n_i}$, and
  clearly $g_{\psi'}\ge g_{I'}$, whence $\psi'\in I'$. We have proved
  $I=I'\,\prod_{i>r}(\phi_i)^{n_i}=\prod_i I_{\nu_i}^{n_i}$.
  
  This concludes the proof.
\end{proof}
%
%
\subsection{Integral closures}
The mapping $I\mapsto\rho_I$ is not injective in general.  For
instance, the ideals $\langle x^2,y^2\rangle$ and $\langle
x^2,xy,y^2\rangle$ both have tree measure $2\nu_\fm$.  However, the
lack of injectivity can be well understood.  Recall that the
\emph{integral closure}\index{ideal!integral closure of} $\bI$ of $I$ is the
set of $\phi\in R$ such that $\phi^n+a_1\phi^{n-1}+\dots+a_n=0$ for
some $n\ge1$ and $a_i\in I^i$. Then $\bar{\bI}=\bI$ and $I$ is
\emph{integrally closed} if $\bI=I$.  We have the following classical
result (see~\cite[p.~350]{ZS}), rephrased in our language.
\begin{proposition}\label{P702-ideal}
  For any ideal $I\subset R$ we have
  $\bI=\{\phi\in R\ ;\ g_\phi\ge g_I\}$.
\end{proposition}
\begin{remark}\label{R703-ideal}
  Fix an ideal $I$ and $\phi\in\fm$.  Suppose we want to show that
  $\phi\in\bI$.  By Proposition~\ref{P702-ideal} we must show that
  $g_\phi\ge g_I$.  Since $g_I$ is locally constant outside
  $\cT_I=\supp g_I$ it suffices to show that $g_\phi\ge g_I$ on
  $\cT_I$.  Write $\rho_I=\sum n_ib_i\nu_i$ as before.  Assume that
  $\nu_i$ is divisorial for $i\le r$ and a curve valuation for $r<i\le s$.  
  Then $g_I$ is affine and $g_\phi$ concave on any
  segment in $\cT_I\setminus\supp\rho_I$.  Hence it suffices to check
  that $\nu_i(\phi)\ge\nu_i(I)$ for $i\le r$ and
  $\lim\a(\nu)^{-1}(\nu(\phi)-\nu(I))\ge0$ as $\nu\to\nu_i$ for $i>r$.
  
  In particular, if $I$ is primary, then it suffices to check that
  $g_\phi\ge g_I$ at the Rees valuations of $I$ 
  (see \eg~\cite[p.333]{teissier-mult}).
\end{remark}
\begin{corollary}\label{C701-ideal}
  If $\rho\in\cM^+_\cI$, then $I_\rho$ is integrally closed.
\end{corollary}
Let $\cI$ be the set of integrally closed ideals in $R$.
Theorem~\ref{Tcharideal} and Corollary~\ref{C701-ideal} imply that the mapping
$\rho\mapsto I_\rho$ gives a bijection between $\cM^+_\cI$ and $\cI$, with
inverse is given by $I\mapsto\rho_I$.  Now $g_{IJ}=g_I+g_J$ for any
ideals $I$, $J$, hence $\rho_{IJ}=\rho_I+\rho_J$.  We get from this
that if $I$ is integrally closed and $\rho=\rho_I=\sum n_ib_i\nu_i$,
then $I=I_\rho= I_{\nu_1}^{n_1}\cdots I_{\nu_s}^{n_s}$ by
Proposition~\ref{P0-ideal}.  This leads to the following result that
incorporates Zariski's celebrated decomposition of complete 
\index{ideal!complete}
(=integrally closed) ideals (see~\cite{ZS}).
\begin{theorem}\label{T702-ideal}
  The set $\cI$ of integrally closed ideals in $R$ is a semigroup
  under multiplication and the mappings
  \begin{equation*}
    \cI\ni I\mapsto\rho_I\in\cM^+_\cI
    \qquad\text{and}\qquad
    \cM^+_\cI\ni\rho\mapsto I_\rho\in\cI
  \end{equation*}
  define inverse semigroup isomorphisms.

  In particular, every $I\in\cI$ has a unique
  factorization $I=I_{\nu_1}^{n_1}\cdots I_{\nu_s}^{n_s}$, where 
  $n_i$ are positive integers and $\nu_i$ are
  divisorial or curve valuations.
\end{theorem}
%
%
\subsection{Multiplicities}
Suppose $I$ and $J$ are primary ideals in $R$.  We define the
\emph{mixed multiplicity}
\index{multiplicity!mixed (of an ideal)}
$e(I,J)\in\N^*$ following~\cite[Chapter~I,~\S2]{teissier0} (when $I=J$
this gives the \emph{multiplicity} $e(I)$ of $I$ in the classical
sense).  For all $n,m\ge0$, the vector space $R/(I^nJ^m)$ is finite
dimensional over $\C$, and we write
$\dim_k(R/(I^nJ^m))=\frac12n^2e(I)+nme(I,J)+\frac12m^2e(J)+O(n+m)$.
The multiplicity $e(I,J)$ coincides with the multiplicity of their
respective integral closures (\eg~\cite[Chapter~0,~\S0.6]{teissier0}).

Here we show how to compute mixed multiplicities in
terms of the associated tree measures. 
Recall from Section~\ref{S416} the definition of the
inner product on the cone $\cM^+$ of positive Borel measures on $\cV$.
We then have
\begin{theorem}\label{cor-mult-ideal}
  For any primary ideals $I,J$ we have
  \begin{equation}\label{e34-ideal}
    e(I,J)
    =\rho_I\cdot\rho_J
    =\iint\limits_{\cV\times\cV}\mu\cdot\nu\,d\rho_I(\mu)d\rho_J(\nu)
    =\sum_{i,j}n_ib_im_jc_j\,\mu_i\cdot\nu_j,
  \end{equation}
  where $\rho_I=\sum_1^pn_ib_i\mu_i$ 
  and $\rho_J=\sum_1^qm_jc_j\nu_j$ 
  are the measures associated to $I$ and $J$, respectively
  (see Proposition~\ref{P0-ideal}), and $b_i, c_j$ are the generic 
  multiplicities of $\nu_i$ and $\mu_j$ respectively.
\end{theorem}
\begin{remark}
  In~\cite{pshfnts} we shall prove an analytic version of
  Theorem~\ref{cor-mult-ideal} where $I$ and $J$
  are replaced by plurisubharmonic functions $u$ and $v$, and
  where $e(I,J)$ is replaced by the mass of the mixed
  Monge-Amp\`ere measure $dd^cu\wedge dd^cv$ at the origin. 
\end{remark}
\begin{proof}
  We may assume that $I$ and $J$ are integrally closed.
  Fix a finite set of generators for $I$ and $J$, say
  $I=\langle\phi_i\rangle$, $J=\langle\psi_j\rangle$. 
    By~\cite[Chapter~I,~\S2]{teissier0} 
  $e(I,J)$ is equal to the intersection
  multiplicity of $\{\phi_\a=\sum\a_i\phi_i=0\}$ with 
  $\{\psi_\b=\sum \b_j\psi_j=0\}$ if the 
  coefficients $(\a_i)$, $(\b_j)$ are sufficiently generic.
  
Introduce a composition of blow-ups $\pi$ dominating the blow-ups of
both ideals $I$ and $J$, \ie such that all Rees valuations of $I$ and
$J$ are determined by some exceptional divisor of $\pi$. Decompose
the strict transform of $\{\phi_\a=0\}$ into irreducible components
$V^{is}_\a$, $1\le i\le p$, $1\le s\le n_i$ as in the proof of
Proposition~\ref{P0-ideal}. Each $V^{is}_\a$ corresponds to a curve
valuation $\mu^{is}_\a$ dominating $\mu_i$, of multiplicity
$b_i=b(\mu_i)$, and the tangent vectors defined by all $\mu^{is}_\a$
at $\mu_i$ are distinct.  We do the same decomposition of the strict
transform of $\{\psi_\b=0\}$ as a union of irreducible germs
$W^{jt}_\b$, with associated valuations $\nu^{jt}_\b$.  We note that
if $\mu_i=\nu_j$ for some $i,j$, then for sufficiently generic
$\a,\b$, the tangent vectors of $\mu^{jt}_\b$ and $\nu^{is}_\a$ at
$\mu_i=\nu_j$ are all distinct. We infer $V^{is}_\a\cdot
W^{jt}_\b=b_ic_j\,\mu_i\cdot\nu_j$.  Formula~\eqref{e34-ideal} now
follows immediately by bilinearity since
  \begin{equation*}
    \{\phi_\a=0\}\cdot\{\psi_\b=0\}
    =\sum_{i,j,s,t}V^{is}_\a\cdot W^{jt}_\b
    =\sum_{i,j}b_ic_jn_im_j\,\mu_i\cdot\nu_j.
  \end{equation*}
  The proof is complete.
\end{proof}
%
%
%
%
\section{The vo\^ute \'etoil\'ee}\label{sec-voute}
We now turn to the vo\^ute \'etoil\'ee $\fX$.
Our objective is to analyze its cohomology
$H^2(\fX,\C)$ in terms of complex Borel
measures on the valuative tree.
We start by giving the definition and basic properties of
$\fX$, of its cohomology $H^2(\fX,\C)$, and of the inner product
on cohomology induced from the cup product.
In Section~\ref{sec-ve-ident} we then show 
how each cohomology class $\omega$ defines a 
function $g_\omega:\cVqm\to\C$.
The function $g_\omega$ turns out to be
a complex tree potential in the sense of Chapter~\ref{part-potent},
so its Laplacian $\rho_\omega=\Delta g_\omega$ 
is a well-defined complex Borel measure on $\cV$.
In fact, the map $\omega\to\rho_\omega$ gives
an isomorphism between $H^2(\fX,\C)$ and the set
of complex atomic measure supported on divisorial
valuations.
Moreover, as we prove in Section~\ref{sec-isomembedd},
this isomorphism preserves the inner product in the sense
that $-\omega\cdot\omega'=\rho_\omega\cdot\rho_\omega$.
We end the chapter by describing the images of various
subsets of $H^2(\fX,\C)$ under the isomorphism: 
see Theorem~\ref{T286}.

Throughout the section, we shall make essential use of
the fundamental isomorphism from Chapter~\ref{A3}
between the valuative tree and the universal dual graph.
%
%
\subsection{Definition}\label{sec-ve-def}
We shall denote by $\fB$
\index{$\fB$ (modifications above the origin)} 
the set of all blow-ups above the origin in $\C^2$. 
We let $X_\pi$ be the total
space of a fixed element $\pi\in\fB$ so that 
$\pi:X_\pi\to(\C^2,0)$. 
We saw in Chapter~\ref{A3} that $\fB$ forms an inverse system.
\begin{definition}
  The vo\^ute \'etoil\'ee\index{vo\^ute \'etoil\'ee}
  \index{$\fX$ (vo\^ute \'etoil\'ee)} 
  is the projective limit
  \begin{equation*}
    \fX \= \projlim_{\pi\in\fB} X_\pi.
  \end{equation*}
\end{definition}
Each $X_\pi$ is an algebraic variety, hence $\fX$ is naturally a
pro-algebraic variety. We endow it with the topology induced by the
product topology from the natural embedding of $\fX$ into the product
$\prod X_\pi$. There is a natural proper projection map
$\fX\to(\C^2,0)$.  The space $\fX$ is not algebraic (nor even a
topological manifold): we shall see that its second cohomology group
is infinite dimensional.

Let us quickly indicate why our definition is equivalent to the usual
 one, given for instance in~\cite{hironaka}. Any
 "\'etoile"\footnote{An \'etoile is a collection of finite composition
 of local blow-ups, an element is hence a map $\varpi : U \to
 (\C^2,0)$, where $U$ is some analytic space.} in the sense of
 Hironaka has a well-defined "center" in $X_\pi$ for any $\pi\in\fB$:
 take the intersection of all images $\varpi(U)$ over all $(\varpi,U)$
 in the \'etoile.  As we are in dimension two, this center is always a
 (closed) point. This gives a natural map from the set of "\'etoiles"
 to $\fX$, which is easily seen to be bijective, and also
 bicontinuous.

%
%
\subsection{Cohomology}\label{sec-ve-co}
Consider $\pi,\pi'\in\fB$ such that $\pi'=\pi\circ\varpi$ for some
modification $\varpi$. The map $\varpi$ induces a map between
cohomology groups $\varpi^*:H^2(X_{\pi},\C)\to
H^2(X_{\pi'},\C)$.\footnote{When a class $\omega$ is represented by a
smooth form $\a$, $\varpi^*\omega$ is the class of $\varpi^*\a$.}
By~\cite[Corollary 14.6]{bredon} (see also~\cite{kiyek}), 
the sheaf cohomology
\index{vo\^ute \'etoil\'ee!cohomology of} 
of the constant sheaf $\C$ on the projective limit $\fX$ is equal to
\begin{equation}\label{e303}
  H^2(\fX,\C)=\injlim_{\pi\in\fB} H^2(X_\pi,\C).
\end{equation}
The choice of the base field is essentially irrelevant in the sequel.
One may replace $\C$ by $\R$ or $\Q$ or even $\Z$. 
We restrict our attention to $H^2$ 
as all the other cohomology groups are easy to compute.

\medskip 
In order to proceed further, we need to describe more precisely the
cohomology $H^2(X_{\pi},\C)$ and the pullback $\varpi^*$ 
discussed above.
We shall use some of the notation from Chapter~\ref{A3}.
Specifically, we let $\Gast_\pi$ be the (partially ordered)
set consisting of irreducible components of the
exceptional divisor $\pi^{-1}(0)$.

The following result is classical; 
we refer to~\cite[p.473-474]{GH} for a proof.
\begin{proposition}{~}\label{P98}
  \begin{itemize}
  \item
    Each $E\in\Gast_\pi$ induces a natural class 
    $[E]\in H^2(X_{\pi},\C)$.
    The cohomology group $H^2(X_{\pi},\C)$ is isomorphic to the 
    direct sum of $\C[E]$ over $E\in\Gast_\pi$.
 \item
   Suppose $\varpi:X_{\pi'}\to X_\pi$ is the blow-up at 
   \emph{one point} $p\in\pi^{-1}(0)$, 
   and let $E\=\varpi^{-1}\{p\}$. 
   Then $\varpi^*:H^2(X_\pi,\C)\to H^2(X_{\pi'},\C)$ is an
   injective map and we have
   \begin{equation*}
     H^2(X_{\pi'},\C)=\varpi^*H^2(X_{\pi},\C)\oplus\C[E].
   \end{equation*}
 \end{itemize}
\end{proposition}
\begin{remark}\label{R302}
  If $\pi\in\fB$ and $E\in\Gast_\pi$
  (\ie $E$ is an irreducible component of $\pi^{-1}(0)$) 
  then there are three natural objects associated to $E$: 
  a cohomology class $[E]\in H^2(X_\pi,\C)$;
  an element $E$ in the universal dual graph $\Gast$;
  and a divisorial valuation $\nu_E\in\cVdiv$.
  The last two of these are independent of $\pi$ 
  (as long as $E\in\Gast_\pi$) but the cohomology class
  \emph{does} depend on $\pi$. 
  Indeed, if $\pi'=\pi\circ\varpi$ for some modification $\varpi$,
  then the image of $[E]$ in $H^2(X_{\pi'},\C)$ corresponds to
  the \emph{total} transform of $E$ 
  whereas the image of $E$ in $\Gast_{\pi'}$ is the
  \emph{strict} transform.
\end{remark}
From Proposition~\ref{P98} and~\eqref{e303} we infer:
\begin{itemize}
\item
  The natural map $\imath_\pi$ from $H^2(X_\pi,\C)$ to
  $H^2(\fX,\C)$ is injective for any $\pi\in\fB$.
\item
  For each element $\omega\in H^2(\fX,\C)$, there exists
  $\pi\in\fB$, and $\omega_\pi\in H^2(X_\pi,\C)$ such that
  $\imath_\pi\omega_\pi=\omega$. Moreover,
  $\omega_\pi=\sum_{E\in\Gast_\pi}a(E)[E]$, where $a(E)\in\C$.
\item
  Two elements $\omega_1\in H^2(X_{\pi_1},\C)$ and $\omega_2\in
  H^2(X_{\pi_2},\C)$ determine the same element $\omega\in
  H^2(\fX,\C)$ iff there exists $\pi\in\fB$ and
  modifications $\varpi_1 : X_\pi \to X_{\pi_1}$, $\varpi_2 : X_\pi \to
  X_{\pi_2}$ such that $\varpi^*_1\omega_1 = \varpi^*_2 \omega_2$.
\end{itemize}
\begin{corollary}\label{C99}
  The set $H^2(\fX,\C)$ is an infinite
  dimensional vector space. 
\end{corollary}
%
%
\subsection{Intersection product}\label{sec-ve-int}
Each complex vector space $H^2(X_\pi,\C)$ is endowed with a natural
hermitian form, the cup product. 
If $E,E'\in\Gast_\pi$, then 
$[E]\cdot[E']$ is by definition the intersection product of the 
curves $E$ and $E'$.
For two arbitrary elements in $H^2(X_\pi,\C)$, 
we have 
$(\sum a_i[E_i])\cdot(\sum b_j[E_j])\=\sum a_i\overline{b}_j\,E_i\cdot E_j$.

The map $\varpi$ also induces a push-forward map 
$\varpi_*:H^2(X_{\pi'},\C)\to H^2(X_\pi,\C)$.\footnote{As $\varpi$ 
  is proper, $\varpi_* T$ is even defined for any current 
  $T$ on $X_{\pi'}$.} 
It is a basic fact that
$\varpi^*\omega\cdot\omega'=\omega\cdot\varpi_*\omega'$ 
for $\omega\in H^2(X_\pi,\C)$, $\omega'\in H^2(X_{\pi'},\C)$. 
Moreover $\varpi_*\varpi^* =\id$ as $\varpi$ is birational. 
From these two facts one immediately deduces
\begin{proposition}\label{P100}
  Suppose $\varpi:X_{\pi'}\to X_\pi$ is the blow-up at 
  \emph{one point} $p\in\pi^{-1}(0)$, and let 
  $E\=\varpi^{-1}(p)$. 
  Then $\varpi^*:H^2(X_{\pi},\C)\to H^2(X_{\pi'},\C)$ 
  is an isometric embedding and
  we have the orthogonal direct sum decomposition
  \begin{equation*}
    H^2(X_{\pi'},\C)
    =\varpi^*H^2(X_{\pi},\C)
    \perp\C[E].
  \end{equation*}
  Moreover, $[E]\cdot[E]=-1$.
\end{proposition}
\begin{corollary}{~}\label{C100}
  \begin{itemize}
  \item
    The cup product is a negative definite hermitian form on 
    $H^2(X_\pi,\C)$.
  \item
    For any modification $\varpi:X_{\pi'}\to X_\pi$,
    the map $\varpi^*:H^2(X_{\pi},\C)\to H^2(X_{\pi'},\C)$ 
    is an isometric embedding.
  \end{itemize}
\end{corollary}
\begin{remark}\label{R301}
  It follows from Proposition~\ref{P98}, and the definition of
  $H^2(\fX,\C)$ as an injective limit, that
  $H^2(\fX,\C)$ is in fact generated by classes $[E]$ of
  exceptional divisors with self-intersection $-1$. 
  More precisely, it is generated by classes $\omega$
  of the following form: $\omega=\imath_\pi[E]$, where
  $E\in\Gast$ and $\pi\in\fB$ is minimal such that
  $E\in\Gast_\pi$.
\end{remark}
\begin{remark}
  Proposition~\ref{P100} is also true when replacing 
  $\C$ by the ring of integers $\Z$. 
  By decomposing $\pi\in\fB$ into a sequence of point
  blow-ups, we infer the existence
  of a basis $F_1,\dots,F_n$ of $H^2(X_{\pi},\Z)$ (as a $\Z$-module)
  such that $F_i\cdot F_j=-\delta_{ij}$.
\end{remark}

Now we can define the intersection product on 
$H^2(\fX,\C)$. 
Pick two elements $\omega,\omega'\in H^2(\fX,\C)$. 
By Proposition~\ref{P98}, we can find $\pi\in\fB$ and
$\omega_\pi,\omega'_\pi\in H^2(X_\pi,\C)$ such that 
$\imath_\pi\omega_\pi=\omega$ and $\imath_\pi\omega'_\pi=\omega'$.
We set
\index{inner product!on $H^2(\fX)$}
\begin{equation*}
  \omega\cdot\omega'\=\omega_\pi\cdot\omega'_\pi.
\end{equation*}
By Corollary~\ref{C100}, 
this number does not depend on the choice of $\pi\in\fB$. 
\begin{corollary}{~}\label{C101}
  The intersection product is a negative definite hermitian 
  form on $H^2(\fX,\C)$. 
\end{corollary}
%
%
\subsection{Associated complex tree potentials}\label{sec-ve-ident}
Let us associate a function $g_\omega:\cVqm\to\C$ to
each cohomology class $\omega\in H^2(\fX,\C)$.

Fix $\pi\in\fB$, and irreducible components $E,F\subset\pi^{-1}\{0\}$
(\ie $E,F\in\Gast_\pi$). 
We define a function $g_{[E]}:\Phi(\Gast_\pi)\to\Q$ by
\begin{equation*}
  g_{[E]}(\nu_F)=
  \begin{cases}
    b(\nu_E)^{-1} &\text{if $F=E$}\\
    0 &\text{otherwise}.
  \end{cases}
\end{equation*}
Here $b(\nu_E)$ is the generic multiplicity
of the divisorial valuation $\nu_E$ and
$\Phi:\Gast\to\cVdiv$ denotes the isomorphism
between the universal dual graph and the valuative tree
as in Theorem~\ref{thm-universal}.

We then define $g_{\omega_\pi}:\Phi(\Gast_\pi)\to\C$ 
for any $\omega_\pi\in H^2(X_\pi,\C)$ by linearity. 
If $\pi'\in\fB$ and $\pi'=\pi\circ\varpi$ 
for some modification $\varpi$, then the function
$g_{\varpi^*\omega_\pi}$ is defined on
$\Phi(\Gast_{\pi'})\supset\Phi(\Gast_\pi)$ and
restricts to $g_{\omega_\pi}$ on $\Phi(\Gast_\pi)$.

Now fix a cohomology class $\omega\in H^2(\fX,\C)$, and a divisorial
valuation $\nu$. Pick $\pi\in\fB$ such that 
$\imath_\pi\omega_\pi=\omega$ 
with $\omega_\pi\in H^2(X_\pi,\C)$ and $\nu=\nu_E$ 
for some $E\in\Gast_\pi$. We set
\begin{equation*}
  g_{\omega}(\nu_E)\=g_{\omega_\pi}(\nu_E).
\end{equation*}
By what precedes, this does not depend on the choice of $\pi$.
\begin{theorem}\label{T345}
  The function $g_\omega:\cVdiv\to\C$ extends (uniquely)
  to a complex tree potential $g_\omega:\cVqm\to\C$
  whose associated Laplacian $\rho_\omega=\Delta g_\omega$
  is a complex atomic measure on $\cV$ supported on divisorial
  valuations. 
\end{theorem}
\begin{remark}
  As before, we are using the parameterization of $\cV$
  by skewness when talking about complex tree potentials.
\end{remark}
\begin{remark}
  We shall later show that any  
  complex atomic measure $\rho$ supported
  on divisorial valuations is of the form
  $\rho=\rho_\omega$ for a unique $\omega\in H^2(\fX,\C)$.
\end{remark}
In view of Remark~\ref{R301}, Theorem~\ref{T345} follows by 
linearity from the following more precise result:
\begin{proposition}\label{P307}
  Consider a cohomology class $\omega\in H^2(\fX,\C)$
  of the form 
  $\omega=\imath_\pi[E]$, where $\pi\in\fB$,
  $E\in\Gast_\pi$.
  Assume that $E$ has self-intersection -1.
  Then $g_\omega:\cVdiv\to\C$ extends to a complex 
  tree potential whose Laplacian $\rho=\Delta g_\omega$ 
  is an atomic measure given as follows:
  \begin{itemize}
  \item[(i)]
    if $E=E_0$ is the exceptional divisor obtained by blowing
    up the origin once, then $b(E)=1$ and
    \begin{equation*}
      \rho=\nu_\fm;
    \end{equation*}
  \item[(ii)]
    if $E$ intersects a unique $E'\in\Gast_\pi$, \ie if
    $E$ is obtained by blowing up a free point on $E'$, then 
    $b(E)=b(E')$ and 
    \begin{equation*}
      \rho=b(\nu_E)\nu_E-b(\nu_{E'})\nu_{E'};
    \end{equation*}
  \item[iii)]
    if $E$ intersects $E',E''\in\Gast_\pi$, \ie if
    $E$ is obtained by blowing up the intersection
    point $E'\cap E''$, then $b(E)=b(E')+b(E'')$ and
    \begin{equation*}
      \rho=b(\nu_E)\nu_E-b(\nu_{E'})\nu_{E'}-b(\nu_{E''})\nu_{E''}.
    \end{equation*}
  \end{itemize}
\end{proposition}
In order to prove this proposition, 
we shall use the following two results, whose proofs are 
postponed until the end of the section.
\begin{lemma}\label{L846}
  Fix $\pi\in\fB$ and a point $p\in\pi^{-1}(0)$.
  Let $\varpi$ be the blow-up of $p$ with exceptional divisor
  $F=\varpi^{-1}(p)$.
  Consider $\omega\in H^2(\fX,\C)$ and assume that
  $\omega=\imath_\pi\omega_\pi$ for some 
  $\omega_\pi\in H^2(X_\pi,\C)$.
  Then the following hold:
  \begin{itemize}
  \item[(a)]
    if $p$ is a free point, 
    \ie $p$ belongs to a unique $F'\in\Gast_\pi$, then 
    \begin{equation}\label{e786}
      g_\omega(\nu_F)=g_\omega(\nu_{F'});
    \end{equation}
  \item[(b)]
    if $p$ is a satellite point, 
    \ie $p=F'\cap F''$ for $F',F''\in\Gast_\pi$, then
    \begin{equation}\label{e796}
      g_\omega(\nu_F)
      =\frac{b_{F'}}{b_{F'}+b_{F''}}g_\omega(\nu_{F'})
      +\frac{b_{F''}}{b_{F'}+b_{F''}}g_\omega(\nu_{F''}).
    \end{equation}
  \end{itemize}
\end{lemma}
\begin{lemma}\label{L305}
  Consider $\omega\in H^2(\fX,\C)$ and pick
  $\pi\in\fB$ such that 
  $\omega=\imath_\pi\omega_\pi$ 
  for some $\omega_\pi\in H^2(X_\pi,\C)$.
  Assume that $\Gast_\pi$ has more than one element 
  and pick two adjacent elements $F',F''\in\Gast_\pi$ 
  (\ie $F'$ intersects $F''$). 
  Then $g_\omega$ is an affine function of skewness on the segment
  $[\nu_{F'},\nu_{F''}]$ in $\cVqm$.
\end{lemma}
\begin{proof}[Proof of Proposition~\ref{P307}]
  We only need to prove the formulas for $\rho$ as 
  the expressions for the generic multiplicities 
  $b(\nu_E)$ are known from Chapter~\ref{A3}.

  In~(i) we need to show that $g_\omega(\nu)=1$ for every
  divisorial valuation $\nu=\nu_F$. 
  This is clear for $F=E_0$. 
  We now proceed by induction on the ``length'' of $\nu$, 
  \ie the number of elements in $\Gast_\pi$, where $\pi\in\fB$
  is minimal such that $F\in\Gast_\pi$. If this length is
  one, then $F=E_0$ and we are done. Otherwise we may apply
  Lemma~\ref{L846}. The inductive assumption gives
  $g_\omega(\nu_{F'})=1$ in case~(a) and
  $g_\omega(\nu_{F'})=g_\omega(\nu_{F''})=1$ in case~(b).
  In both cases we get $g_\omega(\nu_F)=1$, completing
  the proof of~(i).

  The proof in cases~(ii) and~(iii) is similar to that
  in case~(i).
  By Lemma~\ref{L305}, $g_\omega$ is an affine function
  of skewness on all segments $[\nu_{F'},\nu_{F''}]$,
  where $F'$ and $F''$ are adjacent vertices in $\Gast_\pi$.
  Moreover, $g_\omega(\nu_E)=1$ and $g_\omega(\nu_F)=0$
  for every $F\in\Gast_\pi$, $F\ne E$.
  This determines the value of $g_\omega$ on the smallest 
  subtree $\cS\subset\cVqm$ containing all valuations 
  $\nu_F$, $F\in\Gast_\pi$. 
  Proving the formulas for $\rho$ in~(ii) or~(iii)
  then amounts to showing that $g_\omega$ is
  locally constant outside $\cS$.

  Pick a divisorial valuation $\nu\not\in\cS$ 
  and define $\nu_0=\max\{\mu\in\cS\ ;\ \mu\le\nu\}$. 
  We need to prove that $g_\omega(\nu)=g_\omega(\nu_0)$.  
  The valuations $\nu_0$ and $\nu$ are both divisorial, 
  say $\nu_0=\nu_{F_0}$ and $\nu=\nu_F$ 
  for some $F_0,F\in\Gast$. Moreover, $F$ is the last 
  exceptional divisor obtained by blowing up a sequence
  of infinitely nearby points $p_1,\dots,p_n$, starting with
  a point $p_1\in F_0$.
  Let $F_i$ be the exceptional divisor obtained by blowing up
  $p_i$, and write $\nu_i=\nu_{F_i}$.
  The key point is now that $p_1$ is a \emph{free} point on $F_0$.
  Thus $g_\omega(\nu_1)=g_\omega(\nu_0)$ by Lemma~\ref{L846}.
  Inductively, we obtain from~(a) or~(b) in the same lemma
  that $g_\omega(\nu_i)=g_\omega(\nu_0)$ for
  $1\le i\le n$. But $\nu_n=\nu$ so we are done.
\end{proof}
\begin{proof}[Proof of Lemma~\ref{L846}]
  Both proofs are analogous; we only treat~(b).
  In the basis of $H^2(X_\pi,\C)$ consisting of classes 
  of irreducible components of $\pi^{-1}(0)$, 
  write $\omega_\pi=c'[F']+c''[F'']+\dots$
  with $c',c''\in\C$. 
  By definition, $g_\omega(\nu_{F'})=c'/b_{F'}$
  and $g_\omega(\nu_{F''})=c''/b_{F''}$.
  In $H^2(X_{\varpi\circ\pi},\C)$, we have
  $\varpi^*\omega_\pi=(c'+c'')[F]+\dots$. 
  Together with $b_F=b_{F'}+b_{F''}$, this gives~\eqref{e796}.
\end{proof}
\begin{proof}[Proof of Lemma~\ref{L305}]
  Assume that $\nu_{F'}<\nu_{F''}$.
  The divisorial valuations in the segment 
  $]\nu_{F'},\nu_{F''}[$
  are of the form $\nu_F$, where $F\in\Gast$ 
  is obtained by a finite sequence
  of blowups at satellite points, starting with the point
  $F'\cap F''$. By induction it therefore suffices to show that
  $g_\omega$ is an affine function on the
  totally ordered set $\{\nu_{F'},\nu_F,\nu_{F''}\}$,
  where $F$ is the exceptional divisor obtained by blowing
  up the point $F'\cap F''$ once.
  
  If $(a',b')$ and $(a'',b'')$ are the Farey weights
  of $F'$ and $F''$, respectively, then the Farey weight of
  $F$ is $(a'+a'',b'+b'')$.
  We thus obtain from Lemma~\ref{L846}:
  \begin{multline*}
    \frac{g_\omega(\nu_F)-g_\omega(\nu_{F'})}
    {\a(\nu_F)-\a(\nu_{F'})}
    =\frac{g_\omega(\nu_F)-g_\omega(\nu_{F'})}
    {\frac{A(\nu_F)-A(\nu_{F'})}{b'}}
    =\frac{
      \frac{b'g_\omega(\nu_{F'})+b''g_\omega(\nu_{F''})}{b'+b''}
      -g_\omega(\nu_{F'})}
    {\frac1{b'}\left(\frac{a'+a''}{b'+b''}-\frac{a'}{b'}\right)}\\
    =\frac{g_\omega(\nu_{F''})-g_\omega(\nu_{F'})}
    {\frac1{b'}\left(\frac{a''}{b''}-\frac{a'}{b'}\right)}
    =\frac{g_\omega(\nu_{F''})-g_\omega(\nu_{F'})}
    {\frac{A(\nu_{F''})-A(\nu_{F'})}{b'}}
    =\frac{g_\omega(\nu_{F''})-g_\omega(\nu_{F'})}
    {\a(\nu_{F''})-\a(\nu_{F'})}.
  \end{multline*}
  Here we have used that the multiplicity is constant, equal
  to $b'$ on $]\nu_{F'},\nu_{F''}[$.
  This completes the proof.
\end{proof}
%
%
\subsection{Isometric embedding}\label{sec-isomembedd}
As we saw above, we have a natural intersection product
on the cohomology $H^2(\fX,\C)$.
In Section~\ref{S416} we showed that there is a subspace 
$\cM_0$ of the spaces of complex Borel measures on which
we have a well-defined inner product.

If $\omega\in H^2(\fX,\C)$, then the measure 
$\rho_\omega$ is atomic and supported on $\cVdiv$.
This implies that $\rho_\omega\in\cM_0$: 
see Remark~\ref{R428}.
\begin{theorem}\label{T347}
  The map $\omega\mapsto\rho_\omega$ gives an isometric
  embedding of $H^2(\fX,\C)$ into $\cM_0$ in the following sense:
  for any classes $\omega,\omega'\in H^2(\fX,\C)$, we have
  \begin{equation}\label{E348}
    -\omega\cdot\omega'=\rho_\omega\cdot\rho_{\omega'}.
  \end{equation}
  Here $\omega\cdot\omega'$ denotes the intersection product
  on $H^2(\fX,\C)$ as in Section~\ref{sec-ve-int} and 
  $\rho_\omega\cdot\rho_{\omega'}$ the inner product 
  on $\cM_0$ as in Section~\ref{S416}.
\end{theorem}
\begin{proof}
  Instead of proving~\eqref{E348} we shall prove the equivalent
  formula
  \begin{equation}\label{e304}
    -\omega\cdot\omega'
    =g_\omega(\nu_\fm)\,\overline{g_{\omega'}(\nu_\fm)}
    +\int_{\cVqm}
    \frac{dg_\omega}{d\a}\,\overline{\frac{dg_{\omega'}}{d\a}}\,d\l;
  \end{equation}
  see Theorem~\ref{T405}.
  Fix $\pi\in\fB$, and $E,E'\in\Gast_\pi$.
  Let us first prove~\eqref{e304} for $\omega=\imath_\pi[E]$ and
  $\omega'=\imath_\pi[E']$ under suitable assumptions.
  
  \smallskip
  (1) First suppose $E$ and $E'$ do not intersect, \ie
  $E$ and $E'$ are not adjacent elements in $\Gast_\pi$.
  By Proposition~\ref{P307} above, 
  the function $g_\omega$ is supported on the union of segments
  $[\nu_E,\nu_F[$ for all $F\in\Gast_\pi$ with $E\cap F\ne\emptyset$.
  The analogous assertion holds for $g_{\omega'}$. 
  Hence the product $g_\omega\overline{g_{\omega'}}$ 
  is identically zero everywhere on $\cVqm$. 
  This proves~\eqref{e304} in this case.
  
  \smallskip
  (2) Now suppose $E\cap E'\not=\emptyset$ but $E\ne E'$, \ie
  $E$ and $E'$ are adjacent elements in $\Gast_\pi$.
  Then $-\omega\cdot\omega'=-[E]\cdot[E']=-1$. 
  By the same argument as before, 
  the product $g_\omega\,\overline{g_{\omega'}}$ is identically 
  zero except on the segment 
  $I=\,]\nu_E,\nu_{E'}[$. 
  Assume $\nu_E>\nu_{E'}$
  and let $(a_E,b_E)$ and $(a_{E'}, b_{E'})$ be the
  Farey weights of $\nu_E$ and $\nu_{E'}$, respectively.
  
  The function $g_\omega$ is affine on $I$ taking the value
  $b_E^{-1}$ at $\nu_E$ and $0$ at $\nu_{E'}$. 
  The multiplicity on $I$ is constant and is given by 
  $m=|a_Eb_{E'}-a_{E'} b_E|$ by Remark~\ref{rem-compumul}.
  Whence the (left) derivative with respect to
  skewness of $g_\omega$ on $I$ is equal to
  \begin{equation*}
    \frac{dg_\omega}{d\a}
    =\frac{g_\omega(\nu_E)-g_{\omega}(\nu_{E'})}{\a_E-\a_{E'}}
    =\frac{m}{b_E(A_E-A_{E'})}
    =\frac{a_Eb_{E'}-a_{E'}b_E }{b_E(a_E/b_E-a_{E'}/b_{E'})}
    =b_{E'}.
  \end{equation*}
  A similar computation shows that
  $\frac{dg_{\omega'}}{d\a}=-b_{E'}$.
  The $\l$-length of $I$ is given by $1/(b_Eb_{E'})$ 
  so again~\eqref{e304} holds. 
  
  \smallskip 
  (3) If $\pi$ is the blow-up of the origin and $E=E'=E_0$ 
  is the exceptional divisor of $\pi$, then 
  $-\omega\cdot\omega'=+1$, $g_\omega=g_{\omega'}$ 
  is constant equal to $1$ on $\cVqm$, and~\eqref{e304} is
  immediate.
  
  \smallskip
  (4) Suppose $\pi$ contains more than one blow-up, 
  that $E=E'\in\Gast_\pi$ has self-intersection $-1$,
  and that $E$ is obtained by blowing up a free point on a 
  (unique) exceptional component $F$ with Farey weight $(a,b)$. 
  Then $\nu_E>\nu_F$, the Farey weight of $E$ is
  $(a+1,b)$, the multiplicity on $I=\,]\nu_F,\nu_E[$ is 
  constant equal to $b$, and $\l(I)=b^{-2}$.
  The function $g_\omega=g_{\omega'}$ is affine on $I$, 
  sends $\nu_E$ to $b^{-1}$ and $\nu_F$ to $0$, 
  is locally constant outside $I$, and vanishes at $\nu_\fm$. 
  Whence
  \begin{equation*}
    g_\omega(\nu_\fm)\,\overline{g_{\omega}(\nu_\fm)}
    +\int_{\cVqm}
    \frac{dg_\omega}{d\a}\,\overline{\frac{dg_{\omega}}{d\a}}\,d\l 
    =0+\left(\frac{1/b}{1/b^2}\right)^2\,\frac1{b^2}
    =+1.
  \end{equation*}
  This proves~\eqref{e304} in this case.

  \smallskip
  (5) Now suppose $E=E'$ has self-intersection $-1$
  and is obtained by blowing up a satellite point
  lying at the intersection of two divisors $F,F'\in\Gast_\pi$ 
  with Farey weights $(a,b)$ and $(a',b')$, respectively. 
  We may suppose $\nu_F>\nu_{F'}$.  
  The valuation $\nu_E$ has Farey weight
  $(a+a', b+b')$ by definition, and lies in the segment 
  $]\nu_{F'},\nu_F[$ by Lemma~\ref{Lkeylemma1}, 
  where the multiplicity is constant equal to 
  $m=|ab'-ba'|$ (see Remark~\ref{rem-compumul}). 
  The segment $I=\,]\nu_F,\nu_E[$ hence has 
  $\l$-length $\a(\nu_F)-\a(\nu_E)=1/ b (b+b')$
  whereas the segment $I'=\,]\nu_E,\nu_{F'}[$
  has length $1/b'(b+b')$.
  
  The function $g_\omega=g_{\omega'}$ 
  is affine on $I$ and $I'$, locally constant
  outside these two segments, and 
  $g_{\omega}(\nu_\fm)=g_{\omega}(\nu_F)=g_{\omega}(\nu_{F'})=0$, 
  $g_{\omega}(\nu_E)=1/(b+b')$. 
  A direct computation shows
  \begin{equation*}
    g_\omega(\nu_\fm)\,\overline{g_{\omega}(\nu_\fm)}
    +\int_{\cVqm}
    \frac{dg_\omega}{d\a}\,\overline{\frac{dg_{\omega'}}{d\a}}\,d\l 
    =\frac{b^2}{b(b+b')}+\frac{b^{'2}}{b'(b+b')} 
    =+1.
  \end{equation*}
  This gives~\eqref{e304}.
  
  \smallskip
  We can now prove Theorem~\ref{T347} in full generality. 
  Pick $\pi\in\fB$ minimal such that 
  $\omega=\imath_\pi(\omega_\pi)$ and
  $\omega'=\imath_\pi(\omega'_\pi)$ for classes
  $\omega_\pi,\omega'_\pi\in H^2(X_\pi,\C)$. 
  We proceed by induction on the number $N(\pi)$ 
  of elements in $\Gast_\pi$.
  When $N(\pi)=1$,~\eqref{e304} follows from~(3).
  
  For the inductive step, consider $\pi\in\fB$ with 
  $N(\pi)>1$. 
  Then we may write $\pi=\pi'\circ\varpi$, 
  where $N(\pi')=N(\pi)-1$ and $\varpi$ is 
  the blowup at a point on $(\pi')^{-1}(0)$.
  Let $E$ denote the exceptional divisor of $\varpi$. 
  By linearity, and by Corollary~\ref{C100}, 
  it suffices to check~\eqref{e304} for
  $\omega_\pi=[E]$ and 
  $\omega'_\pi\in\varpi^*H^2(X_{\pi'},\C)$; 
  and for $\omega_\pi=\omega'_\pi=[E]$. 
  The last case is taken care of by~(4) and~(5)
  and the first one reduces to~(1) or~(2),
  again using linearity of both sides of~\eqref{e304}.
\end{proof}
%
%
\subsection{Cohomology groups}\label{sec-ve-cogp}
We have shown that the assignment 
$\omega\mapsto\rho_\omega$ gives
an isometric embedding of $H^2(\fX,\C)$ into $\cM_0$.
We now wish to describe the image of $H^2(\fX,\C)$ as
well as of some of its subsets.
\begin{definition}
  We let $H^2_+(\fX)\subset H^2(\fX,\R)$ be the set of 
  cohomology classes $\omega$ 
  which can be written $\omega=\imath_\pi\omega_\pi$, with
  $\omega_\pi=\sum_{E\in\Gast_\pi}a(E)[E]$ where $a(E)\ge0$.  
  Such classes are called
  \emph{pseudo-effective}.
  \index{cohomology class!pseudo-effective}
\end{definition}  
\begin{definition}
  We let $H^2_\mathrm{an}(\fX)\subset H^2(\fX,\R)$ be the dual cone of
  $H^2_+(\fX)$, that is, the set of classes $\omega\in H^2(\fX,\R)$
  for which $\omega\cdot\omega'\le 0$ for all $\omega'\in
  H^2_+(\fX)$. Such classes are called \emph{anti-numerically
  effective}\index{cohomology class!antinef} (or antinef for short).
\end{definition}
\begin{theorem}\label{T286}
  The map $\omega\to\rho_\omega$ induces an isometry between
  \begin{itemize}
  \item[(i)]
    $H^2(\fX,\C)$ and the set of complex atomic measures 
    supported on divisorial valuations;
  \item[(ii)]
    $H^2(\fX,\R)$ and the set of real atomic measures 
    supported on divisorial valuations; 
  \item[(iii)]
    $H^2(\fX,\Z)$ and the set of real atomic measures $\rho$
    supported on divisorial valuations such that
    $\rho\{\nu\}\in b(\nu)\Z$ for every $\nu\in\cVdiv$;
  \item[(iv)]
    $H^2_+(\fX)$ and the set of real atomic measures 
    $\rho$ supported on divisorial valuations
    such that $\rho\cdot\rho'\ge0$ for all positive measures $\rho'$;
  \item[(v)]
    $H^2_\mathrm{an}(\fX)$ and the set of positive atomic 
    measures supported on divisorial valuations.
 \end{itemize}
\end{theorem}
\begin{remark}
  Thus the image of $H^2(\fX;\Z)\cap H^2_\mathrm{an}(\fX)$ is exactly
  the set of tree measures $\rho_I$ for primary ideals $I\subset R$.
  See Theorem~\ref{Tcharideal} and Proposition~\ref{P0-ideal}.  This
  fact is reminiscent of Lefschetz' theorem of realization of
  cohomology classes $H^2(X, \Z) \cap H^{1,1}(X)$ by divisors on a
  projective variety. On the vo\^ute \'etoil\'ee analytic ideals 
  play the role of (effective) divisors.
\end{remark}
\begin{proof}[Proof of Theorem~\ref{T286}]
  We have already seen in Theorem~\ref{T347} that
  $\omega\mapsto\rho_\omega$ is an isometry and in particular
  injective.
  Let us first prove~(iii).
  Then~(i) and~(ii) follow by linearity, 
  whereas~(iv)-(v) will be proved below.

  Denote by $\cM_\Z$ the set of real 
  atomic measures supported on divisorial valuations whose mass 
  on a divisorial valuation $\nu$ is always an
  integer multiple of the generic multiplicity of $\nu$. 
  Then $\cM_\Z$ is an abelian subgroup of $\cM_0$.
  We have to show that $\cM_\Z$ is the image of $H^2(\fX;\Z)$.

  First pick $\omega\in H^2(\fX,\Z)$. 
  Let us show that $\rho_\omega\in\cM_\Z$.
  By linearity it suffices to do this in the case
  $\omega=\imath_\pi[E]$, where 
  $E\in\Gast_\pi$ 
  and $\pi\in\fB$ is minimal such that $E\in\Gast_\pi$:
  see Remark~\ref{R301}.
  But then the conclusion
  is immediate from Proposition~\ref{P307}.

  Conversely, let us show that any measure
  $\rho\in\cM_\Z$ can be obtained as $\rho_\omega$ for some
  $\omega\in H^2(\fX,\Z)$. 
  Again by linearity, it is sufficient to prove this for 
  $\rho=b(\nu)\nu$, where $\nu=\nu_E$ is divisorial.
  Pick $\pi\in\fB$ minimal such that $E\in\Gast_\pi$.
  We proceed by induction on the number of elements in $\Gast_\pi$,
  \ie the minimal number of blow-ups necessary to create the 
  exceptional divisor $E$.
  
  If $E=E_0$ is the exceptional divisor obtained by blowing up 
  the origin once, then $\nu_E=\nu_\fm$, $b(\nu_E)=1$ and
  we conclude by Proposition~\ref{P307}~(i).
  
  Now assume $\Gast_\pi$ has more than one element. 
  Then either $E$ intersects a unique $E'\in\Gast_\pi$
  or two distinct $E',E''\in\Gast_\pi$.
  We will consider only the second of these cases,
  the first one being almost identical.
  We may apply the inductive hypothesis to
  $\nu_{E'}$ and $\nu_{E''}$ and find
  $\omega',\omega''\in H^2(\fX,\Z)$ such that
  $\rho_{\omega'}=b(\nu_{E'})\nu_{E'}$ and
  $\rho_{\omega''}=b(\nu_{E''})\nu_{E''}$.
  Set $\omega=\omega'+\omega''+b(E)\imath_\pi[E]$.
  We get from Proposition~\ref{P307}~(iii) that 
  \begin{equation*}
    \rho_\omega
    =\rho_{\omega'}+\rho_{\omega''}+
    b(\nu_E)\nu_E-b(\nu_{E'})\nu_{E'}-b(\nu_{E''})\nu_{E''}
    =b(\nu_E)\nu_E. 
  \end{equation*}
  This completes the induction step
  and hence the proof of~(iii).

  As noted above~(i) and~(ii) follow by linearity.
  In order to prove~(iv) 
  it suffices to show that $\omega\in H^2_+(\fX)$ iff
  $g_\omega\ge0$ on $\cVdiv$. Indeed,~\eqref{e460} shows that
  $\rho_\omega\cdot\rho'=\int_\cV g_\omega\,d\rho'$ for any
  positive measure $\rho'$. 
  
  That $g_\omega\ge0$ whenever $\omega\in H^2_+(\fX)$ is a
  direct consequence of the definition of $g_\omega$.
  For the converse, consider $\omega\not\in H^2_+(\fX)$. 
  Pick $\pi\in\fB$ and $\omega_\pi\in H^2(X_\pi,\R)$ such that 
  $\omega=\imath_\pi\omega_\pi$. 
  Write $\omega_\pi=\sum a_i[E_i]$ 
  with $E_i\in\Gast_\pi$ and $a_i\in\R$. 
  As $\omega\not\in H^2_+(\fX)$, one of these real numbers is
  negative, say $a_1<0$. 
  Then $g_\omega(\nu_{E_1})=a_1/b(\nu_1)<0$. 

  This completes the proof of~(iv).
  Finally,~(v) follows by duality from~(iv) 
  and from Theorem~\ref{T347}.
\end{proof}

%
%
%
%
%
%
\chapter*{Appendix.}

\renewcommand{\thechapter}{\Alph{section}}
\renewcommand{\thesection}{\Alph{section}}
\setcounter{chapter}{-1}
\setcounter{section}{0}
\setcounter{subsection}{0}
\setcounter{theorem}{0}
\setcounter{equation}{0}
%
%
%
%
We end this monograph with an appendix containing complements to
results already proved in the main body of the text.

Appendix~\ref{sec-inf-sing} is devoted to infinitely singular 
valuations. Specifically it contains a list of properties that each 
characterizes a valuation as being infinitely singular. 
We also give a few constructions of infinitely singular valuations.

In Appendix~\ref{sec-tangent} we give different
characterizations of the tree tangent space at a divisorial
valuation. 

It is a fascinating fact that there are many paths to the valuative 
tree. We summarize in Appendix~\ref{sec-clas} the
classification of Krull valuations on $R$ from 
the different points of view emphasized in the memoir.

In order to help the reader familiar with invariants and terminology
used to describe plane curve singularities (as in Zariski~\cite{Z3}), 
we give a short dictionary between this terminology and ours.
This is done in Appendix~\ref{S501}, where we also explain how 
the Eggers tree of a reduced curve singularity can be 
naturally embedded in the valuative tree.

We conclude in Appendix~\ref{sec-notcomplete} by discussing the
importance of the various assumptions we made on the ring $R$.
Our standing assumption that $R$ be the ring of formal power
series in two complex variables is clearly unnecessarily
restrictive. As the discussion shows, our method applies,
for instance, to the ring of holomorphic germs at the origin
in $\C^2$ and to the local ring at a smooth point on a
surface over an algebraically closed field.
%
%
%
%
\setcounter{theorem}{0}
\setcounter{equation}{0}
\setcounter{table}{0}
\setcounter{figure}{0}
\section{Infinitely singular valuations}\label{sec-inf-sing}
The infinitely singular valuations in $\cV$ are numerous, but 
also the most complicated to describe.
For the convenience of the reader we gather in one place all
the characterizations of infinitely singular valuations that
we have seen so far. Most of what we present here will also
figure in the classification tables in Appendix~\ref{sec-clas}
but we feel it is useful to have the information spelled out
in more detail.
In addition, we present a couple of explicit constructions of 
infinitely singular valuations.

As before we work with the ring $R$ of formal power series in
two complex variables, although as explained in 
Appendix~\ref{sec-notcomplete} this is in fact unnecessarily
restrictive.
%
%
\subsection{Characterizations}
The following result characterizes infinitely singular valuations
from many different points of view. 
\begin{proposition}
For a valuation $\nu\in\cV$ the following properties are equivalent
and characterize $\nu$ as being infinitely singular:
\begin{itemize}
\item[(i)]
  the multiplicity $m(\nu)$ is infinite;
\item[(ii)] 
  the approximating sequence $(\nu_i)_0^g$ of $\nu$ is infinite, 
  \ie $g=\infty$;
\item[(iii)]
  the value semigroup $\nu(R)$ is not finitely generated;
\item[(iv)]
  the numerical invariants of $\nu$ are given as follows:
  $\ratrk\nu=\rk\nu=1$ and $\trdeg\nu=0$;
\item[(v)]
  for some (or, equivalently, any)
  choice of local coordinates $(x,y)$, the SKP
  $[(U_j);(\btilde_j)]$ defined by 
  $\nu$ has infinite length and $n_j\ge2$ infinitely often 
  (or, equivalently, $m(U_j)\to\infty$);
\item[(vi)]
  for any choice of local coordinates $(x,y)$, 
  some (or, equivalently, any) choice of extension $\hnu$
  to a valuation on the ring $\bk[[y]]$, is of the form
  $\hnu=\val[\hphi;\hbeta]$, with 
  $\hphi=\sum_1^\infty a_jx^{\hbeta_j}$, $\hbeta=\lim\hbeta_j$,
  and the $\hbeta_j$'s have unbounded denominators, \ie 
  $\hphi\not\in\hk$;
\item[(vii)]
  the sequence of infinitely nearby points associated to $\nu$
  is of Type~1 in Definition~\ref{D403}, \ie it contains 
  both infinitely many free and infinitely many satellite points.
\end{itemize}
\end{proposition}
\begin{proof}
  As already noted, all the above characterizations can be read off
  from results already proved in the monograph. 
  Specifically, our definition of an infinitely singularly valuation
  was exactly~(v), for a fixed choice of coordinates.
  Theorem~\ref{divis} asserted that this definition is equivalent
  to~(iv) and hence independent on the choice of local coordinates.
  The equivalences of~(i),~(ii) follow from Proposition~\ref{P402}
  and~(iii) from Proposition~\ref{P421}.
  Further~(vi) is a consequence of Theorem~\ref{T602} 
  and~(vii) of Corollary~\ref{C432} and Theorem~\ref{thm-universal}.
\end{proof}
\begin{remark}
  In addition we have seen that the infinitely singular valuations
  are ends in the valuative tree, but this is not a characterizing
  property as the curve valuations are also ends.
\end{remark}
%
%
\subsection{Constructions}
Next we outline a couple of constructions of infinitely singular 
valuations, starting with the 
construction of an infinitely singular valuation with prescribed
skewness or thinness. See also Lemma~\ref{L777}.
\begin{proposition}\label{P-inf401}
  Pick a quasimonomial valuation $\nu$, and any real number $\a>\a(\nu)$
  (resp.\ $A>A(\nu)$). Then there exists an infinitely
  singular valuation $\mu>\nu$ with $\a(\mu)=\a$ 
  (resp.\ with $A(\mu)=A$).
\end{proposition}
\begin{proof}
  We may assume that $\nu$ is divisorial. Indeed, otherwise replace
  $\nu$ by a divisorial valuation $\nu'>\nu$ such that 
  $\a(\nu')<\a$ (resp.\ $A(\nu')<A$).

  Let us first construct an infinitely singular valuation with
  prescribed skewness.  Our strategy is to extend the approximating
  sequence $(\nu_i)_0^g$ of $\nu$ to an infinite approximating sequence
  $(\nu_i)_0^\infty$. By convention, $\nu_{g+1}=\nu$.  Pick a curve
  valuation $\mu_{g+2}>\nu$ with $m(\mu_{g+2})=b(\nu)$, and consider a
  divisorial valuation $\nu_{g+2}$ in the segment $]\nu,\mu_{g+2}[$.  It
  follows from~\eqref{e707} that the set of divisorial valuations in
  $]\nu,\mu_{g+2}[$ with generic multiplicity equal to $b(\nu)$ is a
  discrete set.  We may therefore pick $\nu_{g+2}$ in this segment such
  that $b(\nu_{g+2})>m(\nu_{g+2})=b(\nu)$ and
  $\a-1/2>\a(\nu_{g+2})>\a(\nu)$.
  
  Inductively, given $(\nu_i)_1^{g+k}$ we construct $\nu_{g+k+1}$
  divisorial with $\nu_{g+k+1}>\nu_{g+k}$,
  $b(\nu_{g+k+1})>m(\nu_{g+k+1})=b(\nu_{g+k})$ and 
  $\a-1/2^k>\a(\nu_{g+k+1})>\a(\nu_{g+k})$. 
  Then $(\nu_i)_0^\infty$ defines an
  approximating sequence for an infinitely singular valuation $\mu$,
  satisfying $\mu>\nu$ and $\a(\mu)=\a$.
  
  \smallskip
  For finding an infinitely singular valuation with prescribed thinness we
  may follow the argument in the proof of Lemma~\ref{L777}. 
  However, let us instead show how to use Puiseux series, exploiting 
  the analysis and notation of Chapter~\ref{sec-puis}.
  Fix local coordinates $(x,y)$ such that $\nu(y)\ge\nu(x)=1$
  and pick a valuation $\hnu\in\hcV$ whose image in $\cV$ under the 
  restriction map equals $\nu$.
  Then $\hnu=\val[\hphi;\hbeta]$ for a Puiseux series
  $\hphi=\sum_1^qa_jx^{\hbeta_j}$ with 
  $\hbeta_j<\hbeta_{j+1}<\hbeta$, 
  and all coefficients $a_j\ne0$. 
  By Theorem~\ref{T602}, $\hbeta=A(\nu)-1<A-1$. 
  Since $\nu$ is divisorial, $\hbeta$ is rational.
  Set $\hbeta_{q+1}=\hbeta$ and pick 
  an arbitrary sequence of strictly increasing rational numbers 
  $(\hbeta_j)_{q+2}^\infty$ such that $\hbeta_{q+2}>\hbeta$
  and $\hbeta_j\to A-1$.
  Define $\hpsi=\sum_1^\infty a_jx^{\hbeta_j}$, where, say,
  $a_j=1$ for $j>q$. 
  Set $\hmu=\val[\hpsi;A-1]$. 
  Then $\hmu$ is of special type and has Puiseux parameter $A-1$.
  Moreover, $\hmu>\hnu$. Let $\mu\in\cV$ be the image of $\hmu$
  under the restriction map. 
  Then Theorem~\ref{T602} implies that $\mu$ is infinitely
  singular, $\mu>\nu$ and $A(\mu)=A$.
\end{proof}
\begin{remark}\label{R405}
  By slightly modifying the first part of the proof above 
  we can construct an infinitely singular valuation with prescribed 
  skewness $t\in(1,\infty)$ and
  infinite thinness.  The reason is that we may choose the $\nu_i$'s
  inductively to have very high multiplicity. Thus $m_i$ grows very fast
  and this allows for skewness $\a_i$ to increase to $t\in(1,\infty)$
  whereas thinness $A_i$ tends to infinity. The details are left to the
  reader.
\end{remark}  
\begin{remark}\label{R404}
  Similarly one can construct a sequence $(\nu_n)_0^\infty$ of
  infinitely singular valuations with $\a(\nu_n)\to1$ and $A(\nu_n)\ge
  3$. Here is an outline of the construction.  Fix a large integer
  $n$. Pick a monomial divisorial valuation $\mu'$ with
  $\a(\mu')=1+1/n$. Then $b(\mu')=n$.  Then pick $\mu$ divisorial with
  $\mu>\mu'$, $\mu$ representing a generic tangent vector at $\nu'$
  and $\a(\mu)-\a(\mu') = 1/n$.  Then $\a(\mu) = 1+2/n$ but $A(\mu)=
  3+ 1/n >3$.  Now use Lemma~\ref{L777} to replace $\mu$ by an
  infinitely singular valuation $\nu$.
\end{remark}

%
%
%
%
\setcounter{theorem}{0}
\setcounter{equation}{0}
\setcounter{table}{0}
\setcounter{figure}{0}
\section{The tangent space at a divisorial valuation}\label{sec-tangent}
Recall from Section~\ref{tree-subsec} that in 
a general nonmetric tree $\cT$, the (tree) tangent space $T\tau$ at 
a point $\tau\in\cT$ is the set of equivalence classes 
$\cT\setminus\{\tau\}/\sim$, where $\sigma\sim\sigma'$ 
iff the two segments $]\tau,\sigma]$ and $]\tau,\sigma']$ intersect. 

We saw in Section~\ref{tree-struc} that branch points of the valuative
tree correspond to divisorial valuations. 
Our goal in this appendix is to describe the tangent space at
a divisorial valuation from several points of view.
\begin{theorem}\label{div-tangent}
  \index{tangent space! at a divisorial valuation}
  Let $\nu=\nu_E\in\cV$ be a divisorial valuation, associated to some
  exceptional component $E\subset\pi^{-1}(0)$, where $\pi$ is a
  composition of point blow-ups. Then there exist natural bijections
  between the following four sets:
  \begin{itemize}
  \item
    $\cT_1$: tree tangent vectors at $\nu$ (see Section~\ref{tree-def});
  \item
    $\cT_2$: Krull valuations $\mu$ satisfying $R_\mu\subsetneq
    R_\nu$ (see Sections~\ref{krullval} and~\ref{versus});
  \item
    $\cT_3$: sequences of infinitely nearby points $(p_j)_0^\infty$ 
    of Type~3 (see Section~\ref{S404}), such that the 
    divisorial valuation $\nu_n$ associated to the 
    truncated sequence $(p_j)_0^n$ converges (weakly) to $\nu$
    as $n\to\infty$;
  \item
    $\cT_4$: points on the exceptional component $E$.
  \end{itemize}
\end{theorem}
The Krull valuations in~$\cT_2$ are exactly the exceptional
curve valuations of the form $\nu_{E,p}$: see Lemma~\ref{L601}.

The proof of this result occupies the rest of 
Appendix~\ref{sec-tangent}. 
We shall construct injective maps 
$\Psi_1:\cT_1\to\cT_2$, 
$\Psi_2:\cT_2\to\cT_3$, 
$\Psi_3:\cT_3\to\cT_4$ and 
$\Psi_4:\cT_4\to\cT_1$.
and prove that the composition 
$\Psi_3\circ\Psi_2\circ\Psi_1\circ\Psi_4$
is the identity.

\smallskip
Fix a divisorial valuation $\nu\in\cV$ for the rest of
the proof.
Also fix a composition $\pi\in\fB$ of blowups such that
$\nu=\nu_E$ for some exceptional component 
$E\subset\pi^{-1}(0)$.

\smallskip
The construction of $\Psi_1:\cT_1\to\cT_2$ is interesting as it
allows us to interpret tangent vectors at $\nu$ in terms of 
directional derivatives.  
Let us first prove a preliminary result.
\begin{lemma}\label{deriv}
  Let $\nu\in\cV$ be divisorial, $\e>0$, and
  $[0,\e]\ni t\mapsto\nu_t\in\cV$ a (weakly) 
  continuous map with $\nu_0=\nu$ such that 
  $\frac{d}{dt}\a(\nu_t)$ is a nonzero constant on $]0,\e[$.
  Consider $\R\times\R$ with the lexicographic order.
  Then the function 
  \begin{equation*}
    \mu(\psi)=\left(
      \nu(\psi),\frac{d}{dt} 
      \bigg|_{t=0}\nu_t(\psi)
    \right)
  \end{equation*}
  defines a centered Krull valuation on $R$, with valuation
  ring satisfying $R_\mu\subsetneq R_\nu$.
\end{lemma}
\begin{proof}
  Clearly $\mu(\phi\psi)=\mu(\phi)+\mu(\psi)$ for $\phi,\psi\in R$.
  Further, $\nu_t(\phi+\psi)\ge\min\{\nu_t(\phi),\nu_t(\psi)\}$ for
  any $t$.  If strict inequality holds at $t=0$, then
  $\mu(\phi+\psi)\ge\min\{\mu(\phi),\mu(\psi)\}$ by the lexicographic
  ordering.  Otherwise $t\mapsto\nu_t(\phi+\psi)$, and
  $t\mapsto\min\{\nu_t(\phi),\nu_t(\psi)\}$ are affine functions near
  $t=0$, coinciding at $t=0$.  Then the slope of the latter cannot
  exceed the slope of the former, and this implies that
  $\mu(\phi+\psi)\ge\min\{\mu(\phi),\mu(\psi)\}$.  
  This shows that $\mu$ is a Krull valuation. 
  It is immediate that $R_\mu \subsetneq R_\nu$.
\end{proof} 
Consider a tree tangent vector $\vv\in\cT_1$ at $\nu\in\cV$. This is
represented by a segment that can be parameterized by skewness.
Lemma~\ref{deriv} gives us a Krull valuation $\mu \= \Psi_1(\vv)$ which
satisfies $R_{\mu} \subsetneq R_\nu$, \ie $\mu\in\cT_2$.
It does not depend on the choice of the segment, as two segments
defining the same tree tangent vector intersect in a one-sided
neighborhood of $\nu$, hence define the same Krull valuation.

To show that $\Psi_1$ is injective, take two different tangent vectors
$\vv_1 \ne \vv_2$, and choose $\phi\in\fm$ irreducible such that
$\phi$ represents $\vv_1$ but not $\vv_2$. 
Then it is straightforward to verify that
$\Psi_1(\vv_1)(\phi)\ne\Psi_1(\vv_2)(\phi)$.

\smallskip
The construction of $\Psi_2:\cT_2\to \cT_3$ is done as follows. 
For any Krull valuation $\mu$, we let $\Psi_2(\mu)$ be the sequence of
infinitely nearby points $(p_j)_0^{\infty}$ constructed in
Section~\ref{sec-equivalence}. 
Let us recall the construction.
The point $p_0$ is the origin in $\C^2$. 
We let $\tpi_0:X_0\to(\C^2,0)$ be the blow-up at $p_0$
with exceptional divisor $E_0$. 
Then $p_1\in E_0$ is defined to be the
center of the valuation $\mu$ in $X_0$.
Inductively we construct a sequence of points 
$p_{j+1}\in X_j$, with $\tpi_j:X_j\to X_{j-1}$
being the blow-up at $p_j$, and we write 
$E_j$ for the exceptional divisor of $\tpi_j$. 
The point $p_{j+1}$ is chosen to be the center of $\mu$ inside
$X_j$. In particular, $p_{j+1}\in E_j$.
Write $\pi_j=\tpi_0\circ\dots\circ\tpi_j$.

As $R_\mu\subsetneq R_\nu$, the center of $\mu$ is always included in
the center of $\nu$. As $\nu$ is divisorial, there exists
$j_0\ge1$ such that the center of $\nu$ in $X_j$
is a point for $j=0,\dots,j_0-1$, and an irreducible component
of $\pi_j^{-1}(0)$ for $i\ge j_0$. 
In $X_{j_0}$, it is given by $E_{j_0}$. 
In $X_j$ for $j>j_0$,
it is the strict transform of $E_{j_0}$ by 
$\tpi_j\circ\dots\circ\tpi_{j_0+1}$. 
We shall write $E_{j_0}$ for this strict transform, too.
This shows that for all $j>j_0$, $p_{j+1}$ has to be the
intersection point of $E_{j_0}$ and $E_j$. 
The sequence $(p_j)_0^\infty$ is thus of Type~3.
By Theorem~\ref{Ttype3}, the divisorial valuation $\nu_n$
associated with the truncated sequence $(p_j)_0^n$, converges
to $\nu_{j_0}=\nu$ as $n\to\infty$.
This shows that $\Psi_2(\mu)\in\cT_3$.

It is clear that $\Psi_2$ is injective by Theorem~\ref{blw}.

\smallskip
Now pick a sequence of infinitely nearby points 
$(p_j)_0^\infty\in\cT_3$ of Type~3. 
As before denote by $\tpi_j:X_j\to X_{j-1}$ the
blow-up at $p_j$ and $E_j$ the exceptional divisor of $\tpi_j$. 
The identity map $\id:(\C^2,0)\to (\C^2,0)$ lifts to a rational 
map $\id_j:X_j\to X$ for all $j\ge0$, where $\pi:X\to(\C^2,0)$ was
fixed at the beginning of the proof.
As $(p_j)_0^\infty$ is of Type~3, there
exists an index $j_0$ such that $p_{j+1}$ is the intersection point of
$E_j$ with the strict transform of $E_{j_0}$ (which we shall also
denote by $E_{j_0}$) by $\tpi_j\circ\dots\circ\pi_{j_0+1}$ 
for all $j>j_0$. 
The divisorial valuation associated to $E_j$
converges to $\nu_{j_0}$ (Theorem~\ref{Ttype3}).
As $\nu_{j_0}=\nu$, $\id_j$ sends $E_{j_0}$ bijectively onto $E$ for 
$j\ge j_0$. On the other hand, for $j$ large enough,
$\id_j$ is regular at $p_{j+1}$. 
The point $q\=\id_j(p_{j+1})$ equals $\id_j(E_{j_0}\cap E_{j+1})$ 
hence lies in $\id_j(E_{j_0})=E$. 
It is also clear that 
$\id_{j+1}(p_{j+2})=\id_{j+1}(\pi_{j+1}(p_{j+2}))=\id_j(p_{j+1})$, 
so that $q$ is independent of $j$.
By definition we set $\Psi_3((p_j)_0^\infty)\=q$.

To see that $\Psi_3$ is injective, choose two different sequences of
infinitely nearby points 
$(p_j)_0^\infty, (p'_j)_0^\infty\in\cT_3$.
As the divisorial valuation associated to both sequences is the same,
and equal to $\nu$, it follows that $p_j=p'_j$ until $j=j_0$, 
where $E_{j_0}$ is the exceptional component attached to $\nu$. 
Being both of Type~3, the sequences are determined uniquely by 
$p_{j_0+1}$ and $p'_{j_0+1}$, respectively. 
By assumption these two points are distinct.
Now pick $j>j_0$, and let $\varpi_j:Y_j\to(\C^2,0)$ 
be the composition of blow-ups at the points 
$p_0,\dots,p_{j_0}$ and then at both
$p_{j_0+1},\dots,p_j$ and $p'_{j_0+1},\dots,p'_j$. 
As before, the identity map 
$(\C^2,0)\to(\C^2,0)$ lifts as a map 
$\id_j:Y_j\to X$.
For $j$ large enough, $\id_j$ is
regular at $p_{j+1}$ and $p'_{j+1}$.
These two distinct points
belong to the strict transform of $E_{j_0}$, which we
again denote by $E_{j_0}$. 
As $\id_j$ is a bijection of $E_{j_0}$ onto $E$,
we conclude that $\Psi_3((p_j)_0^\infty)\ne\Psi_3((p'_j)_0^\infty)$, 
\ie that $\Psi_3$ is injective.

\smallskip
Finally, for a point $p\in E$,
we define $\Psi_4(p)$ to be the tree tangent vector at $\nu$ defined
by the divisorial valuation obtained by blowing-up $p$.
Lemma~\ref{Ptangentgamma} and Theorem~\ref{thm-universal}
imply that $\Psi_4$ is injective.

\smallskip
The proof of Theorem~\ref{div-tangent} will be complete, 
if we now show that the composition 
$\Psi_3\circ\Psi_2\circ\Psi_1\circ\Psi_4$ equals the identity.
Pick an element in $\cT_4$, that is, a point $p\in E$.  
By definition, $\Psi_4(p)$ is the tree tangent vector at $\nu$ 
represented by the divisorial valuation $\nu_p$ associated 
to the blow-up at $p$. 
The valuation $\mu\=\Psi_1\circ\Psi_4(p)$ is then
a Krull valuation whose valuation ring is included in
$R_\nu$, hence the center of $\mu$ in $X$ belongs to $E$. 
Pick an irreducible curve $C=\{\phi=0\}$, 
$\phi\in\fm$, whose strict transform $C'$ is smooth 
and intersects $E$ transversely at $p$.
Then $\nu_C$ and $\nu_p$ represent the same
tree tangent vector at $\nu$. 
The segment $[\nu,\nu_C]\cap[\nu,\nu_p]$ is thus nontrivial. 
We parameterize it by skewness: 
$t\mapsto\nu_{C,t+\a(\nu)}$. 
By definition, we have 
$\mu=(\nu,\frac{d}{dt}|_{t=0}\nu_{C,t+\a(\nu)})$. 
In particular, we infer that
$\mu(\phi)=(\nu(\phi),c)$ with $c$ \emph{positive}.

Let $\tilde{\mu}=(\div_E,\frac{d}{d\tau}|_{\tau=0}\nu_{C',\tau})$. 
It defines an exceptional curve valuation
$\tilde{\mu}$ which is centered at $p$. 
Its image by $\pi_*$ is an exceptional curve valuation 
whose first component equals $\pi_*\div_E$ which is 
equivalent to $\nu$. 
Thus $R_{\pi_*\tilde{\mu}}\subsetneq R_\nu$. 
It is also clear by definition that
$(\pi_*\tilde{\mu})(\phi)=(\div_E(\pi^*\phi),c')$ with $c'>0$.

Now pick a generic element $x\in\fm$ such that $\{x=0\}$ is a smooth
curve, $\nu_x\wedge\nu_C=\nu_\fm$, 
and the strict transform of $\{x=0\}$ by $\pi$ does not contain $p$. 
Then $\mu(x)=(\nu(x),0)$, and 
$(\pi_*\tilde{\mu})(x)=(\div_E(\pi^* x),0)$. 
Define $\psi\=\phi^k/x^l$ where 
$\nu(\phi)k=l\nu(x)$ with $k,l\in\N$. 
Then $\psi$ is an element of $\fm_\mu\cap\fm_{\pi_*\tilde{\mu}}$ 
with $\nu(\psi)=0$. 
Remark~\ref{R413} then implies $R_\mu=R_{\pi_*\tilde{\mu}}$. 
In particular, we conclude that the center of $\mu$ in $X$ 
is precisely the center of $\tilde{\mu}$ which is equal to $p$.

At this stage, we have shown that the Krull valuation 
$\Psi_1\circ\Psi_4(p)$ is an exceptional curve valuation 
centered at $p$, whose associated divisorial valuation is equal 
to $\nu$. 
Now let $(p_j)_0^\infty\=\Psi_2\circ\Psi_1\circ\Psi_4(p)$ 
be the sequence of infinitely nearby points associated to $\mu$. 
Consider the sequence of
blow-ups $\tpi_j:X_j\to X_{j-1}$ as above. 
By construction, $\tpi_j$ is the blow-up at $p_j$, 
which is the center of $\mu$ in $X_{j-1}$. 
For $j$ large enough, the identity map 
$(\C^2,0)\to(\C^2,0)$ lifts to a rational map 
$\id_j:X_j\to X$ which is regular at $p_{j+1}$. 
The image of the center of $\mu$ in $X_j$ by $\id_j$ is the
center of $\mu$ in $X$. 
Thus $\id_j(p_{j+1})=p$ for $j$ large enough. 
This proves that $\Psi_3\circ\Psi_2\circ\Psi_1\circ\Psi_4(p)=p$, 
and concludes the proof.
%
%
%
%
\setcounter{theorem}{0}
\setcounter{equation}{0}
\setcounter{table}{0}
\setcounter{figure}{0}
\section{Classification}\label{sec-clas}
As we have shown in this monograph, centered Krull valuations on $R$
can be interpreted in many different ways, each of which leads
to a (full or partial) classification.
Let us review the four approaches that we have considered.

\smallskip
The first (and classical) way to classify Krull valuations is through
their value groups $\nu(K)$, in particular through the
numerical invariants $\rk$, $\ratrk$ and $\trdeg$. In the two-dimensional
case that we are concerned with in this monograph, this is feasible
since Abhyankar's inequalities give strong restrictions on the values of
these invariants.

\smallskip
A second way, and the one emphasized in the monograph at hand, 
is to identify Krull valuations with points or tangent vectors in 
the valuative tree $\cV$. 
The non-metric tree structure on $\cV$ then leads to a classification
of Krull valuations into ends, regular points, branch points
and tangent vectors.
The valuative tree also comes with two natural parameterizations:
skewness and thinness. By checking whether these are 
rational, irrational or infinite we obtain a third classification.
Finally, we can classify valuations according to whether the
multiplicity is finite or infinite.

\smallskip
Thirdly, we can take advantage of the isomorphism of $\cV$ with the
universal dual graph, as worked out in Chapter~\ref{A3}. Specifically,
we can classify Krull valuations according to their associated
sequences of infinitely nearby points, following the 
terminology of Spivakovsky: see Definition~\ref{D403}.

\smallskip
Finally, as explained in Chapter~\ref{sec-puis}, any valuation
in $\cV$ extends to an element of the tree $\hcV$, \ie a valuation on 
the ring of formal power series in one variable with Puiseux series 
coefficients. Even though this extension is not unique, we 
can classify valuations in the valuative tree using the tree
structure on $\hcV$. This classification can then be rephrased
in the terminology of Berkovich.

\smallskip
We list all of these classifications in three tables.
Table~\ref{table8} contains the classification in terms of
the numerical invariants (here we use the 
completeness of $R$---see Appendix~\ref{sec-notcomplete})
and the associated sequences of infinitely nearby points.

\smallskip
In Table~\ref{table2} we instead focus on the tree structure of
the valuative tree, specifically the nonmetric tree structure,
skewness, thinness and multiplicity.
For convenience we define the skewness (thinness) of a tangent
vector in $\cV$ to be the skewness (thinness) of the associated
divisorial valuation.

\smallskip
Finally, in Table~\ref{table6} we present the classification
using the Puiseux series approach. Notice that this
table is a reproduction of Table~\ref{table3} in
Chapter~\ref{sec-puis}, and that we used the relative 
valuative tree $\cV_x$ instead of $\cV$. 
Recall that the terminology for valuations in $\cV$
and in $\cV_x$ coincides except for $\div_x$ which is
quasimonomial in $\cV_x$, and corresponds to the curve 
valuation $\nu_x$ in $\cV$.

\smallskip
A fifth classification is in terms of SKP's. We do not present 
this classification here as it depends on the choice
of local coordinates. Instead we refer to Definition~\ref{D101}.

\begin{table}[ht]
  \begin{center}
    \begin{tabular}{|c|l|c|c|c|c|}
      \hline
      \multicolumn{2}{|c|}{Classification} 
      & rk & rat.rk & tr.deg & Inf nearby pts\\
      \hline\hline
      & Divisorial & & 1 & 1 & Type 0\\
      \cline{2-2}\cline{4-6}
      \raisebox{1.5ex}[0cm][0cm]{Quasim.} & 
      Irrational & 
      \raisebox{1.5ex}[0cm][0cm]{1} &
      2 & 0 & 
      Type 2\\
      \hline
      & Nonexcept. & & & & Type 4\\
      \cline{2-2}\cline{6-6}
      \raisebox{1.5ex}[0cm][0cm]{Curve} & 
      Exceptional & 
      \raisebox{1.5ex}[0cm][0cm]{2} &
      \raisebox{1.5ex}[0cm][0cm]{2} &
      \raisebox{1.5ex}[0cm][0cm]{0} &
      Type 3\\
      \hline
      \multicolumn{2}{|c|}{Infinitely singular} 
      & 1 & 1 & 0 & Type 1\\
      \hline
    \end{tabular}
  \end{center}
  \medskip
  \caption{The table shows the classification of centered Krull 
    valuations on $R$ into five groups,
    and lists the three numerical invariants as well as the type
    of the associated sequence of infinitely nearby points as
    in Definition~\ref{D403}.}
  \label{table8}
\end{table}

\begin{table}[ht]
  \begin{center}
    \begin{tabular}{|c|l|c|c|c|c|}
      \hline
      \multicolumn{2}{|c|}{Classification} 
      & skewness & thinness & mult & Tree terminology\\
      \hline\hline
      & Divisorial & rational & rational & & Branch point\\
      \cline{2-4}\cline{6-6}
      \raisebox{1.5ex}[0cm][0cm]{Quasim.} & 
      Irrational  & irrational & irrational &
      \raisebox{1.5ex}[0cm][0cm]{$<\infty$} & 
      Regular point\\
      \hline
      & Nonexcept. & $\infty$ & $\infty$ & & 
      End\\
      \cline{2-4}\cline{6-6}
      \raisebox{1.5ex}[0cm][0cm]{Curve} & 
      Exceptional & rational & rational &
      \raisebox{1.5ex}[0cm][0cm]{$<\infty$} &
      Tangent vector\\
      \hline
      \multicolumn{2}{|c|}{Infinitely singular} 
      & $\le\infty$ & $\le\infty$ & $\infty$ & End \\
      \hline
    \end{tabular}
  \end{center}
  \medskip
  \caption{Elements of the valuative tree. Here we show the 
    classification of centered Krull valuations on $R$ in terms of
    the tree structure on $\cV$.}
  \label{table2}
\end{table}

\begin{table}[ht]
  \begin{center}
    \begin{tabular}{|c|c|c|c|c|c|}
      \hline
      \multicolumn{3}{|c|}{Valuations in $\cV_x$} & 
      \multicolumn{2}{|c|}{Valuations in $\hcV$} & 
      Berkovich\\
      \hline\hline
      & 
      \multicolumn{2}{|c|}{Divisorial} &
      & 
      Rational & 
      Type 2\\
      \cline{2-3}\cline{5-6}
      \raisebox{1.5ex}[0cm][0cm]{Quasim.} &
      \multicolumn{2}{|c|}{Irrational} &
      \raisebox{1.5ex}[0cm][0cm]{Finite type} & 
      Irrational & 
      Type 3\\ 
      \hline
      & 
      \multicolumn{2}{|c|}{Curve} &
      & 
      $m<\infty$ &
      \\
      \cline{2-3}\cline{5-5}
      \raisebox{2.5ex}[0cm][0cm]{Not quasi-}
      & 
      & 
      $A=\infty$ &
      \raisebox{1.5ex}[0cm][0cm]{Point type} & 
      $m=\infty$ & 
      \raisebox{1.5ex}[0cm][0cm]{Type 1}\\ 

      \cline{3-6}
      \raisebox{2.5ex}[0cm][0cm]{monomial} & 
      \raisebox{1.5ex}[0cm][0cm]{Inf sing} & 
      $A<\infty$ &
      \multicolumn{2}{|c|}{Special type} & 
      Type 4\\ 
      \hline
    \end{tabular}
  \end{center}
  \medskip
  \caption{Here the classification of centered valuations on $R$
    is in terms of their preimages in the tree $\hcV$ and
    in terms of Berkovich's terminology.}
  \label{table6}
\end{table}
%
%
%
%
\setcounter{theorem}{0}
\setcounter{equation}{0}
\setcounter{table}{0}
\setcounter{figure}{0}
\section{Combinatorics of plane curve singularities}\label{S501}
The study of plane curve singularities is a rich and 
well developed subject with its own notation and terminology.
Here we show how to interpret some classical invariants 
in terms of the valuative tree, specifically in terms of 
skewness, thinness and multiplicity.

We divide the study into two parts. In the first, we consider an
irreducible formal curve with its associated invariants introduced
by Zariski. In the second part we extend this study to a 
(possibly reducible but reduced) formal curve and its
associated Eggers tree.

In both cases, we shall see that all the classical invariants can
be understood in terms of the subtree of the valuative tree whose
end points are the curve valuations associated to the irreducible
components of the curve.
This serves to illustrate the idea that the valuative tree provides
an efficient way of encoding singularities (in this case of plane 
curves).\footnote{We do not claim, however, that the valuative tree 
leads to any results on plane curve singularities that could not be 
proved by other means.}
%
%
\subsection{Zariski's terminology for plane curve singularities}
\label{sec-termino}
Let us recall some notation from~\cite{Z3}. 
As before, $R$ denotes the
ring of formal power series in two complex variables.  
Consider an irreducible formal curve 
$C$.\footnote{Zariski actually writes $X$ instead of $C$.}

Zariski writes $n$ for the multiplicity $m(C)$
and gives two equivalent sets of invariants for $C$.
The first set is defined through Puiseux expansions. 
Pick local coordinates $(x,y)$ such that the curve $C$ is 
not tangent to $\{x=0\}$. In these coordinates, 
there is a Puiseux parameterization
$t\to(t^n,\sum a_jt^j)$ of $C$. 
Here $a_j=0$ for $0\le j<n$. Define
\begin{itemize}
\item
  $\b_1=\min\{j\ ;\ a_j\neq0,\ j\not\equiv 0\mod n\}$;
\item 
  $e_1=\gcd\{n,\b_1\}$;
\item 
  $\b_2=\min\{j\ ;\ a_j\ne0,\ j\not\equiv0\mod e_1\}$ (if $e_1>1$);
\item $e_2=\gcd\{n,\b_1,\b_2\}$ etc\dots
\end{itemize}
This process produces two finite sets of integers $(\b_1,\dots,\b_g)$
and $(e_1,\dots,e_g)$, but notice that the $e_i$'s can be recovered
from the $\b_i$'s and $n$.  Set $e_0=\b_0=n$ and $n_i=e_{i-1}/e_i$ for
$i=1,\dots,g$. The set $(\b_1/n,\dots,\b_g/n)$ is usually called the
set of \emph{generic characteristic exponents}\index{characteristic
exponents} of $C$.

The second set of invariants of $C$ is defined through
its value semigroup. This is by definition the 
collection of all intersection products $C\cdot D$ when $D$ 
ranges over all formal curves. 
We let $\bbar_0,\dots,\bbar_g \in\N^*$ be a minimal
set of generators of the semi-group of $C$, that is,
$\bbar_0=\min_D\{C\cdot D\}=n$, and, inductively, 
$\bbar_{i+1}=\min\{a=C\cdot D\ ;\ a\not\in\sum_0^i\N\bbar_l\}$ 
for $i\ge0$. 

The $\b_i$'s (with $n$) determine the $\bbar_i$'s uniquely 
(and vice-versa) as follows:
$\bbar_0=n$, $\bbar_1=\b_1$ and 
\begin{equation}\label{erel}
  \bbar_i=n_{i-1}\bbar_{i-1}+\b_i-\b_{i-1}
  \quad\text{for $i=2,\dots,g$}.
\end{equation}

It is known that $(n;\b_1,\dots,\b_g)$ determines 
the equisingularity class of $C$, 
\ie the topological type of the embedding of $C$ in $\C^2$.

\medskip
We now translate these invariants into the tree language that 
we have developed, specifically using skewness, thinness and 
multiplicity on the valuative tree.
For this we consider the approximating sequence of $\nu_C$ as 
defined in Section~\ref{sec-approx}:
\begin{equation}\label{e431}
  \nu_\fm=\nu_0<\nu_1<\dots<\nu_g<\nu_{g+1}=\nu_C.
\end{equation}
Here $\nu_i$, $1\le i\le g$, are divisorial valuations with
strictly increasing multiplicities $m_i=m(\nu_i)$, such that
the multiplicity is constant, equal to $m_i$m on each segment
$]\nu_{i-1},\nu_i]$. See Figure~\ref{F8}. Thus the generic 
multiplicity of $\nu_i$ is $b(\nu_i)=m_{i+1}$, where 
$m_{g+1}=m(C)$.

Write $\a_i$ and $A_i$ for the skewness and thinness of $\nu_i$,
respectively.
Then we claim that 
we have the identification
given by Table~\ref{table5}.
\begin{table}[ht]
  \begin{center}
    \begin{tabular}{|c|c|}
      \hline
      Classical  & Tree \\
      invariants & language \\
      \hline 
      $n$         & $m(C)$\\
      $g$         & $g$\\
      $\b_i/n$    & $A_i-1$ \\
      $n/e_i$     & $m_{i+1}=b_i$ \\
      $n_i$       & $b_i/m_i$ \\
      $\bbar_i/n$ & $\a_im_i$\\
      \hline     
    \end{tabular}
  \end{center}
  \medskip
  \caption{Dictionary expressing classical invariants of a curve 
    singularity in terms of the multiplicity, skewness and
    thinness of the elements of the approximating sequence of the
    associated curve valuation.}
  \label{table5}
\end{table}

An irreducible formal curve naturally defines an end
in the valuative tree. This provides an embedding
of a fundamentally discrete object (the ultrametric space $\cC$ of
local irreducible curves) into a ``continuous'' object
(the valuative tree $\cV$). 
The dictionary in Table~\ref{table5} shows that 
the \emph{discrete} invariants for curves can be viewed as
special cases of \emph{continuous} invariants for valuations.

\smallskip
Let us now prove the validity of the 
correspondence given in Table~\ref{table5}.
As the invariants $\b_i$ and $e_i$ are defined in terms
of Puiseux series, it is natural to use the analysis in
Chapter~\ref{sec-puis} to verify the first four entries in the 
table. Here we will freely apply the results from that chapter.

Thus define a Puiseux series by $\hphi=\sum a_jt^{j/n}\in\hk$ and let
$\hnu_\hphi=\val[\hphi;\infty]$ be the associated valuation in $\hcV$
of point type. Let $\phi$ be the minimal polynomial of $\hphi$ over $k$.
Then $C=\{\phi=0\}$ and  
$\hnu_\hphi$ is sent to the curve
valuation $\nu_C\in\cV_x$\footnote{By our choice of coordinates,
  $\nu_C$ is a valuation in both $\cV$ and $\cV_x$.}
by the restriction mapping $\Phi:\hcV\to\cV_x$ induced by the 
inclusion $\C[[x,y]]\subset\bk[[y]]$.

For $0\le i\le g$, let
\begin{equation*}
  \hnu_i=\val\left[\sum_0^{\b_i}a_jt^{j/n};\b_i/n\right]=
  \val\left[\sum_0^{\b_i-1}a_jt^{j/n};\b_i/n\right]
\end{equation*}
This defines an increasing sequence of rational valuations in $\hcV$.
The multiplicity of $\hnu_i$ is given by
\begin{equation}\label{e430}
  m(\hnu_i)
  =\lcm\left\{\frac{j}{n}\ ;\ 0<j<\b_i, a_j\ne 0\right\}
  =n/\gcd\{\b_1,\dots,\b_{i-1}\}
  =n/e_{i-1}.
\end{equation}
and the multiplicity is constant in the segment 
$]\hnu_{i-1},\hnu_i]$. By Theorem~\ref{T602} this implies
that if $\nu_i=\Phi(\hnu_i)$ is the restriction of $\hnu_i$
to $R$, then $(\nu_i)_1^g$ is the approximating sequence 
of $\nu_C$.
The first and the second entries in Table~\ref{table5} are thus clear.
The fourth entry follows from~\eqref{e430} and implies the fifth.
Further, as $\hnu_i$ has Puiseux parameter $\hbeta_i/n$, 
$\nu_i$ must have thinness $\hbeta_i/n-1$, which gives the third entry
in the table.

Finally let us prove that $\bbar_i/n=m_i\a_i$ for $1\le i\le g$.
This can be seen from~\eqref{erel} and the relation between
thinness and skewness.
Indeed, we first have $\bbar_1/n=\b_1/n=A_1-1=\a_1m_1$ as $m_1=1$. 
Suppose we have proved that $\bbar_{i-1}/n=\a_{i-1}\,m_{i-1}$. 
Then
\begin{multline*}
  \bbar_{i}
  =n_{i-1}\bbar_{i-1}+\b_i-\b_{i-1}
  =nn_{i-1}\a_{i-1}m_{i-1}+n(A_i-A_{i-1})=\\
  =nn_{i-1}m_{i-1}\a_{i-1}+n m_i(\a_i-\a_{i-1})
  =n\a_i\,m_i.
\end{multline*}
By induction, this completes the verification of Table~\ref{table5}.

An alternative way of showing that $\bbar_i/n=a_im_i$ is to use the
fact that the value group of $\nu_C$ is given by
$\nu_C(R)=\sum_{i=0}^gm_i\a_i\N$
(see Proposition~\ref{P556}) and the formula 
$m(C)\nu_C(\psi)=C\cdot\{\psi=0\}$ for any $\psi\in\fm$.
%
%
\subsection{The Eggers tree}
Let $C$ be a reduced formal curve and let $C_1,\dots,C_r$ be its
branches (\ie irreducible components).  The embedding of $C$ in $\C^2$
is determined up to topological conjugacy\footnote{\ie the
equisingularity type of $C$} by the generic characteristic exponents
of each branch, and by the order of contact between them.

This numerical information can be encoded geometrically in
the \emph{Eggers tree}\index{tree!Eggers}\index{Eggers tree}
$\cT_C$, first introduced by Eggers~\cite{eggers}, and 
subsequently used in studies of plane 
curve singularities~\cite{GB,popescu}.
Here we show that $\cT_C$ has a natural interpretation inside the 
valuative tree $\cV$.

Let us recall the definition of $\cT_C$\index{Eggers tree},
essentially following~\cite{popescu}.
First assume $C$ is irreducible, \ie $r=1$, and let
$\b_1/n,\dots,\b_g/n$ be its generic characteristic exponents
(with the notation in Appendix~\ref{sec-termino}).
Also define $\b_0\=n$ and $\b_{g+1}:=+\infty$.
Then $\cT_C$ is a simplicial tree whose 
vertices are exactly of the form $\{\b_i/n\}_{0\le i\le g+1}$,
and with exactly $g+1$ edges linking $\b_i/n$ to
$\b_{i+1}/n$ for $0\le i\le g$. 
By sending each vertex to the corresponding element in $[1,+\infty]$, 
we may also view $\cT_C$ as the interval $[1,+\infty]$ 
with $g+2$ marked points $\{\b_i/n\}_{0\le i\le g+1}$. 
Alternatively speaking, $\cT_C$ is a rooted, nonmetric tree
with a parameterization onto $[1,+\infty]$, still with
$g+2$ marked points.

When $C=\bigcup C_j$ is reducible, $\cT_C$ is constructed 
by patching together the trees $\cT_{C_j}$ as follows.
First define the coincidence order $K_{jj'}=K(C_j,C_{j'})$ 
between two branches $C_j$ and $C_{j'}$ as follows.
Fix local coordinates $(x,y)$ such that all branches
are transverse to $\{x=0\}$, and let 
$\hphi=\sum a_lx^{\hbeta_l}$, 
$\hphi'=\sum a'_lx^{\hbeta'_l}$ 
be Puiseux series associated to $C_j$ and $C_{j'}$, respectively. 
Let $(\hphi_l)_1^{n_j}$, $(\hphi'_{l'})_1^{n_{j'}}$ 
denote their orbits under the
action of the Galois group $\Gal(\hk/k)$ where 
$k$ and $\hk$ are the fields of Laurent and Puiseux series 
in $x$, respectively (see Section~\ref{Sgalois}). 
Then
\begin{equation*}
  K_{jj'}:=K(C_j,C_{j'})\=\max_{l,l'}\nu_\star(\hphi_l-\hphi'_{l'}),
\end{equation*}
where $\nu_\star$ is the valuation on $\hk$ given by
$\nu_\star(\sum a_jx^{\hbeta_j})\=\min\hbeta_j$.
Note that $K_{jj'}\ge1$.

The trees $\cT_{C_j}\simeq[1,+\infty]$ and 
$\cT_{C_{j'}}\simeq[1,+\infty]$ are now
patched together at the point corresponding to $K_{jj'}$.
In other words, $\cT_C$ is the set of pairs $(C_j,t)$ 
where $1\le j\le r$ and $t\in [1,+\infty]$ 
modulo the relation $(C_j,s)=(C_{j'},t)$ iff 
$s=t\le K_{jj'}$.
Notice the similarity to the construction in 
Section~\ref{def-ultra}.
Thus $\cT_C$ is naturally a rooted, nonmetric tree
with a parameterization onto $[1,+\infty]$.
The root of $\cT_C$ is $(C_j,1)$ and 
the maximal points of $\cT_C$ are exactly of the
form $(C_j,+\infty)$.

To be precise, $\cT_C$ also comes with a finite
number of marked points. 
These are exactly the points coming from marked points
on the trees $T_{C_j}$
(\ie the points of the form $(C_j,\beta_{ji}/n_j)$)
together with the branch points of $\cT_C$
(\ie the points of the form 
$(C_j,K_{jj'})=(C_{j'},K_{jj'})$ for $j\ne j'$).
This marking allows us to recover $\cT_C$ as a
simplicial tree: the vertices are the marked points
and the edges are open segments in $\cT_C$
containing no marked points.
\begin{proposition}\label{P308}
  Let $C$ be a reduced formal curve with irreducible components
  $C_1,\dots C_r$. 
  Let $n_j$ be the multiplicity of $C_j$, 
  and let $\{\b_{ji}/n_j\}_{i=0}^{g_j+1}$ be the set of
  (generic) characteristic exponents of $C_j$ 
  (by convention $1$ and $+\infty$ belong to this set).
  Denote by $\nu_j=\nu_{C_j}$ the curve valuation associated
  to $C_j$ and by $(\nu_{ji})_{i=1}^{g_j}$
  the approximating sequence of $\nu_j$.
  
  Define a map $\Psi:\cT_C\to\cV$ by 
  sending $(C_j,t)$ to the unique valuation in the segment
  $[\nu_\fm,\nu_j]$ with thinness equal to $t+1$. 
  Then $\Psi$ gives an isomorphism of parameterized trees
  from $\cT_C$ onto the subtree 
  $\cV_C=\cup_{1\le j\le r}\{\nu\le\nu_j\}$ 
  of the valuative tree parameterized by thinness (minus one).

  Moreover $\Psi(C_j,\b_{ji}/n_j)=\nu_{ji}$ 
  and $\Psi(C_j,K_{jj'})=\nu_j\wedge\nu_{j'}$.
  In other words, the marked points on $\cT_C$ are sent to
  the points on $\cV_C$ consisting of the root, the ends,
  the branch points, and all regular points where the 
  multiplicity function is not locally constant.
\end{proposition}
\begin{proof}
  For an irreducible curve $C$, we have seen in 
  the previous section that the generic characteristic 
  exponents corresponds exactly to 
  the values $-1+A(\nu_i)$ where $(\nu_i)_1^g$ is the
  approximating sequence of the curve valuation $\nu_C$. 
  When $C_j$, $C_{j'}$ are two irreducible curves, 
  then by Theorem~\ref{T602}, we have
  $K_{jj'}=-1+A(\nu_j\wedge\nu_{j'})$. 
  These two facts immediately imply the proposition.
\end{proof}

In the original definition of the Eggers tree~\cite{eggers}, 
two types of edges were distinguished: dashed or plain; 
as well as two types of vertices: white or black. 
Using the isomorphism above, the different types of 
vertices and edges can be seen inside the tree 
$\cV_C\subset\cV$ as follows: the white vertices are exactly 
the curve valuations, and an edge $]\nu,\nu'[$ is 
dashed iff the multiplicity on the edge is (constant)
equal to that of its left end point $\nu$.

\smallskip
It is also possible to interpret the Eggers tree $\cT_C$ in terms of
the minimal desingularization of $C$.  This was done by Popescu-Pampu:
see~\cite[Th\'eor\`eme 4.4.1. p.152]{popescu}. 
Let us outline how to understand this interpretation 
using the analysis in Section~\ref{Sdualmini}.

Let $\pi:X\to(\C^2,0)$ be the minimal
desingularization\index{desingularization!and Eggers
tree}\index{Eggers tree!and minimal desingularization} of $C$, \ie
$\pi\in\fB$ is minimal with the property that the total transform
$\pi^{-1}(C)$ has normal crossings in $X$. 
Let $\Gast_C$ be the dual graph of $\pi$. 
This is by definition a finite subset of $\Gast\subset\Gamma$
and in particular an $\N$-tree.
Enlarge $\Gast_C$ by adding the elements $C_j\in\Gamma$
corresponding to the irreducible components of $C$.
The resulting set is still an $\N$-tree.
Now remove from this set all points except for the following:
the points $C_j$; the branch points; the root $E_0$.
We obtain a set $X_C$, which by
construction is a subset of the $\R$-tree
$\cS_C:=\bigcup_j[E_0,C_j]\subset\Gamma$.
Let us equip $\cS_C$ with the Farey parameterization, 
from which we subtract a constant 1, and the multiplicity
function induced from $\Gamma$.
Then $X_C$ is exactly the subset consisting of 
the root, the ends, the branch points, and all 
regular points where the multiplicity function 
(induced from $\Gamma$) is not locally constant.

Proposition~\ref{P308} and Theorem~\ref{thm-universal} now
imply that $\cS_C$ is isomorphic to the Eggers tree $\cT_C$.
More precisely, the fundamental isomorphism
$\Phi:\Gamma\to\cV$ restricts to
an isomorphism of parameterized trees $\Phi:\cS_C\to\cV_C$,
where $\cV_C$ is defined in Proposition~\ref{P308}.
Thus the composition $\Psi^{-1}\circ\Phi:\cS_C\to\cT_C$ is
an isomorphism of parameterized trees and maps $X_C$
onto the marked points on $\cT_C$.

We may also recover the minimal desingularization $\pi_C$ of $C$ from
the Eggers tree $\cT_C$ using the algorithm in
Section~\ref{Sdualmini}, as $\cT_C$ gives precisely the
equisingularity type of the curve.

%
%
%
%
\setcounter{theorem}{0}
\setcounter{equation}{0}
\setcounter{table}{0}
\setcounter{figure}{0}
\section{What are the essential assumptions on the ring $R$?}
\label{sec-notcomplete}
Throughout the monograph we have assumed that $R$ is the ring
of formal power series in two complex variables. 
An equivalent way of describing $R$ is as an
\emph{equicharacteristic, complete, two-dimensional, 
  regular local ring with residue field $\C$}.
Here we address the question of which of these conditions are actually
necessary for the analysis to go through.

First of all, nothing changes if we replace the residue field
$\C$ by any algebraically closed field $k$ of characteristic zero.
On the other hand, some assumptions are still crucial 
to our analysis like locality, regularity, and dimension two.
Without these assumptions $\cV$ needs not be a tree 
(see Remark~\ref{dim3} for instance). 

Let us briefly discuss the assumption of completeness, 
and the assumptions on the residue field $k$.
%
%
\subsection{Completeness}
In this text we have always assume $R$ to be a complete ring. Let us
explain why this assumption is not essential.

Recall that Cohen's Theorem~\cite[p.304]{ZS}
asserts that any complete, equicharacteristic, 
regular local ring is isomorphic to the ring of formal power series
with coefficients in the residue field of the ring.
\begin{proposition}
  Suppose $R$ is an equicharacteristic, two-dimensional, 
  regular local ring with algebraically closed residue field, 
  and let $\hat{R}$ be its completion. 
  Denote by $\cV(R)$ the set of normalized centered
  valuations on $R$ with values in $\Rbar$. 
  Then the natural restriction map 
  $\cV(\hat{R}) \to \cV(R)$ is a bijection.
\end{proposition}
\begin{proof}
  If $(x,y)$ is a regular system of parameters for $R$, and $k$ is 
  the residue field of $R$, then the assumptions imply that
  $k[x,y]\subset R\subset k[[x,y]]=\hat{R}$. 
  The assertion then follows immediately from Proposition~\ref{extend}.
\end{proof}

There is still one point where the assumption of completeness is important. 
It concerns the numerical invariants of a curve valuation.
\begin{proposition}
  Let $R$ be a (not necessarily complete) 
  equicharacteristic, two-dimensional, 
  regular local ring with algebraically closed residue field.
  Suppose $\nu\in\cV(R)$ is a curve valuation associated to the
  irreducible curve $V$. Then $\trdeg\nu=0$, and
  \begin{itemize}
  \item
    $\ratrk\nu=\rk\nu=2$ when $V=\{\phi=0\}$ for some $\phi\in R$;  
  \item
    $\ratrk\nu=\rk\nu=1$ otherwise.  
  \end{itemize}
\end{proposition}
In particular when $R$ is complete, we are always in the first case.
When $R$ is the ring of \emph{convergent} power series, the first case
appears exactly when $V$ is a formal curve, and  the second when $V$ is
an analytic curve.
\begin{proof}
  Embed $R$ in its completion $\hat{R}$. 
  Pick $\phi\in\hat{R}$ such that
  the Krull valuation $\nu$ is given by 
  $\nu(\psi)=(\nu_\phi(\tpsi),\div_\phi(\psi))\in\N\times\N$ where 
  $\psi=\phi^{\div_\phi(\tpsi)} \hat{\psi}$ with $\phi$ not dividing
  $\tpsi$, and $\nu_\phi(\tpsi)$ is equal to the intersection
  number between $\{\phi=0\}$ and $\{\tpsi=0\}$.
  
  The set $I\=\{\psi\in R\ ;\ \div_\phi(\psi)>0\}$ is a prime
  ideal in $R$. Its completion $I\cdot\hat{R}$ is included in the
  height one ideal $\{\psi\in\hat{R}\ ;\ \div_\phi(\psi)>0\}$, 
  hence $I$ has height zero or one. 
  As $R$ is a factorial domain~\cite[p.312]{ZS}, we have 
  either $I=(0)$, in which case $\rk\nu=\ratrk\nu=2$; 
  or $I$ is generated by an irreducible element $\phi'\in R$,
  in which case $\rk\nu=\ratrk\nu=1$.
\end{proof}
%
%
\subsection{The residue field}
We assumed the characteristic of the residue field of $k$ to be zero
for sake of convenience. But Weierstrass' preparation theorem holds on
$k[[x,y]]$ without any assumption on the characteristic of $k$
(see~\cite[p.139]{ZS}). Hence $k$ may be taken to be an arbitrary
algebraically closed field.

The fact that $k$ is algebraically closed is used in an essential way
in Corollary~\ref{primes}, hence in Theorem~\ref{approx} which
gives the description of valuations in terms of SKP's. 
Most results in the monograph rely on this description. 
It would be interesting to extend some of our results 
to the case when $k$ is not algebraically closed 
(for instance when $k=\R$).

Finally all methods presented here using SKP's or Puiseux series fail
on a non-equicharacteristic ring. On the other hand, the geometric
method in Chapter~\ref{A3} based on the universal dual graph
is likely to work. Again it would be interesting to extend our results 
to this more general setting.


%
%
%
%
%
%


\printindex
\end{document}